%% file: MAlgebra2.English.tex
\documentclass{amsbook}

\input{MAlgebra2}

\begin{document}
\title{Introduction into Noncommutative Algebra}
\def\subtitle{Volume 2}
\def\subtitleA{Module over Algebra}

\def\PartA{Associative Algebra}
\def\PartB{Non-Associative Algebra}

\ShowEq{contents}
\end{document}

%% file: MAlgebra2.tex
\def\PrintBook{}
\def\UseParacol{}
\def\InputEncoding{on}
\input{Commands}
\def\BookNumber{MAlgebra2}
\xRefDef{0909.0855}
\xRefDef{1003.1544}
\xRefDef{1104.5197}
\xRefDef{1908.04418}

\DefEq
{
\input{\FilePrefix Prolog.\TheLanguage}
\input{\FilePrefix Preliminary.Representation.\TheLanguage}
\input{\FilePrefix Preliminary.D.Module.\TheLanguage}

\ifx\texFuture\Defined
преобразование базиса и геометрический объект;

структурные константы произведения отображений

центр алгебры

Зачем дважды описывать унитальную алгебру для левого и правого модулей?
переместить определение в главу D алгебры.

множество гомоморфизмов модуля столбцов

проверить 8428-0408 in preliminary. i have extra text

стандартное представление отображений сопряжения

я хочу рассмотреть следующие темы:

линейное отображение модуля;

унитальное расширение алгебры Ли;

преобразование базиса и геометрический объект;

Тензор и базис (дневник поиска 2022 18 сентября)

система линейных уравнений;

делители нуля и связь с неоднозначной записью координат вектора;

полилинейное отображение и многочлен;

деление многочленов;

насколько неоднозначно определено частное от деления на многочлен;

алгебра дробей в некоммутативном случае;

алгебра с сопряжением;

двойной определитель в алгебре с сопряжением;

аффинная геометрия;

значение косой формы на группе векторов;

некоммутативная сумма

алгебры счётной размерности

В частности, я долго думал, должен ли я представлять объём
с помощью квазидетерминанта.
Был момент, когда я в этом усомнился.
Но я тут же рассмотрел относительно простой случай.
Я рассмотрел параллелограм, построенный на двух векторах.
В каком порядке я должен рассматривать произведение этих векторов,
если произведение некоммутативно.

С более широкой точки зрения, это вопрос о представлении
дифференциальных форм в модуле над некоммутативной алгеброй.
\fi

\part{\PartA}
\ShowEq{def AModule}
\ShowEq{def left}
\input{\FilePrefix A.Module.\TheLanguage}
\input{\FilePrefix Homomorphism.A.Module.\TheLanguage}
\ShowEq{def right}

\input{\FilePrefix A.Module.\TheLanguage}
\input{\FilePrefix Homomorphism.A.Module.\TheLanguage}
\part{\PartB}
\ShowEq{def left}
\input{\FilePrefix Na.A.Module.\TheLanguage}
\ShowEq{def right}
\input{\FilePrefix Na.A.Module.\TheLanguage}

\input{\FilePrefix Biblio.\TheLanguage}
\input{\FilePrefix Index.\TheLanguage}
\input{\FilePrefix Symbol.\TheLanguage}
}
{contents}

%% file: Commands.tex
\def\Defined{}
\def\ValueOff{off}%
\def\ValueOn{on}%
\scrollmode
\ifx\FilePrefix\undefined
\newcommand{\FilePrefix}{}
\fi
\ifx\UseRussian\Defined
\usepackage{cmap}
\usepackage[T2A,T2B]{fontenc}
\usepackage[cp1251]{inputenc}
\fi

\usepackage{trace}

\ifx\GJSFRA\Defined
\usepackage[top=1.905cm,bottom=1.905cm,inner=1.65cm,outer=1.65cm]{geometry}
\fi
\ifx\Presentation\Defined
\usepackage[margin=1cm]{geometry}
\fi

\ifx\setCACAA\Defined
\usepackage{times,mathptm}
\usepackage{amsmath,amsfonts,amssymb}
\usepackage{graphicx,graphics,epsfig}
\usepackage{fancyhdr,fancybox}
\usepackage{verbatim}
\usepackage{setspace}
\fi

\ifx\CreateSpace\Defined
\usepackage[top=1.905cm,bottom=1.905cm,inner=1.905cm,outer=1.27cm]{geometry}
\def\Publisher{CreateSpace Independent Publishing Platform}
\fi

\ifx\KDP\Defined
\usepackage[top=\MarginTop,bottom=\MarginBottom,inner=\Gutter,outer=\MarginOut]{geometry}
\def\Publisher{Kindle Direct Publishing}
\fi
\ifx\PublishBook\Defined
\def\PrintPaper{}
\usepackage{setspace}
\ifx\UseRussian\undefined
\usepackage{pslatex}
\fi
\usepackage[margin=2cm]{geometry}
\fi
\usepackage{graphicx}
\usepackage[all]{xy}
\usepackage{xcolor}
\ifx\PrintPaper\undefined
\definecolor{CoverColor}{rgb}{.82,.7,.55}
\definecolor{Gray}{rgb}{.94,.95,.94}
\definecolor{UrlColor}{rgb}{.9,0,.3}
\definecolor{SymbColor}{rgb}{.4,0,.9}
\definecolor{IndexColor}{rgb}{1,.3,.6}
\definecolor{BlueColor}{rgb}{.1,.6,.8}
\definecolor{Eq}{rgb}{0.8, 0, 1}
\newcommand\BlueText[1]{\textcolor{BlueColor}{#1}}
\newcommand\EqText[1]{\textcolor{Eq}{#1}}
\newcommand\RedText[1]{\textcolor{red}{#1}}
\else
\definecolor{UrlColor}{rgb}{.1,.1,.1}
\definecolor{SymbColor}{rgb}{.1,.1,.1}
\definecolor{IndexColor}{rgb}{.5,.1,.1}
\newcommand\BlueText[1]{#1}
\newcommand\RedText[1]{#1}
\fi
\usepackage{framed}

\ifx\UseParacol\Defined
\usepackage{paracol}
\globalcounter*
\globalcounter{theorem}
\globalcounter{Index}
\globalcounter{Symbols}
\globalcounter{Symbol}
\globalcounter{enumi}
\footnotelayout{m}
\usepackage{etoolbox}
\makeatletter
\patchcmd{\pcol@zparacol}{\topskip}{0pt}{}{}
\makeatother
\makeatletter
\def\@makefntext{\noindent\@makefnmark}
\makeatother
\fi

\ifx\UseRussian\Defined
\usepackage[english,russian]{babel}
\else
\ifx\InputEncoding\ValueOn
\UseRawInputEncoding
\fi
\fi
\usepackage{xr-hyper}
\ifx\Link\ValueOff
\usepackage[unicode,draft]{hyperref}
\else
\usepackage[unicode]{hyperref}
\fi
\hypersetup{pdfdisplaydoctitle=true}
\hypersetup{colorlinks}
\hypersetup{citecolor=UrlColor}
\hypersetup{urlcolor=UrlColor}
\hypersetup{linkcolor=UrlColor}
\hypersetup{pdffitwindow=true}
\hypersetup{pdfnewwindow=true}
\hypersetup{pdfstartview={FitH}}
\GetTitleStringSetup{expand}

\input{Statement}

\newcommand\Prolog
{
\ifx\PrimaryMSC\undefined
\else
\ifx\SecondaryMSC\undefined
\subjclass[2010]{Primary \PrimaryMSC}
\else
\subjclass[2010]{Primary \PrimaryMSC;
Secondary \SecondaryMSC}
\fi
\fi
\input{Prolog.Parm}\ifx\MoveAbstract\undefined
\input{\FilePrefix Abstract.\BookNumber.\TheLanguage}
\fi
\maketitle
\input{Prolog.Eq}
\ifx\MoveAbstract\undefined
\else
\input{\FilePrefix Abstract.\BookNumber.\TheLanguage}
\fi
\tableofcontents
\input{\FilePrefix Preface.\BookNumber.\TheLanguage}
\ePrints{0767-8264,1305.4547}%
\ifx\Semafor\ValueOn%
\input{\FilePrefix Convention.\TheLanguage}
\fi
}
\newcommand\Epilog[1][0]
{
\input{\FilePrefix Biblio.\TheLanguage}
\input{\FilePrefix Index.\TheLanguage}
\input{\FilePrefix Symbol.\TheLanguage}
\def\TempA{PrintCover}%
\def\TempB{#1}%
\ifx\TempA\TempB%
\newpage
\pagestyle{empty}
\includegraphics[scale=.75]{Cover_Front_\BookNumber_\CurrentLanguage}
\newpage
\includegraphics[scale=.75]{Cover_Back_\BookNumber_\CurrentLanguage}
\fi
}
\newcounter{Index}
\newcounter{Symbol}
\newcounter{Symbols}

\def\hyph{\penalty0\hskip0pt\relax-\penalty0\hskip0pt\relax}
\def\Hyph{-\penalty0\hskip0pt\relax}%
\newcommand{\Basis}[1]{\overline{\overline{#1}}{}}
\newcommand{\Vector}[1]{\overline{#1}{}}
\ifx\PrintPaper\undefined
\newcommand{\gi}[1]{\boldsymbol{\textcolor{IndexColor}{\it #1}}}
\else
\newcommand{\gi}[1]{\boldsymbol{\it #1}}
\fi
\newcommand\gii{\gi i}
\newcommand\giI{\gi I}
\newcommand\gij{\gi j}
\newcommand\giJ{\gi J}
\newcommand\gik{\gi k}
\newcommand\gil{\gi l}
\newcommand\gin{\gi n}
\newcommand\gim{\gi m}
\newcommand\giA{\gi 1}

\def\Items#1{\ItemList#1,LastItem,}%
\def\LastItem{LastItem}%
\def\ItemList#1,{\def\ViewBook{#1}%
\ifx\ViewBook\LastItem%
\else%
\ifx\ViewBook\BookNumber%
\def\Semafor{on}%
\fi%
\expandafter\ItemList%
\fi%
}%

\newcommand{\ePrints}[1]%
{%
\def\Semafor{off}%
\Items{#1}%
}%

\makeatletter
\newcommand{\NameDef}[1]{%
\expandafter\gdef\csname #1\endcsname%
}%
\newcommand{\xNameDef}[1]{%
\expandafter\xdef\csname #1\endcsname%
}%
\newcommand{\ShowSymbol}[2]{%
\@nameuse{ViewSymbol#1,,,#2}%
}%

\newcommand{\symb}[4][]{%
\def\TempA{}%
\def\TempB{#1}%
\ifx\TempA\TempB%
\def\ThisSymbol{#3}%
\else%
\edef\ThisSymbol{#3(#1)}%
\fi%
\@ifundefined{ViewSymbol\ThisSymbol}{%
\addtocounter{Symbols}{1}%
\edef\SymbolId{\arabic{Symbols}}%
\xNameDef{ViewSymbol\ThisSymbol}{\SymbolId}%
\NameDef{ViewSymbol\ThisSymbol:::\SymbolId}{#2}%
\@namedef{RefSymbol}{:}%
}{%
\edef\Symbols{\@nameuse{ViewSymbol\ThisSymbol}}%
\def\aSymbolId{0}%
\@for\Symbol:=\Symbols\do{%
\protected@edef\TempA{#2}%
\protected@edef\TempB{\@nameuse{ViewSymbol\ThisSymbol:::\Symbol}}%
\ifx\TempA\TempB%
\edef\aSymbolId{\Symbol}%
\fi%
}%
\def\Zero{0}%
\ifx\aSymbolId\Zero%
\addtocounter{Symbols}{1}%
\edef\SymbolIds{\@nameuse{ViewSymbol\ThisSymbol},\arabic{Symbols}}%
\xNameDef{ViewSymbol\ThisSymbol}{\SymbolIds}%
\edef\SymbolId{\arabic{Symbols}}%
\NameDef{ViewSymbol\ThisSymbol:::\SymbolId}{#2}%
\else%
\def\SymbolId{\aSymbolId}%
\fi%
\addtocounter{Symbol}{1}%
\@namedef{RefSymbol}{\arabic{Symbol}}%
}%
\@namedef{LabelSymbol}{\label{symbol: \ThisSymbol:\@nameuse{RefSymbol}}}%
\edef\RefIds{RefSymbol\ThisSymbol===\SymbolId}%
\@ifundefined{\RefIds}{%
\xNameDef{\RefIds}{\@nameuse{RefSymbol}}%
}{%
\xNameDef{\RefIds}{\@nameuse{\RefIds},\@nameuse{RefSymbol}}%
}%
\NameDef{ViewSymbol#3,,,#4}{\textcolor{SymbColor}{#2}}%
\def\Temp{#4}%
\def\One{1}%
\def\Two{2}%
\def\Three{3}%
\ifx\Temp\One%
$\@nameuse{ViewSymbol#3,,,#4}$%
\fi%
\ifx\Temp\Two%
\[\@nameuse{ViewSymbol#3,,,#4}\]%
\fi%
\ifx\Temp\Three%
\@nameuse{ViewSymbol#3,,,#4}%
\fi%
\@nameuse{LabelSymbol}%
}%

\newcommand\AddEq[3][0]%
{%
\@ifundefined{ViewEq #2}{%
\expandafter\newcommand\csname ViewEq #2\endcsname[#1]{#3}%
}{%
\errmessage {second entry of AddEq: #2}%
}%
}%

\newcommand\DefEq[2]{%
\@ifundefined{ViewEq #2}{%
\NameDef{ViewEq #2}{#1}%
}{%
\errmessage {second entry of DefEq: #2}%
}%
}%
\newcommand{\DefEquation}[2]{%
\AddEq{#2}%
{%
\begin{equation}%
#1%
\EqLabel{#2}%
\end{equation}%
}%
}%
\newcommand{\AddEquation}[2]{%
\AddEq{#1}{\begin{equation}#2\EqLabel{#1}\end{equation}}%
}%
\def\ViewParm#1{\protect\getParm#1,endParm,}%
\def\endParm{endParm}%
\def\getParm#1,{\def\temp{#1}%
\ifx\temp\endParm%
\else%
\ShowEq{#1}%
\expandafter\getParm%
\fi%
}%

\newcommand{\EqParm}[2]{%
\ViewParm{#2}\ShowEq{#1}%
}%
\newcommand{\EquationParm}[2]{%
\@ifundefined{ViewEq #1[#2]}%
{%
\ViewParm{#2}%
\DefEquation{\ShowEq{#1}}{#1[#2]}%
}{}%
\ShowEq{#1[#2]}%
}%
\newcommand\DrawEqParm[3]{%
\ViewParm{#2}%
\@ifundefined{ViewEq #1(#2)}{%
\DefEq%
{%
\ShowEq{#1}%
}{#1(#2)}%
}{%
}%
\DrawEq{#1(#2)}{#3}%
}%
\newcommand\EqRef[2][]%
{%
\def\Semafor{on}%
\def\Temp{}%
\edef\Tempa{#1}%
\ifx\Tempa\Temp%
\def\Semafor{off}%
\fi%
\@ifundefined{xRefDef#1}{%
\def\Semafor{off}%
}{}%
\ifx\Semafor\ValueOff%
\eqref{eq: #2}%
\else%
\citeBib{#1}-\eqref{#1-\TheLanguage-eq: #2}%
\fi%
}%
\newcommand\eqRef[3][]{\EqRef[#1]{#2(#3)}}%
\newcommand\EqLabel[1]{\label{eq: #1}}%
\newcommand\ShowEq[1]{%
\@ifundefined{ViewEq #1}{%
\message {error: missed ShowEq #1}%
\newline%
\RedText{missed ShowEq #1}%
\newline%
}{%
\csname ViewEq #1\endcsname
}%
}%

\newcommand\DrawEq[3][]{%
\def\Temp{}%
\def\Tempa{#3}%
\ifx\Tempa\Temp%
\[\ShowEq{#2}#1\]%
\else%
\def\Temp{-}%
\ifx\Tempa\Temp%
$\ShowEq{#2}#1$
\else%
\@ifundefined{ViewEq #2(#3)}{%
\AddEquation{#2(#3)}{\ShowEq{#2}#1}%
  }{%
\errmessage {second entry of DrawEq: #2(#3)}%
}%
\ShowEq{#2(#3)}
\fi%
\fi%
}%
\makeatother

\DeclareMathOperator{\Hom}{\mathrm{Hom}} 
\DeclareMathOperator{\End}{\mathrm{End}} 
\DeclareMathOperator{\rank}{\mathrm{rank}} 
\DeclareMathOperator{\id}{\mathrm{id}} 
 
\newcommand{\subs}{${}_*$\Hyph}
\newcommand{\sups}{${}^*$\Hyph}

\newcommand{\CRstar}{{}^*{}_*}
\newcommand{\RCstar}{{}_*{}^*}
\newcommand{\CRcirc}{{}^{\circ}{}_{\circ}}
\newcommand\RCcirc{{}_{\circ}{}^{\circ}}

\newcommand{\RC}{$\RCstar$\Hyph}
\newcommand{\CR}{$\CRstar$\Hyph}
\newcommand{\RCo}{$\RCcirc$\Hyph}
\newcommand{\CRo}{$\CRcirc$\Hyph}
\newcommand{\drc}{$D\RCstar$\Hyph}
\newcommand{\Drc}{$\mathcal D\RCstar$\Hyph}
\newcommand{\dcr}{$D\CRstar$\hyph}
\newcommand{\rcd}{$\RCstar D$\Hyph}
\newcommand{\crd}{$\CRstar D$\Hyph}

\newcommand{\RCPower}[1]{#1\RCstar}
\newcommand{\CRPower}[1]{#1\CRstar}
\newcommand{\RCInverse}{\RCPower{-1}}
\newcommand{\CRInverse}{\CRPower{-1}}
\newcommand{\RCRank}{\rank(\RCstar)}
\newcommand{\CRRank}{\rank(\CRstar)}
\newcommand\RCDet{\det(\RCstar)}
\newcommand\RCdet[2]{\aUD{\RCDet}{#1}{#2}\,}
\newcommand\CRDet{\det(\CRstar)}
\newcommand\CRdet[2]{\aUD{\CRDet}{#1}{#2}\,}
\newcommand{\RCGL}[2]{\ensuremath{GL(\gi #1\RCstar #2)}}
\newcommand{\CRGL}[2]{\ensuremath{GL(\gi #1\CRstar #2)}}
\newcommand\sT[1]{$*#1$\Hyph}%
\newcommand\Ts[1]{$#1*$\Hyph}%
\newcommand\sD{$\star D$\Hyph}%
\newcommand\Ds{$D\star$\Hyph}%

\newcommand\VirtVar{\vphantom{\overset{\rightarrow}{1}^1}}
\newcommand\pC[2]{{}_{#1\cdot #2}}%
\newcommand\DcrPartial[1]%
{%
\def\tempa{}%
\def\tempb{#1}%
\ifx\tempa\tempc%
(\partial\CRstar)%
\else%
(\partial_{\gi{#1}}\CRstar)%
\fi%
}%
\newcommand\rcDPartial[1]%
{%
\def\tempa{}%
\def\tempb{#1}%
\ifx\tempa\tempc%
(\RCstar\partial)%
\else%
(\RCstar\partial_{\gi{#1}})%
\fi%
}%
\newcommand\StandPartial[3]%
{%
\frac{d^{\gi{#3}} #1}{d #2}%
}%

\ifx\setCACAA\undefined
\renewcommand{\uppercasenonmath}[1]{}

\makeatletter
\def\l@chapter{\@tocline{0}{8pt plus1pt}{1pc}{}{}}

\def\part{\newpage\thispagestyle{empty}%
  \null\vfil  \markboth{}{}\secdef\@part\@spart}
\def\@part[#1]#2{%
  \ifnum \c@secnumdepth >-2\relax \refstepcounter{part}%
    \addcontentsline{toc}{part}{\partname\ \thepart.\ #1}%
  \else
    \addcontentsline{toc}{part}{#1}\fi
  \begingroup\centering
  \ifnum \c@secnumdepth >-2\relax
       {\fontsize{\@xviipt}{22}\bfseries
         \partname\ \thepart} \vskip 20\p@ \fi
  \fontsize{\@xxpt}{25}\bfseries
      #1\vfil\vfil\endgroup \newpage\thispagestyle{empty}}

\newcommand\@dotsep{4.5}
\def\dotfill{%
\hbox{$\m@th\mkern \@dotsep mu\hbox{.}\mkern \@dotsep mu$}\hfill}
\def\@tocline#1#2#3#4#5#6#7{\relax
  \ifnum #1>\c@tocdepth 
  \else
    \def\@toclevel{#1}%
    \par \addpenalty\@secpenalty\addvspace{#2}%
    \begingroup 
        \hyphenpenalty\@M
        \@ifempty{#4}{%
          \@tempdima\csname r@tocindent\number#1\endcsname\relax
        }{%
          \@tempdima#4\relax
        }%
		\parindent .5pc
		\leftskip#3\relax \advance\leftskip\@tempdima\relax
        \rightskip\@pnumwidth plus4em \parfillskip-\@pnumwidth
        #5\leavevmode\hskip-\@tempdima #6\relax
        \leaders\dotfill
		\hbox to\@pnumwidth{\@tocpagenum{#7}}\par
        \nobreak
    \endgroup
  \fi
}
\makeatother 

\fi

\ifx\PrintBook\undefined
\def\Chapter{\section}
\def\Section{\subsection}
\ifx\setCACAA\undefined
\makeatletter
\renewcommand{\@indextitlestyle}{%
\twocolumn[\section{\indexname}]%
\def\IndexSpace{off}%
}
\makeatother 
\ifx\PrintPaper\undefined
\thanks{\href{mailto:Aleks\_Kleyn@MailAPS.org}{Aleks\_Kleyn@MailAPS.org}}
\ePrints{1102.1776,1201.4158}
\ifx\Semafor\ValueOff
\thanks{\ \ \ \url{http://AleksKleyn.dyndns-home.com:4080/}\ \ \ \ \ \url{http://arxiv.org/a/kleyn\_a\_1}}
\thanks{\ \ \ \ \url{http://AleksKleyn.blogspot.com/}}
\fi
\fi
\fi
\else
\def\Chapter{\chapter}
\def\Section{\section}

\makeatletter
\def\@maketitle{%
  \cleardoublepage \thispagestyle{empty}%
  \begingroup \topskip\z@skip
  \null\vfil
  \begingroup
  \def\and{\par\medskip}\centering
  \bfseries\authors\par\bigskip
  \par\vspace{24pt}%
  \LARGE\bfseries \centering
  \openup\medskipamount
  \@title
  \par
  \ifx\subtitle\undefined
  \else
  \centerline{\ }
  \centerline{\emph\subtitle}
  \fi
  \ifx\subtitleA\undefined
  \else
  \centerline{\emph\subtitleA}
  \fi
  \ifx\edition\undefined
  \else
  \centerline{\emph\edition}
  \fi
  \endgroup
  \vfill
\noindent
\href{mailto:Aleks\_Kleyn@MailAPS.org}{Aleks\_Kleyn@MailAPS.org}
\newline
\url{http://AleksKleyn.dyndns-home.com:4080/}
\newline
\url{http://arxiv.org/a/kleyn\_a\_1}
\newline
\url{http://AleksKleyn.blogspot.com/}
  \newpage\thispagestyle{empty}
  \begin{center}
    \ifx\@empty\@subjclass\else\@setsubjclass\fi
    \ifx\@empty\@keywords\else\@setkeywords\fi
    \ifx\@empty\@translators\else\vfil\@settranslators\fi
    \ifx\@empty\thankses\else\vfil\@setthanks\fi
  \end{center}
  \vfil
  \@setabstract
\vfil
  \def\Temp{0000}
  \ifx\copyrightyear\Temp
  \else
  \begin{center}
\begin{tabular}{@{}c}
Copyright\ \copyright\ \copyrightyear\ \copyrightholder
\\
All rights reserved.
\end{tabular}
  \end{center}
  \fi
  \ifx\Publisher\undefined%
  \else
  \begin{center}
\begin{tabular}{@{}c}
\Publisher
\end{tabular}
  \end{center}
  \fi
  \ifx\ISBN\undefined%
  \else%
 \begin{center}
\begin{tabular}{@{}r@{\ }l}
ISBN:&\ISBN
\\
ISBN-13:&\ISBNa
\end{tabular}
  \end{center}
  \fi%
  \ifx\titleNote\undefined
  \else
  \par\vspace{24pt}%
  \centerline{\mdseries\titleNote}
	  \centerline{\Title}
	  \ifx\Subtitle\undefined
	  \else
	  \centerline{\emph\Subtitle}
	  \fi
	  \ifx\Edition\undefined
	  \else
	  \centerline{\Edition}
	  \fi
	  \centerline{\Authors}
  \fi
  \endgroup}
\def\chapter{%
	\clearpage
  \thispagestyle{plain}\global\@topnum\z@
  \@afterindenttrue \secdef\@chapter\@schapter}
\renewcommand{\@indextitlestyle}{%
\twocolumn[\chapter{\indexname}]%
\def\IndexSpace{off}%
\let\@secnumber\@empty
\chaptermark{\indexname}%
}
\makeatother 
\email{\href{mailto:Aleks\_Kleyn@MailAPS.org}{Aleks\_Kleyn@MailAPS.org}}
\ePrints{1102.1776,1201.4158}
\ifx\Semafor\ValueOff
\urladdr{\url{http://AleksKleyn.dyndns-home.com:4080/}}
\urladdr{\url{http://arxiv.org/a/kleyn\_a\_1}}
\urladdr{\url{http://AleksKleyn.blogspot.com/}}
\fi
\fi

\newcommand\FrameEqRef[3][]
{
\fbox{%
\eqRef{#2}{#3}%
$\ \ \ $%
\DrawEq[#1]{#2}-
}
}

\newcommand\FrameCiteBib[1]
{
\noindent
\fbox{
\begin{minipage}{350pt}
\setlength{\parindent}{6mm}
\citeBib{#1}
\ShowCiteBib{#1}
\end{minipage}
}
}

\newcommand\arXivOldRef{http://arxiv.org/PS_cache/}
\newcommand\arXivRef{http://arxiv.org/pdf/}
\newcommand\RgRef{https://www.researchgate.net/publication/}
\newcommand\AmazonRef{http://www.amazon.com/s/ref=nb_sb_noss?url=search-alias=aps&field-keywords=aleks+kleyn}
\newcommand\DefRgRef[1]
{
\def\Tempb{#1}
\ifx\Tempa\Tempb
\def\wRef{\RgRef #1}
\fi
}
\newcommand\DefArXivOldRef[2]
{
\def\Tempb{#1}
\ifx\Tempa\Tempb
\def\wRef{\arXivOldRef #2.pdf}
\fi
}
\newcommand\DefArXivRef[2]
{
\def\Tempb{#1}
\ifx\Tempa\Tempb
\def\wRef{\arXivRef #1#2.pdf}
\fi
}
\newcommand\DefAmazonRef[1]
{
\def\Tempb{#1}
\ifx\Tempa\Tempb
\def\wRef{\AmazonRef}
\fi
}
\newcommand\DefCaCaaRef[1]
{
\def\Tempb{CACAA.#1}
\ifx\Tempa\Tempb
\def\wRef{http://www.cliffordanalysis.com/}
\fi
}
\newcommand\wRefDef[2]
{
\def\Tempa{#1}
\DefArXivOldRef{0405.027}{gr-qc/pdf/0405/0405027v3}
\DefArXivOldRef{0405.028}{gr-qc/pdf/0405/0405028v5}
\DefArXivOldRef{0412.391}{math/pdf/0412/0412391v4}
\DefArXivOldRef{0612.111}{math/pdf/0612/0612111v2}
\DefArXivOldRef{0701.238}{math/pdf/0701/0701238v6}
\DefArXivOldRef{0702.561}{math/pdf/0702/0702561v3}
\DefArXivRef{0707.2246}{v2}
\DefArXivRef{0803.3276}{v3}
\DefArXivRef{0812.4763}{v7}
\DefArXivRef{0906.0135}{v3}
\DefArXivRef{0909.0855}{v5}
\DefArXivRef{0912.3315}{v3}
\DefArXivRef{0912.4061}{v2}
\DefArXivRef{1001.4852}{}
\DefArXivRef{1003.3714}{v2}
\DefArXivRef{1003.1544}{v2}
\DefArXivRef{1006.2597}{v2}
\DefArXivRef{1011.3102}{}
\DefArXivRef{1102.1776}{}
\DefArXivRef{1104.5197}{}
\DefArXivRef{1105.4307}{}
\DefArXivRef{1107.1139}{}
\DefArXivRef{1107.5037}{}
\DefArXivRef{1111.6035}{}
\DefArXivRef{1202.6021}{v2}
\DefArXivRef{1211.6965}{}
\DefArXivRef{1302.7204}{v1}
\DefArXivRef{1305.4547}{}
\DefArXivRef{1310.5591}{}
\DefArXivRef{1502.04063}{v2}
\DefArXivRef{1505.03625}{v1}
\DefArXivRef{1506.00061}{}
\DefArXivRef{1601.03259}{v3}
\DefArXivRef{1801.01628}{}
\DefArXivRef{1908.04418}{}
\DefArXivRef{2020.06.01}{}
\DefArXivRef{2021.01.06}{}
\DefArXivRef{2207.06506}{}
\DefRgRef{322019412}
\DefRgRef{323966352}
\DefRgRef{348311191}
\DefAmazonRef{8433-5163}
\DefAmazonRef{8443-0072}
\DefAmazonRef{4776-3181}
\DefAmazonRef{4975-6381}
\DefAmazonRef{4993-2400}
\DefAmazonRef{5059-9176}
\DefAmazonRef{5114-6019}
\DefAmazonRef{5148-4632}
\DefAmazonRef{5410-9916}
\DefAmazonRef{6860-2955}
\DefAmazonRef{0767-8264}
\DefAmazonRef{8428-0408}
\DefCaCaaRef{01.109}
\DefCaCaaRef{01.291}
\DefCaCaaRef{02.097}
\DefCaCaaRef{04.001}
\DefCaCaaRef{05.001}
\DefCaCaaRef{06.121}
 \def\Tempb{GJSFRA.13.1.39}
\ifx\Tempa\Tempb
\def\wRef{http://www.cliffordanalysis.com/}
\fi
\externaldocument[#1-#2-]{\FilePrefix #1.#2}[\wRef]
}

\makeatletter
\newcommand\org@maketitle{}
\let\org@maketitle\maketitle
\def\maketitle{%
\hypersetup{pdftitle={\@title}}%
\hypersetup{pdfauthor={\authors}}%
\hypersetup{pdfsubject=\@keywords}%
\ifx\UseRussian\Defined
\pdfbookmark[1]{\@title}{TitleRussian}
\else
\pdfbookmark[1]{\@title}{TitleEnglish}
\fi
\org@maketitle
}
\def\make@stripped@name#1{%
\begingroup
\escapechar\m@ne
\global\let\newname\@empty
\protected@edef\Hy@tempa{\CurrentLanguage #1}%
\edef\@tempb{%
\noexpand\@tfor\noexpand\Hy@tempa:=%
\expandafter\strip@prefix\meaning\Hy@tempa
}%
\@tempb\do{%
\if\Hy@tempa\else
\if\Hy@tempa\else
\xdef\newname{\newname\Hy@tempa}%
\fi
\fi
}%
\endgroup
}%
\newenvironment{enumBib}{%
\BibTitle
\advance\@enumdepth \@ne
\edef\@enumctr{enum\romannumeral\the\@enumdepth}\list
{\csname biblabel\@enumctr\endcsname}{\usecounter
{\@enumctr}\def\makelabel##1{\hss\llap{\upshape##1}}}
}{%
\endlist
}

\makeatletter

\newcommand{\BiblioItem}[2]
{
\def\Semafor{off}
\@ifundefined{\LanguagePrefix ViewCite#1}{}{%
\def\Semafor{on}%
}%
\ifx\Semafor\ValueOff
\@ifundefined{xRefDef#1}{}{%
\def\Semafor{on}%
}%
\fi
\ifx\Semafor\ValueOn
\ifx\IndexState\ValueOff
\begin{enumBib}
\def\IndexState{on}
\fi
\item \label{\LanguagePrefix bibitem: #1}#2%
\fi
}
\makeatother
\newcommand{\OpenBiblio}
{
\def\IndexState{off}
}
\newcommand{\CloseBiblio}
{
\ifx\IndexState\ValueOn
\end{enumBib}
\def\IndexState{off}
\fi
}

\makeatletter
\def\StartCite{[}%
\def\citeBib#1{\protect\showCiteBib#1,endCite,}%
\def\endCite{endCite}%
\def\showCiteBib#1,{\def\temp{#1}%
\ifx\temp\endCite
]%
\def\StartCite{[}%
\else
\StartCite\LanguagePrefix \ref{\LanguagePrefix bibitem: #1}%
\@ifundefined{\LanguagePrefix ViewCite#1}{%
\NameDef{\LanguagePrefix ViewCite#1}{}%
}{%
}%
\def\StartCite{, }%
\expandafter\showCiteBib%
\fi}%
\makeatother

\newcommand{\arp}{\ar @{-->}}

\newcommand\Bundle[1]{{\mathbb #1}}
\newcommand{\bundle}[4]%
{%
\def\tempa{}%
\def\tempb{#3}%
\def\tempc{#1}%
\ifx\tempa\tempb%
\ifx\tempa\tempc%
#2%
\else%
\xymatrix{#2:#1\arp[r]&#4}%
\fi%
\else%
\ifx\tempa\tempc%
#2[#3]%
\else%
\xymatrix{#2[#3]:#1\arp[r]&#4}%
\fi%
\fi%
}%
\makeatletter
\newcommand{\AddIndex}[2]%
{%
\@ifundefined{RefIndex#2}{%
\xNameDef{RefIndex#2}{:}%
\@namedef{LabelIndex}{\label{index: #2::}}%
}{%
\addtocounter{Index}{1}%
\xNameDef{RefIndex#2}{\@nameuse{RefIndex#2},\arabic{Index}}%
\@namedef{LabelIndex}{\label{index: #2:\arabic{Index}}}%
}%
\@nameuse{LabelIndex}%
{\bf #1}%
}%

\newcommand{\Index}[2]%
{%
\@ifundefined{RefIndex#2}{%
\def\Semafor{off}%
}{%
\def\Semafor{on}%
}%
\ifx\Semafor\ValueOn%
\def\tempa{}%
\def\tempb{#2}%
\ifx\IndexState\ValueOff%
\ifx\setCACAA\undefined
\begin{theindex}%
\else
\section*{\indexname}%
\begin{itemize}
\fi
\def\IndexState{on}%
\fi%
\ifx\IndexSpace\ValueOn%
\indexspace%
\def\IndexSpace{off}%
\fi%
\item #1%
\ifx\tempa\tempb%
\else%
\edef\PageRefs{\@nameuse{RefIndex#2}}
\def\Sep{\ }%
\@for\PageRef:=\PageRefs\do{%
\Sep
\pageref{index: #2:\PageRef}%
\def\Sep{,\ }%
}%
\fi%
\fi%
}%

\newcommand{\Symb}[4]%
{%
\def\Semafor{off}%
\@ifundefined{ViewSymbol#2}{%
\@ifundefined{ViewSymbol#2(#3-#4)}{%
}{%
\def\Semafor{on}
\edef\ThisSymbol{#2(#3-#4)}%
}%
}{%
\def\Semafor{on}%
\edef\ThisSymbol{#2}%
}%
\ifx\Semafor\ValueOn%
\ifx\IndexState\ValueOff%
\ifx\setCACAA\undefined
\begin{theindex}%
\else
\section*{\indexname}%
\begin{itemize}
\fi
\def\IndexState{on}%
\fi%
\ifx\IndexSpace\ValueOn%
\indexspace%
\def\IndexSpace{off}%
\fi%
\edef\Symbols{\@nameuse{ViewSymbol\ThisSymbol}}%
\@for\Symbol:=\Symbols\do{%
\edef\Temp{ViewSymbol\ThisSymbol:::\Symbol}%
\item $\displaystyle\textcolor{SymbColor}{\@nameuse{\Temp}}$
\ \ #1
\edef\PageRefs{\@nameuse{RefSymbol\ThisSymbol===\Symbol}}
\def\Sep{}%
\@for\PageRef:=\PageRefs\do{%
\Sep
\pageref{symbol: \ThisSymbol:\PageRef}%
\def\Sep{,\ }%
}%
}%
\fi
}

\newcommand{\Symba}[2]
{
\def\Semafor{off}
\@ifundefined{ViewSymbol#2}{%
}{%
\def\Semafor{on}
}%
\ifx\Semafor\ValueOn
\ifx\IndexState\ValueOff
\begin{theindex}
\def\IndexState{on}
\fi
\ifx\IndexSpace\ValueOn
\indexspace
\def\IndexSpace{off}
\fi
\item $\displaystyle\@nameuse{ViewSymbol#2}$\ \ #1
\edef\PageRefs{\@nameuse{RefSymbol#2}}
\def\Sep{}%
\@for\PageRef:=\PageRefs\do{%
\Sep
\pageref{symbol: #2:\PageRef}%
\def\Sep{,\ }%
}%
\fi
}

\makeatother

\newcommand{\SetIndexSpace}%
{%
\def\IndexSpace{on}%
}%

\newcommand{\OpenIndex}
{
\def\IndexState{off}
}
\newcommand{\CloseIndex}
{
\ifx\IndexState\ValueOn
\ifx\setCACAA\undefined
\end{theindex}
\else
\end{itemize}
\fi
\def\IndexState{off}
\fi
}

\def\LastMemo{LastMemo}%
\def\MemoList#1//{\def\temp{#1}%
\ifx\temp\LastMemo
\else%
\setlength{\parindent}{5mm}
\par
\BlueText{#1}%
\expandafter\MemoList%
\fi%
}     

%% file: Statement.tex
\ifx\SelectlEnglish\undefined
\ifx\UseRussian\undefined
\def\SelectlEnglish{}
\fi
\fi
\ifx\SelectlEnglish\undefined
\newcommand\TheLanguage{Russian}%
\else
\newcommand\TheLanguage{English}%
\fi

\newcommand\LanguagePrefix{}%
\newcommand\CurrentLanguage{\TheLanguage.}%
\newcommand\xRefDef[1]
{
\wRefDef{#1}{\TheLanguage}
\NameDef{xRefDef#1}{}%
}

\input{\FilePrefix Statement.\TheLanguage}

\theoremstyle{definition}
\theoremstyle{remark}
\renewenvironment{proof}[1][\proofname]
{\par{\sc #1. }}{\qed}%

\ifx\PrintBook\undefined
\numberwithin{footnote}{section}
\numberwithin{Hfootnote}{section}
\else
\numberwithin{section}{chapter}
\numberwithin{footnote}{chapter}
\numberwithin{Hfootnote}{chapter}
\fi

\ifx\Presentation\undefined
\numberwithin{equation}{section}
\numberwithin{figure}{section}
\numberwithin{table}{section}
\numberwithin{Item}{section}
\fi

\def\OneTheorem{the theorem }
\def\ManyTheorems{theorems }

\def\LastRef{LastRef}%
\def\PrintTheorem#1{\PrtTheorem#1||LastRef||}%
\def\PrtTheorem#1||#2||{\def\temp{#2}%
\ifx\temp\LastRef%
\gdef\TheoremName{\OneTheorem}
\else%
\gdef\TheoremName{\ManyTheorems}
\expandafter\DoNothing
\fi%
}
\def\DoNothing#1||
{
\def\temp{#1}%
\ifx\temp\LastRef%
\else%
\expandafter\DoNothing
\fi%
}
\def\Refs#1{\RefList#1||LastRef||}%
\def\RefList#1||{\def\temp{#1}%
\ifx\temp\LastRef%
\else%
\TheoremSep\DoTheoremRef{#1}
\gdef\TheoremSep{, }%
\expandafter\RefList
\fi%
}
\def\DoTheoremRef#1{\DoRefs#1|:LastTRef|:}%
\def\LastTRef{LastTRef}%
\def\DoRefs#1|:#2|:{\def\tempA{#2}%
\ifx\tempA\LastTRef%
\RefTheorem{#1}%
\else%
\refTheorem{#1}{#2}%
\fi%
\expandafter\RefNothing%
}
\def\RefNothing#1
{%
}%
\newenvironment{ProofRefA}[1]{%
\gdef\TheoremSep{}%
{\PrintTheorem{#1}}
{\sc Proof of \TheoremName \Refs{#1}.}
}%
{\qed}

\ifx\texFuture\Defined
\begin{ProofRefA}{Polylinear map is geometric object||%
define homomorphism of vector space|:left||define homomorphism of vector space|:right}
\end{ProofRefA}

\begin{ProofRefA}{define homomorphism of vector space|:left}
\end{ProofRefA}
\fi

\definecolor{TheoremColor}{rgb}{.94,.96,.94}%
\definecolor{LemmaColor}{rgb}{.94,.96,.98}%
\definecolor{DefinitionColor}{rgb}{.94,.94,.95}%
\definecolor{QuestionColor}{rgb}{.93,.93,.96}%
\definecolor{SummaryColor}{rgb}{.92,.93,.95}%
\definecolor{ExampleColor}{rgb}{.94,.95,.90}%
\newcommand{\DefineTheoremColor}{%
\colorlet{shadecolor}{TheoremColor}
}%
\newcommand{\DefineLemmaColor}{%
\colorlet{shadecolor}{LemmaColor}
}%
\newcommand{\DefineDefinitionColor}{%
\colorlet{shadecolor}{DefinitionColor}
}%
\newcommand{\DefineQuestionColor}{%
\colorlet{shadecolor}{QuestionColor}
}%
\newcommand{\DefineSummaryColor}{%
\colorlet{shadecolor}{SummaryColor}
}%
\newcommand{\DefineExampleColor}{%
\colorlet{shadecolor}{ExampleColor}
}%

\newenvironment{Shaded}{%
  \MakeFramed {\FrameRestore}}%
 {\endMakeFramed}

\newenvironment{Convention}
{
\begin{convention}
}
{
\qed
\end{convention}
}

\newenvironment{Definition}
{
\begin{definition}
}
{
\qed
\end{definition}
}

\newenvironment{ShadedDefinition}
{
\DefineDefinitionColor
\begin{Shaded}
\begin{Definition}
}{%
\end{Definition}%
\end{Shaded}%
}

\newenvironment{ShadedTheorem}
{
\DefineTheoremColor%
\begin{Shaded}%
\begin{theorem}%
}{%
\end{theorem}%
\end{Shaded}%
}

\newenvironment{Corollary}
{
\begin{corollary}
}
{
\qed
\end{corollary}
}

\newenvironment{ShadedLemma}
{
\DefineLemmaColor%
\begin{Shaded}%
\begin{lemma}
}{%
\end{lemma}
\end{Shaded}%
}

\newenvironment{Example}
{
\begin{example}
}
{
\qed
\end{example}
}

\newenvironment{Question}
{
\begin{question}
}
{
\qed
\end{question}
}

\newenvironment{Summary}
{
\begin{summary}
}
{
\qed
\end{summary}
}

\newenvironment{Remark}
{
\begin{remark}
}
{
\qed
\end{remark}
}

\def\DefinitionStyle{ShadedDefinition}%

\newcommand{\DefDefinition}[3][]{%
\def\Temp{}%
\def\Tempa{#1}%
\ifx\Tempa\Temp%
\def\TempB{\AddEq}%
\else%
\def\TempB{\AddEq[#1]}%
\fi%
\TempB{definition: #2}
{
\begin{\DefinitionStyle}
\labelDefinition{#2}
#3
\end{\DefinitionStyle}
}
}%

\newcommand{\DefLabeledDefinition}[4][]{%
\def\Temp{}%
\def\Tempa{#1}%
\ifx\Tempa\Temp%
\def\TempB{\AddEq}%
\else%
\def\TempB{\AddEq[#1]}%
\fi%
\TempB{definition: #2}
{
\begin{\DefinitionStyle}
\labelDefinition{#2:: #3}
#4
\end{\DefinitionStyle}
}
}%

\newcommand{\DefDefinitionNote}[4][]{%
\def\Temp{}%
\def\Tempa{#1}%
\ifx\Tempa\Temp%
\def\TempB{\AddEq}%
\else%
\def\TempB{\AddEq[#1]}%
\fi%
\TempB{definition: #2}
{
\begin{\DefinitionStyle}
\labelDefinition{#2}
#3
\end{\DefinitionStyle}
\footnotetext{\,
#4
}
}
}%

\newcommand{\DefLabeledDefinitionNote}[5][]{%
\def\Temp{}%
\def\Tempa{#1}%
\ifx\Tempa\Temp%
\def\TempB{\AddEq}%
\else%
\def\TempB{\AddEq[#1]}%
\fi%
\TempB{definition: #2}
{
\begin{ShadedDefinition}
\labelDefinition{#2:: #3}
#4
\end{ShadedDefinition}
\footnotetext{\,
#5
}
}
}%

\newcommand\ShowDefinition[1]{\ShowEq{definition: #1}}

\def\TheoremStyle{ShadedTheorem}%

\newcommand{\DefTheorem}[3][]{%
\def\Temp{}%
\def\Tempa{#1}%
\ifx\Tempa\Temp%
\def\TempB{\AddEq}%
\else%
\def\TempB{\AddEq[#1]}%
\fi%
\TempB{theorem: #2}
{
\begin{\TheoremStyle}
\labelTheorem{#2}
#3
\end{\TheoremStyle}
}
}%

\newcommand{\DefLabeledTheorem}[4][]{%
\def\Temp{}%
\def\Tempa{#1}%
\ifx\Tempa\Temp%
\def\TempB{\AddEq}%
\else%
\def\TempB{\AddEq[#1]}%
\fi%
\TempB{theorem: #2}
{
\begin{\TheoremStyle}
\labelTheorem{#2:: #3}
#4
\end{\TheoremStyle}
}
}%

\newcommand{\DefTheoremNote}[4][]{%
\def\Temp{}%
\def\Tempa{#1}%
\ifx\Tempa\Temp%
\def\TempB{\AddEq}%
\else%
\def\TempB{\AddEq[#1]}%
\fi%
\TempB{theorem: #2}
{
\begin{\TheoremStyle}
\label{theorem: #2}
#3
\end{\TheoremStyle}
\footnotetext{\,
#4
}
}
}%

\newcommand{\DefLabeledTheoremNote}[5][]{%
\def\Temp{}%
\def\Tempa{#1}%
\ifx\Tempa\Temp%
\def\TempB{\AddEq}%
\else%
\def\TempB{\AddEq[#1]}%
\fi%
\TempB{theorem: #2}
{
\begin{\TheoremStyle}
\labelTheorem{#2:: #3}
#4
\end{\TheoremStyle}
\footnotetext{\,
#5
}
}
}%

\newcommand{\ShowTheorem}[1]{%
\ShowEq{theorem: #1}%
}%

\def\LemmaStyle{ShadedLemma}%

\newcommand\DefLabeledLemma[4]{%
\AddEq{lemma: #1}
{
\begin{\LemmaStyle}
\labelLemma{#1:: #2}
#3
\end{\LemmaStyle}
}
}

\newcommand{\ShowLemma}[1]{\ShowEq{lemma: #1}}%

\def\CorollaryStyle{ShadedCorollary}%

\newcommand{\DefCorollary}[2]{%
\AddEq{corollary: #1}
{
\begin{\CorollaryStyle}
\labelCorollary{#1}
#2
\end{\CorollaryStyle}
}
}

\newcommand{\ShowCorollary}[1]{%
\ShowEq{corollary: #1}%
}%

\newcommand{\DefProof}[3][]{%
\def\Temp{}%
\def\Tempa{#1}%
\ifx\Tempa\Temp%
\def\TempB{\AddEq}%
\else%
\def\TempB{\AddEq[#1]}%
\fi%
\TempB{proof: #2}
{%
\begin{proof}
#3
\end{proof}%
}
}

\newcommand{\DefProofSloppy}[3][]{%
\def\Temp{}%
\def\Tempa{#1}%
\ifx\Tempa\Temp%
\def\TempB{\AddEq}%
\else%
\def\TempB{\AddEq[#1]}%
\fi%
\TempB{proof: #2}
{%
\begin{sloppypar}
\begin{proof}
#3
\end{proof}%
\end{sloppypar}
}
}

\newcommand{\DefProofRef}[4][]{%
\def\Temp{}%
\def\Tempa{#1}%
\ifx\Tempa\Temp%
\def\TempB{\AddEq}%
\else%
\def\TempB{\AddEq[#1]}%
\fi%
\TempB{proof: #2}
{%
\begin{ProofRef}{#2}{#3}
#4
\end{ProofRef}%
}
}

\newcommand{\ShowProof}[1]{\ShowEq{proof: #1}}%

\newcommand{\ProveTheorem}[1]
{
\ShowTheorem{#1}
\ShowProof{#1}
}%

\def\ExampleStyle{ShadedExample}%

\newcommand{\DefExample}[2]{%
\AddEq{example: #1}
{
\begin{\ExampleStyle}
\labelExample{#1}
#2
\end{\ExampleStyle}
}
}%

\newcommand{\DefLabeledExample}[3]{%
\AddEq{example: #1}
{
\begin{\ExampleStyle}
\labelExample{#1:: #2}
#3
\end{\ExampleStyle}
}
}%

\newcommand{\ShowExample}[1]{%
\ShowEq{example: #1}%
}%

\newcommand\DefRemark[2]{%
\AddEq{remark: #1}
{
\begin{Remark}
\labelRemark{#1}
#2
\end{Remark}
}
}%

\newcommand\DefLabeledRemark[3]{%
\AddEq{remark: #1}
{
\begin{Remark}
\labelRemark{#1:: #2}
#3
\end{Remark}
}
}%

\newcommand\DefLabeledSloppyRemark[3]{%
\AddEq{remark: #1}
{
\begin{sloppypar}
\begin{Remark}
\labelRemark{#1:: #2}
#3
\end{Remark}
\end{sloppypar}
}
}%

\newcommand\DefText[3][]{%
\def\Temp{}%
\def\Tempa{#1}%
\ifx\Tempa\Temp%
\def\TempB{\AddEq}%
\else%
\def\TempB{\AddEq[#1]}%
\fi%
\TempB{text: #2}
{
#3
}
}%

\newcommand{\ShowRemark}[1]{%
\ShowEq{remark: #1}%
}%

\newcommand{\ShowText}[1]{%
\ShowEq{text: #1}%
}%

\def\ConventionStyle{ShadedConvention}%

\newcommand\DefConvention[3][]{%
\def\Temp{}%
\def\Tempa{#1}%
\ifx\Tempa\Temp%
\def\TempB{\AddEq}%
\else%
\def\TempB{\AddEq[#1]}%
\fi%
\TempB{convention: #2}
{
\begin{\ConventionStyle}
\labelConvention{#2}
#3
\end{\ConventionStyle}
}
}%

\newcommand\DefLabeledConvention[4][]{%
\def\Temp{}%
\def\Tempa{#1}%
\ifx\Tempa\Temp%
\def\TempB{\AddEq}%
\else%
\def\TempB{\AddEq[#1]}%
\fi%
\TempB{convention: #2}
{
\begin{\ConventionStyle}
\labelConvention{#2:: #3}
#4
\end{\ConventionStyle}
}
}%

\newcommand{\ShowConvention}[1]{%
\ShowEq{convention: #1}%
}%

\def\SummaryStyle{ShadedSummary}%

\newcommand\DefSummary[2]{%
\AddEq{summary: #1}
{
\begin{\SummaryStyle}
\labelSummary{#1}
#2
\end{\SummaryStyle}
}
}%

\newcommand{\DefCiteBib}[2]{\AddEq{citeBib: #1}{#2}}%
\newcommand{\DefBiblioItem}[1]{\BiblioItem{#1}{\ShowCiteBib{#1}}}%

\newcommand\ShowCiteBib[1]{\ShowEq{citeBib: #1}}

\makeatletter
\newcommand\StartLabelItem[1][theorem]%
{
\expandafter \def \csname%
 theenumi\expandafter \endcsname \expandafter {\csname the#1\expandafter%
 \endcsname.\expandafter \@arabic \csname c@enumi\endcsname }%
}
\newcommand\StopLabelItem[1][theorem]%
{
\counterwithout{enumi}{#1}%

}
\makeatother
\newcommand\labelConvention[1]{\label{convention: #1}}%
\newcommand\labelTheorem[1]{\label{theorem: #1}}%
\newcommand\labelLemma[1]{\label{lemma: #1}}%
\newcommand\labelCorollary[1]{\label{corollary: #1}}%
\newcommand\labelStatement[1]{\label{statement: #1}}
\newcommand\labelDefinition[1]{\label{definition: #1}}%
\newcommand\labelExample[1]{\label{example: #1}}%
\newcommand\labelRemark[1]{\label{remark: #1}}%
\newcommand\labelSection[1]{\label{section: #1}}%
\newcommand\labelItem[1]{\label{item: #1}}%

\newcommand\labelSummary[1]{\label{summary: #1}}
\newcommand\labelFootnote[1]{\label{footnote: #1}}
\makeatletter
\newcommand\xRef[2][]%
{%
\def\Semafor{on}%
\def\Temp{}%
\edef\Tempa{#1}%
\ifx\Tempa\Temp%
\def\Semafor{off}%
\fi%
\@ifundefined{xRefDef#1}{%
\def\Semafor{off}%
}{}%
\ifx\Semafor\ValueOff%
\penalty-1\ref{#2}%
\else%
\penalty-1\citeBib{#1}-\ref{#1-\TheLanguage-#2}%
\fi%
}%
\makeatother
\newcommand\RefConvention[2][]{\xRef[#1]{convention: #2}}
\newcommand\RefTheorem[2][]{\xRef[#1]{theorem: #2}}
\newcommand\refTheorem[3][]{\xRef[#1]{theorem: #2:: #3}}
\newcommand\RefLemma[2][]{\xRef[#1]{lemma: #2}}
\newcommand\refLemma[3][]{\xRef[#1]{lemma: #2:: #3}}
\newcommand\RefStatement[2][]{\xRef[#1]{statement: #2}}
\newcommand\RefCorollary[2][]{\xRef[#1]{corollary: #2}}
\newcommand\RefDefinition[2][]{\xRef[#1]{definition: #2}}
\newcommand\refDefinition[3][]{\xRef[#1]{definition: #2:: #3}}
\newcommand\RefExample[2][]{\xRef[#1]{example: #2}}
\newcommand\refExample[3][]{\xRef[#1]{example: #2:: #3}}
\newcommand\RefRemark[2][]{\xRef[#1]{remark: #2}}
\newcommand\refRemark[3][]{\xRef[#1]{remark: #2:: #3}}

\newcommand\RefChapter[2][]{\xRef[#1]{chapter: #2}}
\newcommand\RefItem[2][]{\xRef[#1]{item: #2}}

\newcommand\refFootnote[3][]{\,${}^{\xRef[#1]{footnote: #2:: #3}}$}

\makeatletter
\newcommand\AddFootnote[2]
{
\stepcounter{Hfootnote}
\stepcounter{footnote}
\global \let \Hy@saved@currentHref \@currentHref
\hyper@makecurrent {Hfootnote}
\global \let \Hy@footnote@currentHref \@currentHref 
\footnotetext{
\labelFootnote{#1}
\,#2
}
}

\newcommand\ShowPrevFootnote[1]
{
\addtocounter{Hfootnote}{-1}
\addtocounter{footnote}{-1}
\global \let \Hy@saved@currentHref \@currentHref
\hyper@makecurrent {Hfootnote}
\global \let \Hy@footnote@currentHref \@currentHref 
\footnotetext{
\,#1
}
}

\newcommand\ShowNextFootnote[1]
{
\stepcounter{Hfootnote}
\stepcounter{footnote}
\global \let \Hy@saved@currentHref \@currentHref
\hyper@makecurrent {Hfootnote}
\global \let \Hy@footnote@currentHref \@currentHref 
\footnotetext{
\,#1
}
}
\makeatother

\newcommand\DefFootnote[3][]%
{%
\def\Temp{}%
\def\Tempa{#1}%
\ifx\Tempa\Temp%
\def\TempB{\AddEq}%
\else%
\def\TempB{\AddEq[#1]}%
\fi%
\TempB{footnote: #2}
{
\AddFootnote{#2}{#3}
}
}

\newcommand\DefLabeledFootnote[4][]%
{%
\def\Temp{}%
\def\Tempa{#1}%
\ifx\Tempa\Temp%
\def\TempB{\AddEq}%
\else%
\def\TempB{\AddEq[#1]}%
\fi%
\TempB{footnote: #2}
{
\AddFootnote{#2:: #3}{#4}
}
}

\newcommand{\ShowFootnote}[1]{%
\ShowEq{footnote: #1}%
}%

\newcommand\DefRef[3][]%
{%
\def\Temp{}%
\def\Tempa{#1}%
\ifx\Tempa\Temp%
\def\TempB{\AddEq}%
\else%
\def\TempB{\AddEq[#1]}%
\fi%
\TempB{ref: #2}
{
#3
}
}

\newcommand\ShowRef[1]{\ShowEq{ref: #1}}

\ifx\texFuture\Defined
\newcommand\TwoColText[2]
{
\noindent
\begin{minipage}{0.5\textwidth}
\setlength{\parindent}{4mm}
#1
\end{minipage}
\begin{minipage}{0.5\textwidth}
\setlength{\parindent}{4mm}
#2
\end{minipage}
}
\fi

\ifx\UseParacol\Defined
\newcommand\TwoColText[2]
{
\begin{paracol}{2}
#1
\switchcolumn
#2
\end{paracol}
}

\else
\newcommand\TwoColText[2]
{
\noindent
\begin{tabular}{@{}l|r}
\begin{minipage}{0.49\textwidth}
\setlength{\parindent}{4mm}
#1
\end{minipage}
&
\begin{minipage}{0.49\textwidth}
\setlength{\parindent}{4mm}
#2
\end{minipage}
\end{tabular}
}
\fi

%% file: Statement.English.tex

\newenvironment{ProofRef}[2]{%

{\sc Proof of theorem}
\def\Temp{}%
\edef\Tempa{#2}%
\ifx\Tempa\Temp%
\RefTheorem{#1}.
\else
\refTheorem{#1}{#2}.
\fi
}%
{\qed}

\newcommand{\DefLemmaProof}[3][]{%
\def\Temp{}%
\def\Tempa{#1}%
\ifx\Tempa\Temp%
\def\TempB{\AddEq}%
\else%
\def\TempB{\AddEq[#1]}%
\fi%
\TempB{proof: #2}
{%
{\sc Proof.}
#3
\hfill\(\odot\)
}
}

\author{Aleks Kleyn}
\ifx\setCACAA\undefined
\newtheorem{theorem}{Theorem}[section]
\newtheorem{corollary}[theorem]{Corollary}
\newtheorem{example}[theorem]{Example}
\newtheorem{definition}[theorem]{Definition}
\newtheorem{remark}[theorem]{Remark}
\newtheorem{question}[theorem]{Question}
\newtheorem{summary}[theorem]{Summary of Results}
\newtheorem{lemma}[theorem]{Lemma}
\else
\theoremstyle{definition}
\theoremstyle{remark}
\fi
\newtheorem{Statement}[theorem]{Statement}
\newtheorem{convention}[theorem]{Convention}

\ifx\PrintBook\undefined
\newcommand{\BibTitle}{%
\section{References}%
}
\else
\newcommand{\BibTitle}{%
\chapter{References}%
}
\fi

%% file: Prolog.Parm.tex

\newcommand\Multiply[2]{#1#2}

\AddEq{=DA1DA2}
{
\def\DFDT{D1 D2 }%
\def\MF{r1:D1->D2 }%
\def\DF{1}%
\def\DT{2}%
\def\VF{1}%
\def\VT{2}%
\def\MapE{hgf}%
}

\AddEq{=D1D2}
{
\def\DFDT{D1 D2 }%
\def\MF{r1:D1->D2 }%
\def\DF{1}%
\def\DT{2}%
\def\VF{1}%
\def\VT{2}%
\def\MapE{hf}%
}

\AddEq{=DD}
{
\def\DFDT{D }%
\def\MF{}%
\def\DF{}%
\def\DT{}%
\def\VF{1}%
\def\VT{2}%
\def\MapE{f}%
}

\AddEq{def universal}
{
\def\Cols{*}%
\def\Rows{*}%
\renewcommand\ARow[2]{\ensuremath{##1(\gi{##2})}}%
\renewcommand\ACol[2]{\ensuremath{##1(\gi{##2})}}%
\renewcommand\EBase[2]{\ensuremath{##1(\gi{##2})}}%
\def\Product{*}%
\def\ProductA{}%
\def\ProductS{}%
\def\ProductType{}%
\def\ProductTypeA{}%
\def\GL{\RCGL}%
\def\Group{GL}%
\def\ProductVal{}%
}

\AddEq{def opRow}
{%
\def\Product{*Row}%
\def\ProductA{}%
\def\ProductS{}%
\def\ProductType{}%
\def\ProductTypeA{}%
\def\GL{\RCGL}%
\def\Group{GL}%
\def\ProductVal{}%
}

\AddEq{def opCol}
{%
\def\Product{*Col}%
\def\ProductA{}%
\def\ProductS{}%
\def\ProductType{}%
\def\ProductTypeA{}%
\def\GL{\RCGL}%
\def\Group{GL}%
\def\ProductVal{}%
}

\AddEq{def rc}
{%
\def\Product{rc}%
\def\ProductA{cr}%
\def\ProductS{rc }%
\def\ProductType{\RC}%
\def\ProductTypeA{\CR}%
\def\GL{\RCGL}%
\def\Group{\RCGL}%
\def\ProductVal{\RCstar}%
\def\ProductAVal{\CRstar}%
\def\PowerVal{\RCPower}%
\def\PowerAVal{\CRPower}%
\def\InverseVal{\RCInverse}%
\def\InverseAVal{\CRInverse}%
\def\DetVal{\RCDet}%
\def\detVal{\RCdet}%
\def\RankVal{\RCRank}%
}

\AddEq{def cr}
{%
\def\Product{cr}%
\def\ProductA{rc}%
\def\ProductS{cr }%
\def\ProductType{\CR}%
\def\ProductTypeA{\RC}%
\def\GL{\CRGL}%
\def\Group{\CRGL}%
\def\ProductVal{\CRstar}%
\def\ProductAVal{\RCstar}%
\def\PowerVal{\CRPower}%
\def\PowerAVal{\RCPower}%
\def\InverseVal{\CRInverse}%
\def\InverseAVal{\RCInverse}%
\def\DetVal{\CRDet}%
\def\detVal{\CRdet}%
\def\RankVal{\CRRank}%
}

\AddEq{def DModule}
{
\def\AArg{d}%
\def\DTransf{g_1}%
\def\DArg{n}%
\def\Base{D}%
\def\BaseExt{(1)}%
\def\CBase{Z}%
\def\BaseRings{of ring of rational integers $Z$ and commutative ring $D$ }%
\def\SidePresentation{}%
\def\Algebra{commutative ring}%
\def\ShortAlgebraWS{ring }%
\def\DivAlgebra{field}%
\def\SetRepresentation{Abelian group }%
\def\VectorSet{module }%
\def\VectorSetC{Module }%
\def\VectorSetNS{module}%
\def\VectorSetNSC{Module}%
\def\VectorsSet{modules }%
\def\VectorsSetNS{modules}%
\def\Division{}%
\def\Module{V}%
\def\ModuleA{W}%
\def\ANumber{d}%
\def\VNumber{a}%
\def\VNumberA{b}%
\def\algebraa{ring }%
\def\Free{free }%
\ifx\UseRussian\Defined%
\def\Same{один и тотже }%
\def\To{}%
\def\TO{вых }%
\def\ToV{вое }%
\def\Algebra{ассоциативная $D$\Hyph алгебра}%
\def\algebraa{$D$\Hyph алгебре }%
\def\algebrac{$D$\Hyph алгебра }%
\def\algebrad{$D$\Hyph алгебру }%
\def\algebraD{$D$\Hyph алгебра }%
\def\BaseRings{кольца целых чисел $Z$ и коммутативного кольца $D$ }%
\def\From{}%
\def\VectorSetRuA{модуль }%
\def\VectorSetRuANS{модуль}%
\def\VectorSetRuB{модуля }%
\def\VectorSetRuBNS{модуля}%
\def\VectorSetRuD{модуле }%
\def\VectorSetRuE{модулем }%
\def\VectorsSetRuA{модули }%
\def\VectorSetRuENS{модулем}%
\def\VectorsSetRuB{модулей }%
\def\VectorsSetRuBNS{модулей}%
\def\AlgebraRu{кольцо }%
\def\AlgebraRuNS{коммутативное кольцо}%
\def\ShortAlgebraRuBWS{кольца }%
\def\ShortAlgebraRuCWS{кольцом }%
\def\ShortAlgebraRuDWS{кольце }%
\def\ShortAlgebraRuEWS{кольцо }%
\def\SidePresentationRuWS{}%
\def\SidePresentationRuBWS{}%
\def\SetRepresentationRuA{абелевая группа }%
\def\SetRepresentationRu{абелевой группе }%
\def\DivisionRu{ }%
\def\DivisionRuNS{}%
\def\DivAlgebraRu{поле}%
\def\SideModuleC{}%
\def\FreeRuB{свободного }%
\fi
}

\AddEq{def DAlgebra}
{
\def\AArg{d}%
\def\DTransf{g_1}%
\def\DArg{n}%
\def\Base{D}%
\def\BaseExt{(1)}%
\def\CBase{Z}%
\def\BaseRings{of ring of rational integers $Z$ and commutative ring $D$ }%
\def\SidePresentation{}%
\def\Algebra{commutative ring}%
\def\ShortAlgebraWS{ring }%
\def\DivAlgebra{field}%
\def\SetRepresentation{Abelian group }%
\def\VectorSet{algebra }%
\def\VectorSetC{Algebra }%
\def\VectorSetNS{algebra}%
\def\VectorSetNSC{Algebra}%
\def\VectorsSet{algebras }%
\def\VectorsSetNS{algebras}%
\def\Division{}%
\def\Module{V}%
\def\ModuleA{W}%
\def\ANumber{d}%
\def\VNumber{a}%
\def\VNumberA{b}%
\def\algebraa{ring }%
\def\Free{free }%
\ifx\UseRussian\Defined%
\def\Same{одна и таже }%
\def\To{}%
\def\TO{вых }%
\def\ToV{вое }%
\def\Algebra{ассоциативная $D$\Hyph алгебра}%
\def\algebraa{$D$\Hyph алгебре }%
\def\algebrac{$D$\Hyph алгебра }%
\def\algebrad{$D$\Hyph алгебру }%
\def\algebraD{$D$\Hyph алгебра }%
\def\BaseRings{кольца целых чисел $Z$ и коммутативного кольца $D$ }%
\def\From{}%
\def\VectorSetRuA{алгебра }%
\def\VectorSetRuANS{алгебра}%
\def\VectorSetRuB{алгебры }%
\def\VectorSetRuBNS{алгебры}%
\def\VectorSetRuD{алгебре }%
\def\VectorSetRuE{алгеброй }%
\def\VectorsSetRuA{модули }%
\def\VectorSetRuENS{модулем}%
\def\VectorsSetRuB{модулей }%
\def\VectorsSetRuBNS{модулей}%
\def\AlgebraRu{кольцо }%
\def\AlgebraRuNS{коммутативное кольцо}%
\def\ShortAlgebraRuBWS{кольца }%
\def\ShortAlgebraRuCWS{кольцом }%
\def\ShortAlgebraRuDWS{кольце }%
\def\ShortAlgebraRuEWS{кольцо }%
\def\SidePresentationRuWS{}%
\def\SidePresentationRuBWS{}%
\def\SetRepresentationRuA{абелевая группа }%
\def\SetRepresentationRu{абелевой группе }%
\def\DivisionRu{ }%
\def\DivisionRuNS{}%
\def\DivAlgebraRu{поле}%
\def\SideModuleC{}%
\def\FreeRuB{свободного }%
\fi
}

\AddEq{def DVector}
{
\def\AArg{d}%
\def\DTransf{g_1}%
\def\DArg{n}%
\def\Base{D}%
\def\BaseExt{}%
\def\CBase{}%
\def\BaseRings{of ring of rational integers $Z$ and commutative ring $D$ }
\def\SidePresentation{}
\def\Algebra{field}%
\def\ShortAlgebraWS{field }%
\def\DivAlgebra{field}%
\def\SetRepresentation{Abelian group }%
\def\VectorSet{vector space }%
\def\VectorSetC{Vector Space }%
\def\VectorSetNS{vector space}%
\def\VectorSetNSC{Vector Space}%
\def\VectorsSet{vector spaces }%
\def\VectorsSetNS{vector spaces}%
\def\Division{}%
\def\Module{V}%
\def\ModuleA{W}%
\def\ANumber{d}%
\def\VNumber{v}%
\def\VNumberA{w}%
\def\algebraa{field }%
\def\Free{}
\ifx\UseRussian\Defined%
\def\Same{одно и тоже }%
\def\To{}%
\def\TO{}%
\def\ToV{}%
\def\What{}%
\def\WhatA{ый }%
\def\WhatO{}%
\def\Algebra{кольцо}%
\def\algebraa{кольце }%
\def\algebrac{кольцо }%
\def\algebrad{кольцо }%
\def\algebraD{Кольцо }%
\def\Algebrab{модуль}%
\def\BaseRings{кольца целых чисел $Z$ и коммутативного кольца $D$ }
\def\From{}%
\def\VectorSetRuA{векторное пространство }%
\def\VectorSetRuANS{векторное пространство}%
\def\VectorSetRuB{векторного пространства }%
\def\VectorSetRuBNS{векторного пространства}%
\def\VectorSetRuD{векторном пространстве }%
\def\VectorSetRuE{векторным пространством }%
\def\VectorSetRuENS{векторным пространством}%
\def\VectorsSetRuA{векторные пространства }%
\def\VectorsSetRuB{векторных пространств }%
\def\VectorsSetRuBNS{векторных пространств}%
\def\AlgebraRu{поле }%
\def\AlgebraRuNS{поле}%
\def\ShortAlgebraRuBWS{поля }%
\def\ShortAlgebraRuCWS{полем }%
\def\ShortAlgebraRuDWS{поле }%
\def\ShortAlgebraRuEWS{поле }%
\def\SidePresentationRuWS{}%
\def\SidePresentationRuBWS{}%
\def\SetRepresentationRuA{абелевая группа }%
\def\SetRepresentationRu{абелевой группе }%
\def\DivisionRu{ }%
\def\DivisionRuNS{}%
\def\DivAlgebraRu{поле}%
\def\SideModuleC{}%
\def\FreeRuB{}
\fi
}

\AddEq{def AModule}
{
\def\AArg{a}%
\def\DTransf{g_{1,4}}%
\def\DArg{d}%
\def\Base{A}%
\def\BaseExt{(1)}%
\def\CBase{D}%
\def\BaseRings{of commutative ring $D$ and $D$\Hyph algebra $A$ }
\def\SidePresentation{\Hyph side }
\def\Algebra{associative $D$\Hyph algebra}%
\def\ShortAlgebraWS{$D$\Hyph algebra }%
\def\DivAlgebra{associative division $D$\Hyph algebra}%
\def\SetRepresentation{$D$\Hyph module }%
\def\VectorSet{module }%
\def\VectorSetC{Module }%
\def\VectorSetNS{module}%
\def\VectorSetNSC{Module}%
\def\VectorsSet{modules }%
\def\VectorsSetNS{modules}%
\def\Division{}%
\def\Module{A}%
\def\ModuleA{B}%
\def\ANumber{a}%
\def\VNumber{v}%
\def\VNumberA{w}%
\def\algebraa{$D$\Hyph algebra }%
\def\Free{free }
\ifx\UseRussian\Defined%
\def\Same{один и тотже }%
\def\To{вый }%
\def\TO{вых }%
\def\ToV{вое }%
\def\Algebra{ассоциативная $D$\Hyph алгебра}%
\def\algebraa{$D$\Hyph алгебре }%
\def\algebrac{$D$\Hyph алгебра }%
\def\algebrad{$D$\Hyph алгебру }%
\def\algebraD{$D$\Hyph алгебра }%
\def\BaseRings{коммутативного кольца $D$ и $D$\Hyph алгебры $A$ }
\def\From{вого }%
\def\VectorSetRuA{модуль }%
\def\VectorSetRuANS{модуль}%
\def\VectorSetRuB{модуля }%
\def\VectorSetRuBNS{модуля}%
\def\VectorSetRuD{модуле }%
\def\VectorSetRuE{модулем }%
\def\VectorSetRuENS{модулем}%
\def\VectorsSetRuA{модули }%
\def\VectorsSetRuB{модулей }%
\def\VectorsSetRuBNS{модулей}%
\def\AlgebraRu{ассоциативная $D$\Hyph алгебра }%
\def\AlgebraRuNS{ассоциативная $D$\Hyph алгебра}%
\def\ShortAlgebraRuBWS{$D$\Hyph алгебры }%
\def\ShortAlgebraRuCWS{$D$\Hyph алгеброй }%
\def\ShortAlgebraRuDWS{$D$\Hyph алгебре }%
\def\ShortAlgebraRuEWS{$D$\Hyph алгебру }%
\def\SidePresentationRuWS{востороннее }%
\def\SidePresentationRuBWS{восторонним }%
\def\SetRepresentationRuA{$D$\Hyph модуль }%
\def\SetRepresentationRu{$D$\Hyph модуле }%
\def\DivisionRu{ }%
\def\DivisionRuNS{}%
\def\DivAlgebraRu{ассоциативная $D$\Hyph алгебра с делением}%
\def\SideModuleC{вым }%
\def\FreeRuB{свободного }
\fi
}

\AddEq{def AVector}
{
\def\AArg{a}%
\def\DTransf{g_{1,4}}%
\def\DArg{d}%
\def\Base{A}%
\def\BaseExt{}%
\def\CBase{D}%
\def\BaseRings{of commutative ring $D$ and $D$\Hyph algebra $A$ }
\def\SidePresentation{\Hyph side }
\def\Algebra{associative division $D$\Hyph algebra}%
\def\ShortAlgebraWS{$D$\Hyph algebra }%
\def\DivAlgebra{associative division $D$\Hyph algebra}%
\def\SetRepresentation{$D$\Hyph vector space }%
\def\VectorSet{vector space }%
\def\VectorSetC{Vector Space }%
\def\VectorSetNS{vector space}%
\def\VectorSetNSC{Vector Space}%
\def\VectorsSet{vector spaces }%
\def\VectorsSetNS{vector spaces}%
\def\Division{division }%
\def\Module{V}%
\def\ModuleA{W}%
\def\ANumber{a}%
\def\VNumber{v}%
\def\VNumberA{w}%
\def\algebraa{$D$\Hyph algebra }%
\def\Free{}
\ifx\UseRussian\Defined%
\def\Same{одно и тоже }%
\def\To{вое }%
\def\TO{}%
\def\ToV{}%
\def\What{}%
\def\WhatA{вый }%
\def\WhatO{вое }%
\def\Algebra{кольцо}%
\def\algebraa{кольце }%
\def\algebrac{кольцо }%
\def\algebrad{кольцо }%
\def\algebraD{Кольцо }%
\def\Algebrab{модуль}%
\def\BaseRings{коммутативного кольца $D$ и $D$\Hyph алгебры $A$ }
\def\From{вого }%
\def\VectorSetRuA{векторное пространство }%
\def\VectorSetRuANS{векторное пространство}%
\def\VectorSetRuD{векторном пространстве }%
\def\VectorSetRuB{векторного пространства }%
\def\VectorSetRuBNS{векторного пространства}%
\def\VectorSetRuE{векторным пространством }%
\def\VectorSetRuENS{векторным пространством}%
\def\VectorsSetRuA{векторные пространства }%
\def\VectorsSetRuB{векторных пространств }%
\def\VectorsSetRuBNS{векторных пространств}%
\def\AlgebraRu{ассоциативная $D$\Hyph алгебра }%
\def\AlgebraRuNS{ассоциативная $D$\Hyph алгебра}%
\def\ShortAlgebraRuBWS{$D$\Hyph алгебры }%
\def\ShortAlgebraRuCWS{$D$\Hyph алгеброй }%
\def\ShortAlgebraRuDWS{$D$\Hyph алгебре }%
\def\ShortAlgebraRuEWS{$D$\Hyph алгебру }%
\def\SidePresentationRuWS{востороннее }%
\def\SidePresentationRuBWS{восторонним }%
\def\SetRepresentationRuA{$D$\Hyph векторное пространство }%
\def\SetRepresentationRu{$D$\Hyph векторном пространстве }%
\def\DivisionRu{ с делением }%
\def\DivisionRuNS{ с делением}%
\def\DivAlgebraRu{ассоциативная $D$\Hyph алгебра с делением}%
\def\SideModuleC{вым }%
\def\FreeRuB{}
\fi
}

\AddEq{def left}
{%
\def\SideNS{left}%
\def\SideWS{left }%
\def\SideWSC{Left }%
\def\DefRow{rc-rows}%
\def\DefCol{cr-cols}%
\def\HSide{\Hyph side }%
\def\ATransf{g_{3,4}}%
\def\OtherSideWS{right }%
\def\OtherSideNS{right}%
\def\RefDistributive##1{(b1+b2).##1.a=}%
\renewcommand\Multiply[2]{##1##2}%
\ifx\UseRussian\Defined%
\def\HSide{востороннее }%
\def\SideRu{ле}%
\def\OtherSideRu{пра}%
\def\SideRuC{Ле}%
\def\What{вым }%
\def\WhatTo{вом }%
\def\Which{вого }%
\fi%
}

\AddEq{def right}
{%
\def\SideNS{right}%
\def\SideWS{right }%
\def\SideWSC{Right }%
\def\DefRow{cr-rows}%
\def\DefCol{rc-cols}%
\def\HSide{\Hyph side }%
\def\ATransf{g_{3,4}}%
\def\OtherSideWS{left }%
\def\OtherSideNS{left}%
\def\RefDistributive##1{a.##1.(b1+b2)=}%
\renewcommand\Multiply[2]{##2##1}
\ifx\UseRussian\Defined%
\def\HSide{востороннее }%
\def\SideRu{пра}%
\def\OtherSideRu{ле}%
\def\SideRuC{Пра}%
\def\What{вым }%
\def\WhatTo{вом }%
\def\Which{вого }%
\fi%
}

\AddEq{def ab=ba}
{
\def\SideNS{}%
\def\SideWS{}%
\def\DefRow{-rows}%
\def\DefCol{-cols}%
\def\HSide{}%
\def\ATransf{g_2}%
\ifx\UseRussian\Defined%
\def\SideRu{}%
\def\SideRuC{}%
\def\What{}%
\def\WhatTo{}%
\def\Which{}%
\fi%
}

\AddEq{=Omega}
{
\def\Base{\Omega}%
\def\Algebrab{group}%
\def\Algebrac{Omega group}%
\def\Module{A}%
\ifx\UseRussian\Defined%
\def\WhatA{ая }%
\def\Algebrab{группа}%
\fi%
}

\DefEq%
{%
\def\Pt{.}%
}%
{=.}%

\DefEq%
{%
\def\Pt{;}%
}%
{=.c}%

\DefEq%
{%
\def\Pt{,}%
}%
{=c}%

\DefEq%
{%
\def\Pt{}%
}%
{=z}%

\AddEq{D= F= A=}
{
\def\PD{}%
\def\PF{}%
\def\PA{}%
}

\DefEq
{
\def\PX{X}
}
{X}

\AddEq{f=f}
{%
\def\Pf{f}%
}%

\AddEq{f=g}
{%
\def\Pf{g}%
}%

\DefEq%
{%
\def\Pf{m}%
}%
{f=m}%

\DefEq
{
\def\Pf{I}%
}
{f=I}

\DefEq
{
\def\pD{D}%
\def\pA{A}%
\def\pB{A}%
}
{A=A}

\DefEq%
{%
\def\pD{D}%
\def\pA{A}%
\def\pB{B}%
}%
{A=AB}

\DefEq
{
\def\pD{D}%
\def\pA{B}%
\def\pB{C}%
}
{A=BC}

\DefEq%
{%
\def\pD{R}%
\def\pA{C}%
\def\pB{C}%
}%
{A=C}%

\DefEq%
{%
\def\pD{C}%
\def\pA{C}%
\def\pB{C}%
}%
{A=CC}%

\DefEq
{
\def\pD{D}%
\def\pA{A_1}%
\def\pB{A_2}%
}
{A=12}

\DefEq
{
\def\pD{D}%
\def\pA{A_2}%
\def\pB{A_2}%
}
{A=2}

\DefEq
{
\def\pD{D}%
\def\pA{A_1\rightarrow A_2}%
\def\pB{A_3}%
}
{A=123}

\DefEq
{
\def\pD{D}%
\def\pA{A_1}%
\def\pB{A_3}%
}
{A=13}

\DefEq
{
\def\pD{D}%
\def\pA{A}%
\def\pB{C}%
}
{A=AC}

\DefEq
{%
\def\Pn{0}%
}%
{n=0}

\DefEq
{%
\def\Pn{n}%
}%
{n=n}

\DefEq
{%
\def\Pn{k}%
}%
{n=k}

\DefEq
{%
\def\Pn{1}%
\def\Pp{p}%
}%
{n=1}

\DefEq
{%
\def\Pn{2}%
\def\Pp{q}%
}%
{n=2}

\DefEq
{%
\def\Pn{3}%
\def\Pp{r}%
}%
{n=3}

\DefEq
{%
\def\tM{1}%
}%
{m=1}

\DefEq
{%
\def\tM{2}%
}%
{m=2}

%% file: Prolog.Eq.tex

\def\iI{\ensuremath{i\in I}}%
\def\iIg{\ensuremath{\gii\in\giI}}%
\newcommand\jJg[2]{\ensuremath{\gi{#1}\in\gi{#2}}}%
\def\Times{A_1\times...\times A_n}%
\newcommand\LAA[2]{\ensuremath{\mathcal L(#1;#2\rightarrow #2)}}%
\newcommand\LAB[3]{\ensuremath{\mathcal L(#1;#2\rightarrow #3)}}%
\newcommand\Kn[2][k]{$#1=#2$, ..., $n$}%
\newcommand\Kb[2][k]{$(#1)=(#2)$, ..., $(n)$}%
\newcommand\Ki[2][i]{\ensuremath{\gi{#1}=\gi 1}, ..., \ensuremath{\gi{#2}}}%
\def\ATwo{A_2\otimes A_2}%
\newcommand\AxA[1]{\ensuremath{#1\times #1}}%
\newcommand\AoxA[1]{\ensuremath{#1\otimes #1}}%
\newcommand\BoxB[1]{\AoxA{#1}\Hyph}%
\newcommand\AoxB[2]{\ensuremath{#1\otimes #2}}%

\newcommand\FC[2]{f^{\gi{#1#2}}}%
\newcommand\TensorBasis[1]%
{e_{1\cdot\gi{#1_1}}\otimes...\otimes e_{n\cdot\gi{#1_n}}}%
\newcommand\re{\mathrm{Re}\,}%
\newcommand\im{\mathrm{Im}\,}%

\ePrints{309618526,5284-0163,1801.01628,322019412,323236904,330532713,323966352}
\Items{8428-0408,2207.06506,1908.04418,6860-2955,367250917,1506.00061}%
\Items{2021.01.06,1601.03259,2022.01.05,2021.11.26,2112.00613,348311191,2020.06.01}%
\Items{357606033,MAlgebra2,2204.06320,2306.00880,2022.11.26,2023.01.07}%
\ifx\Semafor\ValueOn%
\input{Prolog.Cite}\fi%

\ePrints{MAlgebra2,1908.04418,6860-2955,1601.03259,4975-6381,5410-9916,1305.4547,1502.04063,1310.5591}%
\Items{309618526,CACAA.06.121,5059-9176,323236904,322019412,5284-0163,1801.01628,9835-2163,9856-6693,0767-8264}%
\Items{8428-0408,2207.06506,CACAA.07,323966352,0921-2370,5148-4632,1506.00061,7287-9339,330532713}%
\Items{2020.06.01,2021.01.06,2204.06320,2022.01.05,348311191,2306.00880,2022.11.26,2023.01.07}%
\Items{2023.10.28}
\ifx\Semafor\ValueOn%
\input{\FilePrefix Stmt.Representation.\TheLanguage}\fi%

\ePrints{MAlgebra2,309618526,CACAA.06.121,1601.03259,4975-6381,5410-9916,1502.04063,5114-6019}%
\Items{323236904,1908.04418,6860-2955,322019412,5284-0163,1801.01628,9835-2163,9856-6693,0767-8264,MTensor}%
\Items{323966352,0921-2370,CACAA.07,1506.00061,7287-9339,5148-4632}%
\Items{8428-0408,2207.06506,0803.3276,7100-9423,2020.06.01,1211.6965,MAlgebra4,2021.01.06}%
\Items{2021.01.06,2021.11.26,2112.00613,2022.01.05,2204.06320,2306.00880,2022.11.26}%
\Items{2023.01.07,367250917,2023.10.28}%
\ifx\Semafor\ValueOn%
\input{\FilePrefix Stmt.Module.\TheLanguage}\fi%

\ePrints{1801.01628,2020.06.01,339295394,348311191,5284-0163,1506.00061}%
\Items{8428-0408,2207.06506,2204.06320,2022.01.05,1908.04418,6860-2955}%
\Items{MAlgebra2,5148-4632,2306.00880,2022.11.26,2023.01.07,357606033,367250917,2023.10.28}%
\ifx\Semafor\ValueOn%
\input{\FilePrefix Vector.Space.2020.Stmt.\TheLanguage}\fi%

\ePrints{1502.04063,323236904,0803.3276,1801.01628}%
\Items{0767-8264,2020.06.01,0412.391,4827-2437,7100-9423,5284-0163}%
\ifx\Semafor\ValueOn%
\input{\FilePrefix D.Module.Stmt.\TheLanguage}%
\fi%

\ePrints{1801.01628,5284-0163}%
\Items{2204.06320,2022.01.05,348311191,357606033,2020.06.01}%
\Items{2023.01.07}%
\ifx\Semafor\ValueOn%
\input{\FilePrefix DiffUr.3.Stmt.\TheLanguage}%
\fi%

\ePrints{8428-0408,2207.06506,367250917}%
\ifx\Semafor\ValueOn%
\input{DiffUr.3.Stmt.Eq}%
\fi%

\ePrints{1908.04418,6860-2955,1801.01628,2020.06.01,5284-0163,0767-8264}%
\Items{8428-0408,2207.06506,2204.06320,2022.01.05}%
\ifx\Semafor\ValueOn%
\input{\FilePrefix Biring.Stmt.\TheLanguage}\fi%

\ePrints{1908.04418,6860-2955,2207.06506,8428-0408,1305.4547,5284-0163,1801.01628}%
\Items{6860-2955,5148-4632,MAlgebra2,5410-9916,1601.03259,309618526,CACAA.06.121,4975-6381,322019412}%
\ifx\Semafor\ValueOn%
\input{\FilePrefix Stmt.Omega.\TheLanguage}\fi%

\ePrints{309618526,CACAA.06.121,4975-6381,1601.03259,323236904}%
\Items{5284-0163,1801.01628,9835-2163,0767-8264,322019412,1211.6965,MAlgebra4}%
\Items{8428-0408,2207.06506,9856-6693,CACAA.07,330532713,2020.06.01,2021.01.06,348311191}%
\ifx\Semafor\ValueOn%
\input{\FilePrefix Stmt.Derivative.\TheLanguage}%
\fi%

\ePrints{323966352,0921-2370,1908.04418,6860-2955,0803.3276,7100-9423,FAlgebra,0405.027}%
\ifx\Semafor\ValueOn%
\input{\FilePrefix Stmt.Gravity.\TheLanguage}%
\fi%

\ePrints{324329551}%
\ifx\Semafor\ValueOn%
\input{Stmt.Gravity.Eq}%
\fi%

\ePrints{MAlgebra3,1601.03259,4975-6381,322019412,5284-0163,1801.01628}%
\ifx\Semafor\ValueOn%
\input{Algebra.Complex.2016.Eq}%
\fi%

\ePrints{309618526,323966352,0921-2370,CACAA.06.121,MAlgebra3,CACAA.07,322019412}%
\Items{1506.00061,7287-9339,5284-0163,1801.01628,9835-2163,0767-8264,330532713,2021.01.06}%
\Items{339295394,348311191}%
\ifx\Semafor\ValueOn%
\input{Algebra.Quaternion.2016.Eq}%
\fi%

\ePrints{309618526,323966352,0921-2370,MAlgebra3}%
\ifx\Semafor\ValueOn%
\input{Algebra.Octonion.2016.Eq}%
\fi%

\ePrints{330532713,339295394,348311191,357606033}%
\ifx\Semafor\ValueOn%
\input{Stmt.Module.Eq}\fi%

\ePrints{323966352,0921-2370,339295394,2204.06320,2022.01.05,367250917}%
\ifx\Semafor\ValueOn%
\input{Stmt.Derivative.Eq}%
\fi%

\ePrints{339295394,357606033,367250917}%
\ifx\Semafor\ValueOn%
\input{Stmt.Representation.Eq}
\fi%

\ePrints{2021.11.26,2112.00613,8428-0408,2207.06506,357606033}%
\ifx\Semafor\ValueOn%
\input{\FilePrefix Stmt.2021.11.26.\TheLanguage}%
\fi%

\ePrints{2022.01.05,2306.00880,2022.11.26,2204.06320,2023.01.07}%
\ifx\Semafor\ValueOn%
\input{\FilePrefix Stmt.2022.01.05.\TheLanguage}%
\fi%

\input{Prolog.Ref}
\AddEq{L D->A=A}
{
$\mathcal L(D;D\rightarrow A)=A$
\qed
}

\AddEq{|x'-x|<delta}
{
\[
\|x'-x\|_1<\delta
\]
}

\AddEquation{fd=da}
{
f\circ d=da
}

\AddEq{fd=df1}
{
\[f\circ d=d(f\circ 1)\]
}

\AddEq{ab=ba}
{
ab=ba
}

\AddEq[1]{a=b}
{
$a_1=b_1$#1
}

\DefEq
{
$\Basis e_i$
}
{tensor product of algebras, basis i}

\AddEq[1]{set a1n}
{
$\{a_1,...,a_n\}$#1
}

\AddEquation{x2+ax+xa=}
{
x^2+(1\otimes a+a\otimes 1)\circ x
=x\left(\frac 12 x+a\right)+\left(\frac 12x+a\right)x
}

\AddEq[1]{Quaternion a}
{
#1=\aU{#1}0+\aU{#1}1i+\aU{#1}2j+\aU{#1}3k
}

\DefEq
{
\symb{a^{\gi{i_1...i_n}}}{standard component of tensor}{}
}
{standard component of tensor}

\AddEquation{module L(A;A) is algebra}
{
\circ:(g,f)\in\mathcal L(D;A\rightarrow A)\times\mathcal L(D;A\rightarrow A)
\rightarrow g\circ f\in\mathcal L(D;A\rightarrow A)
}

\AddEq [3]{f:A->B}
{
\[#1:#2\rightarrow #3\]
}

\DefEq
{
\[
\xymatrix{
&A\ar[rd]^g
\\
A\ar[ru]^f\ar[rr]^{g\circ f}&&A
}
\]
}
{product of linear map, algebra 1}

\DefEquation
{
a=
\ShowSymbol{standard component of tensor}{}
\TensorBasis i
}
{tensor canonical representation, algebra}

\AddEq{k,1n}
{
$k$, $k=1$, ..., $n$,
}

\AddEq{n,1.}
{
$n$, $n=1$, ...,
}

\newcommand{\Tensor}[1]{#1_1\otimes...\otimes #1_n}

\AddEq[1]{set Ai}
{
$#1_i$, \iI,
}

\AddEq{coproduct in category diagram}
{
\[
\xymatrix{
P\ar[d]_h&B_i\ar[l]_{f_i}\ar[dl]^{g_i}&h\circ f_i=g_i
\\
R
}
\]
}

\AddEq[2]{A=Ai iI}
{
\(#1=(#1_i,\iI)\)#2
}

\DefEq
{
\begin{align*}
\re&:A\rightarrow A
\\
\im&:A\rightarrow A
\end{align*}
}
{Re Im A->A}

\DefEq
{
\symb{\re d}{scalar of algebra}{}
\symb{\im d}{vector of algebra}{}
}
{Re Im A->F 0}

\AddEquation{Re Im A->F 1}
{
\begin{matrix}
\ShowSymbol{scalar of algebra}{}=d^{\gi 0}
&\ShowSymbol{vector of algebra}{}=d-d^{\gi 0}
&d\in D
&d=d^{\gii}e_{\gii}
\end{matrix}
}

\DefEq
{
$\ShowSymbol{scalar of algebra}{}$
}
{Re Im A->F 2}

\DefEq
{
$\ShowSymbol{vector of algebra}{}$
}
{Re Im A->F 3}

\DefEq
{
\[
F=\{d\in A:\re d=d\}
\]
}
{Re Im A->F 4}

\DefEq
{
\symb{\re A}{scalar algebra of algebra}1
}
{scalar algebra of algebra}

\DefEq
{
\symb{\im A}{vector module of algebra}{}
\begin{equation}
\EqLabel{vector module of algebra}
\ShowSymbol{vector module of algebra}{}=\{d\in A:\re d=0\}
\end{equation}
}
{vector module of algebra}

\DefEquation
{
C_{\gi{0k}}^{\gil}=C^{\gil}_{\gi{k0}}=\delta^{\gik}_{\gil}
}
{structural constant re algebra, 1}

\DefEquation
{
d=\re d+\im d
}
{d Re Im}

\DefEq
{
\symb{d^*}{conjugation in algebra}{}
}
{conjugation in algebra}

\DefEquation
{
\ShowSymbol{conjugation in algebra}{}=\re d-\im d
}
{conjugation in algebra, 0}

\DefEquation
{
(cd)^*=d^*\,c^*
}
{conjugation in algebra, 1}

\DefEquation
{
\begin{matrix}
C^{\gi 0}_{\gi{kl}}=C^{\gi 0}_{\gi{lk}}
&
C^{\gi p}_{\gi{kl}}=-C^{\gi p}_{\gi{lk}}
\end{matrix}
}
{conjugation in algebra, 5}

\AddEq{p=p circ x}
{
$p(x)=p_1\circ x+p_0$
}

\AddEq{x in C=>fx=gx}
{
$x\in C\Rightarrow f(x)=g(x)$.
}

\AddEq[3]{A in B}
{
$#1\subseteq #2$#3
}

\AddEq{(a+x)2=a2 1}
{
x^2+(1\otimes a+a\otimes 1)\circ x=0
}

\DefEquation
{
p(x)=p_0+p_1\circ x
}
{p=p+p circ x}

\DefEq
{
\[
r(x)=r_0+r_1\circ x+...+r_k\circ x^k
\]
}
{r power k}

\DefEquation
{
\begin{split}
r(x)&=r_0
-((r_{1\cdot 0\cdot s}\otimes r_{1\cdot 1\cdot s})\circ
p^{-1}_1)\circ p_0
\\
&-(((r_{2\cdot 0\cdot s}\circ x)\otimes r_{2\cdot 1\cdot s})\circ
p^{-1}_1)\circ p_0
\\
&-...-(((r_{k\cdot 0\cdot s}\circ x^{k-1})
\otimes r_{k\cdot 1\cdot s})\circ
p^{-1}_1)\circ p_0
\\
&+((r_{1\cdot 0\cdot s}\otimes r_{1\cdot 1\cdot s})\circ
p^{-1}_1)\circ p(x)
\\
&+(((r_{2\cdot 0\cdot s}\circ x)\otimes r_{2\cdot 1\cdot s})\circ
p^{-1}_1)\circ p(x)
\\
&+...+(((r_{k\cdot 0\cdot s}\circ x^{k-1})
\otimes r_{k\cdot 1\cdot s})\circ
p^{-1}_1)\circ p(x)
\\
&=r_0
-((r_{1\cdot 0\cdot s}\otimes r_{1\cdot 1\cdot s}
+(r_{2\cdot 0\cdot s}\circ x)\otimes r_{2\cdot 1\cdot s}
\\
&+...+(r_{k\cdot 0\cdot s}\circ x^{k-1})
\otimes r_{k\cdot 1\cdot s})\circ
p^{-1}_1)\circ p_0
\\
&+((r_{1\cdot 0\cdot s}\otimes r_{1\cdot 1\cdot s}
+(r_{2\cdot 0\cdot s}\circ x)\otimes r_{2\cdot 1\cdot s}
\\
&+...+(r_{k\cdot 0\cdot s}\circ x^{k-1})
\otimes r_{k\cdot 1\cdot s})\circ
p^{-1}_1)\circ p(x)
\end{split}
}
{r=+q circ p 1}

\DefEq
{
\[
ax-xa=1
\]
}
{ax-xa=1}

\DefEquation
{
r(x)=q_{i\cdot 0}(x)p(x)q_{i\cdot 1}(x)
=(q_{i\cdot 0}(x)\otimes q_{i\cdot 1}(x))\circ p(x)
}
{r=qpq(x)}

\DefEq
{
h=f_1...f_n\omega
}
{h=f1...fn omega}

\DefEq
{
\(f_X(g_{1\cdot n})(g_{2\cdot n})\)
}
{f(g1n)(g2n)}

\DefEq
{
\[b_i=f_i(a_i)\]
}
{b=f(a)i}

\AddEq [2]{f->g}
{
\[#1\rightarrow #2\]
}

\AddEq [3]{A->B}
{
$#1\rightarrow #2$#3
}

\DefEq
{
\[h:S_1\rightarrow S_2\]
}
{f->g 1}

\AddEq{polylinear maps category, diagram}
{
\[
\xymatrix{
&S_1\ar[dd]^h
\\
\Times
\ar[ru]^f\ar[rd]_g
\\
&S_2
}
\]
}

\AddEq [3]{map r,R}
{
$(#1,#2)$#3
}

\DefEq
{
\symb{A^{\otimes n}}{tensor power of algebra}{}
\[
\begin{matrix}
\ShowSymbol{tensor power of algebra}{}
=\Tensor A
&
A_1=...=A_n=A
\end{matrix}
\]
}
{tensor power of algebra}

\AddEq{basis of complex field}
{
\[
\begin{array}{rr}
e_{\gi 0}=1&e_{\giA}=i
\end{array}
\]
}

\AddEquation{product of complex field}
{
e_{\giA}^2=-e_{\gi 0}
}

\DefEq
{
\[f:H\rightarrow H\ \ \ y^{\gii}=f^{\gii}_{\gij}x^{\gij}\]
}
{f:H->H ij}

\DefEq
{
\[HE=\{aE:a\in H\}\]
}
{HE set}

\AddEquation{structural constants of complex field}
{
\begin{array}{r@{}rr@{}r}
\aUD C0{00}=&1&\aUD C1{01}=&1
\\
\VirtVar
\aUD C1{10}=&1&\aUD C0{11}=&-1
\end{array}
}

\DefEq
{
\[
x=x^{\gi 0}+x^{\giA}i\rightarrow
f=f^{\gi 0}(x^{\gi 0},x^{\giA})+f^{\giA}(x^{\gi 0},x^{\giA})i
\]
}
{map of complex variable}

\AddEq [1]{Cauchy Riemann equations, complex field, 1}
{
\begin{split}
\frac{\partial #1^{\giA}}{\partial x^{\gi 0}}
&=-\frac{\partial #1^{\gi 0}}{\partial x^{\giA}}
\\[5pt]
\frac{\partial #1^{\gi 0}}{\partial x^{\gi 0}}
&=\hphantom{-\,}\frac{\partial #1^{\giA}}{\partial x^{\giA}}
\end{split}
}

\AddEq [1]{Cauchy Riemann equations, complex field, 2, 1}
{
\frac{\partial #1^{\gi 0}}{\partial x^{\gi 0}}
+i\frac{\partial #1^{\giA}}{\partial x^{\gi 0}}
+i\left(\frac{\partial #1^{\gi 0}}{\partial x^{\giA}}
+i\frac{\partial #1^{\giA}}{\partial x^{\giA}}\right)=
\frac{\partial #1^{\gi 0}}{\partial x^{\gi 0}}
+i\frac{\partial #1^{\giA}}{\partial x^{\gi 0}}
+i\frac{\partial #1^{\gi 0}}{\partial x^{\giA}}
-\frac{\partial #1^{\giA}}{\partial x^{\giA}}=
0
}

\DefEq
{
\[
\begin{matrix}
a\circ E:H\rightarrow H&
a=a^{\gi 0}+a^{\gi 1}i+a^{\gi 2}j+a^{\gi 3}k\in H
\end{matrix}
\]
}
{aE, quaternion}

\AddEq{df=...dx 03}
{
\begin{split}
\frac{d f}{d x}=&-\frac 12
\left(
\frac{\partial f}{\partial x^{\gi 0}}
+\frac{\partial f}{\partial x^{\gi 1}}i
+\frac{\partial f}{\partial x^{\gi 2}}j
+\frac{\partial f}{\partial x^{\gi 3}}k
\right)\bullet E
\\&
+\frac 12\left(
\frac{\partial f}{\partial x^{\gi 0}}
+\frac{\partial f}{\partial x^{\gi 1}}i
\right)\bullet\aU I1
+\frac 12\left(
\frac{\partial f}{\partial x^{\gi 0}}
+\frac{\partial f}{\partial x^{\gi 2}}j
\right)\bullet\aU I2
\\&
+\frac 12\left(
\frac{\partial f}{\partial x^{\gi 0}}
+\frac{\partial f}{\partial x^{\gi 3}}k
\right)\bullet\aU I3
\end{split}
}

\AddEq{tensor product of algebras}
{
\symb{\Tensor A}{tensor product}1
}

\DefEq
{
\[
(f+g)(a)=f(a)+g(a)
\]
}
{(f+g)(a)=}

\AddEquation{sum of transformations of Abelian group, 1}
{
f(a+b)(x)=f(a)(x)+f(b)(x)
}

\ePrints{1502.04063,5114-6019,1211.6965,MAlgebra4}
\ifx\Semafor\ValueOff
\DefEquation
{
f(a)(m)=f(b)(m)
}
{representation of ring, 1}

\AddEquation{representation of ring, 2}
{
f(a-b)(m)=0
}
\fi

\AddEq{polynomials p,r}
{
$p\circ F_{[1]}\circ x^n$, $r\circ F_{[2]}\circ x^m$
}

\DefEq
{
\[
(p\circ F_{[1]}\circ x^n)(r\circ F_{[2]}\circ x^m)
=(p\underline{\otimes}r)\circ
(F_{[1](1)},...,F_{[1](n)},F_{[2](1)},...,F_{[2](m)})
\circ x^{n+m}
\]
}
{pr=p circ r}

\DefEquation
{
p(x)=
(a_{k\cdot 0}\circ x^{k-1})
((1\otimes a_{k\cdot 1})\circ x)
}
{p(x)=a k-1 k 1}

\DefEquation
{
p(x)=
((a_{k\cdot 0}\circ x^{k-1})
\otimes a_{k\cdot 1})\circ x
}
{p(x)=a k-1 k 2}

\DefEq
{
$A_{ij}$, $i=1$, ..., $n$, $j=1$, ..., $m$,
}
{Aij}

\DefEq
{
$\Omega_{ij}$\Hyph%
}
{Omegaij}

\AddEq{polylinear maps category}
{
\[
\xymatrix@C=15pt{
f:\Times\ar[r]&S_1
&
g:\Times\ar[r]&S_2
}
\]
}

\AddEq[1]{R1 X->}
{
\[R_1:#1\rightarrow #1'\]
}

\DefEq
{
\[
R\circ m=R_1(m)
\]
}
{Rm=R1m}

\DefEq
{
$r_1=q_1\circ p_1$, $r_2=q_2\circ p_2$.
}
{r=q*p}

\DefEq
{
\[
r_2:A_2\rightarrow A_2
\]
}
{R:A2->A2}

\DefEq
{
\[
a\pC i0xa\pC i1
\]
}
{Sum over repeated index}

\DefEq
{
\symb{\mathcal L(D;A_1\times...\times A_n\rightarrow S)}{set polylinear maps}1
}
{set polylinear maps}

\DefEquation
{
(A_1\otimes A_2)\otimes A_3=A_1\otimes(A_2\otimes A_3)
=A_1\otimes A_2\otimes A_3
}
{A1xA2xA3}

\DefEq
{
\[
(v_1, ..., v_n)\in V_1\times...\times V_n
\rightarrow v_1\otimes...\otimes v_n\in V_1\otimes...\otimes V_n
\]
}
{V times->V otimes}

\DefEquation
{
\begin{split}
&\,a_1\otimes...\otimes(a_i+b_i)\otimes...\otimes a_n
\\
=&\,a_1\otimes...\otimes a_i\otimes...\otimes a_n+
a_1\otimes...\otimes b_i\otimes...\otimes a_n
\\
&\,a_i, b_i\in A_i
\end{split}
}
{tensors 1, tensor product}

\DefEquation
{
\begin{split}
&a_1\otimes...\otimes(ca_i)\otimes...\otimes a_n=
c(a_1\otimes...\otimes a_i\otimes...\otimes a_n)
\\
&a_i\in A_i\ \ \ c\in D
\end{split}
}
{tensors 2, tensor product}

\AddEq [4]{f: A->B}
{
#1:#2\rightarrow #3#4
}

\AddEq [4]{vf: A->B}
{
$\Vector{#1}:#2\rightarrow #3$#4
}

\AddEq [5]{f:a in A->b in B}
{
\[#1:#2\in #3\rightarrow #4\in #5\]
}

\DefEquation
{
f:A^n\rightarrow A,
a=f\circ(a_1,...,a_n)
}
{polylinear map, algebra}

\DefEquation
{
a=f\pC{s}{0}^n\ \sigma_s(I_{(s\cdot 1)}\circ a_1)
\ f\pC{s}{1}^n\ ...\ \sigma_s(I_{(s\cdot n)}\circ a_n)\ f\pC{s}{n}^n
}
{polylinear map, algebra, canonical morphism}

\DefEq
{
$\{a_1,...,a_n\}$
\[
\sigma_s=
\begin{pmatrix}
a_1&...&a_n
\\
\sigma_s(a_1)&...&\sigma_s(a_n)
\end{pmatrix}
\]
}
{transposition of set of variables, algebra}

\DefEq
{
$I_{(s\cdot 1)}$, ..., $I_{(s\cdot n)}\in\mathcal L(D;A\rightarrow A)$
}
{I1n in L(A;A)}

\DefEq
{
\labelItem{monomial of power 0}
$p_0(x)=a_0$, $a_0\in A$.
}
{p0(x)=a0}

\AddEq{monomial of power k}
{
\labelItem{monomial of power k}
\[
p_k(x)=p_{k-1}(x)(F\circ x)a_k
\]
}

\AddEq{p1(x)=aFxa}
{
$p_1(x)=a_0(F\circ x)a_1$.
}

\AddEq{Basis F=...}
{
$\Basis F=(F_1,...,F_n)$.
}

\AddEq{Fi=Fj}
{
$F_{(i)}=F_{(j)}$
}

\AddEq{F=...}
{
$F=(F_{(1)},...,F_{(k)})$
}

\AddEq[1]{set of homogeneous polynomials}
{
\symb{A_{#1}[x]}{set of homogeneous polynomials}1
}

\DefEq
{
$a\in A^{\otimes (n+1)}$\Pt
}
{a in Aoxn+1}

\DefEq
{
$x_1=...=x_n=x$,
}
{x1=xn=x}

\DefEq
{
\[
a\circ F\circ x^n=a\circ(F_{(1)},...,F_{(n)})\circ(x_1\otimes...\otimes x_n)
\]
}
{a xn=}

\AddEq{F1n in F}
{
$F_{(1)}$, ..., $F_{(n)}\in\Basis F$.
}

\AddEq{F[k]=...}
{
\[
F_{[k]}=(F_{[k](1)},...,F_{[k](n)})
\]
}

\DefEq
{
\[
a=a_{i\cdot 0}\otimes a_{i\cdot 1}\otimes...\otimes a_{i\cdot n}
\ \ \ i\in I
\]
}
{a=oxi}

\DefEq
{
$\sigma=\{\sigma_i\in S(n):i\in I\}$
}
{si in Sn}

\DefEq
{
\[
(a,\sigma):A^{\times n}\rightarrow A
\]
}
{ox:An->A}

\DefEq
{
\begin{align*}
(a,\sigma)\circ (b_1,...,b_n)&=
(a_{i\cdot 0}\otimes a_{i\cdot 1}\otimes...\otimes a_{i\cdot n},\sigma_i)\circ (b_1,...,b_n)
\\&=
a_{i\cdot 0}\sigma_i(b_1)a_{i\cdot 1}...\sigma_i(b_n)a_{i\cdot n}
\end{align*}
}
{ox circ =}

\AddEq{p(x)=a circ xk}
{
\[
\begin{matrix}
p(x)=a_{[s]}\circ F_{[s]}\circ x^k
&a_{[s]}\in A^{\otimes(k+1)}
\end{matrix}
\]
}

\DefEq
{
\symb{A[x]}{algebra of polynomials over algebra}{}
\[
\ShowSymbol{algebra of polynomials over algebra}{}
=\bigoplus_{n=0}^{\infty}A_n[x]
\]
}
{algebra of polynomials over algebra}

\AddEq[3]{a1n}
{
$#1_1$, ..., $#1_{#2}$#3
}

\AddEq[3]{A1n}
{
$#1^1$, ..., $#1^{#2}$#3
}

\AddEq[3]{aU A1n}
{
$\aU{#1}1$, ..., $\aU{#1}{#2}$#3
}

\DefEq
{
$a_0\in U$.
}
{a0 in U}

\DefEq
{
\[
f:U\rightarrow \mathcal L(D;A^p\rightarrow B)
\]
}
{U->L(Ap,B)}

\DefEquation
{
\begin{array}{r@{\ }l}
&((a_0,...,a_n,\sigma)\circ(f_1,...,f_n))\circ(x_1,...,x_n)
\\=&
(a_0\sigma(f_1)a_1...a_{n-1}\sigma(f_n)a_n)\circ(x_1,...,x_n)
\\=&
a_0\sigma(f_1\circ x_1)a_1...a_{n-1}\sigma(f_n\circ x_n)a_n
\end{array}
}
{n linear map A LA}

\AddEq[3]{set polylinear maps An}
{
\symb{\mathcal L(#1;#2^n\rightarrow #3)}{set polylinear maps An}1
}

\DefEquation
{
\begin{array}{r@{\,}l@{\ \ \ }r@{\,}l@{\ \ \ }r@{\,}l}
a_k&=a_{k\cdot 0}\underline{\otimes}(1\otimes a_{k\cdot 1})
&a_{k\cdot 0}&\in A^{\otimes k}&a_{k\cdot 1}&\in A
\end{array}
}
{p(x)=a k-1 k 3}

\AddEq [2]{M M1 Mn}
{
\[#1=\max(#1_1,...,#1_{#2})\]
}

\AddEq [2]{M max M1 Mn}
{
\[#1=\max(\aD{#1}1,...,\aD{#1}{#2})\]
}

\DefEq
{
\[N=\max(N_1,N_2)\]
}
{N N1 N2}

\DefEq
{
f_i(x)=\lim_{m\rightarrow\infty}f_{i\cdot m}(x)
}
{fi(x)=lim}

\DefEq
{
\delta_1\in R,\ \delta_1>0
}
{delta1 in R}

\DefEq
{
\(m>M_i\)
}
{m>Mi}

\DefEq
{
\[h_m=f_{1\cdot m}...f_{n\cdot m}\omega\]
}
{hm=f1m...fnm omega}

\DefEq
{
\[
\epsilon_1(\delta)<\epsilon
\]
}
{e(d)<e}

\DefEq
{
\[\delta_1=0\Rightarrow\epsilon_1=0\]
}
{d1=0=>e1=0}

\DefEq
{
F_i\ge 0
}
{F>0}

\DefEq
{
F_i=\sup\|f_i(x)\|
}
{F=sup|f|}

\AddEq[3]{a in AA}
{
$#1\in #2\otimes #2$#3
}

%% file: Prolog.Cite.tex

\input{\FilePrefix Prolog.Cite.\TheLanguage}

\DefCiteBib{1812.03397}
{
Ivan Kyrchei,
Linear differential systems over the quaternion skew field,
eprint \href{http://arxiv.org/abs/1812.03397}{arXiv:1812.03397} (2018)
}

\DefCiteBib{math.QA-0208146}
{
I. Gelfand, S. Gelfand, V. Retakh, R. Wilson,
Quasideterminants,\\
eprint \href{http://arxiv.org/abs/math.QA/0208146}{arXiv:math.QA/0208146} (2002)
}

\DefCiteBib{q-alg-9705026}
{
I. Gelfand, V. Retakh,
Quasideterminants, I,\\
eprint \href{http://arxiv.org/abs/q-alg/9705026}{arXiv:q-alg/9705026} (1997)
}

\DefCiteBib{1510.02224}
{
Kit Ian Kou, Yong-Hui Xia.
Linear Quaternion Differential Equations: Basic Theory and Fundamental Results.\\
eprint \href{http://arxiv.org/abs/1510.02224}{arXiv:1510.02224} (2017)
}

\DefCiteBib{Cohn: Skew Fields}
{
Paul M. Cohn,
Skew Fields,
Cambridge University Press, 1995
}

\DefCiteBib{Sudbery Quaternionic Analysis}
{
A. Sudbery,
Quaternionic Analysis,\\
Math. Proc. Camb. Phil. Soc. (1979), {\bf 85}, 199 - 225
}

%% file: Prolog.Cite.English.tex

\DefCiteBib{1302.7204}
{
Aleks Kleyn,
Polynomial over Associative $D$-Algebra,\\
eprint \href{http://arxiv.org/abs/1302.7204}{arXiv:1302.7204} (2013)
}

\DefCiteBib{1801.01628}
{
Aleks Kleyn,
Differential Equation over Banach Algebra,\\
eprint \href{http://arxiv.org/abs/1801.01628}{arXiv:1801.01628} (2018)
}

\DefCiteBib{1908.04418}
{
Aleks Kleyn,
Diagram of Representations of Universal Algebras,\\
eprint \href{http://arxiv.org/abs/1908.04418}{arXiv:1908.04418} (2019)
}

\DefCiteBib{Korn}
{
Granino A. Korn, Theresa M. Korn,
Mathematical Handbook for Scientists and Engineer,
McGraw-Hill Book Company, New York, San Francisco,
Toronto, London, Sydney, 1968
}

\DefCiteBib{Rashevsky}
{
P. K. Rashevsky,
Riemann Geometry and Tensor Calculus,\\
Moscow, Nauka, 1967
}

\DefCiteBib{1502.04063}
{
Aleks Kleyn,
Linear Map of $D$\Hyph Algebra,\\
eprint \href{http://arxiv.org/abs/1502.04063}{arXiv:1502.04063} (2015)
}

\DefCiteBib{0812.4763}
{
Aleks Kleyn,
Introduction into Calculus over Division Ring,\\
eprint \href{http://arxiv.org/abs/0812.4763}{arXiv:0812.4763} (2010)
}

\DefCiteBib{1601.03259}
{
Aleks Kleyn,
Introduction into Calculus over Banach Algebra,\\
eprint \href{http://arxiv.org/abs/1601.03259}{arXiv:1601.03259} (2016)
}

\DefCiteBib{2207.06506}
{
Aleks Kleyn,
Introduction into Noncommutative Algebra,
Volume 1, Division Algebra\\
eprint \href{http://arxiv.org/abs/2207.06506}{arXiv:2207.06506} (2022)
}

%% file: Stmt.Representation.English.tex
\input{Stmt.Representation.Eq}

\DefDefinition{morphism of representations of universal algebra}
{
Let
\ShowEq{f:A->*B}f{A_1}{A_2}
be representation of $\Omega_1$\Hyph algebra $A_1$
in $\Omega_2$\Hyph algebra $A_2$ and
\ShowEq{f:A->*B}g{B_1}{B_2}
be representation of $\Omega_1$\Hyph algebra $B_1$
in $\Omega_2$\Hyph algebra $B_2$.
For
\ShowEq{i=1,2},
let the map
\ShowEq{ri:A->B}
be homomorphism of $\Omega_i$\Hyph algebra.
The tuple of maps
\ShowEq{map r12}r{}
such, that
\ShowEq{morphism of representations of universal algebra, definition, 2}
is called
\AddIndex{morphism of representations from $f$ into $g$}
{morphism of representations from f into g}.
We also say that
\AddIndex{morphism of representations of $\Omega_1$\Hyph algebra
in $\Omega_2$\Hyph algebra}
{morphism of representations of Omega1 algebra in Omega2 algebra} is defined.
}

\DefRemark{morphism of representations of universal algebra}
{
There are two ways
to interpret
\eqRef{morphism of representations of universal algebra, 2m}{representation}
\begin{itemize}
\item Let transformation $\BlueText{f(a)}$ map $m\in A_2$
into $\BlueText{f(a)}(m)$.
Then transformation
\ShowEq{g(r1(a))}
maps
\ShowEq{r2(m)in B2}
into
\ShowEq{r2(f(a,m))}
\item We represent morphism of representations from $f$ into $g$
using diagram
\DrawEq{morphism of representations of universal algebra, 2m 1}{ShadedDefinition}
From \EqRef{morphism of representations of universal algebra, definition, 2},
it follows that diagram $(1)$ is commutative.
\end{itemize}
We also use diagram
\ShowEq{morphism of representations of universal algebra, definition, 2m 2}
instead of diagram
\eqRef{morphism of representations of universal algebra, 2m 1}{ShadedDefinition}.
}

\DefRemark{morphism of representations of universal algebra as map}
{
We may consider a pair of maps $r_1$, $r_2$ as map
\ShowEq{F:A1+A2->B1+B2}
such that
\ShowEq{F:A1+A2->B1+B2 1}
Therefore, hereinafter the tuple of maps
\ShowEq{map r12}r{}
also is called map
and we will use map
\ShowEq{f:A->B}rfg
Let
\ShowEq{a=a12}
be tuple of $A$\Hyph numbers.
We will use notation
\ShowEq{r(a)12}
for image of tuple of $A$\Hyph numbers
with respect to morphism of representations $r$.
}

\DefTheorem{unique morphism of representations of universal algebra}
{
Let the representation
\ShowEq{f:A->*B}f{A_1}{A_2}
of $\Omega_1$\Hyph algebra $A_1$ be single transitive representation
and the representation
\ShowEq{f:A->*B}g{B_1}{B_2}
of $\Omega_1$\Hyph algebra $B_1$ be single transitive representation.
Given homomorphism of $\Omega_1$\Hyph algebra
\ShowEq{f:A->B}{r_1}{A_1}{B_1}
consider a homomorphism of $\Omega_2$\Hyph algebra
\ShowEq{f:A->B}{r_2}{A_2}{B_2}
such that
\ShowEq{map r12}r{}
is morphism
of representations from $f$ into $g$.
The map $H$ is unique up to
choice of image
\ShowEq{n=r2(m)}
of given element $m\in A_2$.
}

\DefLabeledFootnote[1]{iso end aut morphism}{#1}
{
I follow the definition on
page \citeBib{Cohn: Universal Algebra}\Hyph 49.
}

\DefConvention[2]{A number}
{
Element of $#1$\Hyph #2 $A$ is called
\AddIndex{$A$\Hyph number}{A number}.
For instance, complex number is also called
$C$\Hyph number, and quaternion is called $H$\Hyph number.
}

\DefTheorem{map is reduced morphism of representations iff}
{
Let
\ShowEq{f:A->*B}f{A_1}{A_2}
be representation of $\Omega_1$\Hyph algebra $A_1$
in $\Omega_2$\Hyph algebra $A_2$ and
\ShowEq{f:A->*B}g{A_1}{B_2}
be representation of $\Omega_1$\Hyph algebra $A_1$
in $\Omega_2$\Hyph algebra $B_2$.
The map
\ShowEq{f:A->B}{r_2}{A_2}{B_2}
is reduced morphism of representations iff
\DrawEq{reduced morphism representation =}{definition}
}

\DefDefinition{morphism of representation f}
{
If representation $f$ and $g$ coincide, then morphism of representations
\ShowEq{map r12}r{}
is called
\AddIndex{morphism of representation $f$}{morphism of representation f}.
}

\DefRemark{notation for morphism of representations}
{
Consider morphism of representations
\ShowEq{r12:A->B}rAB
We denote elements of the set $B_1$ by letter using pattern $b\in B_1$.
However if we want to show that $b$ is image of element
\ShowEq{a in A1},
we use notation $\RedText{r_1(a)}$.
Thus equation
\ShowEq{r1(a)=r1(a)}
means that $\RedText{r_1(a)}$ (in left part of equation)
is image
\ShowEq{a in A1}{}
(in right part of equation).
Using such considerations, we denote
element of set $B_2$ as $\BlueText{r_2(m)}$.
We will follow this convention when we consider correspondences
between homomorphisms of $\Omega_1$\Hyph algebra
and maps between sets
where we defined corresponding representations.
}

\DefTheorem{Tuple of maps is morphism of representations iff}
{
Let
\ShowEq{f:A->*B}f{A_1}{A_2}
be representation of $\Omega_1$\Hyph algebra $A_1$
in $\Omega_2$\Hyph algebra $A_2$ and
\ShowEq{f:A->*B}g{B_1}{B_2}
be representation of $\Omega_1$\Hyph algebra $B_1$
in $\Omega_2$\Hyph algebra $B_2$.
The map
\ShowEq{r12:A->B}rAB
is morphism of representations iff
\DrawEq{morphism of representations of universal algebra, 2m}{representation}
}

\DefDefinition{tower of representations}
{
Consider set of $\Omega_k$\Hyph algebras $A_k$, \Kn1.
Let $A=(A_1,...,A_n)$.
Let $f=(f_{1\,2},...,f_{n-1\,n})$.
Set of representations $f_{k\,k+1}$, \Kn1,
of $\Omega_k$\Hyph algebra $A_k$ in
$\Omega_{k+1}$\Hyph algebra $A_{k+1}$
is called
\AddIndex{tower $(f,A)$ of representations of $\Omega$\Hyph algebras}
{tower of representations of algebras}.
}

\DefDefinition{diagram of representations}
{
\AddIndex{Diagram $(f,A)$ of representations of universal algebras}
{diagram of representations of algebras}
is oriented graph such that
\StartLabelItem
\begin{enumerate}
\item
the vertex of $A_k$, \Kn1, is $\Omega_k$\Hyph algebra;
\item
the edge $f_{kl}$ is representation
of $\Omega_k$\Hyph algebra $A_k$ in
$\Omega_l$\Hyph algebra $A_l$;
\end{enumerate}
We require that this graph
is connected graph and does not have loops.
Let $A_{[0]}$ be set of initial vertices of the graph.
Let $A_{[k]}$ be set of vertices of the graph
for which the maximum path from the initial vertices is $k$.
}

\DefRemark{diagram of representations}
{
Since different vertices of the graph
can be the same algebra,
then we denote
\ShowEq{A=A1n}A{}
the set of universal algebras which are distinct.
From the equality
\ShowEq{A=A(1n)1n}A
it follows that, for any index $(i)$,
there exists at least one index $i$ such that
\ShowEq{A(i)=Ai}i.
If there are two sets of sets
\ShowEq{A=A1n}A,
\ShowEq{A=A1n}B
and there is a map
\ShowEq{hi:Ai->Bi}
for an index $(i)$,
then also there is a map
\ShowEq{f:A->B}{h_i}{A_i}{B_i}
for any index $i$ such that
\ShowEq{A(i)=Ai}i{}
and in this case $h_i=h_{(i)}$.
}

\DefTheorem{diagram of representations}
{
{\bf(}\AddIndex{Induction over diagram of representations}{induction over diagram of representations}{\bf)}.
Let the theorem $\mathcal T$ be true for the set
of universal algebras $A_{[0]}$
of diagram $(f,A)$ of representations of universal algebras.
Let the statement that the theorem $\mathcal T$ is true for the set
of universal algebras $A_{[k]}$
of diagram $(f,A)$ of representations imply the statement
that the theorem $\mathcal T$ is true for the set
of universal algebras $A_{[k+1]}$
of diagram $(f,A)$ of representations.
Then the theorem $\mathcal T$ is true for the set
of universal algebras
of diagram $(f,A)$ of representations.
}

\DefDefinitionNote{commutative diagram of representations}
{
Diagram $(f,A)$ of representations of universal algebras
is called
\AddIndex{commutative}{commutative diagram of representations}
when diagram meets the following requirement.
for each pair of representations
\ShowEq{f:A->*B}{f_{ik}}{A_i}{A_k}
\ShowEq{f:A->*B}{f_{jk}}{A_j}{A_k}
the following equality is true\,\footnotemark
\ShowEq{fik fjk = fjk fik}
}
{
Metaphorically speaking,
representations $f_{ik}$ and $f_{jk}$
are transparent to each other.
}

\DefTheoremNote{diagram of representations, define map fik}
{
Let
\ShowEq{f:i->*j}ij
be representation of $\Omega_i$\Hyph algebra $A_i$
in $\Omega_j$\Hyph algebra $A_j$.
Let
\ShowEq{f:i->*j}jk
be representation of $\Omega_j$\Hyph algebra $A_j$
in $\Omega_k$\Hyph algebra $A_k$.
We represent the fragment\,\footnotemark
\DrawEq[ijk]{Ai->*Aj->*Ak}{}
of the diagram of representations using the diagram
\ShowEq{tower of representations, 1 2 3}
The map
\ShowEq{f:A->End 2 3}
is defined by the equality
\ShowEq{define map f13}
where
\ShowEq{a in A}i,
\ShowEq{a in A}j.
If the representation $f_{jk}$ is effective
and the representation $f_{ij}$ is free,
then the map $f_{ijk}$ is free representation
\ShowEq{f:1->*3}
of $\Omega_i$\Hyph algebra $A_i$
in $\Omega_j$\Hyph algebra
\ShowEq{End Ak}Ak.
}
{
The theorem
\RefTheorem{diagram of representations, define map fik}
states that transformations in diagram of representations are coordinated.
}

\DefDefinition{Morphism of Diagram of Representations}
{
Let
$(f,A)$
be the diagram of representations where
\ShowEq{A=A1n}A{}
is the set of universal algebras.
Let
$(B,g)$
be the diagram of representations where
\ShowEq{A=A1n}B{}
is the set of universal algebras.
The set of maps
\ShowEq{A=A1n}h{}
\ShowEq{hi:Ai->Bi}
is called \AddIndex{morphism from diagram of representations
$(f,A)$ into diagram of representations $(B,g)$}
{morphism from diagram of representations into diagram of representations},
if for any indexes $(i)$, $(j)$, $i$, $j$ such that
\ShowEq{A(i)=Ai}i,
\ShowEq{A(i)=Ai}j{}
and for any representation
\ShowEq{f:A->*B}{f_{ji}}{A_j}{A_i}
the tuple of maps $(h_j\ \ h_i)$ is
morphism of representations from $f_{ji}$ into $g_{ji}$.
}

\AddEq{remark: Morphism of Diagram of Representations}
{
For any representation $f_{ij}$,
\Kn[i]1, \Kn[j]1,
we have diagram
\ShowEq{morphism of diagram of representations of F algebra, level k, diagram}
Equalities
\ShowEq{morphism of diagram of representations, level k}
\ShowEq{morphism of diagram of representations, levels k k+1}
express commutativity of diagram (1).
}

\DefTheorem{reduced polymorphism of representations}
{
Let the map $r_2$ be reduced polymorphism of
effective representations $f_1$, ..., $f_n$ into effective representation $f$.
\begin{itemize}
\item
For any
\ShowEq{k,1n}
the map $r_2$ satisfies to the equality
\ShowEq{reduced polymorphism of representation, ak}
\item
For any
\ShowEq{kl,1n}
the map $r_2$ satisfies to the equality
\ShowEq{reduced polymorphism of representation, akl}
\item
Let $\omega_2\in\Omega_2(p)$.
For any
\ShowEq{k,1n}
the map $r_2$ satisfies to the equality
\ShowEq{reduced polymorphism of representation, omega2}
\end{itemize}
}

\AddEq{remark: morphism of representation of Omega group}
{
\begin{remark}
\labelRemark{morphism of \SideWS representation of Omega group}
Let the map
\ShowEq{f:A->*B}f{A_1}{A_2}
be the \SideNS\Hyph side representation
of multiplicative $\Omega$\Hyph group $A_1$
in $\Omega_2$\Hyph algebra $A_2$.
Let the map
\ShowEq{f:A->*B}g{B_1}{B_2}
be the \SideNS\Hyph side representation
of multiplicative $\Omega$\Hyph group $B_1$
in $\Omega_2$\Hyph algebra $B_2$.
Let the map
\ShowEq{r12:A->B}rAB
be morphism of representations.
We use notation
\ShowEq{r2a=r2 o a}
for image of $A_2$\Hyph number $a_2$
with respect to the map $r_2$.
Then we can write the equality
\eqRef{morphism of representations of universal algebra, 2m}{representation}
in the following form
\ShowEq{morphism of \SideWS representation of Omega group}
\qed
\end{remark}
}

\DefTheorem{proper definition of orbit}
{
Let the map
\ShowEq{f:A->*B}f{A_1}{A_2}
be the left\Hyph side representation
of multiplicative $\Omega$\Hyph group $A_1$.
Let
\ShowEq{orbit, proposition}
Then
\ShowEq{A1a2=A1b2}
}

\AddEq{definition: orbit of representation}
{
\begin{ShadedDefinition}
\labelDefinition{orbit of \SideNS-side representation}
Let $A_1$ be $\Omega$\Hyph groupoid with product
\ShowEq{(ab)->ab}
Let the map
\ShowEq{f:A->*B}f{A_1}{A_2}
be the \SideNS\Hyph side representation
of $\Omega$\Hyph groupoid $A_1$
in $\Omega_2$\Hyph algebra $A_2$.
For any
\ShowEq{a in A}2,
we define
\AddIndex{orbit of representation}{orbit of representation}
of the $\Omega$\Hyph groupoid $A_1$ as set
\ShowEq{orbit of \SideNS-side representation}
\end{ShadedDefinition}
}

\AddEq{Let A->*AB2 be representations}
{
Let
\ShowEq{f:A->*B}f{A_1}{A_2}
be representation of $\Omega_1$\Hyph algebra $A_1$
in $\Omega_2$\Hyph algebra $A_2$.
Let
\ShowEq{f:A->*B}g{A_1}{B_2}
be representation of $\Omega_1$\Hyph algebra $A_1$
in $\Omega_2$\Hyph algebra $B_2$.
}

\AddEq{Let XY be generating sets}
{
Let $X$ be the generating set of the representation
\ShowEq{f:A->*B}f{A_1}{A_2}
of $\Omega_1$\Hyph algebra $A_1$
in $\Omega_2$\Hyph algebra $A_2$.
Let $Y$ be the generating set of the representation
\ShowEq{f:A->*B}g{A_1}{B_2}
of $\Omega_1$\Hyph algebra $A_1$
in $\Omega_2$\Hyph algebra $B_2$.
}

\DefDefinitionNote{product in category}
{
Let $\mathcal A$ be a category.
Let
\ShowEq{set Bi}B
be the set of objects of $\mathcal A$.
Object
\ShowEq{product in category}
and set of morphisms
\ShowEq{set f:A->B}fP{B_i}
is called a
\AddIndex{product of set of objects
\ShowEq{set Bi}B
in category $\mathcal A$}
{product in category}\,\footnotemark
if for any object $R$
and set of morphisms
\ShowEq{set f:A->B}gR{B_i}
there exists a unique morphism
\ShowEq{f:A->B}hRP
such that diagram
\ShowEq{product in category diagram}
is commutative for all $i\in I$.

If $|I|=n$, then we also will use notation
\ShowEq{product in category, 1 n}
for product of set of objects
$\{B_i,\iI\}$ in $\mathcal A$.
}
{
I made definition according to \citeBib{Serge Lang}, page 58.
}

\DefExample{Cartesian product of sets}
{
Let \(\mathcal S\) be the category of sets.\,\footnote{See
also the example in
\citeBib{Serge Lang},
page 59.
}
According to the definition
\ShowEq{ref product in category}
Cartesian product
\ShowEq{Cartesian product of sets}
of family of sets
\ShowEq{Ai iI}A{}
and family of projections on the \(i\)\Hyph th factor
\ShowEq{projection on i factor}
are product in the category \(\mathcal S\).
}

\DefTheorem{product of effective representations}
{
In category
\ShowEq{A1(mA2)}
there exists product
of effective representations of $\Omega_1$\Hyph algebra $A_1$
in $\Omega_2$\Hyph algebra
and the product is
effective representation of $\Omega_1$\Hyph algebra $A_1$.
}

\DefTheoremNote{structure of subrepresentations}
{
Let\,\footnotemark
\ShowEq{f:A->*B}g{A_1}{A_2}
be representation of $\Omega_1$\Hyph algebra $A_1$
in $\Omega_2$\Hyph algebra $A_2$.
Let $X\subset A_2$.
Define a subset $X_k\subset A_2$ by induction on $k$.
\ShowEq{structure of subrepresentations}
Then
\ShowEq{structure of subrepresentations, 1}
}{
The statement of theorem is similar to the
statement of theorem 5.1, \citeBib{Cohn: Universal Algebra}, page 79.
}

\DefDefinitionNote{direct sum of Abelian groups}
{
Coproduct in category of Abelian groups $Ab$ is called
\AddIndex{direct sum}{direct sum}.\,\footnotemark
We will use notation
\ShowEq{direct sum of Abelian groups}
for direct sum of Abelian groups $A$ and $B$.
}{See also definition in \citeBib{Serge Lang}, pages 36, 37.}

\DefDefinitionNote{coproduct in category}
{
Let $\mathcal A$ be a category.
Let
\ShowEq{set Bi}B
be the set of objects of $\mathcal A$.
Let
\ShowEq{category A[Bi]}
be a category
whose objects are tuples $(P,f)$
where $P$ is object of category $\mathcal A$
and $f$ is set of morphisms
\ShowEq{set f:A->B}f{B_i}P
Universally repelling object of category
\ShowEq{category A[Bi]}
\ShowEq{coproduct in category}
is called a
\AddIndex{coproduct of set of objects
\ShowEq{set Bi}B
in category $\mathcal A$}
{coproduct in category}.\,\footnotemark

If $|I|=n$, then we also will use notation
\ShowEq{coproduct in category, 1 n}
for coproduct of set of objects
\ShowEq{set Bi}B
in $\mathcal A$.
}
{
I made definition according to the definition on page
\citeBib{Serge Lang}\Hyph 59.
}

\AddEq{definition: free Abelian group}
{
\begin{ShadedDefinition}
\labelDefinition{free Abelian group}
Let $S$ be a set and\,\footnotemark
\ShowEq{category AbS}
be category objects of which are maps
\ShowEq{f:A->B}fSG
of the set $S$ into Abelian groups.
If
\ShowEq{f:A->B}fSG
\ShowEq{f:A->B}{f'}S{G'}
are objects of the category
\ShowEq{category AbS},
then morphism from $f$ to $f'$
is the homomorphism of groups
\ShowEq{f:A->B}gG{G'}
such that the diagram
\DrawEq{AbS ff'g}{}
is commutative.
\end{ShadedDefinition}
\footnotetext{\,
See also the definition in \citeBib{Serge Lang}, pages 37.
}
}

\AddEq{theorem: free Abelian group}
{
\begin{ShadedTheorem}
\labelTheorem{free Abelian group}
There exists\,\footnotemark
universally repelling object $G$ of the category
\ShowEq{category AbS}
\StartLabelItem
\begin{enumerate}
\item
\ShowEq{S in G}
\item
The set $S$ generates Abelian group $G$.
\labelItem{set generates Abelian group}
\end{enumerate}
Abelian group $G$ is called
\AddIndex{free Abelian group}{free Abelian group}
generated by the set $S$.
\end{ShadedTheorem}
\footnotetext{\,
See also similar statement in \citeBib{Serge Lang}, pages 38.
}
}

\DefProof{free Abelian group}
{
Let $G$
be set of maps
\ShowEq{f:A->B}hSZ
such that the set
\ShowEq{x in -> h(x)}
is finite.

\begin{ShadedLemma}
{\it
The set $G$ is Abelian group.
}
\end{ShadedLemma}
{\sc Proof.}
Let
\ShowEq{h12 in G}
According to the definition, sets
\ShowEq{H12=}
are finite. Then the set
\ShowEq{H1vH2}
is finite.
Therefore,
\ShowEq{h1+h2 in G}
Properties of Abelian group are evident.
\hfill\(\odot\)

Let $k\in Z$, $x\in S$.
We denote by $k*x$ the map $h$ defined by the equality
\ShowEq{k*x=k}
The map
\ShowEq{f:S->G}
is injective and allows us to identify $S$ and image $f(S)\subseteq G$
(the statement
\RefItem{S in G})
and to use notation
\ShowEq{kx=k*x}
According to the definition of the group $G$,
we can write any map $h\in G$ as
\ShowEq{h=k*x}
where
\ShowEq{k in Z x in S}
and the set
\ShowEq{|kx ne 0|}
is finite.
The statement
\RefItem{set generates Abelian group}
follows from the theorem
\RefTheorem{structure of Abelian group}.

\begin{ShadedLemma}
{\it
The representation
\EqRef{h=k*x}
of the map $h$ is unique.
}
\end{ShadedLemma}
{\sc Proof.}
If the map $h$ admits representations
\ShowEq{h=k12*x i}
then the equality
\ShowEq{k1-k2 *x}
follows from the equality
\EqRef{h=k12*x i}.
\ShowEq{k1=k2 i}
follows from the equality
\EqRef{k1-k2 *x}.
\hfill\(\odot\)

Let
\ShowEq{f:A->B}{f'}S{G'}
be the map of the set $S$ into Abelian group $G'$.
We define the homomorphism
\ShowEq{f:A->B}gG{G'}
with request that
following diagram is commutative
\DrawEq{AbS ff'g}{1}
From the diagram
\eqRef{AbS ff'g}{1},
it follows that
\ShowEq{g(1x)=f'(x)}
According to the theorem
\RefTheorem{homomorphism f(na)=nf(a)},
since the map $g$ is homomorphism of Abelian group, then the equality
\ShowEq{g(x)=}
for any map $h\in G$
follows from equalities
\EqRef{g(1x)=f'(x)},
\EqRef{h=k*x}.
Therefore, homomorphism $g$ is defined uniquely.
According to definitions
\RefDefinition{universally repelling object of category},
\RefDefinition{free Abelian group},
Abelian group $G$ is free group.
}

\DefTheorem{homomorphism f(na)=nf(a)}
{
The homomorphism
\ShowEq{f:A->B}fG{G'}
holds the equality
\ShowEq{f(na)=nf(a)}
}

\DefTheoremNote{coproduct in category}
{
Let $\mathcal A$ be a category.
Let
\ShowEq{set Bi}B
be the set of objects of $\mathcal A$.
Object
\ShowEq{coproduct in category}
and set of morphisms
\ShowEq{set f:A->B}f{B_i}P
is called a
\AddIndex{coproduct of set of objects
\ShowEq{set Bi}B
in category $\mathcal A$}
{coproduct in category}\,\footnotemark
if for any object $R$ and set of morphisms
\ShowEq{set f:A->B}g{B_i}R
there exists a unique morphism
\ShowEq{f:A->B}hPR
such that diagram
\ShowEq{coproduct in category diagram}
is commutative for all $i\in I$.
}{I made definition according to \citeBib{Serge Lang}, page 59.}

\DefProof{coproduct in category}
{
The theorem follows from definitions
\RefDefinition{universally repelling object of category},
\RefDefinition{coproduct in category}.
}

\DefTheoremNote{direct sum of Abelian groups}
{
Let
\ShowEq{set Bi}A
be set of Abelian groups.
Let
\ShowEq{A in xAi}
be such set that
\ShowEq{(xi)in A},
if
$x_i\ne 0$
for finite number of indices $i$.
Then\,\footnotemark
\ShowEq{A=o+Ai}A
}{See also proposition \citeBib{Serge Lang}-7.1, page 37.}

\DefProof{direct sum of Abelian groups}
{
According to the theorem
\RefTheorem[\RefRepresentation]{product exists in category of Omega algebras},
there exists Abelian group $B=\prod A_i$.
According to the statement
\RefItem[\RefRepresentation]{tuple represent A number},
$B$\Hyph number $a$ can be represented as tuple
\ShowEq{set Bi}a
$A_i$\Hyph numbers.
According to the statement
\RefItem[\RefRepresentation]{operation is defined componentwise},
the sum in group $B$ is defined by the equality
\ShowEq{(ai)+(bi)=(ai+bi)}

According to construction,
$A$ is subgroup of Abelian group $B$.
The map
\ShowEq{f:A->B}{\lambda_j}{A_j}A
defined by the equality
\ShowEq{lj(x)=0x0}
is an injective homomorphism.

Let
\ShowEq{set f:A->B}f{A_i}B
be set of homomorphisms into Abelian group $B$.
We define the map
\ShowEq{f:A->B}fAB
by the equality
\DrawEq{fxi,i=sum fxi}{group}
The sum in the right side of the equality
\eqRef{fxi,i=sum fxi}{group}
is finite, since all summands, except for a finite number, equal $0$.
From the equality
\ShowEq{fi(x+y)=}
and the equality
\eqRef{fxi,i=sum fxi}{group},
it follows that
\ShowEq{f(x+y)i=}
Therefore, the map $f$ is homomorphism of Abelian group.
The equality
\ShowEq{flj=fj}
follows from equalities
\EqRef{lj(x)=0x0},
\eqRef{fxi,i=sum fxi}{group}.
Since the map $\lambda_i$ is injective,
then the map $f$ is unique.
Therefore, the theorem folows from definitions
\RefDefinition{coproduct in category},
\RefDefinition{direct sum of Abelian groups}.
}

\AddEq{remark: Non-commutative Sum}
{
Since our childhood, we know that sum does not depend on order of summands.
So, when I learned about non-commutative addition
(\citeBib{Alain Connes 1994}),
I was little bit confused.
However I was lucky to see geometry where sum of vectors is non-commutative.

At school I liked many subjects: physics, chemistry, drafting.
I did not think about mathematics seriously.
I did not like arithmetic, but I liked algebra and geometry.
I did not yet know that all roads lead to mathematics.

Drafting was my favorite subject. My teacher,
Gabis Sergey Aleksandrovich, gave me problem book in drafting;
I solved every problem from this problem book.
After this I introduced my problem and solved it.
To my surprise, the first theorem which we proved on
geometry class next year was exactly the solution to this problem.

I knew that it will be hard to enter college after school.
So I decided to enter technical school after eight class.
Bela Markovna offered me few lessons during summer;
she wanted to see how far I was prepared for the exam on math.
I remember one lesson.
When we considered line and parabola graphs,
Bela Markovna said: it is easy to draw straight line or parabola;
however, what can we do with cardiogram?
I decided to find an analytical expression 
for cardiogram.

My friend offered me to study calculus
together during summer.
Boys who studied at a technical school lived in our building;
they gave me calculus textbook for summer time.
At this moment my friend lost interest to his idea;
and I began to study calculus independently.

The first chapter of the textbook was dedicated to analytic geometry.
So I realized that this is the tool that I need to analyze curves of cardiogram.
analytic geometry determined my future.
I realized that I will prepare myself for math department of university.

I have read a lot of books dedicated to calculus.
However, at some point of time, I got the impression
that one author copied off from another.
I started to look for different books.
I do not remember where did I get the first volume of Fihtengolts.
But this book impressed me very much.
Closer to spring, I felt that the book becomes difficult;
however I continued to read this book.
In March, I had the flu with an unusually high temperature.
But I doubt that it was a flu,
because it became easier to read the first volume of Fihtengolts.

The study of calculus reflected
in my attempts to restore functional
dependence on graph. I did not yet understand
that not every map could be represented analytically.
I drew graphs, made tables of dependence of function increment
on increment of argument.
However soon this idea bothered me;
and I changed direction of my learning.

When I started to read the book dedicated to tensor calculus,
I discovered that simple theorems of linear algebra
substituted complex calculations of analytic geometry
related to classification of second-order curves.
I lost my interest to analytic geometry
after this.

Almost half a century has passed.
And suddenly, analytic geometry brought me
two pleasant surprises.
I discovered that analytic geometry
has interesting structure from the point of view
of the theory of representation of universal algebra.
I identify analytic and affine geometries;
however, I think that this is not big mistake.

The second surprise is related to the fact that I realized
that I have model of affine geometry
in Riemann space.
Since I studied metric\hyph affine manifold for a very long time,
then there was one step left to see model of affine geometry
where sum is non\Hyph commutative.
Henceforth non\Hyph commutative addition
ceased to be abstract concept for me.
}

\DefText{side representation of group}
{
From the equality
\eqRef{* \SideNS-side representation of group}{def},
it follows that
\DrawEq{o \SideNS-side representation of group}{def}
We consider the equality
\eqRef{o \SideNS-side representation of group}{def}
as
\AddIndex{associative law}{associative law}.

\begin{sloppypar}
\begin{example}
Let
\ShowEq{f:A->B}f{A_2}{A_2}
\ShowEq{f:A->B}g{A_2}{A_2}
be endomorphisms of $\Omega_2$\Hyph algebra $A_2$.
Let product
in multiplicative $\Omega$\Hyph group
\ShowEq{End O2}{A_2}{}
is composition of endomorphisms.
Since the product of maps $f$ and $g$ is defined in the same order
as these maps act on $A_2$\Hyph number,
then we consider the equality
\ShowEq{fga= o \SideNS}
as
associative law.
This allows writing of equality
\EqRef{fga= o \SideNS}
without using of brackets
\ShowEq{fga= 1 \SideNS}
as well it allows writing of equality
\eqRef{\SideNS-side representation of group}{def}
in the following form
\DrawEq{\SideNS-side 1 representation of group}{def}
\qed
\end{example}
\end{sloppypar}
}

\AddEq[4]{theorem: select second basis of representation}
{
\begin{ShadedTheorem}
\labelTheorem{select second basis of representation, #4}
Let $#2$ be passive transformation of the
basis manifold of the #3 $#4$.
Let $\Basis e_1$ be the basis of the #3 $#4$,
\ShowEq{e2=e1 o S}{#2}21
For basis $\Basis e_3$, let there exists an active transformation $#1$ such that
\ShowEq{e3=R o e1}{#1}31
Let
\ShowEq{e3=R o e1}{#1}42
Then
\ShowEq{e2=e1 o S}{#2}43
\end{ShadedTheorem}
}

\AddEq[3]{proof: select second basis of representation}
{
\begin{proof}
According to the equality
\eqRef{active transformation}{#3},
active transformation of coordinates of basis $\Basis e_3$ has form
\DrawEq{active transformation, 1}{#3}
Let
\ShowEq{e2=e1 o S}{#2}53
From the equality
\eqRef{passive transformation}{#3},
it follows that
\DrawEq{passive transformation, 1}{#3}
From match of expressions in equalities
\eqRef{active transformation, 1}{#3},
\eqRef{passive transformation, 1}{#3},
it follows that
\ShowEq{Basis e4=Basis e5}
Therefore, the diagram
\ShowEq{passive and active maps}{#1}{#2}
is commutative.
\end{proof}
}

\AddEq[3]{definition: basis manifold}
{
The set
\ShowEq{basis manifold}
of bases of #1 $#2$ is called
\AddIndex{basis manifold}{basis manifold}
of #1 $#2$.
}

\AddEq{ref active representation, representation}
{
According to theorems
\RefTheorem{superposition of coordinates, representation},
\RefTheorem{Automorphism of representation maps a basis into basis},
}

\AddEq{ref active representation, diagram of representations}
{
According to the theorem
\RefTheorem{Automorphism of diagram of representations maps a basis into basis}
and to the definition
\RefDefinition{superposition of coordinates and words, diagram of representations},
}

\AddEq[4]{definition: active representation in basis manifold}
{
\begin{ShadedDefinition}
\labelDefinition{active representation in basis manifold, #1}
\ShowEq{ref active representation, #1}
automorphism $#4$ of the #2 $#3$
generates transformation
\DrawEq[{#4}]{active transformation}{#1}
of the basis manifold of #2.
This transformation is called
\AddIndex{active}{active transformation of basis}.
According to the theorem
\RefTheorem{group of automorphisms of #1},
we defined left\Hyph side representation
\ShowEq{active representation in basis manifold}
\ShowEq{active representation in basis manifold def}
of group $GA(f)$
in basis manifold $\mathcal B[f]$.
Representation $A(f)$ is called
\AddIndex{active representation}{active representation in basis manifold}.
According to the corollary
\RefCorollary{automorphism uniquely defined by image of basis, #1},
this representation is single transitive.
\qed
\end{ShadedDefinition}
}

\AddEq[3]{remark: active representation}
{
\begin{remark}
\labelRemark{active representation, #1}
According to remark
\RefRemark{X is quasibasis of #1},
it is possible that there exist bases of #2 $#3$
such that there is no active transformation between them.
Then we consider the orbit of selected basis
as basis manifold.
Therefore, it is possible that the #2 $#3$ has different basis manifolds.
We will assume that we have chosen a basis manifold.
\end{remark}
}

\AddEq[5]{theorem: Coordinate transformation does not depend}
{
\begin{ShadedTheorem}
\labelTheorem{Coordinate transformation does not depend, #1}
Let passive transformation $#2\in GA(f)$ maps
basis $\Basis e_1\in\mathcal{B}[f]$
into basis $\Basis e_2\in\mathcal{B}[f]$
\DrawEq[{#2}]{passive transformation S}{#1}
Let #4 $#3$ has #5
\DrawEq[1{#3}]{Omega_2 word, basis e}{1 #1}
relative to basis $\Basis e_1$ and has #5
\DrawEq[2{#3}]{Omega_2 word, basis e}{2 #1}
relative to basis $\Basis e_2$.
Coordinate transformation
\DrawEq[{#3}{#2}]{coordinate transformation}{#1}
does not depend on #4 $#3$ or basis $\Basis e_1$, but is
defined only by coordinates of #4 $#3$
relative to basis $\Basis e_1$.
\end{ShadedTheorem}
}

\AddEq[3]{proof: Coordinate transformation does not depend}
{
\begin{proof}
From \eqRef{passive transformation S}{#1}
and
\eqRef{Omega_2 word, basis e}{2 #1},
it follows that
\DrawEq[{#2}{#3}]{Omega word in representation, basis Y, 1}{#1}
Comparing \eqRef{Omega_2 word, basis e}{1 #1}
and
\eqRef{Omega word in representation, basis Y, 1}{#1}
we get
\DrawEq[{#2}{#3}]{coordinate transformation in representation, 1}{#1}
Since $#2$ is automorphism, then the equality
\eqRef{coordinate transformation}{#1}
follows from
\eqRef{coordinate transformation in representation, 1}{#1}
and the theorem
\RefTheorem{coordinates of inverse transformation, #1}.
\end{proof}
}

\AddEq[2]{remark: active and passive transformation}
{
An active transformation changes a basis of the #1
and #2 uniformly
and coordinates of $\Omega_2$\Hyph number relative basis do not change.
A passive transformation changes only the basis and it leads to change
of coordinates of #2 relative to the basis.
}

\AddEq[2]{theorem: Coordinate transformations form right-side representation}
{
\begin{ShadedTheorem}
\labelTheorem{Coordinate transformations form right-side representation, #1}
Coordinate transformations
\eqRef{coordinate transformation}{#1}
form effective contravariant
right\Hyph side representation of group $GA(f)$ which is called
\AddIndex{coordinate representation}{coordinate representation}
in #2.
\end{ShadedTheorem}
}

\AddEq[8]{proof: Coordinate transformations form right-side representation}
{
\begin{proof}
According to corollary
\RefCorollary{map of coordinates of #1},
the transformation
\eqRef{coordinate transformation}{#1}
is the endomorphism of #4\,\footnote{This transformation
does not generate an endomorphism of the #4 $#6$. Coordinates change because
basis relative which we determinate coordinates changes. However,
#7, coordinates of which we are considering, does not change.}
#5

Suppose we have two consecutive passive transformations
$#2$ and $#8$. Coordinate transformation
\DrawEq[{#3}{#2}]{coordinate transformation}{#2 #1}
corresponds to passive transformation $#2$.
Coordinate transformation
\DrawEq[{#3}{#8}]{coordinate transformation}{#8 #1}
corresponds to passive transformation $#8$.
According to the theorem
\RefTheorem{passive representation in basis manifold, diagram of representations},
product of coordinate transformations
\eqRef{coordinate transformation}{#2 #1}
and
\eqRef{coordinate transformation}{#8 #1}
has form
\ShowEq{coordinate transformation in representation, TS}{#3}{#2}{#8}
and is coordinate transformation
corresponding to passive transformation
\ShowEq{coordinate transformation in representation, TS 1}{#8}{#2}
According to theorems
\RefTheorem{coordinates of inverse transformation, #1},
\RefTheorem{coordinates of product of inverse transformation, #1}
and to the definition
\RefDefinition{Left side contravariant representation},
coordinate transformations
form
right\Hyph side contravariant representation of group $GA(f)$.

Suppose coordinate transformation does not change coordinates of selected basis.
Then unit of group $GA(f)$ corresponds to it because representation
is single transitive. Therefore,
coordinate representation is effective.
\end{proof}
}

\AddEq[5]{remark: geometric object}
{
Passive representation $P(g)$ is coordinated
with passive representation $P(f)$,
if there exists homomorphism $h$ of group $GA(f)$ into group $GA(g)$.
Consider diagram
\ShowEq{passive representation coordinated with passive representation}
Since maps $P(f)$, $P(g)$ are isomorphisms of group,
then map $H$ is homomorphism of groups.
Therefore, map $f'$ is representation of group $GA(f)$
in basis manifold $\mathcal B(g)$.
According to design, passive transformation $H(#2)$ of basis manifold $\mathcal B(g)$
corresponds to passive transformation $#2$ of basis manifold $\mathcal B(f)$
\DrawEq[{#2}]{passive transformation of representation g}{#1}
Then coordinate transformation in #3 $#4$ gets form
\DrawEq[{#5}{#2}]{coordinate transformation g}{#1}
}

\AddEq[5]{definition: coordinate of geometric object}
{
\begin{ShadedDefinition}
\labelDefinition{coordinate of geometric object, #1}
Orbit
\ShowEq{coordinates of geometric object}{#3}%
\ShowEq{show coordinates of geometric object}{#2}{#3}
is called
\AddIndex{coordinates of geometric object}{coordinates of geometric object}
defined in the #4 $#5$.
For any basis
\ShowEq{geometric object 1}{#2}
corresponding point 
\eqRef{coordinate transformation g}{#1}
of orbit defines
\AddIndex{coordinates of geometric object}{coordinates of geometric object}
relative basis
\ShowEq{geometric object 2}
\qed
\end{ShadedDefinition}
}

\AddEq[1]{remark: 2 views on geometric object}
{
Since a geometric object is an orbit of representation, we see that
according to the theorem
\RefTheorem{proper definition of orbit}
the definition of the geometric object is a proper definition.

Definition
\RefDefinition{coordinate of geometric object, #1}
introduces a geometric object in coordinate space.
We assume in definition \RefDefinition{geometric object, #1}
that we selected a basis of representation $g$.
This allows using a representative of the geometric object
instead of its coordinates.
}

\AddEq[6]{definition: geometric object}
{
\begin{ShadedDefinition}
\labelDefinition{geometric object, #1}
Orbit
\ShowEq{geometric object}{#3}%
\ShowEq{show geometric object}{#2}{#3}%
is called
\AddIndex{geometric object}{geometric object} 
defined in the #4 $#5$.
We also say that $#3$ is
a \AddIndex{geometric object of type $H$}{type of geometric object}.
For any basis
\ShowEq{geometric object 1}{#2}
corresponding point
\eqRef{coordinate transformation g}{#1}
of orbit defines #6
\ShowEq{geometric object 2, g}{#3}%
called
\AddIndex{representative of geometric object}{representative of geometric object}
in the #4 $#5$.
\qed
\end{ShadedDefinition}
}

\AddEq[1]{theorem: invariance principle}
{
\begin{ShadedTheorem}[invariance principle]
\AddIndex{}{invariance principle}
\labelTheorem{invariance principle, #1}
Representative of geometric object does not depend on selection
of basis $\Basis e_f$.
\end{ShadedTheorem}
}

\AddEq[6]{proof: invariance principle}
{
\begin{proof}
To define representative of geometric object,
we need to select basis $\Basis e_f$ of #2 $#3$,
basis $\Basis e_g$ of #2 $#4$
and coordinates of geometric object
\ShowEq{invariance principle 4}{#5}
Corresponding representative of geometric object
has form
\ShowEq{invariance principle 1}{#5}
Suppose we map basis $\Basis e_f$ to basis $\Basis e_{f1}$
by passive transformation
\ShowEq{invariance principle 2 passive}{#6}
According building this forms passive transformation
\eqRef{passive transformation of representation g}{#1}
and coordinate transformation
\eqRef{coordinate transformation g}{#1}.
Corresponding
representative of geometric object has form
\ShowEq{invariance principle 3 algebra}{#5}{#6}
Therefore representative of geometric object
is invariant relative selection of basis.
\end{proof}
}

\AddEq[1]{theorem: passive representation in basis manifold}
{
\begin{ShadedTheorem}
\labelTheorem{passive representation in basis manifold, #1}
There exists single transitive right\Hyph side representation
\ShowEq{passive representation in basis manifold}
\ShowEq{passive representation in basis manifold def}
of group $GA(f)$
in basis manifold $\mathcal B[f]$.
Representation $P(f)$ is called
\AddIndex{passive representation}{passive representation in basis manifold}.
\end{ShadedTheorem}
}

\DefProof{passive representation in basis manifold}
{
Since $A(f)$ is single transitive left\Hyph side representation of group $GA(f)$,
then single transitive right\Hyph side representation $P(f)$
is uniquely defined
according to the theorem
\RefTheorem{two representations of group}.
}

\AddEq[1]{theorem: passive transformation is automorphism of A(f)}
{
\begin{ShadedTheorem}
\labelTheorem{passive transformation is automorphism of A(f), #1}
Passive transformation of the basis manifold
is automorphism of representation $A(f)$.
\end{ShadedTheorem}
}

\DefProof{passive transformation is automorphism of A(f)}
{
The theorem follows from the theorem
\RefTheorem{twin representations of group, automorphism}.
}

\AddEq[3]{theorem: passive transformation}
{
\begin{ShadedTheorem}
\labelTheorem{passive transformation, #1}
Transformation of representation $P(f)$
is called
\AddIndex{passive transformation of the basis manifold}
{passive transformation of basis}
of #2.
We also use notation
\ShowEq{passive transformation of basis}{#3}
\ShowEq{passive transformation of basis def}{#3}
to denote the image of basis $\Basis e$ under passive transformation $#3$.
Passive transformation of basis has form
\DrawEq[#3]{passive transformation}{#1}
\end{ShadedTheorem}
}

\AddEq[1]{proof: passive transformation}
{
\begin{proof}
According to the equality
\eqRef{active transformation}{#1},
active transformation acts from left on coordinates
of basis.
The equality
\eqRef{passive transformation}{#1}
follows from theorems
\RefTheorem{shift is automorphism of representation},
\RefTheorem{two representations of group},
\RefTheorem{twin representations of group, automorphism};
according to these theorems,
passive transformation acts from right on coordinates
of basis.
\end{proof}
}

\DefLabeledDefinitionNote{side representation of group}{\SideNS}
{
Let
\ShowEq{End O2}{A_2}{}
be a multiplicative $\Omega$\Hyph group with product\,\footnotemark
\ShowEq{(fg)->fg}
Let an endomorphism $f$ act on $A_2$\Hyph number $a$ on the \SideNS.
We will use notation
\ShowEq{f(a)=o \SideNS}
Let $A_1$ be multiplicative $\Omega$\Hyph group with product
\ShowEq{(ab)->ab}
We call a homomorphism of multiplicative $\Omega$\Hyph group
\DrawEq{representation of group, map}{\SideNS-side}
\AddIndex{\SideNS\Hyph side representation}{\SideNS-side representation}
of multiplicative $\Omega$\Hyph group $A_1$
or
\AddIndex{\SideNS\Hyph side $A_1$\Hyph representation}{\SideNS-side A representation}
in $\Omega_2$\Hyph algebra $A_2$ if the map $f$ holds
\DrawEq{\SideNS-side representation of group}{def}
We identify
an $A_1$\Hyph number $a_1$ and its image $f(a_1)$
and write \SideNS\Hyph side transformation
caused by $A_1$\Hyph number $a_1$
as
\ShowEq{\SideNS-side representation}
In this case, the equality
\eqRef{\SideNS-side representation of group}{def}
gets following form
\DrawEq{* \SideNS-side representation of group}{def}
The map
\ShowEq{A12->A2 1 \SideNS}
generated by \SideNS\Hyph side representation $f$ is called
\AddIndex{\SideNS\Hyph side product}{\SideNS-side product}
of $A_2$\Hyph number $a_2$ over $A_1$\Hyph number $a_1$.
}
{
Very often a product
in multiplicative $\Omega$\Hyph group
\ShowEq{End O2}{A_2}{}
is superposition of endomorphisms
\ShowEq{Act=circ}
However, a product
in multiplicative $\Omega$\Hyph group
\ShowEq{End O2}{A_2}{}
may be different from superposition of endomorphisms.
}

\DefDefinition{reduced morphism of representations}
{
Let
\ShowEq{f:A->*B}f{A_1}{A_2}
be representation of $\Omega_1$\Hyph algebra $A_1$
in $\Omega_2$\Hyph algebra $A_2$ and
\ShowEq{f:A->*B}g{A_1}{B_2}
be representation of $\Omega_1$\Hyph algebra $A_1$
in $\Omega_2$\Hyph algebra $B_2$.
Let
\ShowEq{id:A1->A1 A2->B2}
be morphism of representations.
In this case we identify morphism
\ShowEq{map r,R}{\id}{r_2}{}
of representations of $\Omega_1$\Hyph algebra and corresponding homomorphism $r_2$ of $\Omega_2$\Hyph algebra
and the homomorphism $r_2$ is called
\AddIndex{reduced morphism of representations}
{reduced morphism of representations}.
We will use diagram
\ShowEq{morphism id,R of representations}
to represent reduced morphism $r_2$ of representations of
$\Omega_1$\Hyph algebra.
From diagram it follows
\ShowEq{morphism of representations of universal algebra}
We also use diagram
\ShowEq{morphism id,R of representations 2}
instead of diagram
\EqRef{morphism id,R of representations}.
}

\DefDefinition{action of rational integers in Abelian group}
{
The action of ring of rational integers $Z$
in Abelian group $G$
is defined using following rules
\ShowEq{action of rational integers in Abelian group}
}

\DefTheorem{action of ring of rational integers in Abelian group}
{
The action of ring of rational integers $Z$
in Abelian group $G$
defined in the definition
\RefDefinition{action of rational integers in Abelian group}
is representation.
The following equalities are true
\ShowEq{representation of Z in G}
}

\DefProof{action of ring of rational integers in Abelian group}
{
The equality
\EqRef{1a=a}
follows from the equality
\EqRef{0g=0}
and from the equality
\EqRef{(n+1)g=}
when $n=0$.

\ShowEq{step of proof: action in Abelian group}{(m+n)a=ma+na}

The equality
\ShowEq{(k+n)a-na=ka}
follows from the equality
\EqRef{(m+n)a=ma+na}.
The equality
\EqRef{(m-n)a=ma-na}
follows from the equality
\EqRef{(k+n)a-na=ka},
if we assume
$m=k+n$, $k=m-n$.

\ShowEq{step of proof: action in Abelian group}{(nm)a=n(ma)}

\ShowEq{step of proof: action in Abelian group}{n(a+b)=na+nb}

From the equality
\EqRef{n(a+b)=na+nb},
it follows that the map
\ShowEq{phi(n):G->G}
is an endomorphism of Abelian group $G$.
From equalities
\EqRef{(m+n)a=ma+na},
\EqRef{(nm)a=n(ma)},
it follows that the map
\ShowEq{phi:Z->End(G)}
is a homomorphism of the ring $Z$.
According to the definition
\RefDefinition[\RefRepresentation]{representation of algebra},
the map $\varphi$ is representation
of ring of rational integers $Z$
in Abelian group $G$.
}

\DefDefinitionNote{tensor product of representations}
{
Let $A$ be Abelian multiplicative $\Omega_1$\Hyph group.
Let
\ShowEq{a1n}An{}
be $\Omega_2$\Hyph algebras.\,\footnotemark
Let, for any
\ShowEq{k,1n}
\ShowEq{representation A Ak 1}
be effective representation of multiplicative $\Omega_1$\Hyph group $A$
in $\Omega_2$\Hyph algebra $A_k$.
Consider category $\mathcal A$ whose objects are
reduced polymorphisms of representations $f_1$, ..., $f_n$
\ShowEq{polymorphisms category}
where $S_1$, $S_2$ are $\Omega_2$\Hyph algebras and
\ShowEq{representation of algebra in S1 S2}
are effective representations of multiplicative $\Omega_1$\Hyph group $A$.
We define morphism
\ShowEq{r1->r2}
to be reduced morphism of representations $h:S_1\rightarrow S_2$
making following diagram commutative
\ShowEq{polymorphisms category, diagram}
Universal object
\ShowEq{tensor product of representations}
of category $\mathcal A$ is called
\AddIndex{tensor product}{tensor product}
of representations
\ShowEq{a1n}An.
}
{
I give definition
of tensor product of representations of universal algebra
following to definition in \citeBib{Serge Lang}, p. 601 - 603.
}

\DefEq
{
\begin{ShadedTheorem}
\labelTheorem{tensor product of representations}
Since there exists polymorphism of representations,
then there exists tensor product of representations.
\end{ShadedTheorem}
}
{theorem: tensor product of representations}

\DefEq
{
\begin{ShadedTheorem}
\labelTheorem{representation, tensor product}
Let
\ShowEq{equivalence, 1, representation, tensor product}
Tensor product is distributive over operation $\omega$
\ShowEq{tensors 1, representation, tensor product}
The representation of multiplicative $\Omega_1$\Hyph group $A$
in tensor product is defined by equality
\ShowEq{tensors 2, representation, tensor product}
\end{ShadedTheorem}
}
{theorem: representation, tensor product}

\DefEq
{
\begin{ShadedTheorem}
\labelTheorem{tensor product and polymorphism}
Let $B_1$, ..., $B_n$ be
$\Omega_2$\Hyph algebras.
Let
\ShowEq{map f, 1, representation, tensor product}
be reduced polymorphism defined by equality
\ShowEq{map f, representation, tensor product}
Let
\ShowEq{map g, representation, tensor product}
be reduced polymorphism into $\Omega$\Hyph algebra $V$.
There exists morphism of representations
\ShowEq{map h, representation, tensor product}
such that the diagram
\ShowEq{map gh, representation, tensor product}
is commutative.
\end{ShadedTheorem}
}
{theorem: tensor product and polymorphism}

\DefEq
{
\begin{ShadedTheorem}
\labelTheorem{B times->B otimes}
The map
\ShowEq{B times->B otimes}
is polymorphism.
\end{ShadedTheorem}
}
{theorem: B times->B otimes}

\DefTheorem{|a-b|>|a|-|b|}
{
Let $A$ be normed $\Omega$\Hyph group.
Then
\ShowEq{|a-b|>|a|-|b|}
}

\DefDefinitionNote{free representation of algebra}
{
The representation
\ShowEq{f:A->*B}f{A_1}{A_2}
of the $\Omega_1$\Hyph algebra $A_1$ is called
\AddIndex{free}{free representation},\,\footnotemark
if the statement
\ShowEq{faa=fba}
for any
\ShowEq{a in A}2{}
implies that
\ShowEq{a=b}.
}
{
See similar definition of free representation of group in
\citeBib{Postnikov: Differential Geometry}, page 16.
\ePrints{1908.04418,6860-2955}
\ifx\Semafor\ValueOn
See also the theorem
\RefTheorem{free representation of group}.
\fi
}

\DefDefinition{continuous map, Omega group}
{
A map
\DrawEq[f{A_1}{A_2}{}]{f: A->B}{}
of normed $\Omega_1$\Hyph group $A_1$ with norm $\|x\|_1$
into normed $\Omega_2$\Hyph group $A_2$ with norm $\|y\|_2$
is called \AddIndex{continuous}{continuous map}, if
for every as small as we please $\epsilon>0$
there exist such $\delta>0$, that
\ShowEq{|x'-x|<delta}
implies
\ShowEq{|f(x)-f(x')|<epsilon}
}

\DefTheorem{continuous map, open set}
{
A map
\DrawEq[f{A_1}{A_2}{}]{f: A->B}{}
of normed $\Omega_1$\Hyph group $A_1$ with norm $\|x\|_1$
into normed $\Omega_2$\Hyph group $A_2$ with norm $\|y\|_2$
is continuous, iff
preimage of an open set
is the open set.
}

\DefTheorem{image of interval is interval, real field}
{
Let
\ShowEq{f:A->B}fRR
be continuous map of real field.
Then image of interval is interval.
}

\DefEq
{
\begin{ShadedDefinition}
\labelDefinition{norm of operation}
Let $A$
be normed $\Omega$\Hyph group.
For $n$\Hyph ary operation $\omega$, the value
\ShowEq{norm of operation}
\ShowEq{norm of operation, definition}
is called
\AddIndex{norm of operation}{norm of operation} $\omega$.
\qed
\end{ShadedDefinition}
}
{definition: norm of operation}

\DefTheorem{|fx|<|f||x|1n}
{
Let $A$ be normed $\Omega$\Hyph group.
For $n$\Hyph ary operation $\omega$,
\ShowEq{|a omega|<|omega||a|1n}
}

\DefEq
{
\begin{ShadedDefinition}
\labelDefinition{norm of representation}
Let
\ShowEq{f:A->*B}g{A_1}{A_2}
be representation
\ePrints{1305.4547}
\ifx\Semafor\ValueOn
\footnote{
See the definition
\ShowEq{ref definition: left-side representation of algebra}{}
of the representation of universal algebra.
According to the definition
\xRef[1211.6965]{definition: module over ring},
module is representation of ring in Abelian group.
Since ring and Abelian group are $\Omega$\Hyph groups,
then module is representation of $\Omega$\Hyph group.
}
\fi
of $\Omega_1$\Hyph group $A_1$ with norm $\|x\|_1$
in $\Omega_2$\Hyph group $A_2$ with norm $\|x\|_2$.
The value
\ShowEq{norm of representation}
\ShowEq{norm of representation, definition}
is called
\AddIndex{norm of representation}{norm of representation} $f$.
\qed
\end{ShadedDefinition}
}
{definition: norm of representation}

\DefTheorem{|fab|<|f||a||b|}
{
Let
\ShowEq{f:A->*B}g{A_1}{A_2}
be representation
of $\Omega_1$\Hyph group $A_1$ with norm $\|x\|_1$
in $\Omega_2$\Hyph group $A_2$ with norm $\|x\|_2$.
Then
\ShowEq{|fab|<|f||a||b|}
}

\DefEq
{
\begin{ShadedDefinition}
\labelDefinition{open set}
Let $A$ be normed $\Omega$\Hyph group.
A set
$U\subset A$
is called
\AddIndex{open}{open set},\,\footnote{
In topology,
we usually define an open set before we
define base of topology.
In the case of a metric or normed space,
it is more convenient to define an open set,
based on the definition of base of topology.
In such case, the definition is based on one of the properties of
base of topology. An immediate proof
allows us to see that defined such
an open set satisfies the basic properties.
}
if, for any $A$\Hyph number
\ShowEq{a in U},
there exists
\ShowEq{epsilon in R}
such that
\ShowEq{B(a) subset U}
\qed
\end{ShadedDefinition}
}
{definition: open set}

\DefEq
{
\begin{ShadedDefinition}
\labelDefinition{compact set}
A set $T$ of topological space is called
\AddIndex{compact}{compact set},
if every open cover of $T$ has
finite subcover.\,\footnote{
See also the definition
\citeBib{Kolmogorov Fomin}-1, page 92.
}
\qed
\end{ShadedDefinition}
}
{definition: compact set}

\DefTheorem{c1+c2 in B}
{
Let $A$ be normed $\Omega$\Hyph group.
For
\ShowEq{c1 c2 in A}
let
\ShowEq{c1 c2 in B}
Then
\ShowEq{c1+c2 in B}
}

\DefLabeledDefinition{limit of sequence}{normed \Algebrac}
{
Let $\Module$
be normed $\Base$\Hyph \Algebrab.
$\Module$\Hyph number $a$ is called 
\AddIndex{limit of a sequence}{limit of sequence}
\ShowEq{[an] in A}
\ShowEq{limit of sequence}
if for any
\ShowEq{epsilon in R}
there exists positive integer $n_0$ depending on $\epsilon$ and such, that
\DrawEq{|an-a|<epsilon}{-}
for every $n>n_0$.
We also say that
\AddIndex{sequence $a_n$ converges}{sequence converges}
to $a$.
}

\AddEq{theorem: limit of sequence}
{
\begin{ShadedTheorem}
\labelTheorem{limit of sequence, normed \Algebrac}
Let $\Module$
be normed $\Base$\Hyph \Algebrab.
$\Module$\Hyph number $a$ is
limit of a sequence
\ShowEq{[an] in A}
\DrawEq{a=lim an}{}
if for any
\ShowEq{epsilon in R}
there exists positive integer $n_0$ depending on $\epsilon$ and such, that
\DrawEq{an in B(a)}{}
for every $n>n_0$.
\end{ShadedTheorem}
}

\DefProof{limit of sequence}
{
The theorem follows from definitions
\RefDefinition{open ball},
\refDefinition{limit of sequence}{normed \Algebrac}.
}

\DefLabeledDefinition{fundamental sequence}{normed \Algebrac}
{
Let $\Module$
be normed $\Base$\Hyph \Algebrab.
The sequence
\ShowEq{[an] in A}
is called 
\AddIndex{fundamental}{fundamental sequence}
or \AddIndex{Cauchy sequence}{Cauchy sequence},
if for every
\ShowEq{epsilon in R}
there exists positive integer $n_0$ depending on $\epsilon$ and such, that
\DrawEq{|ap-aq|<epsilon}{-}
for every $p$, $q>n_0$.
}

\AddEq{theorem: fundamental sequence}
{
\begin{ShadedTheorem}
\labelTheorem{fundamental sequence, normed \Algebrac}
Let $\Module$
be normed $\Base$\Hyph \Algebrab.
The sequence
\ShowEq{[an] in A}
is fundamental sequence,
iff for every
\ShowEq{epsilon in R}
there exists positive integer $n_0$ depending on $\epsilon$ and such, that
\ShowEq{aq in B(ap)}
for every $p$, $q>n_0$.
\end{ShadedTheorem}
}

\DefProof{fundamental sequence}
{
The theorem follows from definitions
\RefDefinition{open ball},
\refDefinition{fundamental sequence}{normed \Algebrac}.
}

\DefTheorem{norm of compact set is bounded}
{
Let $C$ be compact set of
normed $\Omega$\Hyph group $A$.
Then the norm
\ShowEq{|x in C|}
is bounded from both sides.
}

\DefTheorem{lim a=lim b}
{
Let $A$ be normed $\Omega$\Hyph group.
Let $a_n$, $b_n$, $n=1$, ..., be fundamental sequences.
Let
\DrawEq{lim a-b=0}{limit}
If the sequence $a_n$ converges, then
the sequence $b_n$ converges and
\ShowEq{lim a=lim b}
}

\DefEq
{
\begin{ShadedDefinition}
\labelDefinition{complete Omega group}
Normed $\Omega$\Hyph group $A$ is called
\AddIndex{complete}{complete Omega group}
if any fundamental sequence of elements
of $\Omega$\Hyph group $A$ converges, i.e.
has limit in $\Omega$\Hyph group $A$.
\qed
\end{ShadedDefinition}
}
{definition: complete Omega group}

\DefEq
{
\begin{ShadedDefinition}
\labelDefinition{limit of sequence, map to Omega group}
Let
\ShowEq{fn M(X,A)}
be sequence of maps into normed
$\Omega$\Hyph group $A$.
The map
\ShowEq{f M(X,A)}
is called
\AddIndex{limit of sequence}{limit of sequence}
$f_n$, if for any $x\in X$
\DrawEq{f(x)=lim}{}
We also say that
\AddIndex{sequence $f_n$ converges}{sequence converges}
to the map $f$.
\qed
\end{ShadedDefinition}
}
{definition: limit of sequence, map to Omega group}

\DefEq
{
\begin{ShadedDefinition}
\labelDefinition{sequence converges uniformly}
Let
\ShowEq{fn M(X,A)}
be sequence of maps into normed
$\Omega$\Hyph group $A$.
\AddIndex{Sequence $f_n$ converges uniformly}
{sequence converges uniformly}
to the map $f$,
if for any
\ShowEq{epsilon in R}
there exists $N$ such that
\ShowEq{fn(x) - f(x)}
for any $n>N$.
\qed
\end{ShadedDefinition}
}
{definition: sequence converges uniformly}

\DefTheorem{sequence converges uniformly, fn-fm}
{
Sequence of maps
\ShowEq{fn M(X,A)}
into normed
$\Omega$\Hyph group $A$ converges uniformly
to the map $f$,
if for any
\ShowEq{epsilon in R}
there exists $N$ such that
\ShowEq{|fn(x)-fm(x)|<e}
for any \(n\), \(m>N\).
}

\DefDefinition{quasibasis of representation}
{
Let
\ShowEq{f:A->*B}f{A_1}{A_2}
be representation
of $\Omega_1$\Hyph algebra $A_1$
in $\Omega_2$\Hyph algebra $A_2$
and
\ShowEq{Gen f =}
If, for the set $X\subset A_2$, it is true that
\ShowEq{in Gen}Xf,
then for any set $Y$, $X\subset Y\subset A_2$,
also it is true that
\ShowEq{in Gen}Yf.
If there exists minimal set
\ShowEq{in Gen}Xf,
then the set $X$ is called
\AddIndex{quasibasis}{quasibasis}
of representation $f$.
}

\DefTheorem{X is quasibasis of representation}
{
If the set $X$ is the quasibasis of the representation $f$,
then, for any $m\in X$,
the set $X\setminus\{m\}$ is not
generating set of the representation $f$.
}

\DefRemark{X is quasibasis of representation}
{
The proof of the theorem
\RefTheorem{X is quasibasis of representation}
gives us effective method for constructing the quasibasis of the representation $f$.
Choosing an arbitrary generating set, step by step, we
remove from set those elements which have coordinates
relative to other elements of the set.
If the generating set of the representation is infinite,
then this construction may not have the last step.
If the representation has finite generating set,
then we need a finite number of steps to construct a quasibasis of this representation.
}

\AddEq{remark: representation in form of Omega2-word is ambiguous}
{
We introduced $\Omega_2$\Hyph word of $x\in A_2$ relative generating set $X$
in the definition
\RefDefinition{word of element relative to set, representation}.
From the theorem
\RefTheorem{X is quasibasis of representation},
it follows that if the generating set $X$ is not an quasibasis,
then a choice of $\Omega_2$\Hyph word relative generating set $X$ is ambiguous.
However, even if the generating set $X$ is an quasibasis,
then a representation of $m\in A_2$ in form of $\Omega_2$\Hyph word is ambiguous.
}

\DefDefinitionNote{word of element relative to set, representation}
{
Let $X\subset A_2$.
For each
\ShowEq{x in JX}m{}
there exists
\AddIndex{$\Omega_2$\Hyph word}
{word of element relative to generating set, representation}
defined according to following rules.
\ShowEq{word of element relative to generating set, representation}
\StartLabelItem
\begin{enumerate}
\item If $m\in X$, then $m$ is $\Omega_2$\Hyph word.
\labelItem{word of element relative to set, representation, m in X}
\item If $m_1$, ..., $m_n$ are $\Omega_2$\Hyph words and
$\omega\in\Omega_2(n)$, then $m_1...m_n\omega$
is $\Omega_2$\Hyph word.
\labelItem{word of element relative to set, representation, omega}
\item If $m$ is $\Omega_2$\Hyph word and
\ShowEq{a in A1},
then $f(a)(m)$
is $\Omega_2$\Hyph word.
\labelItem{word of element relative to set, representation, am}
\end{enumerate}
We will identify an element
\ShowEq{x in JX}m{}
and corresponding it $\Omega_2$\Hyph word using equation
\ShowEq{identify element jf(X) and word}
Similarly, for an arbitrary set
\ShowEq{subset of representation}{}
we consider the set of $\Omega_2$\Hyph words\,\footnotemark
\ShowEq{subset of words of representation}
We also use notation
\ShowEq{subset of words of representation, 1}
Denote
\ShowEq{word set of representation}
\AddIndex{the set of $\Omega_2$\Hyph words of representation $J[f,X]$}
{word set of representation}.
}{
The expression $\wXm$
is a special case
of the expression $\wXm[B]$, namely
\ShowEq{subset of words of representation, 2}
}

\DefRemark{three reasons of ambiguity in Omega2-word}
{
There are three reasons of ambiguity in notation of
$\Omega_2$\Hyph word.
\StartLabelItem
\begin{enumerate}
\item
In $\Omega_i$\Hyph algebra $A_i$, $i=1$, $2$,
equalities may be defined.
For instance, if $e$ is unit of multiplicative group $A_i$,
then the equality
\[ae=a\]
is true for any $a\in A_i$.
\item
\begin{sloppypar}
Ambiguity of choice of $\Omega_2$\Hyph word
may be associated with properties of representation.
For instance, if $m_1$, ..., $m_n$ are $\Omega_2$\Hyph words,
$\omega\in\Omega_2(n)$ and
\ShowEq{a in A1},
then\,\footnote{
For instance, let
$\{e_1,e_2\}$
be the basis of vector space over field $k$.
The equation \EqRef{ambiguity of coordinates 1}
has the form of distributive law
\ShowEq{ambiguity of coordinates 1, vector}
}
\ShowEq{ambiguity of coordinates 1}
At the same time, if $\omega$ is operation of
$\Omega_1$\Hyph algebra $A_1$ and operation of
$\Omega_2$\Hyph algebra $A_2$, then we require
that $\Omega_2$\Hyph words
\ShowEq{a1n omega x}
and
\ShowEq{a1n x omega}
describe the same
element of $\Omega_2$\Hyph algebra $A_2$.\,\footnote{For vector space,
this requirement has the form of distributive law
\ShowEq{(a+b)e=ae+be}
}
\DrawEq{ambiguity of coordinates 2}{representation}
\end{sloppypar}
\item
Equalities like
\EqRef{ambiguity of coordinates 1},
\eqRef{ambiguity of coordinates 2}{representation}
persist under morphism of representation.
Therefore we can ignore this
form of ambiguity of $\Omega_2$\Hyph word.
However, a fundamentally different form of ambiguity is possible.
We can see an example of such ambiguity
in theorems
\refTheorem[\RefRepresentation]{linear dependence between vectors of basis}{\SideWS module},
\refTheorem[\RefRepresentation]{coordinates of vector with linear dependence}{\SideWS module}.
\end{enumerate}
So we see that we can define different equivalence relations
on the set of $\Omega_2$\Hyph words.\,\footnote{
Evidently each of the equalities
\EqRef{ambiguity of coordinates 1},
\eqRef{ambiguity of coordinates 2}{representation}
generates some equivalence relation.
}
Our goal is to find a maximum equivalence
on the set of $\Omega_2$\Hyph words which
persist under morphism of representation.
}

\DefTheoremNote{equivalence generated by basis}
{
Let $X$ be quasibasis of the representation
\ShowEq{f:A->*B}f{A_1}{A_2}
For an generating set $X$ of representation $f$,
there exists equivalence
\ShowEq{l fX in w fX}
which is generated exclusively by the following statements.
\StartLabelItem
\begin{enumerate}
\item
If in $\Omega_2$\Hyph algebra $A_2$ there is an equality
\ShowEq{w1[]=w2[]}
defining structure of $\Omega_2$\Hyph algebra, then
\ShowEq{w1 w2 in l fX}
\item
If in $\Omega_1$\Hyph algebra $A_1$ there is an equality
\ShowEq{w1[]=w2[]}
defining structure of $\Omega_1$\Hyph algebra, then
\ShowEq{fw1m fw2m in l fX}
\item
For any operation
\ShowEq{omega n ari}{}1n,
\ShowEq{(f(a1n)a2,f(a)a2) in l}
\item
For any operation
\ShowEq{omega n ari}{}2n,
\ShowEq{(fa1(a21n),fa1(a2)) in l}
\item
Let
\ShowEq{omega in O1AO2}12.
If the representation $f$ satisfies equality\,\footnotemark
\ShowEq{f(a1n)a2=f(a)a2}
then we can assume that the following equality is true
\ShowEq{(f(a1n)a2,f(a)a2) in l 1}
\end{enumerate}
}
{
Consider a representation of commutative ring $D$
in $D$\Hyph algebra $A$. We will use notation
\ShowEq{f(a)(v)=av}
The operations of addition and multiplication are defined in both algebras.
However the equality
\ShowEq{f(a+b)v=}
is true, and the equality
\ShowEq{f(ab)v=}
is wrong.
}

\DefDefinition{basis of representation}
{
Quasibasis $\Basis e$ of the representation $f$ such that
\ShowEq{rho=lambda}
is called
\AddIndex{basis of representation}{basis of representation} $f$.
}

\DefDefinition{equivalence on the set of Omega2-words}
{
Generating set $X$ of representation $f$
generates equivalence
\ShowEq{rho fX=}
on the set of $\Omega_2$\Hyph words.
}

\DefTheorem{h=f+g, converges uniformly}
{
Let sequence of maps
\ShowEq{fn M(X,A)}
into complete $\Omega$\Hyph group $A$
converge uniformly
to the map $f$.
Let sequence of maps
\ShowEq{gn M(X,A)}
into complete $\Omega$\Hyph group $A$
converge uniformly
to the map $g$.
Then sequence of maps
\ShowEq{hn=fn+gn}
into complete $\Omega$\Hyph group $A$
converges uniformly
to the map
\DrawEq{h=f+g}{converges uniformly}
}

\DefTheorem{representation f generates representation fX}
{
The representation
\ShowEq{f:A->*B}g{A_1}{A_2}
of $\Omega_1$\Hyph group $A_1$ with norm $\|x\|_1$
in $\Omega_2$\Hyph group $A_2$ with norm $\|x\|_2$
generates representation
\ShowEq{f:A->*B}{f_X}{M(X,A_1)}{M(X,A_2)}
of $\Omega_1$\Hyph group
\ShowEq{M(X,A1)}
in $\Omega_2$\Hyph group
\ShowEq{M(X,A2)}
where (
\ShowEq{g1 in M(X,A1)},
\ShowEq{g2 in M(X,A2)}
)
\ShowEq{f(x)g(x):X->A2}
}

\DefTheorem{fX(g1)(g2), converges uniformly}
{
Let
\ShowEq{f:A->*B}g{A_1}{A_2}
be representation
of complete $\Omega_1$\Hyph group $A_1$ with norm $\|x\|_1$
in complete $\Omega_2$\Hyph group $A_2$ with norm $\|x\|_2$.
Let sequence of maps
\ShowEq{g1n M(X,A1)}
converge uniformly
to the map $g_1$.
Let sequence of maps
\ShowEq{g2n M(X,A2)}
converge uniformly
to the map $g_2$.
Let range of the map $g_i$, $i=1$, $2$,
be compact set.
Then the sequence of maps
\ShowEq{f(g1n)(g2n)}
converge uniformly
to the map
\ShowEq{fX(g1)(g2)}.
}

\DefTheorem{set of maps to Omega group}
{
Let
\ShowEq{set of maps to Omega group}
be set of maps of a set $X$ to $\Omega$\Hyph group $A$.
We can define the structure of $\Omega$\Hyph group on the set
\EqParm{M(X,A)}{=.}
}

\DefEq
{
Since $X$ is an arbitrary set, we cannot define
norm in $\Omega$\Hyph group
\EqParm{M(X,A)}{=.}
However we can define the convergence of the sequence in
\EqParm{M(X,A)}{=.c}
therefore, we can define a topology in
\EqParm{M(X,A)}{=.}
}
{remark: set of maps to Omega group}

\DefDefinition{closed ball}
{
Let $\Module$
be normed $\Base$\Hyph \Algebrab.
Let $a\in \Module$.
The set
\ShowEq{closed ball}
is called
\AddIndex{closed ball}{closed ball}
with center at $a$.
}

\DefDefinition{open ball}
{
Let $\Module$
be normed $\Base$\Hyph \Algebrab.
Let $a\in \Module$.
The set
\ShowEq{open ball}
is called
\AddIndex{open ball}{open ball}
with center at $a$.
}

\DefDefinition{reduced polymorphism of representations}
{
Let
\ShowEq{set of universal algebras 1}
be universal algebras.
Let, for any
\ShowEq{k,1n}
\ShowEq{f:A->*B}{f_k}A{B_k}
be effective representation of $\Omega_1$\Hyph algebra $A$
in $\Omega_2$\Hyph algebra $B_k$.
Let
\ShowEq{f:A->*B}fAB
be effective representation of $\Omega_1$\Hyph algebra $A$
in $\Omega_2$\Hyph algebra $B$.
The map
\ShowEq{reduced polymorphism of representation}
is called
\AddIndex{reduced polymorphism of representations}
{reduced polymorphism of representations}
$f_1$, ..., $f_n$ into representation $f$,
if, for any
\ShowEq{k,1n}
provided that all variables except the variable $x_k\in B_k$
have given value, the map $r_2$
is a reduced morphism of representation $f_k$ into representation $f$.

If $f_1=...=f_n$, then we say that the map $r_2$
is reduced polymorphism of representation $f_1$ into representation $f$.

If $f_1=...=f_n=f$, then we say that the map $r_2$
is reduced polymorphism of representation $f$.
}

\DefConvention{Einstein summation}
{
We will use Einstein summation convention.
When an index is present in
an expression twice (one above and one below) and a set of index is known,
we have the sum with respect to repeated index.
In this case we assume that we know the set
of summation index and do not use summation symbol
\ShowEq{av=sum av}av
If needed to clearly
show set of index, I will do it.
}

\DefTheorem{structure of Abelian group}
{
Let $G$ be Abelian group.
The set of $G$\Hyph numbers generated by the set
\ShowEq{S=si}
has form
\ShowEq{J(S)=gi si}
where the set
\ShowEq{|gi ne 0|}g
is finite.
}

\DefProof{structure of Abelian group}
{
We prove the theorem by induction based on the theorems
\citeBib{Cohn: Universal Algebra}\Hyph 5.1, page 79 and
\RefTheorem[\RefRepresentation]{structure of subrepresentations}.

For any $s_k\in S$,
let
\ShowEq{gi=dik}
Then
\ShowEq{sk=sum si}
\ShowEq{sk in JS}
follows from
\EqRef{J(S)=gi si},
\EqRef{sk=sum si}.

Let
\ShowEq{w12 in Jv}gS
Since $G$ is Abelian group,
then, according to the statement
\RefItem[\RefRepresentation]{x1n omega in Xk+1},
\ShowEq{w1+w2 in Jv}gS
According to the equality
\EqRef{J(S)=gi si},
there exist $Z$\Hyph numbers
\ShowEq{gi12}g
such that
\DrawEq[g]{w12= }{group}
where sets
\DrawEq[g]{|ci12 ne 0|}{group}
are finite.
From the equality
\eqRef{w12= }{group},
it follows that
\DrawEq[g]{w1+w2= }{group}
The equality
\DrawEq[g]{w1+w2= 1 }{group}
follows from equalities
\EqRef{(m+n)a=ma+na},
\eqRef{w1+w2= }{group}.
From the equality
\eqRef{|ci12 ne 0|}{group},
it follows that
the set
\ShowEq{|gi 1+2 ne 0|}g
is finite.
From the equality
\eqRef{w1+w2= 1 }{group},
it follows that
\ShowEq{w1+w2 in Jv}gS
}

\AddEq [1]{step of proof: action in Abelian group}
{
From the equality
\EqRef{0g=0},
it follows that the equality
\EqRef{#1}
is true when $n=0$.
\def\TempA{(nm)a=n(ma)}
\def\TempB{#1}
\begin{itemize}
\item
Let the equality
\EqRef{#1}
is true when $n=k\ge 0$.
Then
\ShowEq{#1,n=k}
The equality
\ShowEq{#1,n=k+1}
follows from
\ifx\TempA\TempB
equalities
\EqRef{(n+1)g=},
\EqRef{(m+n)a=ma+na}.
\else
the equality
\EqRef{(n+1)g=}.
\fi
Therefore, the equality
\EqRef{#1}
is true when $n=k+1$.
According to mathematical induction,
the equality
\EqRef{#1}
is true for any $n\ge 0$.
\item
Let the equality
\EqRef{(m+n)a=ma+na}
is true when $n=k\le 0$.
Then
\ShowEq{#1,n=k}
The equality
\ShowEq{#1,n=k-1}
follows from
\ifx\TempA\TempB
equalities
\EqRef{(n-1)g=},
\EqRef{(m-n)a=ma-na}.
\else
the equality
\EqRef{(n-1)g=}.
\fi
Therefore, the equality
\EqRef{#1}
is true when $n=k-1$.
According to mathematical induction,
the equality
\EqRef{#1}
is true for any $n\le 0$.
\item
Therefore, the equality
\EqRef{#1}
is true for any $n\in Z$.
\end{itemize}
}

\DefDefinition{category of representations A1(mA2)}
{
Let $A_1$ be $\Omega_1$\Hyph algebra.
Let $\mathcal A_2$ be category of $\Omega_2$\Hyph algebras.
We define \AddIndex{category
\ShowEq{A1(mA2) symb}
of representations}
{category of representations}
of $\Omega_1$\Hyph algebra $A_1$ in $\Omega_2$\Hyph algebra.
Representations of $\Omega_1$\Hyph algebra $A_1$ in $\Omega_2$\Hyph algebra
are objects of this category.
Reduced morphisms of corresponding representations are morphisms of this category.
}

\DefProofRef{cr transpose}{}
{
We can prove \eqRef{cr transpose}{Theorem}
in case of matrices the same way as we
proved \eqRef{rc transpose}{Theorem}. However it is
more important for us to show that \eqRef{cr transpose}{Theorem}
follows directly from \eqRef{rc transpose}{Theorem}.

Applying \EqRef{double transpose} to each term in left side of
\eqRef{cr transpose}{Theorem} we get
\ShowEq{cr transpose, 1}
From \EqRef{cr transpose, 1} and \eqRef{rc transpose}{Theorem}
it follows that
\ShowEq{cr transpose, 2}
\eqRef{cr transpose}{Theorem} follows from \EqRef{cr transpose, 2}
and \EqRef{double transpose}.
}

\DefProofRef{rc transpose}{}
{
The chain of equalities
\ShowEq{rc transpose, 1}
follows from \EqRef{transpose of matrix, 1},
\EqRef{rc-product of matrices} and
\EqRef{cr-product of matrices}.
The equality \eqRef{rc transpose}{Theorem} follows from
\EqRef{rc transpose, 1}.
}

\DefText{matrix operations 1}
{
According to the custom the product of
matrices $a$ and $b$ is defined as product of
rows of the matrix $a$ and columns of the matrix $b$.
\ePrints{1908.04418,6860-2955,2020.06.01,2204.06320,2022.01.05}
\Items{8428-0408,2207.06506}%
\ifx\Semafor\ValueOn
In non\Hyph commutative algebra, this product
is not enough to solve some problems.
\ePrints{8428-0408,2207.06506}
\ifx\Semafor\ValueOff
\ShowExample{matrix in right vector space}

\ShowExample{matrix in left vector space}

From examples
\RefExample{matrix in right vector space},
\RefExample{matrix in left vector space},
\else
\TwoColText
{
\ShowEq{def left}
\ShowEq{def AVector}%
\ShowEq{\DefCol}
\ShowExample{matrix in vector space}
}
{
\ShowEq{def right}
\ShowEq{def AVector}%
\ShowEq{\DefCol}
\ShowExample{matrix in vector space}
}

From examples
\refExample{matrix in vector space}{left},
\refExample{matrix in vector space}{right}
\fi
it follows that
\else
If the product in
\ePrints{5284-0163,1801.01628}
\ifx\Semafor\ValueOn
$D$\Hyph algebra
\else
$\Omega$\Hyph ring
\fi
is non\Hyph commutative,
\fi
we cannot confine ourselves to traditional product of matrices
and we need to define two products of matrices.
To distinguish between these products we introduced
a new notation.
}

\DefText{matrix operations 2}
{
We also consider following operations on the set of matrices.
}

\DefExample{matrix in left vector space}
{
We represent the basis $\Basis e$
of left vector space $V$
\ePrints{2020.06.01,2204.06320,2022.01.05}
\ifx\Semafor\ValueOff
over $D$\Hyph algebra $A$
(see the definition
\ShowEq{def left}%
\ShowEq{def AVector}%
\refDefinition{module over associative algebra}{\SideWS \VectorSetNS}
and the theorem
\refTheorem{basis over division algebra}{left})
\fi
as row of matrix
\DrawEq[en]{a=(a1.n row)}{left basis}
We represent coordinates of vector
\ShowEq{w=w*e left-cols}v{}
as vector column
\DrawEq[vn]{a=(a1.n col)}{left basis}
However, we cannot represent the vector $\Vector v$
as product of matrices
\ShowEq{v=(v)(e)}
because this product is not defined.
\ePrints{2020.06.01,2204.06320,2022.01.05}
\ifx\Semafor\ValueOff
We represent homomorphism
of left vector space $V$
using matrix
\ShowEq{v'=vf}
We cannot express the equality
\EqRef{v'=vf}
as traditional product of matrices $v$ and $f$.
\fi
}

\AddEq{define basis for example}
{
We represent the set of vectors
\ShowEq{\RowN() 1n}en{}
of basis $\Basis e$
of \SideWS vector space $V$
over $D$\Hyph algebra $A$
(see the definition
\ShowRef{\SideWS module over algebra}
and the theorem
\refTheorem{basis over division algebra}{left})
as
\RowNWS of matrix
\DrawEq[en]{a=(a1.n \RowN)}{\SideWS basis}
We represent coordinates
\ShowEq{\ColNS() 1n}vn{}
of vector
\ShowEq{w=w*e left-cols}v{}
as \ColWS of matrix
\DrawEq[vn]{a=(a1.n \ColNS)}{vn \SideNS-\Cols}
}

\AddEq{vector as matrix product left}
{
Therefore, we can represent the vector $\Vector v$
as product of matrices
\ShowEq{v=(e)(v)}
}

\AddEq{vector as matrix product right}
{
However, we cannot represent the vector $\Vector v$
as product of matrices
\ShowEq{v=(v)(e)}
because this product is not defined.
}

\DefLabeledExample{matrix in vector space}{\SideNS}
{
\ShowEq{define basis for example}
\ShowEq{vector as matrix product \SideNS}
}

\DefExample{matrix in right vector space}
{
We represent the basis $\Basis e$
of right vector space $V$
\ePrints{2020.06.01,2204.06320,2022.01.05}
\ifx\Semafor\ValueOff
over $D$\Hyph algebra $A$
(see the definition
\ShowEq{def right}%
\ShowEq{def AVector}%
\refDefinition{module over associative algebra}{\SideWS \VectorSetNS}
and the theorem
\refTheorem{basis over division algebra}{right})
\fi
as row of matrix
\DrawEq[en]{a=(a1.n row)}{right basis}
We represent coordinates of vector
\ShowEq{w=w*e right-cols}v{}
as vector column
\DrawEq[vn]{a=(a1.n col)}{right basis}
Therefore, we can represent the vector $\Vector v$
as product of matrices
\ShowEq{v=(e)(v)}
\ePrints{2020.06.01,2204.06320,2022.01.05}
\ifx\Semafor\ValueOff
We represent homomorphism
of right vector space $V$
using matrix
\ShowEq{v'=fv}
The equality
\EqRef{v'=fv}
expresses a traditional product of matrices $f$ and $v$.
\fi
}

\DefDefinition{row over column product}
{
Let the nubmer of columns of the matrix $a$ equal the number of rows of the matrix $b$.
\AddIndex{\RC product}{rc-product}
of matrices $a$ and $b$ has form
\ShowEq{rc-product}
\ShowEq{rc-product of matrices}
\ShowEq{entry rc-product of matrices}
\ShowEq{rc-product, matrices}
\RC product can be expressed as product of a row of the matrix $a$
over a column of the matrix $b$.
}

\DefDefinition{column over row product}
{
Let the nubmer of rows of the matrix $a$ equal the number of columns of the matrix $b$.
\AddIndex{\CR product}{cr-product}
of matrices $a$ and $b$ has form
\ShowEq{cr-product}
\ShowEq{cr-product of matrices}
\ShowEq{entry cr-product of matrices}
\ShowEq{cr-product, matrices}
\CR product can be expressed as product of a column of the matrix $a$
over a row of the matrix $b$.
}

\DefDefinition{transpose of matrix}
{
The transpose $a^T$ of the matrix $a$
exchanges rows and columns
\ShowEq{transpose of matrix, 1}
}

\DefDefinition{sum of matrices}
{
The sum of matrices $a$ and $b$ is defined by the equality
\ShowEq{(a+b)=}
}

\DefDefinition{rc power}
{
We introduce \AddIndex{\RC power}{rc power} of
\ePrints{2020.06.01}
\ifx\Semafor\ValueOff
$\mathcal{A}$\Hyph number
\else
the matrix
\fi
$a$
using recursive definition
\ShowEq{rc power}
\DrawEq{rc power, 0}1
\DrawEq{rc-power, n}1
}

\DefDefinition{cr power}
{
We introduce \AddIndex{\CR power}{cr power} of
$\mathcal{A}$\Hyph number $a$
using recursive definition
\ShowEq{cr power}
\DrawEq{cr power, 0}1
\DrawEq{cr-power, n}1
}

\DefDefinitionNote{Hadamard inverse of matrix}
{
Let
\DrawEq[ann]{a=(a11.nm matrix)}{}
be a matrix of $A$\Hyph numbers.
We call matrix\,\footnotemark
\ShowEq{Hadamard inverse of matrix}
\ShowEq{Hadamard inverse of matrix 2}
\ShowEq{Hadamard inverse of matrix 2 entry}
\AddIndex{Hadamard inverse of matrix}{Hadamard inverse of matrix} $a$
(\citeBib{q-alg-9705026}-\href{https://arxiv.org/pdf/math/9705026.pdf\#Page=4}{page 4}).
}{
The notation in the equality
\EqRef{Hadamard inverse of matrix 2 entry}
means that we exchange rows and columns
in Hadamard inverse.
\label{footnote: index of inverse element}
We can
formally write right site of the equality
\EqRef{Hadamard inverse of matrix 2 entry}
in the following form
\ShowEq{Hadamard inverse of matrix 3}
}

\DefDefinition{transformation coordinated with equivalence}
{
Let $S$ be equivalence on the set $A_2$.
Transformation $f$ is called
\AddIndex{coordinated with equivalence}{transformation coordinated with equivalence} $S$,
when
\ShowEq{fm1 fm2 modS}{}
follows from condition
\ShowEq{m1 m2 modS}.
}

\DefTheorem{transformation correlated with equivalence}
{
Consider equivalence $S$ on set $A_2$.
Consider $\Omega_1$\Hyph algebra
on the set
\ShowEq{End O2}{A_2}.
If any transformation
\ShowEq{f in End A2}{}
is coordinated with equivalence $S$,
then we can define the structure of $\Omega_1$\Hyph algebra
on the set
\ShowEq{End A2/S}S.
}

\DefTheorem{decompositions of morphism of representations}
{
Let
\ShowEq{f:A->*B}f{A_1}{A_2}
be representation of $\Omega_1$\Hyph algebra $A_1$,
\ShowEq{f:A->*B}g{B_1}{B_2}
be representation of $\Omega_1$\Hyph algebra $B_1$.
Let
\ShowEq{r12:A->B}tAB
be morphism of representations from $f$ into $g$.
Then there exist decompositions of $t_1$ and $t_2$,
which we describe using diagram
\ShowEq{decompositions of morphism of representations, diagram}
\StartLabelItem
\begin{enumerate}
\item
The kernel of homomorphism
\ShowEq{kernel of homomorphism i}
is a congruence on $\Omega_i$\Hyph algebra $A_i$,
\ShowEq{i=1,2}.
\labelItem{kernel of homomorphism i}
\item
There exists decomposition of homomorphism $t_i$,
\ShowEq{i=1,2},
\ShowEq{morphism of representations of algebra, homomorphism, 1}
\labelItem{exists decomposition of homomorphism}
\item
Maps
\ShowEq{morphism of representations of algebra, p1=}
\ShowEq{morphism of representations of algebra, p2=}
are natural homomorphisms.
\labelItem{p12 are natural homomorphisms}
\item
Maps
\ShowEq{morphism of representations of algebra, q1=}
\ShowEq{morphism of representations of algebra, q2=}
are isomorphisms.
\labelItem{q12 are isomorphisms}
\item
Maps
\ShowEq{morphism of representations of algebra, r1=}
\ShowEq{morphism of representations of algebra, r2=}
are monomorphisms.
\labelItem{r12 are monomorphisms}
\item $F$ is representation of $\Omega_1$\Hyph algebra $A_1/s$ in $A_2/s_2$
\item $G$ is representation of $\Omega_1$\Hyph algebra $t_1A_1$ in $t_2A_2$
\item The map
\ShowEq{map r12}p{}
is the morphism of representations $f$ and $F$
\labelItem{(p12) is morphism of representations}
\item The map
\ShowEq{map r12}q{}
is the isomorphism of representations $F$ and $G$
\labelItem{(q12) is isomorphism of representations}
\item The map
\ShowEq{map r12}r{}
is the morphism of representations $G$ and $g$
\labelItem{(r12) is morphism of representations}
\item There exists decompositions of morphism of representations
\labelItem{exists decompositions of morphism of representations}
\ShowEq{decompositions of morphism of representations}
\end{enumerate}
}

\DefDefinition{stable set of representation}
{
Let
\ShowEq{f:A->*B}f{A_1}{A_2}
be representation of $\Omega_1$\Hyph algebra $A_1$
in $\Omega_2$\Hyph algebra $A_2$.
The set
\ShowEq{B2 subset A2}
is called
\AddIndex{stable set of representation}
{stable set of representation} $f$,
if
\ShowEq{fam in B2}
for each
\ShowEq{a in A1},
$m\in B_2$.
}

\DefTheorem{subrepresentation of representation}
{
Let
\ShowEq{f:A->*B}f{A_1}{A_2}
be representation of $\Omega_1$\Hyph algebra $A_1$
in $\Omega_2$\Hyph algebra $A_2$.
Let set
\ShowEq{B2 subset A2}
be subalgebra of $\Omega_2$\Hyph algebra $A_2$
and stable set of representation $f$.
Then there exists representation
\ShowEq{representation of algebra A in algebra B}
such that
\ShowEq{fB2(a)=}
Representation $f_{B_2}$ is called
\AddIndex{subrepresentation}{subrepresentation}
of representation $f$.
}

\DefTheoremNote{lattice of subrepresentations}
{
The set\,\footnotemark
\ShowEq{lattice of subrepresentations}
of all subrepresentations of representation $f$ generates
a closure system on $\Omega_2$\Hyph algebra $A_2$
and therefore is a complete lattice.
}{
This definition is
similar to definition of the lattice of subalgebras
(\citeBib{Cohn: Universal Algebra}, p. 79, 80).
\ePrints{1908.04418,6860-2955}
\ifx\Semafor\ValueOn%
In general, in this and subsequent theorems of this chapter,
it is necessary to consider the structure of universal algebras
$A_1$ and $A_2$.
Because the main task of this chapter is
is the study of the structure of the representation,
I deliberately simplified the theorems so
that the details do not obscure the basic statements.
This topic will be discussed in more details in the chapter
\RefChapter{Basis of Diagram of Representations of Universal Algebra},
where theorems will be formulated in general form.
\fi
}

\AddEq{remark: closure operator, representation}
{
We denote the corresponding closure operator by
\ShowEq{closure operator, representation}
Thus
\ShowEq{subrepresentation generated by set}
is the intersection of all subalgebras
of $\Omega_2$\Hyph algebra $A_2$
containing $X$ and stable with respect to representation $f$.
}

\DefDefinition{generating set of representation}
{
\ShowEq{show closure operator, representation}
is called
\AddIndex{subrepresentation}{subrepresentation}
generated by set $X$,
and $X$ is a
generating set
of subrepresentation $J[f,X]$.
In particular, a
\AddIndex{generating set}{generating set}
of representation $f$
is a subset $X\subset A_2$ such that
\ShowEq{generating set of representation}
}

\DefTheoremNote{map of words of representation}
{
\ShowEq{Let A->*AB2 be representations}
Let $X$ be the generating set of representation $f$.
Let
\ShowEq{f:A->B}R{A_2}{B_2}
be reduced morphism of representation\,\footnotemark
and $X'=R(X)$.
Reduced morphism $R$ of representation
generates the map of $\Omega_2$\Hyph words
\ShowEq{map of words of representation}
such that
\StartLabelItem
\begin{enumerate}
\item If
\ShowEq{map of words of representation, m in X, 1}
then
\ShowEq{map of words of representation, m in X}
\item If
\ShowEq{map of words of representation, omega, 1}
then for operation
$\omega\in\Omega_2(n)$ holds
\ShowEq{map of words of representation, omega}
\item If
\ShowEq{map of words of representation, am, 1}
then
\ShowEq{map of words of representation, am}
\end{enumerate}
}{
I considered morphism of representation in the theorem
\RefTheorem[\RefRepresentation]{map of words of diagram of representations}.
}

\DefTheorem{Representation is single transitive iff}
{
Representation is single transitive iff for any $a, b \in A_2$
exists one and only one $g\in A_1$ such that $a=f(g)(b)$.
}

\DefProof{Representation is single transitive iff}
{
The theorem follows from definitions \RefDefinition{free representation of algebra}
and \RefDefinition{transitive representation of algebra}.
}

\DefTheorem{single transitive representation generates algebra}
{
Let
\ShowEq{f:A->*B}f{A_1}{A_2}
be a single transitive representation
of $\Omega_1$\Hyph algebra $A_1$
in $\Omega_2$\Hyph algebra $A_2$.
There is the structure of $\Omega_1$\Hyph algebra on the set $A_2$.
}

\AddEq{theorem: automorphism uniquely defined by image of basis}
{
\begin{ShadedTheorem}
\labelTheorem{automorphism uniquely defined by image of basis}
Let $\Basis e$ be the basis of the representation $f$.
Let
\ShowEq{R1 X->}{\Basis e}
be arbitrary map of the set $\Basis e$.
Consider the map of $\Omega_2$\Hyph words
\ShowEq{map of words of representation, W}{\Basis e}
that satisfies conditions
\RefItem{map of words of representation, m in X},
\RefItem{map of words of representation, omega},
\RefItem{map of words of representation, am}
and such that
\ShowEq{map of words of representation, W in X}e{\Basis e}
There exists unique endomorphism of representation $f$\,\footnotemark
\ShowEq{R:A2->A2}
defined by rule
\ShowEq{endomorphism based on map of words of representation}{\Basis e}
\end{ShadedTheorem}
\footnotetext{\,
This statement is
similar to the theorem \citeBib{Serge Lang}-4.1, p. 135.
}
}

\DefDefinition{parity of permutation}
{
The map
\ShowEq{sigma->+-}
defined by the equality
\def\permutation{permutation }
\def\even{even}
\def\odd{odd}
\ShowEq{parity of permutation}
is called
\AddIndex{parity of permutation}{parity of permutation}.
}

\DefText{passive representation of vector space 1}
{
\ShowText{coordinates of vector with respect to basis (1,2) \Cols}v3
Then
\ShowEq{v=vi ei 123, \SideNS-\Cols}
Let passive transformation $g_1$ map the basis
\ShowEq{Basis e in BVG}1
into the basis
\ShowEq{Basis e in BVG}2
\DrawEq[12g1]{e2i=aij e1j \Product-\Cols}{12}
Let passive transformation $g_2$ map the basis
\ShowEq{Basis e in BVG}2
into the basis
\ShowEq{Basis e in BVG}3
\DrawEq[23g2]{e2i=aij e1j \Product-\Cols}{23}
}

\DefText{passive representation of vector space 2}
{
The transformation
\ShowEq{e3=g12 cr e1 \SideNS-\Cols}
is the product of transformations
\eqRef{e2i=aij e1j \Product-\Cols}{12},
\eqRef{e2i=aij e1j \Product-\Cols}{23}.
The coordinate transformation
\DrawEq[v2v1g1{\SideNS}{\Cols}]{v2=v1g-}{1 \SideNS-\Cols}
corresponds to passive transformation $g_1$
and follows from equalities
\EqRef{v=vi ei 123, \SideNS-\Cols},
\eqRef{e2i=aij e1j \Product-\Cols}{12}
and the theorem
\refTheorem{passive transformation of vector space}{\SideNS-\Cols}.
}

\DefText{passive representation of vector space 3}
{
Similarly,
\ePrints{1908.04418,6860-2955}%
\ifx\Semafor\ValueOff%
according to the theorem
\refTheorem{passive transformation of vector space}{\SideNS-\Cols},
\fi
the coordinate transformation
\DrawEq[v3v2g2{\SideNS}{\Cols}]{v2=v1g-}{23 \SideNS-\Cols}
corresponds to passive transformation $g_2$
and follows from equalities
\EqRef{v=vi ei 123, \SideNS-\Cols},
\eqRef{e2i=aij e1j \Product-\Cols}{23}.
The product of coordinate transformations
\eqRef{v2=v1g-}{1 left-cols},
\eqRef{v2=v1g-}{23 left-cols}
has form
\ShowEq{v3=v1 cr g1- cr g2- \SideNS-\Cols}
and is coordinate transformation
corresponding to passive transformation
\EqRef{e3=g12 cr e1 \SideNS-\Cols}.
The equality
\ShowEq{v3=v1 cr g21- \SideNS-\Cols}
follows from equalities
\EqRef{cr-inverse matrix of product},
\EqRef{v3=v1 cr g1- cr g2- \SideNS-\Cols}.
}

\DefText{passive representation of left vector space}
{
\ShowText{passive representation of vector space 1}
\ShowText{passive representation of vector space 2}
\ShowText{passive representation of vector space 3}
}

\DefDefinition{Left side contravariant representation}
{
Left\Hyph side representation
\ShowEq{f:A->*B}f{A_1}{A_2}
is called
\AddIndex{contravariant}{contravariant representation}
if the equality
\ShowEq{Left side contravariant representation}
is true.
}

\DefDefinition{Left side covariant representation}
{
Left\Hyph side representation
\ShowEq{f:A->*B}f{A_1}{A_2}
is called
\AddIndex{covariant}{covariant representation}
if the equality
\ShowEq{Left side covariant representation}
is true.
}

\DefDefinition{cr-inverse matrix}
{
\ePrints{2020.06.01}
\ifx\Semafor\ValueOff
$\mathcal{A}$\Hyph number
\else
The matrix
\fi
\ShowEq{cr-inverse element}
is \AddIndex{\CR inverse element}{cr-inverse element}
\ePrints{2020.06.01}
\ifx\Semafor\ValueOff
of $\mathcal{A}$\Hyph number
\else
of the matrix
\fi
$a$ if
\DrawEq{cr-inverse matrix}1
\ePrints{2020.06.01}
\ifx\Semafor\ValueOn
Matrix $a$ is called
\CR regular,
if there exists \CR inverse matrix.
\fi
}

\DefDefinition{delta=Matrix}
{
Matrix
\ShowEq{delta=Matrix}
is identity for both products.
}

\DefRemark{pronunciation of product}
{
We will use symbol \RC\,\, or \CR\,\, in name of properties of each product
and in the notation.
We can read the symbol
$\RCstar$ as
$rc$\hyph product (product of row over column)
and the symbol $\CRstar$ as
$cr$\hyph product (product of column over row).
In order to keep this
notation consistent with the existing one we assume that
we have in mind \RC product when
no clear notation is present.
\ePrints{2020.06.01}
\ifx\Semafor\ValueOff
This rule we extend to following terminology.
To draw symbol of product,
we put two symbols of product in the place
of index which participate in sum.
For instance, if product of $A$\Hyph numbers has form
$a\circ b$,
then \RC product
of matrices $a$ and $b$ has form
$a\RCcirc b$
and \CR product
of matrices $a$ and $b$ has form
$a\CRcirc b$.
\fi
}

\DefTheorem{quasideterminant rc cr, matrix 2x2}
{
\ShowEq{quasideterminant matrix 2x2 def}
Consider matrix
\ShowEq{inverse matrix 2x2, 0}
Then
\ShowEq{quasideterminant rc, matrix 2x2}
\ShowEq{quasideterminant cr, matrix 2x2}
\ShowEq{rc inverse matrix 2x2}
}

\DefDefinition{RC-quasideterminant}
{
\AddIndex{$(\aUD{}ji)$\hyph \RC quasideterminant}{j i RC-quasideterminant}
of \nTimes matrix $a$
is formal expression
\ShowEq{j i RC-quasideterminant definition}
\ShowEq{j i RC-quasideterminant =}
We consider $\aUD (ji)$\hyph \RC quasideterminant
as an entry of the matrix
\ShowEq{RC-quasideterminant definition}
\ShowEq{RC-quasideterminant definition=}
which is called
\AddIndex{\RC quasideterminant}{RC-quasideterminant}.
}

\DefDefinition{biring}
{
The set $\mathcal A$ is a \AddIndex{biring}{biring}
if we defined on $\mathcal A$ an unary operation, say transpose,
and three binary operations,
say \RC product, \CR product and sum, such that
\begin{itemize}
\item \RC product and sum define structure of ring on $\mathcal A$
\item \CR product and sum define structure of ring on $\mathcal A$
\item both products have common identity $\delta$
\item products satisfy equation
\DrawEq{rc transpose}{}
\item transpose of identity is identity
\DrawEq{transpose of identity}E
\item double transpose is original element
\ShowEq{double transpose}
\end{itemize}
}

\DefTheorem{double transpose is original matrix}
{
Double transpose is original matrix
\ShowEq{double transpose}
}

\DefProof{double transpose is original matrix}
{
The theorem follows from the definition
\RefDefinition{transpose of matrix}.
}

\DefEq
{
\begin{ShadedTheorem}[duality principle for biring]
\AddIndex{}{duality principle for biring}
\labelTheorem{duality principle for biring}
Let $\mathcal{A}$ be true statement about
biring $A$.
If we exchange the same time
\begin{itemize}
\item $a\in A$ and $a^T$
\item \RC product and \CR product
\end{itemize}
then we soon get true statement.
\end{ShadedTheorem}
}
{theorem: duality principle for biring}

\DefEq
{
\begin{ShadedTheorem}[duality principle for biring of matrices]
\AddIndex{}{duality principle for biring of matrices}
\labelTheorem{duality principle for biring of matrices}
Let $A$ be biring of matrices.
Let $\mathcal{A}$ be true statement about matrices.
If we exchange the same time
\begin{itemize}
\item rows and columns of all matrices
\item \RC product and \CR product
\end{itemize}
then we soon get true statement.
\end{ShadedTheorem}
\begin{proof}
This is the immediate consequence
of the theorem \RefTheorem{duality principle for biring}.
\end{proof}
}
{theorem: duality principle for biring of matrices}

\DefEq
{
\begin{remark}
\labelRemark{reducible biring}
If product in $\Omega$\Hyph $A$ ring is commutative, then
\ShowEq{reducibility of products}
\AddIndex{Reducible biring}{reducible biring} is
the biring which holds
\AddIndex{condition of reducibility of products}
{condition of reducibility of products}
\EqRef{reducibility of products}.
So, in reducible biring, it is enough to consider only
\RC product. However in case when
the order of factors is essential
we will use \CR product also.
\qed
\end{remark}
}
{remark: reducible biring}

\DefEq
{
Since left\Hyph side representation and right\Hyph side representation
are based on homomorphism of $\Omega$\Hyph group,
then the following statement is true.

\begin{ShadedTheorem}[duality principle for representation of multiplicative $\Omega$\Hyph group]
\labelTheorem{duality principle, algebra representation}
Any statement which holds for
left\Hyph side representation of multiplicative $\Omega$\Hyph group $A_1$
holds also for right\Hyph side representation of multiplicative $\Omega$\Hyph group $A_1$, if
we will use right\Hyph side product over $A_1$\Hyph number $a_1$ instead of
left\Hyph side product over $A_1$\Hyph number $a_1$.
\qed
\end{ShadedTheorem}
}
{theorem: duality principle, algebra representation}

\DefTheorem{left right shift}
{
The product
\ShowEq{(ab)->ab}
in multiplicative $\Omega$\Hyph group $A$
determines two different representations.
\begin{itemize}
\EqParm{associative shift}{def left}
\EqParm{associative shift}{def right}
\end{itemize}
}

\AddEq{non associative shift}
{
\item
The \SideWS shift
\ShowEq{\SideWS shift}
is representation
of $\Omega$\Hyph algebra $A$
in $\Omega$\Hyph algebra $A$.
}

\AddEq{associative shift}
{
\item
the \AddIndex{\SideWS shift}{\SideWS shift}
\ShowEq{\SideWS shift =}
\ShowEq{\SideWS shift}
is \SideNS\Hyph side representation
of multiplicative $\Omega$\Hyph group $A$
in $\Omega$\Hyph algebra $A$
\ShowEq{\SideWS shift, product}
}

\DefDefinition{representation of algebra}
{
Let the set $A_2$ be $\Omega_2$\Hyph algebra.
Let
the set of transformations
\ShowEq{End O2}{A_2}{}
be  $\Omega_1$\Hyph algebra.
The homomorphism
\ShowEq{representation of algebra}
of $\Omega_1$\Hyph algebra $A_1$ into
$\Omega_1$\Hyph algebra
\ShowEq{End O2}{A_2}{}
is called
\AddIndex{representation of $\Omega_1$\Hyph algebra}
{representation of algebra}
$A_1$ or
\AddIndex{$A_1$\Hyph representation}{A representation of algebra}
in $\Omega_2$\Hyph algebra $A_2$.
}

\DefEq
{
If multiplicative $\Omega$\Hyph group $A_1$ is Abelian,
then there is no difference between left\Hyph side and right\Hyph side representations.

\begin{ShadedDefinition}
\labelDefinition{representation of Abelian group}
Let $A_1$ be Abelian multiplicative $\Omega$\Hyph group.
We call a homomorphism of multiplicative $\Omega$\Hyph group
\DrawEq{representation of group, map}{Abelian group}
\AddIndex{representation}{representation of algebra}
of multiplicative $\Omega$\Hyph group $A_1$
or
\AddIndex{$A_1$\Hyph representation}{A representation of algebra}
in $\Omega_2$\Hyph algebra $A_2$ if the map $f$ holds
\DrawEq{left-side representation of group}{Abelian}
\end{ShadedDefinition}

Usually we identify a representation of the Abelian multiplicative $\Omega$\Hyph group $A_1$
and a left\Hyph side representation of the multiplicative $\Omega$\Hyph group $A_1$.
However, if it is necessary for us,
we identify a representation of the Abelian multiplicative $\Omega$\Hyph group $A_1$
and a right\Hyph side representation of the multiplicative $\Omega$\Hyph group $A_1$.
}
{definition: representation of Abelian group}

\DefTheorem{two representations of group}
{
Let $A_1$ be multiplicative $\Omega$\Hyph group.
Let left\Hyph side $A_1$\Hyph representation $f$
on $\Omega_2$\Hyph algebra $A_2$ be single transitive.
Then we can uniquely define a single transitive right\Hyph side $A_1$\Hyph representation $h$
on $\Omega_2$\Hyph algebra $A_2$
such that diagram
\ShowEq{two representations of group}
is commutative for any $a_1$, $b_1\in A_1$.
}

\AddEq{definition: twin representations}
{
Representations $f$ and $h$ are called
\AddIndex{twin representations}{twin representations}
of the multiplicative $\Omega$\Hyph group $A_1$.
}

\DefRemark{one representation of group}
{
It is clear that transformations $L(a)$ and $R(a)$
are different until the multiplicative $\Omega$\Hyph group $A_1$ is nonabelian.
However they both are maps onto.
Theorem \RefTheorem{two representations of group} states that if both
right and left shift presentations exist on the set $A_2$,
then we can define two commuting representations on the set $A_2$.
The right shift or the left shift only cannot represent both types of representation.
To understand why it is so let us change diagram
\EqRef{Diagram: two representations of group}
and assume
\ShowEq{haa=Laa}
instead of
\ShowEq{haa=aRa}
and let us see what expression $h(a_1)$ has at
the point $c_2$. The diagram
\ShowEq{Diagram: two representations of group 1}
is equivalent to the diagram
\ShowEq{Diagram: two representations of group 2}
and we have
\ShowEq{d2=...c2 1}
Therefore
\ShowEq{d2=...c2}
We see that the representation of $h$
depends on its argument.
}

\DefTheorem{single transitive representation of group}
{
If there exists single transitive representation
\ShowEq{f:A->*B}f{A_1}{A_2}
of multiplicative $\Omega$\Hyph group $A_1$
in $\Omega_2$\Hyph algebra $A_2$,
then we can uniquely define coordinates on $A_2$ using $A_1$\Hyph numbers.

If $f$ is left\Hyph side single transitive representation
then $f(a)$ is equivalent to the left shift $L(a)$ on the group $A_1$.
If $f$ is right\Hyph side single transitive representation
then $f(a)$ is equivalent to the right shift $R(a)$ on the group $A_1$.
}

\DefDefinition{algebra of sets}
{
A nonempty system of sets $\mathcal R$ is called
\AddIndex{ring of sets}{ring of sets},\,\footnote{
See also the definition
\citeBib{Kolmogorov Fomin}\Hyph 1,  page 31.}
if condition
\ShowEq{A, B in R}
imply
\ShowEq{A, B in R 1}
A set
\ShowEq{E in R}
is called
\AddIndex{unit of ring of sets}{unit of ring of sets}
if
\ShowEq{A cap E=A}
A ring of sets with unit is called
\AddIndex{algebra of sets}{algebra of sets}.
}

\DefEq
{
\begin{remark}
\labelRemark{A cup/minus B}
For any
\ShowEq{A cup/minus B}
Therefore, if
\ShowEq{A, B in R},
then
\ShowEq{A cup/minus B R}
\qed
\end{remark}
}
{remark: A cup/minus B}

\DefDefinition{sigma algebra of sets}
{
The ring of sets $\mathcal R$ is called
\AddIndex{\(\sigma\)\Hyph ring of sets}{sigma ring of sets},\,\footnote{
See similar definition in
\citeBib{Kolmogorov Fomin}, page 35, definition 3.}
if condition
\ShowEq{Ai in R}
imply
\ShowEq{Ai in R 1}
\(\sigma\)\Hyph Ring of sets with unit is called
\AddIndex{\(\sigma\)\Hyph algebra of sets}{sigma algebra of sets}.
}

\DefDefinition{Borel algebra}
{
Minimal \(\sigma\)\Hyph algebra
\ShowEq{Borel algebra}
generated by the set of all open balls of normed $\Omega$\Hyph group $A$
is called
\AddIndex{Borel algebra}{Borel algebra}.\,\footnote{
See remark in
\citeBib{Kolmogorov Fomin}, p. 36.
According to the remark
\ref{remark: A cup/minus B},
the set of closed balls also generates
Borel algebra.
}
A set
belonging to Borel algebra
is called
\AddIndex{Borel set}{Borel set}
or
\AddIndex{$B$\Hyph set}{B set}.
}

\DefDefinition{C_X C_Y measurable}
{
Let $\mathcal C_X$ be
\(\sigma\)\Hyph algebra of sets of set $X$.
Let $\mathcal C_Y$ be
\(\sigma\)\Hyph algebra of sets of set $Y$.
The map
\ShowEq{f:A->B}fXY
is called
$(\mathcal C_X,\mathcal C_Y)$\Hyph measurable\,\footnote{
See similar definition in
\citeBib{Kolmogorov Fomin}, page 284.}
if for any set
\ShowEq{CX CY measurable map}
}

\DefEq
{
\begin{example}
\labelExample{measurable map}
Let $\mu$ be a \(\sigma\)\Hyph additive measure defined on the set $X$.
Let
\ShowEq{algebra of sets measurable with respect to measure}
be \(\sigma\)\Hyph algebra of sets measurable with respect to measure $\mu$.
Let
$\mathcal B(A)$
be Borel algebra of
normed $\Omega$\Hyph group $A$.
The map
\ShowEq{f:A->B}fXA
is called
\AddIndex{$\mu$\Hyph measurable}{mu measurable map}\,\footnote{
See similar definition in
\citeBib{Kolmogorov Fomin},  pp. 284, 285, definition 1.
If the measure $\mu$ is defined on the set $X$ by context,
then we also call the map
\ShowEq{f:A->B}fXA
\AddIndex{measurable}{measurable map}.
}
if for any set
\ShowEq{measurable map}
\qed
\end{example}
}
{example: measurable map}

\DefDefinition{simple map}
{
Let $\mu$ be a \(\sigma\)\Hyph additive measure defined on the set $X$.
The map
\ShowEq{f:A->B}fXA
into normed $\Omega$\Hyph group $A$ is called
\AddIndex{simple map}{simple map},
if this map is $\mu$\Hyph measurable and
its range is finite or countable set.
}

\DefTheorem{measurable simple map}
{
Let range
\ShowEq{range f}
of the map
\ShowEq{f:A->B}fXA
be finite or countable set.
The map $f$ is $\mu$\Hyph measurable iff
all sets
\DrawEq{Fn=}{}
are $\mu$\Hyph measurable.\,\footnote{
See similar theorem in
\citeBib{Kolmogorov Fomin},  page 286, theorem 4.}
}

\DefDefinition{measure}
{
Let \(\mathcal C_m\) be semiring of sets.\,\footnote{
See also the definition
\citeBib{Kolmogorov Fomin}\Hyph 1 on page 270.}
The map
\ShowEq{f:A->B}f{\mathcal C_m}R
is called
\AddIndex{measure}{measure},
if
\StartLabelItem
\begin{enumerate}
\item
\ShowEq{m(A)>0}
\item
The map \(m\) is additive map.
If a set
\ShowEq{A in Cm}
has finite expansion
\DrawEq{A=+Ai Cm}{}
where
\ShowEq{Ai*Aj=0},
then
\ShowEq{m A=+m Ai}
\end{enumerate}
}

\DefDefinition{complete measure}
{
A measure \(\mu\) is called
\AddIndex{complete measure}{complete measure},
if conditions
\ShowEq{B in A, mu(A)=0}
imply
that \(B\) is measurable set.
}

\DefDefinition{sigma-additive measure}
{
Let $\mathcal C_{\mu}$ be
\(\sigma\)\Hyph algebra of sets of set $F$.\,\footnote{
See similar definitions in
\citeBib{Kolmogorov Fomin}, definition 1 on page 270
and definition 2 on page 272.}
The map
\ShowEq{f:A->B}{\mu}{\mathcal C_{\mu}}R
into real field $R$ is called
\AddIndex{\(\sigma\)\Hyph additive measure}{sigma-additive measure},
if, for any set
\ShowEq{X in Cmu},
following conditions are true.
\StartLabelItem
\begin{enumerate}
\item
\ShowEq{muX>0}
\item
Let
\DrawEq{X=+Xi}{}
be finite or countable union of sets
\ShowEq{Xn in Cmu}
Then
\DrawEqParm{mu(X=+Xi)}X{}
where series on the right converges absolutely.
\labelItem{mu(X=+Xi)}
\end{enumerate}
}

\DefEq
{
Let $\mu$ be a \(\sigma\)\Hyph additive measure defined on the set $X$.
Let effective representation of real field \(R\)
in complete Abelian \(\Omega\)\Hyph group \(A\) be defined.\,\footnote{In
other words,
\(\Omega\)\Hyph group \(A\) is
\(R\)\Hyph vector space.}
}
{remark - measure - effective representation of real field}

\DefDefinition{series converges normally}
{
Let \(a_i\) be a sequence of \(A\)\Hyph numbers.
If
\ShowEq{sum |ai|<infty}
then we say that the {\bf series}
\ShowEq{sum ai}
\AddIndex{converges normally}
{series converges normally}.\,\footnote{See also
the definition of normal convergence of the series
on page \citeBib{Cartan differential form}-12.}
}

\DefDefinition{Integral of Map over Set of Finite Measure}
{
\(\mu\)\Hyph measurable map
\ShowEq{f:A->B}fXA
is called
\AddIndex{integrable map}{integrable map}
over the set \(X\),\,\footnote{
See also the definition in
\citeBib{Kolmogorov Fomin}, page 296.
}
if there exists a sequence
of simple integrable over the set \(X\) maps
\ShowEq{f:A->B}{f_n}XA
converging uniformly to \(f\).
Since the map \(f\) is integrable map,
then the limit
\ShowEq{integral of map}
\ShowEq{integral of measurable map}
is called
\AddIndex{Lebesgue integral}{Lebesgue integral}
of map \(f\) over the set \(X\).
}

\DefDefinition{integrable map}
{
For simple map
\ShowEq{f:A->B}fXA
consider series
\DrawEq{sum mu(F) f}{integral}
where
\begin{itemize}
\item
The set
\ShowEq{f1...}
is domain of the map \(f\)
\item
Since \(n\ne m\), then
\ShowEq{fn ne fm}
\item
\DrawEq{Fn=}{-}
\end{itemize}
Simple map
\ShowEq{f:A->B}fXA
is called
\AddIndex{integrable map}{integrable map}
over the set \(X\)
if series
\eqRef{sum mu(F) f}{integral}
converges normally.\,\footnote{
See similar definition in
\citeBib{Kolmogorov Fomin}, p. 294.}
Since the map \(f\) is integrable map,
then sum of series
\eqRef{sum mu(F) f}{integral}
is called
\AddIndex{Lebesgue integral}{Lebesgue integral}
of map \(f\) over the set \(X\)
\ShowEq{integral of map}
\ShowEq{integral of simple map}
}

\DefTheorem{|int f|<int|f|}
{
Let
\ShowEq{f:A->B}fXA
be measurable map.
Integral
\ShowEq{int f X}
exists
iff
integral
\ShowEq{int |f| X}
exists.
Then
\DrawEq{|int f|<int|f|}{measurable}
}

\DefTheorem{int |f|<M mu}
{
Let
\ShowEq{f:A->B}fXA
be $\mu$\Hyph measurable map such that
\DrawEq{|f(x)|<M}{}
Since the measure of the set \(X\) is finite, then
\DrawEq{int |f|<M mu}{measurable}
}

\DefTheorem{int f+g X}
{
Let
\ShowEq{f:A->B}fXA
\ShowEq{f:A->B}gXA
be $\mu$\Hyph measurable maps
with compact range.
Since there exist integrals
\ShowEq{int f X}
\ShowEq{int g X}
then there exists integral
\ShowEq{int f+g X}
and
\DrawEq{int f+g X =}{measurable}
}

%% file: Stmt.Representation.Eq.tex

\def\Act{\bullet}
\newcommand\wXm[1][m]{w[f,X,#1]}
\newcommand\wXR{w[f\rightarrow g,X,R]}
\def\UnitNMatrix{\aD En}

\AddEq[1]{i=1,2}
{
$i=1$, $2$#1
}

\AddEquation{morphism of representations of algebra, homomorphism, 1}
{
t_i=r_i\circ q_i\circ p_i
}

\AddEq{morphism of representations of algebra, p1=}
{
\[
\RedText{p_1(a)}=a^{\mathrm{ker}\,t_1}
\]
}

\AddEq{morphism of representations of algebra, p2=}
{
\[
\BlueText{p_2(a)}=a^{\mathrm{ker}\,t_2}
\]
}

\AddEquation{morphism of representations of algebra, q1=}
{
q_1(\RedText{p_1(a)})=\RedText{t_1(a)}
}

\AddEquation{morphism of representations of algebra, q2=}
{
q_2(\BlueText{p_2(a)})=\BlueText{t_2(a)}
}

\AddEq{morphism of representations of algebra, r1=}
{
\[
r_1:\RedText{t_1(a)}\in f(A_1)\rightarrow \RedText{t_1(a)}\in B_1
\]
}

\AddEq{morphism of representations of algebra, r2=}
{
\[
r_2:\BlueText{t_2(a)}\in f(A_2)\rightarrow \BlueText{t_2(a)}\in B_2
\]
}

\AddEq [2]{map r12}
{
$#1=(#1_1,#1_2)$#2
}

\AddEquation{decompositions of morphism of representations}
{
(t_1,t_2)
=
(r_1,r_2)
\circ
(q_1,q_2)
\circ
(p_1,p_2)
}

\AddEq{r1(a)=r1(a)}
{
\[
\RedText{r_1(a)}=r_1(a)
\]
}

\AddEq{g(r1(a))}
{
$g(\RedText{r_1(a)})$
}

\AddEq{r2(m)in B2}
{
$\BlueText{r_2(m)}\in B_2$
}

\AddEq{r2(f(a,m))}
{
$\BlueText{r_2(f(a)(m))}$.
}

\AddEquation{morphism of representations of universal algebra, definition, 2m 2}
{
\xymatrix{
A_2\ar[rr]^{r_2}&&B_2\\
A_1\ar[rr]^{r_1}\ar[u]^f|{*}&&B_1\ar[u]^g|{*}
}
}

\AddEq{n=r2(m)}
{
$n=r_2(m)\in B_2$
}

\AddEq{morphism of representations of universal algebra, 2m}
{
\BlueText{r_2(f(a)(m))}=g(\RedText{r_1(a)})(\BlueText{r_2(m)})
}

\AddEq{morphism of representations of universal algebra, 2m 1}
{
\xymatrix{
&A_2\ar[dd]_(.3){f(a)}\ar[rr]^{r_2}&&B_2\ar[dd]^(.3){g(\RedText{r_1(a)})}\\
&\ar @{}[rr]|{(1)}&&\\
&A_2\ar[rr]^{r_2}&&B_2\\
A_1\ar[rr]^{r_1}\ar@{=>}[uur]^(.3)f&&B_1\ar@{=>}[uur]^(.3)g
}
}

\AddEq[1]{a in A1}
{
$a\in A_1$#1
}

\AddEq{fam in B2}
{
$f(a)(m)\in B_2$
}

\AddEq{B2 subset A2}
{
$B_2\subset A_2$
}

\AddEq{fB2(a)=}
{
$f_{B_2}(a)=f(a)|_{B_2}$.
}

\AddEq{lattice of subrepresentations}
{
\symb{\mathcal B_f}{lattice of subrepresentations}1
}

\AddEq{representation of algebra A in algebra B}
{
\[
\xymatrix
{
f_{B_2}:A_1\ar[r]|{*}&B_2
}
\]
}

\AddEq{rc transpose}
{
(a\RCstar b)^T=a^T\CRstar b^T
}

\AddEq{cr transpose}
{
(a\CRstar b)^T=(a^T)\RCstar (b^T)
}

\DefTheorem{rc transpose}
{
\DrawEq{rc transpose}{Theorem}
}

\DefTheorem{cr transpose}
{
\DrawEq{cr transpose}{Theorem}
}

\AddEquation{rc transpose, 1}
{
\aUD{((a\RCstar b)^T)}ji
=\aUD{(a\RCstar b)}ij
=\aUD aik\aUD bkj
=\aUD{(a^T)}ki\aUD{(b^T)}jk
=\aUD{((a^T)\CRstar (b^T))}ji
}

\AddEquation{L(b)a=}
{
L(b)\circ a=[b,a]
}

\AddEquation{L(c)L(b)a=}
{
L(c)\circ L(b)\circ a=L(c)\circ(L(b)\circ a)=[c,[b,a]]
}

\AddEq{Lie product on set of left shifts}
{
\[
[L(c),L(b)]\circ a=L(c)\circ L(b)\circ a-L(b)\circ L(c)\circ a
\]
}

\AddEquation{[L(c)L(b)]a= 2}
{
[L(c),L(b)]\circ a=L([c,b])\circ a
}

\AddEquation{[L(c)L(b)]a= 1}
{
L(c)\circ L(b)\circ a-L(b)\circ L(c)\circ a=L([c,b])\circ a
}

\AddEquation{Lee identity 1}
{
[c,[b,a]]+[b,[a,c]]=-[a,[c,b]]=[[c,b],a]
}

\AddEquation{Lee identity}
{
[c,[b,a]]+[b,[a,c]]+[a,[c,b]]=0
}

\AddEquation{[L(c)L(b)]a=}
{
\begin{split}
L(c)\circ L(b)\circ a-L(b)\circ L(c)\circ a&=[c,[b,a]]-[b,[c,a]]
\\
&=[c,[b,a]]+[b,[a,c]]
\end{split}
}

\AddEq{direct sum of Abelian groups}
{
\symb{A\oplus B}{direct sum}1
}

\AddEquation{fik fjk = fjk fik}
{
f_{ik}(a_i)(f_{jk}(a_j)(a_k))
=
f_{jk}(a_j)(f_{ik}(a_i)(a_k))
}

\AddEquation{reduced polymorphism of representation, akl}
{
\begin{array}{r@{\,}l}
&r_2(m_1,...,\BlueText{f_k(a)(m_k)},...,m_l,...,m_n)
\\
=&r_2(m_1,...,m_k,...,\BlueText{f_l(a)(m_l)},...,m_n)
\end{array}
}

\ePrints{1502.04063,5114-6019}
\ifx\Semafor\ValueOff
\AddEquation{reduced polymorphism of representation, omega2}
{
\begin{array}{r@{\,}l}
&\,r_2(m_1,...,m_{k\cdot 1}...m_{k\cdot p}\omega_2,...,m_n)
\\
=&\,
r_2(m_1,...,m_{k\cdot 1},...,m_n)
...
r_2(m_1,...,m_{k\cdot p},...,m_n)
\omega_2
\end{array}
}
\fi

\AddEq{kl,1n}
{
$k$, $l$, $k=1$, ..., $n$, $l=1$, ..., $n$,
}

\AddEq{reduced morphism representation =}
{
r_2(f(a)(m))=g(a)(r_2(m))
}

\AddEquation{morphism of representations of universal algebra, definition, 2m-1}
{
\BlueText{r_2(f(a')(r_2^{-1}(m')))}
=f(\RedText{a'})(\BlueText{m'})
}

\AddEq{A in xAi}
{
\[A\subseteq\prod_{\iI}A_i\]
}

\AddEq[1]{A=o+Ai}
{
\[
#1=\bigoplus_{\iI}#1_i
\]
}

\AddEq[1]{m1 m2 modS}
{
$m_1\equiv m_2(\mathrm{mod} S)$#1
}

\AddEq[1]{fm1 fm2 modS}
{
$f(m_1)\equiv f(m_2)(\mathrm{mod}\ S)$#1
}

\AddEq[2]{End A2/S}
{
$\End(\Omega_2;A_2/#1)$#2
}

\AddEq[1]{f in End A2}
{
$f\in\End(\Omega_2;A_2)$#1
}

\AddEq[2]{End O2}
{
$\End(\Omega_2,#1)$#2
}

\AddEq [3]{r12:A->B}
{
\[
(#1_1:#2_1\rightarrow #3_1,\ #1_2:#2_2\rightarrow #3_2)
\]
}

\AddEq{decompositions of morphism of representations, diagram}
{
\[
\xymatrix{
&&A_2/s_2\ar[rrrrr]^{q_2}\ar@{}[drrrrr]|{(5)}&&\ar@{}[dddddll]|{(4)}&
\ar@{}[dddddrr]|{(6)}&&t_2A_2\ar[ddddd]^{r_2}\\
&&&&&&&\\
A_1/s_1\ar[r]^{q_1}\ar@/^2pc/@{=>}[urrr]^F&
t_1A_1\ar[d]^{r_1}\ar@{=>}[urrrrr]^(.4)G&&&
A_2/s_2\ar[r]^{q_2}\ar[lluu]_{F(\RedText{p_1(a)})}&
t_2A_2\ar[d]^{r_2}\ar[rruu]_{G(\RedText{t_1(a)})}\\
A_1\ar[r]_{t_1}\ar[u]^{p_1}\ar@{}[ur]|{(1)}\ar@/_2pc/@{=>}[drrr]^f&
B_1\ar@{=>}[drrrrr]^(.4)g&&&
A_2\ar[r]_{t_2}\ar[u]^{p_2}\ar@{}[ur]|{(2)}\ar[ddll]^{f(a)}&
B_2\ar[ddrr]^{g(\RedText{t_1(a)})}\\
&&&&&&&\\
&&A_2\ar[uuuuu]^{p_2}\ar[rrrrr]_{t_2}\ar@{}[urrrrr]|{(3)}&&&&&B_2
}
\]
}

\AddEq{kernel of homomorphism i}
{
$\mathrm{ker}\,t_i
=t_i\circ t_i^{-1}$
}

\AddEq{fxi,i=sum fxi}
{
f\left(\bigoplus_{\iI}x_i\right)=\sum_{\iI}f_i(x_i)
}

\AddEq{fi(x+y)=}
{
\[
f_i(x_i+y_i)=f_i(x_i)+f_i(y_i)
\]
}

\AddEq{f(x+y)i=}
{
\begin{align*}
f(\bigoplus_{\iI}(x_i+y_i))&=
\sum_{\iI}f_i(x_i+y_i)=
\sum_{\iI}(f_i(x_i)+f_i(y_i))
\\&=\sum_{\iI}f_i(x_i)+\sum_{\iI}f_i(y_i)
\\&=f(\bigoplus_{\iI}x_i)+f(\bigoplus_{\iI}y_i)
\end{align*}
}

\AddEq{structure of subrepresentations}
{
\StartLabelItem
\begin{enumerate}
\item
$X_0=X$
\labelItem{X0=X}
\item
$x\in X_k=>x\in X_{k+1}$
\labelItem{Xk in Xk+1}
\item
$x_1\in X_k$, ..., $x_n\in X_k$, $\omega\in\Omega_2(n)
=>x_1...x_n\omega\in X_{k+1}$
\labelItem{x1n omega in Xk+1}
\item
$x\in X_k$, $a\in A=>f(a)(x)\in X_{k+1}$
\labelItem{ax in Xk+1}
\end{enumerate}
}

\AddEquation{structure of subrepresentations, 1}
{
\bigcup_{k=0}^{\infty}X_k=J[f,X]
}

\AddEquation{lj(x)=0x0}
{
\lambda_j(x)=(\delta^i_jx,\iI)
}

\AddEq{flj=fj}
{
\[
f\circ\lambda_j(x)=\sum_{\iI}f_i(\delta_j^ix)=f_j(x)
\]
}

\AddEq[1]{(xi)in A}
{
$(x_i,\iI)\in A$#1
}

\DefEq
{
$\|x\|$, $x\in C$,
}
{|x in C|}

\DefEquation
{
\|a-b\|\ge|\,\|a\|-\|b\|\,|
}
{|a-b|>|a|-|b|}

\DefEq
{
\[
\|f(x')-f(x)\|_2<\epsilon
\]
}
{|f(x)-f(x')|<epsilon}

\AddEq[1]{End A(O2) B}
{
$\End(A(\Omega_2);B)$#1
}

\AddEq{open ball}
{
\symb{B_o(a,\rho)}{open ball}{}
\[
\ShowSymbol{open ball}{}=
\{
b\in A:\|b-a\|<\rho
\}
\]
}

\AddEq{closed ball}
{
\symb{B_c(a,\rho)}{closed ball}{}
\[
\ShowSymbol{closed ball}{}=
\{
b\in A:\|b-a\|\le\rho
\}
\]
}

\DefEq
{
$f_n\in M(X,A)$, $n=1$, ...,
}
{fn M(X,A)}

\DefEquation
{
\|f_n(x)-f_m(x)\|<\epsilon
}
{|fn(x)-fm(x)|<e}

\DefEq
{
\|a_p-a_q\|<\epsilon
}
{|ap-aq|<epsilon}

\AddEq{[an] in A}
{
$\{a_n\}$, $a_n\in \Module$,
}

\AddEq{aq in B(ap)}
{
\[a_q\in B_o(a_p,\epsilon)\]
}

\AddEq{an in B(a)}
{
a_n\in B_o(a,\epsilon)
}

\AddEq{limit of sequence}
{
\symb{\lim_{n\rightarrow\infty}a_n}{limit of sequence}{}
\[
a=\ShowSymbol{limit of sequence}{}
\]
}

\AddEq{a=lim an}
{
a=\lim_{n\rightarrow\infty}a_n
}

\DefEq
{
\|a_n-a\|<\epsilon
}
{|an-a|<epsilon}

\DefEquation
{
\lim_{n\rightarrow\infty}a_n=\lim_{n\rightarrow\infty}b_n
}
{lim a=lim b}

\DefEq
{
\lim_{n\rightarrow\infty}(a_n-b_n)=0
}
{lim a-b=0}

\DefEq
{
$g_n\in M(X,A)$, $n=1$, ...,
}
{gn M(X,A)}

\DefEq
{
\[h_n=f_n+g_n\]
}
{hn=fn+gn}

\DefEq
{
h=f+g
}
{h=f+g}

\DefEq
{
$c_1$, $c_2\in A$,
}
{c1 c2 in A}

\DefEq
{
$c_1\in B_c(a_1,R_1)$, $c_2\in B_c(a_2,R_2)$.
}
{c1 c2 in B}

\DefEq
{
$g_{1\cdot n}\in M(X,A_1)$, $n=1$, ...,
}
{g1n M(X,A1)}

\DefEq
{
$g_{2\cdot n}\in M(X,A_2)$, $n=1$, ...,
}
{g2n M(X,A2)}

\DefEq
{
$f_X(g_1)(g_2)$
}
{fX(g1)(g2)}

\DefEq
{
$M(X,A_1)$
}
{M(X,A1)}

\DefEq
{
$M(X,A_2)$
}
{M(X,A2)}

\DefEq
{
$g_1\in M(X,A_1)$
}
{g1 in M(X,A1)}

\DefEq
{
$g_2\in M(X,A_2)$
}
{g2 in M(X,A2)}

\DefEquation
{
\begin{split}
f_X(g_1)(g_2)&:X->A_2
\\
(f_X(g_1)(g_2))(x)&=f(g_1(x))(g_2(x))
\end{split}
}
{f(x)g(x):X->A2}

\DefEq
{
\symb{\|\omega\|}{norm of operation}{}
}
{norm of operation}

\DefEquation
{
\ShowSymbol{norm of operation}{}=
\text{sup}\frac{\| a_1...a_n\omega\|}{\|a_1\|...\|a_n\|}
}
{norm of operation, definition}

\DefEquation
{
\|a_1...a_n\omega\|\le\|\omega\|\|a_1\|...\|a_n\|
}
{|a omega|<|omega||a|1n}

\DefEq
{
\symb{\|f\|}{norm of representation}{}
}
{norm of representation}

\DefEquation
{
\ShowSymbol{norm of representation}{}=
\text{sup}\frac{\|f(a_1)(a_2)\|_2}{\|a_1\|_1\|a_2\|_2}
}
{norm of representation, definition}

\DefEquation
{
\|f(a_1)(a_2)\|_2\le\|f\|\|a_1\|_1\|a_2\|_2
}
{|fab|<|f||a||b|}

\DefEquation
{
c_1+c_2\in B_c(a_1+a_2,R_1+R_2)
}
{c1+c2 in B}

\DefEq
{
$f\in M(X,A)$
}
{f M(X,A)}

\DefEq
{
f(x)=\lim_{n\rightarrow\infty}f_n(x)
}
{f(x)=lim}

\DefEq
{
\[\|f_n(x)-f(x)\|<\epsilon\]
}
{fn(x) - f(x)}

\AddEq[1]{a in U}
{
$a \in U$#1
}

\AddEq{sum mu(F) f}
{
\displaystyle\sum_n\mu(F_n)f_n
}

\AddEq{f1...}
{
\(\{f_1,f_2,...\}\)
}

\AddEq{fn ne fm}
{
\(f_n\ne f_m\)
}

\AddEq{Fn=}
{
F_n=\{x\in X:f(x)=f_n\}
}

\AddEq{mu(X=+Xi)}
{
\mu(\PX)=\sum_i\mu(\PX_i)
}

\AddEq{m A=+m Ai}
{
\[
m(A)=\sum_{i=1}^nm(A_i)
\]
\labelItem{m A=+m Ai}
}

\DefEq
{
\(A\), \(B\in\mathcal R\)
}
{A, B in R}

\DefEq
{
$A\cup B\in\mathcal R$, $A\setminus B\in\mathcal R$.
}
{A cup/minus B R}

\DefEq
{
\symb{\mathcal B(A)}{Borel algebra}1
}
{Borel algebra}

\DefEq
{
$C\in\mathcal C_Y$
\[f^{-1}(C)\in\mathcal C_X\]
}
{CX CY measurable map}

\DefEq
{
$C\in\mathcal B(A)$
\[f^{-1}(C)\in\mathcal C_{\mu}\]
}
{measurable map}

\AddEq{range f}
{
\(y_1\), \(y_2\), ...
}

\DefEq
{
\(m(A)\ge 0\)
\labelItem{m(A)>0}
}
{m(A)>0}

\DefEq
{
\(A\in\mathcal C_m\)
}
{A in Cm}

\DefEq
{
A=\bigcup_{i=1}^nA_i\ \ \ A_i\in\mathcal C_m
}
{A=+Ai Cm}

\DefEq
{
\symb{\mathcal C_{\mu}}
{algebra of sets measurable with respect to measure}1
}
{algebra of sets measurable with respect to measure}

\DefEq
{
\(A\triangle B\), \(A\cap B\in\mathcal R\).
}
{A, B in R 1}

\DefEq
{
$E\in\mathcal R$
}
{E in R}

\DefEq
{
\[A\cap E=A\]
}
{A cap E=A}

\DefEq
{
$A_i\in\mathcal R$, $i=1$, ..., $n$, ...,
}
{Ai in R}

\DefEq
{
$A$, $B$
\begin{align*}
A\cup B&=(A\Delta B)\Delta (A\cap B)
\\
A\setminus B&=A\Delta (A\cap B)
\end{align*}
}
{A cup/minus B}

\DefEq
{
\[\bigcup_nA_n\in\mathcal R\]
}
{Ai in R 1}

\DefEq
{
$B_o(a,\epsilon)\subset U$.
}
{B(a) subset U}

\DefEq
{
$M(X,A)$\Pt
}
{M(X,A)}

\DefEq
{
\symb{M(X,A)}{set of maps to Omega group}1
}
{set of maps to Omega group}

\AddEq{transpose of identity}
{
\UnitNMatrix^T=\UnitNMatrix
}

\AddEquation{double transpose}
{
(a^T)^T=a
}

\DefEq
{
\[
\xymatrix
{
f_k:A\ar[r]|{*}&A_k
}
\]
}
{representation A Ak 1}

\DefEq
{
\[
\xymatrix{
r_1:\Times\ar[r]&S_1
&
r_2:\Times\ar[r]&S_2
}
\]
}
{polymorphisms category}

\DefEq
{
\[
\begin{matrix}
\xymatrix
{
g_1:A\ar[r]|{*}&S_1
}
&
\xymatrix
{
g_2:A\ar[r]|{*}&S_2
}
\end{matrix}
\]
}
{representation of algebra in S1 S2}

\DefEq
{
$r_1\rightarrow r_2$
}
{r1->r2}

\DefEq
{
\[
\xymatrix{
&S_1\ar[dd]^h
\\
\Times
\ar[ru]^{r_1}\ar[rd]_{r_2}
\\
&S_2
}
\]
}
{polymorphisms category, diagram}

\DefEq
{
\symb{\Tensor B}{tensor product}1
}
{tensor product of representations}

\DefEq
{
\[
\begin{matrix}
b_k\in B_k&k=1,...,n&
b_{i\cdot 1},...,b_{i\cdot p}\in B_i&\omega\in\Omega_2(p)&a\in A
\end{matrix}
\]
}
{equivalence, 1, representation, tensor product}

\DefEquation
{
\begin{split}
&b_1\otimes...\otimes(b_{i\cdot 1}...b_{i\cdot p}\omega)\otimes...\otimes b_n
\\
=&(b_1\otimes...\otimes b_{i\cdot 1}\otimes...\otimes b_n)...
(b_1\otimes...\otimes b_{i\cdot p}\otimes...\otimes b_n)
\omega
\end{split}
}
{tensors 1, representation, tensor product}

\DefEquation
{
b_1\otimes...\otimes(f_i(a)\circ b_i)\otimes...\otimes b_n=
f(a)\circ(b_1\otimes...\otimes b_i\otimes...\otimes b_n)
}
{tensors 2, representation, tensor product}

\DefEq
{
\[
f:\Times\rightarrow\Tensor B
\]
}
{map f, 1, representation, tensor product}

\DefEquation
{
f\circ(b_1,...,b_n)=\Tensor b
}
{map f, representation, tensor product}

\DefEq
{
\[
g:\Times\rightarrow V
\]
}
{map g, representation, tensor product}

\DefEq
{
\[
h:\Tensor B\rightarrow V
\]
}
{map h, representation, tensor product}

\AddEq{tower of representations}
{
\[
\xymatrix
{
A_1\ar[r]|{*}^{f_{1\,2}}&A_2\ar[r]|{*}^{f_{2\,3}}
&...\ar[r]|{*}^{f_{n-1\,n}}&A_n
}
\]
}

\AddEquation{morphism of diagram of representations, level k}
{
h_j\circ\BlueText{f_{ij}(a_i)}
=g_{ij}(\RedText{h_i(a_i)})\circ h_j
}

\AddEq{hi:Ai->Bi}
{
\ShowEq{f:A->B}{h_{(i)}}{A_{(i)}}{B_{(i)}}
}

\AddEq[2]{A=A1n}
{
$#1=(#1_{(1)}\ \ ...\ \ #1_{(n)})$#2
}

\AddEq[1]{A=A(1n)1n}
{
\[
#1=(#1_{(1)}\ \ ...\ \ #1_{(n)})=(#1_1\ \ ...\ \ #1_n)
\]
}

\AddEq[2]{A(i)=Ai}
{
$A_{(#1)}=A_{#1}$#2
}

\AddEq{closure operator, representation}
{
\symb{J[f]}{closure operator, representation}1.
}

\AddEq{show closure operator, representation}
{
$\ShowSymbol{subrepresentation generated by set}1$
}

\AddEq[1]{map of words of representation, W}
{
\[
w[f\rightarrow g,#1,#1',R_1]:w[f,#1]\rightarrow w[g,#1']
\]
}

\AddEq[2]{map of words of representation, W in X}
{
\[
#1\in #2=>w[f\rightarrow g,#2,#2',R_1](#1)=R_1(#1)
\]
}

\AddEq[1]{endomorphism based on map of words of representation}
{
\[
R(m)=w[f\rightarrow g,#1,#1',R_1](w[f,#1,m])
\]
}

\AddEq{map of words of representation 2}
{
$m'\in J[g,X']$.
}

\AddEq{map of words of representation, m in X, 1}
{
$m\in X$, $m'=R(m)$,
}

\AddEq{map of words of representation}
{
\[
\wXR:w[f,X]\rightarrow w[g,X']
\]
}

\AddEq{a=a12}
{
$a=(a_1,a_2)$
}

\AddEq{r(a)12}
{
\[
r(a)=(r_1(a_1),r_2(a_2))
\]
}

\AddEq{map of words of representation, m in X}
{
\labelItem{map of words of representation, m in X}
\[
\wXR(m)=m'
\]
}

\AddEq{map of words of representation, omega, 1}
{
\[
m_1, ..., m_n\in w[f,X]
\]
\[
\begin{matrix}
m'_1=\wXR(m_1)&...&m'_n=\wXR(m_n)
\end{matrix}
\]
}

\AddEq{map of words of representation, omega}
{
\labelItem{map of words of representation, omega}
\[
\wXR(m_1...m_n\omega)=m_1'...m_n'\omega
\]
}

\AddEq{map of words of representation, am, 1}
{
\[
\begin{matrix}
m\in w[f,X]& m'=\wXR(m)
&a\in A_1
\end{matrix}
\]
}

\AddEq{map of words of representation, am}
{
\labelItem{map of words of representation, am}
\[
\wXR(f(a)(m))=g(a)(m')
\]
}

\AddEq{generating set of representation}
{
$J[f,X]=A_2$.
}

\AddEq{subrepresentation generated by set}
{
\symb{J[f,X]}{subrepresentation generated by set}1
}

\AddEquation{morphism of diagram of representations, levels k k+1}
{
\BlueText{h_j(f_{ij}(a_i)(a_j))}
=g_{ij}(\RedText{h_i(a_i)})(\BlueText{h_j(a_j)})
}

\AddEq[2]{f:i->*j}
{
\ShowEq{f:A->*B}{f_{#1#2}}{A_{#1}}{A_{#2}}
}

\AddEq{f:ij->*k}
{
\ShowEq{f:A->*B}{f_{ij,k}}{A_i\times A_j}{A_k}
}

\AddEq[3]{Ai->*Aj->*Ak}
{
\xymatrix
{
A_{#1}\ar[r]|{*}^{f_{#1#2}}&A_{#2}\ar[r]|{*}^{f_{#2#3}}&A_{#3}
}
}

\AddEquation{tower of representations, 1 2 3}
{
\xymatrix{
\\
A_k\ar@/_3pc/[rrrr]^{f_{jk}(a_j)}="fkaj"
\ar@/^3pc/[rrrr]^{f_{jk}(\RedText{f_{ij}(a_i)(a_j)})}="fkbj"&&&&A_k
\\
\\
&A_j\ar[rr]^{f_{ij}(a_i)}\ar @{}[rl]="aj"&&A_j\ar @{}[rl]="bj"
\\
\\
&&A_i\ar @{}[rl]="ai"\ar@{=>}[uu]^{f_{ij}}&&
\ar @{=>}^{f_{jk}} "aj";"fkaj"
\ar @{=>}@/_3pc/^{f_{jk}} "bj";"fkbj"
\ar @{=>} "fkaj";"fkbj"_{f_{ijk}(a_i)}="ajbj"
\ar @{=>}@/^6pc/ "ai";"ajbj"^{f_{ijk}}
}
}

\AddEquation{tower of representations, ij,k}
{
\xymatrix{
&&A_i\times A_j\ar@{=>}[d]_{f_{ij,k}}
\\
A_k\ar[rrrr]_{f_{jk}(\RedText{f_{ij}(a_i)(a_j)})}="fkbj"&&&&A_k
\\
\\
&A_j\ar[rr]^{f_{ij}(a_i)}\ar @{}[rl]="aj"&&A_j\ar @{}[rl]="bj"
\\
\\
&&A_i\ar @{}[rl]="ai"\ar@{=>}[uu]^{f_{ij}}&&
\ar @{=>}@/_3pc/^{f_{jk}} "bj";"fkbj"
}
}

\AddEq{f:A->End 2 3}
{
\ShowEq{f:A->B}{f_{ijk}}{A_i}{
\ShowEq{End 2 3}
}
}

\AddEq{End 2 3}
{
\ensuremath{\End(\Omega_j,\End(\Omega_k,A_k))}
}

\AddEquation{define map f13}
{
f_{ijk}(a_i)(\BlueText{f_{jk}(a_j)})
=\BlueText{f_{jk}(\RedText{f_{ij}(a_i)(a_j)})}
}

\AddEq{f:1->*3}
{
\ShowEq{f:A->*B}{f_{ijk}}{A_i}{\End(\Omega_k,A_k)}
}

\AddEq[3]{End Ak}
{
$\End(\Omega_{#2},#1_{#2})$#3
}

\AddEquation{morphism of diagram of representations of F algebra, level k, diagram}
{
\xymatrix{
&A_j\ar[dd]_(.3){f_{ij}(a_i)}
\ar[rr]^{h_j}
&&
B_j\ar[dd]^(.3){g_{ij}(\RedText{h_i(a_i)})}\\
&\ar @{}[rr]|{(1)}&&\\
&A_j\ar[rr]^{h_j}&&B_j\\
A_i\ar[rr]^{h_i}\ar@{=>}[uur]^(.3){f_{ij}}
&&
B_i\ar@{=>}[uur]^(.3){g_{ij}}
}
}

\DefEq
{
\[
\xymatrix{
&\Tensor B\ar[dd]^h
\\
\Times
\ar[ru]^f\ar[rd]_g
\\
&V
}
\]
}
{map gh, representation, tensor product}

\DefEq
{
$A$, $B_1$, ..., $B_n$, $B$
}
{set of universal algebras 1}

\DefEq
{
\[
r_2:B_1\times...\times B_n\rightarrow B
\]
}
{reduced polymorphism of representation}

\AddEq{rc-product}
{
\symb{a\RCstar b}{rc-product}{}
}

\AddEquation{rc-product of matrices}
{
\ShowSymbol{rc-product}{}=
\begin{pmatrix}
\aUD aik\aUD bkj
\end{pmatrix}
}

\AddEquation{entry rc-product of matrices}
{
\aUD {(\ShowSymbol{rc-product}{})}ij=
\aUD aik\aUD bkj
}

\AddEquation{rc-product, matrices}
{
\begin{split}
\PMatrix apn\RCstar\PMatrix bmp
&=
\begin{pmatrix}
\aUD a1k\aUD bk1&...&\aUD a1k\aUD bkm\\
...&...&...\\
\aUD ank\aUD bk1&...&\aUD ank\aUD bkm
\end{pmatrix}
\\ &=
\PMatrix{(\ShowSymbol{rc-product}{})}mn
\end{split}
}

\AddEq{cr-product}
{
\symb{a\CRstar b}{cr-product}{}
}

\AddEquation{cr-product of matrices}
{
\ShowSymbol{cr-product}{}=
\begin{pmatrix}
\aUD aki\aUD bjk
\end{pmatrix}
}

\AddEquation{entry cr-product of matrices}
{
\aUD{(\ShowSymbol{cr-product}{})}ij=
\aUD aki\aUD bjk
}

\AddEquation{cr-product, matrices}
{
\begin{split}
\PMatrix amp\CRstar\PMatrix bpn
&=
\begin{pmatrix}
\aUD ak1\aUD b1k&...&\aUD akm\aUD b1k\\
...&...&...\\
\aUD ak1\aUD bnk&...&\aUD akm\aUD bnk
\end{pmatrix}
\\ &=
\PMatrix{(\ShowSymbol{cr-product}{})}mn
\end{split}
}

\AddEquation{(a+b)=}
{
\aUD{(a+b)}ij=\aUD aij+\aUD bij
}

\AddEquation{transpose of matrix, 1}
{
\aUD {(a^T)}ij=\aUD aji
}

\AddEquation{v=(e)(v)}
{
v=
\RowMatrix en
\ColMatrix vn
=\aD ei\aU vi
}

\AddEquation{v'=fv}
{
\aU {v'}i=\aUD fij\aU vj
}

\AddEquation{v'=vf}
{
\aU {v'}i=\aU vj\aUD fij
}

\AddEquation{v=(v)(e)}
{
v=\ColMatrix vn
\ \ \ \ e=\RowMatrix en
}

\AddEquation{reducibility of products}
{
a\RCstar b
=(\aUD aki\aUD bjk)
=(\aUD bjk\aUD aki)
=b\CRstar a
}

\AddEq{(m+n)a=ma+na,n=k}
{
\[(m+k)a=ma+ka\]
}

\AddEq{(m+n)a=ma+na,n=k+1}
{
\begin{align*}
(m+(k+1))a&=((m+k)+1)a=(m+k)a+a=ma+ka+a\\&=ma+(k+1)a
\end{align*}
}

\AddEq{(m+n)a=ma+na,n=k-1}
{
\begin{align*}
(m+(k-1))a&=((m+k)-1)a=(m+k)a-a=ma+ka-a\\&=ma+(k-1)a
\end{align*}
}

\AddEq{n(a+b)=na+nb,n=k}
{
\[k(a+b)=ka+kb\]
}

\AddEq{n(a+b)=na+nb,n=k+1}
{
\begin{align*}
(k+1)(a+b)&=k(a+b)+a+b=ka+kb+a+b\\&=ka+a+kb+b\\&=(k+1)a+(k+1)b
\end{align*}
}

\AddEq{n(a+b)=na+nb,n=k-1}
{
\begin{align*}
(k-1)(a+b)&=k(a+b)-(a+b)=ka+kb-a-b\\&=ka-a+kb-b\\&=(k-1)a+(k-1)b
\end{align*}
}

\AddEq{(nm)a=n(ma),n=k}
{
\[(km)a=k(ma)\]
}

\AddEq{(nm)a=n(ma),n=k+1}
{
\begin{align*}
((k+1)m)a&=(km+m)a=(km)a+ma=k(ma)+ma\\&=(k+1)(ma)
\end{align*}
}

\AddEq{(nm)a=n(ma),n=k-1}
{
\begin{align*}
((k-1)m)a&=(km-m)a=(km)a-ma=k(ma)-ma\\&=(k-1)(ma)
\end{align*}
}

\AddEq{f(na)=nf(a),n=k}
{
\[f(ka)=kf(a)\]
}

\AddEq{f(na)=nf(a),n=k+1}
{
\begin{align*}
f((k+1)a)&=f(ka+a)=f(ka)+f(a)=kf(a)+f(a)\\&=(k+1)f(a)
\end{align*}
}

\AddEq{f(na)=nf(a),n=k-1}
{
\begin{align*}
f((k-1)a)&=f(ka-a)=f(ka)-f(a)=kf(a)-f(a)\\&=(k-1)f(a)
\end{align*}
}

\AddEq{A1(mA2) symb}
{%
\symb{A_1(\mathcal A_2)}{category of representations}1
}

\AddEq[2]{x in JX}
{
$#1\in J[f,X]$#2
}

\AddEq{word of element relative to generating set, representation}
{
\symb{\wXm}
{word of element relative to generating set, representation}1
}

\AddEq{identify element jf(X) and word}
{
\[
m=\wXm
\]
}

\AddEq[1]{subset of representation}
{
$B\subset J_{#1}[f,X]$
}

\AddEq{subset of words of representation, 2}
{
\[
\wXm[\{m\}]=\{\wXm\}
\]
}

\AddEq{subset of words of representation}
{
\symb{\wXm[B]}{subset of words of representation}{}
\[
\ShowSymbol{subset of words of representation}{}
=\{\wXm:m\in B\}
\]
}

\AddEq{subset of words of representation, 1}
{
\[
\wXm[B]=(\wXm,m\in B)
\]
}

\AddEq{word set of representation}
{
\symb{w[f,X]}{word set of representation}1
}

\AddEquation{ambiguity of coordinates 1}
{
f(a)(m_1...m_n\omega)=(f(a)(m_1))...(f(a)(m_n))\omega
}

\AddEq{ambiguity of coordinates 1, vector}
{
\[
a(b^1e_1+b^2e_2)=(ab^1) e_1+(ab^2) e_2
\]
}

\AddEq{a1n omega x}
{
$f(a_1...a_n\omega)(x)$
}

\AddEq{a1n x omega}
{
$(f(a_1)(x))...(f(a_n)(x))\omega$
}

\AddEq{(a+b)e=ae+be}
{
\[(a+b) e_1=ae_1+b e_1\]
}

\AddEq{l fX in w fX}
{
\[
\lambda[f,X]\subseteq w[f,X]\times w[f,X]
\]
}

\AddEq{w1 w2 in l fX}
{
\[
(w_1[f,X,m],w_2[f,X,m])\in\lambda[f,X]
\]
}

\AddEq{fw1m fw2m in l fX}
{
\[
(f(w_1)(w[f,X,m]),f(w_2)(w[f,X,m]))\in\lambda[f,X]
\]
}

\AddEq{(f(a1n)a2,f(a)a2) in l}
{
\[
(f(a_{11}...a_{1n}\omega)(a_2),(f(a_{11})...f(a_{1n})\omega)(a_2))\in\lambda[f,X]
\]
}

\AddEq{(fa1(a21n),fa1(a2)) in l}
{
\[
(f(a_1)(a_{21}...a_{2n}\omega),f(a_1)(a_{21})...f_(a_1)(a_{2n})\omega)
\in\lambda[f,X]
\]
}

\AddEq[3]{omega in O1AO2}
{
$\omega\in\Omega_{#1}(n)\cap\Omega_{#2}(n)$#3
}

\AddEq{f(a)(v)=av}
{
\[f(a)(v)=av\]
}

\AddEq{f(a+b)v=}
{
\[f(a+b)(v)=f(a)(v)+f(b)(v)\]
}

\AddEq{f(ab)v=}
{
\[f(ab)(v)=f(a)(v)f(b)(v)\]
}

\AddEq{f(a1n)a2=f(a)a2}
{
\[
f(a_{11}...a_{1n}\omega)(a_2)=(f(a_{11})(a_2))...(f(a_{1n})(a_2))\omega
\]
}

\AddEq{(f(a1n)a2,f(a)a2) in l 1}
{
\[
(f(a_{11}...a_{1n}\omega)(a_2),(f(a_{11})(a_2))...(f(a_{1n})(a_2))\omega)\in\lambda[f,X]
\]
}

\AddEq{rho=lambda}
{
\[
\rho[f,\Basis e]=\lambda[f,\Basis e]
\]
}

\AddEq{w1[]=w2[]}
{
\[
w_1[f,X,m]=w_2[f,X,m]
\]
}

\AddEq{rho fX=}
{
\[
\rho[f,X]=
\{
(w[f,X,m],w_1[f,X,m]):m\in A_2
\}
\]
}

\AddEq{ambiguity of coordinates 2}
{
f(a_1...a_n\omega)(x)=(f(a_1)(x))...(f(a_n)(x))\omega
}

\AddEq{A1(mA2)}
{
$A_1(\mathcal A_2)$
}

\AddEq[1]{A1(mA2,x)}
{
$\mathcal A(
\xymatrix
{
A_1\ar[r]|{*}&\mathcal A(\Omega_2,X)
}
)$#1
}

\AddEq{X in A2}
{
$X\subset A_2$,
}

\AddEq{product in category}
{
\symb[Pi-1]{\prod_{\iI}B_i}{product in category}{}
\[P=\ShowSymbol{product in category}{}\]
}

\AddEq{coproduct in category}
{
\symb[Pi-1]{\coprod_{\iI}B_i}{coproduct in category}{}
\[P=\ShowSymbol{coproduct in category}{}\]
}

\AddEq{coproduct in category, 1 n}
{
\symb[Pi-1]{\coprod_{i=1}^nB_i}{coproduct in category}{i 1 n}
\symb[B-0]{B_1\coprod...\coprod B_n}{coproduct in category}{1 n}
\[
P=\ShowSymbol{coproduct in category}{i 1 n}
=\ShowSymbol{coproduct in category}{1 n}
\]
}

\AddEq [2]{av=sum av}
{
\[#1^i#2_i=\sum_{\iI}#1^i#2_i\]
}

\AddEquation{f(a)=o left}
{
f(a_2)=f\Act a_2
}

\AddEquation{f(a)=o right}
{
f(a_2)=a_2\Act f
}

\AddEquation{f(a)=* right}
{
f(a)=a\Act f
}

\AddEq{A12->A2 left}
{
\[
(a_1,a_2)\in A_1\times A_2\rightarrow a_1*a_2\in A_2
\]
}

\AddEq{A12->A2 right}
{
\[
(a_1,a_2)\in A_1\times A_2\rightarrow a_2*a_1\in A_2
\]
}

\AddEq{A12->A2 1 left}
{
\[
(a_1,a_2)\in A_1\times A_2\rightarrow a_1*a_2\in A_2
\]
}

\AddEq{A12->A2 1 right}
{
\[
(a_1,a_2)\in A_1\times A_2\rightarrow a_2*a_1\in A_2
\]
}

\AddEquation{fga= o left}
{
(f\circ g)\circ a=f\circ (g\circ a)
}

\AddEquation{fga= o right}
{
a\circ (g\circ f)=(a\circ g)\circ f
}

\DefEq
{
\[f\circ g\circ a=f\circ (g\circ a)=(f\circ g)\circ a\]
}
{fga= 1 left}

\DefEq
{
\[a\circ g\circ f=(a\circ g)\circ f=a\circ (g\circ f)\]
}
{fga= 1 right}

\AddEquation{right shift, product}
{
R(b*c)=R(b)\circ R(c)
}

\AddEquation{left shift, product}
{
L(c*b)=L(c)\circ L(b)
}

\AddEq{right shift =}
{
\symb{R(b)}{right shift}{}
}

\AddEq{right shift}
{
\[a'=a\circ\ShowSymbol{right shift}{}=a*b\]
}

\AddEq{left shift =}
{
\symb{L(b)}{left shift}{}
}

\AddEq{left shift}
{
\[a'=\ShowSymbol{left shift}{}\circ a=b*a\]
}

\DefEq
{
\[
(x_1, ..., x_n)\in B_1\times...\times B_n
\rightarrow x_1\otimes...\otimes x_n\in B_1\otimes...\otimes B_n
\]
}
{B times->B otimes}

\DefEq
{
f:A_1\rightarrow\End(\Omega_2,A_2)
}
{representation of group, map}

\DefEq
{
\[
(\id:A_1\rightarrow A_1,\ r_2:A_2\rightarrow B_2)
\]
}
{id:A1->A1 A2->B2}

\DefEquation
{
\xymatrix{
&A_2\ar[dd]_(.3){f(a)}\ar[rrr]^{r_2}&&&B_2\ar[dd]^(.3){g(a)}\\
&&&&\\
&A_2\ar[rrr]_(.65){r_2}&&&B_2\\
A_1\ar@{=>}[uur]^(.3)f\ar@{=>}[uurrrr]^(.7)g
}
}
{morphism id,R of representations}

\AddEquation{morphism of representations of universal algebra}
{
r_2\circ f(a)=g(a)\circ r_2
}

\DefEq
{
\[
\xymatrix{
A_2\ar[rr]^{r_2}&&B_2\\
&A_1\ar[lu]|{*}^f\ar[ru]|{*}_g
}
\]
}
{morphism id,R of representations 2}

\AddEq{Gen f =}
{
\[
Gen[f]=\{X\subseteq A_2:J[f,X]=A_2\}
\]
}

\AddEq[3]{in Gen}
{
$#1\in Gen[#2]$#3
}

\AddEq[3]{f:A->*B}
{
\[
\xymatrix@C=30pt
{
#1:#2\ar[r]|-{*}&#3
}
\]
}

\AddEq[3]{A->*B}
{
\xymatrix
{
#1\ar[r]|-{*}&#2
}
#3
}

\AddEq{action of rational integers in Abelian group}
{
\begin{align}
\EqLabel{0g=0}
0g&=0
\\
\EqLabel{(n+1)g=}
(n+1)g&=ng+g
\\
\EqLabel{(n-1)g=}
(n-1)g&=ng-g
\end{align}
}

\AddEq{representation of Z in G}
{
\begin{align}
1a&=a
\EqLabel{1a=a}
\\
\EqLabel{(nm)a=n(ma)}
(nm)a&=n(ma)
\\
\EqLabel{(m+n)a=ma+na}
(m+n)a&=ma+na
\\
\EqLabel{(m-n)a=ma-na}
(m-n)a&=ma-na
\\
\EqLabel{n(a+b)=na+nb}
n(a+b)&=na+nb
\end{align}
}

\AddEquation{(k+n)a-na=ka}
{
(k+n)a-na=ka
}

\AddEq{phi(n):G->G}
{
\[
\varphi(n):a\in G\rightarrow na\in G
\]
}

\AddEq{S=si}
{
$S=\{s_i:\iI\}$
}

\AddEquation{J(S)=gi si}
{
J(S)=\left\{g:g=\sum_{\iI}g^is_i, g^i\in Z\right\}
}

\AddEq[1]{|gi ne 0|}
{
$\{\iI:#1^i\ne 0\}$
}

\AddEq{category A(Bi)}
{
$\mathcal A(B_i,\iI)$
}

\AddEq{category A[Bi]}
{
$\mathcal A[B_i,\iI]$
}

\AddEq{gi=dik}
{
$g^i=\delta^i_k$.
}

\AddEquation{sk=sum si}
{
s_k=\sum_{\iI}g^is_i
}

\AddEq{sk in JS}
{%
$s_k\in J(S)$
}

\AddEquation{(ai)+(bi)=(ai+bi)}
{
(a_i,\iI)+(b_i,\iI)=(a_i+b_i,\iI)
}

\AddEq{phi:Z->End(G)}
{
\[
\varphi:Z\rightarrow \End(Ab,G)
\]
}

\AddEq{right-side representation}
{
\[a_2'=a_2\Act f(a_1)=a_2*a_1\]
}

\AddEq{left-side representation}
{
\[a_2'=f(a_1)\Act a_2=a_1*a_2\]
}

\AddEq{right-side representation of group}
{
a_2\Act f(a_1*b_1)=a_2\Act(f(a_1)\Act f(b_1))
}

\AddEq{left-side representation of group}
{
f(a_1*b_1)\Act a_2=(f(a_1)\Act f(b_1))\Act a_2
}

\AddEq{* right-side representation of group}
{
a_2\Act f(a_1*b_1)=a_2*(a_1*b_1)
}

\AddEq{* left-side representation of group}
{
f(a_1*b_1)\Act a_2=(a_1*b_1)*a_2
}

\AddEq{o right-side representation of group}
{
a_2*(a_1*b_1)=(a_2* a_1)* b_1
}

\AddEq{o left-side representation of group}
{
(a_1*b_1)* a_2=a_1*( b_1* a_2)
}

\AddEq{left-side 1 representation of group}
{
f(a_1*b_1)\circ a_2=f(a_1)\circ f(b_1)\circ a_2
}

\AddEq{right-side 1 representation of group}
{
a_2\circ f(a_1*b_1)=a_2\circ f(a_1)\circ f(b_1)
}

\AddEq{morphism of left representation of Omega group}
{
\[
\BlueText{r_2\circ(a_1*a_2)}=\RedText{r_1(a_1)}*(\BlueText{r_2\circ a_2})
\]
}

\AddEq{morphism of right representation of Omega group}
{
\[
\BlueText{r_2\circ(a_2*a_1)}=(\BlueText{r_2\circ a_2})*\RedText{r_1(a_1)}
\]
}

\AddEq{orbit of left-side representation}
{
\symb{A_1*a_2}{orbit of representation}{}
\[
\ShowSymbol{orbit of representation}{}=\{b_2=a_1*a_2:a_1\in A_1\}
\]
}

\AddEq{orbit of right-side representation}
{
\symb{a_2*A_1}{orbit of representation}{}
\[
\ShowSymbol{orbit of representation}{}=\{b_2=a_2*a_1:a_1\in A_1\}
\]
}

\AddEq{space of orbits of representation}
{
\symb{A_2/A_1}{space of orbits of representation}1
}



\AddEq{r2a=r2 o a}
{
\[r_2(a_2)=r_2\circ a_2\]
}

\AddEq[2]{passive and active maps}
{
\[
\xymatrix{
\Basis e_1\in\mathcal B[f]\ar[d]_{#2}\ar[rr]^{#1}&&\Basis e_3\in\mathcal B[f]\ar[d]_{#2}
\\
\Basis e_2\in\mathcal B[f]\ar[rr]^{#1}&&\Basis e_4\in\mathcal B[f]
}
\]
}

\AddEq{sigma->+-}
{
\symb{|\sigma|}{parity of permutation}{}
\[
\ShowSymbol{parity of permutation}{}
:\sigma\in S(n)\rightarrow \{-1,1\}
\]
}

\AddEquation{parity of permutation}
{
\ShowSymbol{parity of permutation}{}=
\left\{
\begin{matrix}
1&\textrm{\permutation}\sigma\textrm{ \even}
\\
-1&\textrm{\permutation}\sigma\textrm{ \odd}
\end{matrix}
\right.
}

\AddEq{Basis e4=Basis e5}
{
$\Basis e_4=\Basis e_5$.
}

\AddEq{passive transformation, 1}
{
\Basis e_5
=\Basis e_3 \circ W[f,\Basis e_1,\Basis e_2]
=W[f,\Basis e_1,\Basis e_3]\circ \Basis e_1\circ W[f,\Basis e_1,\Basis e_2]
}

\AddEq[2]{Omega_2 word, basis e}
{
#2=\Basis e_{#1}\circ W[f,\Basis e_{#1},#2]
}

\AddEq[2]{coordinate transformation in representation, 1}
{
W[f,\Basis e_1,#2]=W[f,\Basis e_1,\Basis e_1\circ #1]\circ W[f,\Basis e_2,#2]
}

\AddEq[2]{Omega word in representation, basis Y, 1}
{
\begin{split}
\Basis e_1\circ W[f,\Basis e_1,#2]&=\Basis e_2\circ W[f,\Basis e_2,#2]
=\Basis e_1\circ W[f,\Basis e_1,\Basis e_2]\circ W[f,\Basis e_2,#2]
\\
&=\Basis e_1\circ W[f,\Basis e_1,\Basis e_1\circ #1]\circ W[f,\Basis e_2,#2]
\end{split}
}

\AddEq[1]{passive transformation S}
{
\Basis e_2=\Basis e_1\circ #1=\Basis e_1\circ W[f,\Basis e_1,\Basis e_1\circ #1]
}

\AddEq[2]{coordinate transformation}
{
W[f,\Basis e_2,#1]=W[f,\Basis e_1\circ #2,\Basis e_1]\circ W[f,\Basis e_1,#1]
}

\AddEq[2]{coordinate transformation in representation, TS 1}
{
$#2\circ #1$.
}

\AddEq[3]{coordinate transformation in representation, TS}
{
\begin{align*}
W[f,\Basis e_3,#1]
&=W[f,\Basis e_1\circ #3,\Basis e_1]\circ W[f,\Basis e_1\circ #2,\Basis e_1]\circ W[f,\Basis e_1,#1]
\\
&=W[f,\Basis e_1\circ #3\circ #2,\Basis e_1]\circ W[f,\Basis e_1,#1]
\end{align*}
}

\AddEq[1]{coordinates of geometric object}
{
\symb[O-0]{\mathcal O(f,g,\Basis e_g,#1)}{coordinates of geometric object}{}
}

\AddEq[2]{show coordinates of geometric object}
{
\begin{align*}
\ShowSymbol{coordinates of geometric object}{}
&=H(GA(f))\circ W[g,\Basis e_g,#2]
\\
&=(W[g,\Basis e_g\circ H(#1),\Basis e_g]\circ W[g,\Basis e_g,#2],\Basis e_f\circ #1,#1\in GA(f))
\end{align*}
}

\AddEq[1]{geometric object 2, g}
{
\[
#1=\Basis e_g\circ W[g,\Basis e_g,#1]
\]
}

\AddEq[2]{invariance principle 3 algebra}
{
\begin{align*}
#1'&=\Basis e_{g1}\circ W[g,\Basis e_{g1},#1']
\\
&=\Basis e_g\circ W[g,\Basis e_g,\Basis e_g\circ H(#2)]
\circ W[g,\Basis e_g\circ H(#2),\Basis e_g]\circ W[g,\Basis e_g,#1]
\\
&=\Basis e_g\circ W[g,\Basis e_g,#1]=#1
\end{align*}
}

\AddEq[1]{invariance principle 2 passive}
{
\[\Basis e_{f1}=\Basis e_f\circ #1\]
}

\AddEq[1]{geometric object}
{
\symb[O-0]{\mathcal O(f,g,#1)}{geometric object}{}%
}

\AddEq[2]{show geometric object}
{
\[
\ShowSymbol{geometric object}{}=
(W[g,\Basis e_g\circ H(#1),\Basis e_g]\circ W[g,\Basis e_g,#2],
\Basis e_g\circ H(#1),\Basis e_f\circ #1,#1\in GA(f))
\]
}

\AddEq[1]{geometric object 1}
{
$\Basis e_{f1}=\Basis e_f\circ #1$
}

\AddEq{geometric object 2}
{
$\Basis e_{f1}$.
}

\AddEq[1]{passive transformation of representation g}
{
\Basis e_{g1}=\Basis e_g\circ H(#1)
}

\AddEq[2]{coordinate transformation g}
{
W[g,\Basis e_{g1},#1]=W[g,\Basis e_g\circ H(#2),\Basis e_g]\circ W[g,\Basis e_g,#1]
}

\AddEq[1]{invariance principle 4}
{
$W[g,\Basis e_g,#1]$.
}

\AddEq{passive representation coordinated with passive representation}
{
\[
\xymatrix{
\mathcal \End(B[f])\ar[rr]^H
&&
\mathcal \End(B[g])
\\
GA(f)\ar[u]^{P(f)}\ar[rr]_h\ar[urr]_{f'}
&&
GA(g)\ar[u]_{P(g)}
}
\]
}

\AddEq{active transformation, 1}
{
\Basis e_4
=W[f,\Basis e_1,\Basis e_3]\circ \Basis e_2
=W[f,\Basis e_1,\Basis e_3]\circ \Basis e_1\circ W[f,\Basis e_1,\Basis e_2]
}

\AddEq{basis manifold}
{
\symb{\mathcal B[f]}{basis manifold}1
}

\AddEq{passive representation in basis manifold def}
{
\ShowEq{f:A->*B}{\ShowSymbol{passive representation in basis manifold}{}}
{GA(f)}{\mathcal B[f]}
}

\AddEq{passive representation in basis manifold}
{
\symb{P(f)}{passive representation in basis manifold}{}
}

\AddEq[1]{passive transformation of basis def}
{
\[
#1(\Basis e)=
\ShowSymbol{passive transformation of basis}{}
\]
}

\AddEq[1]{passive transformation of basis}
{
\symb{\ecS[#1]}{passive transformation of basis}{}
}

\AddEq[1]{passive transformation}
{
\begin{split}
#1:\Basis h&\rightarrow\Basis h\circ #1
\\
\Basis h\circ #1&=\Basis h\circ W[f,\Basis e,\ecS[#1]]
\end{split}
}

\AddEq[1]{active transformation}
{
\begin{split}
#1:\Basis h&\rightarrow #1\circ\Basis h
\\
#1\circ\Basis h&
=W[f,\Basis e,\Rce[#1]]\circ\Basis h
\end{split}
}

\AddEq[3]{e2=e1 o S}
{
$\Basis e_{#2}=\Basis e_{#3}\circ #1$.
}

\AddEq[3]{e3=R o e1}
{
$\Basis e_{#2}=\Rce[#1]_{#3}$.
}

\DefText{two matrix products}
{
\ShowDefinition{row over column product}

\ShowDefinition{column over row product}
}

\DefText{matrix operations}
{
\ShowText{matrix operations 1}

\ShowText{two matrix products}

\ePrints{2020.06.01,2204.06320,2022.01.05}
\ifx\Semafor\ValueOff
\ShowText{matrix operations 2}

\ShowDefinition{transpose of matrix}

\ePrints{1908.04418,6860-2955}
\ifx\Semafor\ValueOff
\ShowTheorem{double transpose is original matrix}
\ShowProof{double transpose is original matrix}
\fi

\ShowDefinition{sum of matrices}
\fi

\ShowEq{remark: pronunciation of product}

\ePrints{5284-0163,1801.01628,2207.06506,8428-0408,MAlgebra2,5148-4632}
\ifx\Semafor\ValueOn
\ShowEq{product of matrices is distributive over addition}
\fi

\ePrints{2020.06.01,2022.01.05}
\ifx\Semafor\ValueOff
\TwoColText
{
\ShowTheorem{rc transpose}
}
{
\ShowTheorem{cr transpose}
}
\ShowProof{rc transpose}
\ShowProof{cr transpose}
\fi
}

\AddEquation{cr transpose, 1}
{
(a\CRstar b)^T=((a^T)^T\CRstar (b^T)^T)^T
}

\AddEquation{cr transpose, 2}
{
(a\CRstar b)^T=((a^T\RCstar b^T)^T)^T
}

\DefText{matrix over Omega ring}
{
\ShowDefinition{matrix over Omega ring}

\ShowConvention{Einstein summation}

\ShowText{matrix operations}

\ShowDefinition{biring}

\ShowTheorem{duality principle for biring}

\ShowTheorem{duality principle for biring of matrices}

\ShowEq{power matrix}

\ShowText{matrix product}

\ShowRemark{reducible biring}
}

\DefText{matrix product}
{
\TwoColText
{
\ShowEq{def rc}
\ShowTheorem{two products equal}
\begin{sloppypar}
\ShowProof{two products equal}
\end{sloppypar}
}
{
\ShowEq{def cr}
\ShowTheorem{two products equal}
\begin{sloppypar}
\ShowProof{two products equal}
\end{sloppypar}
}

\TwoColText
{
\ShowEq{def rc}
\ProveTheorem{inverse matrix of product}
}
{
\ShowEq{def cr}
\ProveTheorem{inverse matrix of product}
}
}

\AddEq{delta=Matrix}
{
\gdef\UnitNMatrix{\delta}
$\UnitNMatrix=(\aUD {\delta}ij)$
}

\AddEq{rc-inverse matrix}
{
a\RCstar a^{\RCInverse}=\UnitNMatrix
}

\AddEq{cr-inverse matrix}
{
a\CRstar a^{\CRInverse}=\UnitNMatrix
}

\ePrints{1908.04418,6860-2955}%
\ifx\Semafor\ValueOff%
\AddEquation{v=vi ei 123, left-cols}
{
\begin{aligned}
\Vector v&=v_1\CRstar\Basis e_1=v_2\CRstar\Basis e_2\\&=v_3\CRstar\Basis e_3
\end{aligned}
}
\else
\AddEquation{v=vi ei 123, left-cols}
{
\Vector v=v_1\CRstar\Basis e_1=v_2\CRstar\Basis e_2=v_3\CRstar\Basis e_3
}
\fi

\AddEquation{v=vi ei 123, right-cols}
{
\begin{aligned}
\Vector v&=\Basis e_1\RCstar v_1=\Basis e_2\RCstar v_2\\&=\Basis e_3\RCstar v_3
\end{aligned}
}

\AddEquation{v=vi ei 123, left-rows}
{
\begin{aligned}
\Vector v&=v_1\RCstar\Basis e_1=v_2\RCstar\Basis e_2\\&=v_3\RCstar\Basis e_3
\end{aligned}
}

\AddEquation{v=vi ei 123, right-rows}
{
\begin{aligned}
\Vector v&=\Basis e_1\CRstar v_1=\Basis e_2\CRstar v_2\\&=\Basis e_3\CRstar v_3
\end{aligned}
}

\AddEquation{v=vi ei 12, left-cols}
{
\Vector v=v_1\CRstar\Basis e_1=v_2\CRstar\Basis e_2
}

\AddEquation{v=vi ei 12, right-cols}
{
\Vector v=\Basis e_1\RCstar v_1=\Basis e_2\RCstar v_2
}

\AddEquation{v=vi ei 12, left-rows}
{
\Vector v=v_1\RCstar\Basis e_1=v_2\RCstar\Basis e_2
}

\AddEquation{v=vi ei 12, right-rows}
{
\Vector v=\Basis e_1\CRstar v_1=\Basis e_2\CRstar v_2
}

\AddEq{Basis e i=12}
{
$\Basis e_i$, $i=1$, $2$.
}

\AddEq{Basis e i=13}
{
$\Basis e_i$, $i=1$, $2$, $3$.
}

\AddEq[4]{e2i=aij e1j cr-cols}
{
\Basis e_{#2}=#3_{#4}\CRstar\Basis e_{#1}
}

\AddEq[4]{e2i=aij e1j rc-cols}
{
\Basis e_{#2}=\Basis e_{#1}\RCstar #3_{#4}
}

\AddEq[4]{e2i=aij e1j cr-rows}
{
\Basis e_{#2}=\Basis e_{#1}\CRstar #3_{#4}
}

\AddEq[4]{e2i=aij e1j rc-rows}
{
\Basis e_{#2}=#3_{#4}\RCstar\Basis e_{#1}
}

\AddEq[1]{Basis e in BVG}
{
$\Basis e_{#1}\in
\ShowEq{basis manifold of V, \SideNS-\Cols}eG{}$
}

\AddEquation{e3=g12 cr e1 left-cols}
{
\Basis e_3=(g_2\CRstar g_1)\CRstar\Basis e_1
}

\AddEquation{e3=g12 cr e1 right-cols}
{
\Basis e_3=\Basis e_1\RCstar(g_1\RCstar g_2)
}

\AddEquation{e3=g12 cr e1 left-rows}
{
\Basis e_3=(g_2\RCstar g_1)\RCstar\Basis e_1
}

\AddEquation{e3=g12 cr e1 right-rows}
{
\Basis e_3=\Basis e_1\CRstar(g_1\CRstar g_2)
}

\AddEquation{v1 cr e1=v2 cr g cr e1, left-cols}
{
v_1\CRstar\Basis e_1=v_2\CRstar g\CRstar\Basis e_1
}

\AddEquation{v1 cr e1=v2 cr g cr e1, right-cols}
{
\Basis e_1\RCstar v_1=\Basis e_1\RCstar g\RCstar v_2
}

\AddEquation{v1 cr e1=v2 cr g cr e1, left-rows}
{
v_1\RCstar\Basis e_1=v_2\RCstar g\RCstar\Basis e_1
}

\AddEquation{v1 cr e1=v2 cr g cr e1, right-rows}
{
\Basis e_1\CRstar v_1=\Basis e_1\CRstar g\CRstar v_2
}

\AddEq[6]{v1=v2*a left-cols}
{
#1_{#2}=#3_{#4}\CRstar #5_{#6}
}

\AddEq[6]{v1=v2*a right-cols}
{
#1_{#2}=#5_{#6}\RCstar #3_{#4}
}

\AddEq[6]{v1=v2*a left-rows}
{
#1_{#2}=#3_{#4}\RCstar #5_{#6}
}

\AddEq[6]{v1=v2*a right-rows}
{
#1_{#2}=#5_{#6}\CRstar #3_{#4}
}

\AddEq[4]{v1g- left-cols}
{
#1_{#2}\CRstar #3_{#4}^{\CRInverse}
}

\AddEq[4]{v1g- right-cols}
{
#3_{#4}^{\RCInverse}\RCstar #1_{#2}
}

\AddEq[4]{v1g- left-rows}
{
#1_{#2}\RCstar #3_{#4}^{\RCInverse}
}

\AddEq[4]{v1g- right-rows}
{
#3_{#4}^{\CRInverse}\CRstar #1_{#2}
}

\AddEq[8]{v2=v1g-}
{
#1_{#2}=
\ShowEq{v1g- #7-#8}{#3}{#4}{#5}{#6}
}

\AddEquation{v3=v1 cr g21- left-cols}
{
v_3=v_1\CRstar (g_2\CRstar g_1)^{\CRInverse}
}

\AddEquation{v3=v1 cr g21- right-cols}
{
v_3=(g_2\RCstar g_1)^{\RCInverse}\RCstar v_1
}

\AddEquation{v3=v1 cr g21- left-rows}
{
v_3=v_1\RCstar (g_2\RCstar g_1)^{\RCInverse}
}

\AddEquation{v3=v1 cr g21- right-rows}
{
v_3=(g_2\CRstar g_1)^{\CRInverse}\CRstar v_1
}

\AddEquation{v3=v1 cr g1- cr g2- left-cols}
{
v_3=v_1\CRstar g_1^{\CRInverse}\CRstar g_2^{\CRInverse}
}

\AddEquation{v3=v1 cr g1- cr g2- right-cols}
{
v_3=g_2^{\RCInverse}\RCstar g_1^{\RCInverse}\RCstar v_1
}

\AddEquation{v3=v1 cr g1- cr g2- left-rows}
{
v_3=v_1\RCstar g_1^{\RCInverse}\RCstar g_2^{\RCInverse}
}

\AddEquation{v3=v1 cr g1- cr g2- right-rows}
{
v_3=g_2^{\CRInverse}\CRstar g_1^{\CRInverse}\CRstar v_1
}

\AddEq{rc power, 0}
{
a^{\RCPower 0}=\UnitNMatrix
}

\AddEq{rc power}
{
\symb{a^{\RCPower k}}{rc power}{}
}

\AddEq{rc-power, n}
{
\ShowSymbol{rc power}{}=a^{\RCPower{k-1}}\RCstar a
}

\AddEq{cr power, 0}
{
a^{\CRPower 0}=\UnitNMatrix
}

\AddEq{cr power}
{
\symb{a^{\CRPower k}}{cr power}{}
}

\AddEq{cr-power, n}
{
\ShowSymbol{cr power}{}=a^{\CRPower{k-1}}\CRstar a
}

\AddEq{Hadamard inverse of matrix 3}
{
\[
(\aUD aji)^{-1}=\frac 1{\aUD aij}
\]
}

\AddEquation{Hadamard inverse of matrix 2}
{
\ShowSymbol{Hadamard inverse of matrix}{}
=
\begin{pmatrix}
(\aUD a11)^{-1}&...&(\aUD an1)^{-1}\\
...&...&...\\
(\aUD a1n)^{-1}&...&(\aUD ann)^{-1}
\end{pmatrix}
}

\AddEquation{Hadamard inverse of matrix 2 entry}
{
\aUD{(\ShowSymbol{Hadamard inverse of matrix}{})}ij
=(\aUD {(a^T)}ij)^{-1}=(\aUD aji)^{-1}
}

\AddEq{Hadamard inverse of matrix}
{
\symb{\mathcal{H}a}{Hadamard inverse of matrix}{}
}

\AddEq{RC-quasideterminant definition}
{
\symb{\RCDet\,a}{RC-quasideterminant definition}{}
}

\AddEquation{RC-quasideterminant definition=}
{
\begin{split}
\ShowSymbol{RC-quasideterminant definition}{}&=
\begin{pmatrix}
\RCdet 11a&...&\RCdet 1na\\
...&...&...\\
\RCdet n1a&...&\RCdet nna
\end{pmatrix}
\\ &=
\begin{pmatrix}
(\aUD{(a^{\RCInverse})}11)^{-1}&...&(\aUD{(a^{\RCInverse})}n1)^{-1}\\
...&...&...\\
(\aUD{(a^{\RCInverse})}1n)^{-1}&...&(\aUD{(a^{\RCInverse})}nn)^{-1}
\end{pmatrix}
\end{split}
}

\AddEq{j i RC-quasideterminant definition}
{
\symb{\RCdet jia}{j i RC-quasideterminant definition}{}
}

\AddEquation{j i RC-quasideterminant =}
{
\ShowSymbol{j i RC-quasideterminant definition}{}
=\minorD
-\minorC\RCstar
\left(\minorA\right)^{\RCInverse}\RCstar
\minorB
=(\aUD{(a^{\RCInverse})}ij)^{-1}
}

\AddEq{inverse matrix 2x2, 0}
{
\[
a=
\begin{pmatrix}
\aUD a11&\aUD a12
\\
\aUD a21&\aUD a22
\end{pmatrix}
\]
}

\AddEq{quasideterminant matrix 2x2 def}
{
\def\aAA{(\aUD a11)^{-1}}
\def\aBA{(\aUD a21)^{-1}}
\def\aAB{(\aUD a12)^{-1}}
\def\aBB{(\aUD a22)^{-1}}
\def\BAA{\aUD a11-\aUD a12\aBB\aUD a21}
\def\BAB{\aUD a12-\aUD a11\aBA\aUD a22}
\def\BBA{\aUD a21-\aUD a22\aAB\aUD a11}
\def\BBB{\aUD a22-\aUD a21\aAA\aUD a12}
\def\CAA{\aUD a11-\aUD a21\aBB\aUD a12}
\def\CAB{\aUD a12-\aUD a22\aBA\aUD a11}
\def\CBA{\aUD a21-\aUD a11\aAB\aUD a22}
\def\CBB{\aUD a22-\aUD a12\aAA\aUD a21}
\def\ABA{(\BAB)^{-1}}
\def\DAA{(\CAA)^{-1}}
\def\DBA{(\CAB)^{-1}}
\def\DAB{(\CBA)^{-1}}
\def\DBB{(\CBB)^{-1}}
\def\UMBA{\aUD{(a^{\RCInverse})}21}
\def\UMAB{\aUD{(a^{\RCInverse})}12}
}

\AddEquation{quasideterminant rc, matrix 2x2}
{
\RCDet a=
\begin{pmatrix}
\BAA&\BAB
\\
\BBA&\BBB
\end{pmatrix}
\ifx\texFuture\Defined
\aUD a11-\aUD a12(\aUD a22)^{-1}\aUD a21&\aUD a12-\aUD a11(\aUD a21)^{-1}\aUD a22
\\
\aUD a21-\aUD a22(\aUD a12)^{-1}\aUD a11&\aUD a22-\aUD a21(\aUD a11)^{-1}\aUD a12
\fi
}

\AddEquation{diagram of representations, left module nonassociative 1}
{
\begin{matrix}
\vcenter
{
\xymatrix
{
A\ar[rr]|{*}^{g_{34}}&&V
\\
&D\ar[ur]|{*}_{g_{14}}\ar[ul]|(.4){*}^{g_{12}}
}
}
&
\begin{array}{r@{\,}l}
g_{12}(d):a&\rightarrow d\,a
\\
g_{34}(a):v&\rightarrow a\,v
\\
g_{14}(d):v&\rightarrow d\,v
\end{array}
\end{matrix}
}

\AddEquation{quasideterminant cr, matrix 2x2}
{
\CRDet a=
\begin{pmatrix}
\CAA&\CAB
\\
\CBA&\CBB
\end{pmatrix}
}

\AddEquation{rc inverse matrix 2x2}
{
a^{\RCInverse}=
\begin{pmatrix}
\DAA&\DAB
\\
\DBA&\DBB
\end{pmatrix}
}

\AddEquation{orbit, proposition}
{
b_2\in A_1*a_2
}

\AddEquation{A1a2=A1b2}
{
A_1*a_2=A_1*b_2
}

\AddEq{AbS ff'g}
{
\xymatrix
{
G\ar[rr]^g&&G'\\
&S\ar[ul]^f\ar[ur]_{f'}
}
}

\AddEq{S in G}
{
$S\subseteq G$.
\labelItem{S in G}
}

\AddEq{h12 in G}
{
$h_1$, $h_2\in G$.
}

\AddEq{H12=}
{
\begin{align*}
H_1&=\{x\in S:h_1(x)\ne 0\}
\\
H_2&=\{x\in S:h_2(x)\ne 0\}
\end{align*}
}

\AddEq{H1vH2}
{
\[
\{x\in S:h_1(x)+h_2(x)\ne 0\}\subseteq H_1\cup H_2
\]
}

\AddEq{h1+h2 in G}
{
$h_1+h_2\in G$.
}

\AddEq{k*x=k}
{
\[
h(y)=\left\{
\begin{matrix}
k&y=x\\
0&y\ne x
\end{matrix}
\right.
\]
}

\AddEq{f:S->G}
{
\[f:x\in S\rightarrow 1*x\in G\]
}

\AddEquation{h=k*x}
{
h=\sum_{x\in S}k_xx
}

\AddEq{kx=k*x}
{
\[kx=k(1*x)=k*x\]
}

\AddEq{|kx ne 0|}
{
$\{x:k_x\ne 0\}$
}

\AddEq{k in Z x in S}
{
$k_x\in Z$, $x \in S$,
}

\AddEquation{h=k12*x i}
{
\sum_{i=1}^nk_1^i*x_i=\sum_{i=1}^nk_2^i*x_i
}

\AddEquation{k1-k2 *x}
{
\sum_{i=1}^n(k_1^i-k_2^i)*x_i=0
}

\AddEq{k1=k2 i}
{
$k_1^i=k_2^i$, $i=1$, ..., $n$,
}

\AddEquation{g(1x)=f'(x)}
{
g(1x)=f'(x)
}

\AddEquation{f(na)=nf(a)}
{
f(na)=nf(a)
}

\AddEquation{g(x)=}
{
g(h)=
g\left(\sum_{x\in S}k_xx\right)=\sum_{x\in S}g(k_xx)=\sum_{x\in S}k_xg(1x)=\sum_{x\in S}k_xf'(x)
}

\AddEq{x in -> h(x)}
{
$\{x\in S:h(x)\ne 0\}$
}

\AddEq{category AbS}
{
$Ab(S)$
}

\DefText{representation of group}
{
\ShowEq{def left}
\ShowDefinition{side representation of group}
\ShowText{side representation of group}
\ShowRemark{morphism of representation of Omega group}
\ShowEq{def right}
\ShowDefinition{side representation of group}
\ShowText{side representation of group}
\ShowRemark{morphism of representation of Omega group}
\ShowDefinition{representation of Abelian group}
}

\AddEquation{Left side contravariant representation}
{
a_1^{-1}*(b_1^{-1}*a_2)=(b_1*a_1)^{-1}*a_2
}

\AddEq{Left side covariant representation}
{
\[
a_1*(b_1*a_2)=(a_1*b_1)*a_2
\]
}

\AddEq{two representations of group}
{
\[
\xymatrix{
A_2\ar[rr]^{h(a_1)}\ar[d]^{f(b_1)} & & A_2\ar[d]^{f(b_1)}\\
A_2\ar[rr]_{h(a_1)}& &A_2
}
\]
}

\AddEq{d2=...c2 1}
{
$d_2=b_1b_2=b_1a_1a_2=b_1a_1b_1^{-1}c_2$.
}

\AddEq{d2=...c2}
{
\[h(a_1)\Act c_2=(b_1a_1b_1^{-1})c_2\]
}

\AddEq{Diagram: two representations of group 1}
{
\[\xymatrix{
a_2\ar[rrr]^{h(a_1)=L(a_1) }\ar[d]^{f(b_1)=L(b_1)}
&&& b_2\ar[d]^{f(b_1)=L(b_1)}\\
c_2\ar[rrr]_{h(a_1)}&&&d_2
}\]
}

\AddEq{Diagram: two representations of group 2}
{
\[\xymatrix{
a_2\ar[rrr]^{h(a_1)=L(a_1)} &&& b_2\ar[d]^{f(b_1)=L(b_1)}\\
c_2\ar[rrr]_{h(a_1)}\ar[u]_{(f(b_1))^{-1}=L(b^{-1}_1)}&&&d_2
}\]
}

\AddEq{haa=aRa}
{
\[a_2\Act h(a_1)=a_2\circ R(a_1)=b_2\]
}

\AddEq{haa=Laa}
{
\[h(a_1)\Act a_2=L(a_1)\circ a_2=b_2\]
}

\AddEq{(fg)->fg}
{
\[
(f,g)\rightarrow f\Act g
\]
}

\AddEq{Act=circ}
{
\[f\Act g=f\circ g\]
}

\AddEq{(ab)->ab}
{
\[
(a,b)\rightarrow a*b
\]
}

\DefEq
{
\[\sum_{i=1}^{\infty}\|a^i\|<\infty\]
}
{sum |ai|<infty}

\DefEq
{
\[\sum_{i=1}^{\infty}a^i\]
}
{sum ai}

\AddEq{integral of map}
{
\symb{\int_X d\mu(x)f(x)}{Lebesgue integral}{}
}

\AddEquation{integral of simple map}
{
\ShowSymbol{Lebesgue integral}{}=
\sum_n \mu(F_n)f_n
}

\AddEquation{integral of measurable map}
{
\ShowSymbol{Lebesgue integral}{}=
\lim_{n\rightarrow\infty}\int_X d\mu(x)f_n(x)
}

\DefEq
{
\[\int_Xd\mu(x)f(x)\]
}
{int f X}

\DefEq
{
\[\int_Xd\mu(x)g(x)\]
}
{int g X}

\DefEq
{
\[\int_Xd\mu(x)(f(x)+g(x))\]
}
{int f+g X}

\DefEq
{
\int_Xd\mu(x)(f(x)+g(x))=
\int_Xd\mu(x)f(x)+
\int_Xd\mu(x)g(x)
}
{int f+g X =}

\DefEq
{
\[\int_Xd\mu(x)\|f(x)\|\]
}
{int |f| X}

\DefEq
{
\left\|\int_Xd\mu(x)f(x)\right\|\le\int_Xd\mu(x)\|f(x)\|
}
{|int f|<int|f|}

\DefEq
{
\|f(x)\|\le M
}
{|f(x)|<M}

\DefEq
{
\int_Xd\mu(x) \|f(x)\|\le M \mu(X)
}
{int |f|<M mu}

%% file: Stmt.Module.English.tex
\input{Stmt.Module.Eq}

\DefExample{Abelian Group is Z module}
{
From theorems
\RefTheorem{action of ring of rational integers in Abelian group},
\RefTheorem{Abelian group is free}
and the definition
\refDefinition{module over algebra}{\SideWS \VectorSetNS},
it follows that free Abelian group $G$ is module over ring of integers $Z$.
}

\DefText[2]{ix-xi-1 roots}
{
the polynomial
\DrawEq[1]{ix-xi-1}{#1}
have no root
and the polynomial
\DrawEq[k]{ix-xi-1}{#2}
has the set of roots
\DrawEq{x=c12+j}{#2}
}

\DefEq
{
\begin{ShadedDefinition}
\labelDefinition{module over ring}
Let commutative ring $D$ has unit $1$.
Effective representation
\DrawEq{D->*V}{def}
of ring $D$
in an Abelian group $V$
is called
\AddIndex{module over ring}{module over ring} $D$
or
\AddIndex{$D$\Hyph module}{D module}.
Effective representation
\eqRef{D->*V}{def}
of commutative ring $D$ in an Abelian group $V$
is called
\AddIndex{vector space over field}{vector space over field} $D$
or
\AddIndex{$D$\Hyph vector space}{D vector space}.
\end{ShadedDefinition}
}
{... definition: module over ring}

\DefTheorem{tensor product of D-algebras is D-algebra}
{
Let
\ShowEq{A1...n}
be $D$\Hyph algebras.
Tensor product
$\Tensor A$
of $D$\Hyph modules
\ShowEq{A1...n}
is $D$\Hyph algebra,
if we define product by the equality
\ShowEq{xA1n*xA1n->oxA1n=}
}

\DefTheorem{representation of algebra H2 in LH}
{
Let product in algebra
\AoxA H
be defined according to rule
\DrawEq[pq]{product in algebra AA}{RH}
A representation
\DrawEq[RH]{h:AoxA->L(A)}H
of $R$\Hyph algebra
\AoxA H
in Abelian group
\ShowEq{L(A->B)}RHH{}
defined by the equality
\DrawEq[RH]{representation AA in LA}{RH}
is left \BoxB{H}module.
If we put $g=E$
where
\ShowEq{product in algebra AA 3}RHE
is identity map,
then we indentify the linear map
\DrawEq[fHH{}]{f: A->B}-
of quaternion algebra and tensor
\DrawEq{linear map f, basis E}1
using the following equality
\ShowEq{value of linear map f, basis E}
}

\DefProof{representation of algebra H2 in LH}
{
The theorem follows from the statement
that the map
\ShowEq{Eox=x}
generates any linear map of quaternion algebra.
}

\DefDefinition{square root}
{
The root
\ShowEq{square root}
of the equation
\DrawEq{x2=a}{root}
in $D$\Hyph algebra $A$
is called
\AddIndex{square root}{square root}
of $A$\Hyph number $a$.
}

\DefDefinition{polynomial*polynomial}
{
Bilinear map
\ShowEq{polynomial*polynomial}
is defined by the equality
\ShowEq{polynomial*polynomial, 1}
}

\DefTheorem{(a*b)x=(ax)(bx)}
{
For any tensors
\ShowEq{a in Aoxn}a{n+1},
\ShowEq{a in Aoxn}b{m+1}{},
product of homogeneous polynomials
\ShowEq{polynomials a o x,b o x}
is defined by the equality
\ShowEq{(a*b)x=(ax)(bx)}
}

\DefTheorem{x2=a quaternion root}
{
Let $H$ be quaternion algebra and $a$ be $H$\Hyph number.
\StartLabelItem
\begin{enumerate}
\item
Since
\ShowEq{Re sqrt a ne 0}
then the equation
\DrawEq{x2=a}{}
has roots
\ShowEq{x=x12}
such that
\ShowEq{x2=-x1}
\labelItem{Re sqrt a ne 0}
\item
Since $a=0$, then the equation
\DrawEq{x2=a}{}
has root $x=0$ with multiplicity $2$.
\labelItem{a=0}
\item
Since conditions
\RefItem{Re sqrt a ne 0},
\RefItem{a=0}
are not true,
then the equation
\DrawEq{x2=a}{}
has infinitly many roots such that
\ShowEq{|x|=sqrt a}
\labelItem{Re sqrt a=0}
\end{enumerate}
}

\DefDefinition{polylinear map of algebras, property}
{
The polylinear map of $D$\Hyph algebra $A$
\ShowEq{f:A->B}f{A^n}A
satisfies to equalities
\DrawEq[f]{f(ai+bi)=fai+fbi}{1}
\DrawEq[f]{f(pai)=pfai}{1}
\ShowEq{polylinear map of algebras, 1}AD
Let us denote
\ShowEq{set polylinear maps An}DAA
set of $n$\hyph linear maps
of $D$\Hyph algebra $A$.
}

\DefRemark{monomial of power k}
{
In the theorem
\RefTheorem{monomial of power k},
I considered recursive representation of a monomial in associative $D$\Hyph algebra.
Since product is independent of the way
in which brackets are placed,
recursive representation of a monomial is not unique
in associative $D$\Hyph algebra.
For instance, I can use any of the following forms
\ShowEq{ax^2bxcx^3d}
to represent monomial $ax^2bxcx^3d$.
I chose the equality
\EqRef{monomial of power k, F=1}
as the most simple for algorithm of division of polynomials.
}

\DefTheorem{monomial of power k}
{
Let $p_k(x)$ be
\AddIndex{monomial of power}{monomial of power} $k$
over associative $D$\Hyph algebra $A$.
Then
\StartLabelItem
\begin{enumerate}
\item
Monomial of power $0$ has form
\ShowEq{p0(x)=a0}
\item
If $k>0$, then
\labelItem{monomial of power k}
\ShowEq{monomial of power k, F=1}
where $a_k\in A$.
\end{enumerate}
}

\DefProof{monomial of power k}
{
We prove the theorem by induction over power
$n$ of monomial.

Let $n=0$.
We get the statement
\RefItem{monomial of power 0}
since monomial $p_0(x)$ is constant.

Let $n=k$. Last factor of monomial $p_k(x)$ is either $a_k\in A$,
or has form $x^l$, $l\ge 1$.
In the later case we assume $a_k=1$.
Factor preceding $a_k$ has form $x^l$, $l\ge 1$.
We can represent this factor as $x^{l-1}x$.
Therefore, we proved the statement.
}

\DefTheorem{r=+q circ p}
{
Let
\ShowEq{p=p circ x}
be polynomial of power $1$ and $p_1$ be nonsingular tensor.
Let
\ShowEq{r power k}
be polynomial of power $k$.
Then
\ShowEq{r=+q circ p}
where
\ShowEq{qij i j}
is homogeneous polynomial of power $i$.
}

\DefLabeledDefinition{linear map of A vector space}{\SideNS}
{
\ShowText{algebra over ring ()}
Let $V$, $W$ be \SideWS vector spaces over $D$\Hyph algebra $A$.
Linear map
\ShowEq{f:A->B}fVW
of $D$\Hyph vector space $V$
into $D$\Hyph vector space $W$
is called
\AddIndex{linear map}{linear map}
of \SideWS $A$\Hyph vector space $V$
into \SideWS $A$\Hyph vector space $W$.
Let us denote
\ShowEq{set linear maps, VW module}
set of linear maps
of \SideWS $A$\Hyph vector space $V$
into \SideWS $A$\Hyph vector space $W$.
}

\DefDefinition{component of linear map, basis E}
{
Expression
\ShowEq{component of linear map, basis E}
in equality
\eqRef{linear map f, basis E}1
is called
\AddIndex{component of linear map}
{component of linear map} $f$.
}

\DefTheorem{representation of algebra A2 in LA}
{
Let $A$ be $D$\Hyph algebra.
Let product in $D$\Hyph module
\AoxA A
be defined according to rule
\DrawEq[pq]{product in algebra AA}{}
A representation
\DrawEq[DA]{h:AoxA->L(A)}{}
of $D$\Hyph algebra
\AoxA A
in module
\ShowEq{L(A->B)}DAA{}
defined by the equality
\DrawEq[DA]{representation AA in LA}{}
allows us to identify tensor
\ShowEq{d in AxoA}A
and linear map
\ShowEq{product in algebra AA 2}DA{\delta}{}
where
\ShowEq{product in algebra AA 3}DA{\delta}
is identity map.
Linear map
\ShowEq{a ox b}{\delta}
has form
\ShowEq{a ox b c=}
}

\AddEq[1]{theorem: matrices fUD I}
{
\begin{ShadedTheorem}
\labelTheorem{#1, matrices fUD I}
Maps of conjugation
have
standard representation
\ShowEq{#1. standard representation I}
\end{ShadedTheorem}
}

\AddEq[2]{proof: matrices fUD I}
{
According to the equality
\EqRef{#2x= #1},
the map of conjugation
\ShowEq{Map of Conjugation}{#2}
has the matrix
\ShowEq{#1. matrix I#2}
The equality
\DrawEq[{#2}{}]{#1.fUD->fUU I}{#1#2}
follows from equalities
\EqRef{#1.fUD->fUU},
\EqRef{#1. matrix I#2}.
}

\DefTheorem{unit of algebra and ring}
{
Let $D$ be commutative ring.
Let $A$ be associative $D$\Hyph algebra with unit \(1_A\).
\StartLabelItem
\begin{enumerate}
\item
The map
\ShowEq{f:a in A->b in B}fdD{d1_A}A
is isomorphism of the ring $D$ into the center of the algebra $A$
and allows to identify \(D\)\Hyph number \(d\)
and \(A\)\Hyph number \(d\,1_A\).
\labelItem{isomorphism d->d1}
\item
The ring \(D\) is subalgebra of algebra $A$.
\labelItem{D is subalgebra of A}
\item
\ShowEq{D in Z(A)}
\end{enumerate}
}

\DefProof{unit of algebra and ring}
{
Let
\ShowEq{d12 in D}
According to the statement
\RefItem{distributive law, module},
\ShowEq{d1+d2 in A}
According to the definition
\RefDefinition{algebra over ring},
\ShowEq{d1 d2 in A}
The statement
\RefItem{isomorphism d->d1}
follows from equalities
\EqRef{d1+d2 in A},
\EqRef{d1 d2 in A}.
The statement
\RefItem{D is subalgebra of A}
follows from the statement
\RefItem{isomorphism d->d1}.

According to the definition
\RefDefinition{algebra over ring},
\ShowEq{ad=da in A}
The statement
\RefItem{D in Z(A)}
follows from the equation
\EqRef{ad=da in A}.
}

\DefTheorem{multiplication in algebra is distributive over addition}
{
The multiplication in the algebra $A$ is distributive over addition
\ShowEq{(a+b)c=.}
\ShowEq{a(b+c)=.}
}

\DefDefinition{module over commutative ring}
{
Effective representation of commutative ring $D$
in an Abelian group $V$
\DrawEq{D->*V}{module}
is called
\AddIndex{module over ring}{module over ring} $D$
or
\AddIndex{$D$\Hyph module}{D module}.
$V$\Hyph number is called
\AddIndex{vector}{vector}.
\ePrints{309618526,CACAA.06.121}
\ifx\Semafor\ValueOn
If $D$ is a field, then the Abelian group $V$
is called
\AddIndex{vector space}{vector space} over field $D$
or
\AddIndex{$D$\Hyph vector space}{D vector space}.
\fi
}

\DefDefinition{linear map of algebra}
{
\ShowEq{def DModule}
\ShowEq{=DD}
Linear map
\ShowEq{f:A->B}fAA
of $D$\Hyph algebra $A$
satisfies to equalities
\DrawEq[fab]{f(a+b)=...}{D }
\DrawEq[{}fa]{f(da)=h... module}{module}
\ShowEq{D ab in A, d in D}DA
Let us denote
\ShowEq{set linear maps, module}DAA
set of linear maps
of $D$\Hyph algebra $A$.
}

\DefLabeledDefinitionNote{direct sum of D modules}{\SideWS \Base-\VectorsSetNS}
{
Coproduct in category of \SideWS $\Base$\Hyph \VectorsSet is called
\AddIndex{direct sum}{direct sum}.\,\footnotemark
We will use notation
\ShowEq{direct sum of modules}{\Module}{\ModuleA}
for direct sum of \SideWS $\Base$\Hyph \VectorsSet $V$ and $W$.
}{
See also the definition
of direct sum of modules in
\citeBib{Serge Lang},
page 128.
On the same page, Lang proves the existence of direct sum of modules.
}

\DefLabeledTheorem{direct sum of D modules}{\SideWS \Base-\VectorsSetNS}
{
Let
\ShowEq{set Bi}{\Module}
be set of \SideWS $\Base$\Hyph \VectorsSetNS.
Then the representation
\ShowEq{\SideWS D->*o+A}{\Base}{\Module}{\ANumber}{\VNumber}
of the \algebraa $\Base$
in direct sum of Abelian groups
\ShowEq{o+Ai}{\Module}
is direct sum of \SideWS $\Base$\Hyph modules
\ShowEq{o+Ai}{\Module}
}

\DefProof{direct sum of D modules}
{
Let
\ShowEq{set f:A->B}f{\Module_i}W
be set of linear maps into \SideWS $\Base$\Hyph module $W$.
We define the map
\ShowEq{f:A->B}f{\Module}{\ModuleA}
by the equality
\DrawEq{fxi,i=sum fxi}{\SideWS \Base-\VectorSetNS}
The sum in the right side of the equality
\eqRef{fxi,i=sum fxi}{\SideWS \Base-\VectorSetNS}
is finite, since all summands, except for a finite number, equal $0$.
From the equality
\ShowEq{fi(x+y)=}
and the equality
\eqRef{fxi,i=sum fxi}{\SideWS \Base-\VectorSetNS},
it follows that
\ShowEq{f(x+y)i=}
From the equality
\ShowEq{fi(dx)=}
and the equality
\eqRef{fxi,i=sum fxi}{\SideWS \Base-\VectorSetNS},
it follows that
\ShowEq{f(dx)i=}
Therefore, the map $f$ is linear map.
The equality
\ShowEq{flj=fj}
follows from equalities
\EqRef{lj(x)=0x0},
\eqRef{fxi,i=sum fxi}{\SideWS \Base-\VectorSetNS}.
Since the map $\lambda_i$ is injective,
then the map $f$ is unique.
Therefore, the theorem folows from definitions
\RefDefinition{coproduct in category},
\RefDefinition{direct sum of Abelian groups}.
}

\DefProof{cr-power, 1}
{
The equality
\ShowEq{cr-power, 1 1}
follows from equalities
\eqRef{cr power, 0}1,
\eqRef{cr-power, n}1.
The equality
\EqRef{cr-power, 1}
follows from the equality
\EqRef{cr-power, 1 1}.
}

\DefProof{rc-power, 1}
{
The equality
\ShowEq{rc-power, 1 1}
follows from equalities
\eqRef{rc power, 0}1,
\eqRef{rc-power, n}1.
The equality
\EqRef{rc-power, 1}
follows from the equality
\EqRef{rc-power, 1 1}.
}

\DefTheorem{a in ZAb=>b in ZAa}
{
Since
\ShowEq{c in ZAa}ab,
then
\ShowEq{c in ZAa}ba.
}

\DefProof{a in ZAb=>b in ZAa}
{
Let
\ShowEq{c in ZAa}ab.
The equality
\DrawEq{ab=ba}Z
follows from the statement
\eqRef{c in S=>bc=cb}{Z(A,b)}.
The theorem follows from the equality
\eqRef{ab=ba}Z.
}

\DefTheorem{a in ZAb, c in ZA=>a in ZAbc}
{
Since
\ShowEq{c in ZAa}ab{}
and
\ShowEq{c in ZA}c
then
\ShowEq{c in ZAa}a{bc}.
}

\DefProof{a in ZAb, c in ZA=>a in ZAbc}
{
The equality
\ShowEq{a(bc)=...=(bc)a}
follows from definitions
\RefDefinition{nucleus of algebra},
\RefDefinition{center of algebra}
and from the theorem
\RefTheorem{c in S=>bc=cb submodule of A}.
The theorem follows from the equality
\EqRef{a(bc)=...=(bc)a}
and from the theorem
\RefTheorem{c in S=>bc=cb submodule of A}.
}

\DefTheorem{c in Za=> p(c) in Za}
{
Let $A$ be non\Hyph commutative $D$\Hyph algebra.
For any $a\in A$, if
\ShowEq{c in ZAa}ca,
then
\DrawEq[{p(c)}a{}{}]{c in ZAa}1
for any polynomial
\DrawEq{p(c) p in D}1
}

\DefProof{c in Za=> p(c) in Za}
{
}

\DefTheorem{c in S=>bc=cb submodule of A}
{
Let $A$ be non\Hyph commutative $D$\Hyph algebra.
For any $b\in A$, there exists subalgebra
\ShowEq{center of A number}
of $D$\Hyph algebra $A$ such that
\DrawEq[{Z(A,b)}{}]{c in S=>bc=cb}{Z(A,b)}
$D$\Hyph algebra
\ShowEq{center Ac}{}
is called
\AddIndex{center of $A$\Hyph number}{center of A number}
$b$.
}

\DefProof{c in S=>bc=cb submodule of A}
{
The subset
\ShowEq{center Ac}{}
is not empty because
\ShowEq{c in ZAa}0b,
\ShowEq{c in ZAa}bb.
Let
\ShowEq{d in D c in S 12}
Then
\ShowEq{dcb=bdc}
And this implies
\ShowEq{dc in S}
Therefore, the set
\ShowEq{center Ac}{}
is $D$\Hyph module.

Let
\ShowEq{c in ZAb}
Then
\ShowEq{c1c2 in ZAb}
From the equality
\EqRef{c1c2 in ZAb},
it follows that
\ShowEq{c1 c2 in ZAb}
Therefore, $D$\Hyph module
\ShowEq{center Ac}{}
is $D$\Hyph algebra.
}

\DefTheorem{set of maps B->A is Abelian group}
{
Let $A$ be algebra over commutative ring $D$.
For any set $B$,
the set of maps $A^B$ is Abelian group with respect to operation
\ShowEq{(f+g)x=}x
}

\DefProof{set of maps B->A is Abelian group}
{
The theorem follows from definitions
\ShowRef{set of maps B->A is Abelian group}
}

\DefTheorem{set of endomorphisms - D module}
{
The set
\ShowEq{End DV}{}
of endomorphisms of $D$\Hyph module $V$
is $D$\Hyph module.
}

\DefProof{set of endomorphisms - D module}
{
The theorem follows from the definition
\RefDefinition{D module, endomorphism}
and the corollary
\RefCorollary{product of linear map over scalar, D module}.
}

\DefDefinition{D module, endomorphism}
{
Linear map
\DrawEq[fVV]{f: A->B}{}
of $D$\Hyph module $V$
is called
\AddIndex{endomorphism}{endomorphism}
of $D$\Hyph module $V$.
We use notation
\ShowEq{set of endomorphisms D module}
set of endomorphisms
of $D$\Hyph module $V$.
}

\DefLabeledDefinition{module over algebra}{\SideWS \VectorSetNS}
{
Let $\Base$ be \Algebra.
Effective \SideNS\SidePresentation representation
\DrawEq{\SideWS \Base->*V}{\SideWS \VectorSetNS}
of \ShortAlgebraWS $\Base$
in \SetRepresentation $V$
is called
\AddIndex{\SideWS \VectorSet}{\SideWS \VectorSetNS}
over \ShortAlgebraWS $\Base$.
We will also say that \SetRepresentation $V$ is
\AddIndex{\SideWS $\Base$\Hyph \VectorSetNS}{\SideWS A \VectorSetNS}.
$V$\Hyph number is called
\AddIndex{vector}{vector}.
Bilinear map
\DrawEq{\SideWS \Base*V->V}{\VectorSet}
generated by \SideNS\SidePresentation representation
\DrawEq{\Base*V->V to \SideWS product}{\VectorSetNS}
is called
\SideNS\SidePresentation product
of vector over scalar.
}

\DefLabeledDefinition{module over associative algebra}{\SideWS \VectorSetNS}
{
Let $\Base$ be \Algebra.
Let $V$ be $\CBase$\Hyph \VectorSetNS.
Let in $\CBase$\Hyph \VectorSet
\ShowEq{End DV}{}
the product of endomorphisms is defined as composition of maps.
Let there exist homomorphism
\DrawEq[{g_{34}}A{\End(V)}{}]{f: A->B}{}
of $\CBase$\Hyph algebra $A$
into $\CBase$\Hyph algebra
\ShowEq{End DV}.

Effective \SideNS\SidePresentation representation
\DrawEq{\SideWS \Base->*V}{\SideWS \VectorSetNS}
of \ShortAlgebraWS $\Base$
in \SetRepresentation $V$
is called
\AddIndex{\SideWS \VectorSet}{\SideWS \VectorSetNS}
over \ShortAlgebraWS $\Base$.
We will also say that \SetRepresentation $V$ is
\AddIndex{\SideWS $\Base$\Hyph \VectorSetNS}{\SideWS A \VectorSetNS}.
$V$\Hyph number is called
\AddIndex{vector}{vector}.
Bilinear map
\DrawEq{\SideWS \Base*V->V}{\VectorSet}
generated by \SideNS\SidePresentation representation
\DrawEq{\Base*V->V to \SideWS product}{\VectorSetNS}
is called
\SideNS\SidePresentation product
of vector over scalar.
}

\DefLabeledDefinition{module over non-associative algebra}{\SideWS \VectorSetNS}
{
Let $\Base$ be non\Hyph\Algebra.
Let $V$ be $\CBase$\Hyph \VectorSetNS.
Let in $\CBase$\Hyph \VectorSet
\ShowEq{End DV}{}
the product of endomorphisms is defined
such way that
Let there exist homomorphism
\DrawEq[{g_{34}}A{\End(V)}{}]{f: A->B}{}
of $\CBase$\Hyph algebra $A$
into $\CBase$\Hyph algebra
\ShowEq{End DV}.
We will use the symbol $*$ to denote the product
in $D$\Hyph algebra $A$
and $D$\Hyph algebra
\ShowEq{End DV}.

Effective \SideNS\SidePresentation representation
\DrawEq{\SideWS \Base*->*V}{\SideWS \VectorSetNS}
of \ShortAlgebraWS $\Base$
in \SetRepresentation $V$
is called
\AddIndex{\SideWS \VectorSet}{\SideWS \VectorSetNS}
over \ShortAlgebraWS $\Base$.
We will also say that \SetRepresentation $V$ is
\AddIndex{\SideWS $\Base$\Hyph \VectorSetNS}{\SideWS A \VectorSetNS}.
$V$\Hyph number is called
\AddIndex{vector}{vector}.
Bilinear map
\DrawEq{\SideWS \Base**V->V}{\VectorSetNS}
generated by \SideNS\SidePresentation representation
\DrawEq{\Base**V->V to \SideWS product}{\VectorSetNS}
is called
\SideNS\SidePresentation product
of vector over scalar.
}

\DefLabeledTheorem{module over algebra}{\SideNS}
{
The following diagram of representations describes \SideWS $\Base$\Hyph module $V$
\DrawEq[{}{}{}{}]{diagram of representations, \SideWS module}1
The diagram of representations
\eqRef{diagram of representations, \SideWS module}1
holds
\AddIndex{commutativity of representations}{commutativity of representations}
\BaseRings
in Abelian group $V$
\ShowEq{\SideWS module, a d v}
}

\DefProof{module over algebra}
{
The diagram of representations
\eqRef{diagram of representations, \SideWS module}1
follows from the definition
\ShowRef{\SideWS module over algebra}
and from the theorem
\def\Temp{}
\ifx\SideNS\Temp
\RefTheorem{action of ring of rational integers in Abelian group}.
\else
\RefTheorem{Free Algebra over Ring}.
\fi
Since \SideNS\HSide transformation $\ATransf(a)$
is endomorphism
of $\CBase$\Hyph module $V$,
we obtain the equality
\EqRef{\SideWS module, a d v}.
}

\DefLabeledTheorem{module over non-associative algebra}{\SideNS}
{
The following diagram of representations describes \SideWS $\Base$\Hyph module $V$
\DrawEq[{}{}{}{}]{diagram of representations, \SideWS nonA module}1
The diagram of representations
\eqRef{diagram of representations, \SideWS nonA module}1
holds
\AddIndex{commutativity of representations}{commutativity of representations}
\BaseRings
in Abelian group $V$
\ShowEq{\SideWS module, a* d v}
}

\DefProof{module over non-associative algebra}
{
The diagram of representations
\eqRef{diagram of representations, \SideWS nonA module}1
follows from the definition
\refDefinition{module over non-associative algebra}{\SideWS \VectorSetNS}
and from the theorem
\RefTheorem{Free Algebra over Ring}.
Since \SideNS\HSide transformation $\ATransf(a)$
is endomorphism
of $\CBase$\Hyph module $V$,
we obtain the equality
\EqRef{\SideWS module, a* d v}.
}

\DefProof{Free Algebra over Ring}
{
The structure of $D$\Hyph module $A$ is generated by effective representation
\ShowEq{f:A->*B}{g_{12}}DA
of ring $D$ in Abelian group $A$.

\begin{ShadedLemma}
\labelLemma{structure of D algebra is generated by product}
{\it
Let the structure of $D$\Hyph algebra $A$
defined in $D$\Hyph module $A$,
be generated by product
\DrawEq{product in D algebra}{}
\AddIndex{Left shift of $D$\Hyph module $A$}{left shift of module}
defined by equation
\ShowEq{l(v):w->vw}
generates the representation
\ShowEq{endomorphism of module from product, 1}
of $D$\Hyph module $A$
in $D$\Hyph module $A$
}
\end{ShadedLemma}

{\sc Proof.}
According to definitions
\RefDefinition{algebra over ring}
and
\RefDefinition{polylinear map of modules},
left shift of $D$\Hyph module $A$
is linear map.
According to the definition
\refDefinition{linear map of D module}1,
the map \(l(v)\)
is endomorphism of $D$\Hyph module $A$.
The equation
\ShowEq{l(v1+v2)w}
follows from the definition
\RefDefinition{polylinear map of modules}
and from the equation
\EqRef{l(v):w->vw}.
\ShowEq{def sum of linear maps}
\ShowEq{ref sum of linear maps}
the equation
\ShowEq{l(v1+v2)}
follows from equation
\EqRef{l(v1+v2)w}.
The equation
\ShowEq{l(dv)w}
follows from the definition
\RefDefinition{polylinear map of modules}
and from the equation
\EqRef{l(v):w->vw}.
\ShowEq{ref sum of linear maps}
the equation
\ShowEq{l(dv)}
follows from equation
\EqRef{l(dv)w}.
The lemma follows from equalities
\EqRef{l(v1+v2)},
\EqRef{l(dv)}.
\hfill\(\odot\)

\begin{ShadedLemma}
\labelLemma{representation of D module in D module determines the product}
{\it
The representation
\ShowEq{endomorphism of module from product, 1}
of $D$\Hyph module $A$ in $D$\Hyph module $A$
determines the product in
$D$\Hyph module $A$ according to rule
\ShowEq{endomorphism of module from product, 8}
}
\end{ShadedLemma}

{\sc Proof.}
Since map $g_{23}\circ v$ is endomorphism of $D$\Hyph module $A$, then
\ShowEq{endomorphism of module from product, 3}
Since the map $g_{23}$ is
linear map
\ShowEq{g23:A->L}
then,
\ShowEq{ref sum and product over scalar, linear map}
\ShowEq{endomorphism of module from product, 4}
\ShowEq{endomorphism of module from product, 7}
From equations
\EqRef{endomorphism of module from product, 3},
\EqRef{endomorphism of module from product, 4},
\EqRef{endomorphism of module from product, 7}
and the definition
\RefDefinition{polylinear map of modules},
it follows that the map $g_{23}$ is bilinear map.
Therefore, the map $g_{23}$ determines the product in
$D$\Hyph module $A$ according to rule
\ShowEq{endomorphism of module from product, 8}
\hfill\(\odot\)

The theorem follows from lemmas
\RefLemma{structure of D algebra is generated by product},
\RefLemma{representation of D module in D module determines the product}.
}

\DefLabeledTheorem{A module -> algebra is associative}{\SideNS}
{
Let $g$ be effective left\Hyph side representation of $D$\Hyph algebra $A$
in $D$\Hyph module $V$.
Then $D$\Hyph algebra $A$ is associative.
}

\DefProof{A module -> algebra is associative}
{
Let
\ShowEq{abc in A, v in v}
Since \SideNS\Hyph side representation
$g$ is \SideNS\Hyph side representation
of the multiplicative group
of $D$\Hyph algebra $A$,
we obtain the equality
\DrawEq{\SideWS module, associative law}0
The equality
\ShowEq{a(b(cv))=(a(bc))v \SideNS}
follows from the equality
\eqRef{\SideWS module, associative law}0.
Since
\ShowEq{cv \SideNS}
the equality
\ShowEq{a(b(cv))=((ab)c)v \SideNS}
follows from the equality
\eqRef{\SideWS module, associative law}0.
The equality
\ShowEq{(a(bc))v=((ab)c)v \SideNS}
follows from equalities
\EqRef{a(b(cv))=(a(bc))v \SideNS},
\EqRef{(a(bc))v=((ab)c)v \SideNS}.
Since $v$ is any vector of $A$\Hyph module $V$,
the equality
\ShowEq{a(bc)=(ab)c \SideNS}
follows from the equality
\EqRef{(a(bc))v=((ab)c)v \SideNS}.
Therefore, $D$\Hyph algebra $A$ is associative.
}

\DefLabeledTheoremNote{definition of A module}{\SideWS \VectorSetNS}
{
Let $V$ be \SideWS $\Base$\Hyph module.
For any vector $v\in V$,
vector generated by the diagram of representations
\eqRef{diagram of representations, \SideWS module}1
has the following form
\ShowEq{\SideWS module, (a+n)v=av+nv}
\StartLabelItem
\begin{enumerate}
\item
The set of maps
\ShowEq{\SideWS module, a+n:V->V}
generates\,\footnotemark
\algebraa $\Base_{(1)}$
where the sum is defined by the equality
\DrawEq{(a+n)+(b+m)=}{\SideWS module}
and the product is defined by the equality
\DrawEq{(a+n)(b+m)=}{\SideWS module}
The \algebraa $\Base_{(1)}$ is called
\AddIndex{unital extension}{unital extension}
of the \algebraa $\Base$.
\begin{table}[h]
\begin{tabular}{|l|l|l|}
\hline
If \algebraa $\Base$ has unit, then
\ShowEq{A1=A unital extension}
If \algebraa $\Base$ is ideal of $\CBase$, then
\ShowEq{A1=D unital extension}
Otherwise
\ShowEq{A1=A+D unital extension}
\end{tabular}
\end{table}
\item
The \algebraa $\Base$ is \SideWS ideal of \algebraa $\Base_{(1)}$.
\labelItem{Algebra is \SideWS ideal of algebra (1)}
\item
The set of transormations
\EqRef{\SideWS module, (a+n)v=av+nv}
is \SideNS\HSide representation of \algebraa $\Base_{(1)}$ in Abelian group $V$.
\labelItem{\SideWS representation of D1 in Abelian group}
\end{enumerate}
We use the notation
\ShowEq{set of vectors generated by vector \Base}
for the set of vectors generated by vector $v$.
}
{See the definition of unital extension also on the pages
\citeBib{McCrimmon: Jordan Algebras}\Hyph 52,
\citeBib{Zharinov: foundation of mathematical physics}\Hyph 64.}

\DefLabeledTheoremNote{definition of non-associative A module}{\SideWS \VectorSetNS}
{
Let $V$ be \SideWS $\Base$\Hyph module.
For any vector $v\in V$,
vector generated by the diagram of representations
\eqRef{diagram of representations, \SideWS nonA module}1
has the following form
\ShowEq{\SideWS module, (a+d)*v=a*v+dv}
\StartLabelItem
\begin{enumerate}
\item
The set of maps
\ShowEq{\SideWS module, a*+d:V->V}
generates\,\footnotemark
\algebraa $\Base_{(1)}$
where the sum is defined by the equality
\DrawEq{(a+n)+(b+m)=}{\SideWS nonA module}
and the product is defined by the equality
\DrawEq{(a+n)*(b+m)=}{\SideWS nonA module}
The \algebraa $\Base_{(1)}$ is called
\AddIndex{unital extension}{unital extension}
of the \algebraa $\Base$.
\begin{table}[h]
\begin{tabular}{|l|l|l|}
\hline
If \algebraa $\Base$ has unit, then
\ShowEq{A1=A unital extension}
If \algebraa $\Base$ is ideal of $\CBase$, then
\ShowEq{A1=D unital extension}
Otherwise
\ShowEq{A1=A+D unital extension}
\end{tabular}
\end{table}
\item
The \algebraa $\Base$ is \SideWS ideal of \algebraa $\Base_{(1)}$.
\labelItem{Algebra is \SideWS ideal of nonA algebra (1)}
\item
The set of transormations
\EqRef{\SideWS module, (a+d)*v=a*v+dv}
is \SideNS\HSide representation of \algebraa $\Base_{(1)}$ in Abelian group $V$.
\labelItem{\SideWS representation of D1 in Abelian group, nonA}
\end{enumerate}
We use the notation
\ShowEq{set of vectors generated by vector \Base}
for the set of vectors generated by vector $v$.
}
{See the definition of unital extension also on the pages
\citeBib{McCrimmon: Jordan Algebras}\Hyph 52,
\citeBib{Zharinov: foundation of mathematical physics}\Hyph 64.}

\DefLabeledTheorem{definition of A module, property}{\SideNS}
{
Let $V$ be \SideWS $A$\Hyph\VectorSetNS.
Following conditions hold for $V$\Hyph numbers:
\StartLabelItem
\begin{enumerate}
\item 
\AddIndex{commutative law}{commutative law}
\DrawEq{commutative law}{\SideWS vector space}
\item 
\AddIndex{associative law}{associative law}
\labelItem{associative law, \SideWS module}
\DrawEq{associative law, \SideWS module}1
\item 
\AddIndex{distributive law}{distributive law}
\labelItem{distributive law, \SideWS module}
\DrawEq{distributive law, \SideWS module, 1}1
\DrawEq{distributive law, \SideWS module, 2}1
\item
\AddIndex{unitarity law}{unitarity law}
\labelItem{unitarity law, \SideWS \Base-module}
\ShowEq{unitarity law, \SideWS \Base-module}
\end{enumerate}
for any
\ShowEq{p,q in D, v,w in V}
}

\DefLabeledTheorem{definition of non-associative A module, property}{\SideNS}
{
Let $V$ be \SideWS $A$\Hyph\VectorSetNS.
Following conditions hold for $V$\Hyph numbers:
\StartLabelItem
\begin{enumerate}
\item 
\AddIndex{commutative law}{commutative law}
\DrawEq{commutative law}{\SideWS nonA vector space}
\item 
\AddIndex{associative law}{associative law}
\labelItem{associative law, nonA \SideWS module}
\DrawEq{associative law*, \SideWS module}1
\item 
\AddIndex{distributive law}{distributive law}
\labelItem{distributive law, nonA \SideWS module}
\DrawEq{distributive law*, \SideWS module, 1}1
\DrawEq{distributive law*, \SideWS module, 2}1
\item
\AddIndex{unitarity law}{unitarity law}
\labelItem{unitarity law, nonA \SideWS \Base-module}
\ShowEq{unitarity law*, \SideWS \Base-module}
\end{enumerate}
for any
\ShowEq{m,n in D}
\ShowEq{p,q in D, v,w in V}
}

\AddEq{proof: definition of A module}
{%
{\sc Proof of theorems
\refTheorem{definition of A module}{\SideWS \VectorSetNS},
\refTheorem{definition of A module, property}{\SideNS}.}
Let $v\in V$.

\begin{ShadedLemma}
\labelLemma{\SideWS module, map a+n:V->V is endomorphism of Abelian group}
{\it
Let
\ShowEq{module, d in D}
\ShowEq{module, a in A}
The map
\EqRef{\SideWS module, a+n:V->V}
is endomorphism of Abelian group $V$.
}
\end{ShadedLemma}

{\sc Proof.}
Statements
\ShowEq{\SideWS module, dv in V}
\ShowEq{\SideWS module, av in V}
follow from the theorems
\RefTheorem[\RefRepresentation]{structure of subrepresentations},
\refTheorem{module over algebra}{\SideNS}.
Since $V$ is Abelian group, then
\ShowEq{\SideWS module, dv+av in V}
Therefore,
for any $\CBase$\Hyph number $\DArg$
and for any $\Base$\Hyph number $a$,
we defined the map
\EqRef{\SideWS module, a+n:V->V}.
Since transformation $\DTransf(\DArg)$
and \SideNS\HSide transformation $\ATransf(a)$
are endomorphisms of Abelian group $V$,
then the map
\EqRef{\SideWS module, a+n:V->V}
is endomorphism of Abelian group $V$.
\hphantom{aaaa}\hfill\(\odot\)

Let $\Base_{(1)}$ be the set of maps
\EqRef{\SideWS module, a+n:V->V}.
The equality
\eqRef{distributive law, \SideWS module, 1}1
follows from the lemma
\RefLemma{\SideWS module, map a+n:V->V is endomorphism of Abelian group}.

Let
\ShowEq{p,q in D1}
According to the statement
\RefItem{representation of D1 in Abelian group},
we define the sum of $\Base_{(1)}$\Hyph numbers $p$ and $q$ by the equality
\eqRef{distributive law, \SideWS module, 2}1.
The equality
\ShowEq{\SideWS module, (a+n)+(b+m)=, 1}
follows from the equality
\eqRef{distributive law, \SideWS module, 2}1.
Since representation
$\DTransf$ is homomorphism of the aditive group of ring $\CBase$,
we obtain the equality
\ShowEq{distributive law, \CBase, \SideWS module, 2}
Since \SideNS\HSide representation
$\ATransf$ is homomorphism of the aditive group of \algebraa $\Base$,
we obtain the equality
\ShowEq{distributive law, \Base, \SideWS module, 2}
Since $V$ is Abelian group, then the equality
\ShowEq{\SideWS module, (a+n)+(b+m)=}
follows from equalities
\EqRef{\SideWS module, (a+n)+(b+m)=, 1},
\EqRef{distributive law, \CBase, \SideWS module, 2},
\EqRef{distributive law, \Base, \SideWS module, 2}.
From the equality
\EqRef{\SideWS module, (a+n)+(b+m)=},
it follows that the definition
\eqRef{(a+n)+(b+m)=}{\SideWS module}
of sum on the set $\Base_{(1)}$ does not depend on vector $v$.

Equalities
\eqRef{associative law, \SideWS module}1,
\EqRef{unitarity law, \SideWS \Base-module}
follow from the statement
\RefItem{\SideWS representation of D1 in Abelian group}.
Let
\ShowEq{p,q in D1}
\def\Temp{}
\ifx\SideNS\Temp
Since representation $\DTransf$ is representation
of the multiplicative group of ring $\CBase$,
we obtain the equality
\DrawEq{associative law, \CBase, \SideWS module}1
Since representation $g_2$ is representation
of the multiplicative group of ring $D$,
we obtain the equality
\DrawEq{\SideWS module, associative law}1
Since the ring $D$
is Abelian group,
we obtain the equality
\DrawEq{associative law, \CBase\Base, \SideWS module}1
The equality
\ShowEq{module, (a+n)(b+m)=}
follows from equalities
\ShowEq{ref module, (a+n)(b+m)=}
The equality
\eqRef{(a+n)(b+m)=}{\SideWS module}
follows from the equality
\EqRef{module, (a+n)(b+m)=}.
\else%
Since the product in $D$\Hyph algebra $A$ can be non associative,
then, based on the theorem
\refTheorem{definition of A module, property}{\SideNS},
we consider product
of $\Base_{(1)}$\Hyph numbers $p$ and $q$
as bilinear map
\ShowEq{f:D1xD1->D1}
such that following equalities are true
\DrawEq{fab=ab}{\SideWS module}
\DrawEq{f1p=p}{\SideWS module}
The equality
\DrawEq{pq=fpq}{\SideWS module}
follows from equalities
\eqRef{fab=ab}{\SideWS module},
\eqRef{f1p=p}{\SideWS module}.
The equality
\eqRef{(a+n)(b+m)=}{\SideWS module}
follows from the equality
\eqRef{pq=fpq}{\SideWS module}.
\fi%

The statement
\RefItem{Algebra is \SideWS ideal of algebra (1)}
follows from the equality
\eqRef{(a+n)(b+m)=}{\SideWS module}.
\qed
}

\DefLabeledDefinition{submodule}{\SideWS module}
{
Subrepresentation of \SideWS $\Base$\Hyph module $V$ is called
\AddIndex{submodule}{submodule}
of \SideWS $\Base$\Hyph module $V$.
}

\AddEq{theorem: submodule}
{
\begin{ShadedTheorem}
\labelTheorem{submodule, \SideWS module}
Let
\ShowEq{vi V}{}
be set of vectors of \SideWS $\Base$\Hyph module $V$.
If vectors
\ShowEq{set vi cols}v
belongs submodule $V'$ of \SideWS $\Base$\Hyph module $V$,
then linear combination of vectors
\ShowEq{set vi cols}v
belongs submodule $V'$.
\end{ShadedTheorem}
}

\DefProof{submodule}
{
The theorem follows from
the theorem
\refTheorem{set of vectors generated by set of vectors}{\SideWS module}
and definitions
\refDefinition{linear combination of vectors}{\SideNS},
\refDefinition{submodule}{\SideWS module}.
}

\AddEq{ref in item cvk module}
{
theorems
\RefTheorem[\RefRepresentation]{structure of subrepresentations},
\refTheorem{definition of A \VectorSetNS}{\SideNS}
}

\AddEq{ref in item cvk vector space}
{
the theorem
\RefTheorem[\RefRepresentation]{structure of subrepresentations},
the definition
\refDefinition{module over associative algebra}{\SideWS \VectorSetNS}.
}

\DefLabeledTheorem{set of automorphisms}{\SideNS}
{
The set $GL(V)$
of automorphisms
of \SideWS $A$\Hyph vector space $V$
is group.
}

\DefLabeledTheorem{active transformations in module}{\SideNS-\Cols}
{
\ShowEq{Let V be vector space of and basis}
Any automorphism $\Vector f$ of \SideWS $A$\Hyph vector space $V$
has form
\DrawEq[{v'}{}v{}f{}]{v1=v2*a \SideNS-\Cols}{automorphism}
where $f$ is a \ProductType nonsingular matrix.
Matrices of automorphisms of \SideWS $A$\Hyph vector space $V$ of \ColN s form a group
$GL(V_*)$
isomorphic to group $GL(V)$.
Automorphisms of \SideWS $A$\Hyph vector space of \ColsWS form
a \OtherSideNS\Hyph side linear
effective representation
\DrawEq[{GL(V_*)}{V_*}{}]{A->*B}{\SideNS-\Cols}
of the group $GL(V_*)$
in \SideWS $A$\Hyph vector space $V_*$.
}

\DefLabeledTheoremNote{set of vectors generated by set of vectors}{\SideWS \VectorSetNS}
{
Let $V$ be \SideWS $\Base$\Hyph \VectorSetNS.
The set of vectors generated by the set of vectors
\ShowEq{v(i) V}{}
has form\,\footnotemark
\DrawEq[{\Base}{\BaseExt}]{w=sum vi, \SideWS module}{\Base}
}
{For a set $A$,
we denote by $|A|$ the cardinal number of the set $A$.
The notation $|A|<\infty$ means that the set $A$ is finite.}

\AddEq{ref module structure of subrepresentations}
{
and theorems
\RefTheorem[\RefRepresentation]{structure of subrepresentations},
\refTheorem{definition of A module}{\SideWS \VectorSetNS},
}

\AddEq{ref vector space structure of subrepresentations}
{
and the theorems
\RefTheorem[\RefRepresentation]{structure of subrepresentations},
}

\DefLabeledLemma{w=w12 in Xk-1}{\SideWS \VectorSetNS}
{\leftskip=25pt
Let
\ShowEq{w=w12 in Xk-1}.
Then $w\in J(v)$.
}

\DefLabeledLemma{w=aw in Xk-1}{\SideWS \VectorSetNS}
{\leftskip=25pt
Let
\ShowEq{w=aw in Xk-1,\SideNS}
Then $w\in J(v)$.
}

\DefProof{set of vectors generated by set of vectors}
{
We prove the theorem by induction based on the theorem
\RefTheorem[\RefRepresentation]{structure of subrepresentations}.

For any
\ShowEq{vk in Jv},
let
\ShowEq{c(i)=dik}
Then
\DrawEq[{\Base_{\BaseExt}}{}]{v(k)=sum v(i), \SideWS module}{\Base}
From equalities
\eqRef{w=sum vi, \SideWS module}{\Base},
\eqRef{v(k)=sum v(i), \SideWS module}{\Base},
it follows that
the theorem is true on the set
\ShowEq{X0=v in J}

Let
\ShowEq{Xk in Jv}{k-1}
Acording to the definition
\ShowRef{\SideWS module over algebra}
\ShowEq{ref \VectorSet structure of subrepresentations}
if $w\in X_k$, then or
\ShowEq{w=w12 in Xk-1},
either
\ShowEq{w=aw in Xk-1,\SideNS}

{\leftskip=25pt
\ShowLemma{w=w12 in Xk-1}

{\sc Proof.}
According to the equality
\eqRef{w=sum vi, \SideWS module}{\Base},
there exist $\Base_{\BaseExt}$\Hyph numbers
\ShowEq{g(i)12}w
such that
\DrawEq[w]{w()12= \SideNS}{\VectorSetNS}
where sets
\DrawEq[w]{|c(i)12 ne 0|}{\SideWS \VectorSetNS}
are finite.
Since $V$ is \SideWS $\Base$\Hyph \VectorSetNS,
then from equalities
\eqRef{distributive law, \SideWS module, 2}1,
\eqRef{w()12= \SideNS}{\VectorSetNS},
it follows that
\DrawEq[w]{w()1+w()2= \SideNS}{\VectorSetNS}
From the equality
\eqRef{|c(i)12 ne 0|}{\SideWS \VectorSetNS},
it follows that
the set
\ShowEq{|gi() 1+2 ne 0|}w
is finite.
\hfill\(\odot\)

\ShowLemma{w=aw in Xk-1}

{\sc Proof.}
According to the statement
\RefItem[\RefRepresentation]{ax in Xk+1},
for any $\Base_{\BaseExt}$\Hyph number $a$,
\ShowEq{aw in Xk \SideWS module}
According to the equality
\eqRef{w=sum vi, \SideWS module}{\Base},
there exist $\Base_{\BaseExt}$\Hyph numbers
\ShowEq{set au v(i)}w
such that
\ShowEq{w()= \SideWS module}
where
\DrawEq{|w(i) ne 0|}{\SideWS module}
From the equality
\EqRef{w()= \SideWS module}
it follows that
\ShowEq{aw()= \SideWS module}
From the statement
\eqRef{|w(i) ne 0|}{\SideWS module},
it follows that
the set
\ShowEq{\SideWS module, |aw(i) ne 0|}
is finite.
\hfill\(\odot\)

\,
}

From lemmas
\refLemma{w=w12 in Xk-1}{\SideWS \VectorSetNS},
\refLemma{w=aw in Xk-1}{\SideWS \VectorSetNS},
it follows that
\ShowEq{Xk in Jv}k
}

\DefLabeledTheorem{da in DA->da in D1A}{\SideWS module}
{
The map
\ShowEq{da in DA->da in D1A}
is homomorphism of \SideWS $\Base$\Hyph module $V$
into \SideWS $\Base_{(1)}$\Hyph module $V$.
}

\DefProof{da in DA->da in D1A}
{
The equality
\DrawEq[{\AArg}]{(d1+0)+(d2+0)}{\SideWS module}
follows from the equality
\eqRef{(a+n)+(b+m)=}{\SideWS module}.
The equality
\DrawEq[{\AArg}]{(d1+0)(d2+0)}{\SideWS module}
follows from the equality
\eqRef{(a+n)(b+m)=}{\SideWS module}.
Equalities
\eqRef{(d1+0)+(d2+0)}{\SideWS module},
\eqRef{(d1+0)(d2+0)}{\SideWS module}
imply that the map $I_{(1)}$
is homomorphism of \algebraa $\Base$
into \algebraa $\Base_{(1)}$.

The equality
\ShowEq{da in DA->da in D1A 1, \SideWS module}
follows from the equality
\ShowRef{da in DA->da in D1A}
The theorem follows from the equality
\EqRef{da in DA->da in D1A 1, \SideWS module}
and the definition
\RefDefinition[\RefRepresentation]{morphism of representations of universal algebra}.
}

\DefLabeledConvention{da in DA->da in D1A}{\SideNS}
{
According to theorems
\refTheorem{set of vectors generated by set of vectors}{\SideWS module},
\refTheorem{da in DA->da in D1A}{\SideWS module},
the set of words generating
\SideWS $\Base$\Hyph module $V$
is the same as the set of words generating
\SideWS $\Base_{(1)}$\Hyph module $V$.
Therefore, without loss of generality,
we will assume that \algebraa $\Base$
has unit.
}

\DefLabeledConvention{linear combination of vectors}{\SideNS}
{
If it is necessary to explicitly show
that we multiply vector $v(\gii)$
over $\Base_{\BaseExt}$\Hyph number $c(i)$ on the \SideNS,
then we will call the expression
\ShowEq{the \SideWS linear combination()}
\AddIndex{\SideWS linear combination}{\SideWS linear combination}.
We will use this convention to similar terms.
For instance, we will say that a vector
\ShowEq{w=c(i)v(i).\SideNS}
linearly depends
on vectors
\ShowEq{set vi *}v
on the \SideNS.
}

\DefLabeledDefinition{linear combination of vectors}{\SideNS}
{
Let $\Module$ be \SideNS\HSide \VectorSetNS.
Let
\ShowEq{v(i) V}{}
be set of vectors.
\ShowEq{\SideWS linear combination()}%
The expression
\ShowEq{A linear combination() =}
is called
\AddIndex{linear combination}{linear combination} of vectors
\ShowEq{v(i)}
A vector
\ShowEq{w=c(i)v(i).\SideNS}
is called
\AddIndex{linearly dependent}{linearly dependent}
on vectors
\ShowEq{v(i)}
}

\DefLabeledDefinition{generating set of module}{\SideNS}
{
$J(v)$
is called
\AddIndex{submodule generated by set}
{submodule generated by set} $v$,
and $v$ is a
\AddIndex{generating set}{generating set}
of submodule $J(v)$.
In particular, a
\AddIndex{generating set}{generating set}
of \SideWS $\Base$\Hyph module $V$
is a subset $X\subset V$ such that
\ShowEq{generating set of module}
}

\DefLabeledDefinition{generating set of vector space}{\SideNS}
{
$J(v)$
is called
vector subspace generated by set $v$,
and $v$ is a
\AddIndex{generating set}{generating set}
of vector subspace $J(v)$.
In particular, a
\AddIndex{generating set}{generating set}
of \SideWS $\Base$\Hyph vector space $V$
is a subset $X\subset V$ such that
\ShowEq{generating set of module}
}

\DefText{generating set of module}
{
The following definition follows from the theorems
\refTheorem{set of vectors generated by set of vectors}{\SideWS module},
\RefTheorem[\RefRepresentation]{structure of subrepresentations}
and from the definition
\RefDefinition[\RefRepresentation]{generating set of representation}.
}

\DefText{quasibasis of module}
{
The following definition follows from the theorems
\refTheorem{set of vectors generated by set of vectors}{\SideWS module},
\RefTheorem[\RefRepresentation]{structure of subrepresentations}
and from the definition
\RefDefinition[\RefRepresentation]{basis of representation}.
}

\DefLabeledDefinition{basis of module}{\SideNS}
{
If the set $X\subset V$ is generating set of \SideWS $\Base$\Hyph \VectorSet
$V$, then any set $Y$, $X\subset Y\subset V$
also is generating set of \SideWS $\Base$\Hyph \VectorSet $V$.
If there exists minimal set $X$ generating
the \SideWS $\Base$\Hyph \VectorSet $V$, then the set $X$ is called
\AddIndex{quasibasis}{quasibasis} of \SideWS $\Base$\Hyph \VectorSet $V$.
}

\DefLabeledTheorem{linearly depends on rest of vectors}{\SideWS module}
{
Let $\Base$ be \DivAlgebra.
Since the equation
\DrawEq[wv]{\SideWS wi vi=0}{1 \Cols}
implies existence of index
\ShowEq{i=j}
such that
\ShowEq{wj ne 0}w,
then the vector
\ARow vj{}
linearly depends on rest of vectors $v$.
}

\DefProof{linearly depends on rest of vectors}
{
The theorem follows from the equality
\ShowEq{\SideWS vj=sum vi}vw
and from the definition
\refDefinition{linear combination of vectors}{\SideNS}.
}

\DefText{0=0vi}
{
It is evident that for any set of vectors
\ARow vi{}
\ShowEq{\SideWS 0=0vi}
}

\DefLabeledDefinitionNote{linearly independent vectors}{\SideNS}
{
The set of vectors\,\footnotemark
\ShowEq{set vi *}v
of \SideWS $\Base$\Hyph \VectorSet $V$ is
\AddIndex{linearly independent}{linearly independent set}
if
\ShowEq{w(i)=0}
follows from the equation
\DrawEq{\SideWS w(i) v(i)=0}{}
Otherwise the set of vectors
\ShowEq{set vi *}v
is \AddIndex{linearly dependent}{linearly dependent set}.
}
{I follow to the definition in
\citeBib{Serge Lang}, page 130.}

\DefLabeledTheorem{quasibasis of module}{\SideNS}
{
The set of vectors
\ShowEq{basis, module}
is quasibasis of \SideWS $\Base$\Hyph module
$V$, if following statements are true.
\StartLabelItem
\begin{enumerate}
\item
\labelItem{vector is linear combination of set, \SideWS module}
Arbitrary vector $v\in V$
is linear combination of
vectors of the set $\Basis e$.
\item
\labelItem{cannot be represented as a linear combination, \SideWS module}
Vector $e(\gii)$
cannot be represented as a linear combination
of the remaining vectors of the set $\Basis e$.
\end{enumerate}
}

\DefProof{quasibasis of module}
{
According to the statement
\RefItem{vector is linear combination of set, \SideWS module},
the theorem
\refTheorem{set of vectors generated by set of vectors}{\SideWS module}
and the definition
\refDefinition{linear combination of vectors}{\SideNS},
the set $\Basis e$ generates \SideWS $\Base$\Hyph module $V$
(the definition
\refDefinition{generating set of module}{\SideNS}).
According to the statement
\RefItem{cannot be represented as a linear combination, \SideWS module},
the set $\Basis e$ is minimal set
generating \SideWS $\Base$\Hyph module $V$.
According to the definitions
\refDefinition{basis of  module}{\SideNS},
the set $\Basis e$ is a quasibasis of \SideWS $\Base$\Hyph module $V$.
}

\DefLabeledTheorem{basis over division algebra}{\SideNS}
{
Let $\Base$ be \Algebra.
The set of vectors
\ShowEq{basis, module}
is a
\AddIndex{basis of \SideWS $\Base$\Hyph vector space}{basis, vector space} $V$
if vectors
\ARow ei
are linearly independent and any vector $v\in V$
linearly depends on vectors
\ARow ei.
}

\DefLabeledTheorem{vector space is free module}{\SideNS}
{
The \SideWS $\Base$\Hyph vector space
is free $\Base$\Hyph module.
}

\DefProof{basis over division algebra}
{
Let the set of vectors
\ShowEq{set vi \Cols}e
be linear dependent. Then the equation
\DrawEq[we]{\SideWS wi vi=0}{}
implies existence of index $\gii=\gij$ such that
\ShowEq{wj ne 0}w.
According to the theorem
\refTheorem{linearly depends on rest of vectors}{\SideWS module},
the vector
\ARow ej{}
linearly depends on rest of vectors of the set $\Basis e$.
According to the definition
\refDefinition{basis of module}{\SideNS},
the set of vectors
\ShowEq{set vi \Cols}e
is not a basis for \SideWS $\Base$\Hyph vector space $V$.

Therefore, if the set of vectors
\ShowEq{set vi \Cols}e
is a basis, then these vectors
are linearly independent.
Since an arbitrary vector $v\in V$
is linear combination of vectors
\ShowEq{set vi \Cols}e
then the set of vectors $v$,
\ShowEq{set vi \Cols}e
is not linearly independent.
}

\DefLabeledConvention{unit of algebra in basis}{\Cols}
{
Let \eV be the basis of $D$\Hyph algebra of \ColsWS $A$.
If $D$\Hyph algebra $A$ has unit,
then we assume that \ARow e0 is the unit of $D$\Hyph algebra $A$.
}

\DefTheoremNote{algebra A2 representation in LA}
{
Let $A_1$ and $A_2$ be $D$\Hyph algebras.
Let product in algebra \AoxA A
be defined according to rule
\DrawEq[pq]{product in algebra AA}{A2}
A linear map
\ShowEq{representation A2 in LA}
defined by the equality
\ShowEq{representation A2 in LA, 1}
is representation\,\footnotemark
of algebra $\ATwo$
in module
\ShowEq{L(A->B)}D{A_1}{A_2}.
}
{See the definition of representation
of $\Omega$\Hyph algebra in the definition
\ShowEq{ref definition: representation of algebra}.}

\DefLabeledTheorem[4]{matrix generates D module homomorphism}{\Cols(#1#2)}
{
\ShowText{map be homomorphism of ring (#1)}
\ShowText{matrix of numbers}D{#3}fiIjJ
The map\refFootnote{homomorphism of D module}{\Cols(#1#2)}
\newline
\FrameEqRef[{\Vector f}V]{homomorphism D algebra #1#2}{Vector module, coordinates \Cols}
defined by the equality
\ShowText{define homomorphism D module by matrix(#1)}
is homomorphism
of $D_{#1}$\Hyph module of \ColsWS $V_{#2}$
into $D_{#3}$\Hyph module of \ColsWS $V_{#4}$.
The homomorphism
\eqRef{homomorphism D algebra #1#2}{Vector module, coordinates \Cols}
which has the given%
\ShowText{define homomorphism by given matrix()}%
is unique.
}

\DefProof[5]{matrix generates D module homomorphism}
{
The equality
\ShowEq{fo(v+w) \Cols(#1)}
\begin{sloppypar}
\noindent
follows from equalities
\ShowRef{homomorphism of D module}{#1}{#2}{#5}
From the equality
\EqRef{fo(v+w) \Cols(#1)},
it follows that the map $\Vector f$ is homomorphism
of Abelian group.
The equality
\end{sloppypar}
\ShowEq{fo(va) \Cols(#1)}
follows from the equality
\eqRef{left homomorphism, f av=}{f D module #1#2}.
From the equality
\EqRef{fo(va) \Cols(#1)},
and definitions
\RefDefinition{Morphism of Diagram of Representations},
\refDefinition{linear map of D module}{#1#2},
it follows that the map
\newline
\FrameEqRef[{\Vector f}V]{homomorphism D algebra #1#2}{Vector module, coordinates \Cols}
\newline
\begin{sloppypar}
\noindent
is homomorphism
of $D_{#1}$\Hyph module of \ColN s $V_{#2}$
into $D_{#3}$\Hyph module of \ColN s $V_{#4}$.
\end{sloppypar}

Let $f$ be
matrix of homomorphisms $\Vector f$, $\Vector g$
relative to bases
\ShowEq{Bases eVW}12{}.%
The equality
\ShowEq{fov=gov \Cols(#1)}
follows from the theorem
\refTheorem{linear map of D module}{#1#2}.
Therefore, $\Vector f=\Vector g$.
}

\DefLabeledTheorem[4]{matrix generates D algebra homomorphism}{\Cols(#1#2)}
{
\ShowText{map be homomorphism of ring (#1)}
\ShowText{matrix of numbers and C}D{#3}fiIjJ{#1}{#2}
The map\refFootnote{homomorphism of d algebra}{#1#2\Cols}
\newline
\FrameEqRef[{\Vector f}V]{homomorphism D algebra #1#2}{Vector module, coordinates \Cols}
defined by the equality
\ShowText{define homomorphism D module by matrix(#1)}
is homomorphism
of $D_{#1}$\Hyph algebra of \ColsWS $A_{#2}$
into $D_{#3}$\Hyph algebra of \ColsWS $A_{#4}$.
The homomorphism
\eqRef{homomorphism D algebra #1#2}{Vector module, coordinates \Cols}
which has the given%
\ShowText{define homomorphism by given matrix()}%
is unique.
}

\DefProof[5]{matrix generates D algebra homomorphism}
{
According to the theorem
\refTheorem{matrix generates D module homomorphism}{\Cols(#1#2)},
the map
\eqRef{homomorphism D algebra #1#2}{Vector module, coordinates \Cols}
is homomorphism
of $D_{#1}$\Hyph module $A_{#2}$
into $D_{#3}$\Hyph module $A_{#4}$
and the homomorphism
\eqRef{homomorphism D algebra #1#2}{Vector module, coordinates \Cols}
is unique.

Let
\ShowEq{a=ai ei}a1i,
\ShowEq{a=ai ei}b1i.
The equality
\DrawEq{algebra, homomorphism and product abe #1#2}{\Cols}
follows from the equality
\eqRef{algebra, homomorphism and product (#1#2)}{\Cols(\SideNS)}.
\ShowText{matrix generates D algebra homomorphism (#1)}
The equality
\DrawEq{algebra, homomorphism and product abe 3(#1)}{\Cols}
follows from the equality
\eqRef{algebra, homomorphism and product abe 2(#1)}{\Cols}
and the equality
\ShowRef{homomorphism of d algebra, coordinates 1}{#1}{#2}
The equality
\DrawEq{algebra, homomorphism and product abe 4(#1)}{\Cols}
follows from the equality
\ShowEq{algebra, homomorphism and product, 1}
The equality
\DrawEq{algebra, homomorphism and product abe 5(#1)}{\Cols}
follows from the equality
\eqRef{algebra, homomorphism and product abe 4(#1)}{\Cols}
and the equality
\ShowRef{homomorphism of d algebra, coordinates 1}{#1}{#2}
The equality
\ShowRef{homomorphism, f vw=}{#1}{#2}
follows from equalities
\ShowRef{algebra, homomorphism and product 1}{#1}{#2}
According to the theorem
\refTheorem{homomorphism from A1 to A2, D algebra}{#1#2},
the map
\eqRef{homomorphism D algebra #1#2}{Vector module, coordinates \Cols}
is homomorphism
of $D_{#1}$\Hyph algebra $A_{#2}$
into $D_{#3}$\Hyph algebra $A_{#4}$.
}

\DefText{matrix generates D algebra homomorphism (1)}
{
Since the map $h$ is homomorphism
of the ring $D_1$ into the ring $D_2$, then
\DrawEq{algebra, homomorphism and product abe 2}{\Cols}
The equality
\DrawEq{algebra, homomorphism and product abe 2(1)}{\Cols}
follows from the equality
\eqRef{algebra, homomorphism and product abe 2}{\Cols}
and the equality
\ShowRef{product in algebra}
}

\DefText{matrix generates D algebra homomorphism ()}
{
The equality
\DrawEq{algebra, homomorphism and product abe 2()}{\Cols}
follows from the equality
\ShowRef{product in algebra}
}

\DefLabeledFootnote[4]{homomorphism of D module}{\Cols(#1#2)}
{
In theorems
\refTheorem{linear map of D module, coordinates}{#1#2\Cols},
\refTheorem{matrix generates D module homomorphism}{\Cols(#1#2)},
we use the following convention.
\ShowEq{Let be basis of module}1{}iI{}D{#1}V{#2}
\ShowEq{Let be basis of module}2{}jJ{}D{#3}V{#4}
}

\DefLabeledFootnote[6]{homomorphism of A module}{\SideNS-\Cols(#1#2#3)}
{
In theorems
\refTheorem{homomorphism A module}{\SideNS-\Cols(#1#2#3)},
\refTheorem{matrix generates A module homomorphism}{\SideNS-\Cols(#1#2#3)},
we use the following convention.
\ShowEq{prolog homomorphism of vector space(#1#2#3)}{#1}{#2}{#3}{#4}{#5}{#6}
}

\DefLabeledDefinition[4]{homomorphism from A1 to A2, D algebra}{#1#2}
{
\ShowText{algebra over ring (#1#2)}
Let diagram of representations
\DrawEq[{#1}{1.}{#2}h]{diagram of representations of D algebra}{->1(#1#2)}
describe $D_1$\Hyph algebra $A_1$.
Let diagram of representations
\DrawEq[{#3}{2.}{#4}h]{diagram of representations of D algebra}{->2(#1#2)}
describe $D_2$\Hyph algebra $A_2$.
Morphism
\DrawEq[fA]{homomorphism D algebra #1#2}1
of diagram of representations
\eqRef{diagram of representations of D algebra}{->1(#1#2)}
into diagram of representations
\eqRef{diagram of representations of D algebra}{->2(#1#2)}
is called
\AddIndex{homomorphism}{homomorphism}
of $D_{#1}$\Hyph algebra $A_{#2}$
into $D_{#3}$\Hyph algebra $A_{#4}$.
Let us denote
\ShowEq{set homomorphisms, D algebra #1#2}
set of homomorphisms
of $D_{#1}$\Hyph algebra $A_{#2}$
into $D_{#3}$\Hyph algebra $A_{#4}$.
}

\DefLabeledTheorem[5]{homomorphism from A1 to A2, D algebra}{#1#2}
{
Homomorphism
\newline
\FrameEqRef[fA]{homomorphism D algebra #1#2}1
\newline
of $D_{#1}$\Hyph algebra $A_{#2}$
into $D_{#3}$\Hyph algebra $A_{#4}$
is linear map
\eqRef[fA]{homomorphism D algebra #1#2}1
of $D_{#1}$\Hyph module $A_{#2}$
into $D_{#3}$\Hyph module $A_{#4}$
such that
\DrawEq[fab]{homomorphism, f vw=}{(#1#2)g}
and satisfies following equalities
\ShowEq{define homomorphism of D algebra #1#2}{#1}{#2}{#5}
}

\DefProof[4]{homomorphism from A1 to A2, D algebra}
{
According to definitions
\ShowRef{Morphism of Diagram of Representations}
the map
\eqRef{homomorphism D algebra #1#2}1
is morphism of representation $h_{1.12}$
describing $D_{#1}$\Hyph module $A_{#2}$
into representation $h_{2.12}$
describing $D_{#3}$\Hyph module $A_{#4}$.
Therefore, the map $f$ is linear map of $D_{#1}$\Hyph module $A_{#2}$
into $D_{#3}$\Hyph module $A_{#4}$ and equalities
\ShowEq{ref homomorphism of D algebra(#1)}{#1}{#2}
follow from the theorem
\refTheorem{linear map of D module}{#1#2}.

According to definitions
\ShowRef{Morphism of Diagram of Representations}
the map $f$
is morphism of representation $h_{1.23}$
into representation $h_{2.23}$.
The equality
\DrawEq{morphism of representation, algebra, 23}{#1#2}
follows from the equality
\eqRef{morphism of representations of universal algebra, 2m}{representation}.
From equalities
\eqRef{diagram of representations of D algebra}{->1(#1#2)},
\eqRef{diagram of representations of D algebra}{->2(#1#2)},
\eqRef{morphism of representation, algebra, 23}{#1#2},
it follows that
\DrawEq{algebra, morphism of representation 23, 1}{#1#2}
The equality
\eqRef{homomorphism, f vw=}{(#1#2)g}
follows from equalities
\eqRef{product in D algebra}{definition},
\eqRef{algebra, morphism of representation 23, 1}{#1#2}.
}

\AddEq{definition: unital algebra}
{
\begin{ShadedDefinition}
\labelDefinition{unital algebra}
If the product in $D$\Hyph algebra $A$ has unit element,
then $D$\Hyph algebra $A$ is called
\AddIndex{unital algebra}{unital algebra}\,\footnotemark
\end{ShadedDefinition}
\footnotetext{\,
See the definition of unital $D$\Hyph algebra also on the pages
\citeBib{McCrimmon: Jordan Algebras}\Hyph 137.
}
}

\DefDefinition{division algebra}
{
\(D\)\Hyph algebra \(A\) is called
\AddIndex{division algebra}{division algebra},
if for any \(A\)\Hyph number \(a\ne 0\)
there exists \(A\)\Hyph number \(a^{-1}\).
}

\DefTheorem{division algebra}
{
Let $A$ be associative division $D$\Hyph algebra.
The statement
\ShowEq{ab=ac a ne 0}
implies $b=c$.
}

\DefProof{division algebra}
{
The equality
\ShowEq{ab=ac a ne 0 1}
follows from the statement
\EqRef{ab=ac a ne 0 1}.
}

\DefLabeledDefinition[4]{linear map of D module}{#1#2}
{
Morphism of representations
\DrawEq[fV]{homomorphism D algebra #1#2}{module}
of $D_{#1}$\Hyph module $V_{#2}$
into $D_{#3}$\Hyph module $V_{#4}$
is called homomorphism or
\AddIndex{linear map}{linear map}
of $D_{#1}$\Hyph module $V_{#2}$
into $D_{#3}$\Hyph module $V_{#4}$.
Let us denote
\ShowEq{set linear maps, module (#1#2)}DV
set of linear maps
of $D_{#1}$\Hyph module $V_{#2}$
into $D_{#3}$\Hyph module $V_{#4}$.
}

\DefTheorem[3]{module of skew symmetric polylinear maps}
{
The set
\ShowEq{module of skew symmetric polylinear maps}{#1}{#2}{#3}
of skew symmetric polylinear maps
is $D$\Hyph module.
}

\AddEq[3]{remark: module of skew symmetric polylinear maps}
{
Without loss of generality, we assume
\ShowEq{L(A1,A0,B)=}{#1}{#2}{#3}
}

\DefTheorem{norm in D module A->B}
{
Let $A$ be normed $D$\Hyph module with norm $|x|_A$.
Let $B$ be normed $D$\Hyph module with norm $|y|_B$.
The map
\ShowEq{f in BA->|f|}
defined by the equality
\ShowEq{norm of map}
\ShowEq{norm of map, algebra}
is the norm in $D$\Hyph module $B^A$
and the value
\ShowEq{show|f|}
is called
\AddIndex{norm of map $f$}{norm of map}.
}

\DefTheorem{set of A->B is D module}
{
Let $A$ be Banach $D$\Hyph module with norm $|x|_A$.
Let $B$ be Banach $D$\Hyph module with norm $|y|_B$.
\StartLabelItem
\begin{enumerate}
\item
The set
$B^A$
of maps
\ShowEq{f:A->B}fAB
is $D$\Hyph module.
\labelItem{set of A->B is D module}
\item
The map
\ShowEq{f in BA->|f|}
defined by the equality
\ShowEq{norm of map}
\ShowEq{norm of map, algebra}
is the norm in $D$\Hyph module $B^A$
and the value
\ShowEq{show|f|}
is called
\AddIndex{norm of map $f$}{norm of map}.
\labelItem{norm of map}
\end{enumerate}
}

\DefTheorem{symmetrization of polylinear map}
{
Let
\ShowEq{f in L(A->B)}D{A^n}B{}
be a polylinear map.
Then the map
\ShowEq{symmetrization of polylinear map}
\EqParm{<f>}{=z}
defined by the equality
\ShowEq{symmetrization of polylinear map =}
is symmetric polylinear map
and is called
\AddIndex{symmetrization of polylinear map}{symmetrization of polylinear map}
$f$.
}

\DefTheorem{exterior product = skew symmetric}
{
Exterior product satisfies the following equation
\ShowEq{exterior product = skew symmetric}
}

\DefProof{exterior product = skew symmetric}
{
The equality
\EqRef{exterior product = skew symmetric}
follows from equalities
\EqRef{alternation of polylinear map =},
\EqRef{exterior product =}.
}

\DefDefinition{exterior product polylinear map}
{
The skew symmetric polylinear map
\ShowEq{exterior product}
\ShowEq{exterior product =}
is called
\AddIndex{exterior product}{exterior product}.
}

\AddEq{remark: product of polylinear maps}
{
If $B_1=B_2$, then, in the right side of the equality
\EqRef{hpq(fg)=},
we consider product of $B_1$\Hyph numbers
\ShowEq{f()},
\ShowEq{g()}.
According to the theorem
\RefTheorem{Free Algebra over Ring},
this definition is compatible with the definition
\RefDefinition{product of polylinear maps}.
}

\AddEq{definition: product of polylinear maps}
{
\begin{ShadedDefinition}
\labelDefinition{product of polylinear maps}
Let $A$, $B_2$ be free algebras over commutative ring $D$.\,\footnotemark
Let
\ShowEq{h:B1->*B2}
be left\Hyph side representation of
free associative $D$\Hyph algebra $B_1$ in $D$\Hyph module $B_2$.
The map
\ShowEq{hpq:Lpq->Lp+q}
is defined by the equality
\ShowEq{hpq(fg)=}
where, in the right side of the equality
\EqRef{hpq(fg)=},
we consider left\Hyph side product
of $B_2$\Hyph number
\ShowEq{g()}{}
over $B_1$\Hyph number
\ShowEq{f()}.
\end{ShadedDefinition}
\footnotetext{\,
To define product of skew symmetric polylinear maps,
I follow definition in section
\citeBib{Cartan differential form}-1.4 of chapter 1,
pages 12 - 14.
}
}

\DefDefinition{symmetric polylinear map}
{
Let $A$, $B$ be algebras over commutative ring $D$.
A polylinear map
\ShowEq{f in L(A->B)}D{A^n}B{}
is called
\AddIndex{symmetric}{symmetric polylinear map},
if
\ShowEq{fa=fsa}
for any permutation $\sigma$ of the set
\ShowEq{set a1n}.
}

\DefTheoremNote{fx1n=0, xi=xi1}
{
A polylinear map
\ShowEq{f in L(A->B)}D{A^n}B{}
is skew symmetric iff
\ShowEq{fx1n=0}
as soon as
$a_i=a_{i+1}$
for at list one\,\footnotemark
\ShowEq{1<=i<n}
}{
In the book
\citeBib{Cartan differential form},
page 9,
Henri Cartan considered the theorem
\RefTheorem{fx1n=0, xi=xi1}
as definition of skew symmetric map.
}

\DefTheorem[3]{alternation of polylinear map}
{
Let
\ShowEq{f in L(A->B)}{#1}{#2^n}{#3}{}
be a polylinear map.
Then the map
\ShowEq{alternation of polylinear map}
\ShowEq{[f]}{}
defined by the equality
\ShowEq{alternation of polylinear map =}
is skew symmetric polylinear map
and is called
\AddIndex{alternation of polylinear map}{alternation of polylinear map}
$f$.
}

\DefDefinition[3]{skew symmetric polylinear map}
{
A polylinear map
\ShowEq{f in L(A->B)}{#1}{#2^n}{#3}{}
is called
\AddIndex{skew symmetric}{skew symmetric polylinear map},
if
\ShowEq{fa=sfsa}
for any permutation $\sigma$ of the set
\ShowEq{set a1n}.
}

\DefLabeledTheorem{linear dependence between vectors of basis}{\SideWS module}
{
Let $\Base$ be \Algebra.
Let $\Basis e$ be quasibasis of \SideWS $\Base$\Hyph \VectorSet $V$.
Let
\DrawEq[ce]{\SideWS wi vi=0}{we \Cols}
be linear dependence of vectors of the quasibasis $\Basis e$.
Then
\StartLabelItem
\begin{enumerate}
\item
$\Base_{\BaseExt}$\Hyph number
\ShowEq{ci i in I}
does not have inverse element
in \ShortAlgebraWS $\Base_{\BaseExt}$.
\item
The set $\Base'$ of matrices
\ShowEq{c=ci i in I}
generates \SideWS $\Base$\Hyph module $\Base'$.
\end{enumerate}
}

\DefProof{linear dependence between vectors of basis}
{
Let there exist matrix
\ShowEq{c=ci i in I}
such that the equality
\eqRef{\SideWS wi vi=0}{we \Cols}
is true and there exist index
\ShowEq{i=j}
such that
\ShowEq{wj ne 0}c.
If we assume that $\Base$\Hyph number
\ACol cj
has inverse one, then the equality
\ShowEq{\SideWS vj=sum vi}ec
follows from the equality
\eqRef{\SideWS wi vi=0}{we \Cols}.
Therefore, the vector
\ARow ej{}
is linear combination of other vectors of the set $\Basis e$
and the set $\Basis e$ is not quasibasis.
Therefore, our assumption is false,
and $\Base_{\BaseExt}$\Hyph number
\ACol cj
does not have inverse.

Let matrices
\ShowEq{b in D'}b,
\ShowEq{b in D'}c.
From equalities
\DrawEq[be]{\SideWS wi vi=0}{}
\DrawEq[ce]{\SideWS wi vi=0}{}
it follows that
\ShowEq{(b+c)*e, \SideWS module}
Therefore, the set $\Base'$ is Abelian group.

Let matrix
\ShowEq{b in D'}c{}
and $a\in\Base$.
From the equality
\DrawEq[ce]{\SideWS wi vi=0}{}
it follows that
\ShowEq{(ac)*e, \SideWS module}
Therefore, Abelian group $\Base'$ is \SideWS $\Base$\Hyph module.
}

\DefLabeledTheorem{coordinates of vector with linear dependence}{\SideWS module}
{
Let \SideWS $\Base$\Hyph module $V$
have the quasibasis $\Basis e$ such that in the equality
\DrawEq[ce]{\SideWS wi vi=0}{2 we \Cols}
there exists index
\ShowEq{i=j}
such that
\ShowEq{wj ne 0}c.
Then
\StartLabelItem
\begin{enumerate}
\item
The matrix
\ShowEq{c=ci i in I}
determines coordinates of vector $0\in V$ with respect to basis $\Basis e$.
\labelItem{coordinates of vector 0 with linear dependence, \SideWS module}
\item
Coordinates of vector $\Vector v$ with respect to basis $\Basis e$
are uniquely determined up to a choice of coordinates of vector $0\in V$.
\labelItem{coordinates of vector with linear dependence, \SideWS module}
\end{enumerate}
}

\DefProof{coordinates of vector with linear dependence}
{
The statement
\RefItem{coordinates of vector 0 with linear dependence, \SideWS module}
follows from the equality
\eqRef{\SideWS wi vi=0}{2 we \Cols}
and from the definition
\refDefinition{coordinates of vector}{\SideWS \VectorSetNS}.

Let vector $\Vector v$ have expansion
\DrawEq{vv=ve \SideWS module}{2}
with respect to quasibasis $\Basis e$.
The equality
\ShowEq{v=v+0, \SideWS module}
follows from equalities
\eqRef{\SideWS wi vi=0}{2 we \Cols},
\eqRef{vv=ve \SideWS module}{2}.
The statement
\RefItem{coordinates of vector with linear dependence, \SideWS module}
follows from equalities
\eqRef{vv=ve \SideWS module}{2},
\EqRef{v=v+0, \SideWS module}
and from the definition
\refDefinition{coordinates of vector}{\SideWS \VectorSetNS}.
}

\DefLabeledTheorem{quasibasis of module is basis}{\SideWS \VectorSetNS}
{
Let \SideWS $\Base$\Hyph module $V$
have the quasibasis $\Basis e$ such that the equality
\DrawEq[ce]{\SideWS wi vi=0}{3 ce \Cols}
implies that
\ShowEq{ci=0 1n}
Then quasibasis $\Basis e$ is
basis of \SideWS $\Base$\Hyph module $V$.
}

\DefProof{quasibasis of module is basis}
{
The theorem follows from the definition
\RefDefinition{basis of representation}
and the theorem
\refTheorem{coordinates of vector with linear dependence}{\SideWS module}.
}

\DefLabeledDefinitionNote{free module over ring}{\SideNS}
{
The \SideWS $\Base$\Hyph module $V$ is
\AddIndex{free \SideWS $\Base$\Hyph module}{free module},\,\footnotemark
if \SideWS $\Base$\Hyph module $V$ has basis.
}
{
I follow to the
definition in \citeBib{Serge Lang}, page 135.
}

\DefProof{coordinates of vector of free module}
{
The theorem follows from the theorem
\refTheorem{coordinates of vector with linear dependence}{\SideWS module}
and from definitions
\refDefinition{linearly independent vectors}{\SideNS},
\refDefinition{free module over ring}{\SideNS}.
}

\DefTheoremNote{standard representation of map A->A, associative algebra}
{
Let $A$ be finite dimensional associative $D$\Hyph algebra.
Let $\Basis e$ be basis of $D$\Hyph module $A$.
Let $\Basis F$
be the basis\,\footnotemark
of left \BoxB{A}module
\ShowEq{L(A->B)}DAA.
\StartLabelItem
\begin{enumerate}
\item
The linear map
\ShowEq{f in L(A->B)}DAA{}
has the following expansion
\labelItem{map f generated by basis F, associative algebra}
\DrawEq{map f generated by basis F}{associative algebra}
where
\ShowEq{fk= in AxA}
\item
The linear map $f$ has the standard representation
\labelItem{standard representation of map A->A, associative algebra}
\ShowEq{standard representation of map A->A, associative algebra}
\ShowEq{standard representation of map A->A, o, associative algebra}
\end{enumerate}
}{
If $D$\Hyph module $A$
is not free $D$\Hyph nodule,
then we may consider the set
\ShowEq{Ik 1n}
of linear independent linear maps. The theorem is true for any linear map
\ShowEq{f:A->B}fAA
generated by the set of linear maps $\Basis F$.
}

\DefProof{standard representation of map A->A, associative algebra}
{
Since $\Basis F$ is the basis of left \BoxB{A}module
\ShowEq{L(A->B)}DAA,
then according to the definition
\refDefinition{free module over ring}{\SideNS},
there exists expansion
\DrawEq[A]{expansion of linear map with respect to basis}A
of the linear map $f$ with respect to the basis $\Basis F$.
According to the definition
\ShowEq{ref map j, representation, tensor product}
\DrawEq{f=fkxfk}{associative algebra}
The equality
\eqRef{map f generated by basis F}{associative algebra}
follows from equalities
\eqRef{expansion of linear map with respect to basis}A,
\eqRef{f=fkxfk}{associative algebra}.
According to theorem
\RefTheorem{standard component of tensor, algebra},
the standard representation of the tensor $f^k$ has form
\ShowEq{standard representation of map A1 A2, 3, associative algebra}
The equation
\EqRef{standard representation of map A->A, associative algebra}
follows from equations
\eqRef{map f generated by basis F}{associative algebra},
\EqRef{standard representation of map A1 A2, 3, associative algebra}.
}

\DefTheoremNote{set FoG generates left module L(A;B), n<m}
{
Let $A$ be $D$\Hyph module,
\ShowEq{n=dim A}nA.
Let $B$ be associative $D$\Hyph algebra,
\ShowEq{n=dim A}mB.
Let $\Basis F$ be basis of left \BoxB{B}module
\ShowEq{L(A->B)}DBB.
Let
\ShowEq{gi n<m}
Let
\ShowEq{f:A->B}GAB
be linear map of maximal rank.
The set
\DrawEq{F o G}1
generates left \BoxB{B}module\,\footnotemark
\ShowEq{L(A->B)}DAB.
}{
I do not claim that this set is a basis,
because maps
\ShowEq{FijG}
can be linearly dependent.
}

\DefProof{set FoG generates left module L(A;B), n<m}
{
Let
\ShowEq{f:A->B}gAB
be a linear map.
Let $\Basis e_A$ be the basis of $D$\Hyph module $A$.
Let $\Basis e_B$ be the basis of $D$\Hyph module $B$.
According to the theorem
\ShowEq{ref linear map of D1 D2 module, coordinates}
the linear map $G$ has coordinates
\ShowEq{matrix amn}Gmn
with respect to bases $\Basis e_A$, $\Basis e_B$
and the linear map $g$ has coordinates
\ShowEq{matrix amn}gmn
with respect to bases $\Basis e_A$, $\Basis e_B$.
A row $G_{\gii}$ of the matrix $G$,
as well a row $g_{\gii}$ of the matrix $g$,
is coordinates of linear form
\ShowEq{A->B}AD.
Since the matrix $G$ has maximal rank,
then rows of the matrix $G$ generate $D$\Hyph module
\ShowEq{L(A->B)}DAD{}
and rows of the matrix $g$ are linear combination of rows of the matrix $G$
\ShowEq{g=CG}
Since we can consider the matrix $C$
as coordinates of linear map
\ShowEq{f:A->B}CBB
then the equality
\ShowEq{g=(cF)G}
follows from the equality
\EqRef{g=CG}
and from the equality
\ShowEq{C=cF}
Since
\ShowEq{Gkj in D}
then the equality
\ShowEq{g=c(FG)}
follows from the equality
\EqRef{g=(cF)G}.
Therefore, the map $g$ belongs to
linear span of the set of maps
\eqRef{F o G}1.
}

\DefTheorem{set FoG generates left module L(A;B), n>m}
{
Let
\ShowEq{n=dim A}nA,
\ShowEq{n=dim A}mB.
Let $\Basis F$ be basis of left \BoxB{B}module
\ShowEq{L(A->B)}DBB.
Let
\ShowEq{gi n>m}
Let
\ShowEq{f:A->B}GAB
be linear map of maximal rank.
The set
\DrawEq{F o G}{}
generates the set of maps
\ShowEq{set ker G in ker g}
}

\DefProof{set FoG generates left module L(A;B), n>m}
{
The proof of the theorem is similar to the proof of the theorem
\RefTheorem{set FoG generates left module L(A;B), n<m}.
However, since number of rows of the matrix $G$ less then dimention
of $D$\Hyph module $A$, then rows of the matrix $G$ do not generate $D$\Hyph module
\ShowEq{L(A->B)}DAD{}
and the map $G$ has non trivial kernel.
In particular, rows of the matrix $g$ linearly depend on rows of the matrix $G$
iff
\ShowEq{ker G in ker g}g.
}

\DefTheorem{set FoG generates left module L(A;B)}
{
Let $\Basis F$ be basis of left \BoxB{B}module
\ShowEq{L(A->B)}DBB.
Let
\ShowEq{f:A->B}GAB
be linear map of maximal rank.
The set
\DrawEq{F o G}{}
generates the set of maps
\ShowEq{set ker G in ker g}
}

\DefProof{set FoG generates left module L(A;B)}
{
It is easy to see that theorem
\RefTheorem{set FoG generates left module L(A;B), n<m}
is particular case of the theorem
\RefTheorem{set FoG generates left module L(A;B), n>m},
because, in the theorem
\RefTheorem{set FoG generates left module L(A;B), n>m},
$\ker G=\emptyset$.
}

\AddEq{remark: set FoG generates left module L(A;B), n>m}
{
From the theorem
\RefTheorem{set FoG generates left module L(A;B), n>m},
it follows that
choice of the map $G$ depends on the map $g$.
}

\DefTheorem{basis D module L(D->A)}
{
Let $\Basis e$ be the basis of $D$\Hyph module $A$.
Then the set of maps
\ShowEq{d->dei}i
is the basis of $D$\Hyph module
\ShowEq{L(A->B)}DDA.
}

\DefProof{basis D module L(D->A)}
{
The theorem follows from the equality
\ShowEq{f(t)=fit ei}
}

\DefDefinition{component of linear map}
{
Expression
\ShowEq{component of linear map}
in equality
\EqRef{fk= in AxA}
is called
\AddIndex{component of linear map}
{component of linear map} $f$.
Expression
\ShowEq{standard component of linear map}
in the equality
\EqRef{standard representation of map A->A, associative algebra}
is called
\AddIndex{standard component of linear map}
{standard component of linear map} $f$.
}

\DefTheorem{endomorphism of module from product}
{
The representation
\ShowEq{endomorphism of module from product}
of $D$\Hyph module $A$ in $D$\Hyph module $A$
is equivalent to structure of $D$\Hyph algebra $A$.
}

\DefLabeledConvention{sum av() convention}{\SideNS}
{
We will use summation convention
in which repeated index
in linear combination
implies summation with respect to repeated index.
In this case we assume that we know the set
of summation index and do not use summation symbol
\ShowEq{av=sum av()\SideNS}
If needed to clearly show a set of indices, I will do it.
}

\DefConvention{av=sum av convention}
{
We will use Einstein summation convention
in which repeated index (one above and one below)
implies summation with respect to repeated index.
In this case we assume that we know the set
of summation index and do not use summation symbol
\ShowEq{av=sum av convention}
If needed to clearly show a set of indices, I will do it.
}

\DefTheorem{effective representation of the ring}
{
Representation
\ShowEq{f:A->*B}fDA
of the ring $D$
in an Abelian group $A$ is
\AddIndex{effective}{effective representation}
iff
$a=0$ follows from equation $f(a)=0$.
}

\DefProof{effective representation of the ring}
{
Suppose $a$, $b\in R$
cause the same transformation. Then
\ShowEq{representation of ring, 1}
for any $m\in A$.
From equalities
\EqRef{sum of transformations of Abelian group, 1},
\EqRef{representation of ring, 1},
it follows that
\ShowEq{representation of ring, 2}
The equality
\ShowEq{f(a-b)=0}
follows from equalities
\EqRef{0ov=0},
\EqRef{representation of ring, 2}.
Therefore, the representation $f$ is effective
iff $a=b$.
}

\DefTheorem{Representation of ring f(0)=v0}
{
Representation
\ShowEq{f:A->*B}fDA
of the ring $D$
in an Abelian group $A$
satisfies the equality
\ShowEq{f(0)=v0}
where
\ShowEq{0:A->A}
such that
\ShowEq{0ov=0}
}

\DefProof{Representation of ring f(0)=v0}
{
The equality
\ShowEq{fa=fa+f0}
follows from the equality
\EqRef{sum of transformations of Abelian group, 1}.
The equality
\ShowEq{f0x=0}
follows from the equality
\EqRef{fa=fa+f0}.
The equality
\EqRef{0ov=0}
follows from the equalities
\EqRef{f(0)=v0},
\EqRef{f0x=0}.
}

\DefText{sum of transformations of Abelian group}
{
We define the sum of transformations $f$ and $g$ of an Abelian group
according to rule
\ShowEq{(f+g)(a)=}
Therefore, considering the representation
\ShowEq{f:A->*B}fDA
of the ring $D$ in the Abelian group $A$, we assume
\ShowEq{f(a)+f(b)=}
According to the definition
\RefDefinition{representation of algebra},
the map $f$ is homomorphism of the ring $D$.
Therefore
\ShowEq{f(a+b)=f(a)+f(b)}
The equalty
\ShowEq{sum of transformations of Abelian group, 1}
follows from equalities
\EqRef{f(a)+f(b)=},
\EqRef{f(a+b)=f(a)+f(b)}.
}

\AddEq{remark: we do not have definition of determinant for division algebra}
{
According to
\citeBib{q-alg-9705026},
\citeBib{math.QA-0208146},
we do not have an appropriate definition
of a determinant for a division algebra.\,\footnote{
Professor Kyrchei
uses double determinant
(see the definition in the section
\citeBib{1812.03397}\Hyph 2.2)
to solve system of linear equations in quaternion algebra
and to solve eigenvalues problem 
(see the section
\citeBib{1812.03397}\Hyph 2.5).
I confine myself by consideration of quasideterminant,
because I am interested in a wider set of algebras.
\ePrints{2020.06.01,2204.06320,2022.01.05}
\ifx\Semafor\ValueOn

\FrameCiteBib{1812.03397}
\fi
}
However, we can define a quasideterminant which finally gives a
similar picture.
In definition
\RefDefinition{RC-quasideterminant},
I follow the definition
\citeBib{math.QA-0208146}-\href{http://arxiv.org/PS_cache/math/pdf/0208/0208146.pdf\#Page=9}{1.2.2}.
\ePrints{2020.06.01,2204.06320,2022.01.05}
\ifx\Semafor\ValueOn

\FrameCiteBib{q-alg-9705026}
\FrameCiteBib{math.QA-0208146}
\fi
}

\DefConvention{basis of algebra as basis of module}
{
According to definitions
\refDefinition{free module over ring}{},
\RefDefinition{algebra over ring},
free $D$\Hyph algebra $A$ is
$D$\Hyph module $A$,
which has a basis \eV[][.]
However, generally speaking, \eV is not a
basis of $D$\Hyph algebra $A$,
because there is product in $D$\Hyph algebra $A$.
For instance, in quaternion algebra, any quaternion
\ShowEq{any quaternion}
can be represented as
\ShowEq{any quaternion i j ij}
However, for most problems
using a basis of $D$\Hyph module $A$ is easier than
using a basis of $D$\Hyph algebra $A$.
For instance, it is easier to determine coordinates of $H$\Hyph number with respect to the basis
$(1,i,j,k)$
of $R$\Hyph vector space $H$,
than to determine coordinates of $H$\Hyph number with respect to the basis
$(1,i,j)$
of $R$\Hyph algebra $H$.
Therefore the phrase "we consider basis of $D$\Hyph algebra $A$"
means that we consider $D$\Hyph algebra $A$
and basis of $D$\Hyph module $A$.
}

\DefConvention{we use separate color for index of element}
{
Let $A$ be free algebra
with finite or countable basis.
Considering expansion of element of algebra $A$ relative basis $\Basis e$
we use the same root letter to denote this element and its coordinates.
In expression $a^2$, it is not clear whether this is component
of expansion of element
$a$ relative basis, or this is operation $a^2=aa$.
To make text clearer we use separate color for index of element
of algebra. For instance,
\ShowEq{Expansion relative basis in algebra}
}

\DefTheorem{direct sum of n Abelian groups}
{
Direct sum of Abelian groups
\ShowEq{a1n}An{}
coincides with their Cartesian product
\ShowEq{A1o+An=A1xAn}A
}

\AddEq{remark: direct sum of n Abelian groups}
{
Let
\ShowEq{a1o+.n}An
be direct sum of Abelian groups
\ShowEq{a1n}An.
According to the proof of the theorem
\RefTheorem{direct sum of Abelian groups},
any $A$\Hyph number $a$ has form
$(a_1,...,a_n)$
where $a_i\in A_i$.
We also will use notation
\ShowEq{a=a1 o+ an}
}

\AddEq{theorem: direct sum of n D modules}
{
\begin{ShadedTheorem}
\labelTheorem{direct sum of n \SideWS \Base-\VectorsSetNS}
Direct sum of \SideWS $\Base$\Hyph \VectorsSet
\ShowEq{a1n}{\Module}n{}
coincides with their Cartesian product
\ShowEq{A1o+An=A1xAn}{\Module}
\end{ShadedTheorem}
}

\DefLabeledTheorem{foa=fi ai}{\SideWS \Base-\VectorsSetNS}
{
Let
\ShowEq{A 1n}{\Module}n{}
be $\Base$\Hyph \VectorsSet and
\ShowEq{A 1o+.n}{\Module}n
Let us represent $\VNumber$\Hyph number
\ShowEq{A 1o+.n}{\VNumber}n
as column vector
\DrawEq[{\VNumber}n]{a=(a1.n col)}{}
Let us represent a linear map
\ShowEq{f:A->B}f{\Module}{\ModuleA}
as row vector
\DrawEq[fn]{a=(a1.n row)}{}
\ShowEq{f:A->B}{\aD fi}{\aU {\Module}i}{\ModuleA}
Then we can represent value of the map $f$ in $\Module$\Hyph number $\VNumber$
as product of matrices
\DrawEq[{\VNumber}{}]{foa=fi ai}{\SideWS \Base-\VectorsSetNS}
}

\DefProof{foa=fi ai}
{
The theorem follows from the definition
\eqRef{fxi,i=sum fxi}{\SideWS \Base-\VectorSetNS}.
}

\DefLabeledTheorem{foa=fi a}{\SideWS \Base-\VectorsSetNS}
{
Let
\ShowEq{A 1n}{\ModuleA}m{}
be $\Base$\Hyph \VectorsSet and
\ShowEq{A 1o+.n}{\ModuleA}m
Let us represent $\ModuleA$\Hyph number
\ShowEq{A 1o+.n}{\VNumberA}m
as column vector
\DrawEq[{\VNumberA}m]{a=(a1.n col)}{}
Then the linear map
\ShowEq{f:A->B}f{\Module}{\ModuleA}
has representation as column vector of maps
\DrawEq[fm]{a=(a1.n col)}{}
such way that, if
\ShowEq{b=foa}{\VNumberA}{\VNumber}
then
\DrawEq[{\VNumberA}{\VNumber}]{foa=fi a}{}
}

\DefProof{foa=fi a}
{
The theorem follows from the theorem
\RefTheorem{map from product into product}.
}

\DefTheorem{standard component of tensor, module}
{
Tensor product $\Tensor A$ of free
finite dimensional modules
$A_1$, ..., $A_n$ over the commutative ring $D$ is free
finite dimensional module.

Let
\ShowEq{tensor product of algebras, basis i}
be the basis of module $A_i$ over ring $D$.
We can represent any tensor $a\in\Tensor A$ in the following form
\ShowEq{standard component of tensor}
\ShowEq{tensor canonical representation, algebra}
The expression
$\ShowSymbol{standard component of tensor}{}$
is defined uniquely and
is called
\AddIndex{standard component of tensor}{standard component of tensora}.
}

\DefLabeledTheorem{map of direct sum of modules}{\SideWS \Base-\VectorsSetNS}
{
Let
\ShowEq{A 1n}{\Module}n,
\ShowEq{A 1n}{\ModuleA}m{}
be $\Base$\Hyph \VectorsSet and
\ShowEq{A 1o+.n}{\Module}n
\ShowEq{A 1o+.n}{\ModuleA}m
Let us represent $\Module$\Hyph number
\ShowEq{A 1o+.n}{\VNumber}n
as column vector
\DrawEq[{\VNumber}n]{a=(a1.n col)}{v,\SideWS \Base-\VectorsSetNS}
Let us represent $\ModuleA$\Hyph number
\ShowEq{A 1o+.n}{\VNumberA}m
as column vector
\DrawEq[{\VNumberA}m]{a=(a1.n col)}{w,\SideWS \Base-\VectorsSetNS}
Then the linear map
\ShowEq{f:A->B}fVW
has representation as a matrix of maps
\DrawEq[fnm]{a=(a11.nm matrix)}{f,\SideWS \Base-\VectorsSetNS}
such way that, if
\ShowEq{b=foa}{\VNumberA}{\VNumber}
then
\DrawEq[f{\VNumberA}{\VNumber}nm]{b=f rco a}{\SideWS \Base-\VectorsSetNS}
The map
\ShowEq{fij:->}{\Module}{\ModuleA}
is a linear map and is called
\AddIndex{partial linear map}{partial linear map}.
}

\DefProof{map of direct sum of modules}
{
According to the theorem
\refTheorem{foa=fi a}{\SideWS \Base-\VectorsSetNS},
there exists the set of linear maps
\ShowEq{fi:->}
such that
\DrawEq[{\VNumberA}{\VNumber}]{foa=fi a}{\SideWS \Base-\VectorsSetNS}
According to the theorem
\refTheorem{foa=fi ai}{\SideWS \Base-\VectorsSetNS},
for every $\gii$,
there exists the set of linear maps
\ShowEq{fij:->}{\Module}{\ModuleA}
such that
\DrawEq{fioa=fij aj}{\SideWS \Base-\VectorsSetNS}
If we identify matrices
\ShowEq{fij=(fi)j}
then the equality
\eqRef{b=f rco a}{\SideWS \Base-\VectorsSetNS}
follows from equalities
\eqRef{foa=fi a}{\SideWS \Base-\VectorsSetNS},
\eqRef{fioa=fij aj}{\SideWS \Base-\VectorsSetNS}.
}

\DefDefinition{Matrix of tensors}
{
Let $a$ be a matrix and
\ShowEq{a in Aoxn}{\aUD aij}n.
The matrix $a$ is called matrix of tensors
\ShowEq{Aoxn}n.
}

\AddEq{remark: map of direct sum of modules}
{
Let
\ShowEq{aU A1n}Bm{}
be $D$\Hyph algebras.
Then we can represent linear map $\aUD fij$
using \BoxB{\aU Bi}number.
}

\DefText{Free Algebra over Ring}
{
Let $D$ be commutative ring and $A$ be Abelian group.
The diagram of representations
\DrawEq[{}{}{}g]{diagram of representations of D algebra}{}
generates the structure of $D$\Hyph algebra $A$.
}

\DefText[4]{linear map of A module, coordinates 111}
{
\ShowText{Let be module of}V{#2}{}D_{#1}{algebra}
\ShowText{Let be module of}V{#4}{}D_{#3}{algevra}
\ShowEq{Let be basis of module}1{}iI{}D{#1}V{#2}
\ShowEq{Let be basis of module}2{}jJ{}D{#3}V{#4}
}

\DefText[1]{linear map of A module, coordinates 11}
{
}

\DefText[1]{linear map of A module, coordinates 1}
{
}

\DefLabeledTheorem[7]{linear map of A module, coordinates}{\SideWS \Base\DF\Base\DT}
{
\ShowText{linear map of A module, coordinates #1#2#3}
\ShowEq{Let be basis of module}A{#2}kK{}D{#1}A{#2}
\ShowEq{Let be basis of module}A{#5}lL{}D{#4}A{#5}
\ShowEq{Let be basis of module}V{#3}iI{\SideWS}A{#2}V{#3}
\ShowEq{Let be basis of module}V{#6}jJ{\SideWS}A{#5}V{#6}
Then linear map
\newline
\FrameEqRef[{\Vector g}{\Vector f}{}]{homomorphism A module #1#2#3}{\SideNS-\Cols}
\newline
has presentation
\DrawEq[#7]{f:A1->A2, A module}{#1#2#3}
\ShowEq{linear map in \SideWS A module}
\DrawEq{linear map in A module}{\SideNS}
relative to selected bases. Here
\begin{itemize}
\ShowText{homomorphism of vector space, algebra(#2)}{#1}{#2}{#3}{#4}{#5}{6}
\ShowText{homomorphism of vector space, algebra 1}{#1}{#2}{#3}{#4}{#5}{6}
\[===\]
\end{itemize}
The map
\ShowEq{fij:-> \SideWS A}
is a linear map and is called
\AddIndex{partial linear map}{partial linear map}.
}

\AddEq [9]{coordinates of the linear map}
{
\item $#1$ is coordinate matrix of $#2_1$\Hyph number
$\Vector #1$
relative the basis
\ShowEq{basis e}{#2_1}{}
\DrawEq [#1{#2_1}]{va=ae1, #8module}{#1 \DFDT\SideWS module}
\def\Temp{D1 D2 }
\ifx\DFDT\Temp
\item
\ShowEq{h(a)=...}{#1}{#4}{#5}{#6}
is a matrix of $#7$\Hyph numbers.
\fi
\item $#9$ is coordinate matrix of vector
\DrawEq[#1#9#3]{vb=f(va)}{#1 \DFDT\SideWS \VectorSetNS}
relative the basis
\ShowEq{basis e}{#2_2}{}
\DrawEq [#9{#2_2}]{va=ae1, #8module}{#9 \DFDT\SideWS module}
\item $#3$ is coordinate matrix of set of vectors
\ShowEq{Vector f(e1) module}#3#2#5#6
relative the basis
\ShowEq{basis e}{#2_2}.
The matrix $#3$ is
called \AddIndex{matrix of linear map}{matrix of linear map}
$\Vector #3$ relative bases
\ShowEq{basis e}{#2_1}{}
and
\ShowEq{basis e}{#2_2}.
}

\DefProof{linear map of A module, coordinates}
{
The equality
\eqRef{r2:A1->A2, \DFDT module}{1 \SideWS g}
follows from the theorem
\RefTheorem{linear map of D1 D2 module, coordinates}.

\ShowEq{module as direct sum}1kK
\ShowEq{module as direct sum}2lL
The equality
\eqRef{linear map in A module}{\SideNS}
follows from the theorem
\RefTheorem{map of direct sum of modules}.
}

\AddEq[3]{module as direct sum}
{
$A_{#1}$\Hyph module $V_{#1}$ is direct sum
\ShowEq{V=o+Ae}#1#2#3
}

\DefTheorem{representation of algebra An in LAnA}
{
Consider $D$\Hyph algebra $A$.
A representation
\ShowEq{representation An in LAnA}DA
of algebra \AOn
in module \LAnA
defined by the equality
\DrawEq[DA]{representation An in LAnA, 1}{}
allows us to identify tensor
\ShowEq{product in algebra An 1}
and transposition $\sigma\in S^n$
with map
\DrawEq[DA]{product in algebra An 2}{DA}
where
\ShowEq{product in algebra AA 3}DA{\delta}
is identity map.
}

\DefDefinition{linear map from A1 to A2, algebra}
{
Let $A_1$ and
$A_2$ be algebras over commutative ring $D$.
The linear map
of the $D_{\DF}$\Hyph module $A_\VF$
into the $D_{\DT}$\Hyph module $A_\VT$
is called
\AddIndex{linear map}{linear map}
of $D_{\DF}$\Hyph algebra $A_\VF$ into $D_{\DT}$\Hyph algebra $A_\VT$.

Let us denote
\ShowEq{set linear maps, module (1)}DA
set of linear maps
of $D$\Hyph algebra
$A_\VF$
into $D$\Hyph algebra
$A_\VT$.
}

\DefText[4]{notation for linear map}
{
If the map
\eqRef{homomorphism D algebra #1#2}{module}
is linear map of $D_{#1}$\Hyph module $V_{#2}$ into $D_{#3}$\Hyph module $V_{#4}$,
then, according to the definition
\refDefinition[\RefRepresentation]{side representation of group}{left},
I use notation
\DrawEq{f circ a =}{}
for image of the map $f$.
}

\DefTheorem{linear map times constant, algebra}
{
Let map
\DrawEq[f{A_1}{A_2}{}]{f: A->B}{}
be linear map of $D$\Hyph algebra $A_1$ into $D$\Hyph algebra $A_2$.
Then maps
\ShowEq{linear map times constant, algebra}
defined by equalities
\ShowEq{linear map times constant, 0, algebra}
are linear.
}

\DefTheorem{linear map AA LAA}
{
\ePrints{1502.04063}
\ifx\Semafor\ValueOn
Consider $D$\Hyph algebras $A_1$ and $A_2$.
For given map
\ShowEq{f in L(A->B)}D{A_1}{A_2},
\else
For given map
\ShowEq{f in L(A->B)}D{A_1}{A_2}{}
of $D$\Hyph algebra $A_1$ into $D$\Hyph algebra $A_2$,
\fi
there exists linear map
\ShowEq{linear map AA LAA}
defined by the equality
\ShowEq{linear map AA LAA, 1}
\ShowEq{linear map AA LAA, 2}
}

\DefEq
{
\begin{ShadedTheorem}
\labelTheorem{conjugation transformation}
Let $A_1$ be free $D$\Hyph module.
Let $A_2$ be free associative $D$\Hyph algebra.
Let $\Basis F$ be the basis of left \BoxB{A_2}module
\ShowEq{L(A->B)}D{A_1}{A_2}.
For any map
\ShowEq{Ik in I}
there exists set of linear maps
\ShowEq{conjugation transformation}
\ShowEq{conjugation transformation:}
of $D$\Hyph module
$A_1\otimes A_1$
into $D$\Hyph module
$A_2\otimes A_2$
such that
\ShowEq{conjugation transformation =}
The map
\ShowEq{show conjugation transformation}
is called
\AddIndex{conjugation transformation}{conjugation transformation}.
\end{ShadedTheorem}
}
{theorem: conjugation transformation}

\DefProof{conjugation transformation}
{
According to the theorem
\RefTheorem{product of linear map, algebra},
for any tensor
$a\in A_1\otimes A_1$,
the map
\ShowEq{x->Ik a x}
is linear.
According to the statement
\RefItem{map f generated by basis F},
there exists expansion
\ShowEq{expansion Ik a x}
Let
\ShowEq{b=Ikl a}
The equality
\EqRef{conjugation transformation =}
follows from equalities
\EqRef{expansion Ik a x},
\EqRef{b=Ikl a}.
From equalities
\ShowEq{Ikl a1+a2}
\ShowEq{Ikl d a}
it follows that the map $I_k^l$ is linear map.
}

\DefEq
{
\begin{ShadedTheorem}
\labelTheorem{representation of composition of linear maps}
Let $A_1$ be free $D$\Hyph module.
Let $A_2$, $A_2$ be free associative $D$\Hyph algebras.
Let $\Basis F$ be the basis of left \BoxB{A_2}module
\ShowEq{L(A->B)}D{A_1}{A_2}.
Let $\Basis G$ be the basis of left \BoxB{A_3}module
\ShowEq{L(A->B)}D{A_2}{A_3}.
\StartLabelItem
\begin{enumerate}
\item
The set of maps
\labelItem{basis of composition of linear maps}
\DrawEq{JlIk}{linear map}
is the basis of left \BoxB{A_3}module
\ShowEq{L(A->B)}D{A_1\rightarrow A_2}{A_3}.
\item
Let
\ShowEq{expansion of f with respect to basis I}
be expansion of linear map
\DrawEq[f{A_1}{A_2}{}]{f: A->B}{}
with respect to the basis $\Basis I$.
Let
\ShowEq{expansion of g with respect to basis J}
be expansion of linear map
\ShowEq{f:A->B}g{A_2}{A_3}
with respect to the basis $\Basis J$.
Then linear map
\DrawEq{h=g o f}{123}
has expansion
\labelItem{hlk=glfk}
\ShowEq{expansion of h with respect to basis K}
with respect to the basis $\Basis K$ where
\ShowEq{hlk=glfk}
\end{enumerate}
\end{ShadedTheorem}
}
{theorem: representation of composition of linear maps}

\DefProof{representation of composition of linear maps}
{
The equality
\ShowEq{h(a)1}
follows from equalities
\EqRef{expansion of f with respect to basis I},
\EqRef{expansion of g with respect to basis J},
\eqRef{h=g o f}{123}.
The equality
\ShowEq{h(a)2}
follows from equalities
\eqRef{JlIk}{linear map},
\EqRef{h(a)1}
and from the theorem
\RefTheorem{conjugation transformation}.
From the equality
\EqRef{h(a)2}
it follows that set of maps $\Basis K$ generates
left \BoxB{A_3}module\newline
\ShowEq{L(A->B)}D{A_1\rightarrow A_2}{A_3}.
From the equality
\ShowEq{aK=aJI}
it follows that
\ShowEq{aJ=0}
and, therefore, $a^{lk}=0$.
Therefore, the set $\Basis K$ is the basis of
left \BoxB{A_3}module
\ShowEq{L(A->B)}D{A_1\rightarrow A_2}{A_3}.
}

\DefDefinition{similar A numbers}
{
$A$\Hyph numbers $a$ and $b$ are called similar,
if there exists $A$\Hyph number $c$ such that
\ShowEq{similar A numbers}
}

\DefEq
{
\begin{ShadedTheorem}
\labelTheorem{representation of composition of linear maps A->A}
Let $A$ be free associative $D$\Hyph algebra.
Let left \BoxB{A}module
\ShowEq{L(A->B)}DAA{}
is generated by the identity map $F_0=\delta$.
Let
\ShowEq{expansion of f A->A}
be expansion of linear map
\ShowEq{f:A->B}fAA
Let
\ShowEq{expansion of g A->A}
be expansion of linear map
\ShowEq{f:A->B}gAA
Then linear map
\DrawEq{h=g o f}{A}
has expansion
\ShowEq{expansion of h A->A}
where
\ShowEq{hlk=glfk 01}
\end{ShadedTheorem}
}
{theorem: representation of composition of linear maps A->A}

\DefEq
{
\begin{proof}
The equality
\ShowEq{h(a)}
follows from equalities
\EqRef{expansion of f A->A},
\EqRef{expansion of g A->A},
\eqRef{h=g o f}{A}.
The equality
\EqRef{hlk=glfk 01}
follows from the equality
\EqRef{h(a)}.
\end{proof}
}
{proof: representation of composition of linear maps A->A}

\DefTheorem{h generated by f, associative algebra}
{
Consider $D$\Hyph algebra $A_1$
and associative $D$\Hyph algebra $A_2$.
Consider the representation of algebra $\ATwo$
in the module
\ShowEq{L(A->B)}D{A_1}{A_2}.
The map
\ShowEq{h:A1->A2}
generated by the map
\DrawEq[f{A_1}{A_2}{}]{f: A->B}{}
has form
\ShowEq{h generated by f, associative algebra}
}

\DefTheorem{coordinates of map A->A, algebra}
{
Let $A$ be free finite dimensional associative $D$\Hyph algebra.
Let $\Basis e$ be basis of $D$\Hyph module $A$.
Let
\ShowEq{structural constants, algebra}
be structural constants of algebra $A$.
Let $\Basis F$ be the basis
of left \BoxB{A}module
\ShowEq{L(A->B)}DAA{}
and
\ShowEq{coordinates of map Ik}
be coordinates of map $F_k$ with respect to basis $\Basis e$.
Coordinates
\ShowEq{Coordinates of map f}
of the map
\ShowEq{f in L(A->B)}DAA{}
and its standard components
\ShowEq{standard components of map f}k
are connected by the equation
\DrawEq{coordinates of map A->A, 2}{associative algebra}
}

\DefProof{coordinates of map A->A, algebra}
{
Relative to basis
$\Basis e$, linear maps $f$ and $I_k$ have form
\ShowEq{coordinates of map f, associative algebra}
\ShowEq{coordinates of map Fk, associative algebra}
The equality
\ShowEq{coordinates of map A->A, 3, associative algebra}
follows from equalities
\EqRef{standard representation of map A->A, associative algebra},
\EqRef{coordinates of map f, associative algebra},
\EqRef{coordinates of map Fk, associative algebra}.
Since vectors $\aD ek$
are linear independent and $x^{\gi i}$ are arbitrary,
then the equation
\eqRef{coordinates of map A->A, 2}{associative algebra}
follows from the equation
\EqRef{coordinates of map A->A, 3, associative algebra}.
}

\DefTheorem{coordinates of map A1 A2 FoG, algebra}
{
Let $\Basis e_1$ be basis of the finite dimensional
$D$\Hyph module $A_1$.
Let $\Basis e_2$ be basis of the finite dimensional associative
$D$\Hyph algebra $A_2$.
Let
\ShowEq{f in L(A->B)}D{A_1}{A_2}.
Let
\ShowEq{structural constants, algebra}
be structural constants of algebra $A_2$.
Let $\Basis F$ be the basis
of left \BoxB{A_2}module
\ShowEq{L(A->B)}D{A_2}{A_2}{}
and
\ShowEq{coordinates of map Ik}
be coordinates of map $F_k$ with respect to basis $\Basis e_2$.
Let
\ShowEq{f:A->B}GAB
be linear map of maximal rank such that
\ShowEq{ker G in ker g}f{}
and
\ShowEq{coordinates of map G}
be coordinates of map $G$ with respect to bases $\Basis e_1$ and $\Basis e_2$.
Coordinates
\ShowEq{Coordinates of map f}
of the map $f$
and its standard components
\ShowEq{standard components of map f}k
are connected by the equation
\DrawEq[f-]{coordinates of map A1 A2, FoG, associative algebra}f
}

\AddEq[1]{theorem: maps of conjugation as basis}
{
\begin{ShadedTheorem}
\labelTheorem{maps of conjugation as basis, algebra #1}
The set of linear maps
\ShowEq{vI=EI}
is the basis
of left $#1$\Hyph vector space
\ShowEq{L(A->B)}R{#1}{#1}
and a linear map
\ShowEq{f in L(A->B)}R{#1}{#1}{}
have expansion
\ShowEq{linear map of algebra #1, structure, 1}
\ShowEq{linear map of algebra #1, structure, 2}
where
$#1$\Hyph numbers
\ShowEq{a... #1}
are defined by the equality
\ShowEq{a...= #1}
\end{ShadedTheorem}
}

\AddEq[1]{theorem: linear map, maps of conjugation, algebra}
{
\begin{ShadedTheorem}
\labelTheorem{linear map, maps of conjugation, algebra #1}
A linear map
\ShowEq{f in L(A->B)}R{#1}{#1}{}
have expansion
\ShowEq{linear map of algebra #1, structure, 1}
\ShowEq{linear map of algebra #1, structure, 2}
where
\ShowEq{vI=EI}
and $#1$\Hyph numbers
\ShowEq{a... #1}
are defined by the equality
\ShowEq{a...= #1}
\end{ShadedTheorem}
}

\AddEq [1]{proof: L is left vector space}
{
\begin{proof}
According to the theorem
\RefTheorem{linear map, maps of conjugation, algebra #1},
the expansion
\EqRef{linear map of algebra #1, structure, 1}
of the linear map $f$
exists and is unique.
Therefore, the set
\ShowEq{I=(E,I)}{#1}{}
is basis of left $#1$\Hyph vector space
\ShowEq{L(A->B)}R{#1}{#1}.
\end{proof}
}

\DefTheorem{f=.E+.I+.J+.K}
{
Let
\ShowEq{f:H->H ij}
be linear map of quaternion algebra.
Let
\DrawEq{fi=fij ej}H
be quaternions represented by rows of the matrix $f$
\ShowEq{fi=fij ej H}
Then
\ShowEq{f=.E+.I+.J+.K}
}

\DefProof{f=.E+.I+.J+.K}
{
Equalities
\ShowEq{a0=f03}
\ShowEq{a1=f03}
\ShowEq{a2=f03}
\ShowEq{a3=f03}
follow from the equality
\EqRef{a...= H}.
The equalitiy
\EqRef{f=.E+.I+.J+.K}
follows from equalities
\ShowEq{f=.E+.I+.J+.K ref}
}

\DefTheorem{f=.I0.7}
{
Let
\ShowEq{f:H->H ij}
be linear map of quaternion algebra.
Let
\DrawEq{fi=fij ej}O
Then
\ShowEq{f=.I0.7}
}

\DefProof{f=.I0.7}
{
Equalities
\ShowEq{a0=f0.7}
follow from the equality
\EqRef{a...= H}.
The equalitiy
\EqRef{f=.I0.7}
follows from equalities
\ShowEq{f=.I0.7 ref}
}

\DefDefinition{maps of conjugation, complex field}
{
\ShowEq{C maps of conjugation}
Complex field has following
\AddIndex{maps of conjugation}{map of conjugation}
\ShowEq{C list maps of conjugation}
}

\DefTheorem{HE is algebra isomorphic to quaternion algebra}
{
The set
\ShowEq{HE set}
is $R$\Hyph algebra isomorphic to quaternion algebra.
}

\DefProof{HE is algebra isomorphic to quaternion algebra}
{
The theorem follows from equalities
\ShowEq{aE+bE=(a+b)E}
\ShowEq{aE o bE=(ab)E}
based on the theorem
\RefTheorem{linear map, maps of conjugation, algebra H}.
}

\DefDefinition{algebra, left and right action on L}
{
Let $A$ be $D$\Hyph module
and $B$ be $D$\Hyph algebra.
For any map
\ShowEq{f in L(A->B)}DAB{}
and $b\in B$,
we define left side transformation of the map $f$ using equality
\ShowEq{b o f =}
and right side transformation of the map $f$ using equality
\ShowEq{b * f =}
}

\DefDefinition{projection maps, complex field}
{
Following projection maps are defined in complex field
\ShowEq{projection maps, complex field}
}

\DefTheorem{Expansion of projection maps relative E,I}
{
Expansion of projection maps relative basis
\ShowEq{e=(E,I)}{}
has form
\ShowEq{projection mappings 1, complex field}
}

\DefProof{Expansion of projection maps relative E,I}
{
The theorem follows from the theorem
\RefTheorem{linear map, maps of conjugation, algebra C}
and from the definition
\RefDefinition{projection maps, complex field}.
}

\DefTheorem{f=.E+.I}
{
Let
\DrawEq[fCC{}]{f: A->B}{}
be linear map of complex field.
Let
\DrawEq{fi=fij ej}C
Then
\ShowEq{f=.E+.I}
}

\DefProof{f=.E+.I}
{
Equalities
\ShowEq{a0=f01}
\ShowEq{a1=f01}
follow from the equality
\EqRef{a...= C}.
The equalitiy
\EqRef{f=.E+.I}
follows from equalities
\ShowEq{f=.E+.I ref}
}

\DefDefinition{polylinear map of algebras}
{
Let $A_1$, ..., $A_n$, $S$ be $D$\Hyph algebras.
Polylinear map
\ShowEq{f:A->B}f{\Times}S
of $D$\Hyph modules
$A_1$, ..., $A_n$
into $D$\Hyph module $S$
is called
\AddIndex{polylinear map}{polylinear map} of $D$\Hyph algebras
$A_1$, ..., $A_n$
into $D$\Hyph algebra $S$.
Let us denote
\ShowEq{set polylinear maps}
set of polylinear maps
of $D$\Hyph modules
$A_1$, ..., $A_n$
into $D$\Hyph module
$S$.
Let us denote
\ShowEq{set polylinear maps An}DAS
set of $n$\hyph linear maps
of $D$\Hyph module $A$ ($A_1=...=A_n=A$)
into $D$\Hyph module
$S$.
}

\DefTheorem{sum of linear maps, D module}
{
Let $A_1$, $A_2$ be $D$\Hyph modules.
The map
\ShowEq{sum of maps,,D module}
\ShowEq{sum of maps, 1, D module}
defined by equation
\ShowEq{sum of maps, D module}
is called
\AddIndex{sum of maps}{sum of maps}
$f$ and $g$
and is linear map.
The set
$\mathcal L(D;A_1;A_2)$
is an Abelian group
relative sum of maps.
}

\DefTheorem{sum of polylinear maps, module}
{
Let $D$ be the commutative ring.
Let $A_1$, ..., $A_n$, $S$ be $D$\Hyph modules.
The map
\ShowEq{sum of maps,,polylinear}
\ShowEq{sum of maps, 1, polylinear}
defined by the equality
\ShowEq{sum of  maps, polylinear}
is called
\AddIndex{sum of polylinear maps}{sum of maps}
$f$ and $g$
and is polylinear map.
The set
\ShowEq{module of polylinear maps}
is an Abelian group
relative sum of maps.
}

\DefCorollary{sum of linear maps, D module}
{
Let $A_1$, $A_2$ be $D$\Hyph modules.
The map
\ShowEq{sum of maps,,D module}
\ShowEq{sum of maps, 1, D module}
defined by equation
\ShowEq{sum of maps, D module}
is called
\AddIndex{sum of maps}{sum of maps}
$f$ and $g$
and is linear map.
The set \LAB D{A_1}{A_2}
is an Abelian group
relative sum of maps.
}

\AddEq{definition: linear map A module}
{
\begin{ShadedDefinition}
\labelDefinition{linear map \SideWS A module}
Let
\ShowEq{Ai, i=}A
be algebra over commutative ring $D_i$.
Let
\ShowEq{Ai, i=}V
be
\ShowEq{\SideWS column module}{A_i}
module.
Morphism of diagram of representations
\ShowEq{A module diagram for linear}1
into diagram of representations
\ShowEq{A module diagram for linear}2
is called
\AddIndex{linear map}{linear map}
of
\ShowEq{\SideWS column module}{A_1}
module $V_1$ into
\ShowEq{\SideWS column module}{A_2}
module $V_2$.
Let us denote
\ShowEq{set linear maps, \SideWS A12 module}
set of linear maps
of
\ShowEq{\SideWS column module}{A_1}
module $V_1$ into
\ShowEq{\SideWS column module}{A_2}
module $V_2$.
\qed
\end{ShadedDefinition}
}

\DefText[7]{homomorphism of A module 1}
{
\subsection{Homomorphism of \SideWSC \texorpdfstring{$A$}{A}-\VectorSetC of \ColTNS}
\ShowText{homomorphism of A module 2}{#1}{#2}{#3}{#4}{#5}{#6}{#7}
}

\DefText{homomorphism A module}
{
\section{General Definition}

\ShowText{def homomorphism A module}111222{h(p)}{(g\circ a)}

\section{Homomorphism When Rings
\ShowEq{A1=A2 pdf}D}

\ShowText{def homomorphism A module}{}11{}22p{(g\circ a)}

\def\Temp{module}%
\ifx\Temp\VectorSetNS
\ProveTheorem{da in DA->da in D1A}

\ShowConvention{da in DA->da in D1A}
\fi

\section{Homomorphism When \texorpdfstring{$D$}{D}-algebras
\ShowEq{A1=A2 pdf}A}

\ShowText{def homomorphism A module}{}{}1{}{}2pa
}

\DefText{algebra over ring (11)}
{
Let
\ShowEq{Ai, i=}A
be \Division algebra over commutative ring $D_i$.%
}

\DefText{algebra over ring (1)}
{
Let
\ShowEq{Ai, i=}A
be \Division algebra over commutative ring $D$.%
}

\DefText{algebra over ring ()}
{
Let $A$
be \Division algebra over commutative ring $D$.%
}

\DefText{module over algebra (11)}
{
Let
\ShowEq{Ai, i=}V
be \SideWS $A_i$\Hyph \VectorSetNS.
}

\DefText{module over algebra (1)}
{
Let
\ShowEq{Ai, i=}V
be \SideWS $A$\Hyph \VectorSetNS.
}

\DefLabeledDefinition[6]{homomorphism A module}{\SideNS(#1#2#3)}
{
Let diagram of representations
\DrawEq[{#1}{#2}{#3}{1.}]{diagram of representations, \SideWS module}{->1(#1#2#3)}
describe
\SideWS $A_{#2}$\Hyph \VectorSet $V_{#3}$.
Let diagram of representations
\DrawEq[{#4}{#5}{#6}{2.}]{diagram of representations, \SideWS module}{->2(#1#2#3)}
describe
\SideWS $A_{#5}$\Hyph \VectorSet $V_{#6}$.
Morphism
\DrawEq[gf]{homomorphism A module #1#2#3}{\SideNS}
of diagram of representations
\eqRef{diagram of representations, \SideWS module}{->1(#1#2#3)}
into diagram of representations
\eqRef{diagram of representations, \SideWS module}{->2(#1#2#3)}
is called
\AddIndex{homomorphism}{homomorphism}
of \SideWS $A_{#2}$\Hyph \VectorSet $V_{#3}$
into \SideWS $A_{#5}$\Hyph \VectorSet $V_{#6}$.
Let us denote
\ShowEq{set homomorphisms, A module #1#2#3 \SideNS}
set of homomorphisms
of \SideWS $A_{#2}$\Hyph \VectorSet $V_{#3}$
into \SideWS $A_{#5}$\Hyph \VectorSet $V_{#6}$.
}

\DefLabeledTheorem[8]{define homomorphism A module}{\SideNS(#1#2#3)}
{
The homomorphism
\newline
\FrameEqRef[gf]{homomorphism A module #1#2#3}{\SideNS}
\newline
of \SideWS $A_{#2}$\Hyph \VectorSet $V_{#3}$
into \SideWS $A_{#5}$\Hyph \VectorSet $V_{#6}$
satisfies following equalities
\ShowEq{define homomorphism of A module #1#2#3}{#1}{#2}{#3}{#7}{#8}
}

\DefProof[6]{define homomorphism A module}
{
\ShowText{define homomorphism of vector space(#2)}{#1}{#2}{#3}{#4}
The equality
\eqRef{homomorphism, f v+w=}{f(#1#2)\SideWS A module}
follows from the definition
\refDefinition{homomorphism A module}{\SideNS(#1#2#3)},
since, accorfing to the definition
\RefDefinition{morphism of representations of universal algebra},
the map $f$ is homomorpism of Abelian group.
The equality
\eqRef{\SideWS homomorphism, f av=}{(#1#2)\SideWS A module}
follows from the equality
\EqRef{morphism of representations of universal algebra}
because the map
\DrawEq[gf]{homomorphism A module #2#3}{}
is morphism of representation $g_{1.34}$
into representation $g_{2.34}$.
}

\DefTheorem{L(An;B) is free D module}
{
Let $A_1$, ..., $A_n$, $B$ be free modules over commutative ring $D$.
$D$\Hyph module
\ShowEq{L(A->B)}D{A_1\times...\times A_n}B{}
is free $D$\Hyph module.
}

\DefLabeledTheoremNote[5]{linear map of D module}{#1#2}
{
Linear map
\newline
\FrameEqRef[fV]{homomorphism D algebra #1#2}{module}
\newline
of $D_{#1}$\Hyph module $V_{#2}$
into $D_{#3}$\Hyph module $V_{#4}$
satisfies to equalities\,\footnotemark
\ShowEq{define homomorphism of D module #1#2}{#1}{#2}{#5}
}
{\,
In some books
(for instance, on page \citeBib{Serge Lang}\Hyph 119) the theorem
\refTheorem{linear map of D module}{#1#2}
is considered as a definition.
}

\DefText{linear map of D module (11)}
{
According to definitions
\RefDefinition[\RefRepresentation]{morphism of representations of universal algebra},
the map $h$
is homomorphism of ring $D_1$
into ring $D_2$.
Equalities
\ShowEq{ref linear map of D module 1}
follow from this statement.
}

\DefText{linear map of D module (1)}
{
}

\DefProof[2]{linear map of D module}
{
\ShowText{linear map of D module (#1#2)}
From the definition
\RefDefinition[\RefRepresentation]{morphism of representations of universal algebra},
it follows that
the map $f$ is a homomorphism of the Abelian group $V_1$
into the Abelian group $V_2$ (the equality
\eqRef{homomorphism, f v+w=}{f D module #1#2}).
The equality
\eqRef{left homomorphism, f av=}{f D module #1#2}
follows from the equality
\EqRef[\RefRepresentation]{morphism of representations of universal algebra}.
}

\DefDefinition{polylinear map of modules}
{
Let $D$ be the commutative ring.
Reduced polymorphism of $D$\Hyph modules
$A_1$, ..., $A_n$ into $D$\Hyph module $S$
\ShowEq{f:A->B}f{\Times}S
is called
\AddIndex{polylinear map}{polylinear map} of $D$\Hyph modules
$A_1$, ..., $A_n$
into $D$\Hyph module $S$.

We denote
\ShowEq{set polylinear maps}
the set of polylinear maps
of $D$\Hyph modules
$A_1$, ..., $A_n$
into $D$\Hyph module
$S$.
Let us denote
\ShowEq{set polylinear maps An}DAS
set of $n$\hyph linear maps
of $D$\Hyph module $A$ ($A_1=...=A_n=A$)
into $D$\Hyph module
$S$.
}

\DefTheorem{polylinear map of modules}
{
Let $D$ be the commutative ring.
The polylinear map of $D$\Hyph modules
$A_1$, ..., $A_n$
into $D$\Hyph module $S$
\ShowEq{f:A->B}f{\Times}S
satisfies to equalities
\DrawEq[f]{f(ai+bi)=fai+fbi}{}
\DrawEq[f]{f(pai)=pfai}{}
\ShowEq{polylinear map of algebras, 1}{A_i}D
}

\DefTheorem{there exists tensor product of modules}
{
Let $A_1$, ..., $A_n$ be
modules over commutative ring $D$.
There exists  and unique \AddIndex{tensor product}{tensor product}
\ShowEq{tensor product of modules}
\ShowEq{f:xA->oxA}
of $D$\Hyph modules $A_1$, ..., $A_n$.
We use notation
\ShowEq{fxa=oxa}
for the image of the map $f$.
}

\DefTheorem{tensor product and polylinear map}
{
Let $A_1$, ..., $A_n$ be
modules over commutative ring $D$.
Let
\ShowEq{map f, 1, tensor product}
be
polylinear map defined by the equality
\ShowEq{map f, tensor product}
Let
\ShowEq{map g, tensor product}
be polylinear map into $D$\Hyph module $V$.
There exists a linear map
\ShowEq{map h, tensor product}
such that the diagram
\ShowEq{map gh, tensor product}
is commutative.
The map \(h\) is defined by the equality
\ShowEq{g=h, tensor product}
}

\DefProof{polylinear map of modules}
{
The theorem follows from definitions
\RefDefinition[\RefRepresentation]{reduced polymorphism of representations},
\refDefinition{linear map of D module}1,
\RefDefinition{polylinear map of modules}
and from the theorem
\refTheorem{linear map of D module}1.
}

\DefProof{sum of polylinear maps, module}
{
According to the theorem
\RefTheorem{polylinear map of modules}
\DrawEq[f]{f(ai+bi)=fai+fbi}{sum of maps f}
\DrawEq[f]{f(pai)=pfai}{sum of maps f}
\DrawEq[g]{f(ai+bi)=fai+fbi}{sum of maps g}
\DrawEq[g]{f(pai)=pfai}{sum of maps g}
The equality
\ShowEq{sum of maps, 31, polylinear}
follows from the equalities
\EqRef{sum of maps, polylinear},
\eqRef{f(ai+bi)=fai+fbi}{sum of maps f},
\eqRef{f(ai+bi)=fai+fbi}{sum of maps g}.
The equality
\ShowEq{sum of maps, 32, polylinear}
follows from the equalities
\EqRef{sum of maps, polylinear},
\eqRef{f(pai)=pfai}{sum of maps f},
\eqRef{f(pai)=pfai}{sum of maps g}.
From equalities
\EqRef{sum of maps, 31, polylinear},
\EqRef{sum of maps, 32, polylinear}
and from the theorem
\RefTheorem{polylinear map of modules},
it follows that the map
\EqRef{sum of maps, 1, polylinear}
is linear map of $D$\Hyph modules.

Let
\ShowEq{module of polylinear maps, 1}
For any
\ShowEq{module of polylinear maps, 2}
Therefore, sum of polylinear maps is commutative and associative.

From the equality
\EqRef{sum of maps, polylinear},
it follows that the map
\ShowEq{0:A1n->S}
is zero of addition
\ShowEq{0+f=f 1n}
From the equality
\EqRef{sum of maps, polylinear},
it follows that the map
\ShowEq{-f:A1n->S}
is map inversed to map \(f\)
\ShowEq{f-f=0}
because
\ShowEq{f-f=0 1n}
From the equality
\ShowEq{sum of maps, 4, polylinear}
it follows that sum of maps is commutative.
Therefore, the set
\ShowEq{module of polylinear maps}
is an Abelian group.
}

\DefTheorem{module of polylinear maps}
{
Let $D$ be the commutative ring.
Let $A_1$, ..., $A_n$, $S$ be $D$\Hyph modules.
The map
\ShowEq{product of map over scalar,,polylinear}
\ShowEq{product of map over scalar, 1, polylinear}
defined by equality
\ShowEq{product of map over scalar, polylinear}
is polylinear map
and is called
\AddIndex{product of map $f$ over scalar}
{product of map over scalar} $d$.
The representation
\ShowEq{a:LA1n->LA1n}
of ring $D$ in Abelian group
\ShowEq{module of polylinear maps}
generates structure of $D$\Hyph module.
}

\DefProof{module of polylinear maps}
{
According to the theorem
\RefTheorem{polylinear map of modules}
\DrawEq[f]{f(ai+bi)=fai+fbi}{product of map over scalar}
\DrawEq[f]{f(pai)=pfai}{product of map over scalar}
The equality
\ShowEq{product of map over scalar, 31, polylinear}
follows from equalities
\EqRef{product of map over scalar, polylinear},
\eqRef{f(ai+bi)=fai+fbi}{product of map over scalar}.
The equality
\ShowEq{product of map over scalar, 32, polylinear}
follows from equalities
\EqRef{product of map over scalar, polylinear},
\eqRef{f(pai)=pfai}{product of map over scalar}.
From equalities
\EqRef{product of map over scalar, 31, polylinear},
\EqRef{product of map over scalar, 32, polylinear}
and from the theorem
\RefTheorem{polylinear map of modules},
it follows that the map
\EqRef{product of map over scalar, 1, polylinear}
is polylinear map of $D$\Hyph modules.

The equality
\DrawEq{(p+q)f=pf+qf}{polylinear}
follows from the equality
\ShowEq{(p+q)f=pf+qf 1}
The equality
\DrawEq{p(qf)=(pq)f}{polylinear}
follows from the equality
\ShowEq{p(qf)=(pq)f 1}
From equalities
\eqRef{(p+q)f=pf+qf}{polylinear}
\eqRef{p(qf)=(pq)f}{polylinear}
it follows that the map
\EqRef{a:LA1n->LA1n}
is representation of ring $D$
in Abelian group
\ShowEq{module of polylinear maps}.
Since specified representation is effective,
then, according to the definition
\refDefinition{module over algebra}{\SideWS \VectorSetNS}
and the theorem
\RefTheorem{sum of polylinear maps, module},
Abelian group
\ShowEq{L(A->B)}D{A_1}{A_2}{}
is $D$\Hyph module.
}

\DefCorollary{product of linear map over scalar, D module}
{
Let $A_1$, $A_2$ be $D$\Hyph modules.
The map
\ShowEq{product of map over scalar,,D module}
\ShowEq{product of map over scalar, 1, D module}
defined by the equality
\ShowEq{product of map over scalar, D module}
is linear map
and is called
\AddIndex{product of map $f$ over scalar}
{product of map over scalar} $d$.
The representation
\ShowEq{a:LA12->LA12}
of ring $D$ in Abelian group
\ShowEq{L(A->B)}D{A_1}{A_2}{}
generates structure of $D$\Hyph module.
}

\DefLabeledTheoremNote{linear map of A module}{\Base\DF\Module\VF->\Base\DT\Module\VT, \SideNS}
{
Linear map
\ShowEq{show linear map \MapE}gf
of \SideWS $\Base_\DF$\Hyph module $\Module_\VF$
into \SideWS $\Base_\DT$\Hyph module $\Module_\VT$
satisfies to equalities\,\footnotemark
\ShowEq{property linear map \MapE}
}
{
In some books
(for instance, on page \citeBib{Serge Lang}\Hyph 119) the theorem
\refTheorem{linear map of A module}{\Base\DF\Module\VF->\Base\DT\Module\VT, \SideNS}
is considered as a definition.
}

\DefProof[2]{linear map of module}
{
From definitions
\RefDefinition[\RefRepresentation]{morphism of representations of universal algebra},
\refDefinition{linear map of D module}{#1#2},
it follows that
the map $h$ is a homomorphism of the ring $D_1$
into the ring $D_2$ (the equalities
\eqRef{h(d1+d2)=...}{\SideWS module},
\eqRef{h(d1d2)=...}{\SideWS module})
and the map $f$ is a homomorphism of the Abelian group $A_1$
into the Abelian group $A_2$ (the equality
\eqRef{f(a+b)=...}{D1 D2 \SideWS module}).
The equality
\eqRef{f(da)=h... \SideWS module}{\DFDT \SideWS module}
follows from the equality
\eqRef[\RefRepresentation]{morphism of representations of universal algebra, 2m}{representation}.
}

\DefProof{linear map of A module}
{
From definitions
\RefDefinition[\RefRepresentation]{morphism of representations of universal algebra},
\RefDefinition{linear map \SideWS A module},
it follows that
\begin{itemize}
\item
the map $h$ is a homomorphism of the ring $D_1$
into the ring $D_2$ (the equalities
\eqRef{h(d1+d2)=...}{\SideWS module},
\eqRef{h(d1d2)=...}{\SideWS module});
\item
the map $g$ is homomorphism of the Abelian group $A_1$
into the Abelian group $A_2$ (the equality
\eqRef{f(a+b)=...}{g \DFDT \SideWS module});
\item
the map $f$ is homomorphism of the Abelian group $V_1$
into the Abelian group $V_2$ (the equality
\eqRef{f(a+b)=...}{f \DFDT \SideWS module}).
\end{itemize}
Equalities
\eqRef{f(da)=h... \SideWS module}{g \DFDT},
\eqRef{f(da)=h... \SideWS module}{f \DFDT}
follow from the equality
\EqRef[\RefRepresentation]{morphism of representations of F algebra, definition, 2m}.
}

\DefTheorem{complex field over real field}
{
Consider complex field $C$ as two-dimensional algebra over real field.
Let
\ShowEq{basis of complex field}
be the basis of algebra $C$.
Then in this basis product has form
\ShowEq{product of complex field}
and structural constants have form
\ShowEq{structural constants of complex field}
}

\DefProof{complex field over real field}
{
Equalities
\EqRef{product of complex field} and
\EqRef{structural constants of complex field}
follow from the equality $i^2=-1$.
}

\DefLabeledTheorem[1]{maps of conjugation antilinear}{#1}
{
Maps of conjugation
\ShowEq{I... #1}
are \DoVerb antilinear homomorphisms.
}

\DefProof[1]{maps of conjugation antilinear}
{
The product of $#1$\Hyph numbers
\ShowEq{#1 number Ea}a
and
\ShowEq{#1 number Ea}b
has form
\DrawEq[ab]{#1 product aEx}{ab}
\ShowText{maps of conjugation antilinear #1}
}

\DefText{maps of conjugation antilinear H}
{
\ShowEq{H items maps of conjugation antilinear}
}

\DefText{maps of conjugation antilinear O}
{
To prove the theorem, it is enough to consider map $\aU I1$,
because the proof is the same for other maps.
\ShowEq{item maps of conjugation antilinear}O1
}

\DefTheorem{coordinates of map A1 A2, algebra}
{
Let $\Basis e_1$ be basis of the free finite dimensional
$D$\Hyph module $A_1$.
Let $\Basis e_2$ be basis of the free finite dimensional associative
$D$\Hyph algebra $A_2$.
Let
\ShowEq{structural constants, algebra}
be structural constants of algebra $A_2$.
Let $\Basis F$ be the basis
of left \BoxB{A_2}module
\ShowEq{L(A->B)}D{A_1}{A_2}{}
and
\ShowEq{coordinates of map Ik}
be coordinates of map $F_k$ with respect to bases $\Basis e_1$ and $\Basis e_2$.
Coordinates
\ShowEq{Coordinates of map f}
of the map
\ShowEq{f in L(A->B)}D{A_1}{A_2}{}
and its standard components
\ShowEq{standard components of map f}k
are connected by the equation
\DrawEq[f-]{coordinates of map A1 A2, 2, associative algebra}f
}

\DefProof{coordinates of map A1 A2, algebra}
{
Relative to bases
$\Basis e_1$ and $\Basis e_2$, linear maps $f$ and $I_k$ have form
\ShowEq{coordinates of map f 18, associative algebra}
\ShowEq{coordinates of map Ik, associative algebra}
The equality
\ShowEq{coordinates of map A1 A2, 3, associative algebra}
follows from equalities
\EqRef{standard representation of map A1 A2, associative algebra},
\EqRef{coordinates of map f 18, associative algebra},
\EqRef{coordinates of map Ik, associative algebra}.
Since vectors $e_{2\cdot\gik}$
are linear independent and $x^{\gi i}$ are arbitrary,
then the equality
\eqRef{coordinates of map A1 A2, 2, associative algebra}f
follows from the equation
\EqRef{coordinates of map A1 A2, 3, associative algebra}.
}

\DefTheorem{linear map in L(A,A), associative algebra}
{
Let $A_1$ be $D$\Hyph algebra.
Let $A_2$ be free finite dimensional associative $D$\Hyph algebra.
Let $\Basis e$ be basis of $D$ module $A_2$.
Left \BoxB{A_2}module
\ShowEq{L(A->B)}D{A_1}{A_2}{}
has finite
\AddIndex{basis}{basis of algebra L(A,A)} $\Basis I$.
\StartLabelItem
\begin{enumerate}
\item
The linear map
\ShowEq{f in L(A->B)}D{A_1}{A_2}{}
has form
\labelItem{f in L(A,A), 1, associative algebra}
\ShowEq{f in L(A,A), 1, associative algebra}
\item
Its standard representation has form
\labelItem{f in L(A,A), 2, associative algebra}
\ShowEq{f in L(A,A), 2, associative algebra}
\end{enumerate}
}

\DefTheoremNote{standard representation of map A1 A2, associative algebra}
{
Let $A_1$ be free $D$\Hyph module.
Let $A_2$ be free finite dimensional associative $D$\Hyph algebra.
Let $\Basis e$ be basis of $D$\Hyph module $A_2$.
Let $\Basis F$
be the basis of left \BoxB{A_2}module
\ShowEq{L(A->B)}D{A_1}{A_2}.\,\footnotemark
\StartLabelItem
\begin{enumerate}
\item
The map
\DrawEq[f{A_1}{A_2}{}]{f: A->B}{}
has the following expansion
\labelItem{map f generated by basis F}
\DrawEq{map f generated by basis F}{expansion}
where
\ShowEq{fk= in A2xA2}
\item
The map $f$ has the standard representation
\labelItem{standard representation of map A1 A2, associative algebra}
\ShowEq{standard representation of map A1 A2, associative algebra}
\end{enumerate}
}{
If $D$\Hyph module $A_1$ or $D$\Hyph module $A_2$
is not free $D$\Hyph nodule,
then we may consider the set
\ShowEq{Ik 1n}
of linear independent linear maps. The theorem is true for any linear map
\DrawEq[f{A_1}{A_2}{}]{f: A->B}{}
generated by the set of linear maps $\Basis F$.
}

\DefProof{standard representation of map A1 A2, associative algebra}
{
Since $\Basis F$ is the basis of left \BoxB{A_2}module
\ShowEq{L(A->B)}D{A_1}{A_2},
then according to the definition
\ShowEq{ref definition: basis of representation}
and the theorem
\ShowEq{RefTheorem set of vectors generated by set of vectors, module}
there exists expansion
\DrawEq[{A_2}]{expansion of linear map with respect to basis}{A2}
of the linear map $f$ with respect to the basis $\Basis I$.
According to the definition
\ShowEq{ref map j, representation, tensor product}
\DrawEq{f=fkxfk}{module}
The equality
\eqRef{map f generated by basis F}{expansion}
follows from equalities
\eqRef{expansion of linear map with respect to basis}{A2},
\eqRef{f=fkxfk}{module}.
According to theorem
\RefTheorem{standard component of tensor, algebra},
the standard representation of the tensor $f^k$ has form
\ShowEq{standard representation of map A1 A2, 3, associative algebra}
The equation
\EqRef{standard representation of map A1 A2, associative algebra}
follows from equations
\eqRef{map f generated by basis F}{expansion},
\EqRef{standard representation of map A1 A2, 3, associative algebra}.
}

\DefTheorem{map > bullet map}
{
Let $A$ be $D$\Hyph algebra.
For any $A$\Hyph number $a$, the map
\ShowEq{f->ao f}
defined by the equality
\ShowEq{ao f=a f o}
is endomorphism of $D$\Hyph module
\ShowEq{L(A->B)}DAA.
}

\DefProof{map > bullet map}
{
}

\DefTheorem{a > a bullet}
{
$D$\Hyph module
\ShowEq{L(A->B)}DAA
is left $A$\Hyph module generated by the representation
\ShowEq{a > a bullet}
}

\DefProof{a > a bullet}
{
}

\AddEq[3]{theorem: ao Jacobian matrix}
{
\begin{ShadedTheorem}
\labelTheorem{a#1, #2, Jacobian matrix}
The map
\DrawEq[{#1}{#2}]{map ao}{#1#2}
has matrix
\ShowEq{maps of conjugation, Jacobian matrix}{#1}{#2}{#3}
\DrawEq[{#1}{#2}]{a o, Jacobian matrix}{#1#2}
\end{ShadedTheorem}
}

\AddEq[2]{proof: ao Jacobian matrix}
{
\begin{proof}
The product of $#2$\Hyph numbers
\ShowEq{#2 number Ea}a
and
\ShowEq{#2 number #1a}x
has form
\DrawEq[ax]{#2 product a#1x}{}
Therefore, function
\eqRef{map ao}{#1#2}
has Jacobian matrix
\eqRef{a o, Jacobian matrix}{#1#2}.
\end{proof}%
}

\AddEq [2]{item maps of conjugation antilinear}
{

The equality
\ShowEq{#1#2 product ab}
follows from the equalities
\EqRef{#2x= #1},
\eqRef{#1 product aEx}{ab}.
From the equality
\eqRef{#1 product aEx}{ab},
it follows that product of $#1$\Hyph numbers
\ShowEq{#1 number #2a}b
and
\ShowEq{#1 number #2a}a
has form
\ShowEq{#1 #2b*#2a}
From equalities
\EqRef{#1#2 product ab},
\EqRef{#1 #2b*#2a},
it follows that the map $\aU I#2$
is \DoVerb linear antihomomorphism.
}

\AddEq[2]{theorem: ao Jacobian matrix 1}
{
\begin{ShadedTheorem}
\labelTheorem{a#1, #2, Jacobian matrix 1}
\DrawEq[#1]{ao, #2 matrix 1}{#1}
\end{ShadedTheorem}
}

\AddEq[2]{proof: ao Jacobian matrix 1}
{
\begin{proof}
The equality
\eqRef{ao, #2 matrix 1}{#1}
follows from the chain of equalities
\ShowEq{a#1, #2 matrix 2}
\end{proof}%
}

\AddEq[2]{theorem: a* Jacobian matrix}
{
\begin{ShadedTheorem}
\labelTheorem{a*#1, #2, Jacobian matrix}
The map
\DrawEq[{#1}{#2}]{map a*}{#1#2}
has matrix
\ShowEq{maps of conjugation, Jacobian matrix}{#1}{#2}r
\DrawEq[{#1}{#2}]{a *, Jacobian matrix}{#1#2}
\end{ShadedTheorem}
}

\AddEq[2]{proof: a* Jacobian matrix}
{
\begin{proof}
The product of $#2$\Hyph numbers
\ShowEq{#2 number #1a}x
and
\ShowEq{#2 number Ea}a
has form
\ShowEq{#2 product #1xa}
Therefore, function
\eqRef{map a*}{#1#2}
has Jacobian matrix
\eqRef{a *, Jacobian matrix}{#1#2}.
\end{proof}%
}

\AddEq[2]{theorem: a* Jacobian matrix 1}
{
\begin{ShadedTheorem}
\labelTheorem{a*#1, #2, Jacobian matrix 1}
\DrawEq[#1]{a*, #2 matrix 1}{#1}
\end{ShadedTheorem}
}

\AddEq[2]{proof: a* Jacobian matrix 1}
{
\begin{proof}
The equality
\eqRef{a*, #2 matrix 1}{#1}
follows from the chain of equalities
\ShowEq{#1a, #2 matrix 2}
\end{proof}%
}

\DefDefinition{antilinear homomorphism}
{
The map
\ShowEq{f in L(A->B)}DAA{}
is called
\AddIndex{antilinear homomorphism}{antilinear homomorphism}
if the map $f$ satisfies the equality
\ShowEq{fab=fbfa}
}

\DefDefinition{quaternion maps of conjugation}
{
\ShowEq{H maps of conjugation}
Quaternion algebra has following
\AddIndex{maps of conjugation}{map of conjugation}
\ShowEq{H list maps of conjugation}
We also use notation
\ShowEq{I0=E}
}

\DefDefinition{octonion maps of conjugation}
{
\ShowEq{O maps of conjugation}
Octonion algebra has following
\AddIndex{maps of conjugation}{map of conjugation}
\ShowEq{O list maps of conjugation}
We also use notation
\ShowEq{I0=E}
}

\AddEq[1]{theorem: L is left vector space}
{
\begin{ShadedTheorem}
\labelTheorem{L(#1->#1) is left #1-vector space}
$#1\otimes #1$\Hyph module
\ShowEq{L(A->B)}R{#1}{#1}
is left $#1$\Hyph vector space
and has the basis
\ShowEq{I=(E,I)}{#1}.
\end{ShadedTheorem}
}

\AddEq [9]{Let be basis of vector space}
{%
Let the set of vectors
\DrawEq[{#1}{#2}{#3}]{basis e of module #8}{#9}
be a basis of \SideWS $#4_{#5}$\Hyph \VectorSet $#6_{#7}$.%
}%

\AddEq [3]{Let be Basis of vector space}
{%
Let \eV[#1]
be a basis of \SideWS $#2$\Hyph vector space $#3$.%
}%

\AddEq [9]{Let be basis of module}
{%
Let
\DrawEq[{#1_{#2}}{#3}{#4}]{basis e of module \Cols}-
be a basis of #5 $#6_{#7}$\Hyph \VectorSet of \ColsWS $#8_{#9}$.
}

\AddEq [9]{Let be basis of algebra}
{%
Let the set of vectors
\DrawEq[{#1_{#2}}{#3}{#4}]{basis e of module \Cols}{#9}
be a basis of $#5_{#6}$\Hyph algebra of \ColsWS $#7_{#8}$.%
}

\AddEq [9]{Let be basis of algebra and C}
{%
Let the set of vectors
\DrawEq[{#1}{#3}{#4}]{basis e of module \Cols}{#9}
be a basis and
\ShowEq{structural constants of algebra}{#2}kij,
\ShowEq{kij in I}{#4}
be structural constants of $#5_{#6}$\Hyph algebra of \ColsWS $#7_{#8}$.%
}

\AddEq [5]{Let e be basis 1n}
{
Let
\ShowEq{basis e of module 1n}{#1}{#2}
be a basis of #3 $#4$\Hyph module $#5$.
}

\DefText[5]{coordinates of the linear map 1}
{
\item $#1$ is coordinate matrix of $#2_1$\Hyph number
$\Vector #1$
relative the basis \eV[#3]
\DrawEq [{#1}{#3}{}]{va=ae1, module (#5)(\Cols)}{\SideNS-\VectorSet #1 #4}
}

\DefText[5]{coordinates of the linear map 2(1)}
{
\item
\ShowEq{h(a)=...}{#1}{#2}{#3}{#4}
is a matrix of $#5$\Hyph numbers.
}

\DefText[5]{coordinates of the linear map 2()}
{
}

\DefText[7]{coordinates of the linear map 3}
{
\item $#6$ is coordinate matrix of $#2_2$\Hyph number
\DrawEq[#1#6#4]{vb=f(va)}{#1 #7(\Cols)\SideNS-\VectorSetNS}
relative the basis \eV[#3]
\DrawEq [#6{#3}]{va=ae1, module (#5)(\Cols)}{#6 #7\SideNS-\VectorSetNS}
}

\DefText[6]{coordinates of the linear map 4}
{
\item $#1$ is coordinate matrix of set of $#2_2$\Hyph numbers
\ShowEq{Vector f(e1) module}#1{#3}#5#6
relative the basis \eV[#4][.]
}

\DefText[6]{Let be module of}
{
Let $#1_{#2}$ be #3 $#4_{#5}$\Hyph #6 of \ColTNS.
}

\DefText[2]{we identify linear map and matrix}
{
On the basis of theorems
\refTheorem{linear map of D module, coordinates}{#1#2cols},
\refTheorem{matrix generates D module homomorphism}{cols(#1#2)},
as well on the basis of theorems
\refTheorem{linear map of D module, coordinates}{#1#2rows},
\refTheorem{matrix generates D module homomorphism}{rows(#1#2)},
we identify the linear map
\DrawEq[{\Vector f}V]{homomorphism D algebra #1#2}{}
and the matrix $f$ of its presentation.
}

\DefLabeledFootnote[4]{homomorphism of d algebra}{#1#2\Cols}
{
In theorems
\refTheorem{homomorphism of d algebra, coordinates}{#1#2\Cols},
\refTheorem{matrix generates D algebra homomorphism}{\Cols(#1#2)},
we use the following convention.
\ShowEq{Let be basis of algebra and C}11iID{#1}A{#2}-
\ShowEq{Let be basis of algebra and C}22jJD{#3}A{#4}-
}
 
\DefLabeledTheorem[5]{homomorphism of d algebra, coordinates}{#1#2\Cols}
{
The homomorphism\refFootnote{homomorphism of d algebra}{#1#2\Cols}
\DrawEq[{\Vector f}A]{homomorphism D algebra #1#2}{}
of $D_{#1}$\Hyph algebra of \ColsWS $A_{#2}$
into $D_{#3}$\Hyph algebra of \ColsWS $A_{#4}$ has presentation
\DrawEq[f{#5}{}]{f:V1->V2, D module \Cols}{algebra(#1#2)}
\DrawEq[hf{#2}{#4}a]{f o ea=efa i (#1#2)}{(\Cols)algebra}
\DrawEq[hf12]{f o ea=efa (#1#2)(\Cols)}{algebra}
relative to selected bases. Here
\begin{itemize}
\ShowText{coordinates of the linear map 1}aA1{(#1)(#2)}{}
\ShowText{coordinates of the linear map 2(#1)}ahiI{D_{#3}}
\ShowText{coordinates of the linear map 3}aA2f{}b{(#1)(#2)}
\ShowText{coordinates of the linear map 4}fA{#2}{#4}iI
\ShowText{Morphism of D algebra}{#1}{#2}f
\end{itemize}
}

\DefProof[4]{homomorphism of d algebra, coordinates}
{
Equalities
\ShowRef{homomorphism of d algebra, coordinates}{#1}{#2}{#4}
follow from theorems
\refTheorem{linear map of D module, coordinates}{#1#2\Cols},
\refTheorem{homomorphism from A1 to A2, D algebra}{#1#2}.

The equality
\DrawEq{algebra, homomorphism and product eiej(#1#2)}{\Cols}
follows from equalities
\ShowRef{algebra, homomorphism and product eiej}{#1}{#2}{#3}
The equality
\DrawEq{algebra, homomorphism and product eiej 4(#1#2)}{\Cols}
follows from the equality
\ShowRef{algebra, homomorphism and product eiej 1}{#1}{#2}
and the equality
\eqRef{algebra, homomorphism and product eiej(#1#2)}{\Cols}.
The equality
\DrawEq{algebra, homomorphism and product eiej 5}{(#1#2)\Cols}
follows from the equality
\ShowRef{homomorphism, f vw=}{#1}{#2}
The equality
\DrawEq{algebra, homomorphism and product eiej 2(#1#2)}{\Cols}
follows from the equality
\ShowEq{algebra, homomorphism and product, 1}
and the equality
\eqRef{algebra, homomorphism and product eiej 5}{(#1#2)\Cols}.
The equality
\DrawEq{algebra, homomorphism and product eiej 3(#1#2)}{\Cols}
follows from equalities
\eqRef{algebra, homomorphism and product eiej 4(#1#2)}{\Cols},
\eqRef{algebra, homomorphism and product eiej 2(#1#2)}{\Cols}.
The equality
\ShowRef{algebra, homomorphism and product}{#1}{#2}f
follows from the equality
\eqRef{algebra, homomorphism and product eiej 3(#1#2)}{\Cols},
and and the theorem
\refTheorem{coordinates of vector}{\SideNS-\Cols}.
}

\DefText{Representation of ring f(0)=v0}
{
\begin{sloppypar}
The equality
\ShowEq{fa=+0}
follows from the equalities
\EqRef{f(a+b)=f(a)+f(b)},
\EqRef{f(0)=v0}.
Therefore, the map $\Vector 0$
is zero of ring
\ShowEq{End +A}
\end{sloppypar}
}

\DefLabeledTheorem[5]{linear map of D module, coordinates}{#1#2\Cols}
{
Linear map\refFootnote{homomorphism of D module}{\Cols(#1#2)}
\DrawEq[{\Vector f}V]{homomorphism D algebra #1#2}{Vector module, coordinates \Cols}
of $D_{#1}$\Hyph module of \ColsWS $A_{#2}$
into $D_{#3}$\Hyph module of \ColsWS $A_{#4}$ has presentation
\DrawEq[hf{#2}{#4}{}]{f o ea=efa (#1#2)(\Cols)}{module}
\DrawEq[hf{#2}{#4}a]{f o ea=efa i (#1#2)}{(\Cols)module}
\DrawEq[f{#5}{}]{f:V1->V2, D module \Cols}{#1#2}
relative to selected bases. Here
\begin{itemize}
\ShowText{coordinates of the linear map 1}aV1{(#1)(#2)}{}
\ShowText{coordinates of the linear map 2(#1)}ahiI{D_{#3}}
\ShowText{coordinates of the linear map 3}aV2f{}b{(#1)(#2)}
\ShowText{coordinates of the linear map 4}fV{#2}{#4}iI
\end{itemize}
}

\DefProofRef[5]{linear map of D module, coordinates}{#1#2\Cols}
{
Since
\newline
\FrameEqRef[{\Vector f}V]{homomorphism D algebra #1#2}{Vector module, coordinates \Cols}
\newline
is a linear map, then the equality
\ShowEq{vb=h(e1)a (#1)(#2)(\Cols)\SideNS}
follows from equalities
\ShowRef{vb=h(e1)a}{#1}{#2}{#5}
$V_2$\Hyph number
\ShowEq{f(e1)(\Cols)\SideNS}
has expansion
\DrawEq{f(e1)=(\Cols)\SideNS}{(#1)(#2)}
relative to basis $\Basis e_2$.
Combining \EqRef{vb=h(e1)a (#1)(#2)(\Cols)\SideNS}
and \eqRef{f(e1)=(\Cols)\SideNS}{(#1)(#2)}, we get
\ShowEq{vb=a f e (#1)(#2)(\Cols)\SideNS}
\eqRef{f:V1->V2, D module \Cols}{#1#2}
follows from comparison of
\newline
\FrameEqRef[b2]{va=ae1, module ()(\Cols)}{b (#1)(#2)\SideNS-\VectorSetNS}
\newline
and \EqRef{vb=a f e (#1)(#2)(\Cols)\SideNS} and
the theorem
\refTheorem{coordinates of vector}{\SideNS-\Cols}.
}

\DefDefinitionNote{algebra over ring}
{
Let $D$ be commutative ring.
$D$\Hyph module $A$ is called
\AddIndex{algebra over ring}{algebra over ring} $D$
or
\AddIndex{$D$\Hyph algebra}{D algebra},
if we defined product\,\footnotemark
in $A$
\DrawEq{product in D algebra}{definition}
where $C$ is bilinear map
\ShowEq{product in algebra, definition 1}
If $A$ is free
$D$\Hyph module, then $A$ is called
\AddIndex{free algebra over ring}{free algebra over ring} $D$.
}
{
I follow the definition
given in
\citeBib{Richard D. Schafer}, page 1,
\citeBib{0105.155}, page 4. The statement which
is true for any $D$\Hyph module,
is true also for $D$\Hyph algebra.
}

\DefText{algebra inherits type of module}
{
$D$\Hyph algebra $A$ inherits type of $D$\Hyph module $A$.
}

\DefLabeledDefinition{type of algebra}{\Cols}
{
$D$\Hyph algebra $A$ is called
$D$\Hyph algebra of \ColTNS,
if $D$\Hyph module $A$
is $D$\Hyph module of \ColTNS.
}

\DefLabeledTheorem{associative product in algebra}{\Cols}
{
Since the algebra $A$ is commutative, then
\DrawEq{commutative product in algebra, 1}{\Cols}
Since the algebra $A$ is associative, then
\DrawEq{associative product in algebra, 1}{\Cols}
}

\DefProof{associative product in algebra}
{
For commutative algebra,
the equation
\eqRef{commutative product in algebra, 1}{\Cols}
follows from equation
\ShowEq{commutative product in algebra}
\begin{sloppypar}
\noindent
For associative algebra,
the equation
\eqRef{associative product in algebra, 1}{\Cols}
follows from equation
\end{sloppypar}
\ShowEq{associative product in algebra}
}

\DefLabeledTheorem{product in algebra}{\Cols}
{
\ShowEq{Let be basis of algebra}{}{}iID{}A{}-
Let
\ShowEq{a b in basis of algebra}
We can get the product of $a$, $b$ according to rule
\DrawEq{product in algebra}{\Cols}
where
\ShowEq{structural constants of algebra symb}
are \AddIndex{structural constants}{structural constants}
of algebra $A$ over ring $D$.
The product of basis vectors in the algebra $A$ is defined according to rule
\DrawEq{product of basis vectors, algebra}{\Cols}
}

\DefProof{product in algebra}
{
The equality
\eqRef{product of basis vectors, algebra}{\Cols}
is corollary of the statement that \eV
is the basis of $D$\Hyph algebra $A$.
Since the product in the algebra is a bilinear map,
then we can write the product of $a$ and $b$ as
\DrawEq{product in algebra, 1}{\Cols}
From equalities
\eqRef{product of basis vectors, algebra}{\Cols},
\eqRef{product in algebra, 1}{\Cols},
it follows that
\DrawEq{product in algebra, 2}{\Cols}
Since \eV is a basis of the algebra $A$, then the equality
\eqRef{product in algebra}{\Cols}
follows from the equality
\eqRef{product in algebra, 2}{\Cols}.
}

\DefText{D-module type}
{
Organization of coordinates of vector in matrix
is called $D$\Hyph module type.

In this section, we consider column vector
and row vector.
it is evident that there exist other forms of representation of vector.
For instance, we can represent coordinates of vector as
\nmTimes nm matrix or as triangular matrix.
Format of representation depends on considered problem.
}

\DefText{D-module type 1}
{
If a statement depends on the format of representation of vector,
we will specify either type of $D$\Hyph module,
or format of representation of basis.
Both ways of specifying the type of $D$\Hyph module are equivalent.
}

\DefText{Morphism of D algebra, injection 11}
{
(and therefore,
\ShowEq{Morphism of D algebra, injection}
}

\DefText{Morphism of D algebra, injection 1}
{
}

\DefText[3]{Morphism of D algebra}
{
\item There is relation between the matrix of homomorphism
and structural constants
\DrawEq[{#3}]{algebra, homomorphism and product (#1#2)}{\Cols(\SideNS)}
}

\DefText{algebra, homomorphism and product 11}
{
Since the map $h$ is homomorphism of rings,
then the equality
\DrawEq{algebra, homomorphism and product, 4}{\Cols}
follows from the equality
\eqRef{algebra, homomorphism and product eiej(11)}{\Cols}.
The equality
\DrawEq{algebra, homomorphism and product, 5 11}{\Cols}
follows from the theorem
\refTheorem{homomorphism of d algebra, coordinates}{11\Cols}
and the equality
\eqRef{algebra, homomorphism and product, 4}{\Cols}.
}

\DefText{algebra, homomorphism and product 1}
{
The equality
\DrawEq{algebra, homomorphism and product, 5 1}{\Cols}
follows from the theorem
\refTheorem{homomorphism of d algebra, coordinates}{1\Cols}
and the equality
\eqRef{algebra, homomorphism and product eiej(1)}{\Cols}.
}

\DefDefinition{commutator of algebra}
{
The \AddIndex{commutator}{commutator of algebra}
\ShowEq{commutator of algebra}
measures commutativity in $D$\Hyph algebra $A$.
$D$\Hyph algebra $A$ is called
\AddIndex{commutative}{commutative D algebra},
if
\ShowEq{commutative D algebra}
}

\DefDefinitionNote{nucleus of algebra}
{
The set\,\footnotemark
\ShowEq{nucleus of algebra}
is called the
\AddIndex{nucleus of an $D$\Hyph algebra $A$}{nucleus of algebra}.
}
{The definition is based on
the similar definition in
\citeBib{Richard D. Schafer}, p. 13.}

\DefDefinition{associator of algebra}
{
The \AddIndex{associator}{associator of algebra}
\ShowEq{associator of algebra}
\ShowEq{associator of algebra =}
measures associativity in $D$\Hyph algebra $A$.
$D$\Hyph algebra $A$ is called
\AddIndex{associative}{associative D algebra},
if
\ShowEq{associative D algebra}
}

\DefTheorem{associator of algebra}
{
Let $\Basis e$ be basis of $D$\Hyph algebra $A$. Then
\ShowEq{associator of algebra = A}
where
\AddIndex{coordinates of associator}{coordinates of associator}
are defined by the equality
\ShowEq{coordinates of associator}
\ShowEq{coordinates of associator =}
}

\DefProof{associator of algebra}
{
The equality
\ShowEq{associator of algebra = Ae}
follows from the equality
\EqRef{associator of algebra =}.
The equality
\EqRef{coordinates of associator =}
follows from the equality
\EqRef{associator of algebra = Ae}.
}

\DefDefinitionNote{center of algebra}
{
The set\,\footnotemark
\ShowEq{center of algebra}
is called the
\AddIndex{center of an $D$\Hyph algebra $A_1$}{center of algebra}.
}{
The definition is based on
the similar definition in
\citeBib{Richard D. Schafer}, page 14.
}

\DefDefinition{otimes -}
{
Bilinear map
\ShowEq{otimes -}
is defined by the equality
\ShowEq{otimes -, 1}
}

%% file: Stmt.Module.Eq.tex
\ePrints{1506.00061,7287-9339,0803.3276,7100-9423,357606033,2023.01.07}
\ifx\Semafor\ValueOff
\input{Stmt.Module.Ref}
\fi

\newcommand\aUD[3]{\ensuremath{{#1}^{\gi{#2}}_{\gi{#3}}}}%
\newcommand\Entry[5][]{#2_{#3}^{}\AUD{{}}{#4}{#5}{}_{#1}}%
\newcommand\aU[3][]{\ensuremath{#2^{#1\gi{#3}}}}%
\newcommand\aD[3][]{\ensuremath{#2_{#1\gi{#3}}}}%
\newcommand\eV[1][]
{
\def\argA{#1}%
\eVA%
}
\newcommand\eVA[1][]{\ensuremath{\Basis e_{\argA}}{#1} }%
\newcommand\EV[2]{e_{#1\cdot\gi{#2}}}
\newcommand\ECol[2]{\ensuremath{e_{#1\gi{#2}}}}
\newcommand\ERow[2]{\ensuremath{\aU{e_{#1}}{#2}}}
\newcommand\AOn[1][A]{\ensuremath{#1^{n+1\otimes}}}
\newcommand\LAnA[1][A]{\ensuremath{\mathcal L(D;#1^n\rightarrow #1)}}
\newcommand\nTimes[1][n]{$\gi{#1}\times\gi{#1}$ }
\newcommand\nmTimes[2]{$\gi{#1}\times\gi{#2}$ }
\def\Eij{\ARow ei\ARow ej}

\newcommand{\LAVW}{\mathcal L(A\CRcirc;V\rightarrow W)}

\newcommand\Ei
{
\begin{pmatrix}
1&0&0&0\\
0&-1&0&0\\
0&0&1&0\\
0&0&0&1
\end{pmatrix}
}

\newcommand\Ej
{
\begin{pmatrix}
1&0&0&0\\
0&1&0&0\\
0&0&-1&0\\
0&0&0&1
\end{pmatrix}
}

\newcommand\Ek
{
\begin{pmatrix}
1&0&0&0\\
0&1&0&0\\
0&0&1&0\\
0&0&0&-1
\end{pmatrix}
}

\newcommand\KR[1]
{
\left(
\begin{array}{rr@{}rr@{}rr@{}r}
\end{array}
\right)
}

\newcommand\PMatrix[4][]
{
\begin{pmatrix}
\aUD {#2}11#1&...&\aUD {#2}1{#3}#1\\
...&...&...\\
\aUD {#2}{#4}1#1&...&\aUD {#2}{#4}{#3}#1
\end{pmatrix}
}

\newcommand\ColMatrix[2]
{
\begin{pmatrix}
\aU {#1}1\\...\\ \aU {#1}{#2}
\end{pmatrix}
}

\newcommand\RowMatrix[2]
{
\begin{pmatrix}
\aD {#1}1&...&\aD {#1}{#2}
\end{pmatrix}
}

\AddEq{basis e of V}
{
$\eV=(\aD e1,...,\aD en)$
}

\AddEq{basis e=1}
{
$\eV=\{1\}$.
}

\AddEq[1]{Let aij rco bij}
{
\ShowEq{a in Aoxn}{\aUD aij}2,
\ShowEq{a in Aoxn}{\aUD bij}2#1
}

\AddEquation{aij rco bij}
{
\PMatrix anm
\RCcirc
\PMatrix bkn
=
\begin{pmatrix}
\aUD a1i\circ\aUD bi1&...&\aUD a1i\circ\aUD bik\\
...&...&...\\
\aUD ami\circ\aUD bi1&...&\aUD ami\circ\aUD bik
\end{pmatrix}
}

\AddEq{a rco b ij}
{
\aUD{(a\RCcirc b)}ij=\aUD aik\circ\aUD bkj
}

\AddEq{matrix vx}
{
\[
a=
\begin{pmatrix}
\aU v1&\aU x1\\
\aU v2&\aU x2
\end{pmatrix}
\]
}

\AddEquation{qdet21 vx}
{
\RCdet 21a=
\aU v2-\aU x2(\aU x1)^{-1}\aU v1
}

\AddEq{any quaternion}
{
\[
a=\aU a0+\aU a1i+\aU a2j+\aU a3k
\]
}

\AddEq{any quaternion i j ij}
{
\[
a=\aU a0+\aU a1i+\aU a2j+\aU a3ij
\]
}

\DefText{D-Algebra Type}
{
\ShowConvention{basis of algebra as basis of module}

\ShowConvention{we use separate color for index of element}

\ShowText{algebra inherits type of module}

\ShowEq{def ab=ba}%
\TwoColText
{
\ShowEq{\DefCol}%
\ShowDefinition{type of algebra}
\ShowConvention{unit of algebra in basis}
}
{
\ShowEq{\DefRow}%
\ShowDefinition{type of algebra}
\ShowConvention{unit of algebra in basis}
}
}

\DefText{Preliminary product in algebra}
{
\TwoColText
{
\ShowEq{\DefCol}%
\ShowTheorem{product in algebra}
\proofTheorem{\RefLinearMap}{product in algebra}{\Cols}
}
{
\ShowEq{\DefRow}%
\ShowTheorem{product in algebra}
\proofTheorem{\RefLinearMap}{product in algebra}{\Cols}
}
}

\AddEq{spectrum of matrix}
{
\ProductType $\mathrm{spec}(a)$
}

\AddEq[2]{left map a En}
{
\ensuremath{\Vector{\eV[#2]#1}}
}

\AddEq[2]{right map a En}
{
\ensuremath{\Vector{#1\eV[#2]}}
}

\AddEq[2]{fov=b e o v}
{
\Vector f\circ{#2}=
\ShowEq{\SideWS map a En}{#1}{}
\circ{#2}
}

\AddEq [3]{L(A->B)}
{
\ensuremath{\mathcal L(#1;#2\rightarrow #3)}%
}

\AddEq{conjugation in algebra, 4}
{
\[
\begin{matrix}
\giA\le\gik\le\gin&\giA\le\gil\le\gin&\giA\le\gi p\le\gin
\end{matrix}
\]
}

\AddEquation{structural constants Cikl}
{
\begin{matrix}
\aUD Ck{k0}=\aUD Ck{0k}=1
&
\aUD Cp{kl}=-\aUD Cp{lk}
\end{matrix}
}

\AddEquation{structural constants C000}
{
\begin{matrix}
\hphantom{\aUD C0{00}=}
\aUD C0{00}=1
&\aUD C0{kl}=\hphantom{-}\aUD C0{lk}
\end{matrix}
}

\AddEq{x=c12+j}
{
x=C_1+C_2i+\frac 12 j
}

\AddEq[1]{ix-xi-1}
{
p_{#1}(x)=ix-xi-#1
}

\AddEq{characteristic polynomial}
{
\det(a-b\aD En)
}

\AddEq{projection maps, quaternion}
{
\begin{matrix}
P^{\gi 0}\circ a=a^{\gi 0}&P^{\gi 0\cdot}_{}{}^{\gi 0}_{\gi 0}=1&
R^{\gi 0}\circ a=a^{\gi 0}&R^{\gi 0\cdot}_{}{}^{\gi 0}_{\gi 0}=1
\\
P^{\giA}\circ a=a^{\giA}i&P^{\giA\cdot}_{}{}^{\giA}_{\giA}=1&
R^{\giA}\circ a=a^{\giA}&R^{\giA\cdot}_{}{}^{\gi 0}_{\giA}=1
\\
P^{\gi 2}\circ a=a^{\gi 2}j&P^{\gi 2\cdot}_{}{}^{\gi 2}_{\gi 2}=1&
R^{\gi 2}\circ a=a^{\gi 2}&R^{\gi 2\cdot}_{}{}^{\gi 0}_{\gi 2}=1
\\
P^{\gi 3}\circ a=a^{\gi 3}k&P^{\gi 3\cdot}_{}{}^{\gi 3}_{\gi 3}=1&
R^{\gi 3}\circ a=a^{\gi 3}&R^{\gi 3\cdot}_{}{}^{\gi 0}_{\gi 3}=1
\end{matrix}
}

\AddEq{P0z=}
{
P^{\gi 0}\circ x=\frac 14(1\otimes 1-i\otimes i-j\otimes j-k\otimes k)\circ x
=\frac 14(x-ix i-jx j-kx k)
}

\AddEq{P1z=}
{
P^{\gi 1}\circ x=\frac 14(1\otimes 1-i\otimes i+j\otimes j+k\otimes k)\circ x
=\frac 14(x-ix i+jx j+kx k)
}

\AddEq{P2z=}
{
P^{\gi 2}\circ x=\frac 14(1\otimes 1+i\otimes i-j\otimes j+k\otimes k)\circ x
=\frac 14(x+ix i-jx j+kx k)
}

\AddEq{P3z=}
{
P^{\gi 3}\circ x=\frac 14(1\otimes 1+i\otimes i+j\otimes j-k\otimes k)\circ x
=\frac 14(x+ix i+jx j-kx k)
}

\AddEquation{polynomial*polynomial}
{
*:A^{n\otimes}\times A^{m\otimes}
\rightarrow A^{n+m-1\otimes}
}

\AddEquation{polynomial*polynomial, 1}
{
(a_1\otimes...\otimes a_n)*(b_1\otimes...\otimes b_n)
=
a_1\otimes...\otimes a_{n-1}\otimes a_nb_1\otimes b_2\otimes...\otimes b_n
}

\AddEq{polynomials a o x,b o x}
{
$a\circ x^n$, $b\circ x^m$
}

\AddEquation{(a*b)x=(ax)(bx)}
{
(a*b)\circ x^{n+m}=(a\circ x^n)(b\circ x^m)
}

\AddEq{x=i,j}
{
$x=i$, $x=j$.
}

\AddEquation{p=a1 o ij+a2 o ji}
{
p(x)=a_1\circ((x-i)(x-j))+a_2\circ((x-j)(x-i))
}

\AddEq{r=xx+1}
{
r(x)=x^2+1
}

\AddEquation{a1 o ij+a2 o ji=xx+1}
{
a_1\circ(x^2-ix-xj+k)+a_2\circ(x^2-jx-xi-k)=x^2+1
}

\AddEq{square root}
{
\symb{\sqrt a}{square root}{}%
$x=\ShowSymbol{square root}{}$
}

\AddEq{x2=a}
{
x^2=a
}

\AddEq{Re sqrt a ne 0}
{
$\re\sqrt a\ne 0$,
}

\AddEq{x=x12}
{
$x=x_1$, $x=x_2$
}

\AddEquation{x2=-x1}
{
x_2=-x_1
}

\AddEquation{|x|=sqrt a}
{
\begin{matrix}
x\in\im H&a\in\re H&|x|=\sqrt{-a}
\end{matrix}
}

\AddEq[1]{aij=.ox.}
{
\Entry a{#1}ij=\AoxB{\Entry[s0]a{#1}ij}{\Entry[s1]a{#1}ij}
}

\AddEquation{vi=aij o wj}
{
\aU vi=(\AoxB{\Entry[s1]a{}ij}{\Entry[s2]a{}ij})\circ\aU wj
=\Entry[s1]a{}ij\aU wj\Entry[s2]a{}ij
}

\AddEq[1]{a1 o+a2 o=xx+1 (ij) x x2=}
{
a_1\circ #1+a_2\circ #1=#1
}

\AddEquation{a1 o+a2 o=xx+1 xx}
{
a_1\circ x^2+a_2\circ x^2=x^2
}

\AddEquation{a1 o+a2 o=xx+1 x}
{
a_1\circ (ix)+a_1\circ (xj)+a_2\circ (jx)+a_2\circ (xi)=0
}

\AddEquation{a1 o+a2 o=xx+1 1}
{
a_1\circ k-a_2\circ k=1
}

\AddEquation{a1 o+a2 o=xx+1 x x=1}
{
a_1\circ i+a_1\circ j+a_2\circ j+a_2\circ i=0
}

\AddEq{D*V->V to product}
{
(a,v)\rightarrow av
}

\AddEq{A*V->V to left product}
{
(a,v)\rightarrow av
}

\AddEq{A*V->V to right product}
{
(v,a)\rightarrow va
}

\AddEq{A**V->V to left product}
{
(a,v)\rightarrow a*v
}

\AddEq{A**V->V to right product}
{
(v,a)\rightarrow v*a
}

\AddEq{a rc x=b pdf}
{
\texorpdfstring{$a\RCstar x=b$}{a rc x=b}
}

\AddEquation{coordinates of map f 18, associative algebra}
{
f\circ x=f^{\gii}_{\gij}x^{\gij}e_{2\cdot\gii}
}

\AddEquation{coordinates of map Ik, associative algebra}
{
F_k\circ x=F_{k\cdot}^{}{}^{\gii}_{\gij}x^{\gij}e_{2\cdot\gii}
}

\AddEquation{coordinates of map A1 A2, 3, associative algebra}
{
\begin{split}
f^{\gik}_{\gil}x^{\gil}e_{2\cdot\gik}
&=
f^{k\cdot\gi{ij}} e_{2\cdot\gii}
F_{k\cdot}^{}{}^{\gim}_{\gil}x^{\gil}e_{2\cdot\gim}
e_{2\cdot\gij}
\\
&=
f^{k\cdot\gi{ij}}
F_{k\cdot}^{}{}^{\gim}_{\gil}x^{\gil}
C^{\gi p}_{\gi{im}}
C^{\gik}_{\gi{pj}}
e_{2\cdot\gik}
\end{split}
}

\AddEq [2]{coordinates of map A1 A2, 2, associative algebra}
{
\def\Cdot{}
\def\Temp{-}%
\edef\Tempa{#2}%
\ifx\Tempa\Temp%
\else
\def\Cdot{#2\cdot}
\fi
#1^{\gik}_{\gil}=#1^{k\cdot\gi{ij}}
F_{k\cdot}^{}{}^{\gim}_{\gil}
C_{\Cdot}{}^{\gi p}_{\gi{im}}C_{\Cdot}{}^{\gik}_{\gi{pj}}
}

\AddEq{1 ox i-i ox 1}
{
\[
i\otimes 1-1\otimes i\in H^{2\otimes}
\]
}

\AddEquation{1 ox i-i ox 1 matrix}
{
\begin{pmatrix}
0&-1&0&0\\
1&0&0&0\\
0&0&0&-1\\
0&0&1&0
\end{pmatrix}
-
\begin{pmatrix}
0&-1&0&0\\
1&0&0&0\\
0&0&0&1\\
0&0&-1&0
\end{pmatrix}
=
\begin{pmatrix}
0&0&0&0\\
0&0&0&0\\
0&0&0&-2\\
0&0&2&0
\end{pmatrix}
}

\AddEquation{ax^2bxcx^3d}
{
\begin{split}
ax^2bxcx^3d&=ax(xbxcx^3)=(ax)x(bxcx^3d)=(ax^2b)x(cx^3d)
\\
&=(ax^2bxc)x(x^2d)=(ax^2bxcx)x(xd)=(ax^2bxcx^2)xd
\end{split}
}

\AddEquation{f in L(A,A), 2, associative algebra}
{
f=
a^{k\cdot \gi{ij}}(e_{\gii}\otimes e_{\gij})\circ F_k
=
a^{k\cdot \gi{ij}}e_{\gii}I_ke_{\gij}
}

\AddEquation{standard representation of map A1 A2, associative algebra}
{
f=
f^{k\cdot\gi{ij}}(\aU ei\otimes \aU ej)\circ F_k
=
f^{k\cdot\gi{ij}}\aU eiF_k\aU ej
}

\AddEquation{f(a)+f(b)=}
{
(f(a)+f(b))(x)=f(a)(x)+f(b)(x)
}

\AddEquation{f(a+b)=f(a)+f(b)}
{
f(a+b)=f(a)+f(b)
}

\AddEquation{f(a-b)=0}
{
f(a-b)=\Vector 0
}

\AddEquation{f(0)=v0}
{
f(0)=\Vector 0
}

\AddEq{0:A->A}
{
\ShowEq{f:A->B}{\Vector 0}AA
}

\AddEquation{0ov=0}
{
\Vector 0\circ v=0
}

\AddEquation{fa=fa+f0}
{
f(a)(x)=f(a+0)(x)=f(a)(x)+f(0)(x)
}

\AddEquation{f0x=0}
{
f(0)(x)=0
}

\AddEq[3]{a in Aoxn}
{
$#1\in A^{#2\otimes}$#3
}

\AddEq[2]{Aoxn}
{
$A^{#1\otimes}$#2
}

\AddEq{set of endomorphisms D module}
{
\symb{\End(D,V)}{set of endomorphisms}1
}

\AddEq[1]{End DV}
{
$\End(D,V)$#1
}

\DefText{definition of A module}
{
\ShowEq{def AModule}
\ShowEq{def universal}

\ShowDefinition{module over associative algebra}

\ShowEq{def AVector}
\ShowDefinition{module over algebra}

\ShowEq{def AModule}
\ShowTheorem{module over algebra}
\ShowProof{module over algebra}

\def\TheoremStyle{ParacolTheorem}%
\ShowTheorem{definition of A module}
\def\TheoremStyle{ShadedTheorem}%
\ShowTheorem{definition of A module, property}
\ShowProof{definition of A module}
}

\AddEq{X0=v in J}
{
$X_0=v\subseteq J(v)$.
}

\DefText{linear combination of vectors}
{
\ProveTheorem{set of vectors generated by set of vectors}

\ShowDefinition{linear combination of vectors}

\ShowConvention{sum av() convention}

\ProveTheorem{linearly depends on rest of vectors}

\ShowText{0=0vi}

\ShowDefinition{linearly independent vectors}
}

\DefText{basis of module}
{
\ShowText{generating set of module}

\ShowDefinition{generating set of module}

\ShowText{quasibasis of module}

\ShowDefinition{basis of module}

\ShowDefinition{coordinates of vector}

\ProveTheorem{quasibasis of module}

\ProveTheorem{linear dependence between vectors of basis}

\ProveTheorem{coordinates of vector with linear dependence}

\ProveTheorem{quasibasis of module is basis}

\ShowDefinition{free module over ring}

\ShowTheorem{vector space is free module}
\ShowProof{basis over division algebra}

\ShowTheorem{coordinates of vector}
\ShowProof{coordinates of vector of free module}
}

\AddEq{fa=+0}
{
\[
f(a)=f(a+0)=f(a)+f(0)=f(a)+\Vector 0
\]
}

\AddEq{End +A}
{
$\End(\{+\},A)$.
}

\AddEq{ci=0 1n}
{
$\ACol ci=0$, $\gii=\gi 1$, ..., $\gin$.
}

\AddEq[2]{homomorphism D algebra 11}
{
(
\ShowEq{f: A->B}h{D_1}{D_2}\ \ \,
\ShowEq{f: A->B}{#1}{#2_1}{#2_2}{}
)
}

\AddEq[2]{homomorphism D algebra 1}
{
\ShowEq{f: A->B}{#1}{#2_1}{#2_2}{}
}

\AddEq[2]{homomorphism A module 111}
{
(
\ShowEq{f: A->B}h{D_1}{D_2}{}\ \ \,
\ShowEq{f: A->B}{#1}{A_1}{A_2}{}\ \ \,
\ShowEq{f: A->B}{#2}{V_1}{V_2}{}
)
}

\AddEq[2]{homomorphism A module 11}
{
(
\ShowEq{f: A->B}{#1}{A_1}{A_2}{}\ \ \,
\ShowEq{f: A->B}{#2}{V_1}{V_2}{}
)
}

\AddEq[2]{homomorphism A module 1}
{
\ShowEq{f: A->B}{#2}{V_1}{V_2}{}
}

\AddEq[2]{homomorphism- A module 111}
{
(
\ShowEq{f: A->B}{h^{-1}}{D_2}{D_1}{}\ \ \,
\ShowEq{f: A->B}{#1^{-1}}{A_2}{A_1}{}\ \ \,
\ShowEq{f: A->B}{#2^{-1}}{V_2}{V_1}{}
)
}

\AddEq[2]{homomorphism- A module 11}
{
(
\ShowEq{f: A->B}{#1^{-1}}{A_2}{A_1}{}\ \ \,
\ShowEq{f: A->B}{#2^{-1}}{V_2}{V_1}{}
)
}

\AddEq[2]{homomorphism- A module 1}
{
\ShowEq{f: A->B}{#2^{-1}}{V_2}{V_1}{}
}

\AddEq{ref linear map of D module 1}
{
\eqRef{homomorphism, f(p+q)=}{D module},
\eqRef{homomorphism, f(pq)=}{D module}.
}

\AddEq[2]{ref homomorphism of D algebra(1)}
{
\eqRef{homomorphism, f(p+q)=}1,
\eqRef{homomorphism, f(pq)=}1,
\ShowEq{ref homomorphism of D algebra()}{#1}{#2}
}

\AddEq[2]{ref homomorphism of D algebra()}
{
\eqRef{homomorphism, f v+w=}{(#1#2)g},
\eqRef{left homomorphism, f av=}{(#1#2)g}
}

\AddEq[3]{ref homomorphism of A vector space(1)}
{
\eqRef{homomorphism, f(p+q)=}{\SideWS A module},
\eqRef{homomorphism, f(pq)=}{\SideWS A module},
\ShowEq{ref homomorphism of A vector space()}{#1}{#2}{#3}
}

\AddEq[3]{ref homomorphism of A vector space()}
{
\eqRef{homomorphism, f v+w=}{g(#1#2)\SideWS A module},
\eqRef{homomorphism, f vw=}{(#1#2)\SideWS A module},
\eqRef{\SideWS homomorphism, f av=}{(#1#2)}
}

\AddEq[3]{define homomorphism of D algebra 11}
{
\DrawEq[hpq]{homomorphism, f(p+q)=}1
\DrawEq[hpq]{homomorphism, f(pq)=}1
\ShowEq{define homomorphism of D algebra 1}{#1}{#2}{#3}
}

\AddEq[3]{define homomorphism of D algebra 1}
{
\DrawEq[fab]{homomorphism, f v+w=}{(#1#2)g}
\FrameEqRef[fab]{homomorphism, f vw=}{(#1#2)g}
\DrawEq[fp{#3}a]{left homomorphism, f av=}{(#1#2)g}
\ShowEq{ap AD #1#2}
}

\AddEq[5]{define homomorphism of module 111}
{
\DrawEq[fpq]{homomorphism, f(p+q)=}{\SideNS}
\DrawEq[fpq]{homomorphism, f(pq)=}{\SideNS}
\ShowEq{define homomorphism of module 11}{#1}{#2}{#3}{#4}{#5}
}

\AddEq[5]{define homomorphism of module 11}
{
\DrawEq[gab]{homomorphism, f v+w=}{({#1}{#2}{#3})g \SideNS}
\DrawEq[gab]{homomorphism, f vw=}{({#1}{#2}{#3})g \SideNS}
\DrawEq[gp{#4}a]{left homomorphism, f av=}{({#1}{#2}{#3})g \SideNS}
\ShowEq{define homomorphism of module 1}{#1}{#2}{#3}{#4}{#5}
}

\AddEq[5]{define homomorphism of module 1}
{
\DrawEq[hvw]{homomorphism, f v+w=}{({#1}{#2}{#3})h \SideNS}
\DrawEq[ha{#5}v]{\SideWS homomorphism, f av=}{({#1}{#2}{#3})h \SideNS}
\ShowEq{vap VAD #1#2#3}
}

\AddEq[5]{define homomorphism of A module 111}
{
\DrawEq[hpq]{homomorphism, f(p+q)=}{\SideWS A module}
\DrawEq[hpq]{homomorphism, f(pq)=}{\SideWS A module}
\ShowEq{define homomorphism of A module 11}{#1}{#2}{#3}{#4}{#5}
}

\AddEq[5]{define homomorphism of A module 11}
{
\DrawEq[gab]{homomorphism, f v+w=}{g(#1#2)\SideWS A module}
\DrawEq[gab]{homomorphism, f vw=}{(#1#2)\SideWS A module}
\DrawEq[gp{#4}a]{\SideWS homomorphism, f av=}{(#1#2)}
\ShowEq{define homomorphism of A module 1}{#1}{#2}{#3}{#4}{#5}
}

\AddEq[5]{define homomorphism of A module 1}
{
\DrawEq[fuv]{homomorphism, f v+w=}{f(#1#2)\SideWS A module}
\DrawEq[fa{#5}v]{\SideWS homomorphism, f av=}{(#1#2)\SideWS A module}
\ShowEq{vap VAD #1#2#3}
}

\AddEq[3]{define homomorphism of D module 11}
{
\DrawEq[hpq]{homomorphism, f(p+q)=}{D module}
\DrawEq[hpq]{homomorphism, f(pq)=}{D module}
\ShowEq{define homomorphism of D module 1}{#1}{#2}{#3}
}

\AddEq[3]{define homomorphism of D module 1}
{
\DrawEq[fab]{homomorphism, f v+w=}{f D module #1#2}
\DrawEq[fp{#3}a]{left homomorphism, f av=}{f D module #1#2}
\ShowEq{vp VD #1#2}
}

\AddEq{morphism of representation, algebra, 23}
{
\BlueText{f\circ (h_{1.23}(a)(b))}
=h_{2.23}(\RedText{f\circ a})
(\BlueText{f\circ b})
}

\AddEq{algebra, morphism of representation 23, 1}
{
f\circ (C_1\circ (a,b))=C_2\circ (f\circ a,f\circ b)
}

\DefText[7]{homomorphism of A module 2}
{
\,
\ShowFootnote{homomorphism of A module}{#1}{#2}{#3}{#4}{#5}{#6}
\ShowTheorem{homomorphism A module}{#1}{#2}{#3}{#4}{#5}{#6}
\AddIndex{}{matrix of homomorphism}
\AddIndex{}{coordinates of homomorphism}
\ShowProof{homomorphism A module}{#1}{#2}{#3}{#4}{#5}{#6}

\ShowTheorem{matrix generates A module homomorphism}{#1}{#2}{#3}{#4}{#5}{#6}
\ShowProof{matrix generates A module homomorphism}{#1}{#2}{#3}{#4}{#5}{#6}{#7}

\ShowText{homomorphism of vector space 2}{#1}{#2}{#3}
}

\DefText[7]{homomorphism of A module (1)}
{
\ShowEq{\DefCol}
\ShowText{homomorphism of A module 1}{#1}{#2}{#3}{#4}{#5}{#6}{#7}

\ShowEq{\DefRow}
\ShowText{homomorphism of A module 1}{#1}{#2}{#3}{#4}{#5}{#6}{#7}
}

\DefText[7]{homomorphism of A module ()}
{
\ShowEq{\DefCol}
\ShowFootnote{homomorphism of A module}{#1}{#2}{#3}{#4}{#5}{#6}
\ShowEq{\DefRow}
\ShowFootnote{homomorphism of A module}{#1}{#2}{#3}{#4}{#5}{#6}

\TwoColText
{
\ShowEq{\DefCol}
\ShowTheorem{homomorphism A module}{#1}{#2}{#3}{#4}{#5}{#6}
}
{
\ShowEq{\DefRow}
\ShowTheorem{homomorphism A module}{#1}{#2}{#3}{#4}{#5}{#6}
}
\AddIndex{}{matrix of homomorphism}
\AddIndex{}{coordinates of homomorphism}
\TwoColText
{
\ShowEq{\DefCol}
\ShowProof{homomorphism A module}{#1}{#2}{#3}{#4}{#5}{#6}
}
{
\ShowEq{\DefRow}
\ShowProof{homomorphism A module}{#1}{#2}{#3}{#4}{#5}{#6}
}

\TwoColText
{
\ShowEq{\DefCol}
\ShowTheorem{matrix generates A module homomorphism}{#1}{#2}{#3}{#4}{#5}{#6}
\ShowProof{matrix generates A module homomorphism}{#1}{#2}{#3}{#4}{#5}{#6}{#7}
}
{
\ShowEq{\DefRow}
\ShowTheorem{matrix generates A module homomorphism}{#1}{#2}{#3}{#4}{#5}{#6}
\ShowProof{matrix generates A module homomorphism}{#1}{#2}{#3}{#4}{#5}{#6}{#7}
}

\TwoColText
{
\ShowEq{\DefCol}
\ShowText{homomorphism of vector space 2}{#1}{#2}{#3}
}
{
\ShowEq{\DefRow}
\ShowText{homomorphism of vector space 2}{#1}{#2}{#3}
}
}

\DefText[8]{def homomorphism A module}
{
\ShowText{algebra over ring (#1#2)}%
\ShowText{module over algebra (#2#3)}%

\ShowDefinition{homomorphism A module}{#1}{#2}{#3}{#4}{#5}{#6}

\ShowText{notation for homomorphism of module}

\ShowTheorem{define homomorphism A module}{#1}{#2}{#3}{#4}{#5}{#6}{#7}{#8}
\ShowProof{define homomorphism A module}{#1}{#2}{#3}{#4}{#5}{#6}

\ShowFootnote{iso end aut morphism}{\SideNS(#1#2#3)}

\ShowDefinition{isomorphism module}{#1}{#2}{#3}{#4}{#5}{#6}

\def\TempA{}%
\def\TempB{#2}%
\ifx\TempA\TempB%
\ShowDefinition{end aut morphism module}

\ShowTheorem{set of automorphisms}
\fi

\ShowText{homomorphism of A module (#1)}{#1}{#2}{#3}{#4}{#5}{#6}{#8}
}

\AddEq[1]{A1=A2 pdf}
{
\texorpdfstring{$#1_1=#1_2=#1$}{#1 1=#1 2=#1}
}

\DefText[6]{linear map of D module}
{
\,

\ShowDefinition{linear map of D module}{#1}{#2}{#3}{#4}

\ShowText{notation for linear map}{#1}{#2}{#3}{#4}

\ShowTheorem{linear map of D module}{#1}{#2}{#3}{#4}{#5}
\ShowProof{linear map of D module}{#1}{#2}

\ShowEq{\DefCol}
\ShowFootnote{homomorphism of D module}{#1}{#2}{#3}{#4}
\ShowEq{\DefRow}
\ShowFootnote{homomorphism of D module}{#1}{#2}{#3}{#4}

\TwoColText
{
\ShowEq{\DefCol}
\ShowTheorem{linear map of D module, coordinates}{#1}{#2}{#3}{#4}{#6}
}
{
\ShowEq{\DefRow}
\ShowTheorem{linear map of D module, coordinates}{#1}{#2}{#3}{#4}{#6}
}

\ShowEq{\DefCol}
\ShowProof{linear map of D module, coordinates}{#1}{#2}{#3}{#4}{#5}

\ShowEq{\DefRow}
\ShowProof{linear map of D module, coordinates}{#1}{#2}{#3}{#4}{#5}

\TwoColText
{
\ShowEq{\DefCol}
\ShowTheorem{matrix generates D module homomorphism}{#1}{#2}{#3}{#4}
}
{
\ShowEq{\DefRow}
\ShowTheorem{matrix generates D module homomorphism}{#1}{#2}{#3}{#4}
}

\TwoColText
{
\ShowEq{\DefCol}
\ShowProof{matrix generates D module homomorphism}{#1}{#2}{#3}{#4}{#5}
}
{
\ShowEq{\DefRow}
\ShowProof{matrix generates D module homomorphism}{#1}{#2}{#3}{#4}{#5}
}

\ShowText{we identify linear map and matrix}{#1}{#2}
}

\AddEq{linear map, 0, D algebra}
{
\[f\circ 0=0\]
}

\AddEq[4]{D->*o+A}
{
\[
\xymatrix
{
#1\ar[r]|(.4){*}&\displaystyle\bigoplus_{\iI}#2_i
}
\ \ \ \ \ \,
#3\left(\bigoplus_{\iI}#4_i\right)=\bigoplus_{\iI}#3 #4_i
\]
}

\AddEq[4]{left D->*o+A}
{
\[
\xymatrix
{
#1\ar[r]|(.4){*}&\displaystyle\bigoplus_{\iI}#2_i
}
\ \ \ \ \ \,
#3\left(\bigoplus_{\iI}#4_i\right)=\bigoplus_{\iI}#3 #4_i
\]
}

\AddEq[4]{right D->*o+A}
{
\[
\xymatrix
{
#1\ar[r]|(.4){*}&\displaystyle\bigoplus_{\iI}#2_i
}
\ \ \ \ \ \,
\left(\bigoplus_{\iI}#4_i\right)#3=\bigoplus_{\iI}#4_i #3
\]
}

\AddEq[1]{f(pai)=pfai}
{
#1\circ(
a_1, ...,
pa_i, ...,
a_n)
=
p#1\circ(
a_1, ...,
a_i, ...,
a_n)
}

\AddEq[4]{a=ai ei}
{
\ensuremath{#1=\ACol{#1}{#3}\EBase {#2}{#3}#4}
}

\AddEquation{sum of maps, 31, polylinear}
{
\begin{split}
&(f+g)\circ(x_1,...,x_i+y_i,...,x_n)
\\
=&\,f\circ(x_1,...,x_i+y_i,...,x_n)
+g\circ(x_1,...,x_i+y_i,...,x_n)
\\
=&\,f\circ (x_1,...,x_i,...,x_n)+f\circ (x_1,...,y_i,...,x_n)
\\
+&g\circ (x_1,...,x_i,...,x_n)+g\circ (x_1,...,y_i,...,x_n)
\\
=&(f+g)\circ (x_1,...,x_i,...,x_n)
+(f+g)\circ (x_1,...,y_i,...,x_n)
\end{split}
}

\AddEquation{sum of maps, 32, polylinear}
{
\begin{split}
&(f+g)\circ(x_1,...,px_i,...,x_n)
\\
=&f\circ(x_1,...,px_i,...,x_n)+g\circ(x_1,...,px_i,...,x_n)
\\
=&pf\circ (x_1,...,x_i,...,x_n)+pg\circ (x_1,...,x_i,...,x_n)
\\
=&p(f\circ (x_1,...,x_i,...,x_n)+g\circ (x_1,...,x_i,...,x_n))
\\
=&p(f+g)\circ (x_1,...,x_i,...,x_n)
\end{split}
}

\AddEq{module of polylinear maps, 1}
{
$f$, $g$, $h\in\mathcal L(D;A_1\times...\times A_2\rightarrow S)$.
}

\AddEq{module of polylinear maps, 2}
{
$a=(a_1,...,a_n)$, $a_1\in A_1$, ..., $a_n\in A_n$,
\begin{align*}
(f+g)\circ a=&f\circ a+g\circ a=g\circ a+f\circ a
\\
=&(g+f)\circ a
\\
((f+g)+h)\circ a=&
(f+g)\circ a+h\circ a=(f\circ a+g\circ a)+h\circ a
\\
=&f\circ a+(g\circ a+h\circ a)=f\circ a+(g+h)\circ a
\\
=&(f+(g+h))\circ a
\end{align*}
}

\AddEq{0:A1n->S}
{
\[
0:v\in A_1\times...\times A_n\rightarrow 0\in S
\]
}

\AddEq{0+f=f 1n}
{
\[
(0+f)\circ (a_1,...,a_n)=0\circ (a_1,...,a_n)+f\circ (a_1,...,a_n)=f\circ (a_1,...,a_n)
\]
}

\AddEq{f-f=0 1n}
{
\begin{align*}
(f+(-f))\circ (a_1,...,a_n)&=f\circ (a_1,...,a_n)+(-f)\circ (a_1,...,a_n)
\\&=f\circ (a_1,...,a_n)-f\circ (a_1,...,a_n)
\\&=0
=0\circ (a_1,...,a_n)
\end{align*}
}

\AddEq{sum of maps, 4, polylinear}
{
\begin{align*}
(f+g)\circ (a_1,...,a_n)
&=f\circ (a_1,...,a_n)+g\circ (a_1,...,a_n)\\
&=g\circ (a_1,...,a_n)+f\circ (a_1,...,a_n)\\
&=(g+f)\circ (a_1,...,a_n)
\end{align*}
}

\AddEq{f-f=0}
{
\[
f+(-f)=0
\]
}

\AddEq{-f:A1n->S}
{
\[
-f:(a_1,...,a_n)\in A_1\times...\times A_n\rightarrow -(f\circ(a_1,...,a_n))\in S
\]
}

\AddEq{(p+q)f=pf+qf 1}
{
\begin{align*}
((p+q)f)\circ (x_1,...,x_n)=&(p+q)(f\circ (x_1,...,x_n))
\\
=&p(f\circ (x_1,...,x_n))+q(f\circ (x_1,...,x_n))
\\
=&(pf)\circ (x_1,...,x_n)+(qf)\circ (x_1,...,x_n)
\end{align*}
}

\AddEq{(p+q)f=pf+qf}
{
(p+q)f=pf+qf
}

\AddEq{p(qf)=(pq)f}
{
p(qf)=(pq)f
}

\AddEq{p(qf)=(pq)f 1}
{
\begin{align*}
(p(qf))\circ (x_1,...,x_n)=&p\ (qf)\circ (x_1,...,x_n)
=p\ (q\ f\circ (x_1,...,x_n))
\\
=&(pq)\ f\circ (x_1,...,x_n)=((pq)f)\circ (x_1,...,x_n)
\end{align*}
}

\AddEquation{product of map over scalar, 1, D module}
{
\begin{matrix}
\ShowSymbol{product of map over scalar}{D module}:A_1\rightarrow A_2
&d\in D&f\in\mathcal L(D;A_1\rightarrow A_2)
\end{matrix}
}

\AddEq{endomorphism of module from product, 1}
{
\[
\begin{matrix}
\vcenter
{
\xymatrix
{
A\ar[r]|{*}^{g_{23}}&A
}
}
&
\begin{matrix}
g_{23}:v\rightarrow l\circ v
\\
g_{23}\circ v:w\rightarrow (l\circ v)\circ w
\end{matrix}
\end{matrix}
\]
}

\AddEquation{a:LA12->LA12}
{
a:f\in\mathcal L(D;A_1\rightarrow A_2)\rightarrow af\in\mathcal L(D;A_1\rightarrow A_2)
}

\AddEquation{product of map over scalar, D module}
{
(\ShowSymbol{product of map over scalar}{D module})\circ x=d(f\circ x)
}

\AddEq[2]{b=foa}
{
$#1=f\circ #2$,
}

\AddEq[1]{f(ai+bi)=fai+fbi}
{
#1\circ(
a_1, ...,
a_i+ b_i, ...,
a_n)
=
#1\circ(
a_1, ...,
a_i, ...,
a_n)
+
#1\circ(
a_1, ...,
b_i, ...,
a_n)
}

\AddEq[1]{w=sum vi 2015}
{
J(v)=\left\{w:w=\sum_{\iIg}\aU ci\aD vi,\aU ci\in #1,
|\{\gii:\aU ci\ne 0\}|<\infty\right\}
}

\AddEq{p,q in D, v,w in V}
{
$p$, $q \in \Base_{\BaseExt}$, $v$, $w \in V$.
}

\AddEq{m,n in D}
{
$m$, $n \in D$,
}

\AddEq[1]{o+Ai}
{
\[#1=\bigoplus_{\iI}#1_i\]
}

\AddEq{a=a1 o+ an}
{
\[a=a_1\oplus...\oplus a_n\]
}

\AddEq[2]{A1o+.n}
{
\[
#1=#1^1\oplus...\oplus #1^#2
\]
}

\AddEq[2]{A 1o+.n}
{
\[
#1=\aU{#1}1\oplus...\oplus \aU{#1}#2
\]
}

\AddEq[2]{a1o+.n}
{
\[
#1=#1_1\oplus...\oplus #1_#2
\]
}

\AddEq[1]{A1o+An=A1xAn}
{
\[
#1_1\oplus...\oplus #1_n=#1_1\times...\times #1_n
\]
}

\AddEq{rcd linear span, 0}
{
\[
\Vector a=\RowMatrix{\Vector a}n
=(\aD{\Vector a}i,\iIg)
\]
}

\AddEq{rcd linear span, 2}
{
\Vector b=\Vector a\RCstar x
}

\AddEq{rcd linear span, 6}
{
$\Vector b$, $\aD{\Vector a}i$
}

\AddEq{rcd linear span, 7}
{
\begin{align}
\EqLabel{rcd linear span, b}
\Vector b&=e\RCstar b\\
\EqLabel{rcd linear span, ai}
\aD{\Vector a}i&=e\RCstar\aD ai
\end{align}
}

\AddEq{rcd linear span, 8}
{
e\RCstar b=e\RCstar a\RCstar x
}

\AddEq{system of rcd linear equations}
{
a\RCstar x=b
}

\AddEq{inverce quaternion}
{
\[
x^{-1}=|x|^{-2}x^*
\]
}

\AddEq[2]{a=(a1.n set)}
{
#1=\AcMatrix{#1}{#2}
}

\AddEq[2]{a=(a1.n col)}
{
#1=\ColMatrix{#1}{#2}
}

\AddEq[2]{a=(a1.n row)}
{
#1=\RowMatrix{#1}{#2}
}

\AddEq{Expansion relative basis in algebra}
{
\[
a=a^{\gi i}e_{\gi i}
\]
}

\AddEq{av=sum av convention}
{
\[\aU ci\aD vi=\sum_{\gii\in\giI}\aU ci\aD vi\]
}

\AddEq{av=sum av()}
{
\[c(\gii)v(\gii)=\sum_{\gii\in\giI}c(\gii)v(\gii)\]
}

\AddEq{av=sum av()left}
{
\[c(\gii)v(\gii)=\sum_{\gii\in\giI}c(\gii)v(\gii)\]
}

\AddEq{av=sum av()right}
{
\[v(\gii)c(\gii)=\sum_{\gii\in\giI}v(\gii)c(\gii)\]
}

\AddEq[1]{foa=fi ai}
{
f\circ #1=
\RowMatrix fn
\RCcirc
\ColMatrix {#1}n
=\aD fi\circ\aU {#1}i
}

\AddEq{fioa=fij aj}
{
\aU fi\circ\VNumber=
\RowMatrix{\aU fi}n
\RCcirc
\ColMatrix{\VNumber}n
=\aUD fij\circ \aU{\VNumber}j
}

\AddEq{fij=(fi)j}
{
\[
\begin{pmatrix}
\RowMatrix{\aU f1}n
\\...\\
\RowMatrix{\aU fm}n
\end{pmatrix}
=
\PMatrix fnm
\]
}

\AddEq{vI=EI}
{
$\Vector I=(E,\aU I1,\aU I2,\aU I3)$
}

\AddEq[2]{foa=fi a}
{
\ColMatrix {#1}m
=
\ColMatrix fm
\circ #2=
\begin{pmatrix}
\aU f1\circ #2\\...\\ \aU fm\circ #2
\end{pmatrix}
}

\AddEq[3]{a=(a11.nm matrix)}
{
#1=\PMatrix {#1}{#2}{#3}
}

\AddEq[5]{b=f rco a}
{
\ColMatrix{#2}{#5}
=
\PMatrix {#1}{#4}{#5}
\RCcirc
\ColMatrix {#3}{#4}
=
\begin{pmatrix}
\aUD {#1}1i\circ\aU {#3}i\\
...\\
\aUD {#1}{#5}i\circ\aU {#3}i
\end{pmatrix}
}

\AddEq{fi:->}
{
\ShowEq{f:A->B}{\aU fi}{\Module}{\aU{\ModuleA}i}
}

\AddEquation{ft=tf1}
{
f(t)=tf(1)
}

\AddEq[2]{fij:->}
{
\ShowEq{f:A->B}{\aUD fij}{\aU {#1}j}{\aU {#2}i}
}

\AddEq{fij:-> left A}
{
\ShowEq{f:A->B}{f^{\gii}_{\gij}}{A_2(g\circ e_{V_1\cdot\gij})}{A_2e_{\gii}}
}

\AddEq{fij:-> right A}
{
\ShowEq{f:A->B}{f^{\gii}_{\gij}}{(g\circ e_{V_1\cdot\gij})A_2}{e_{\gii}A_2}
}

\DefProof{direct sum of n Abelian groups}
{
\TheoremFollows
\RefTheorem{direct sum of Abelian groups}.
}

\AddEq[2]{direct sum of modules}
{
\symb{#1\oplus #2}{direct sum}1
}

\AddEq[1]{c in ZA}
{
$#1\in Z(A)$,
}

\AddEq[1]{ca=ac}
{
c_{#1}a=ac_{#1}
}

\AddEquation{a(bc)=...=(bc)a}
{
a(bc)=(ab)c=(ba)c=b(ac)=b(ca)=(bc)a
}

\AddEq{c in ZAb}
{
$c_1$, $c_2\in Z(A,b)$.
}

\AddEq{c1 c2 in ZAb}
{
$c_1c_2\in Z(A,b)$.
}

\AddEquation{c1c2 in ZAb}
{
b(c_1c_2)=(bc_1)c_2=(c_1b)c_2=c_1(bc_2)=c_1(c_2b)=(c_1c_2)b
}

\AddEq[1]{(f+g)x=}
{
\[(f+g)(#1)=f(#1)+g(#1)\]
}

\AddEq[1]{(fp)x=}
{
\[(fp)(#1)=f(#1)p\]
}

\DefTheorem{cr-power, 1}
{
\ShowEq{cr-power, 1}
}

\AddEquation{cr-power, 1}
{
a^{\CRPower 1}=a
}

\AddEquation{cr-power, 1 1}
{
a^{\CRPower 1}=a^{\CRPower 0}\CRstar a=\UnitNMatrix\CRstar a=a
}

\DefTheorem{rc-power, 1}
{
\ShowEq{rc-power, 1}
}

\AddEquation{rc-power, 1}
{
a^{\RCPower 1}=a
}

\AddEquation{rc-power, 1 1}
{
a^{\RCPower 1}=a^{\RCPower 0}\RCstar a=\UnitNMatrix\RCstar a=a
}

\AddEq{center of A number}
{
\symb{Z(A,b)}{center of A number}1
}

\AddEq[1]{center Ac}
{
$Z(A,b)$#1
}

\AddEq[1]{c in S=>bc=cb}
{
c\in #1\Leftrightarrow cb=bc
}

\AddEq[3]{c in ZAa}
{
\ensuremath{#1\in Z(A,#2)}#3
}

\AddEq{p(c) p in D}
{
\begin{aligned}
p(x)&=p_0+p_1x+...+p_nx^n\\ &p_0,\ ...,\ p_n\in D 
\end{aligned}
}

\AddEq{c=p(b)}
{
\[
c=p(b)
\]
}

\AddEq{d in D c in S 12}
{
$d_1$, $d_2\in D$, $c_1$, $c_2\in Z(A,b)$.
}

\AddEq{dcb=bdc}
{
\begin{align*}
(d_1c_1+d_2c_2)b&=d_1c_1b+d_2c_2b=bd_1c_1+bd_2c_2\\&=b(d_1c_1+d_2c_2)
\end{align*}
}

\AddEq{dc in S}
{
$d_1c_1+d_2c_2\in Z(A,b)$.
}

\DefProof{direct sum of n D modules}
{
\TheoremFollows
\refTheorem{direct sum of D modules}{\SideWS \Base-\VectorsSetNS}.
}

\AddEq{v=veV}
{
$v=v^{\gi i}\EV V i$.
}

\AddEquation{linear map in left A module}
{
f\circ(v\CRstar e_V )=(f_{\gi i}^{\gi j}\circ g\circ v^{\gi i})\ \EV W j
}

\AddEquation{linear map in right A module}
{
f\circ(e_V\RCstar v)=\EV W j(f_{\gi i}^{\gi j}\circ g\circ v^{\gi i})\ 
}

\AddEq{linear map in A module}
{
w=f\RCcirc (g\circ v)
}

\AddEq [1]{Ai, i=}
{
$#1_i$, $i=1$, $2$,
}

\AddEquation{module, (a+n)v=av+nv}
{
(a+n)v=av+nv\ \ \ a\in D\ \ \ n\in Z
}

\AddEquation{left module, (a+n)v=av+nv}
{
(a+d)v=av+dv\ \ \ a\in A\ \ \ d\in D
}

\AddEquation{right module, (a+n)v=av+nv}
{
v(a+d)=va+vd\ \ \ a\in A\ \ \ d\in D
}

\AddEquation{left module, (a+d)*v=a*v+dv}
{
(a+d)*v=a*v+dv\ \ \ a\in A\ \ \ d\in D
}

\AddEquation{right module, (a+d)*v=a*v+dv}
{
v*(a+d)=v*a+vd\ \ \ a\in A\ \ \ d\in D
}

\AddEq{(a+n)+(b+m)=}
{
(a+n)+(b+m)=(a+b)+(n+m)
}

\AddEq[1]{(d1+0)+(d2+0)}
{
(#1_1+0)+(#1_2+0)=(#1_1+#1_2)+(0+0)=(#1_1+#1_2)+0
}

\AddEquation{module, (a+n)+(b+m)=}
{
\begin{split}
((a+n)+(b+m))v&=av+nv+bv+mv=av+bv+nv+mv\\
&=(a+b)v+(n+m)v=((a+b)+(n+m))v
\end{split}
}

\AddEquation{g=CG}
{
g^{\gii}_{\gik}=C^{\gii}_{\gij}G^{\gij}_{\gik}
}

\AddEq [3]{matrix amn}
{
\[
#1=
\begin{pmatrix}
\aUD{#1}11&...&\aUD{#1}1{#3}\\
...&...&...\\
\aUD{#1}{#2}1&...&\aUD{#1}{#2}{#3}
\end{pmatrix}
\]
}

\AddEquation{g=(cF)G}
{
g^{\gii}_{\gik}=(c_{1.k.l}F_{k\cdot}^{}{}^{\gii}_{\gij}c_{2.k.l})G^{\gij}_{\gik}
}

\AddEq{C=cF}
{
\[
C=c_{1.k.l}F_kc_{2.k.l}
\]
}

\AddEq{Gkj in D}
{
$G^{\gij}_{\gik}\in D$,
}

\AddEquation{g=c(FG)}
{
g^{\gii}_{\gik}=c_{1.k.l}(F_{k\cdot}^{}{}^{\gii}_{\gij}G^{\gij}_{\gik})c_{2.k.l}
}

\AddEquation{left module, (a+n)+(b+m)=}
{
\begin{split}
((a+n)+(b+m))v&=av+nv+bv+mv=av+bv+nv+mv\\
&=(a+b)v+(n+m)v=((a+b)+(n+m))v
\end{split}
}

\AddEquation{right module, (a+n)+(b+m)=}
{
\begin{split}
v((a+n)+(b+m))&=va+vn+vb+vm=va+vb+vn+vm\\
&=v(a+b)+v(n+m)=v((a+b)+(n+m))
\end{split}
}

\AddEquation{module, (a+n)+(b+m)=, 1}
{
((a+n)+(b+m))v=(a+n)v+(b+m)v
}

\AddEquation{left module, (a+n)+(b+m)=, 1}
{
((a+n)+(b+m))v=(a+n)v+(b+m)v
}

\AddEquation{right module, (a+n)+(b+m)=, 1}
{
v((a+n)+(b+m))=v(a+n)+v(b+m)
}

\AddEq{(a+n)(b+m)=}
{
(a+n)(b+m)=(ab+ma+nb)+(nm)
}

\AddEq{(a+n)*(b+m)=}
{
(a+n)*(b+m)=(a*b+ma+nb)+(nm)
}

\AddEq[1]{(d1+0)(d2+0)}
{
(#1_1+0)(#1_2+0)=(#1_1#1_2+0#1_1+0#1_2)+(0\,0)=(#1_1#1_2)+0
}

\AddEq{da in DA->da in D1A}
{
\[
(
I_{(1)}:d\in \Base\rightarrow d+0\in \Base_{(1)}\ \ \ \,
\id:v\in V\rightarrow v\in V
)
\]
}

\AddEquation{da in DA->da in D1A 1, module}
{
(d+0)v=dv+0v=dv+0=dv
}

\AddEquation{da in DA->da in D1A 1, left module}
{
(d+0)v=dv+0v=dv+0=dv
}

\AddEquation{da in DA->da in D1A 1, right module}
{
v(d+0)=vd+v0=vd+0=vd
}

\AddEq [1]{A module diagram for linear}
{
\[
\xymatrix
{
A_{#1}&&V_{#1}\\
&D_{#1}\ar[lu]|{*}\ar[ru]|{*}
}
\]
}

\DefEquation
{
f^k=f^{k\cdot\gi{ij}}e_{\gii}\otimes e_{\gij}
}
{standard representation of map A1 A2, 3, associative algebra}

\AddEquation{A=re+im}
{
A=\re A\oplus\im A
}

\AddEquation{standard representation of map A->A, associative algebra}
{
f=
f^{k\cdot\gi{ij}}(\aD ei\otimes \aD ej)\circ F_k
}

\AddEq{standard representation of map A->A, o, associative algebra}
{
\[
f\circ x=
f^{k\cdot\gi{ij}}\aU ei(F_k\circ x)\aU ej
\]
}

\AddEq{product of map over scalar,,D module}
{
\symb{d\,f}{product of map over scalar}{D module}
}

\AddEq{product of map over scalar,,polylinear}
{
\symb{d\,f}{product of map over scalar}{polylinear}
}

\AddEquation{product of map over scalar, 1, polylinear}
{
\begin{matrix}
\ShowSymbol{product of map over scalar}{polylinear}{}:
A_1\times...\times A_n\rightarrow S
&d\in D&f\in\mathcal L(D;A_1\times...\times A_n\rightarrow S)
\end{matrix}
}

\AddEquation{product of map over scalar, 31, polylinear}
{
\begin{split}
&\,(pf)\circ(x_1,...,x_i+y_i,...,x_n)
\\
=&\,p\ f\circ(x_1,...,x_i+y_i,...,x_n)
\\
=&\,p\ (f\circ (x_1,...,x_i,...,x_n)+f\circ (x_1,...,y_i,...,x_n))
\\
=&\,p(f\circ (x_1,...,x_i,...,x_n))+p(f\circ (x_1,...,y_i,...,x_n))
\\
=&\,(pf)\circ (x_1,...,x_i,...,x_n)+(pf)\circ (x_1,...,y_i,...,x_n)
\end{split}
}

\AddEquation{product of map over scalar, 32, polylinear}
{
\begin{split}
& \,(pf)\circ(x_1,...,qx_i,...,x_n)
\\=& \,p(f\circ(x_1,...,qx_i,...,x_n))
=pq(f\circ (x_1,...,x_i,...,x_n))
\\=& \,qp(f\circ (x_1,...,x_n))
=q(pf)\circ (x_1,...,x_n)
\end{split}
}

\AddEquation{product of map over scalar, polylinear}
{
(\ShowSymbol{product of map over scalar}{polylinear}{})\circ (a_1,...,a_n)
=d(f\circ (a_1,...,a_n))
}

\DefCorollary{Free Algebra over Ring}
{
\ShowText{Free Algebra over Ring}
}

\AddEquation{a:LA1n->LA1n}
{
a:f\in\mathcal L(D;A_1\times...\times A_n\rightarrow S)\rightarrow
af\in\mathcal L(D;A_1\times...\times A_n\rightarrow S)
}

\DefTheorem{Free Algebra over Ring}
{
\ShowText{Free Algebra over Ring}
}

\AddEq{D*V->V}
{
(a,v)\in D\times V\rightarrow av\in V
}

\AddEq{left A*V->V}
{
(a,v)\in A\times V\rightarrow av\in V
}

\AddEq{right A*V->V}
{
(v,a)\in V\times A\rightarrow va\in V
}

\AddEq{left A**V->V}
{
(a,v)\in A\times V\rightarrow a*v\in V
}

\AddEq{right A**V->V}
{
(v,a)\in V\times A\rightarrow v*a\in V
}

\AddEquation{l(v1+v2)w}
{
(l\circ (v_1+v_2))\circ w=(v_1+v_2)w=v_1w+v_2w
=(l\circ v_1)\circ w+(l\circ v_2)\circ w
}

\AddEquation{l(v1+v2)}
{
l\circ (v_1+v_2)=l\circ v_1+l\circ v_2
}

\AddEquation{l(dv)w}
{
(l\circ (dv))\circ w=(dv)w=d(vw)
=d((l\circ v)\circ w)
}

\AddEquation{l(dv)}
{
l\circ (dv)=d(l\circ v)
}

\AddEq{g23:A->L}
{
\[g_{23}:A\rightarrow\mathcal L(D;A\rightarrow A)\]
}

\AddEquation{coordinates of map A->A, 3, associative algebra}
{
\begin{split}
\aUD fkl\aU xl\aD ek
&=
f^{k\cdot\gi{ij}} \aD ei
\aUD {F_{k\cdot}^{}{}}ml\aU xl\aD em
\aD ej
\\
&=
f^{k\cdot\gi{ij}}
\aUD {F_{k\cdot}^{}{}}ml\aU xl
\aUD Cp{im}\aUD Ck{pj}
\aD ek
\end{split}
}

\DefEq
{
f=f^k\circ F_k
}
{map f generated by basis F}

\AddEquation{fk= in A2xA2}
{
f^k=f^k_{s_k\cdot 0}\otimes f^k_{s_k\cdot 1}\ \ \ \ f^k\in A_2\otimes A_2
}

\AddEq [3]{n=dim A}
{
$\gi #1=\dim #2$#3
}

\AddEq{gi n<m}
{
$\gi n\le\gi m$.
}

\AddEq{gi n>m}
{
$\gi n>\gi m$.
}

\AddEquation{set ker G in ker g}
{
\{g\in\mathcal L(D;A\rightarrow B):\ker G\subseteq\ker g\}
}

\AddEquation{fk= in AxA}
{
f^k=f^k_{s_k\cdot 0}\otimes f^k_{s_k\cdot 1}\ \ \ \ f^k\in A\otimes A
}

\AddEquation{fk= in (Ax)A}
{
f^k=(f^k_{s_k\cdot 0}\otimes) f^k_{s_k\cdot 1}\ \ \ \ f^k\in A\otimes A
}

\DefEquation
{
f=f_{s\cdot 0}\otimes f_{s\cdot 1}
}
{expansion of f A->A}

\DefEquation
{
f=f^k\circ F_k
}
{expansion of f with respect to basis I}

\DefEquation
{
g=g_{t\cdot 0}\otimes g_{t\cdot 1}
}
{expansion of g A->A}

\DefEquation
{
g=g^l\circ G_l
}
{expansion of g with respect to basis J}

\DefEquation
{
h=h_{ts\cdot 0}\otimes h_{ts\cdot 1}
}
{expansion of h A->A}

\DefEquation
{
h=h^{lk}\circ K_{lk}
}
{expansion of h with respect to basis K}

\DefEquation
{
h\circ a=g\circ f\circ a
=g^l\circ J_l\circ f^k\circ I_k\circ a 
}
{h(a)1}

\DefEquation
{
\begin{split}
h\circ a=g\circ f\circ a
&=g^l\circ(J^m_l\circ f^k)\circ
J_m\circ I_k\circ a
\\&=g^l\circ(J^m_l\circ f^k)\circ K_{mk}\circ a
\end{split}
}
{h(a)2}

\DefEquation
{
\begin{split}
h_{ts\cdot 0}&=g_{t\cdot 0}f_{s\cdot 0}
\\
h_{ts\cdot 1}&=f_{s\cdot 1}g_{t\cdot 1}
\end{split}
}
{hlk=glfk 01}

\DefEq
{
\[a^{lk}K_{lk}=(a^{lk}\circ J_l)\circ I_k=0\]
}
{aK=aJI}

\DefEquation
{
\begin{split}
h\circ a&=g\circ f\circ a
\\&=(g_{t\cdot 0}\otimes g_{t\cdot 1})\circ(f_{s\cdot 0}\otimes f_{s\cdot 1})\circ a
\\&=(g_{t\cdot 0}\otimes g_{t\cdot 1})\circ(f_{s\cdot 0}a f_{s\cdot 1})
\\&=g_{t\cdot 0}f_{s\cdot 0}a f_{s\cdot 1} g_{t\cdot 1}
\end{split}
}
{h(a)}

\DefEq
{
\[a^{lk}\circ J_l=0\]
}
{aJ=0}

\DefEquation
{
h^{lk}=g^l\circ (G^k_m\circ f^m)
}
{hlk=glfk}

\DefEq
{
\Basis H=\{H_{lk}:H_{lk}=G_l\circ F_k, G_l\in\Basis G, F_k\in\Basis F\}
}
{JlIk}

\DefEq
{
h=g\circ f
}
{h=g o f}

\AddEq[2]{h:AoxA->L(A)}
{
\xymatrix
{
h:#2\otimes #2\ar[r]|{*}&\mathcal L(#1;#2\rightarrow #2)
}
\ \ \ h(p):g\rightarrow p\circ g
}

\AddEq[2]{representation AA in LA}
{
\begin{matrix}
(a\otimes b)\circ g=agb
&a,b\in #2&g\in\mathcal L(#1;#2\rightarrow #2)
\end{matrix}
}

\ePrints{5114-6019}
\ifx\Semafor\ValueOff
\AddEq{component of linear map}
{
\symb{f^k_{s_k\cdot p}}{component of linear map, division ring}1,
$p=0$, $1$,
}

\AddEq{standard component of linear map}
{
\symb{f^{k\cdot\gi{ij}}}{standard component of linear map}1
}
\fi

\AddEq[2]{product in algebra AA}
{
(#1_0\otimes #1_1)\circ(#2_0\otimes #2_1)=(#1_0 #2_0)\otimes(#2_1 #1_1)
}

\AddEquation{representation A2 in LA}
{
\xymatrix
{
h:\ATwo\ar[r]|(.4){*}&\mathcal L(D;A_1\rightarrow A_2)
}
}

\AddEquation{representation A2 in LA, 1}
{
\begin{matrix}
(a\otimes b)\circ f=afb
&a,b\in A_2&f\in\mathcal L(D;A_1\rightarrow A_2)
\end{matrix}
}

\AddEq{xA1n*xA1n->oxA1n=}
{
\[
(a_1,...,a_n)*(b_1,...,b_n)=(a_1b_1)\otimes...\otimes(a_nb_n)
\]
}

\DefEquation
{
(a\otimes b)\circ c=acb
}
{a ox b c=}

\AddEq[1]{a ox b}
{
$(a\otimes b)\circ #1$
}

\AddEq[1]{d in AxoA}
{
$d\in\AxA{#1}$
}

\AddEq[4]{product in algebra AA 2}
{
$d\circ #3\in\mathcal L(#1;#2\rightarrow #2)$#4
}

\DefEq
{
$A_1$, ..., $A_n$
}
{A1...n}

\AddEq{diagram of representations, module}
{
\xymatrix
{
D\ar[r]|{*}^{g_2}&V\\
&Z\ar[u]|{*}_{g_1}
}
}

\AddEq[8]{diagram of representations, left right module}
{
\begin{matrix}
\vcenter
{
\xymatrix
{
A_{#6}\ar[r]|{*}^{g_{#8 23}}&
A_{#6}\ar[r]|{*}^{g_{#8 34}}&V_{#7}
\\
&D_{#5}\ar[u]|{*}_{g_{#8 12}}\ar@/_1pc/[ur]|{*}_{g_{#8 14}}\ar[ul]|(.4){*}^{g_{#8 12}}
}
}
&
\begin{array}{r@{\,}l}
g_{12}(d):a&\rightarrow d\,a
\\
g_{23}(v):w&\rightarrow
C(w, v)
\\
C&\in\mathcal L(A^2\rightarrow A)
\\
g_{34}(a):v&\rightarrow #3\,#4
\\
g_{14}(d):v&\rightarrow #1\,#2
\end{array}
\end{matrix}
}

\AddEq[8]{diagram of representations, left right nonA module}
{
\begin{matrix}
\vcenter
{
\xymatrix
{
A_{#6}\ar[r]|{*}^{g_{#8 23}}&
A_{#6}\ar[r]|{*}^{g_{#8 34}}&V_{#7}
\\
&D_{#5}\ar[u]|{*}_{g_{#8 12}}\ar@/_1pc/[ur]|{*}_{g_{#8 14}}\ar[ul]|(.4){*}^{g_{#8 12}}
}
}
&
\begin{array}{r@{\,}l}
g_{12}(d):a&\rightarrow d\,a
\\
g_{23}(v):w&\rightarrow
C(w, v)
\\
C&\in\mathcal L(A^2\rightarrow A)
\\
g_{34}(a):v&\rightarrow #3*#4
\\
g_{14}(d):v&\rightarrow #1\,#2
\end{array}
\end{matrix}
}

\AddEq{diagram of representations, left module nonassociative}
{
\[
\begin{matrix}
\vcenter
{
\xymatrix
{
A\ar[rd]|{*}^{g_{23}}\ar[rr]|{*}^{g_{34}}&&V
\\
&A
\\
&D\ar[u]|{*}_{g_{12}}\ar@/_1pc/[uur]|{*}_{g_{14}}\ar[uul]|(.4){*}^{g_{12}}
}
}
&
\begin{array}{r@{\,}l}
g_{12}(d):a&\rightarrow d\,a
\\
g_{23}(v):w&\rightarrow
C(w, v)
\\
C&\in\mathcal L(A^2\rightarrow A)
\\
g_{34}(a):v&\rightarrow v\,a
\\
g_{14}(d):v&\rightarrow d\,v
\end{array}
\end{matrix}
\]
}

\AddEq[4]{diagram of representations, left module}
{
\ShowEq{diagram of representations, left right module}dvav{#1}{#2}{#3}{#4}
}

\AddEq[4]{diagram of representations, right module}
{
\ShowEq{diagram of representations, left right module}vdva{#1}{#2}{#3}{#4}
}

\AddEq[4]{diagram of representations, left nonA module}
{
\ShowEq{diagram of representations, left right nonA module}dvav{#1}{#2}{#3}{#4}
}

\AddEq[4]{diagram of representations, right nonA module}
{
\ShowEq{diagram of representations, left right nonA module}vdva{#1}{#2}{#3}{#4}
}

\AddEq{A1=A unital extension}
{
&$\CBase\subseteq\Base$&$\Base_{(1)}=\Base$\\
\hline
}

\AddEq{A1=D unital extension}
{
&$\Base\subseteq\CBase$&$\Base_{(1)}=\CBase$\\
\hline
}

\AddEq{A1=A+D unital extension}
{
&\multicolumn{2}{c|}{$\Base_{(1)}=\Base\oplus\CBase$}\\
\hline
}

\AddEq{set of vectors generated by vector A}
{
\symb[A-0]{A_{(1)}v}{set of vectors generated by vector}1
}

\AddEq{set of vectors generated by vector D}
{
\symb[D-0]{D_{(1)}v}{set of vectors generated by vector}1
}

\AddEq{D->*V}
{
\begin{matrix}
\xymatrix
{
f:D\ar[r]|{*}&V
}
&
f(d):v\rightarrow d\,v
\end{matrix}
}

\AddEq{right A->*V}
{
\begin{matrix}
\xymatrix
{
g_{34}:A\ar[r]|{*}&V
}
&
g_{34}(a):v\in V\rightarrow va\in V
&a\in A
\end{matrix}
}

\AddEq{left A->*V}
{
\begin{matrix}
\xymatrix
{
g_{34}:A\ar[r]|{*}&V
}
&
g_{34}(a):v\in V\rightarrow av\in V
&a\in A
\end{matrix}
}

\AddEq{right A*->*V}
{
\begin{matrix}
\xymatrix
{
g_{34}:A\ar[r]|{*}&V
}
&
g_{34}(a):v\in V\rightarrow v*a\in V
&a\in A
\end{matrix}
}

\AddEq{left A*->*V}
{
\begin{matrix}
\xymatrix
{
g_{34}:A\ar[r]|{*}&V
}
&
g_{34}(a):v\in V\rightarrow a*v\in V
&a\in A
\end{matrix}
}

\AddEq{module, associative law}
{
(ab)v=a(bv)
}

\AddEquation{module, (a+n)(b+m)=}
{
\begin{split}
((a+n)(b+m))v&=(a+n)((b+m)v)=(a+n)(bv+mv)\\
&=a(bv+mv)+n(bv+mv)\\
&=a(bv)+a(mv)+n(bv)+n(mv)\\
&=(ab)v+m(av)++n(bv)+(nm)v\\
&=(ab)v+(ma)v++(nb)v+(nm)v\\
&=((ab+ma+nb)+nm)v
\end{split}
}

\AddEq{fab=ab}
{
f(a,b)=ab\ \ \ a,b\in\Base
}

\AddEq{f1p=p}
{
f(1,p)=f(p,1)=p\ \ \ p\in\Base_{(1)}\ \ \ 1\in\CBase_{(1)}
}

\AddEq{pq=fpq}
{
\begin{split}
(a+n)(b+m)&=f(a+n,b+m)\\
&=f(a,b)+f(a,m)+f(n,b)+f(n,m)\\
&=f(a,b)+mf(a,1)+nf(1,b)+nf(1,m)\\
&=ab+ma+nb+nm
\end{split}
}

\AddEq{f:D1xD1->D1}
{
\[
f:\Base_{(1)}\times\Base_{(1)}\rightarrow\Base_{(1)}
\]
}

\AddEq{distributive law, module, 1}
{
p(v+w)=pv+pw
}

\AddEq{distributive law, module, 2}
{
(p+q)v=pv+qv
}

\AddEq{p,q in D1}
{
$p=a+n\in\Base_{(1)}$, $q=b+m\in\Base_{(1)}$.
}

\AddEq [4]{h:D1->D2 f:A1->A2}
{
\[
\begin{pmatrix}
#1:#2_1\rightarrow #2_2&
#3:#4_1\rightarrow #4_2
\end{pmatrix}
\]
}

\AddEq{set linear maps, D12 module}
{
\symb{\mathcal L(\Base_1\rightarrow \Base_2;\Module_1\rightarrow \Module_2)}{set linear maps}1
}

\AddEq{set linear maps, VW module}
{
\symb{\mathcal L(A;V\rightarrow W)}{set linear maps}1
}

\AddEq{set linear maps, left A12 module}
{
\symb{\mathcal L(D_1\rightarrow D_2;A_1\CRstar\rightarrow A_2\CRstar;V_1\rightarrow V_2)}{set linear maps}1
}

\AddEq{set homomorphisms, D algebra 11}
{
\symb{\Hom(D_1\rightarrow D_2;A_1\rightarrow A_2)}{set homomorphisms}1
}

\AddEq{set homomorphisms, D algebra 1}
{
\symb{\Hom(D;A_1\rightarrow A_2)}{set homomorphisms}1
}

\AddEq{set homomorphisms, A module 111 left}
{
\symb{\Hom(D_1\rightarrow D_2;A_1\rightarrow A_2;*V_1\rightarrow *V_2)}{set homomorphisms}1
}

\AddEq{set homomorphisms, A module 11 left}
{
\symb{\Hom(D;A_1\rightarrow A_2;*V_1\rightarrow *V_2)}{set homomorphisms}1
}

\AddEq{set homomorphisms, A module 1 left}
{
\symb{\Hom(D;A;*V_1\rightarrow *V_2)}{set homomorphisms}1
}

\AddEq{set homomorphisms, A module 111 right}
{
\symb{\Hom(D_1\rightarrow D_2;A_1\rightarrow A_2;V_1*\rightarrow V_2*)}{set homomorphisms}1
}

\AddEq{set homomorphisms, A module 11 right}
{
\symb{\Hom(D;A_1\rightarrow A_2;V_1*\rightarrow V_2*)}{set homomorphisms}1
}

\AddEq{set homomorphisms, A module 1 right}
{
\symb{\Hom(D;A;V_1*\rightarrow V_2*)}{set homomorphisms}1
}

\AddEq{set homomorphisms, left A12 module}
{
\symb{\Hom(D_1\rightarrow D_2;A_1\CRstar\rightarrow A_2\CRstar;V_1\rightarrow V_2)}{set homomorphisms}1
}

\AddEq{set linear maps, right A12 module}
{
\symb{\mathcal L(D_1\rightarrow D_2;\RCstar A_1\rightarrow\RCstar A_2;V_1\rightarrow V_2)}{set linear maps}1
}

\AddEq{set homomorphisms, right A12 module}
{
\symb{\Hom(D_1\rightarrow D_2;\RCstar A_1\rightarrow\RCstar A_2;V_1\rightarrow V_2)}{set homomorphisms}1
}

\AddEq [1]{left column module}
{
$#1\CRstar$\Hyph
}

\AddEq [1]{right column module}
{
$\RCstar #1$\Hyph
}

\AddEquation{unitarity law, D-module}
{
1v=v
}

\AddEq{associative law, module}
{
(pq)v=p(qv)
}

\AddEq{associative law, left module}
{
(pq)v=p(qv)
}

\AddEq{associative law, right module}
{
v(pq)=(vp)q
}

\AddEq{associative law*, left module}
{
(mn)v=m(nv)
}

\AddEq{associative law*, right module}
{
v(mn)=(vm)n
}

\AddEq{left module, associative law}
{
(ab)v=a(bv)
}

\AddEq{associative law, Z, module}
{
(mn)v=m(nv)
}

\AddEq{associative law, ZD, module}
{
(md)v=m(dv)
}

\AddEq{associative law, D, left module}
{
(cd)v=c(dv)
}

\AddEq{associative law, DA, left module}
{
(da)v=d(av)
}

\AddEq{distributive law, left module, 1}
{
p(v+w)=pv+pw
}

\AddEq{distributive law, left module, 2}
{
(p+q)v=pv+qv
}

\AddEq{distributive law*, left module, 1}
{
p*(v+w)=p*v+p*w
}

\AddEq{distributive law*, left module, 2}
{
(p+q)*v=p*v+q*v
}

\AddEquation{distributive law, Z, module, 2}
{
(n+m)v=nv+mv
}

\AddEquation{distributive law, D, module, 2}
{
(a+b)v=av+bv
}

\AddEquation{distributive law, D, left module, 2}
{
(n+m)v=cn+dm
}

\AddEquation{distributive law, A, left module, 2}
{
(a+b)v=av+bv
}

\AddEquation{distributive law, D, right module, 2}
{
v(n+m)=vn+vm
}

\AddEquation{distributive law, A, right module, 2}
{
v(a+b)=va+vb
}

\AddEq{distributive law, right module, 1}
{
(v+w)p=vp+wp
}

\AddEq{distributive law, right module, 2}
{
v(p+q)=vp+vq
}

\AddEq{distributive law*, right module, 1}
{
(v+w)*p=v*p+w*p
}

\AddEq{distributive law*, right module, 2}
{
v*(p+q)=v*p+v*q
}

\AddEquation{aw in Xk module}
{
aw\in X_{k+1}
}

\AddEquation{aw in Xk left module}
{
aw\in X_{k+1}
}

\AddEquation{aw in Xk right module}
{
wa\in X_{k+1}
}

\AddEq[2]{w1+w2 in Jv}
{
$#1_1+#1_2\in J(#2)$.
}

\AddEquation{v=v+0, module}
{
\Vector v=\Vector v+0=
\ACol vi
\ARow ei+
\ACol ci
\ARow ei=
(\ACol vi+
\ACol ci)
\ARow ei
}

\AddEquation{v=v+0, left module}
{
\Vector v=\Vector v+0=
\ACol vi
\ARow ei+
\ACol ci
\ARow ei=
(\ACol vi+
\ACol ci)
\ARow ei
}

\AddEquation{v=v+0, right module}
{
\Vector v=\Vector v+0=
\ARow ei
\ACol vi+
\ARow ei
\ACol ci=
\ARow ei
(\ACol vi+
\ACol ci)
}

\AddEq{px x-i x-j}
{
\[
p(x)=2(x-i)(x-j)+(x-j)(x-i)
=3x^2-2ix-2xj-jx-xi+k
\]
}

\AddEq{px x-i x-j 1}
{
\[
p(x)=((3-a)(x\otimes 1)+a\otimes x)\circ x-2ix-2xj
-jx-xi+k
\]
}

\AddEq{rx x-i}
{
r(x)=x-i
}

\AddEq{p(x)=sr 2}
{
\begin{align*}
p(x)&=(3-a)(r(x)+i)x+ax(r(x)+i)-2ix-2xj-jx-xi+k
\\
&=(3-a)r(x)x+(3-a)ix+axr(x)+axi-2ix-2xj-jx-xi+k
\\
&=(3-a)r(x)x+axr(x)+(3-a)i(r(x)+i)
\\
&+a(r(x)+i)i-2i(r(x)+i)-2(r(x)+i)j-j(r(x)+i)-(r(x)+i)i+k\\
&=(3-a)r(x)x+axr(x)+3ir(x)-3-air(x)+a
\\
&+ar(x)i-a-2ir(x)+2-2r(x)j-2k-jr(x)+k-r(x)i+1+k
\\
&=(3-a)r(x)x+axr(x)+3ir(x)-air(x)
\\
&+ar(x)i-2ir(x)-2r(x)j-jr(x)-r(x)i
\end{align*}
}

\AddEq{p=s2r}
{
\[p(x)=s_2(a,x)\circ r(x)\]
}

\AddEq{s2(x)=}
{
\begin{align*}
s_2(a,x)&=(3-a)\otimes x+ax\otimes 1+(3-a)i\otimes 1
\\
&+a\otimes i-2i\otimes 1-2\otimes j-j\otimes 1-1\otimes i
\\
&=(3-a)\otimes x+(ax-j-2i+(3-a)i)\otimes 1
+(a-1)\otimes i-2\otimes j
\\
&=(3-a)\otimes x+(ax-j+i-ai)\otimes 1
+(a-1)\otimes i-2\otimes j
\end{align*}
}

\AddEq{s2(1,x)=}
{
\begin{align*}
s_2(1,x)&=2\otimes x+(x-j-2i+2i)\otimes 1
+(1-1)\otimes i-2\otimes j
\\
&=2\otimes (x-j)+(x-j)\otimes 1
\end{align*}
}

\AddEq{s12 in AA}
{
$s_1(x)$, $s_2(x)\in A\otimes A$
}

\AddEq[1]{px=sx o rx}
{
\[
p(x)=s_{#1}(x)\circ r(x)
\]
}

\AddEq{(sx-sx)rx=0}
{
\[
(s_1(x)-s_2(x))\circ r(x)=0
\]
}

\AddEquation{x p0 x*+xQ+Rx*-S=0 3}
{
\begin{split}
p(x)&=(-4\otimes 1\otimes 1-4\otimes i\otimes i-4\otimes j\otimes j
-4\otimes k\otimes k)\circ x^2
\\&+(4\otimes i+2j\otimes 1-2k\otimes i-2\otimes j+2i\otimes k)\circ x+1+2k
\\&=-4x^2-4xixi-4xjxj
-4xkxk 
\\&+4xi+2jx-2kxi-2xj+2ixk+1+2k
\end{split}
}

\AddEquation{q(x)=}
{
q(x)=x-\frac{1}{4}(i+j)
}

\AddEquation{4x=4q+}
{
4x=4q(x)+i+j
}

\AddEquation{rj j-i}
{
r(j)=j-i
}

\AddEq{p(x)=sr}
{
\begin{align*}
p(x)&=3(r(x)+i)x-2ix-2xj-jx-xi+k
\\
&=3r(x)x+ix-2xj-jx-xi+k
\\
&=3r(x)x+i(r(x)+i)-2(r(x)+i)j-j(r(x)+i)-(r(x)+i)i+k
\\
&=3r(x)x+ir(x)-2r(x)j-jr(x)-r(x)i
\\
&=(3\otimes x+i\otimes 1-2\otimes j-j\otimes 1-1\otimes i)\circ r(x)
\end{align*}
}

\AddEq{s(x)=}
{
\[
s(x)=3\otimes x+i\otimes 1-2\otimes j-j\otimes 1-1\otimes i
\]
}

\AddEq{s(j)=}
{
\[
s(j)=1\otimes (j-i)-(j-i)\otimes 1\ne 0\otimes 0
\]
}

\AddEq[1]{Xk in Jv}
{
$X_{#1}\subseteq J(v)$.
}

\AddEq{w1+w2 in V}
{
w_1+w_2\in X_{k+1}
}

\AddEquation{unitarity law, left A-module}
{
1v=v
}

\AddEquation{unitarity law*, left A-module}
{
1*v=v
}

\AddEquation{left module, a d v}
{
a(dv)=d(av)
}

\AddEquation{left module, a* d v}
{
a*(dv)=d(a*v)
}

\AddEquation{module, a d v}
{
a(nv)=n(av)
}

\AddEq{abc in A, v in v}
{
$a$, $b$, $c\in A$, $v\in V$.
}

\AddEquation{a(b(cv))=(a(bc))v left}
{
a(b(cv))=a((bc)v)=(a(bc))v
}

\AddEquation{a(b(cv))=(a(bc))v right}
{
((vc)b)a=(v(cb))a=v((cb)a)
}

\AddEquation{a(b(cv))=((ab)c)v left}
{
a(b(cv))=(ab)(cv)=((ab)c)v
}

\AddEquation{a(b(cv))=((ab)c)v right}
{
((vc)b)a=(vc)(ba)=v(c(ba))
}

\AddEquation{(a(bc))v=((ab)c)v left}
{
(a(bc))v=((ab)c)v
}

\AddEquation{(a(bc))v=((ab)c)v right}
{
v((cb)a)=v(c(ba))
}

\AddEquation{a(bc)=(ab)c left}
{
a(bc)=(ab)c
}

\AddEquation{a(bc)=(ab)c right}
{
(cb)a=c(ba)
}

\AddEq{cv left}
{
$cv\in A$,
}

\AddEq{cv right}
{
$vc\in A$,
}

\AddEq{right module, associative law}
{
v(ab)=(va)b
}

\AddEq{associative law, D, right module}
{
v(dc)=(vd)c
}

\AddEq{associative law, DA, right module}
{
v(ad)=(va)d
}

\AddEquation{unitarity law, right A-module}
{
v1=v
}

\AddEquation{unitarity law*, right A-module}
{
v*1=v
}

\AddEquation{right module, a d v}
{
(vd)a=(va)d
}

\AddEquation{right module, a* d v}
{
(vd)*a=(v*a)d
}

\AddEq{c,d in Z, a,b in D, v,w in V}
{
$m$, $n\in Z$, $a$, $b \in D$, $v$, $w \in V$.
}

\AddEq [1]{vi V}
{
$v=(\aD vi\in V,\iIg)$#1
}

\AddEq [1]{v(i) V}
{
$v=(v(\gii)\in V,\iIg)$#1
}

\AddEq[2]{vv in V}
{
$\Vector{#1}\in #2$
}

\AddEq[3]{polynomial p(x)}
{
\[
#1(x)=\sum_{#2=0}^{#3}#1_{#2}x^{#2}
\]
}

\AddEquation{polynomial p*r(x)}
{
p(x)*r(x)=
(p*r)(x)=\sum_{j=0}^{m+n}(\sum_{i=0}^jp_ir_{j-i})x^j
}

\AddEq{x-=}
{
\[
\tilde{x}=hxh^{-1}=
(i-j)x(i-j)^{-1}=\frac 12(i-j)x(j-i)
\]
}

\AddEquation{L(x-)=}
{
L(\tilde{x})=\frac 12(i-j)x(j-i)-i
}

\AddEquation{p(x)=L(x-)R(x)}
{
x^2-(i+j)x+k=(\frac 12(i-j)x(j-i)-i)(x-j)
}

\AddEquation{p(i+j)=L(i+j-)R(i+j)}
{
(i+j)^2-(i+j)(i+j)+k=(\frac 12(i-j)(i+j)(j-i)-i)(i+j-j)
}

\AddEquation{p(i+j)=L(i+j-)R(i+j) 1}
{
\begin{split}
k&=(\frac 12(i^2+ij-ji-j^2)(j-i)-i)i\\
&=(k(j-i)-i)i=(-i-j-i)i
\end{split}
}

\AddEq{linear map f, basis E}
{
f=f\pC s0\otimes f\pC s1\in\AoxA H
}

\AddEquation{value of linear map f, basis E}
{
f\circ a=(f\pC s0\otimes f\pC s1)\circ a=f\pC s0\,a\,f\pC s1
}

\AddEq{Eox=x}
{
\[E\circ x=x\]
}

\AddEq{component of linear map, basis E}
{
\symb{f\pC sp}{component of linear map, division ring}1,
$p=0$, $1$,
}

\AddEq{vv=ve left module}
{
\Vector v=\ACol vi\ARow ei
}

\AddEq{vv=ve right module}
{
\Vector v=\ARow ei\ACol vi
}

\AddEq{vv=ve module}
{
\Vector v=
\ACol vi
\ARow ei
}

\AddEq{coordinate matrix of vector}
{
$v=(\ACol vi,\iIg)$
}

\AddEquation{w= module}
{
w=\sum_{\iIg}\aU wi\aD vi
}

\AddEquation{w= left module}
{
w=\sum_{\iIg}\aU wi\aD vi
}

\AddEquation{w= right module}
{
w=\sum_{\iIg}\aD vi\aU wi
}

\AddEquation{w()= module}
{
w=\sum_{\iIg}w(\gii)v(\gii)
}

\AddEquation{w()= left module}
{
w=\sum_{\iIg}w(\gii)v(\gii)
}

\AddEquation{w()= right module}
{
w=\sum_{\iIg}v(\gii)w(\gii)
}

\AddEquation{aw= module}
{
aw=a\sum_{\iIg}\aU wi\aD vi
=\sum_{\iIg}a(\aU wi\aD vi)
=\sum_{\iIg}(a\aU wi)\aD vi
}

\AddEquation{aw= left module}
{
aw=a\sum_{\iIg}\aU wi\aD vi
=\sum_{\iIg}a(\aU wi\aD vi)
=\sum_{\iIg}(a\aU wi)\aD vi
}

\AddEquation{aw= right module}
{
wa=\left(\sum_{\iIg}\aD vi\aU wi\right)a
=\sum_{\iIg}(\aD vi\aU wi)a
=\sum_{\iIg}(\aD vi\aU wia)
}

\AddEquation{aw()= module}
{
aw=a\sum_{\iIg}w(\gii)v(\gii)
=\sum_{\iIg}a(w(\gii)v(\gii))
=\sum_{\iIg}(aw(\gii))v(\gii)
}

\AddEquation{aw()= left module}
{
aw=a\sum_{\iIg}w(\gii)v(\gii)
=\sum_{\iIg}a(w(\gii)v(\gii))
=\sum_{\iIg}(aw(\gii))v(\gii)
}

\AddEquation{aw()= right module}
{
wa=\left(\sum_{\iIg}v(\gii)w(\gii)\right)a
=\sum_{\iIg}(v(\gii)w(\gii))a
=\sum_{\iIg}(\aD viw(\gii)a)
}

\AddEq[1]{set au vi}
{
$\aU {#1}i$, \iIg,
}

\AddEq[1]{set au v(i)}
{
$#1(\gii)$, \iIg,
}

\AddEq[1]{set au v1n}
{
$\aU {#1}1$, ..., $\aU {#1}n$
}

\AddEq{basis, module}
{
\symb{\Basis{e}=(
\ARow ei,
\iIg)}{basis for module}1
}

\AddEq[3]{set vi Acols}
{%
$\ACol {#1}{#2}$, \jJg{#2}{#3},
}

\AddEq[1]{set vi cols}
{%
$\aD {#1}i$, \iIg,
}

\AddEq[1]{set vi *}
{%
$#1(\gii)$, \iIg,
}

\AddEquation{coordinate transformation, ba}
{
v''^k=b^{-1}{}^k_ia^{-1}{}^i_jv^j=(ab)^{-1}{}^k_jv^j
}

\AddEq[3]{passive coordinate transformation}
{
v#1^i=#3^{-1}{}^i_jv#2^j
}

\AddEq[3]{w.e, left-cols}
{
#1\CRstar #2=
\ColMatrix{#1}{#3}
\CRstar
\RowMatrix{#2}{#3}
=\ShowEq{w.e 1, cr}{#1}{#2}
}

\AddEq[3]{w.e, right-rows}
{
#2\CRstar #1=
\ColMatrix{#2}{#3}
\CRstar
\RowMatrix{#1}{#3}
=\ShowEq{w.e 1, cr}{#2}{#1}
}

\AddEq[3]{w.e, left-rows}
{
#1\RCstar #2=
\RowMatrix{#1}{#3}
\RCstar
\ColMatrix{#2}{#3}
=\ShowEq{w.e 1, rc}{#1}{#2}
}

\AddEq[3]{w.e, right-cols}
{
#2\RCstar #1=
\RowMatrix{#2}{#3}
\RCstar
\ColMatrix{#1}{#3}
=\ShowEq{w.e 1, rc}{#2}{#1}
}

\AddEq[3]{w.e, -cols}
{
#1 #2=
\ColMatrix{#1}{#3}
\RowMatrix{#2}{#3}
=\ShowEq{w.e 1, -cols}{#1}{#2}
}

\AddEq[3]{w.e, -rows}
{
#2 #1=
\ColMatrix{#2}{#3}
\RowMatrix{#1}{#3}
=\ShowEq{w.e 1, -rows}{#1}{#2}
}

\AddEq[2]{w.e 1, cr}
{
\aU {#1}i\aD {#2}i
}

\AddEq[2]{w.e 1, rc}
{
\aD {#1}i\aU {#2}i
}

\AddEq[2]{w.e 1, -cols}
{
\aU {#1}i\aD {#2}i
}

\AddEq[2]{w.e 1, -rows}
{
\aD {#1}i\aU {#2}i
}

\AddEq[1]{vector w=w*e left-cols}
{
\Vector {#1}=
\ShowEq{w.e, left-cols}{#1}en
}

\AddEq[1]{vector w=w*e left-rows}
{
\Vector {#1}=
\ShowEq{w.e, left-rows}{#1}en
}

\AddEq[1]{vector w=w*e right-cols}
{
\Vector {#1}=
\ShowEq{w.e, right-cols}{#1}en
}

\AddEq[1]{vector w=w*e right-rows}
{
\Vector {#1}=
\ShowEq{w.e, right-rows}{#1}en
}

\AddEq[1]{vector w=w*e -cols}
{
\Vector {#1}=
\ShowEq{w.e, -cols}{#1}en
}

\AddEq[1]{vector w=w*e -rows}
{
\Vector {#1}=
\ShowEq{w.e, -rows}{#1}en
}

\AddEq[3]{w=w*eC left-cols}
{
\ensuremath{\Vector {#1}=
\aU {#1}i\aD [#2]ei}#3
}

\AddEq[2]{w=w*e left-cols}
{
\ensuremath{\Vector {#1}=
\ShowEq{w.e 1, cr}{#1}e}#2
}

\AddEq[2]{w=w*e left-rows}
{
\ensuremath{\Vector {#1}=
\ShowEq{w.e 1, rc}{#1}e}#2
}

\AddEq[2]{w=w*e right-cols}
{
\ensuremath{\Vector {#1}=
\ShowEq{w.e 1, rc}e{#1}}#2
}

\AddEq[2]{w=w*e right-rows}
{
\ensuremath{\Vector {#1}=
\ShowEq{w.e 1, cr}e{#1}}#2
}

\AddEq[2]{w=w*e -cols}
{
\ensuremath{\Vector {#1}=
\ShowEq{w.e 1, cr}{#1}e}#2
}

\AddEq[2]{w=w*e -rows}
{
\ensuremath{\Vector {#1}=
\ShowEq{w.e 1, rc}{#1}e}#2
}

\AddEquation{similar A numbers}
{
a=c^{-1}bc
}

\AddEquation{endomorphism of module from product}
{
\xymatrix
{
g_{23}:A\ar[r]|{*}&A
}
}

\DefEquation
{
l\circ v:w\in A\rightarrow v\,w\in A
}
{l(v):w->vw}

\AddEq[2]{representation An in LAnA}
{
\[
\xymatrix
{
h:\AOn[#2]\times S_n\ar[r]|{*}&\LAnAB {#1}{#2}{#2}
}
\]
}

\AddEq[2]{representation An in LAnA, 1}
{
\begin{matrix}
(a_0\otimes...\otimes a_n,\sigma)\circ(f_1\otimes...\otimes f_n)
=a_0\sigma(f_1)a_1...a_{n-1}\sigma(f_n)a_n
\\
\begin{matrix}
a_0,...,a_n\in #2&\sigma\in S_n&f_1,...,f_n\in\mathcal L(#1;#2\rightarrow #2)
\end{matrix}
\end{matrix}
}

\DefEq
{
$d\in\AOn$
}
{product in algebra An 1}

\AddEq[2]{product in algebra An 2}
{
\begin{matrix}
(d,\sigma)\circ(f_1,...,f_n)&f_i=\delta\in\mathcal L(#1;#2\rightarrow #2)
\end{matrix}
}

\AddEq[3]{product in algebra AA 3}
{
$#3\in\mathcal L(#1;#2\rightarrow #2)$
}

\AddEq[2]{wi vi=0}
{
\ACol {#1}i\ARow {#2}i=0
}

\AddEq[2]{left wi vi=0}
{
\ACol {#1}i\ARow {#2}i=0
}

\AddEq[2]{right wi vi=0}
{
\ARow {#2}i\ACol {#1}i=0
}

\AddEq{w(i)=0}
{
$w(\gii)=0$, \iIg,
}

\AddEq{w(i) v(i)=0}
{
w(\gii)v(\gii)=0
}

\AddEq{left w(i) v(i)=0}
{
w(\gii)v(\gii)=0
}

\AddEq{right w(i) v(i)=0}
{
v(\gii)w(\gii)=0
}

\AddEq{0=0vi}
{
\[
\ACol wi=0
\Rightarrow
\ACol wi
\ARow vi
=0
\]
}

\AddEq{left 0=0vi}
{
\[
\ACol wi=0
\Rightarrow
\ACol wi
\ARow vi
=0
\]
}

\AddEq{right 0=0vi}
{
\[
\ACol wi=0
\Rightarrow
\ARow vi
\ACol wi
=0
\]
}

\AddEq[2]{vj=sum vi}
{
\[
\ARow {#1}j{}
=\sum_{\gii\in \giI\setminus\{\gij\}}
\frac
{\ACol {#2}i{}}
{\ACol {#2}j{}}
\ARow {#1}i{}
\]
}

\AddEq[2]{left vj=sum vi}
{
\[
\ARow {#1}j{}
=\sum_{\gii\in \giI\setminus\{\gij\}}
(\ACol {#2}j{})^{-1}
\ACol {#2}i{}
\ARow {#1}i{}
\]
}

\AddEq[2]{right vj=sum vi}
{
\[
\ARow {#1}j{}
=\sum_{\gii\in \giI\setminus\{\gij\}}
\ARow {#1}i{}
\ACol {#2}i{}
(\ACol {#2}j{})^{-1}
\]
}

\AddEq[2]{b in D'}
{
$#1=(
\ACol {#1}i,
\iIg)\in\Base'$#2
}

\AddEq[2]{wj ne 0}
{
$\ACol {#1}j{}\ne 0$#2
}

\AddEq[2]{w12 in Jv}
{
$#1_1$, $#1_2\in X_k\subseteq J(#2)$.
}

\AddEq[1]{w=w12 in Xk-1}
{
$w=w_1+w_2$, $w_1$, $w_2\in X_{k-1}$#1
}

\AddEq{w=aw in Xk-1,}
{
$w=aw_1$, $a\in\Base_{\BaseExt}$, $w_1\in X_{k-1}$.
}

\AddEq{w=aw in Xk-1,left}
{
$w=aw_1$, $a\in\Base_{\BaseExt}$, $w_1\in X_{k-1}$.
}

\AddEq{w=aw in Xk-1,right}
{
$w=w_1a$, $a\in\Base_{\BaseExt}$, $w_1\in X_{k-1}$.
}

\AddEq{generating set of module}
{
$J(X)=V$.
}

\AddEq[2]{w=sum vi, module}
{
J(v)=\left\{w:w=\sum_{\iIg}c(\gii)v(\gii),c(\gii)\in #1_{#2},
|\{\gii:c(\gii)\ne 0\}|<\infty\right\}
}

\AddEq[2]{w=sum vi, left module}
{
J(v)=\left\{w:w=\sum_{\iIg}c(\gii)v(\gii),c(\gii)\in #1_{#2},
|\{\gii:c(\gii)\ne 0\}|<\infty\right\}
}

\AddEq[2]{w=sum vi, right module}
{
J(v)=\left\{w:w=\sum_{\iIg}v(\gii)c(\gii),c(\gii)\in #1_{#2},
|\{\gii:c(\gii)\ne 0\}|<\infty\right\}
}

\AddEq{vectors generated by set of vectors}
{
\StartLabelItem
\begin{enumerate}
\item
\ShowEq{\SideWS module, sum cvk in Jv}
\item
\ShowEq{v,w in Jv=>v+w in Jv}
\item
\ShowEq{vector space, v in Jv=>av in Jv}
\end{enumerate}
}

\AddEq{vectors generated by set of vectors 2018}
{
\StartLabelItem
\begin{enumerate}
\item
\ShowEq{vk in v}
\labelItem{vector space, vk in Jv}
\item
\ShowEq{vector space, cvk in Jv}
\item
\ShowEq{vector space, sum cvk in Jv}
\item
\ShowEq{v,w in Jv=>v+w in Jv}
\item
\ShowEq{vector space, v in Jv=>av in Jv}
\end{enumerate}
}

\AddEq[1]{vk in Jv}
{%
$v(\gik)\in J(v)$#1
}

\AddEq{module, d in D}
{
$\DArg\in \CBase$,
}

\AddEq{module, dv in V}
{%
$\DArg v\in V$,
}

\AddEq{left module, dv in V}
{%
$\DArg v\in V$,
}

\AddEq{right module, dv in V}
{%
$v\DArg \in V$,
}

\AddEq{module, a in A}
{
$a\in \Base$.
}

\AddEq{module, av in V}
{%
$av\in V$
}

\AddEq{left module, av in V}
{%
$av\in V$
}

\AddEq{right module, av in V}
{%
$va\in V$
}

\AddEq{module, dv+av in V}
{%
\[\DArg v+av\in V\ \ \ \DArg\in \CBase\ \ \ a\in \Base\]
}

\AddEq{left module, dv+av in V}
{%
\[\DArg v+av\in V\ \ \ \DArg\in \CBase\ \ \ a\in \Base\]
}

\AddEq{right module, dv+av in V}
{%
\[v\DArg +va\in V\ \ \ \DArg\in \CBase\ \ \ a\in \Base\]
}

\AddEq{module, cvk in Jv}
{%
$c(\gik)v(\gik)\in J(v)$, $c(\gik)\in \Base_{\BaseExt}$,
\jJg kI
\labelItem{\VectorSetNS, cvk in Jv}
}

\AddEq{vector space, cvk in Jv}
{%
$c^{\gik}v_{\gik}\in J(v)$, $c^{\gik}\in \Base$,
\jJg kI
\labelItem{vector space, cvk in Jv}
}

\AddEquation{module, a+n:V->V}
{
a+n:v\in V\rightarrow (a+n)v\in V
}

\AddEquation{left module, a+n:V->V}
{
a+n:v\in V\rightarrow (a+n)v\in V
}

\AddEquation{right module, a+n:V->V}
{
a+n:v\in V\rightarrow v(a+n)\in V
}

\AddEquation{left module, a*+d:V->V}
{
a+d:v\in V\rightarrow (a+d)*v\in V
}

\AddEquation{right module, a*+d:V->V}
{
a+d:v\in V\rightarrow v*(a+d)\in V
}

\AddEq{module, sum cvk in Jv}
{%
$\displaystyle\sum_{\jJg kI}c(\gik)v(\gik)\in J(v)$,
$c(\gik)\in\Base_{\BaseExt}$, $|\{\gii:c(\gii)\ne 0\}|<\infty$
\labelItem{\VectorSetNS, sum cvk in Jv}
}

\AddEq{vector space, sum cvk in Jv}
{%
$\displaystyle\sum_{\jJg kI}c^{\gik}v_{\gik}\in J(v)$,
$c^{\gik}\in\Base$, $|\{\gii:c^{\gii}\ne 0\}|<\infty$
\labelItem{vector space, sum cvk in Jv}
}

\AddEq{left module, sum cvk in Jv}
{%
$\displaystyle\sum_{\jJg kI}(\gik)v(\gik)\in J(v)$,
$c(\gik)\in\Base_{\BaseExt}$, $|\{\gii:c(\gii)\ne 0\}|<\infty$
\labelItem{left \VectorSetNS, sum cvk in Jv}
}

\AddEq{right module, sum cvk in Jv}
{%
$\displaystyle\sum_{\jJg kI}v(\gik)c(\gik)\in J(v)$,
$c(\gik)\in\Base_{\BaseExt}$, $|\{\gii:c(\gii)\ne 0\}|<\infty$
\labelItem{right \VectorSetNS, sum cvk in Jv}
}

\AddEq[1]{|ci12 ne 0|}
{
H_1=\{\iIg:\aU{#1_1}i\ne 0\}
\ \ \ \,
H_2=\{\iIg:\aU{#1_2}i\ne 0\}
}

\AddEq[1]{|c(i)12 ne 0|}
{
H_1=\{\iIg:#1_1(\gii)\ne 0\}
\ \ \ \,
H_2=\{\iIg:#1_2(\gii)\ne 0\}
}

\AddEq{|wi ne 0|}
{
|\{\iIg:\aU wi\ne 0\}|<\infty
}

\AddEq{|w(i) ne 0|}
{
|\{\iIg:w(\gii)\ne 0\}|<\infty
}

\AddEq{module, |aw(i) ne 0|}
{
$\{\iIg:aw(\gii)\ne 0\}$
}

\AddEq{left module, |aw(i) ne 0|}
{
$\{\iIg:aw(\gii)\ne 0\}$
}

\AddEq{right module, |aw(i) ne 0|}
{
$\{\iIg:w(\gii)a\ne 0\}$
}

\AddEq{module, |awi ne 0|}
{
$\{\iIg:a\aU wi\ne 0\}$
}

\AddEq{left module, |awi ne 0|}
{
$\{\iIg:a\aU wi\ne 0\}$
}

\AddEq{right module, |awi ne 0|}
{
$\{\iIg:\aU wia\ne 0\}$
}

\AddEq{v,w in Jv=>v+w in Jv}
{
$w_1$, $w_2\in J(v)\Rightarrow w_1+w_2\in J(v)$
\labelItem{\SideWS module, v,w in Jv=>v+w in Jv}
}

\AddEq{v in Jv=>av in Jv}
{
$a\in\Base_{(1)}$,
$w\in J(v)\Rightarrow aw\in J(v)$
\labelItem{\SideWS module, v in Jv=>av in Jv}
}

\AddEq{vector space, v in Jv=>av in Jv}
{
$a\in\Base_{\BaseExt}$,
$w\in J(v)\Rightarrow aw\in J(v)$
\labelItem{\SideWS vector space, v in Jv=>av in Jv}
}

\AddEq[1]{gi12}
{
$\aU{#1_1}i$, $\aU{#1_2}i$, \iIg,
}

\AddEq[1]{g(i)12}
{
$#1_1(\gii)$, $#1_2(\gii)$, \iIg,
}

\AddEq[1]{w1+w2= 1 }
{
#1_1+#1_2=\sum_{\iIg}(\aU{#1_1}i+\aU{#1_2}i)v_{\gii}
}

\AddEq[1]{w1+w2= 1 left}
{
#1_1+#1_2=\sum_{\iIg}(\aU{#1_1}i+\aU{#1_2}i)v_{\gii}
}

\AddEq[1]{w1+w2= 1 right}
{
#1_1+#1_2=\sum_{\iIg}v_{\gii}\aU{#1_1}i+\aU{#1_2}i)
}

\AddEq[1]{|gi 1+2 ne 0|}
{
\[
\{\iIg:\aU{#1_1}i+\aU{#1_2}i\ne 0\}\subseteq H_1\cup H_2
\]
}

\AddEq[1]{|gi() 1+2 ne 0|}
{
\[
\{\iIg:#1_1(\gii)+#1_2(\gii)\ne 0\}\subseteq H_1\cup H_2
\]
}

\AddEq[1]{w1+w2= }
{
#1_1+#1_2=\sum_{\iIg}\aU{#1_1}iv_{\gii}+\sum_{\iIg}\aU{#1_2}iv_{\gii}
=\sum_{\iIg}(\aU{#1_1}iv_{\gii}+\aU{#1_2}iv_{\gii})
}

\AddEq[1]{w1+w2= left}
{
#1_1+#1_2=\sum_{\iIg}\aU{#1_1}iv_{\gii}+\sum_{\iIg}\aU{#1_2}iv_{\gii}
=\sum_{\iIg}(\aU{#1_1}iv_{\gii}+\aU{#1_2}iv_{\gii})
}

\AddEq[1]{w1+w2= right}
{
#1_1+#1_2=\sum_{\iIg}v_{\gii}\aU{#1_1}i+\sum_{\iIg}v_{\gii}\aU{#1_2}i
=\sum_{\iIg}(v_{\gii}\aU{#1_1}i+v_{\gii}\aU{#1_2}i)
}

\AddEq[1]{w()1+w()2= }
{
\begin{split}
#1_1+#1_2&=\sum_{\iIg}#1_1(\gii)v(\gii)+\sum_{\iIg}#1_2(\gii)v(\gii)
=\sum_{\iIg}(#1_1(\gii)v(\gii)+#1_2(\gii)v(\gii))
\\ &=\sum_{\iIg}(#1_1(\gii)+#1_2(\gii))v(\gii)
\end{split}
}

\AddEq[1]{w()1+w()2= left}
{
\begin{split}
#1_1+#1_2&=\sum_{\iIg}#1_1(\gii)v(\gii)+\sum_{\iIg}#1_2(\gii)v(\gii)
=\sum_{\iIg}(#1_1(\gii)v(\gii)+#1_2(\gii)v(\gii))
\\ &=\sum_{\iIg}(#1_1(\gii)+#1_2(\gii))v(\gii)
\end{split}
}

\AddEq[1]{w()1+w()2= right}
{
\begin{split}
#1_1+#1_2&=\sum_{\iIg}v(\gii)#1_1(\gii)+\sum_{\iIg}v(\gii)#1_2(\gii)
=\sum_{\iIg}(v(\gii)#1_1(\gii)+v(\gii)#1_2(\gii))
\\ &=\sum_{\iIg}v(\gii)(#1_1(\gii)+#1_2(\gii))
\end{split}
}

\AddEq[1]{w12= }
{
#1_1=\sum_{\iIg}\aU{#1_1}i\aD vi
\ \ \ \,
#1_2=\sum_{\iIg}\aU{#1_2}i\aD vi
}

\AddEq[1]{w12= left}
{
#1_1=\sum_{\iIg}\aU{#1_1}i\aD vi
\ \ \ \,
#1_2=\sum_{\iIg}\aU{#1_2}i\aD vi
}

\AddEq[1]{w12= right}
{
#1_1=\sum_{\iIg}\aD vi\aU{#1_1}i
\ \ \ \,
#1_2=\sum_{\iIg}\aD vi\aU{#1_2}i
}

\AddEq[1]{w()12= }
{
#1_1=\sum_{\iIg}#1_1(\gii)v(\gii)
\ \ \ \,
#1_2=\sum_{\iIg}#1_2(\gii)v(\gii)
}

\AddEq[1]{w()12= left}
{
#1_1=\sum_{\iIg}#1_1(\gii)v(\gii)
\ \ \ \,
#1_2=\sum_{\iIg}#1_2(\gii)v(\gii)
}

\AddEq[1]{w()12= right}
{
#1_1=\sum_{\iIg}v(\gii)#1_1(\gii)
\ \ \ \,
#1_2=\sum_{\iIg}v(\gii)#1_2(\gii)
}

\AddEquation{vk=sum vi, module}
{
v_{\gik}=\sum_{\iIg}\aU ci\aD vi
}

\AddEquation{vk=sum vi, left module}
{
v_{\gik}=\sum_{\iIg}\aU ci\aD vi
}

\AddEquation{vk=sum vi, right module}
{
\aD vk=\sum_{\iIg}\aD vi\aU ci
}

\AddEq[1]{v(k)=sum v(i), module}
{
v(\gik)=\sum_{\iIg}c(\gii)v(\gii)
\ \ \ c(\gii)\in #1
}

\AddEq[1]{v(k)=sum v(i), left module}
{
v(\gik)=\sum_{\iIg}c(\gii)v(\gii)
\ \ \ c(\gii)\in #1
}

\AddEq[1]{v(k)=sum v(i), right module}
{
v(\gik)=\sum_{\iIg}v(\gii)c(\gii)
\ \ \ c(\gii)\in #1
}

\AddEq[1]{ci=dik}
{
$c^{\gii}=\delta^{\gii}_{\gik}\in #1$.
}

\AddEq{c(i)=dik}
{
\[
c(\gii)=
\left\{
\begin{array}{l}
1,\ \gii=\gik\\
0
\end{array}
\right.
\]
}

\AddEq{vk in v}
{
$v_{\gik}\in v$,
}

\AddEq{v(i)}
{
$v(\gii)$.
}

\AddEq{w=c(i)v(i).}
{
$w=\sum c(\gii)v(\gii)$
}

\AddEq{w=c(i)v(i).left}
{
$w=\sum c(\gii)v(\gii)$
}

\AddEq{w=c(i)v(i).right}
{
$w=\sum v(\gii)c(\gii)$
}

\AddEq{w=w cr v}
{
\[\Vector w=\ShowSymbol{linear combination}{\SideWS cr}\]
}

\AddEq{g number linear combination}
{%
\symb{g^is_i}{linear combination}{}%
}

\AddEq{linear combination()}
{%
\symb{\sum c(\gii)v(\gii)}{linear combination}{i}%
}

\AddEq{left linear combination()}
{%
\symb{\sum c(\gii)v(\gii)}{linear combination}{left i}%
}

\AddEq{right linear combination()}
{%
\symb{\sum v(\gii)c(\gii)}{linear combination}{right i}%
}

\AddEq{the left linear combination()}
{%
$c(\gii)v(\gii)$
}

\AddEq{the right linear combination()}
{%
$v(\gii)c(\gii)$
}

\AddEq{A linear combination() =}
{
$\ShowSymbol{linear combination}{\SideWS i}$, $c(\gii)\in\Base_{\BaseExt}$,
}

\AddEq{(b+c)*e, module}
{
\[
(\ACol bi+
\ACol ci)
\ARow ei=0
\]
}

\AddEq{(b+c)*e, left module}
{
\[
(\ACol bi+
\ACol ci)
\ARow ei=0
\]
}

\AddEq{(b+c)*e, right module}
{
\[
\ARow ei
(\ACol bi+
\ACol ci)
=0
\]
}

\AddEq{(ac)*e, module}
{
\[
(a
\ACol ci)
\ARow ei=0
\]
}

\AddEq{(ac)*e, left module}
{
\[
(a
\ACol ci)
\ARow ei=0
\]
}

\AddEq{(ac)*e, right module}
{
\[
\ARow ei(a
\ACol ci)=0
\]
}

\AddEq{sum of maps,,D module}
{
\symb{f+g}{sum of maps}{D module}
}

\AddEquation{sum of maps, 1, D module}
{
\begin{matrix}
\ShowSymbol{sum of maps}{D module}:A_1\rightarrow A_2
&f,g\in\LAB D{A_1}{A_2}
\end{matrix}
}

\DefEquation
{
(\ShowSymbol{sum of maps}{D module})\circ x=f\circ x+g\circ x
}
{sum of maps, D module}

\AddEq [4]{f in L(A->B)}
{
$f\in\LAB{#1}{#2}{#3}$#4
}

\AddEquation{fa=sfsa}
{
f\circ(a_1, ..., a_n)=|\sigma|(f\circ\sigma\circ(a_1,...,a_n))
}

\DefEq
{
}
{}

\AddEq[1]{[f]}
{
$\ShowSymbol{alternation of polylinear map}{}$#1
}

\AddEq{alternation of polylinear map}
{
\symb{[f]}{alternation of polylinear map}{}
}

\AddEq{fx1n=0}
{
\[
f\circ(a_1,...,a_n)=0
\]
}

\AddEq{1<=i<n}
{
$i$, $1\le i<n$.
}

\AddEq{fa=fsa}
{
\[
f\circ(a_1, ..., a_n)=f\circ\sigma\circ(a_1,...,a_n)
\]
}

\AddEq{symmetrization of polylinear map}
{
\symb{<f>}{symmetrization of polylinear map}{}
}

\AddEq{<f>}
{
$\ShowSymbol{symmetrization of polylinear map}{}\circ(a_1,...,a_n)$\Pt
}

\AddEq[1]{f()}
{
$f\circ(a_1,...,a_p)$#1
}

\AddEq{h:B1->*B2}
{
\[
\xymatrix
{
h:B_1\ar[r]|{*}&B_2
}
\]
}

\AddEq[3]{L(A1,A0,B)=}
{
\begin{align}
\mathcal{LA}(#1;#2\rightarrow #3)&=\mathcal L(#1;#2\rightarrow #3)\\
\mathcal{LA}(#1;#2^0\rightarrow #3)&=#3
\end{align}
}

\newcommand\LAnAB[3]{\ensuremath{\mathcal{LA}(#1;#2^n\rightarrow #3)}}

\AddEq[3]{module of skew symmetric polylinear maps}
{
\symb{\LAnAB{#1}{#2}{#3}}{module of skew symmetric polylinear maps}{1}
}

\AddEquation{exterior product = skew symmetric}
{
(f\wedge g)\circ(a_1,...,a_{p+q})=
\frac 1{p!q!}\sum_{\sigma\in S(p+q)}|\sigma|((f\wedge g)\circ\sigma\circ(a_1,...,a_{p+q}))
}

\AddEquation{exterior product =}
{
(\ShowSymbol{exterior product}{1})\circ(a_1,...,a_{p+q})=\frac {(p+q)!}{p!q!}
[fg]\circ(a_1,...,a_{p+q})
}

\AddEq{exterior product}
{
\symb{f\wedge g}{exterior product}{1}
}

\AddEquation{hpq(fg)=}
{
(fg)\circ(a_1,...,a_{p+q})=(f\circ(a_1,...,a_p))(g\circ(a_{p+1},...,a_{p+q}))
}

\AddEq{hpq:Lpq->Lp+q}
{
\[
fg:\mathcal L(D;A^p\rightarrow B)\times\mathcal L(D;A^q\rightarrow C)\rightarrow \mathcal L(D;A^{p+q}\rightarrow C)
\]
}

\AddEq[1]{g()}
{
$g\circ(a_{p+1},...,a_{p+q})$#1
}

\AddEq{symmetrization of polylinear map =}
{
\[
\ShowSymbol{symmetrization of polylinear map}{}\circ(a_1,...,a_n)
=\frac 1{n!}
\sum_{\sigma\in S(n)}f\circ\sigma\circ(a_1,...,a_n)
\]
}

\AddEquation{alternation of polylinear map =}
{
\ShowSymbol{alternation of polylinear map}{}\circ(a_1,...,a_n)
=\frac 1{n!}\sum_{\sigma\in S(n)}|\sigma|(f\circ\sigma\circ(a_1,...,a_n))
}

\AddEq[3]{A 1n}
{
$\aU{#1}1$, ..., $\aU{#1}{#2}$#3
}

\AddEq{fab=fbfa}
{
\[
f(ab)=f(b)f(a)
\]
}

\AddEq{b o f =}
{
\[
(b\circ f)\circ x=b(f\circ x)
\]
}

\AddEq{b * f =}
{
\[
(b*f)\circ x=(f\circ x)b
\]
}

\DefEq
{
$F_k\in\Basis F$,
}
{Ik in I}

\DefEquation
{
x\rightarrow I_k\circ a\circ x
}
{x->Ik a x}

\DefEquation
{
b^l=I^l_k\circ a
}
{b=Ikl a}

\DefEq
{
\begin{align*}
(I_k^l\circ(a_1+\aD a2))\circ I_l\circ x
&=I_k\circ(a_1+\aD a2)\circ x
\\&=I_k\circ a_1\circ x+I_k\circ \aD a2\circ x
\\&=(I_k^l\circ a_1)\circ I_l\circ x+(I_k^l\circ a_1)\circ I_l\circ x
\end{align*}
}
{Ikl a1+a2}

\DefEq
{
\begin{align*}
(I_k^l\circ(da))\circ I_l\circ x
&=I_k\circ(da)\circ x
=I_k\circ(d(a\circ x))
\\&=d(I_k\circ a\circ x)
=d((I_k^l\circ a)\circ I_l\circ x)
\\&=(d(I_k^l\circ a))\circ I_l\circ x
\end{align*}
}
{Ikl d a}

\DefEquation
{
I_k\circ a\circ x=b^l\circ I_l\circ x\ \ \ \ b^l\in A_2\times A_2
}
{expansion Ik a x}

\DefEq
{
\symb{F_k^l}{conjugation transformation}{}
}
{conjugation transformation}

\DefEq
{
\[
\ShowSymbol{conjugation transformation}{}
:A_1\otimes A_1\rightarrow A_2\otimes A_2
\]
}
{conjugation transformation:}

\DefEquation
{
F_k\circ a\circ x=(\ShowSymbol{conjugation transformation}{}\circ a)\circ F_l\circ x
}
{conjugation transformation =}

\DefEq
{
$\ShowSymbol{conjugation transformation}{}$
}
{show conjugation transformation}

\DefEq
{
$\aUD Cp{kl}$
}
{structural constants, algebra}

\AddEq{ci i in I}
{
\ACol ci, \iIg,
}

\AddEq{c=ci i in I}
{
$c=(\ACol ci,\iIg)$
}

\AddEq{i=j}
{
$\gii=\gij$
}

\AddEq{Coordinates of map f}
{
$\aUD fkl$
}

\AddEq[1]{standard components of map f}
{
\def\Cdot{}
\def\Temp{-}%
\edef\Tempa{#1}%
\ifx\Tempa\Temp%
\else
\def\Cdot{#1\cdot}
\fi
$\aU {f^{\Cdot}_{}{}}{ij}$
}

\AddEq{linear map over ring, left side}
{
$\aUD flk$
}

\AddEq{left shift algebra, 2}
{
\[
(a,b)_1\circ x=(a,b,x)
\]
}

\AddEquation{left shift algebra}
{
l(a)\circ x=ax
}

\AddEquation{left shift algebra, 3}
{
\begin{split}
(l(a)\circ l(b))\circ x
&=l(a)\circ(l(b)\circ x)
\\
&=a(bx)=(ab)x-(a,b,x)
\\
&=l(ab)\circ x-(a,b)_1\circ x
\end{split}
}

\AddEquation{left shift algebra, 1}
{
l(a)\circ l(b)=l(ab)-(a,b)_1
}

\AddEquation{right shift algebra}
{
r(a)\circ x=xa
}

\AddEquation{right shift algebra, 3}
{
\begin{split}
(r(a)\circ r(b))\circ x
&=r(a)\circ(r(b)\circ x)
\\
&=(xb)a=x(ba)+(x,b,a)
\\
&=r(ba)\circ x+(x,b,a)
\end{split}
}

\AddEquation{right shift algebra, 1}
{
r(a)\circ r(b)=r(ba)+(b,a)_2
}

\AddEq{right shift algebra, 2}
{
\[
(b,a)_2\circ x=(x,b,a)
\]
}

\AddEq{f over A}
{
\[
f:x\in A\rightarrow (ax)b\in A
\]
}

\AddEq{g over A}
{
\[
\begin{matrix}
g:A\rightarrow A
&&
g=(cf)d
\end{matrix}
\]
}

\AddEquation{linear map, standard representation, nonassociative algebra}
{
f\circ x=\aU f{ij}\ (\aD eix)\aD ej
}

\AddEq{linear map, standard representation, nonassociative algebra, 1}
{
\[
f\circ x=\aU f{ij}\aD ei(x\aD ej)
\]
}

\AddEquation{coordinates of map A1 A2, 2, nonassociative algebra}
{
\aUD gkl
=\aUD fml\aU g{ij}
\aUD Cp{im}\aUD Ck{pj}
}

\AddEquation{coordinates of map A1 A2, 21, nonassociative algebra}
{
\aUD gkl
=\aUD fml \aU g{ij}
\aUD Ck{ip}\aUD Cp{mj}
}

\DefEq
{
$F_{k\cdot}^{}{}_{\gii}^{\gij}$
}
{coordinates of map Ik}

\DefEq
{
$G_{\gii}^{\gij}$
}
{coordinates of map G}

\AddEq{projection maps, complex field}
{
\begin{align*}
P^{\gi 0}\circ a&=\aU a0&P^{\gi 0\cdot}_{}{}^{\gi 0}_{\gi 0}&=1&
R^{\gi 0}\circ a&=\aU a0&R^{\gi 0\cdot}_{}{}^{\gi 0}_{\gi 0}&=1
\\
P^{\giA}\circ a&=\aU a1i&P^{\giA\cdot}_{}{}^{\giA}_{\giA}&=1&
R^{\giA}\circ a&=\aU a1&R^{\giA\cdot}_{}{}^{\gi 0}_{\giA}&=1
\end{align*}
}

\AddEq [1]{e=(E,I)}
{
$\eV=(E,I)$#1
}

\AddEq{projection mappings 1, complex field}
{
\begin{align*}
P^{\gi 0}&=\frac 12\circ E+\frac 12\circ I&
R^{\gi 0}&=\frac 12\circ E+\frac 12\circ I
\\
P^{\giA}&=\frac 12\circ E-\frac 12\circ I&
R^{\giA}&=-\frac i2\circ E+\frac i2\circ I
\end{align*}
}

\AddEquation{df/dx= complex field}
{
\frac{df}{dx}=
\frac 12\left(\frac{\partial f}{\partial \aU x0}
-i\frac{\partial f}{\partial x^{\gi 1}}\right)\circ E
+\frac 12\left(\frac{\partial f}{\partial \aU x0}
+i\frac{\partial f}{\partial x^{\gi 1}}\right)\circ I
}

\DefEquation
{
f=\aD a0\circ E+a_1\circ I+\aD a2\circ J+\aD a3\circ K
}
{expand linear mapping, quaternion}

\AddEquation{fi=fij ej H}
{
\begin{pmatrix}
\aD f0\\ \aD f1\\ \aD f2\\ \aD f3
\end{pmatrix}
=
\begin{pmatrix}
\aUD f00&\aUD f10&\aUD f20&\aUD f30\\
\aUD f01&\aUD f11&\aUD f21&\aUD f31\\
\aUD f02&\aUD f12&\aUD f22&\aUD f32\\
\aUD f03&\aUD f13&\aUD f23&\aUD f33
\end{pmatrix}
\begin{pmatrix}
1\\i\\j\\k
\end{pmatrix}
=
\begin{pmatrix}
\aUD f00+\aUD f10i+\aUD f20j+\aUD f30k\\
\aUD f01+\aUD f11i+\aUD f21j+\aUD f31k\\
\aUD f02+\aUD f12i+\aUD f22j+\aUD f32k\\
\aUD f03+\aUD f13i+\aUD f23j+\aUD f33k
\end{pmatrix}
}

\AddEq{fi=fij ej}
{
\aD fi=\aUD fji \aD ej
}

\AddEquation{f=.E+.I}
{
f=\frac 12(\aD f0-i\aD f1)\circ E+\frac 12(\aD f0+i\aD f1)\circ I
}

\AddEquation{a0=f01}
{
\begin{split}
\aD a0&=\aUD a00+\aUD a10i
=\frac 12(\aUD f00+\aUD f11+(\aUD f10-\aUD f01)i)
\\&=\frac 12(\aUD f00+\aUD f10i+\aUD f11-\aUD f01i)
\\&=\frac 12((\aUD f00+\aUD f10i)-(\aUD f01+\aUD f11i)i)
\end{split}
}

\AddEquation{a1=f01}
{
\begin{split}
a_1&=\aUD a01+\aUD a11i
=\frac 12(\aUD f00-\aUD f11+(\aUD f10+\aUD f01)i)
\\&=\frac 12(\aUD f00+\aUD f10i-\aUD f11+\aUD f01i)
\\&=\frac 12((\aUD f00+\aUD f10i)+(\aUD f01+\aUD f11i)i)
\end{split}
}

\AddEquation{a...= C}
{
\ShowEq{matrix a algebra C}
=
\ShowEq{matrix af algebra C}
\ShowEq{matrix f algebra C}
}

\AddEq [2]{map ao}
{
\def\Map{\aU I#1}
\edef\Tempa{#1}%
\def\Temp{E}
\ifx\Tempa\Temp%
\def\Map{E}
\fi
\def\Temp{I}
\ifx\Tempa\Temp%
\def\Map{I}
\fi
\begin{matrix}
a\bullet \Map:b\in #2\rightarrow a(\Map\circ b)\in #2&a\in #2
\end{matrix}
}

\AddEquation{f->ao f}
{
a\bullet:f\in
\ShowEq{L(A->B)}DAA
\rightarrow a\bullet f\in
\ShowEq{L(A->B)}DAA
}

\AddEquation{ao f=a f o}
{
(a\bullet f)\circ b=a(f\circ b)
}

\AddEquation{f->a* f}
{
a\star:f\in
\ShowEq{L(A->B)}DAA
\rightarrow a\star f\in
\ShowEq{L(A->B)}DAA
}

\AddEquation{a* f=(f o)a}
{
(a\star f)\circ b=(f\circ b)a
}

\AddEq [2]{map a*}
{
\def\Map{E}
\def\Temp{E}
\edef\Tempa{#1}%
\ifx\Tempa\Temp%
\else
\def\Map{\aU I#1}
\fi
\begin{matrix}
\Map\bullet a:b\in #2\rightarrow (\Map\circ b)a\in #2&a\in #2
\end{matrix}
}

\AddEq [1]{ao, H matrix 1}
{
#1_l(a)=E_l(a)\RCcirc #1
}

\AddEq [1]{a*, H matrix 1}
{
#1_r(a)=E_r(a)\RCcirc #1
}

\AddEq {aI, H matrix 2}
{
\begin{align*}
E_l(a)\RCcirc I
&=
\ShowEq{H aE Jacobian matrix}a
\RCcirc\Ei
\\
&=
\ShowEq{H a1 Jacobian matrix}a
=
I_l(a)
\end{align*}
}

\AddEq {aJ, H matrix 2}
{
\begin{align*}
E_l(a)\RCcirc J
&=
\ShowEq{H aE Jacobian matrix}a
\RCcirc\Ej
\\
&=
\ShowEq{H a2 Jacobian matrix}a
=
J_l(a)
\end{align*}
}

\AddEq{aK, H matrix 2}
{
\begin{align*}
E_l(a)\RCcirc K
&=
\ShowEq{H aE Jacobian matrix}a
\RCcirc\Ek
\\
&=
\ShowEq{H a3 Jacobian matrix}a
=
K_l(a)
\end{align*}
}

\AddEq{Ia, H matrix 2}
{
\begin{align*}
E_r(a)\RCcirc I
&=
\ShowEq{H Ea Jacobian matrix}a
\RCcirc\Ei
\\
&=
\ShowEq{H 1a Jacobian matrix}a
=
I_r(a)
\end{align*}
}

\AddEq{Ja, H matrix 2}
{
\begin{align*}
E_r(a)\RCcirc J
&=
\ShowEq{H Ea Jacobian matrix}a
\RCcirc\Ej
\\
&=
\ShowEq{H 2a Jacobian matrix}a
=
J_r(a)
\end{align*}
}

\AddEq{Ka, H matrix 2}
{
\begin{align*}
E_r(a)\RCcirc K
&=
\ShowEq{H Ea Jacobian matrix}a
\RCcirc\Ek
\\
&=
\ShowEq{H 3a Jacobian matrix}a
=
K_r(a)
\end{align*}
}

\AddEq[3]{maps of conjugation, Jacobian matrix}
{
\gdef\Map{\aU I#1}
\gdef\Letter{I}
\def\Temp{E}
\edef\Tempa{#1}%
\ifx\Tempa\Temp%
\gdef\Map{E}
\gdef\Letter{E}
\fi
\def\Temp{I}
\ifx\Tempa\Temp%
\gdef\Map{I}
\fi
\symb[\Letter-0]{\Map_{#3}(a)}{maps of conjugation, Jacobian matrix}{#1#2}
}

\AddEq[2]{a o, Jacobian matrix}
{
\ShowSymbol{maps of conjugation, Jacobian matrix}{#1#2}=
\ShowEq{#2 a#1 Jacobian matrix}a
}

\AddEq[2]{a *, Jacobian matrix}
{
\ShowSymbol{maps of conjugation, Jacobian matrix}{#1#2}=
\ShowEq{#2 #1a Jacobian matrix}a
}

\AddEq[1]{ao, C, Jacobian matrix}
{
\symb[#1-0]{#1(a)}{maps of conjugation, Jacobian matrix}{}
}

\AddEquation{aE, C, matrix}
{
\ShowSymbol{maps of conjugation, Jacobian matrix}{}=
\begin{pmatrix}
\aU a0&-a^{\gi 1}\\
a^{\gi 1}&\hphantom{-}\aU a0
\end{pmatrix}
}

\AddEquation{aI, C, matrix}
{
\ShowSymbol{maps of conjugation, Jacobian matrix}{}=
\begin{pmatrix}
\aU a0&\hphantom{-}a^{\gi 1}\\
a^{\gi 1}&-\aU a0
\end{pmatrix}
}

\AddEquation{f in L(A,A), 2, nonassociative algebra}
{
f
=a^{k\cdot \gi{ij}}(\aU ei\otimes\aU ej)\circ F_k
=
a^{k\cdot \gi{ij}}(\aU eiI_k)\aU ej
}

\AddEquation{system 1 of linear equations, quaternion}
{
\begin{split}
i\aU x1j+j\aU x2k&=i+j
\\
k\aU x1i+i\aU x2j&=i-j
\end{split}
}

\AddEquation{system 1 of linear equations, quaternion, 1}
{
\begin{split}
(i\otimes j)\circ\aU x1+(j\otimes k)\circ\aU x2&=i+j
\\
(k\otimes i)\circ\aU x1+(i\otimes j)\circ\aU x2&=i-j
\end{split}
}

\AddEquation{system 1 of linear equations, quaternion, 2}
{
\begin{pmatrix}
i\otimes j&j\otimes k\\
k\otimes i&i\otimes j
\end{pmatrix}
\RCcirc
\begin{pmatrix}
c\aU x1\\\aU x2
\end{pmatrix}
=
\begin{pmatrix}
i+j\\i-j
\end{pmatrix}
}

\AddEq{system 1 of linear equations, quaternion, 3}
{
\[
a=
\begin{pmatrix}
i\otimes j&j\otimes k\\
k\otimes i&i\otimes j
\end{pmatrix}
\]
}

\AddEq{system 1 of linear equations, quaternion, 4}
{
\ShowEq{system 1 of linear equations, quaternion, 4 1 1}
\ShowEq{system 1 of linear equations, quaternion, 4 1 2}
\ShowEq{system 1 of linear equations, quaternion, 4 2 1}
\ShowEq{system 1 of linear equations, quaternion, 4 2 2}
}

\AddEquation{system 1 of linear equations, quaternion, 4 1 1}
{
\begin{split}
\aUD{\det\left(a,\RCcirc\right)}11&=\aUD a11
-\aUD a1{[1]}\RCcirc
(\aUD a{[1]}{[1]})^{-1\RCcirc}\RCcirc
\aUD a{[1]}1
=\aUD a11
-\aUD a12\circ
(\aUD a22)^{-1}\circ
\aUD a21
\\
&=i\otimes j
-(j\otimes k)\circ
(i\otimes j)^{-1}\circ
(k\otimes i)
\\
&=i\otimes j
-(j\otimes k)\circ
(i\otimes j)\circ
(k\otimes i)
=i\otimes j
-(-1)\otimes (-1)
\\
&=i\otimes j
-1\otimes 1
\end{split}
}

\AddEquation{system 1 of linear equations, quaternion, 4 1 2}
{
\begin{split}
\aUD{\det\left(a,\RCcirc\right)}12&=\aUD a12
-\aUD a1{[2]}\RCcirc
(\aUD a{[1]}{[2]})^{-1\RCcirc}\RCcirc
\aUD a{[1]}2
=\aUD a12
-\aUD a11\circ
(\aUD a21)^{-1}\circ
\aUD a22
\\
&=j\otimes k
-(i\otimes j)\circ
(k\otimes i)^{-1}\circ
(i\otimes j)
\\
&=j\otimes k
-(i\otimes j)\circ
(k\otimes i)\circ
(i\otimes j)
=j\otimes k
-(-k)\otimes i
\\
&=j\otimes k
+k\otimes i
\end{split}
}

\AddEquation{system 1 of linear equations, quaternion, 4 2 1}
{
\begin{split}
\aUD{\det\left(a,\RCcirc\right)}21&=\aUD a21
-\aUD a2{[1]}\RCcirc
(\aUD a{[2]}{[1]})^{-1\RCcirc}\RCcirc
\aUD a{[2]}1
=\aUD a21
-\aUD a22\circ
(\aUD a12)^{-1}\circ
\aUD a11
\\
&=k\otimes i
-(i\otimes j)\circ
(j\otimes k)^{-1}\circ
(i\otimes j)
\\
&=k\otimes i
-(i\otimes j)\circ
(j\otimes k)\circ
(i\otimes j)
=k\otimes i
-(-j)\otimes k
\\
&=k\otimes i
+j\otimes k
\end{split}
}

\AddEquation{system 1 of linear equations, quaternion, 4 2 2}
{
\begin{split}
\aUD{\det\left(a,\RCcirc\right)}22&=\aUD a22
-\aUD a2{[2]}\RCcirc
(\aUD a{[2]}{[2]})^{-1\RCcirc}\RCcirc
\aUD a{[2]}2
=\aUD a22
-\aUD a21\circ
(\aUD a11)^{-1}\circ
\aUD a12
\\
&=i\otimes j
-(k\otimes i)\circ
(i\otimes j)^{-1}\circ
(j\otimes k)
\\
&=i\otimes j
-(k\otimes i)\circ
(i\otimes j)\circ
(j\otimes k)
=i\otimes j
-(-1)\otimes (-1)
\\
&=i\otimes j
-1\otimes 1
\end{split}
}

\AddEquation{fg=1 e o e}
{
\begin{split}
\aU f{ij}\aUD C0{ik}\aUD C0{lj}\aU g{kl}&=1\\
\aU f{ij}\aUD Cp{ik}\aUD Cq{lj}\aU g{kl}&=0\ \ \ p+q>0
\end{split}
}

\AddEq{g o h=1}
{
\[g\circ h=\aD e0\otimes\aD e0\]
}

\AddEq{f1=...}
{
\[f_1=-1\otimes 1+i\otimes j\]
}

\AddEq{f2=...}
{
\[f_2=k\otimes i+j\otimes k\]
}

\AddEq[2]{g=f-1}
{
$g_{#1}=f_{#1}^{-1}$#2
}

\AddEq[3]{row set 1n}
{
$\aD {#1}1$, ..., $\aD {#1}{#2}$#3
}

\AddEq[3]{column set 1n}
{
$\aU {#1}1$, ..., $\aU {#1}{#2}$#3
}

\AddEq[3]{row() 1n}
{
$#1(\gi 1)=\aD {#1}1$, ..., $#1(\gi{#2})=\aD {#1}{#2}$#3
}

\AddEq[3]{col() 1n}
{
$#1(\gi 1)=\aU {#1}1$, ..., $#1(\gi{#2})=\aU {#1}{#2}$#3
}

\AddEq[5]{system of linear equations left column row}
{
\begin{split}
\aU {#1}1\aUD {#2}11+...+\aU {#1}{#4}\aUD {#2}1{#4}&=\aU {#3}1
\\...&=...\\
\aU {#1}1\aUD {#2}{#5}1+...+\aU {#1}{#4}\aUD {#2}{#5}{#4}&=\aU {#3}{#5}
\end{split}
}

\AddEq[5]{system of linear equations linear map}
{
\begin{split}
\aUD {#2}1{11}\aU {#1}1\aUD {#2}1{12}+...+\aUD {#2}1{{#4}1}\aU {#1}{#4}\aUD {#2}1{{#4}2}
&=\aU {#3}1
\\...&=...\\
\aUD {#2}{#5}{11}\aU {#1}1\aUD {#2}{#5}{12}+...+\aUD {#2}{#5}{{#4}1}\aU {#1}{#4}\aUD {#2}{#5}{{#4}2}
&=\aU {#3}{#5}
\end{split}
}

\AddEquation{fg=1 e o e g1}
{
\left\{
\begin{matrix}
-\aUD C0{0k}\aUD C0{l0}\aU g{kl}+
\aUD C0{1k}\aUD C0{l2}\aU g{kl}=1\\
-\aUD C0{0k}\aUD C1{l0}\aU g{kl}+
\aUD C0{1k}\aUD C1{l2}\aU g{kl}=0\\
-\aUD C0{0k}\aUD C2{l0}\aU g{kl}+
\aUD C0{1k}\aUD C2{l2}\aU g{kl}=0\\
-\aUD C0{0k}\aUD C3{l0}\aU g{kl}+
\aUD C0{1k}\aUD C3{l2}\aU g{kl}=0\\
-\aUD C1{0k}\aUD C0{l0}\aU g{kl}+
\aUD C1{1k}\aUD C0{l2}\aU g{kl}=0\\
-\aUD C1{0k}\aUD C1{l0}\aU g{kl}+
\aUD C1{1k}\aUD C1{l2}\aU g{kl}=0\\
-\aUD C1{0k}\aUD C2{l0}\aU g{kl}+
\aUD C1{1k}\aUD C2{l2}\aU g{kl}=0\\
-\aUD C1{0k}\aUD C3{l0}\aU g{kl}+
\aUD C1{1k}\aUD C3{l2}\aU g{kl}=0
\end{matrix}
\ \ \ \ \,
\begin{matrix}
-\aUD C2{0k}\aUD C0{l0}\aU g{kl}+
\aUD C2{1k}\aUD C0{l2}\aU g{kl}=0\\
-\aUD C2{0k}\aUD C1{l0}\aU g{kl}+
\aUD C2{1k}\aUD C1{l2}\aU g{kl}=0\\
-\aUD C2{0k}\aUD C2{l0}\aU g{kl}+
\aUD C2{1k}\aUD C2{l2}\aU g{kl}=0\\
-\aUD C2{0k}\aUD C3{l0}\aU g{kl}+
\aUD C2{1k}\aUD C3{l2}\aU g{kl}=0\\
-\aUD C3{0k}\aUD C0{l0}\aU g{kl}+
\aUD C3{1k}\aUD C0{l2}\aU g{kl}=0\\
-\aUD C3{0k}\aUD C1{l0}\aU g{kl}+
\aUD C3{1k}\aUD C1{l2}\aU g{kl}=0\\
-\aUD C3{0k}\aUD C2{l0}\aU g{kl}+
\aUD C3{1k}\aUD C2{l2}\aU g{kl}=0\\
-\aUD C3{0k}\aUD C3{l0}\aU g{kl}+
\aUD C3{1k}\aUD C3{l2}\aU g{kl}=0
\end{matrix}
\right.
}

\AddEquation{fg=1 e o e g1 1}
{
\left\{
\begin{matrix}
\BlueText{-\aU g{00}+\aU g{12}=1}\\
-\aU g{01}+\aU g{13}=0\\
-\aU g{02}-\aU g{10}=0\\
-\aU g{03}+\aU g{11}=0\\
-\aU g{10}-\aU g{02}=0\\
-\aU g{11}-\aU g{03}=0\\
\BlueText{-\aU g{12}+\aU g{00}=0}\\
-\aU g{13}+\aU g{01}=0
\end{matrix}
\ \ \ \ \,
\begin{matrix}
-\aUD C2{0k}\aUD C0{l0}\aU g{kl}+
\aUD C2{1k}\aUD C0{l2}\aU g{kl}=0\\
-\aUD C2{0k}\aUD C1{l0}\aU g{kl}+
\aUD C2{1k}\aUD C1{l2}\aU g{kl}=0\\
-\aUD C2{0k}\aUD C2{l0}\aU g{kl}+
\aUD C2{1k}\aUD C2{l2}\aU g{kl}=0\\
-\aUD C2{0k}\aUD C3{l0}\aU g{kl}+
\aUD C2{1k}\aUD C3{l2}\aU g{kl}=0\\
-\aUD C3{0k}\aUD C0{l0}\aU g{kl}+
\aUD C3{1k}\aUD C0{l2}\aU g{kl}=0\\
-\aUD C3{0k}\aUD C1{l0}\aU g{kl}+
\aUD C3{1k}\aUD C1{l2}\aU g{kl}=0\\
-\aUD C3{0k}\aUD C2{l0}\aU g{kl}+
\aUD C3{1k}\aUD C2{l2}\aU g{kl}=0\\
-\aUD C3{0k}\aUD C3{l0}\aU g{kl}+
\aUD C3{1k}\aUD C3{l2}\aU g{kl}=0
\end{matrix}
\right.
}

\AddEquation{fg=1 e o e g2}
{
\left\{
\begin{matrix}
\aUD C0{3k}\aUD C0{l1}\aU g{kl}+
\aUD C0{2k}\aUD C0{l3}\aU g{kl}=1\\
\aUD C0{3k}\aUD C1{l1}\aU g{kl}+
\aUD C0{2k}\aUD C1{l3}\aU g{kl}=0\\
\aUD C0{3k}\aUD C2{l1}\aU g{kl}+
\aUD C0{2k}\aUD C2{l3}\aU g{kl}=0\\
\aUD C0{3k}\aUD C3{l1}\aU g{kl}+
\aUD C0{2k}\aUD C3{l3}\aU g{kl}=0\\
\aUD C1{3k}\aUD C0{l1}\aU g{kl}+
\aUD C1{2k}\aUD C0{l3}\aU g{kl}=0\\
\aUD C1{3k}\aUD C1{l1}\aU g{kl}+
\aUD C1{2k}\aUD C1{l3}\aU g{kl}=0\\
\aUD C1{3k}\aUD C2{l1}\aU g{kl}+
\aUD C1{2k}\aUD C2{l3}\aU g{kl}=0\\
\aUD C1{3k}\aUD C3{l1}\aU g{kl}+
\aUD C1{2k}\aUD C3{l3}\aU g{kl}=0
\end{matrix}
\ \ \ \ \,
\begin{matrix}
\aUD C2{3k}\aUD C0{l1}\aU g{kl}+
\aUD C2{2k}\aUD C0{l3}\aU g{kl}=0\\
\aUD C2{3k}\aUD C1{l1}\aU g{kl}+
\aUD C2{2k}\aUD C1{l3}\aU g{kl}=0\\
\aUD C2{3k}\aUD C2{l1}\aU g{kl}+
\aUD C2{2k}\aUD C2{l3}\aU g{kl}=0\\
\aUD C2{3k}\aUD C3{l1}\aU g{kl}+
\aUD C2{2k}\aUD C3{l3}\aU g{kl}=0\\
\aUD C3{3k}\aUD C0{l1}\aU g{kl}+
\aUD C3{2k}\aUD C0{l3}\aU g{kl}=0\\
\aUD C3{3k}\aUD C1{l1}\aU g{kl}+
\aUD C3{2k}\aUD C1{l3}\aU g{kl}=0\\
\aUD C3{3k}\aUD C2{l1}\aU g{kl}+
\aUD C3{2k}\aUD C2{l3}\aU g{kl}=0\\
\aUD C3{3k}\aUD C3{l1}\aU g{kl}+
\aUD C3{2k}\aUD C3{l3}\aU g{kl}=0
\end{matrix}
\right.
}

\AddEquation{fg=1 e o e g2 1}
{
\left\{
\begin{matrix}
\BlueText{\aU g{31}+\aU g{23}=1}\\
-\aU g{30}-\aU g{22}=0\\
-\aU g{33}+\aU g{21}=0\\
\aU g{32}-\aU g{20}=0\\
\aU g{21}-\aU g{33}=0\\
-\aU g{20}+\aU g{32}=0\\
\BlueText{-\aU g{23}-\aU g{31}=0}\\
\aU g{22}+\aU g{30}=0
\end{matrix}
\ \ \ \ \,
\begin{matrix}
\aUD C2{3k}\aUD C0{l1}\aU g{kl}+
\aUD C2{2k}\aUD C0{l3}\aU g{kl}=0\\
\aUD C2{3k}\aUD C1{l1}\aU g{kl}+
\aUD C2{2k}\aUD C1{l3}\aU g{kl}=0\\
\aUD C2{3k}\aUD C2{l1}\aU g{kl}+
\aUD C2{2k}\aUD C2{l3}\aU g{kl}=0\\
\aUD C2{3k}\aUD C3{l1}\aU g{kl}+
\aUD C2{2k}\aUD C3{l3}\aU g{kl}=0\\
\aUD C3{3k}\aUD C0{l1}\aU g{kl}+
\aUD C3{2k}\aUD C0{l3}\aU g{kl}=0\\
\aUD C3{3k}\aUD C1{l1}\aU g{kl}+
\aUD C3{2k}\aUD C1{l3}\aU g{kl}=0\\
\aUD C3{3k}\aUD C2{l1}\aU g{kl}+
\aUD C3{2k}\aUD C2{l3}\aU g{kl}=0\\
\aUD C3{3k}\aUD C3{l1}\aU g{kl}+
\aUD C3{2k}\aUD C3{l3}\aU g{kl}=0
\end{matrix}
\right.
}

\AddEq[1]{Map of Conjugation}
{
\def\Temp{0}%
\def\Tempa{#1}%
\ifx\Tempa\Temp%
\ensuremath{E}
\else
\ensuremath{\aU I{#1}}
\fi
}

\AddEquation{H.fUD<-fUU}
{
\begin{split}
&
\ShowEq{quaternion matrix fUD}
\\=&
\ShowEq{quaternion matrix fUD<-fUU}
\ShowEq{quaternion matrix fUU}f
\end{split}
}

\AddEq[1]{H.fUD->fUU I}
{
\begin{split}
&
\ShowEq{quaternion matrix fUU}f
\\=&
\ShowEq{quaternion matrix fUD->fUU}
\ShowEq{H. matrix fUD I#1}
\end{split}
}

\AddEquation{H.fUD->fUU}
{
\begin{split}
&
\ShowEq{quaternion matrix fUU}f
\\=&
\ShowEq{quaternion matrix fUD->fUU}
\ShowEq{quaternion matrix fUD}
\end{split}
}

\AddEq{quaternion, 3, 1}
{
\[
\ShowEq{quaternion matrix fUD<-fUU}^{-1}
=
\ShowEq{quaternion matrix fUD->fUU}
\]
}

\AddEq{D in Z(A)}
{
\(D\subset Z(A)\).
\labelItem{D in Z(A)}
}

\AddEquation{H.fUD<-fUU f1 1}
{
\begin{split}
&
\ShowEq{quaternion matrix fUD}
\\=&
\ShowEq{quaternion matrix fUD<-fUU}
\left(
\begin{array}{rrrr}
-1&0&0&0\\
0&0&0&-1\\
0&0&0&0\\
0&0&0&0
\end{array}
\right)
\\=&
\left(
\begin{array}{rrrr}
-1&0&0&1\\
-1&0&0&1\\
-1&0&0&-1\\
-1&0&0&-1
\end{array}
\right)
\end{split}
}

\AddEquation{H.fUD<-fUU f2 1}
{
\begin{split}
&
\ShowEq{quaternion matrix fUD}
\\=&
\ShowEq{quaternion matrix fUD<-fUU}
\left(
\begin{array}{rrrr}
0&0&0&0\\
0&0&0&0\\
0&-1&0&0\\
0&0&-1&0
\end{array}
\right)
\\=&
\left(
\begin{array}{rrrr}
0&1&1&0\\
0&-1&-1&0\\
0&1&-1&0\\
0&-1&1&0
\end{array}
\right)
\end{split}
}

\AddEquation{matrix f1}
{
f_1=
\left(
\begin{array}{rrrr}
-1&0&0&1\\
0&-1&-1&0\\
0&-1&-1&0\\
1&0&0&-1
\end{array}
\right)
}

\AddEquation{matrix f2}
{
f_2=
\left(
\begin{array}{rrrr}
0&1&1&0\\
1&0&0&-1\\
1&0&0&-1\\
0&-1&-1&0
\end{array}
\right)
}

\AddEq[1]{f=...ei o ek}
{
#1=\aU {#1}{ij}\aD ei\otimes\aD ej
}

\AddEquation{fg=...ei o ek}
{
f\circ g=\aU {(f\circ g)}{pq}\aD ep\otimes\aD eq
}

\AddEquation{fg=CCfg}
{
\aU {(f\circ g)}{pq}=\aU f{ij}\aU g{kl}\aUD Cp{ik}\aUD Cq{lj}
}

\AddEquation{fg=...e o e}
{
f\circ g=\aU f{ij}\aU g{kl}(\aD ei\aD ek)\otimes(\aD el\aD ej)
}

\AddEquation{fg=...C e o e}
{
f\circ g=\aU f{ij}\aU g{kl}\aUD Cp{ik}\aUD Cq{lj}\aD ep\otimes\aD eq
}

\AddEq[2]{f in AoA}
{
$#1\in A\otimes A$#2
}

\AddEq[2]{texorpdf L(D;A->A)}
{%
\texorpdfstring{%
\ShowEq{L(A->B)}#1#2#2{}%
}{L(#1;#2->#2)}%
}

\AddEq{I0=E}
{
$\aU I0=E$.
}

\AddEquation{a > a bullet}
{
\xymatrix
{
\bullet:A\ar[r]|(.4){*}&
\ShowEq{L(A->B)}DAA
}
\ \ \ \ \ \,
a\bullet:f\rightarrow a\bullet f
}

\AddEq{coordinates of map A->A, 2}
{
\aUD fkl=\aU {f^{k\cdot}_{}{}}{ij}
\aUD {F_{k\cdot}^{}{}}ml
\aUD Cp{im}\aUD Ck{pj}
}

\AddEq [2]{coordinates of map A1 A2, FoG, associative algebra}
{
\def\Cdot{}
\def\Temp{-}%
\edef\Tempa{#2}%
\ifx\Tempa\Temp%
\else
\def\Cdot{#2\cdot}
\fi
#1^{\gik}_{\gil}=#1^{k\cdot\gi{ij}}
F_{k\cdot}^{}{}^{\gim}_{\gi r}G^{\gi r}_{\gil}
C_{\Cdot}{}^{\gi p}_{\gi{im}}C_{\Cdot}{}^{\gik}_{\gi{pj}}
}

\AddEq [2]{coordinates of map A->B}
{
\def\Cdot{}
\def\Temp{-}%
\edef\Tempa{#2}%
\ifx\Tempa\Temp%
\else
\def\Cdot{#2\cdot}
\fi
#1^{\gik}_{\gil}=#1^{k\cdot\gi{ij}}
F_{k\cdot}^{}{}^{\gim}_{\gi r}G^{\gi r}_{\gil}
C_{\Cdot}{}^{\gi p}_{\gi{im}}C_{\Cdot}{}^{\gik}_{\gi{pj}}
}

\AddEquation{a...= H}
{
\begin{split}
&\,
\ShowEq{matrix a algebra H}
\\ = &\,
\ShowEq{matrix af algebra H}
\ShowEq{matrix f algebra H}
\end{split}
}

\AddEquation{a0=f03}
{
\begin{split}
\aD a0&=\aUD a00+\aUD a10i+\aUD a20j+\aUD a30k
\\&=\frac 12
(-\aUD f00+\aUD f11+\aUD f22+\aUD f33
+(-\aUD f10-\aUD f01+\aUD f32-\aUD f23)i
\\&+(-\aUD f20-\aUD f31-\aUD f02+\aUD f13)j
+(-\aUD f30+\aUD f21-\aUD f12-\aUD f03)k)
\\&=\frac 12
((-\aUD f00-\aUD f10i-\aUD f20j-\aUD f30k)
+(-\aUD f01i+\aUD f11+\aUD f21k-\aUD f31j)
\\&+(-\aUD f02j-\aUD f12k+\aUD f22+\aUD f32i)
+(-\aUD f03k+\aUD f13j-\aUD f23i+\aUD f33))
\\&=\frac 12
(-(\aUD f00+\aUD f10i+\aUD f20j+\aUD f30k)
-(\aUD f01+\aUD f11i+\aUD f21j+\aUD f31k)i
\\&-(\aUD f02+\aUD f12i+\aUD f22j+\aUD f32k)j
-(\aUD f03+\aUD f13i+\aUD f23j+\aUD f33k)k)
\end{split}
}

\AddEquation{a1=f03}
{
\begin{split}
\aD a1&=\aUD a01+\aUD a11i+\aUD a21j+\aUD a31k
\\&=\frac 12
(\aUD f00-\aUD f11+(\aUD f10+\aUD f01)i+(\aUD f20+\aUD f31)j+(\aUD f30-\aUD f21)k)
\\&=\frac 12
((\aUD f00+\aUD f10i+\aUD f20j+\aUD f30k)+(\aUD f01i-\aUD f11-\aUD f21k+\aUD f31j))
\\&=\frac 12
((\aUD f00+\aUD f10i+\aUD f20j+\aUD f30k)+(\aUD f01+\aUD f11i+\aUD f21j+\aUD f31k)i)
\end{split}
}

\AddEquation{a2=f03}
{
\begin{split}
\aD a2&=\aUD a02+\aUD a12i+\aUD a22j+\aUD a32k
\\&=\frac 12
(\aUD f00-\aUD f22+(\aUD f10-\aUD f32)i+(\aUD f20+\aUD f02)j+(\aUD f30+\aUD f12)k)
\\&=\frac 12
(\aUD f00-\aUD f22+\aUD f10i-\aUD f32i+\aUD f20j+\aUD f02j+\aUD f30k+\aUD f12k)
\\&=\frac 12
((\aUD f00+\aUD f10i+\aUD f20j+\aUD f30k)+(\aUD f02j+\aUD f12k-\aUD f22-\aUD f32i))
\\&=\frac 12
((\aUD f00+\aUD f10i+\aUD f20j+\aUD f30k)+(\aUD f02+\aUD f12i+\aUD f22j+\aUD f32k)j)
\end{split}
}

\AddEquation{a3=f03}
{
\begin{split}
\aD a3&=\aUD a03+\aUD a13i+\aUD a23j+\aUD a33k
\\&=\frac 12
(\aUD f00-\aUD f33+(\aUD f10+\aUD f23)i+(\aUD f20-\aUD f13)j+(\aUD f30+\aUD f03)k)
\\&=\frac 12
((\aUD f00+\aUD f10i+\aUD f30k+\aUD f20j)+(\aUD f03k-\aUD f13j+\aUD f23i-\aUD f33))
\\&=\frac 12
((\aUD f00+\aUD f10i+\aUD f30k+\aUD f20j)+(\aUD f03+\aUD f13i+\aUD f23j+\aUD f33k)k)
\end{split}
}

\AddEq{a0=f07 1}
{
\begin{align*}
\aD a0&=\aUD a00\aD e0+\aUD a10\aD e1+\aUD a20\aD e2+\aUD a30\aD e3
+\aUD a40\aD e4+\aUD a50\aD e5+\aUD a60\aD e6+\aUD a70\aD e7
\\&=\frac 12
((-5\aUD f00+\aUD f11+\aUD f22+\aUD f33+\aUD f44+\aUD f55+\aUD f66+\aUD f77)\aD e0
\\&+(-5\aUD f10-\aUD f01+\aUD f32-\aUD f23+\aUD f54-\aUD f45-\aUD f76+\aUD f67)\aD e1
\\&+(-5\aUD f20-\aUD f31-\aUD f02+\aUD f13+\aUD f64+\aUD f75-\aUD f46-\aUD f57)\aD e2
\\&+(-5\aUD f30+\aUD f21-\aUD f12-\aUD f03+\aUD f74-\aUD f65+\aUD f56-\aUD f47)\aD e3
\\&+(-5\aUD f40-\aUD f51-\aUD f62-\aUD f73-\aUD f04+\aUD f15+\aUD f26+\aUD f37)\aD e4
\\&+(-5\aUD f50+\aUD f41-\aUD f72+\aUD f63-\aUD f14-\aUD f05-\aUD f36+\aUD f27)\aD e5
\\&+(-5\aUD f60+\aUD f71+\aUD f42-\aUD f53-\aUD f24+\aUD f35-\aUD f06-\aUD f17)\aD e6
\\&+(-5\aUD f70-\aUD f61+\aUD f52+\aUD f43-\aUD f34-\aUD f25+\aUD f16-\aUD f07)\aD e7)
\end{align*}
}

\AddEq{a0=f07 2}
{
\begin{align*}
\aD a0&=\frac 12
(-5\aUD f00\aD e0+\aUD f11\aD e0+\aUD f22\aD e0+\aUD f33\aD e0
+\aUD f44\aD e0+\aUD f55\aD e0+\aUD f66\aD e0+\aUD f77\aD e0
\\&-5\aUD f10\aD e1-\aUD f01\aD e1+\aUD f32\aD e1-\aUD f23\aD e1
+\aUD f54\aD e1-\aUD f45\aD e1-\aUD f76\aD e1+\aUD f67\aD e1
\\&-5\aUD f20\aD e2-\aUD f31\aD e2-\aUD f02\aD e2+\aUD f13\aD e2
+\aUD f64\aD e2+\aUD f75\aD e2-\aUD f46\aD e2-\aUD f57\aD e2
\\&-5\aUD f30\aD e3+\aUD f21\aD e3-\aUD f12\aD e3-\aUD f03\aD e3
+\aUD f74\aD e3-\aUD f65\aD e3+\aUD f56\aD e3-\aUD f47\aD e3
\\&-5\aUD f40\aD e4-\aUD f51\aD e4-\aUD f62\aD e4-\aUD f73\aD e4
-\aUD f04\aD e4+\aUD f15\aD e4+\aUD f26\aD e4+\aUD f37\aD e4
\\&-5\aUD f50\aD e5+\aUD f41\aD e5-\aUD f72\aD e5+\aUD f63\aD e5
-\aUD f14\aD e5-\aUD f05\aD e5-\aUD f36\aD e5+\aUD f27\aD e5
\\&-5\aUD f60\aD e6+\aUD f71\aD e6+\aUD f42\aD e6-\aUD f53\aD e6
-\aUD f24\aD e6+\aUD f35\aD e6-\aUD f06\aD e6-\aUD f17\aD e6
\\&-5\aUD f70\aD e7-\aUD f61\aD e7+\aUD f52\aD e7+\aUD f43\aD e7
-\aUD f34\aD e7-\aUD f25\aD e7+\aUD f16\aD e7-\aUD f07\aD e7)
\end{align*}
}

\AddEq{a0=f07 3}
{
\begin{align*}
\aD a0&=\frac 12
(-5(\aUD f00\aD e0+\aUD f10\aD e1+\aUD f20\aD e2+\aUD f30\aD e3
+\aUD f40\aD e4+\aUD f50\aD e5+\aUD f60\aD e6+\aUD f70\aD e7)
\\&+(-\aUD f01\aD e1+\aUD f11\aD e0+\aUD f21\aD e3-\aUD f31\aD e2
+\aUD f41\aD e5-\aUD f51\aD e4-\aUD f61\aD e7+\aUD f71\aD e6)
\\&+(-\aUD f02\aD e2-\aUD f12\aD e3+\aUD f22\aD e0+\aUD f32\aD e1
+\aUD f42\aD e6+\aUD f52\aD e7-\aUD f62\aD e4-\aUD f72\aD e5)
\\&+(-\aUD f03\aD e3+\aUD f13\aD e2-\aUD f23\aD e1+\aUD f33\aD e0
+\aUD f43\aD e7-\aUD f53\aD e6+\aUD f63\aD e5-\aUD f73\aD e4)
\\&+(-\aUD f04\aD e4-\aUD f14\aD e5-\aUD f24\aD e6-\aUD f34\aD e7
+\aUD f44\aD e0+\aUD f54\aD e1+\aUD f64\aD e2+\aUD f74\aD e3)
\\&+(-\aUD f05\aD e5+\aUD f15\aD e4-\aUD f25\aD e7+\aUD f35\aD e6
-\aUD f45\aD e1+\aUD f55\aD e0-\aUD f65\aD e3+\aUD f75\aD e2)
\\&+(-\aUD f06\aD e6+\aUD f16\aD e7+\aUD f26\aD e4-\aUD f36\aD e5
-\aUD f46\aD e2+\aUD f56\aD e3+\aUD f66\aD e0-\aUD f76\aD e1)
\\&+(-\aUD f07\aD e7-\aUD f17\aD e6+\aUD f27\aD e5+\aUD f37\aD e4
-\aUD f47\aD e3-\aUD f57\aD e2+\aUD f67\aD e1+\aUD f77\aD e0)
\end{align*}
}

\AddEquation{a0=f07}
{
\begin{split}
\aD a0&=-\frac 12
(5(\aUD f00\aD e0+\aUD f10\aD e1+\aUD f20\aD e2+\aUD f30\aD e3
+\aUD f40\aD e4+\aUD f50\aD e5+\aUD f60\aD e6+\aUD f70\aD e7)\aD e0
\\&+(\aUD f01\aD e0+\aUD f11\aD e0+\aUD f21\aD e2+\aUD f31\aD e3
+\aUD f41\aD e4+\aUD f51\aD e5+\aUD f61\aD e6+\aUD f71\aD e7)\aD e1
\\&+(\aUD f02\aD e0+\aUD f12\aD e1+\aUD f22\aD e2+\aUD f32\aD e3
+\aUD f42\aD e4+\aUD f52\aD e5+\aUD f62\aD e6+\aUD f72\aD e7)\aD e2
\\&+(\aUD f03\aD e0+\aUD f13\aD e2+\aUD f23\aD e2+\aUD f33\aD e0
+\aUD f43\aD e4+\aUD f53\aD e5+\aUD f63\aD e6+\aUD f73\aD e7)\aD e3
\\&+(\aUD f04\aD e0+\aUD f14\aD e1+\aUD f24\aD e2+\aUD f34\aD e3
+\aUD f44\aD e4+\aUD f54\aD e5+\aUD f64\aD e6+\aUD f74\aD e7)\aD e4
\\&+(\aUD f05\aD e0+\aUD f15\aD e1+\aUD f25\aD e2+\aUD f35\aD e3
+\aUD f45\aD e4+\aUD f55\aD e5+\aUD f65\aD e3+\aUD f75\aD e7)\aD e5
\\&+(\aUD f06\aD e6+\aUD f16\aD e1+\aUD f26\aD e2+\aUD f36\aD e3
+\aUD f46\aD e4+\aUD f56\aD e5+\aUD f66\aD e6+\aUD f76\aD e7)\aD e6
\\&+(\aUD f07\aD e0+\aUD f17\aD e1+\aUD f27\aD e2+\aUD f37\aD e3
+\aUD f47\aD e4+\aUD f57\aD e5+\aUD f67\aD e6+\aUD f77\aD e7)\aD e7
\end{split}
}

\AddEquation{a1=f07}
{
\begin{split}
\aD a1&=\aUD a01\aD e0+\aUD a11\aD e1+\aUD a21\aD e2+\aUD a31\aD e3
+\aUD a41\aD e4+\aUD a51\aD e5+\aUD a61\aD e6+\aUD a71\aD e7
\\&=\frac 12
((\aUD f00-\aUD f11)\aD e0+(\aUD f10+\aUD f01)\aD e1
+(\aUD f20+\aUD f31)\aD e2+(\aUD f30-\aUD f21)\aD e3
\\ &
+(\aUD f40+\aUD f51)\aD e4+(\aUD f50-\aUD f41)\aD e5
+(\aUD f60-\aUD f71)\aD e6+(\aUD f70+\aUD f61)\aD e7)
\\&=\frac 12
(\aUD f00\aD e0-\aUD f11\aD e0+\aUD f10\aD e1+\aUD f01\aD e1
+\aUD f20\aD e2+\aUD f31\aD e2+\aUD f30\aD e3-\aUD f21\aD e3
\\ &
+\aUD f40\aD e4+\aUD f51\aD e4+\aUD f50\aD e5-\aUD f41\aD e5
+\aUD f60\aD e6-\aUD f71\aD e6+\aUD f70\aD e7+\aUD f61\aD e7)
\\&=\frac 12
((\aUD f00\aD e0+\aUD f10\aD e1+\aUD f20\aD e2+\aUD f30\aD e3+\aUD f40\aD e4
+\aUD f50\aD e5+\aUD f60\aD e6+\aUD f70\aD e7)
\\ &
+(\aUD f01\aD e1-\aUD f11\aD e0-\aUD f21\aD e3+\aUD f31\aD e2
-\aUD f41\aD e5+\aUD f51\aD e4+\aUD f61\aD e7-\aUD f71\aD e6))
\\&=\frac 12
((\aUD f00\aD e0+\aUD f10\aD e1+\aUD f20\aD e2+\aUD f30\aD e3+\aUD f40\aD e4
+\aUD f50\aD e5+\aUD f60\aD e6+\aUD f70\aD e7)
\\ &
+(\aUD f01\aD e0+\aUD f11\aD e1+\aUD f21\aD e2+\aUD f31\aD e3
+\aUD f41\aD e4+\aUD f51\aD e5+\aUD f61\aD e6+\aUD f71\aD e7)\aD e1)
\end{split}
}

\AddEquation{a2=f07}
{
\begin{split}
\aD a2&=\aUD a02\aD e0+\aUD a12\aD e1+\aUD a22\aD e2+\aUD a32\aD e3
+\aUD a42\aD e4+\aUD a52\aD e5+\aUD a62\aD e6+\aUD a72\aD e7
\\&=\frac 12
((\aUD f00-\aUD f22)\aD e0+(\aUD f10-\aUD f32)\aD e1
+(\aUD f20+\aUD f02)\aD e2+(\aUD f30+\aUD f12)\aD e3
\\ &
+(\aUD f40+\aUD f62)\aD e4+(\aUD f50+\aUD f72)\aD e5
+(\aUD f60-\aUD f42)\aD e6+(\aUD f70-\aUD f52)\aD e7)
\\&=\frac 12
(\aUD f00\aD e0-\aUD f22\aD e0+\aUD f10\aD e1-\aUD f32\aD e1
+\aUD f20\aD e2+\aUD f02\aD e2+\aUD f30\aD e3+\aUD f12\aD e3
\\ &
+\aUD f40\aD e4+\aUD f62\aD e4+\aUD f50\aD e5+\aUD f72\aD e5
+\aUD f60\aD e6-\aUD f42\aD e6+\aUD f70\aD e7-\aUD f52\aD e7)
\\&=\frac 12
((\aUD f00\aD e0+\aUD f10\aD e1+\aUD f20\aD e2+\aUD f30\aD e3+\aUD f40\aD e4
+\aUD f50\aD e5+\aUD f60\aD e6+\aUD f70\aD e7)
\\ &
+(\aUD f02\aD e2+\aUD f12\aD e3-\aUD f32\aD e1-\aUD f22\aD e0
-\aUD f42\aD e6-\aUD f52\aD e7+\aUD f62\aD e4+\aUD f72\aD e5))
\\&=\frac 12
((\aUD f00\aD e0+\aUD f10\aD e1+\aUD f20\aD e2+\aUD f30\aD e3+\aUD f40\aD e4
+\aUD f50\aD e5+\aUD f60\aD e6+\aUD f70\aD e7)
\\ &
+(\aUD f02\aD e0+\aUD f12\aD e1+\aUD f22\aD e2+\aUD f32\aD e3
+\aUD f42\aD e4+\aUD f52\aD e5+\aUD f62\aD e6+\aUD f72\aD e7)\aD e2)
\end{split}
}

\AddEquation{a3=f07}
{
\begin{split}
\aD a3&=\aUD a03\aD e0+\aUD a13\aD e1+\aUD a23\aD e2+\aUD a33\aD e3
+\aUD a43\aD e4+\aUD a53\aD e5+\aUD a63\aD e6+\aUD a73\aD e7
\\&=\frac 12
((\aUD f00-\aUD f33)\aD e0+(\aUD f10+\aUD f23)\aD e1
+(\aUD f20-\aUD f13)\aD e2+(\aUD f30+\aUD f03)\aD e3
\\ &
+(\aUD f40+\aUD f73)\aD e4+(\aUD f50-\aUD f63)\aD e5
+(\aUD f60+\aUD f53)\aD e6+(\aUD f70-\aUD f43)\aD e7)
\\&=\frac 12
(\aUD f00\aD e0-\aUD f33\aD e0+\aUD f10\aD e1+\aUD f23\aD e1
+\aUD f20\aD e2-\aUD f13\aD e2+\aUD f30\aD e3+\aUD f03\aD e3
\\ &
+\aUD f40\aD e4+\aUD f73\aD e4+\aUD f50\aD e5-\aUD f63\aD e5
+\aUD f60\aD e6+\aUD f53\aD e6+\aUD f70\aD e7-\aUD f43\aD e7)
\\&=\frac 12
((\aUD f00\aD e0+\aUD f10\aD e1+\aUD f20\aD e2+\aUD f30\aD e3+\aUD f40\aD e4
+\aUD f50\aD e5+\aUD f60\aD e6+\aUD f70\aD e7)
\\ &
+(\aUD f03\aD e3-\aUD f13\aD e2+\aUD f23\aD e1-\aUD f33\aD e0
-\aUD f43\aD e7+\aUD f53\aD e6-\aUD f63\aD e5+\aUD f73\aD e4))
\\&=\frac 12
((\aUD f00\aD e0+\aUD f10\aD e1+\aUD f20\aD e2+\aUD f30\aD e3+\aUD f40\aD e4
+\aUD f50\aD e5+\aUD f60\aD e6+\aUD f70\aD e7)
\\ &
+(\aUD f03\aD e0+\aUD f13\aD e1+\aUD f23\aD e2+\aUD f33\aD e3
+\aUD f43\aD e4+\aUD f53\aD e5+\aUD f63\aD e6+\aUD f73\aD e7)\aD e3)
\end{split}
}

\AddEquation{a4=f07}
{
\begin{split}
\aD a4&=\aUD a04\aD e0+\aUD a14\aD e1+\aUD a24\aD e2+\aUD a34\aD e3
+\aUD a44\aD e4+\aUD a54\aD e5+\aUD a64\aD e6+\aUD a74\aD e7
\\&=\frac 12
((\aUD f00-\aUD f44)\aD e0+(\aUD f10-\aUD f54)\aD e1
+(\aUD f20-\aUD f64)\aD e2+(\aUD f30-\aUD f74\aD e3
\\ &
+(\aUD f40+\aUD f04)\aD e4+(\aUD f50+\aUD f14)\aD e5
+(\aUD f60+\aUD f24)\aD e6+(\aUD f70+\aUD f34)\aD e7)
\\&=\frac 12
(\aUD f00\aD e0-\aUD f44\aD e0+\aUD f10\aD e1-\aUD f54\aD e1
+\aUD f20\aD e2-\aUD f64\aD e2+\aUD f30\aD e3-\aUD f74\aD e3
\\ &
+\aUD f40\aD e4+\aUD f04\aD e4+\aUD f50\aD e5+\aUD f14\aD e5
+\aUD f60\aD e6+\aUD f24\aD e6+\aUD f70\aD e7+\aUD f34\aD e7)
\\&=\frac 12
((\aUD f00\aD e0+\aUD f10\aD e1+\aUD f20\aD e2+\aUD f30\aD e3+\aUD f40\aD e4
+\aUD f50\aD e5+\aUD f60\aD e6+\aUD f70\aD e7)
\\ &
+(\aUD f04\aD e4+\aUD f14\aD e5+\aUD f24\aD e6+\aUD f34\aD e7
-\aUD f44\aD e0-\aUD f54\aD e1-\aUD f64\aD e2-\aUD f74\aD e3))
\\&=\frac 12
((\aUD f00\aD e0+\aUD f10\aD e1+\aUD f20\aD e2+\aUD f30\aD e3+\aUD f40\aD e4
+\aUD f50\aD e5+\aUD f60\aD e6+\aUD f70\aD e7)
\\ &
+(\aUD f04\aD e0+\aUD f14\aD e1+\aUD f24\aD e2+\aUD f34\aD e3
+\aUD f44\aD e4+\aUD f54\aD e5+\aUD f64\aD e6+\aUD f74\aD e7)\aD e4)
\end{split}
}

\AddEquation{a5=f07}
{
\begin{split}
\aD a5&=\aUD a05\aD e0+\aUD a15\aD e1+\aUD a25\aD e2+\aUD a35\aD e3
+\aUD a45\aD e4+\aUD a55\aD e5+\aUD a65\aD e6+\aUD a75\aD e7
\\&=\frac 12
((\aUD f00-\aUD f55)\aD e0+(\aUD f10+\aUD f45)\aD e1
+(\aUD f20-\aUD f75)\aD e2+(\aUD f30+\aUD f65)\aD e3
\\ &
+(\aUD f40-\aUD f15)\aD e4+(\aUD f50+\aUD f05)\aD e5
+(\aUD f60-\aUD f35)\aD e6+(\aUD f70+\aUD f25)\aD e7)
\\&=\frac 12
(\aUD f00\aD e0-\aUD f55\aD e0+\aUD f10\aD e1+\aUD f45\aD e1
+\aUD f20\aD e2-\aUD f75\aD e2+\aUD f30\aD e3+\aUD f65\aD e3
\\ &
+\aUD f40\aD e4-\aUD f15\aD e4+\aUD f50\aD e5+\aUD f05\aD e5
+\aUD f60\aD e6-\aUD f35\aD e6+\aUD f70\aD e7+\aUD f25\aD e7)
\\&=\frac 12
((\aUD f00\aD e0+\aUD f10\aD e1+\aUD f20\aD e2+\aUD f30\aD e3+\aUD f40\aD e4
+\aUD f50\aD e5+\aUD f60\aD e6+\aUD f70\aD e7)
\\ &
+(\aUD f05\aD e5-\aUD f15\aD e4+\aUD f25\aD e7-\aUD f35\aD e6
+\aUD f45\aD e1-\aUD f55\aD e0+\aUD f65\aD e3-\aUD f75\aD e2))
\\&=\frac 12
((\aUD f00\aD e0+\aUD f10\aD e1+\aUD f20\aD e2+\aUD f30\aD e3+\aUD f40\aD e4
+\aUD f50\aD e5+\aUD f60\aD e6+\aUD f70\aD e7)
\\ &
+(\aUD f05\aD e0+\aUD f15\aD e1+\aUD f25\aD e2+\aUD f35\aD e3
+\aUD f45\aD e4+\aUD f55\aD e5+\aUD f65\aD e6+\aUD f75\aD e7)\aD e5)
\end{split}
}

\AddEquation{a6=f07}
{
\begin{split}
\aD a6&=\aUD a06\aD e0+\aUD a16\aD e1+\aUD a26\aD e2+\aUD a36\aD e3
+\aUD a46\aD e4+\aUD a56\aD e5+\aUD a66\aD e6+\aUD a76\aD e7
\\&=\frac 12
((\aUD f00-\aUD f66)\aD e0+(\aUD f10+\aUD f76)\aD e1
+(\aUD f20+\aUD f46)\aD e2+(\aUD f30-\aUD f56\aD e3
\\ &
+(\aUD f40-\aUD f26)\aD e4+(\aUD f50+\aUD f36)\aD e5
+(\aUD f60+\aUD f06)\aD e6+(\aUD f70-\aUD f16)\aD e7)
\\&=\frac 12
(\aUD f00\aD e0-\aUD f66\aD e0+\aUD f10\aD e1+\aUD f76\aD e1
+\aUD f20\aD e2+\aUD f46\aD e2+\aUD f30\aD e3-\aUD f56\aD e3
\\ &
+\aUD f40\aD e4-\aUD f26\aD e4+\aUD f50\aD e5+\aUD f36\aD e5
+\aUD f60\aD e6+\aUD f06\aD e6+\aUD f70\aD e7-\aUD f16\aD e7)
\\&=\frac 12
((\aUD f00\aD e0+\aUD f10\aD e1+\aUD f20\aD e2+\aUD f30\aD e3+\aUD f40\aD e4
+\aUD f50\aD e5+\aUD f60\aD e6+\aUD f70\aD e7)
\\ &
+(\aUD f06\aD e6-\aUD f16\aD e7-\aUD f26\aD e4+\aUD f36\aD e5
+\aUD f46\aD e2-\aUD f56\aD e3-\aUD f66\aD e0+\aUD f76\aD e1))
\\&=\frac 12
((\aUD f00\aD e0+\aUD f10\aD e1+\aUD f20\aD e2+\aUD f30\aD e3+\aUD f40\aD e4
+\aUD f50\aD e5+\aUD f60\aD e6+\aUD f70\aD e7)
\\ &
+(\aUD f06\aD e0+\aUD f16\aD e1+\aUD f26\aD e2+\aUD f36\aD e3
+\aUD f46\aD e4+\aUD f56\aD e5+\aUD f66\aD e6+\aUD f76\aD e7)\aD e6)
\end{split}
}

\AddEquation{a7=f07}
{
\begin{split}
\aD a7&=\aUD a07\aD e0+\aUD a17\aD e1+\aUD a27\aD e2+\aUD a37\aD e3
+\aUD a47\aD e4+\aUD a57\aD e5+\aUD a67\aD e6+\aUD a77\aD e7
\\&=\frac 12
((\aUD f00-\aUD f77)\aD e0+(\aUD f10-\aUD f67)\aD e1
+(\aUD f20+\aUD f57)\aD e2+(\aUD f30+\aUD f47)\aD e3
\\ &
+(\aUD f40-\aUD f37)\aD e4+(\aUD f50-\aUD f27)\aD e5
+(\aUD f60+\aUD f17)\aD e6+(\aUD f70+\aUD f07)\aD e7)
\\&=\frac 12
(\aUD f00\aD e0-\aUD f77\aD e0+\aUD f10\aD e1-\aUD f67\aD e1
+\aUD f20\aD e2+\aUD f57\aD e2+\aUD f30\aD e3+\aUD f47\aD e3
\\ &
+\aUD f40\aD e4-\aUD f37\aD e4+\aUD f50\aD e5-\aUD f27\aD e5
+\aUD f60\aD e6+\aUD f17\aD e6+\aUD f70\aD e7+\aUD f07\aD e7)
\\&=\frac 12
((\aUD f00\aD e0+\aUD f10\aD e1+\aUD f20\aD e2+\aUD f30\aD e3+\aUD f40\aD e4
+\aUD f50\aD e5+\aUD f60\aD e6+\aUD f70\aD e7)
\\ &
+(\aUD f07\aD e7+\aUD f17\aD e6-\aUD f27\aD e5-\aUD f37\aD e4
+\aUD f47\aD e3+\aUD f57\aD e2-\aUD f67\aD e1-\aUD f77\aD e0))
\\&=\frac 12
((\aUD f00\aD e0+\aUD f10\aD e1+\aUD f20\aD e2+\aUD f30\aD e3+\aUD f40\aD e4
+\aUD f50\aD e5+\aUD f60\aD e6+\aUD f70\aD e7)
\\ &
+(\aUD f07\aD e7+\aUD f17\aD e1+\aUD f27\aD e2+\aUD f37\aD e4
+\aUD f61\aD e6+\aUD f47\aD e4+\aUD f57\aD e5+\aUD f77\aD e7)\aD e7)
\end{split}
}

\AddEq{a0=f0.7}
{
\ShowEq{a0=f07 1}
\ShowEq{a0=f07 2}
\ShowEq{a0=f07 3}
\ShowEq{a0=f07}
\ShowEq{a1=f07}
\ShowEq{a2=f07}
\ShowEq{a3=f07}
\ShowEq{a4=f07}
\ShowEq{a5=f07}
\ShowEq{a6=f07}
\ShowEq{a7=f07}
}

\AddEquation{f=.E+.I+.J+.K}
{
\begin{split}
f=&-\frac 12
\left(
\aD f0
+\aD f1i
+\aD f2j
+\aD f3k
\right)\bullet E
\\&
+\frac 12\left(
\aD f0
+\aD f1i
\right)\bullet\aU I1
+\frac 12\left(
\aD f0
+\aD f2j
\right)\bullet\aU I2
+\frac 12\left(
\aD f0
+\aD f3k
\right)\bullet\aU I3
\end{split}
}

\AddEq[2]{ab=a*(.E+.I+.J+.K)b}
{
#1 #2=
(#2^*\bullet E+(i#2^*)\bullet\aU I1+(j#2^*)\bullet\aU I2+(k#2^*)\bullet\aU I3)\circ #1
}

\AddEquation{f=.I0.7}
{
\begin{split}
f=&-\frac 12
\left(
5\aD f0+\aD f1i+\aD f2j+\aD f3k
\aD f4l+\aD f5il+\aD f6jl+\aD f7kl
\right)\circ E
\\&
+\frac 12\left(
\aD f0+\aD f1i
\right)\circ\aU I1
+\frac 12\left(
\aD f0+\aD f2j
\right)\circ\aU I2
+\frac 12\left(
\aD f0+\aD f3k
\right)\circ\aU I3
\\&
+\frac 12\left(
\aD f0+\aD f4l
\right)\circ\aU I4
+\frac 12\left(
\aD f0+\aD f5il
\right)\circ\aU I5
+\frac 12\left(
\aD f0+\aD f6jl
\right)\circ\aU I6
\\&
+\frac 12\left(
\aD f0+\aD f7kl
\right)\circ\aU I7
\end{split}
}

\AddEquation{a...= O}
{
\begin{split}
&\,
\ShowEq{matrix a algebra O}
\\ = &\,
\ShowEq{matrix af algebra O}
\\ * &\,
\ShowEq{matrix f algebra O}
\end{split}
}

\AddEq [2]{I=(E,I)}
{
\def\Set{\aU I1,\aU I2,\aU I3}
\def\Temp{O}%
\edef\Tempa{#1}%
\ifx\Tempa\Temp%
\def\Set{\aU Ii,\gi i=\gi 1,...,\gi 7}
\fi
$\Basis I=(E,\Set)$#2
}

\AddEq{F o G}
{
\Basis F\circ G=\{F_k\circ G:F_k\in \Basis F\}
}

\AddEq{FijG}
{
$F_i\circ G$, \iI,
}

\AddEquation{coordinates of map Fk, associative algebra}
{
F_k\circ x=\aUD {F_{k\cdot}^{}{}}ij\aU xj\aD ei
}

\AddEq [1]{linear map over ring, determinant=0, 1}
{
\[
\rank
\begin{pmatrix}
\mathcal C^{\cdot}{}_{\gim}^{\gik}{}_{\cdot\gi{ij}}
&#1_{\gim}^{\gik}
\end{pmatrix}
=\rank\mathcal C
\]
}

\AddEq[1]{linear map over ring, matrix}
{
\def\Cdot{}
\def\Temp{-}%
\edef\Tempa{#1}%
\ifx\Tempa\Temp%
\else
\def\Cdot{#1\cdot}
\fi
\mathcal C=
\left(
\mathcal C^{\cdot}{}_{\gim}^{\gik}{}_{\cdot\gi{ij}}
\right)
=
\left(
C_{\Cdot}{}^{\gi p}_{\gi{im}}C_{\Cdot}{}^{\gik}_{\gi{pj}}
\right)
}

\AddEq{linear map over ring, row of matrix}
{
${}^{\cdot}{}_{\gim}^{\gik}$
}

\AddEq{linear map over ring, column of matrix}
{
${}_{\cdot\gi{ij}}$
}

\AddEquation{coordinates of map f, associative algebra}
{
f\circ x=\aUD fij\aU xj\aD ei
}

\DefEq
{
\[
\Basis F=\{F_k\in\mathcal L(D;A_1\rightarrow A_2): k=1, ..., n\}
\]
}
{Ik 1n}

\AddEquation{h generated by f, associative algebra}
{
h=
(
a\pC{s}{0}
\otimes
a\pC{s}{1}
)
\circ f
=a\pC{s}{0}
f
a\pC{s}{1}
}

\DefEq
{
\[
h:A_1\rightarrow A_2
\]
}
{h:A1->A2}

\AddEquation{f in L(A,A), 1, associative algebra}
{
f=
(
a_{k\cdot s_k\cdot 0}
\otimes
a_{k\cdot s_k\cdot 1}
)
\circ I_k
=
\sum_{\gik}
a_{k\cdot s_k\cdot 0}
I_k
a_{k\cdot s_k\cdot 1}
}

\DefEq
{
\[
f\in B^A\rightarrow \|f\|\in R
\]
}
{f in BA->|f|}

\AddEq{norm of map}
{
\symb{\|f\|}{norm of map}{}
}

\DefEquation
{
\ShowSymbol{norm of map}{}=
\text{sup}\frac{\|f(x)\|_B}{\|x\|_A}
}
{norm of map, algebra}

\DefEq
{
$\ShowSymbol{norm of map}{}$
}
{show|f|}

\DefEq
{
\symb{f+g}{sum of maps}{polylinear}
}
{sum of maps,,polylinear}

\DefEquation
{
\begin{matrix}
\ShowSymbol{sum of maps}{polylinear}:A_1\times...\times A_n\rightarrow S
&f,g\in\mathcal L(D;A_1\times...\times A_n\rightarrow S)
\end{matrix}
}
{sum of maps, 1, polylinear}

\AddEq{f=fkxfk}
{
f^k=f^k_{s_k\cdot 0}\otimes f^k_{s_k\cdot 1}
}

\AddEq [1]{expansion of linear map with respect to basis}
{
f=f^k\circ I_k\ \ \ f^k\in #1\otimes #1
}

\DefEquation
{
(\ShowSymbol{sum of maps}{polylinear})\circ (a_1,...,a_n)
=f\circ (a_1,...,a_n)+g\circ (a_1,...,a_n)
}
{sum of  maps, polylinear}

\AddEq{module of polylinear maps}
{
$\mathcal L(D;A_1\times...\times A_n\rightarrow S)$
}

\AddEquation{dg h12 xox}
{
dg=(h_2h_1-h_1h_2)x+x(h_1h_2-h_2h_1)
}

\AddEquation{a0=}
{
\aD a0=-\frac 12(\aD f0+\aD f1i+\aD f2j+\aD f3k)
}

\AddEquation{a1=}
{
a_1=\frac 12(\aD f0+\aD f1i)
}

\AddEquation{a2=}
{
\aD a2=\frac 12(\aD f0+\aD f2j)
}

\AddEquation{a3=}
{
\aD a3=\frac 12(\aD f0+\aD f3k)
}

\AddEq[3]{set linear maps, module}
{
\symb{\mathcal L(#1;#2\rightarrow #3)}{set linear maps}1
}

\AddEq[2]{set linear maps, module (1)}
{
\symb{\mathcal L(#1;#2_1\rightarrow #2_2)}{set linear maps}1
}

\AddEq[2]{set linear maps, module (11)}
{
\symb{\mathcal L(#1_1\rightarrow #1_2;#2_1\rightarrow #2_2)}{set linear maps}1
}

\AddEq{set homomorphisms, D algebra}
{
\symb{\Hom(D;A_1\rightarrow A_2)}{set of homomorphisms}1
}

\AddEquation{Am->A(n+m) 1}
{
(p\circ x^n)\circ (q\circ x^m)=(p\underline{\otimes}q)\circ x^{n+m}
}

\AddEq[1]{d->dei}
{
\[
\aD P#1:d\in D\rightarrow d\aD e#1\in A
\]
}

\AddEq{f(t)=fit ei}
{
\[
f\circ t=(\aU fit)\aD ei=\aD Pi\circ(\aU fit)
\]
}

\AddEq{map f, 1, tensor product}
{
\[
f:A_1\times...\times A_n\rightarrow\Tensor A
\]
}

\AddEq{tensor product of modules}
{
\symb{\Tensor A}{tensor product}{}
}

\AddEq{f:xA->oxA}
{
\[f:A_1\times...\times A_n\rightarrow\ShowSymbol{tensor product}{}\]
}

\AddEquation{map f, tensor product}
{
f\circ(d_1,...,d_n)=\Tensor d
}

\AddEq{map g, tensor product}
{
\[
g:A_1\times...\times A_n\rightarrow V
\]
}

\AddEq{map h, tensor product}
{
\[
h:\Tensor A\rightarrow V
\]
}

\AddEquation{map gh, tensor product}
{
\xymatrix{
&\Tensor A\ar[dd]^h
\\
\Times
\ar[ru]^f\ar[rd]_g
\\
&V
}
}

\AddEquation{g=h, tensor product}
{
h\circ(\Tensor a)=g\circ(a_1,...,a_n)
}

\AddEq{fxa=oxa}
{
\[f\circ(a_1,...,a_n)=\Tensor a\]
}

\AddEq[2]{polylinear map of algebras, 1}
{
\[
\begin{matrix}
1\le i\le n
&
a_i, b_i \in #1
&
p\in #2
\end{matrix}
\]
}

\AddEq{aE+bE=(a+b)E}
{
\[(aE)\circ x+(bE)\circ x=ax+bx=(a+b)x=((a+b)E)\circ x\]
}

\AddEq{aE o bE=(ab)E}
{
\[(aE)\circ (bE)\circ x=(aE)\circ (bx)=a(bx)=(ab)x=((ab)E)\circ x\]
}

\AddEq{g=gE}
{
\[g(x)=g^k(x)\circ E_k\ \ \ \ g^k:A\rightarrow B\]
}

\AddEq{linear map times constant, algebra}
{
$a\circ f$, $b\star f$, $a$, $b\in A_2$,
}

\AddEq{linear map times constant, 0, algebra}
{
\begin{align*}
(a\circ f)\circ x&=a(f\circ x)
\\
(b\star f)\circ x&=(f\circ x)b
\end{align*}
}

\AddEq{linear map AA LAA}
{
\[
h:\ATwo\rightarrow \mathcal L(D;A_1\rightarrow A_2)
\]
}

\AddEquation{linear map AA LAA, 1}
{
(a\otimes b)\circ f=b\star (a\circ f)
}

\AddEq{linear map AA LAA, 2}
{
\[
((a\otimes b)\circ f)\circ x=a(f\circ x)b
\]
}

\AddEq{f circ a =}
{
f\circ a=f(a)
}

\AddEq [5]{matrix AIJ}
{
\[
#1=(#1^{\gi #2}_{\gi #4},\gi #2\in\gi #3,\gi #4\in\gi #5)
\]
}

\AddEq{d12 in D}
{
\(d_1\), \(d_2\in D\).
}

\AddEquation{d1 d2 in A}
{
(d_11_A)(d_21_A)=d_1(1_A(d_21_A))=d_1(d_2(1_A1_A))
=d_1(d_21_A)=(d_1d_2)1_A
}

\AddEquation{ad=da in A}
{
a(d1_A)=d(a1_A)=d(1_Aa)=(d1_A)a
}

\AddEquation{d1+d2 in A}
{
d_11_A+d_21_A=(d_1+d_2)1_A
}

\AddEq{fi(dx)=}
{
\[
f_i(dx_i)=df_i(x_i)
\]
}

\AddEq{f(dx)i=}
{
\begin{align*}
f((dx_i,i\in I))&=
\sum_{i\in I}f_i(dx_i)=
\sum_{i\in I}df_i(x_i)
=d\sum_{i\in I}f_i(x_i)
\\&=df((x_i,i\in I))
\end{align*}
}

\AddEquation{linear map, f a+b=}
{
f\circ(a+b)=f\circ a+f\circ b
}

\AddEquation{linear map, f pa=}
{
f\circ(pa)=p(f\circ a)
}

\AddEq [1]{ab V p D}
{
\[
\begin{matrix}
a,b\in #1
&
p\in D
\end{matrix}
\]
}

\AddEq[3]{homomorphism, f v+w=}
{
#1\circ(#2+#3)=#1\circ #2+#1\circ #3
}

\AddEq[3]{homomorphism, f vw=}
{
#1\circ(#2#3)=(#1\circ #2)(#1\circ #3)
}

\AddEq[3]{homomorphism, f(p+q)=}
{
#1(#2+#3)=#1(#2)+#1(#3)
}

\AddEq[3]{homomorphism, f(pq)=}
{
#1(#2#3)=#1(#2)#1(#3)
}

\AddEquation{left linear map, f dv=}
{
f\circ(dv)=d(f\circ v)
}

\AddEquation{right linear map, f dv=}
{
f\circ(vd)=(f\circ v)d
}

\AddEq[4]{left homomorphism, f av=}
{
#1\circ(#2#4)=#3(#1\circ #4)
}

\AddEq[4]{right homomorphism, f av=}
{
#1\circ(#4#2)=(#1\circ #4)#3
}

\AddEq{D1 D2 ab in A, d in D}
{
\[
\begin{matrix}
a,b\in \Module_{\VF}
&d,d_1,d_2\in \Base_{\DF}
\end{matrix}
\]
}

\AddEq{D1 D2 uv in V, ab in A, d in D}
{
\[
\begin{matrix}
u,v\in V_1
&a,b\in A_1
&d,d_1,d_2\in D_1
\end{matrix}
\]
}

\AddEq[2]{D ab in A, d in D}
{
\[
\begin{matrix}
a,b\in #2
&
d\in #1
\end{matrix}
\]
}

\AddEq{ap AD 11}
{
\[
p,q\in D_1
\ \ \ \,
a,b\in A_1
\]
}

\AddEq{ap AD 1}
{
\[
p\in D
\ \ \ \,
a,b\in A_1
\]
}

\AddEq{vp VD 11}
{
\[
p,q\in D_1
\ \ \ \,
v,w\in V_1
\]
}

\AddEq{vp VD 1}
{
\[
p\in D
\ \ \ \,
v,w\in V_1
\]
}

\AddEq{vap VAD 111}
{
\[
p,q\in D_1
\ \ \ \,
a,b\in A_1
\ \ \ \,
u,v\in V_1
\]
}

\AddEq{vap VAD 11}
{
\[
p\in D
\ \ \ \,
a,b\in A_1
\ \ \ \,
u,v\in V_1
\]
}

\AddEq{vap VAD 1}
{
\[
a\in A
\ \ \ \,
u,v\in V_1
\]
}

\AddEq{commutator of algebra}
{
\symb{[a,b]}{commutator of algebra}{}
\[
\ShowSymbol{commutator of algebra}{}=ab-ba
\]
}

\AddEq{commutative D algebra}
{
\[
[a,b]=0
\]
}

\AddEquation{associator of algebra =}
{
\ShowSymbol{associator of algebra}{}=(ab)c-a(bc)
}

\AddEquation{associator of algebra = Ae}
{
\begin{split}
(a,b,c)&=(\aUD Ck{pl}(ab)^{\gi p}c^{\gil}-\aUD Ck{ip}a^{\gii}(bc)^{\gi p})e_{\gik}
\\
&=(\aUD Ck{pl}\aUD Cp{ij}-\aUD Ck{ip}\aUD Cp{jl})a^{\gii}b^{\gij}c^{\gil}e_{\gik}
\end{split}
}

\AddEq{associator of algebra}
{
\symb{(a,b,c)}{associator of algebra}{}
}

\AddEquation{associator of algebra = A}
{
(a,b,c)=A^{\gi k}_{\gi{ijl}}a^{\gii}b^{\gij}c^{\gil}e_{\gik}
}

\AddEq{coordinates of associator}
{
\symb{A^{\gi k}_{\gi{ijl}}}{coordinates of associator}{}
}

\AddEquation{coordinates of associator =}
{
\ShowSymbol{coordinates of associator}{}=
\aUD Ck{pl}\aUD Cp{ij}-\aUD Ck{ip}\aUD Cp{jl}
}

\AddEq{associative D algebra}
{
\[
(a,b,c)=0
\]
}

\AddEq{nucleus of algebra}
{
\symb{N(A)}{nucleus of algebra}{}
\[
\ShowSymbol{nucleus of algebra}{}=
\{
a\in A:
\forall b, c\in A,
(a,b,c)=(b,a,c)=(b,c,a)=0
\}
\]
}

\AddEq{center of algebra}
{
\symb{Z(A)}{center of algebra}{}
\[
\ShowSymbol{center of algebra}{}=
\{
a\in A:
a\in N(A),
\forall b\in A,
ab=ba
\}
\]
}

\AddEq[3]{f(a+b)=...}
{
#1\circ(#2+#3)=#1\circ #2+#1\circ #3
}

\DefEq
{
\[\Vector r_2:A_1\rightarrow A_2\]
}
{r2:A1->A2}

\AddEq{property linear map hf}
{
\DrawEq{h(d1+d2)=...}{\SideWS module}
\DrawEq{h(d1d2)=...}{\SideWS module}
\DrawEq[fab]{f(a+b)=...}{\DFDT \SideWS module}
\DrawEq[hfa]{f(da)=h... \SideWS module}{\DFDT \SideWS module}
\ShowEq{D1 D2 ab in A, d in D}
}

\AddEq{property linear map hgf}
{
\DrawEq{h(d1+d2)=...}{\SideWS module}
\DrawEq{h(d1d2)=...}{\SideWS module}
\DrawEq[gab]{f(a+b)=...}{g \DFDT \SideWS module}
\DrawEq[hga]{f(da)=h... \SideWS module}{g \DFDT}
\DrawEq[fuv]{f(a+b)=...}{f \DFDT \SideWS module}
\DrawEq[hfu]{f(da)=h... \SideWS module}{f \DFDT}
\ShowEq{D1 D2 uv in V, ab in A, d in D}
}

\AddEq[2]{show linear map hgf}
{
\[
\begin{pmatrix}
h:D_1\rightarrow D_2&
#1:A_1\rightarrow A_2&
#2:V_1\rightarrow V_2
\end{pmatrix}
\]
}

\AddEq[2]{show linear map hf}
{
\ShowEq{h:D1->D2 f:A1->A2}{#1} \Base {#2} \Module
}

\AddEq[2]{show linear map f}
{
\ShowEq{f:A->B}{#2}{\Module_1}{\Module_2}
}

\AddEq [3]{V=o+Ae}
{
\[
V_{#1}=\bigoplus_{\gi #2\in\gi #3}A_{#1}\EV{V_{#1}}#2
\]
}

\AddEq{w in Jv}
{
$w\in X_k\subseteq J(v)$.
}

\AddEq[3]{f(da)=h... module}
{
\def\Temp{}%
\edef\Tempa{#1}%
\ifx\Tempa\Temp%
\def\Koef{d}%
\else
\def\Koef{#1(d)}
\fi
#2\circ(d#3)=\Koef(#2\circ #3)
}

\AddEq[3]{f(da)=h... left module}
{
#2\circ(d#3)=#1(d)(#2\circ #3)
}

\AddEq[3]{f(da)=h... right module}
{
#2\circ(d#3)=(#2\circ #3)#1(d)
}

\AddEquation{D , f(da)=... module}
{
f\circ(da)=(h\circ d)(f\circ a)
}

\AddEquation{D , f(da)=... left module}
{
f\circ(da)=(h\circ d)(f\circ a)
}

\AddEquation{D , f(da)=... right module}
{
f\circ(ad)=(f\circ a)(h\circ d)
}

\AddEq{h(d1+d2)=...}
{
h(d_1+d_2)=h(d_1)+h(d_2)
}

\AddEq{h(d1d2)=...}
{
h(d_1d_2)=h(d_1)h(d_2)
}

\AddEq{property linear map f}
{
\DrawEq[fab]{f(a+b)=...}{D }
\DrawEq[{}fa]{f(da)=h... module}{\DFDT \SideWS module}
\ShowEq{D ab in A, d in D}D{A_1}
}

\AddEq{commutative product in algebra, 1}
{
\AUD Cp{ij}=\AUD Cp{ji}
}

\AddEq{associative product in algebra, 1}
{
\AUD Cp{ij}\AUD Cq{pk}
=
\AUD Cq{ip}\AUD Cp{jk}
}

\AddEq{commutative product in algebra}
{
\[
\ARow ei\ARow ej
=
\ARow ej\ARow ei
\]
}

\AddEq{associative product in algebra}
{
\[
(\ARow ei\ARow ej)\ARow ek
=
\ARow ei(\ARow ej\ARow ek)
\]
}

\AddEq [3]{basis ei}
{
\[\eV[#1]=(\ECol{#1}{#2},\jJg {#2}{#3})\]
}

\AddEq [3]{basis e of module cols}
{
\eV[#1]=(\ECol{#1}{#2},\jJg {#2}{#3})
}

\AddEq[3]{basis e of module rows}
{
\eV[#1]=(\ERow{#1}{#2},\jJg {#2}{#3})
}

\AddEq {basis e of module}
{
\eV=(\ARow ej,\jJg jJ)
}

\AddEq [2]{basis e of module 1n}
{
\[\eV[#1]=(\EV{#1}1,...,\EV{#1}{#2})\]
}

\AddEquation{(a+b)c=.}
{
(a+b)c=ac+bc
}

\AddEquation{a(b+c)=.}
{
a(b+c)=ab+ac
}

\AddEq[1]{f:A1->A2, A module}
{
b=#1\CRstar f
}

\AddEq[2]{f:V1->V2, D module cols}
{
b=#1#2
}

\AddEq[2]{f:V1->V2, D module rows}
{
b=#2#1
}

\AddEquation{ab=ac a ne 0}
{
ab=ac\ \ \ \ a\ne 0
}

\AddEquation{ab=ac a ne 0 1}
{
b=(a^{-1}a)b=a^{-1}(\BlueText{ab})=a^{-1}(\BlueText{ac})=(a^{-1}a)c=c
}

\AddEq[2]{ker G in ker g}
{
$\ker G\subseteq\ker #1$#2
}

\AddEq [2]{va=ae1, module (left)(cols)}
{
\Vector #1=#1\CRstar e_{#2}
}

\AddEq [2]{va=ae1, module ()(cols)}
{
\Vector #1= e_{#2}#1
}

\AddEq [2]{va=ae1, module (right)(cols)}
{
\Vector #1=e_{#2}\RCstar #1
}

\AddEq [2]{va=ae1, module (left)(rows)}
{
\Vector #1=#1\RCstar e_{#2}
}

\AddEq [2]{va=ae1, module ()(rows)}
{
\Vector #1=#1 e_{#2}
}

\AddEq [2]{va=ae1, module (right)(rows)}
{
\Vector #1=e_{#2}\CRstar #1
}

\AddEq{Morphism of D algebra, injection}
{
$h(\ACol ai)$, $h(\ACol bi)$)
}

\AddEq[1]{algebra, homomorphism and product (11)}
{
h(\AUD{C_1^{}{}}k{ij})\AUD {#1}lk
=\AUD {#1}pi\AUD {#1}qj\AUD{C_2^{}{}}l{pq}
}

\AddEq[1]{algebra, homomorphism and product (1)}
{
\AUD{C_{1}^{}{}}k{ij}\AUD {#1}lk
=\AUD {#1}pi\AUD {#1}qj\AUD{C_2^{}{}}l{pq}
}

\AddEq[1]{kij in I}
{
$\gik$, $\gii$, \jJg j{#1},
}

\AddEq{algebra, homomorphism and product abe 2}
{
\begin{aligned}
&\,h(\AUD{C_1^{}{}}k{ij})\AUD flkh(\ACol ai)h(\ACol bj)\EBase 2l
\\ =&\,h(\AUD{C_1^{}{}}k{ij}\ACol ai\ACol bj)\AUD flk\EBase 2l
\end{aligned}
}

\AddEq{algebra, homomorphism and product abe 2(1)}
{
\begin{aligned}
&\,h(\AUD{C_1^{}{}}k{ij})\AUD flkh(\ACol ai)h(\ACol bj)\EBase 2l
\\ =&\,h(\ACol {(ab)}k)\AUD flk\EBase 2l
\end{aligned}
}

\AddEq{algebra, homomorphism and product abe 2()}
{
\begin{aligned}
&\,\AUD{C_1^{}{}}k{ij}\AUD flk\ACol ai\ACol bj\EBase 2l
\\ =&\,\ACol {(ab)}k\AUD flk\EBase 2l
\end{aligned}
}

\AddEq{algebra, homomorphism and product abe 3(1)}
{
\begin{aligned}
&\,h(\AUD{C_1^{}{}}k{ij})\AUD flkh(\ACol ai)h(\ACol bj)\EBase 2l
\\ =&\,f\circ(ab)
\end{aligned}
}

\AddEq{algebra, homomorphism and product abe 3()}
{
\begin{aligned}
&\,\AUD{C_1^{}{}}k{ij}\AUD flk\ACol ai\ACol bj\EBase 2l
\\ =&\,f\circ(ab)
\end{aligned}
}

\AddEq{algebra, homomorphism and product abe 11}
{
\begin{aligned}
&\, h(\AUD{C_1^{}{}}k{ij})\AUD flkh(\ACol ai)h(\ACol bj)\EBase 2l
\\ =&\,\AUD fpi\AUD fqj\AUD{C_2^{}{}}l{pq}h(\ACol ai)h(\ACol bj)\EBase 2l
\end{aligned}
}

\AddEq{algebra, homomorphism and product abe 1}
{
\begin{aligned}
&\, \AUD{C_{1}^{}{}}k{ij}\AUD flk\ACol ai\ACol bj\EBase 2l
\\ =&\,\AUD fpi\AUD fqj\AUD{C_2^{}{}}l{pq}\ACol ai\ACol bj\EBase 2l
\end{aligned}
}

\AddEq{algebra, homomorphism and product abe 5(1)}
{
\begin{aligned}
&\, \AUD fpi\AUD fqj\AUD{C_2^{}{}}l{pq}h(\ACol ai)h(\ACol bj)\EBase 2l
\\ =&\,(f\circ a)(f\circ b)
\end{aligned}
}

\AddEq{algebra, homomorphism and product abe 5()}
{
\begin{aligned}
&\, \AUD fpi\AUD fqj\AUD{C_2^{}{}}l{pq}\ACol ai\ACol bj\EBase 2l
\\ =&\,(f\circ a)(f\circ b)
\end{aligned}
}

\AddEq{algebra, homomorphism and product abe 4(1)}
{
\begin{aligned}
&\, \AUD fpi\AUD fqj\AUD{C_2^{}{}}l{pq}h(\ACol ai)h(\ACol bj)\EBase 2l
\\ =&\,(h(\ACol ai)\AUD fpi\EBase 2p)(h(\ACol bj)\AUD fqj\EBase 2q)
\end{aligned}
}

\AddEq{algebra, homomorphism and product abe 4()}
{
\begin{aligned}
&\, \AUD fpi\AUD fqj\AUD{C_2^{}{}}l{pq}\ACol ai\ACol bj\EBase 2l
\\ =&\,(\ACol ai\AUD fpi\EBase 2p)(\ACol bj\AUD fqj\EBase 2q)
\end{aligned}
}

\AddEq{a, b in A1}
{
\[
\begin{matrix}
a,b\in A_1&
a=\ACol ai\EBase 1i&
b=\ACol bi\EBase 1i
\end{matrix}
\]
}

\AddEq{algebra, homomorphism and product, 1}
{
\newline
\FrameEqRef{product of basis vectors, algebra}{\Cols}
\newline
}

\AddEq{algebra, homomorphism and product eiej(11)}
{
\Vector f\circ (\EBase 1i\EBase 1j)
=h(C_1^{}\AUD{{}}k{ij})\Vector f\circ \EBase 1k
}

\AddEq{algebra, homomorphism and product eiej(1)}
{
\Vector f\circ (\EBase 1i\EBase 1j)
=C_1^{}\AUD{{}}k{ij}\Vector f\circ \EBase 1k
}

\AddEq{algebra, homomorphism and product eiej 4(11)}
{
\Vector f\circ (\EBase 1i\EBase 1j)
=h(C_1^{}\AUD{{}}k{ij})\aUD flk\EBase 2l
}

\AddEq{algebra, homomorphism and product eiej 4(1)}
{
\Vector f\circ (\EBase 1i\EBase 1j)
=C_1^{}\AUD{{}}k{ij}\aUD flk\EBase 2l
}

\AddEq{algebra, homomorphism and product, 4}
{
\Vector f\circ (ab)
=h(C_1^{}\AUD{{}}k{ij})
h(\ACol ai)h(\ACol bj)\Vector f\circ \EBase 1k
}

\AddEq{algebra, homomorphism and product, 5 11}
{
\Vector f\circ (ab)
=h(C_1^{}\AUD{{}}k{ij})h(\ACol ai)h(\ACol bj)
\AUD flk\EBase 2l
}

\AddEq{algebra, homomorphism and product, 5 1}
{
\Vector f\circ (ab)
=C_1^{}\AUD{{}}k{ij}\ACol ai\ACol bj
\AUD flk\EBase 2l
}

\AddEq{algebra, homomorphism and product eiej 5}
{
\begin{aligned}
&\,\Vector f\circ (\EBase 1i\EBase 1j)\\ =&\,
(\Vector f\circ \EBase 1i)(\Vector f\circ \EBase 1j)\\ =&\,
\AUD fpi\EBase 2p
\AUD fqj\EBase 2q
\end{aligned}
}

\AddEq{algebra, homomorphism and product eiej 2(11)}
{
\Vector f\circ (\EBase 1i\EBase 2j)=
\AUD fpi
\AUD fqj
C_2^{}\AUD{{}}l{pq}
\EBase 2l
}

\AddEq{algebra, homomorphism and product eiej 2(1)}
{
\Vector f\circ (\EBase 1i\EBase 2j)=
\AUD fpi
\AUD fqj
C_2^{}\AUD{{}}l{pq}
\EBase 2l
}

\AddEq{algebra, homomorphism and product eiej 3(11)}
{
h(C_1^{}\AUD{{}}k{ij})
\AUD flk\EBase 2l
=
\AUD fpi
\AUD fqj
C_2^{}\AUD{{}}l{pq}
\EBase 2l
}

\AddEq{algebra, homomorphism and product eiej 3(1)}
{
C_1^{}\AUD{{}}k{ij}
\AUD flk\EBase 2l
=
\AUD fpi
\AUD fqj
C_2^{}\AUD{{}}l{pq}
\EBase 2l
}

\AddEq{product in D algebra}
{
v\,w=C\circ(v,w)
}

\AddEq{product in algebra, definition 1}
{
\[
C:A\times A\rightarrow A
\]
}

\DefEq
{
\[
\underline{\otimes}:A^{n\otimes }\times A^{m\otimes }
\rightarrow A^{n+m-1\otimes}
\]
}
{otimes -}

\AddEq{a b in basis of algebra}
{
\[
\begin{matrix}
a=\ACol ai\ARow ei&b=\ACol bi\ARow ei&a, b\in A
\end{matrix}
\]
}

\AddEq{product in algebra}
{
\ACol{(ab)}k=\AUD Ck{ij}\ACol ai\ACol bj
}

\AddEq{product of basis vectors, algebra}
{
\ARow ei\ARow ej=
\AUD Ck{ij}\ARow ek
}

\AddEq{product in algebra, 1}
{
ab=\ACol ai\ACol bj\Eij
}

\AddEq{product in algebra, 2}
{
ab=\ACol ai\ACol bj
\AUD Ck{ij}\ARow ek
}

\AddEq[4]{diagram of representations of D algebra}
{
\begin{matrix}
\vcenter
{
\xymatrix
{
A_{#3}\ar[rr]|{*}^{#4_{#2 23}}&&
A_{#3}
\\
&D_{#1}\ar[ur]|{*}_{#4_{#2 12}}\ar[ul]|{*}_{#4_{#2 12}}
}
}
&
\begin{array}{r@{\,}l}
#4_{#2 12}(d):v&\rightarrow d\,v
\\
#4_{#2 23}(v):w&\rightarrow
C_{#3}\circ(v, w)
\\
C_{#3}&\in\mathcal L(D_{#1};A_{#3}^2\rightarrow A_{#3})
\end{array}
\end{matrix}
}

\AddEq{endomorphism of module from product, 8}
{
\[
ab=(g_{23}\circ a)\circ b
\]
}

\AddEquation{endomorphism of module from product, 3}
{
\begin{array}{r@{\,}l}
(g_{23}\circ v)(w_1+w_2)&=(g_{23}\circ v)\circ w_1+(g_{23}\circ v)\circ w_2
\\
(g_{23}\circ v)\circ (dw)&=d((g_{23}\circ v)\circ w)
\end{array}
}

\AddEquation{endomorphism of module from product, 4}
{
(g_{23}\circ (v_1+v_2))\circ w
=(g_{23}\circ v_1+g_{23}\circ v_2)(w)
=(g_{23}\circ v_1)\circ w+(g_{23}\circ v_2)\circ w
}

\AddEquation{endomorphism of module from product, 7}
{
(g_{23}\circ(d\,v))\circ w
=(d\,(g_{23}\circ v))\circ w
=d\,((g_{23}\circ v)\circ w)
}

\AddEquation{r=+q circ p}
{
\begin{split}
r(x)&=s_0+(q_0+q_1(x)+...+q_{k-1}(x))\circ p(x)
\end{split}
}

\AddEq{qij i j}
{
$q_i\in A_i[x]\otimes A$, $i=0$, ..., $k-1$,
}

\AddEquation{monomial of power k, F=1}
{
p_k(x)=p_{k-1}(x)xa_k
}

\AddEq{a circ x=b}
{
a\circ x=b
}

\AddEq[5]{structural constants of algebra}
{
\ensuremath{C_{#1}^{}{}\AUD {}{#2}{#3#4}{#5}}
}

\AddEq{structural constants of algebra symb}
{
\symb{
\ShowEq{structural constants of algebra}{}kij{}
}{structural constants}1
}

\AddEquation{otimes -, 1}
{
(a_1\otimes...\otimes a_n)\underline{\otimes}(b_1\otimes...\otimes b_n)
=
a_1\otimes...\otimes a_{n-1}\otimes a_nb_1\otimes b_2\otimes...\otimes b_n
}

\AddEq[3]{a=oxai}
{
$#1=#1_0\otimes...\otimes #1_{#2}$#3
}

\AddEq[3]{vb=f(va)}
{
\Vector #2=\Vector #3\circ\Vector #1
}

\AddEquation{vb=h(e1)a (1)(1)(cols)}
{
\Vector b=\Vector f\circ\Vector a=\Vector f\circ(e_1a)
= (\Vector f\circ e_1)h(a)
}

\AddEquation{vb=h(e1)a (1)(1)(rows)}
{
\Vector b=\Vector f\circ\Vector a=\Vector f\circ(a e_1)
= h(a)(\Vector f\circ e_1)
}

\AddEquation{vb=h(e1)a ()(1)(cols)}
{
\Vector b=\Vector f\circ\Vector a=\Vector f\circ(e_1a)
= (\Vector f\circ e_1)a
}

\AddEquation{vb=h(e1)a ()(1)(rows)}
{
\Vector b=\Vector f\circ\Vector a=\Vector f\circ(ae_1)
= a(\Vector f\circ e_1)
}

\AddEq[5]{f o ea=efa i (11)}
{
\Vector{#2}\circ(\ACol {#5}i\EBase {#3}i)
= #1(\ACol {#5}i)\AUD {#2}ki\EBase {#4}k
}

\AddEq[5]{f o ea=efa i (1)}
{
\Vector{#2}\circ(\ACol {#5}i\EBase {#3}i)
= \ACol {#5}i\AUD {#2}ki\EBase {#4}k
}

\AddEq[4]{f o ev=efv i (11)}
{
\Vector{#2}\circ(\Multiply{\ACol vi}{\EBase {#3}i})
= \Multiply{(\Vector{#1}\circ \ACol vi)}{\Multiply{\AUD {#2}ki}{\EBase {#4}k}}
}

\AddEq[4]{f o ev=efv i (1)}
{
\Vector{#2}\circ(\Multiply{\ACol vi}\Multiply{\EBase {#3}i})
= \Multiply{\ACol vi}{\Multiply{\AUD {#2}ki}{\EBase {#4}k}}
}

\AddEq[4]{f o ea=efa (11)(cols)}
{
\Vector{#2}\circ(e_{#3}a)
= e_{#4}#2#1(a)
}

\AddEq[4]{f o ea=efa (11)(rows)}
{
\Vector{#2}\circ(a e_{#3})
= #1(a)#2e_{#4}
}

\AddEq[4]{f o ea=efa (1)(cols)}
{
\Vector{#2}\circ(e_{#3}a)
= e_{#4}#2a
}

\AddEq[4]{f o ea=efa (1)(rows)}
{
\Vector{#2}\circ(ae_{#3})
= a#2e_{#4}
}

\DefRef{Morphism of Diagram of Representations}
{
\RefDefinition{morphism of representations of universal algebra},
\RefDefinition{Morphism of Diagram of Representations},
}

\AddEquation{vb=h(e1)a (1)(1)(1)(cols)left}
{
\Vector b=\Vector f\circ\Vector a=\Vector f\circ(a\CRstar e_1)
= h(a)\CRstar(\Vector f\circ e_1)
}

\AddEquation{vb=h(e1)a (1)(1)(1)(cols)right}
{
\Vector b=\Vector a\diamond\Vector f=(e_1\RCstar a)\diamond\Vector f
= (e_1\diamond\Vector f)\RCstar h(a)
}

\AddEquation{vb=h(e1)a (1)(1)(1)(rows)left}
{
\Vector b=\Vector f\circ\Vector a=\Vector f\circ(a\CRstar e_1)
= h(a)\CRstar(\Vector f\circ e_1)
}

\AddEquation{vb=h(e1)a (1)(1)(1)(rows)right}
{
\Vector b=\Vector a\diamond\Vector f=(e_1\RCstar a)\diamond\Vector f
= (e_1\diamond\Vector f)\RCstar h(a)
}

\AddEq{f(e1)(cols)}
{
$\Vector f\circ\aD[1]ei$
}

\AddEq{f(e1)(rows)}
{
$\Vector f\circ\aU{e_1}i$
}

\AddEq{f(e1)(cols)left}
{
$\Vector f\circ\aD[1]ei$
}

\AddEq{f(e1)(cols)right}
{
$\aD[1]ei\diamond\Vector f$
}

\AddEq{f(e1)(rows)left}
{
$\Vector f\circ\aU{e_1}i$
}

\AddEq{f(e1)(rows)right}
{
$\aU{e_1}i\diamond\Vector f$
}

\AddEq{f(e1)=(cols)}
{
\Vector f\circ\aD[1]ei
=e_2\aD fi
=\aD[2]ej\aUD fji
}

\AddEq{f(e1)=(cols)left}
{
\Vector f\circ\aD[1]ei
=\aD fi\CRstar e_2
=\aUD fji\aD[2]ej
}

\AddEq{f(e1)=(cols)right}
{
\aD[1]ei\diamond\Vector f
=e_2\RCstar\aD fi
=\aD[2]ej\aUD fji
}

\AddEq{f(e1)=(rows)}
{
\Vector f\circ\aU{e_1}i
=\aD fie_2
=\aUD fij \aU{e_2}j
}

\AddEq{f(e1)=(rows)left}
{
\Vector f\circ\aU{e_1}i
=\aD fi\RCstar e_2
=\aUD fij\aU{e_2}j
}

\AddEq{f(e1)=(rows)right}
{
\aU{e_1}i\diamond\Vector f
=e_2\CRstar\aD fi
=\aU{e_2}j\aUD fij
}

\AddEquation{vb=a f e (1)(1)(cols)left}
{
\Vector b=h(a)\CRstar f\CRstar e_2
}

\AddEquation{vb=a f e ()(1)(cols)left}
{
\Vector b=a\CRstar f\CRstar e_2
}

\AddEquation{vb=a f e ()(1)(cols)}
{
\Vector b=e_2fa
}

\AddEquation{vb=a f e (1)(1)(cols)}
{
\Vector b=e_2fh(a)
}

\AddEquation{vb=a f e ()(1)(cols)right}
{
\Vector b=e_2\RCstar f\RCstar a
}

\AddEquation{vb=a f e (1)(1)(cols)right}
{
\Vector b=e_2\RCstar f\RCstar h(a)
}

\AddEquation{vb=a f e (1)(1)(rows)left}
{
\Vector b=h(a)\RCstar f\RCstar e_2
}

\AddEquation{vb=a f e ()(1)(rows)left}
{
\Vector b=a\RCstar f\RCstar e_2
}

\AddEquation{vb=a f e ()(1)(rows)}
{
\Vector b=afe_2
}

\AddEquation{vb=a f e (1)(1)(rows)}
{
\Vector b=h(a)fe_2
}

\AddEquation{vb=a f e ()(1)(rows)right}
{
\Vector b=e_2\CRstar f\CRstar a
}

\AddEquation{vb=a f e (1)(1)(rows)right}
{
\Vector b=e_2\CRstar f\CRstar h(a)
}

\AddEq [4]{h(a)=...}
{
$#2(#1)=(#2(\ARow {#1}{#3}),\gi #3\in\gi #4)$
}

\AddEq[4]{Vector f(e1) module}
{
$(\Vector #1\circ\EBase {#2}{#3},\jJg{#3}{#4})$
}

%% file: Stmt.Module.Ref.tex

\AddEq{ref theorem coordinates of map A1 A2, algebra}
{
\ePrints{9835-2163,0767-8264}
\ifx\Semafor\ValueOn
\RefTheorem{coordinates of map A1 A2 FoG, algebra}.
\else
\ePrints{5284-0163,1801.01628,CACAA.07}%
\ifx\Semafor\ValueOn
\RefTheorem{coordinates of map A1 A2, algebra}.
\else
\RefTheorem[\RefPolymorphism]{coordinates of map A1 A2, algebra}.
\fi
\fi
}

\AddEq{RefTheorem set of vectors generated by set of vectors, module}
{
\ePrints{1502.04063,5114-6019}%
\ifx\Semafor\ValueOn%
\RefTheorem{set of vectors generated by the set of vectors},
\else%
\refTheorem{set of vectors generated by set of vectors}{module},
\fi%
}

\AddEq{ref definition: basis of representation}
{
\ePrints{1502.04063,5114-6019,MAlgebra2}%
\ifx\Semafor\ValueOn%
\RefDefinition[0912.3315]{basis of representation}
\else%
\RefDefinition{basis of representation}
\fi
}

\AddEq{f=.I0.7 ref}
{
\EqRef{linear map of algebra O, structure, 1},
\eqRef{fi=fij ej}O,
\EqRef{a0=f07},
\EqRef{a1=f07},
\EqRef{a2=f07},
\EqRef{a3=f07},
\EqRef{a4=f07},
\EqRef{a5=f07},
\EqRef{a5=f07},
\EqRef{a7=f07}.
}

\AddEq{f=.E+.I+.J+.K ref}
{
\EqRef{linear map of algebra H, structure, 1},
\eqRef{fi=fij ej}H,
\EqRef{a0=f03},
\EqRef{a1=f03},
\EqRef{a2=f03},
\EqRef{a3=f03}.
}

\AddEq{f=.E+.I ref}
{
\EqRef{linear map of algebra C, structure, 1},
\eqRef{fi=fij ej}C,
\EqRef{a0=f01},
\EqRef{a1=f01}.
}

\DefRef[3]{homomorphism of d algebra, coordinates}
{
\newline
\FrameEqRef[f{#3}{}]{f:V1->V2, D module \Cols}{algebra(#1#2)}
\newline
\FrameEqRef[hf12a]{f o ea=efa i (#1#2)}{(\Cols)algebra}
\newline
\FrameEqRef[hf12]{f o ea=efa (#1#2)(\Cols)}{algebra}
\newline
}

\DefRef{da in DA->da in D1A}
{
\newline
\FrameEqRef{distributive law, \SideWS module, 2}1
\newline
}

\DefRef[3]{algebra, homomorphism and product}
{
\newline
\FrameEqRef[{#3}{}]{algebra, homomorphism and product (#1#2)}{\Cols(\SideNS)}
\newline
}

\DefRef[3]{matrix generates A module homomorphism()}
{
\newline
\FrameEqRef{\RefDistributive{\Product}}1
\FrameEqRef[fa{#3}v]{\SideWS homomorphism, f av=}{(#1#2)\SideWS A module}
\FrameEqRef[vf{V_1}{V_2}]{f o (ae)=ga o f e, vector space \Product-\Cols}{\SideNS(#1#21)}
\newline
}

\DefRef[3]{matrix generates A module homomorphism(1)}
{
\newline
\FrameEqRef{\RefDistributive{\Product}}1
\FrameEqRef[gab]{homomorphism, f v+w=}{g(#1#2)\SideWS A module}
\FrameEqRef[fa{#3}v]{\SideWS homomorphism, f av=}{(#1#2)\SideWS A module}
\FrameEqRef[vf{V_1}{V_2}]{f o (ae)=ga o f e, vector space \Product-\Cols}{\SideNS(#1#21)}
\newline
}

\DefRef[3]{fo(va)()}
{
\newline
\FrameEqRef[fa{#3}v]{\SideWS homomorphism, f av=}{(#1#2)\SideWS A module}
\newline
\FrameEqRef[vf{V_1}{V_2}]{f o (ae)=ga o f e, vector space \Product-\Cols}{\SideNS(111)}
\newline
}

\DefRef[3]{fo(va)(1)}
{
\newline
\FrameEqRef[fa{#3}v]{\SideWS homomorphism, f av=}{(#1#2)\SideWS A module}
\newline
\FrameEqRef[gab]{homomorphism, f vw=}{(#1#2)g}
\newline
\FrameEqRef[vf{V_1}{V_2}]{f o (ae)=ga o f e, vector space \Product-\Cols}{\SideNS(111)}
\newline
}

\DefRef[2]{homomorphism, f vw=}
{
\newline
\FrameEqRef[fab]{homomorphism, f vw=}{(#1#2)g}
\newline
}

\DefRef[2]{homomorphism of d algebra, coordinates 1}
{
\newline
\FrameEqRef[hf12a]{f o ea=efa i (#1#2)}{(\Cols)algebra}
\newline
}

\DefRef{module over algebra}
{
\refDefinition{module over algebra}{\VectorSetNS}
}

\DefRef{left module over algebra}
{
\refDefinition{module over associative algebra}{left \VectorSetNS}
}

\DefRef{right module over algebra}
{
\refDefinition{module over associative algebra}{right \VectorSetNS}
}

\DefRef{set of maps B->A is Abelian group}
{
\ShowEq{def DModule}%
\refDefinition{module over algebra}{\SideWS \VectorSetNS},
\ShowEq{def DAlgebra}%
\RefDefinition{algebra over ring}.
}

\DefRef[3]{homomorphism of D module}
{
\newline
\FrameEqRef[fab]{homomorphism, f v+w=}{f D module #1#2}
\newline
\FrameEqRef[fp{#3}a]{left homomorphism, f av=}{f D module #1#2}
\newline
}

\DefRef[3]{vb=h(e1)a}
{
\newline
\FrameEqRef[fp{#3}a]{left homomorphism, f av=}{f D module #1#2}
\FrameEqRef[a1{}]{va=ae1, module ()(\Cols)}{\SideNS-\VectorSet a (#1)(#2)}
\newline
}

\DefRef[2]{algebra, homomorphism and product 1}
{
\eqRef{algebra, homomorphism and product abe #1#2}{\Cols},
\eqRef{algebra, homomorphism and product abe 3(#1)}{\Cols},
\eqRef{algebra, homomorphism and product abe 5(#1)}{\Cols}.
}

\DefRef{product in algebra}
{
\newline
\FrameEqRef{product in algebra}{\Cols}
\newline
}

\DefRef[3]{algebra, homomorphism and product eiej}
{
\newline
\FrameEqRef{product of basis vectors, algebra}{\Cols}
\newline
\FrameEqRef[fp{#3}a]{left homomorphism, f av=}{(#1#2)g}
\newline
}

\DefRef[2]{algebra, homomorphism and product eiej 1}
{
\newline
\FrameEqRef{f(e1)=(\Cols)\SideNS}{(#1)(#2)}
\newline
}

\AddEq{ref linear map of D1 D2 module, coordinates}
{
\ePrints{9835-2163,0767-8264}
\ifx\Semafor\ValueOn
\refTheorem{linear map of D module, coordinates}1,
\else
\refTheorem[\RefLinearMap]{linear map of D module}1,
\fi
}

\AddEq{ref module, (a+n)(b+m)=}
{
\EqRef{\SideWS module, a d v},
\EqRef{\SideWS module, (a+n)v=av+nv},
\eqRef{associative law, \SideWS module}1,
\eqRef{associative law, \CBase, \SideWS module}1,
\eqRef{\SideWS module, associative law}1,
\eqRef{associative law, \CBase\Base, \SideWS module}1.
}

\AddEq{def sum of linear maps}
{
\ifx\SelectlEnglish\undefined
\ePrints{1502.04063}
\ifx\Semafor\ValueOn
\def\defTheorem{Согласно теореме}
\else
\def\defTheorem{Согласно следствию}
\fi
\else
\ePrints{1502.04063}
\ifx\Semafor\ValueOn
\def\defTheorem{According to the theorem}
\else
\def\defTheorem{According to the corollary}
\fi
\fi
}

\DefEq
{
\defTheorem
\ePrints{1502.04063}
\ifx\Semafor\ValueOn
\RefTheorem{sum of linear maps, D module},
\else
\RefCorollary{sum of linear maps, D module},
\fi
}
{ref sum of linear maps}

\AddEq [1]{ref definition: representation of algebra}
{
\refDefinition[\RefRepresentation]{side representation of group}{left}#1
}

\AddEq{ref sum and product over scalar, linear map}
{
\ifx\SelectlEnglish\undefined
\ePrints{1502.04063}
\ifx\Semafor\ValueOn
согласно теоремам
\else
согласно следствиям
\fi
\else
\ePrints{1502.04063}
\ifx\Semafor\ValueOn
according to theorems
\else
according to corollarys
\fi
\fi
\ePrints{1502.04063}
\ifx\Semafor\ValueOn
\RefTheorem{sum of linear maps, D module},
\RefTheorem{product of linear map over scalar, D module},
\else
\RefCorollary{sum of linear maps, D module},
\RefCorollary{product of linear map over scalar, D module},
\fi
}

\AddEq{ref map j, representation, tensor product}
{
\ePrints{MAlgebra2}%
\ifx\Semafor\ValueOn%
\EqRef[\RefPolymorphism]{map j, representation, tensor product},
\else%
\ePrints{2207.06506,8428-0408}%
\ifx\Semafor\ValueOn%
\EqRef[\RefLinearMap]{map j, representation, tensor product},
\else
\EqRef{map j, representation, tensor product},
\fi
\fi
}

%% file: Vector.Space.2020.Stmt.English.tex
\input{Vector.Space.2020.Stmt.Eq}

\DefLabeledDefinition{eigenvalue of matrix}{\Product}
{
$A$\Hyph number $b$ is called
{\bf\ProductType eigenvalue}
of the matrix $f$
if the matrix
$f-b\aD En$
is \ProductType singular matrix.
}

\DefLabeledDefinition{eigenvector of matrix}{\Product-\Cols}
{
Let $A$\Hyph number $b$ be
\ProductType eigenvalue
of the matrix $f$.
A \ColWS $v$ is called
{\bf eigen\ColWS}
of matrix $f$ corresponding to \ProductType eigenvalue $b$,
if the following equality is true
\DrawEq{\Product-eigen\Cols}1
}

\DefLabeledTheorem{0 is eigenvalue of multiplicity}{\Product}
{
Let $f$ be $\gin\times\gin$ matrix of $A$\Hyph numbers and
\ShowEq{\Product-rank a=k<n}
Then $0$ is \ProductType eigenvalue of multiplicity
$\gin-\gik$.
}

\DefProof{0 is eigenvalue of multiplicity}
{
The theorem follows from the equality
\[a=a-0E\]
from the definition
\refDefinition{eigenvalue of matrix}{\Product}
and the theorem
\RefTheorem{star rows system of linear equations, solution}.
}

\DefLabeledTheorem{eigenvalues of pair of matrices}{\Product}
{
Let $f$ be $\gin\times\gin$ matrix of $A$\Hyph numbers.
Let $g$ be non\Hyph \ProductType singular $\gin\times\gin$ matrix of $A$\Hyph numbers.
Then non\Hyph zero \ProductType eigenvalues of the pair of matrices $(f,g)$
are roots of any equation
\DrawEq{\Product-det ij f-bg =0}{}
}

\DefProofSloppy{eigenvalues of pair of matrices}
{
The theorem follows from the theorem
\refTheorem{singular matrix and quasideterminant}{\Product}
and from the definition
\refDefinition{eigenvalues of pair of matrices}{\Product}.
}

\DefLabeledTheorem{eigenvalues of matrix}{\Product}
{
Let $f$ be $\gin\times\gin$ matrix of $A$\Hyph numbers.
Then non\Hyph zero \RC eigenvalues of the matrix $f$
are roots of any equation
\DrawEq{\Product-det ij a=0}{}
}

\DefProofSloppy{eigenvalues of matrix}
{
The theorem follows from the theorem
\refTheorem{singular matrix and quasideterminant}{\Product}
and from the definition
\refDefinition{Eigenvalue of Endomorphism}{\SideNS}.
}

\DefLabeledTheorem[6]{product of homomorphisms, A vector space}{\SideNS-\ColN}
{
Let $U$, $V$, $W$ be
\SideWS $A$\Hyph vector spaces of \ColN s.
\ShowEq{Let be Basis of vector space}UAU
\ShowEq{Let be Basis of vector space}VAV
\ShowEq{Let be Basis of vector space}WAW
The homomorphism
\ShowEq{vf: A->B}fUW{}%
is product of homomorphisms
\ShowEq{vf: A->B}hVW,%
\ShowEq{vf: A->B}gUV{}%
iff
matrix $f$ of homomorphism $\Vector f$
relative to bases
\ShowEq{Bases eVW}UW{}{}
is equal to the \ProductType product of matrix $#1$ of the homomorphism $\Vector{#1}$
relative to bases
\ShowEq{Bases eVW}{#2}{#3}{}{}
over matrix $#4$ of the homomorphism $\Vector{#4}$
relative to bases
\ShowEq{Bases eVW}{#5}{#6}{}{}
\DrawEq[f#1{\ProductVal}#4]{f=g*h}{\Product-\Cols}
}

\DefProof[2]{product of homomorphisms, A vector space}
{
According to the theorem
\refTheorem{homomorphism A module}{\SideNS-\Cols(1)}
\DrawEq[ufUW]{f o (ae)=a o f e, vector space \Product-\Cols}{f product \Product-\Cols}
\DrawEq[ugUV]{f o (ae)=a o f e, vector space \Product-\Cols}{g product \Product-\Cols}
\DrawEq[vhVW]{f o (ae)=a o f e, vector space \Product-\Cols}{h product \Product-\Cols}
According to the theorem
\refTheorem{product of homomorphisms is homomorphism, A vector space}{\SideNS},
the equality
\ShowEq{f o v=g o h o v \SideNS-\Cols}
follows from equalities
\eqRef{f o (ae)=a o f e, vector space \Product-\Cols}{g product \Product-\Cols},
\eqRef{f o (ae)=a o f e, vector space \Product-\Cols}{h product \Product-\Cols}.
The equality
\eqRef{f=g*h}{\Product-\Cols}
follows from equalities
\eqRef{f o (ae)=a o f e, vector space \Product-\Cols}{f product \Product-\Cols},
\EqRef{f o v=g o h o v \SideNS-\Cols}
and the theorem
\refTheorem{homomorphism A module}{\SideNS-\Cols(1)}.

Let
\DrawEq[f#1{\ProductVal}{#2.}]{f=g*h}-
According to the theorem
\refTheorem{matrix generates A module homomorphism}{\SideNS-\Cols(1)},
there exist homomorphisms
\ShowEq{vf: A->B}fUV,%
\ShowEq{vf: A->B}gUV,%
\ShowEq{vf: A->B}hUV{}%
such that
\ShowEq{fov=(gh)ov \SideNS-\Cols}
From the equality
\EqRef{fov=(gh)ov \SideNS-\Cols}
and the theorem
\refTheorem{sum of homomorphisms is homomorphism, A vector space}{\SideNS},
it follows that
the homomorphism $\Vector f$
is product of homomorphisms $\Vector h$, $\Vector g$.
}

\DefLabeledTheorem{product of homomorphisms is homomorphism, A vector space}{\SideNS}
{
Let $U$, $V$, $W$ be \SideWS $A$\Hyph vector spaces.
Let diagram of maps
\DrawEq{diagram product of homomorphisms, A vector space}{\SideNS}
\DrawEq{f o u=h o g o u}{\SideNS}
be commutative diagram
where maps $g$, $h$ are
homomorphisms of \SideWS $A$\Hyph vector space.
The map
\ShowEq{f=h o g}
is homomorphism
of \SideWS $A$\Hyph vector space
and is called
\AddIndex{product of homomorphisms}{product of homomorphisms}
$h$, $g$.
}

\DefProof{product of homomorphisms is homomorphism, A vector space}
{
The equality
\eqRef{f o u=h o g o u}{\SideNS}
follows from the commutativity of the diagram
\eqRef{diagram product of homomorphisms, A vector space}{\SideNS}.
The equality
\DrawEq{(h o g)o(u+v)=}{\SideWS vector space}
follows from the equality
\ShowRef{homomorphism, f v+w=}
and the equality
\eqRef{f o u=h o g o u}{\SideNS}.
The equality
\ShowEq{(h o g)o(av)= \SideNS}
follows from the equality
\ShowRef{homomorphism, f av=}
and the equality
\eqRef{sum of homomorphisms f o v=}{\SideWS vector space}.
The theorem follows from the theorem
\refTheorem{define homomorphism A module}{\SideNS(1)}
and equalities
\eqRef{(h o g)o(u+v)=}{\SideWS vector space},
\EqRef{(h o g)o(av)= \SideNS}.
}

\DefLabeledTheorem{sum of homomorphisms is homomorphism, A vector space}{\SideNS}
{
Let $V$, $W$ be \SideWS $A$\Hyph vector spaces
and maps
\DrawEq[gVW,]{f: A->B}-
\DrawEq[hVW{}]{f: A->B}-
be homomorphisms
of \SideWS $A$\Hyph vector space.
Let the map
\DrawEq[fVW{}]{f: A->B}-
be defined by the equality
\DrawEq{sum of homomorphisms f o v=}{\SideWS vector space}
The map $f$ is homomorphism
of \SideWS $A$\Hyph vector space
and is called
\AddIndex{sum
\ShowEq{f=g+h}
of homomorphisms}{sum of maps}
$g$ and $h$.
}

\DefProof{sum of homomorphisms is homomorphism, A vector space}
{
According to the theorem
\refTheorem{definition of A module, property}{\SideNS}
(the equality
\eqRef{commutative law}{\SideWS vector space}),
the equality
\newline
\DrawEq{(g+h)o(u+v)=}{\SideWS vector space}
\begin{sloppypar}
\noindent
follows from the equality
\ShowRef{homomorphism, f v+w=}
and the equality
\eqRef{sum of homomorphisms f o v=}{\SideWS vector space}.
According to the theorem
\refTheorem{definition of A module, property}{\SideNS}
(the equality
\eqRef{distributive law, \SideWS module, 1}1),
the equality
\ShowEq{(g+h)o(av)= \SideNS}
\end{sloppypar}
\noindent
follows from the equality
\ShowRef{homomorphism, f av=}
and the equality
\eqRef{sum of homomorphisms f o v=}{\SideWS vector space}.

The theorem follows from the theorem
\refTheorem{define homomorphism A module}{\SideNS(1)}
and equalities
\eqRef{(g+h)o(u+v)=}{\SideWS vector space},
\EqRef{(g+h)o(av)= \SideNS}.
}

\DefLabeledTheorem{sum of homomorphisms, A vector space}{\SideNS-\Cols}
{
Let $V$, $W$ be
\SideWS $A$\Hyph vector spaces of \ColN s.
\ShowEq{Let be Basis of vector space}VAV
\ShowEq{Let be Basis of vector space}WAW
The homomorphism
\ShowEq{vf: A->B}fVW{}%
is sum of homomorphisms
\ShowEq{vf: A->B}gVW,%
\ShowEq{vf: A->B}hVW{}%
iff
matrix $f$ of homomorphism $\Vector f$
relative to bases
\ShowEq{Bases eVW}VW{}{}
is equal to the sum of the matrix $g$ of the homomorphism $\Vector g$
relative to bases
\ShowEq{Bases eVW}VW{}{}
and the matrix $h$ of the homomorphism $\Vector h$
relative to bases
\ShowEq{Bases eVW}VW{}.
}

\DefProof{sum of homomorphisms, A vector space}
{
According to theorem
\refTheorem{homomorphism A module}{\SideNS-\Cols(1)}
\DrawEq[afUV]{f o (ae)=a o f e, vector space \Product-\Cols}{f sum \Product-\Cols}
\DrawEq[agUV]{f o (ae)=a o f e, vector space \Product-\Cols}{g sum \Product-\Cols}
\DrawEq[ahUV]{f o (ae)=a o f e, vector space \Product-\Cols}{h sum \Product-\Cols}
According to the theorem
\refTheorem{sum of homomorphisms is homomorphism, A vector space}{\SideNS},
the equality
\ShowEq{fov=gov+hov \SideNS-\Cols}
follows from equalities
\eqRef{a.\Product.(b1+b2)=}1,
\eqRef{(b1+b2).\Product.a=}1,
\eqRef{f o (ae)=a o f e, vector space \Product-\Cols}{g sum \Product-\Cols},
\eqRef{f o (ae)=a o f e, vector space \Product-\Cols}{h sum \Product-\Cols}.
The equality
\ShowEq{f=g+h}
follows from equalities
\eqRef{f o (ae)=a o f e, vector space \Product-\Cols}{f sum \Product-\Cols},
\EqRef{fov=gov+hov \SideNS-\Cols}
and the theorem
\refTheorem{homomorphism A module}{\SideNS-\Cols(1)}.

Let
\ShowEq{f=g+h}
According to the theorem
\refTheorem{matrix generates A module homomorphism}{\SideNS-\Cols(1)},
there exist homomorphisms
\ShowEq{vf: A->B}fUV,%
\ShowEq{vf: A->B}gUV,%
\ShowEq{vf: A->B}hUV{}%
such that
\ShowEq{fov=(g+h)ov \SideNS-\Cols}
From the equality
\EqRef{fov=(g+h)ov \SideNS-\Cols}
and the theorem
\refTheorem{sum of homomorphisms is homomorphism, A vector space}{\SideNS},
it follows that
the homomorphism $\Vector f$
is sum of homomorphisms $\Vector g$, $\Vector h$.
}

\AddEq{definition: vector space of maps B->A n}
{
\begin{ShadedDefinition}
\labelDefinition{\SideWS vector space of maps B->A n}
Let $D$ be commutative ring with unit.
Let $A$ be associative division $D$\Hyph algebra.
For any set $B$
and any integer $n>0$,
we introduce \SideWS vector space
\ShowEq{\SideWS vector space of maps B->A n}
using definition by induction
\ShowEq{\SideWS vector space of maps B->A k}
\end{ShadedDefinition}
}

\DefLabeledDefinition{vector space of algebra A n}{\SideNS}
{
Let $D$ be commutative ring with unit.
Let $A$ be associative division $D$\Hyph algebra.
For any integer $n>0$,
we introduce \SideWS vector space
\ShowEq{\SideWS vector space of algebra A n}
using definition by induction
\ShowEq{\SideWS vector space of algebra A k}
}

\DefTheorem{coordinates of vector of vector space}
{
Coordinates of vector $v\in V$ relative to basis $\Basis e$
of left $\Base$\Hyph vector space $V$
are uniquely defined.
}

\DefLabeledTheorem{coordinates of vector}{\SideNS-\Cols}
{
Coordinates of vector $v\in V$ relative to basis $\Basis e$
of \SideWS \Free $\Base$\Hyph \VectorSet $V$
are uniquely defined.
The equality
\DrawEq{e*v=e*w,\SideNS}{\Cols}
implies the equality
\ShowEq{=> v=w}
}

\DefProof{coordinates of vector}
{
According to the theorem
\refTheorem{basis over division algebra}{\SideNS},
the system of vectors
\ShowEq{expansion relative basis, vector space, 0}
is linearly dependent and in equality
\DrawEq{expansion relative basis, vector space, 1}{\SideNS-\Cols}
at least $b$ is different from $0$. Then equality
\DrawEq{expansion relative basis, vector space, 2}{\SideNS-\Cols}
follows from
\eqRef{expansion relative basis, vector space, 1}{\SideNS-\Cols}.
The equality
\DrawEq{expansion relative basis, vector space}{\SideNS-\Cols}
follows from the equality
\eqRef{expansion relative basis, vector space, 2}{\SideNS-\Cols}.

Assume we get another expansion
\DrawEq{expansion relative basis, vector space, 3}{\SideNS-\Cols}
We subtract
\eqRef{expansion relative basis, vector space}{\SideNS-\Cols}
from
\eqRef{expansion relative basis, vector space, 3}{\SideNS-\Cols}
and get
\DrawEq{expansion relative basis, vector space, 4}{\SideNS-\Cols}
Since vectors \ARow ei are linearly independent,
then the equality
\DrawEq{expansion relative basis, vector space, 5}{\SideNS-\Cols}
follows from the equality
\eqRef{expansion relative basis, vector space, 4}{\SideNS-\Cols}.
Therefore, the theorem follows from the equality
\eqRef{expansion relative basis, vector space, 5}{\SideNS-\Cols}.
}

\DefProof{coordinates of vector 2022}
{
According to the theorem
\refTheorem{basis over division algebra}{\SideNS},
the system of vectors
\ShowEq{expansion relative basis, vector space, (0)}
is linearly dependent and in equality
\ShowEq{expansion relative basis, vector space, 1 \SideNS}
at least $b$ is different from $0$. Then equality
\ShowEq{expansion relative basis, vector space, 2 \SideNS}
follows from
\EqRef{expansion relative basis, vector space, 1 \SideNS}.
The equality
\ShowEq{expansion relative basis, vector space \SideNS}
follows from the equality
\EqRef{expansion relative basis, vector space, 2 \SideNS}.

Assume we get another expansion
\ShowEq{expansion relative basis, vector space, 3 \SideNS}
We subtract
\EqRef{expansion relative basis, vector space \SideNS}
from
\EqRef{expansion relative basis, vector space, 3 \SideNS}
and get
\ShowEq{expansion relative basis, vector space, 4 \SideNS}
Since vectors \ARow ei are linearly independent,
then the equality
\DrawEq{expansion relative basis, vector space, 5()}{\SideNS}
follows from the equality
\EqRef{expansion relative basis, vector space, 4 \SideNS}.
Therefore, the theorem follows from the equality
\eqRef{expansion relative basis, vector space, 5()}{\SideNS}.
}

\DefProof{coordinates of vector *}
{
The theorem follows from the theorem
\refTheorem{coordinates of vector}{\SideNS-*}.
}

\DefExample{system of linear equations, right vector space of columns}
{
Let $V$ be a right $A$\Hyph vector space of columns and
row
\ShowEq{rcd linear span, 0}
be set of vectors.
The vector $\Vector b$
linearly dependends on vectors $\aD{\Vector a}i$,
if linear equation
\DrawEq{rcd linear span, 2}f
where
\DrawEq[xn]{a=(a1.n col)}{}
is column of unknown coefficients of expansion
has a solution.
Suppose
\DrawEq[{}jJ]{basis e of module cols}-
is a basis.
Then vectors
\ShowEq{rcd linear span, 6}
have expansion
\ShowEq{rcd linear span, 7}
If we substitute \EqRef{rcd linear span, b}
and \EqRef{rcd linear span, ai}
into
\eqRef{rcd linear span, 2}f
we get
\DrawEq{rcd linear span, 8}f
Applying theorem \RefTheorem{coordinates of vector of vector space}
to
\eqRef{rcd linear span, 8}f
we get
\AddIndex{system of linear equations}
{system of linear equations}
\DrawEq{system of rcd linear equations}f

Let
\ShowEq{a1n= b= column}
Then we can write system of linear equations
\eqRef{system of rcd linear equations}f
in the following form
\ShowEq{star rows system of linear equations 1}
The equality
\ShowEq{star rows system of linear equations}
follows from equalities
\ePrints{339295394}%
\ifx\Semafor\ValueOn%
\ShowEq{rc-product}
\ShowEq{rc-product of matrices}
and
\else%
\EqRef{rc-product of matrices},
\fi%
\EqRef{star rows system of linear equations 1}.
If we write the sum
\EqRef{star rows system of linear equations}
as row, we will get
\DrawEq{star rows system of linear equations 2}2
We see in the system of linear equations
\eqRef{star rows system of linear equations 2}2,
that $A$\Hyph column $b$ is right linear composition of $A$\Hyph columns
\ShowEq{aD1n}
}

\DefLabeledDefinition{rank of matrix}{\Product}
{
If submatrix $\aUD aST$ is \ProductType nonsingular matrix
then we say that
\ProductType rank of matrix $a$ is not less then $\gik$.
{\bf\ProductType rank of matrix} $a$,
\ShowEq{\Product-rank of matrix}
is the maximal value of $\gik$.
We call an appropriate submatrix the
{\bf \ProductType major submatrix}.
}

\DefLabeledTheorem{vector space of maps B->A}{\SideNS}
{
Let $D$ be commutative ring with unit.
Let $A$ be associative division $D$\Hyph algebra.
For any set $B$,
the representation
\ShowEq{representation \SideWS vector space of maps B->A}
is \SideWS vector space
\ShowEq{\SideWS vector space of maps B->A}
over $D$\Hyph algebra $A$.
}

\DefProof{vector space of maps B->A}
{
According to the theorem
\RefTheorem{set of maps B->A is Abelian group},
the set of maps $A^B$ is Abelian group.
From the equality
\ShowEq{\SideWS a(h+g)=}
and definitions
\RefDefinition{homomorphism},
\RefDefinition{endomorphism},
it follows that the map $f(a)$ is endomorphism of Abelian group $A^B$.
From equalities
\ShowEq{\SideWS (a1+a2)h=}
\ShowEq{\SideWS a1a2h=}
and definitions
\RefDefinition{homomorphism},
\RefDefinition{representation of algebra},
it follows that the map $f$ is
representation of $D$\Hyph algebra $A$ in Abelian group $A^B$.
The theorem follows from the definition
\refDefinition{module over associative algebra}{\SideWS \VectorSetNS}.
}

\DefText{linear span, vector space}
{
\subsection{\SideWSC
\texorpdfstring{$A$}{A}-vector space of \ColN s}
\ShowText{linear span in vector space}
}

\DefLabeledDefinition{coordinates of vector}{\SideWS \VectorSetNS}
{
Let $\Basis e$ be the quasibasis of \SideWS $\Base$\Hyph \VectorSet $V$
and vector
\ShowEq{vv in V}vV
has expansion
\DrawEq{vv=ve \SideWS module}{Definition}
with respect to the quasibasis $\Basis e$.
$\Base_{\BaseExt}$\Hyph numbers $\ACol vi$ are called
\AddIndex{coordinates}{coordinates}
of vector $\Vector v$ with respect to the quasibasis $\Basis e$.
Matrix of $\Base_{\BaseExt}$\Hyph numbers
\ShowEq{coordinate matrix of vector}
is called
\AddIndex{coordinate matrix of vector}{coordinate matrix of vector}
$\Vector v$ in quasibasis $\Basis e$.
}

\DefLabeledTheorem{matrix and system of linear equations}{\SideNS-\Cols}
{
Consider the matrix
\DrawEq{matrix a = set ai \Cols}{\SideNS}
Then the system of linear equations
\eqRef{a*x=b \SideNS}{\Cols}
can be represented as
\ShowEq{a*x=b \SideNS-\Cols}
\DrawEq{a*x=b 1 \SideNS-\Cols}1
}

\DefProof{matrix and system of linear equations}
{
The equality
\EqRef{a*x=b \SideNS-\Cols}
follows from equalities
\eqRef{linear span, b \SideNS}{b \Cols},
\eqRef{a=(a1.n set)}{xm \SideNS-\Cols},
\eqRef{a*x=b \SideNS}{\Cols},
\eqRef{matrix a = set ai \Cols}{\SideNS}
and from the definition
\RefDefinition{\RowNWS over \ColNWS product}.
}

\DefLabeledTheoremNote{nonsingular system of linear equations}{\SideNS-\Cols}
{
Solution of nonsingular system of linear equations
\newline
\FrameEqRef{a*x=b 1 \SideNS-\Cols}1
\newline
is determined uniquely and can be presented
in either form\,\footnotemark
\DrawEq{x=a-*b, matrix \SideNS}{\Cols}
\DrawEq{x=a-*b, quasideterminant \SideNS}{\Cols}
}
{
We can see a solution of system
\eqRef{a*x=b 1 \SideNS-\Cols}1
in theorem
\citeBib{math.QA-0208146}-\href{http://arxiv.org/PS_cache/math/pdf/0208/0208146.pdf\#Page=19}{1.6.1}.
I repeat this statement because I slightly changed the notation.
}

\DefProof{nonsingular system of linear equations}
{
Multiplying both sides of the equation
\eqRef{a*x=b 1 \SideNS-\Cols}1
from left by $a^{\InverseVal}$, we get
\eqRef{x=a-*b, matrix \SideNS}{\Cols}.
Using definition
\newline
\FrameEqRef{quasideterminant and inverse}{\Product}
\newline
we get
\eqRef{x=a-*b, quasideterminant \SideNS}{\Cols}.
Based on the theorem
\RefTheorem{two rc-products equal},
the solution is unique.
}

\DefLabeledDefinition{nonsingular system of linear equations}{\SideNS-\Cols}
{
Let $a$
be \ProductType nonsingular matrix. We call appropriate
system of linear equations
\newline
\FrameEqRef{a*x=b 1 \SideNS-\Cols}1
\newline
\AddIndex{nonsingular system of linear equations}
{nonsingular system of linear equations}.
}

\DefLabeledTheorem{linear span and system of equations}{\SideNS-\Cols}
{
Let
\DrawEq{basis e of module}-
be a basis.
Let vectors
\ShowEq{linear span, 6}
have expansion
\DrawEq[b]{linear span, b \SideNS}{b \Cols}
\DrawEq[{\ARow ai}{}]{linear span, b \SideNS}{a \Cols}
The vector $\Vector b$
linearly dependends on vectors $\ARow{\Vector a}i$
\ShowEq{linear span, 1}
if there exist \ColWS of $A$\Hyph numbers
\DrawEq[xm]{a=(a1.n set)}{xm \SideNS-\Cols}
which satisfy the
\AddIndex{system of linear equations}
{system of linear equations}
\DrawEq{a*x=b \SideNS}{\Cols}
}

\DefProof{linear span and system of equations}
{
The equality
\DrawEq{e*b=e*a*x \SideNS}{\Cols}
follows from equalities
\eqRef{linear span, 2 \SideNS}{\SideNS-\Cols},
\eqRef{linear span, b \SideNS}{b \Cols},
\eqRef{linear span, b \SideNS}{a \Cols}.
The equality
\DrawEq{b=a*x \SideNS}{\Cols}
follows from the equality
\eqRef{e*b=e*a*x \SideNS}{\Cols}
and the theorem
\refTheorem{coordinates of vector}{\SideNS-\Cols}.
The equality
\eqRef{a*x=b \SideNS}{\Cols}
follows from the equality
\eqRef{b=a*x \SideNS}{\Cols}.
}

\DefLabeledTheorem{linear span is vector space}{\SideNS-\Cols}
{
$\spanb$ is subspace of \SideWS $A$\Hyph vector space of \ColN s.
}

\DefProof{linear span is vector space}
{
Suppose
\ShowEq{linear span, vector space, 1}
According to the theorem
\refTheorem{vector in linear span}{\SideNS-\Cols}
\DrawEq[bb]{linear span, 2 \SideNS}{b \SideNS-\Cols}
\DrawEq[cc]{linear span, 2 \SideNS}{c \SideNS-\Cols}
Then
\ShowEq{linear span, 3 \SideNS}
This proves the statement.
}

\DefLabeledDefinition{module type}{\SideNS-\Cols}
{
We represented the set of vectors
\ShowEq{\RowN() 1n}vm{}
as
\RowNWS of matrix
\DrawEq[vm]{a=(a1.n \RowN)}{}
and the set of $\Base_{\BaseExt}$\Hyph nummbers
\ShowEq{\ColNS() 1n}cm{}
as
\ColWS of matrix
\DrawEq[cm]{a=(a1.n \ColNS)}{cm \SideNS-\Cols}
Corresponding representation of \SideWS $\Base$\Hyph \VectorSet $V$ is called
\AddIndex{\SideWS $\Base$\Hyph \VectorSet of \ColN s}{\SideWS A-\ColNWS space},
and $V$\Hyph number is called
\AddIndex{\ColWS vector}{\ColNWS vector}.
}

\DefLabeledTheorem{linear combination in module type}{\SideNS-\Cols}
{
We can represent linear combination
\ShowEq{linear combination, \SideNS-\Cols}cvm
of vectors
\ShowEq{\RowNWS set 1n}vm{}
as \ProductType product of matrices
\DrawEq[cvm]{w.e, \SideNS-\Cols}{cv}
}

\DefProof{linear combination in module type}
{
The theorem follows from definitions
\refDefinition{linear combination of vectors}{\SideNS},
\refDefinition{module type}{\SideNS-\Cols}.
}

\DefLabeledTheorem{vector in linear span}{\SideNS-\Cols}
{
Let $V$ be a \SideWS $A$\Hyph vector space of \ColsWS and
\RowNWS
\ShowEq{linear span, 0}
be set of vectors.
The vector $\Vector b$
linearly dependends on vectors $\ARow{\Vector a}i$
\ShowEq{linear span, 1}
if there exist \ColWS of $A$\Hyph numbers
\DrawEq[xm]{a=(a1.n set)}{}
such that the following equality is true
\DrawEq[bx]{linear span, 2 \SideNS}{\SideNS-\Cols}
}

\DefProof{vector in linear span}
{
The equality
\eqRef{linear span, 2 \SideNS}{\SideNS-\Cols}
follows from the definition
\refDefinition{linear span, vector space}{\SideNS-\Cols}
and the theorem
\refTheorem{linear combination in module type}{\SideNS-\Cols}.
}

\DefLabeledDefinition{linear span, vector space}{\SideNS-\Cols}
{
Let $V$ be a \SideWS $A$\Hyph vector space of \ColsWS and
\ShowEq{linear span, set}
be set of vectors.
\AddIndex{Linear span}{linear span, vector space}
in \SideWS $A$\Hyph vector space of \ColsWS is set
\ShowEq{linear span, vector space}
of vectors linearly dependent on vectors
\ARow {\Vector a}i.
}

\DefLabeledTheorem{vector space of algebra A}{\SideNS}
{
Let $D$ be commutative ring with unit.
Let $A$ be associative division $D$\Hyph algebra.
The representation generated by \SideWS shift
\ShowEq{representation \SideWS vector space of algebra A}
is \SideWS vector space
\ShowEq{\SideWS vector space of algebra A}
over $D$\Hyph algebra $A$.
}

\DefProof{vector space of algebra A}
{
According to definitions
\refDefinition{module over associative algebra}{\SideWS \VectorSetNS},
\RefDefinition{algebra over ring},
$D$\Hyph algebra $A$ is Abelian group.
From the equality
\ShowEq{\SideWS a(b+c)=}
and definitions
\RefDefinition{homomorphism},
\RefDefinition{endomorphism},
it follows that the map
\ShowEq{\SideWS a->ab}
is endomorphism of Abelian group $A$.
From equalities
\ShowEq{\SideWS (a1+a2)b=}
\ShowEq{\SideWS a1a2b=}
and definitions
\RefDefinition{homomorphism},
\RefDefinition{representation of algebra},
it follows that the map $f$ is
representation of $D$\Hyph algebra $A$ in Abelian group $A$.
The theorem follows from the definition
\refDefinition{module over associative algebra}{\SideWS \VectorSetNS}.
}

\AddEq{def row text}
{%
\def\ColWS{row }%
\def\ColsWS{rows }%
\def\ColTNS{rows}%
\def\RowsNS{columns}%
\def\RowsRWSA{columns }%
}

\AddEq{def col text}
{%
\def\ColWS{column }%
\def\ColsWS{columns }%
\def\ColTNS{columns}%
\def\RowsNS{rows}%
\def\RowsRWSA{rows }%
}

\DefLabeledDefinition{spectrum of matrix}{\Product}
{
The set
\ShowEq{spectrum of matrix}
of all left and right
\ProductType eigenvalues is called
\ProductType spectrum of the matrix a.
}

\DefLabeledTheorem{vector space over algebra}{\SideNS}
{
The following diagram of representations describes \SideWS $A$\Hyph vector space $V$
\DrawEq[{}{}{}{}]{diagram of representations, \SideWS module}1
The diagram of representations
\eqRef{diagram of representations, \SideWS module}1
holds
\AddIndex{commutativity of representations}{commutativity of representations}
of commutative ring $D$ and $D$\Hyph algebra $A$
in Abelian group $V$
\ShowEq{\SideWS module, a d v}
}

\DefProof{vector space over algebra}
{
The diagram of representations
\eqRef{diagram of representations, \SideWS module}1
follows from the definition
\refDefinition{module over associative algebra}{\SideWS \VectorSetNS}
and from the theorem
\RefTheorem{Free Algebra over Ring}.
Since \SideNS\HSide transformation $\ATransf(a)$
is endomorphism
of $\CBase$\Hyph vector space $V$,
we obtain the equality
\EqRef{\SideWS module, a d v}.
}

\DefProof{definition of A module, property}
{
The equality
\eqRef{commutative law}{\SideWS vector space}
follows from definitions
\refDefinition{module over algebra}{vector space},
\refDefinition{module over associative algebra}{\SideWS \VectorSetNS}.
Since \SideNS-side transformation $p$
is endomorphism of the Abelian group $V$,
we obtain the equation
\eqRef{distributive law, \SideWS module, 1}1.
Since representation
$g_{3,4}$
is homomorphism of the additive group
of $D$\Hyph algebra $A$,
we obtain the equation
\eqRef{distributive law, \SideWS module, 2}1.
Since representation
$g_{3,4}$
is \SideNS\Hyph side representation
of the multiplicative group of $D$\Hyph algebra $A$,
we obtain the equality
\EqRef{unitarity law, \SideWS A-module}.
Since representation
$g_{3,4}$
is \SideNS\Hyph side representation
of the multiplicative group of $D$\Hyph algebra $A$,
we obtain the equality
\eqRef{associative law, \SideWS module}1
when $D$\Hyph algebra $A$ is associative.
}

\DefLabeledTheorem{coordinate matrix of vector}{\SideNS-\Cols}
{
If we write vectors of basis $\Basis e$
as \RowNWS of matrix
\DrawEq[en]{a=(a1.n \RowN)}{type \SideNS-\Cols}
and coordinates of vector
\ShowEq{w=w*e \SideNS-\Cols}w{}
with respect to basis $\Basis e$ as
\ColWS of matrix
\DrawEq[wn]{a=(a1.n \ColNS)}{type \SideNS-\Cols}
then we can represent the vector
$\Vector w$
as \ProductType product of matrices
\DrawEq[w]{vector w=w*e \SideNS-\Cols}w
}

\DefProof{coordinate matrix of vector}
{
The theorem follows from the theorem
\refTheorem{linear combination in module type}{\SideNS-\Cols}.
}

\DefLabeledExample{Vector Space Type}{\SideNS-\Cols}
{
\ePrints{2022.01.05}
\ifx\Semafor\ValueOff
We represented the set of vectors
\ShowEq{\RowNWS set 1n}vm{}
as
\RowNWS of matrix
\DrawEq[vm]{a=(a1.n \RowN)}{}
and the set of $A$\Hyph nummbers
\ShowEq{\ColNWS set 1n}cm{}
as
\ColWS of matrix
\DrawEq[cm]{a=(a1.n \ColNS)}{cm \SideNS-\Cols}
Thus we can represent linear combination of vectors
\ShowEq{\RowNWS set 1n}vm{}
as \ProductType product of matrices
\DrawEq[{\ParmA}{\ParmB}m]{w.e, \SideNS-\Cols}{\ParmA\ParmB}
In particular,
if
\else
If
\fi
we write vectors of basis $\Basis e$ as
\RowNWS of matrix
\DrawEq[en]{a=(a1.n \RowN)}{\SideNS-\Cols}
and coordinates of vector
\ShowEq{w=w*e \SideNS-\Cols}w{}
with respect to basis $\Basis e$ as
\ColWS of matrix
\DrawEq[wn]{a=(a1.n \ColNS)}{\SideNS-\Cols}
then we can represent the vector
$\Vector w$
as \ProductType product of matrices
\DrawEq[w]{vector w=w*e \SideNS-\Cols}w
Corresponding representation of vector space $V$ is called
\AddIndex{\SideWS $A$\Hyph vector space of \ColN s}{\SideWS A-\ColNWS space},
and $V$\Hyph number is called
\AddIndex{\ColWS vector}{\ColNWS vector}.
}

\DefRemark{Left and Right Eigenvalues not equal 0}
{
According to the definition given above,
sets of left and right \ProductType eigenvalues coincide.
However, we considered the most simple case.
Cohn considers definitions

\ShowEq{def right}
\ShowEq{\DefCol}
\FrameEqRef[2{2}{a_1}]{A \ProductS U=U \ProductS D \SideNS}1
\ShowEq{def left}
\ShowEq{\DefRow}
\FrameEqRef[2{2}{a_1}]{A \ProductS U=U \ProductS D \SideNS}1

\noindent
where $a_2$ is non\Hyph \ProductType singular matrix
as starting point for research.
}

\AddEq{remark: Left and Right Eigenvalues not equal}
{
Let multiplicity of all \SideWS \ProductType eigenvalues
of \nTimes matrix $a_2$ equal $1$.
Let \nTimes matrix $a_2$ do not have $\gin$
\SideWS \ProductType eigenvalues,
but have
\ShowEq{gik<gin}
\SideWS \ProductType eigenvalues.
In such case, the matrix $u_2$ in the equality

\FrameEqRef[2{2}{a_1}]{A \ProductS U=U \ProductS D \SideNS}1

\noindent
is
\ShowEq{gin x gik \SideNS}
matrix and the matrix $a_1$
is \nTimes[k] matrix.
However, the equality
\eqRef{A \ProductS U=U \ProductS D \SideNS}1
does not imply that
\ShowEq{rank rc u2 =k}

Let
\ShowEq{rank rc u2 =k}
Then there exists
\ShowEq{gin x gik \OtherSideNS}
matrix
\ShowEq{a in left}w{u_2}{\OtherSideNS}{}
such that
\ShowEq{w rc u2=Ek \SideNS}
The equality
\ShowEq{w rc a2 rc u2=a1 \SideNS}
follows from equalities
\eqRef{A \ProductS U=U \ProductS D \SideNS}1,
\EqRef{w rc u2=Ek \SideNS}.

Let
\ShowEq{rank rc u2 <k}
Then 
there exists linear dependence between
\SideWS eigenvectors.
In such case,
sets of left and right \ProductType eigenvalues do not coincide.
}

\DefText[4]{define homomorphism of vector space(1)}
{
According to definitions
\ShowRef{define homomorphism A module}
the map
\DrawEq[gA]{homomorphism D algebra #1#2}{\SideWS vector}
is homomorphism
of $D_{#1}$\Hyph algebra $A_{#2}$
into $D_{#3}$\Hyph algebra $A_{#4}$.
Therefore, equalities
\ShowEq{ref homomorphism of A vector space(#1)}{#1}{#2}{#3}
follow from the theorem
\refTheorem{linear map of D module}{#1#2}.
}

\DefText[4]{define homomorphism of vector space()}
{
}

\AddEq{remark: colinear vectors}
{
Vectors $v$ and
\ShowEq{av \SideNS}
are called
\OtherSideWS
\AddIndex{colinear}{colinear vectors}.
The statement that vectors $v$ and $w$ are
\OtherSideWS colinear vectors
does not imply that
vectors $v$ and $w$ are \SideWS linearly dependent.
}

\DefText[7]{be extended basis}
{
Let the basis
\DrawEq[{#1}{#2}{#3}{#4}{#5}{#6}]{extended basis}{#7}
be extension of the basis \eV[1][.]
}

\DefLabeledTheorem{Two bases of module}{\SideNS-\Cols}
{
\ShowEq{Let be basis of algebra}{}{}kKD{}A{}-
\ShowEq{Let be basis of module}1{}iI{\SideWS}A{}V{}
Then the set of $V$\Hyph numbers
\DrawEq[2{}{{}}{}ik]{extended basis}{\SideNS-\Cols}
is the basis of $D$\Hyph \VectorSet $V$.
The basis \eV[2] is called extension of the basis \eV[1][.]
}

\DefProofRef{Two bases of module}{\SideNS-\Cols}
{
The theorem follows from the theorem
\refTheorem{coordinates of vector}{-\Cols}
and from lemmas
\refLemma{Two bases of module 1}{\SideNS-\Cols},
\refLemma{Two bases of module 2}{\SideNS-\Cols}.
}

\DefLabeledLemma{Two bases of module 1}{\SideNS-\Cols}
{
The set of $V$\Hyph numbers \EBase 2{ik}
generates $D$\Hyph \VectorSet $V$.
}

\DefLemmaProof{Two bases of module 1}
{
Let
\DrawEq{Two bases of module 1 a}{\SideNS-\Cols}
be any $V$\Hyph number.
For any $A$\Hyph number \ACol ai,
according to the definition
\RefDefinition{algebra over ring}
and the convention
\RefConvention{basis of algebra as basis of module},
there exist $D$\Hyph numbers \ACol a{ik} such that
\DrawEq{Two bases of module 1 ai}{\SideNS-\Cols}
The equality
\DrawEq{Two bases of module 1 a=ai}{\SideNS-\Cols}
follows from equalities
\eqRef{extended basis}{\SideNS-\Cols},
\eqRef{Two bases of module 1 a}{\SideNS-\Cols},
\eqRef{Two bases of module 1 ai}{\SideNS-\Cols}.
According to the theorem
\refTheorem{set of vectors generated by set of vectors}{module},
the lemma
follows from the equality
\eqRef{Two bases of module 1 a=ai}{\SideNS-\Cols}.
}

\DefLabeledLemma{Two bases of module 2}{\SideNS-\Cols}
{
The set of $V$\Hyph numbers \EBase 2{ik}
is linearly independent over the ring $D$.
}

\DefLemmaProof{Two bases of module 2}
{
Let
\DrawEq{extended basis 1}{\SideNS-\Cols}
The equality
\DrawEq{extended basis 2}{\SideNS-\Cols}
follows from equalities
\eqRef{extended basis}{\SideNS-\Cols},
\eqRef{extended basis 1}{\SideNS-\Cols}.
According to the theorem
\refTheorem{coordinates of vector}{\SideNS-\Cols},
the equality
\DrawEq{extended basis 3}{\SideNS-\Cols}
follows from the equality
\eqRef{extended basis 2}{\SideNS-\Cols}.
According to the theorem
\refTheorem{coordinates of vector}{-\Cols},
the equality
\DrawEq{extended basis 4}{\SideNS-\Cols}
follows from the equality
\eqRef{extended basis 3}{\SideNS-\Cols}.
Therefore,
the set of $V$\Hyph numbers \EBase 2{ik}
is linearly independent over the ring $D$.
}

\DefLabeledTheorem{product vector over scalar}{\SideNS}
{
\ShowEq{ab in A,v in V}
\DrawEq{(av)b=a(vb)}{\SideNS}
}

\DefProof{product vector over scalar}
{
According to the theorem
\refTheorem{similarity transformation}{\SideNS-\Cols}
and to the definition
\refDefinition{product vector over scalar}{\SideNS},
the equality
\ShowEq{(av)b= \SideNS}
follows from the equality
\ShowRef{homomorphism, f av=}
The equality
\eqRef{(av)b=a(vb)}{\SideNS}
follows from the equality
\EqRef{(av)b= \SideNS}.
}

\DefLabeledDefinition{product vector over scalar}{\SideNS}
{
\ShowEq{Let V be vector space and basis}
Bilinear map
\DrawEq{\OtherSideWS A*V->V}2
generated by \OtherSideNS\Hyph side representation
\ShowEq{twin to \SideWS product}
is called
\AddIndex{\OtherSideNS\Hyph side product}{\OtherSideNS-side product}
of vector over scalar.
}

\DefLabeledLemma{rcd dcr basis}{\SideNS}
{
We can identify basis manifold
\ShowEq{basis manifold of V, \SideNS-cols}{\aD En}{\GL nA}{}
and the set of \ProductType regular matrices \GL nA.
}

\AddEq{proof: rcd dcr basis 1}
{
\begin{sloppypar}
{\sc Proof.}
According to the remark
\refRemark{identify basis and matrix of coordinates}{\SideNS-\Cols},
we can identify a basis $\Basis e$
of \SideWS $A$\Hyph vector space $V$
and the matrix $e$ of coordinates of basis $\Basis e$
withb respect to the basis $\Basis{\aD En}$
\ShowEq{be=e cr En \SideNS}
\hfill\(\odot\)
\end{sloppypar}
}

\DefLabeledDefinition{dimension of vector space}{\SideNS}
{
We call
\AddIndex{dimension of \SideWS $A$\Hyph vector space}{dimension of vector space}
the number of vectors in a basis.
}

\DefLabeledTheorem{coordinate matrix of basis}{\Product-\Cols}
{
\ShowEq{Let V be vector space of}.
The coordinate matrix of basis $\Basis{g}$ relative basis $\Basis e$ of
\SideWS $A$\Hyph vector space $V$ is \ProductType nonsingular matrix.
}

\DefProof{coordinate matrix of basis}
{
According to the theorem
\refTheorem{rank of matrix}{\Product-\Cols},
\CR rank of the coordinate matrix of basis $\Basis{g}$ relative basis $\Basis e$
equal to the dimension of left $A$\Hyph space of columns.
This proves the statement of the theorem.
}

\DefLabeledTheorem{automorphism of vector space}{\Product-\Cols}
{
\ShowEq{Let V be vector space of}.
Let $\Basis e$ be a basis of \SideWS $A$\Hyph vector space $V$.
Then any automorphism $\Vector f$ of \SideWS $A$\Hyph vector space $V$
has form
\DrawEq[{v'}{}v{}f{}]{v1=v2*a \SideNS-\Cols}{automorphism}
where $f$ is a \ProductType nonsingular matrix.
}

\DefProof{automorphism of vector space}
{
The equality
\eqRef{v1=v2*a \SideNS-\Cols}{automorphism}
follows from theorem
\refTheorem{homomorphism A module}{\SideNS-\Cols(1)}.
Because $\Vector f$ is an isomorphism,
for each vector $\Vector v'$ there exist one and only one vector $\Vector v$ such that
\ShowEq{v'=fov}
Therefore, system of linear equations
\eqRef{v1=v2*a \SideNS-\Cols}{automorphism}
has a unique solution. According to corollary
\RefCorollary{nonsingular rows system of linear equations}
matrix $f$ is a \ProductType nonsingular matrix.
}

\DefLabeledTheorem{Automorphisms of vector space form a group}{\Product-\Cols}
{
Matrices of automorphisms of \SideWS $A$\Hyph space $V$ of \ColN s form a group \Group nA.
}

\DefProof{Automorphisms of vector space form a group}
{
If we have two automorphisms $\Vector f$ and $\Vector g$ then we can write
\ShowEq{Automorphisms of vector space, \Product-\Cols}
Therefore, the resulting automorphism has matrix $f\ProductVal g$.
}

\DefLabeledTheorem{av=bv=>a=b}{\SideNS-\Cols}
{
If \SideWS $A$\Hyph vector space $V$
has finite dimension,
then the statement
\DrawEq{av=bv, v \SideNS}{\Cols}
implies $a=b$.
}

\DefProof{av=bv=>a=b}
{
Let the set of vectors
\DrawEq[{}jJ]{basis e of module cols}-
be the basis of left $A$\Hyph vector space $V$.
According to theorems
\refTheorem{definition of A module, property}{\SideNS},
\refTheorem{coordinates of vector}{\SideNS-\Cols},
the equality
\ShowEq{avi=bvi \SideNS-\Cols}
follows from the equality
\eqRef{av=bv, v \SideNS}{\Cols}.
The theorem follows from the equality
\EqRef{avi=bvi \SideNS-\Cols},
if we assume $v=\aD e1$.
}

\DefText[7]{matrix of numbers}
{
Let
\ShowEq{matrix fIJ}{#3}{#4}{#5}{#6}{#7}
be matrix of $#1_{#2}$\Hyph numbers.
}

\DefText[9]{matrix of numbers and C}
{
Let
\ShowEq{matrix fIJ}{#3}{#4}{#5}{#6}{#7}
be matrix of $#1_{#2}$\Hyph numbers
which satisfy the equality
\ShowRef{algebra, homomorphism and product}{#8}{#9}{#3}
}

\DefText{map be homomorphism of ring (1)}
{
Let the map
\DrawEq[h{D_1}{D_2}{}]{f: A->B}{}
be homomorphism of the ring $D_1$ into the ring $D_2$.
}

\DefText{map be homomorphism of ring ()}
{
}

\DefLabeledTheorem[6]{matrix generates A module homomorphism}{\SideNS-\Cols(#1#2#3)}
{
\ShowText{map be homomorphism of ring (#1)}
\ShowText{matrices of numbers(#2#3)}{#4}{#5}{#1}{#2}
The map\refFootnote{homomorphism of A module}{\SideNS-\Cols(#1#2#3)}
\newline
\FrameEqRef[{\Vector g}{\Vector f}{}]{homomorphism A module #1#2#3}{Vector \SideNS-\Cols}
\newline
defined by the equality
\ShowText{define homomorphism A module by matrix(#1#2#3)}
is homomorphism
of \SideWS $A_{#2}$\Hyph \VectorSet of \ColsWS $V_{#3}$
into \SideWS $A_{#5}$\Hyph \VectorSet of \ColsWS $V_{#6}$.
The homomorphism
\eqRef{homomorphism A module #1#2#3}{Vector \SideNS-\Cols}
which has the given%
\ShowText{define homomorphism by given matrix(#2)}%
is unique.
}

\DefText{define homomorphism by given matrix(1)}
{
set of matrices $(g,f)$
}

\DefText{define homomorphism by given matrix()}
{
matrix $f$
}

\DefText[4]{matrix generates homomorphism (1)}
{
According to the theorem
\refTheorem{matrix generates D algebra homomorphism}{\Cols(#1#2)},
the map
\DrawEq[{\Vector g}A]{homomorphism D algebra #1#2}{(\SideNS-\Cols)algebra}
\begin{sloppypar}
\noindent
is homomorphism
of $D_{#1}$\Hyph algebra $A_{#2}$
into $D_{#3}$\Hyph algebra $A_{#4}$
and the homomorphism
\eqRef{homomorphism D algebra #1#2}{(\SideNS-\Cols)algebra}
is unique.
\end{sloppypar}

}

\DefText[4]{matrix generates homomorphism ()}
{
}

\DefProof[7]{matrix generates A module homomorphism}
{
\ShowText{matrix generates homomorphism (#2)}{#1}{#2}{#4}{#5}
The equality
\DrawEq{fo(v+w) (#2)\SideNS-\Cols}{(#1#2)}
\begin{sloppypar}
\noindent
follows from equalities
\ShowRef{matrix generates A module homomorphism(#2)}{#2}{#3}{#7}
From the equality
\eqRef{fo(v+w) (#2)\SideNS-\Cols}{(#1#2)},
it follows that the map $\Vector f$ is homomorphism
of Abelian group.
The equality
\end{sloppypar}
\DrawEq{fo(va) (#2)\SideNS-\Cols}{(#1#2)}
follows from equalities
\ShowRef{fo(va)(#2)}{#2}{#3}{#7}
From the equality
\eqRef{fo(va) (#2)\SideNS-\Cols}{(#1#2)},
and definitions
\RefDefinition{Morphism of Diagram of Representations},
\refDefinition{homomorphism A module}{\SideNS(#1#2#3)},
it follows that the map
\newline
\FrameEqRef[{\Vector g}{\Vector f}{}]{homomorphism A module #1#2#3}{Vector \SideNS-\Cols}
\newline
\begin{sloppypar}
\noindent
is homomorphism
of $A_{#2}$\Hyph \VectorSet of \ColN s $V_{#3}$
into $A_{#5}$\Hyph \VectorSet of \ColN s $V_{#6}$.
\end{sloppypar}

Let $f$ be
matrix of homomorphisms $\Vector f$, $\Vector g$
relative to bases
\ShowEq{Bases eVW}VW{}.%
The equality
\ShowEq{fov=gov \SideNS-\Cols}
follows from the theorem
\refTheorem{homomorphism A module}{\SideNS-\Cols(1)}.
Therefore, $\Vector f=\Vector g$.
}

\DefLabeledTheorem[6]{homomorphism A module}{\SideNS-\Cols(#1#2#3)}
{
The homomorphism\refFootnote{homomorphism of A module}{\SideNS-\Cols(#1#2#3)}
\DrawEq[{\Vector g}{\Vector f}{}]{homomorphism A module #1#2#3}{Vector \SideNS-\Cols}
of \SideWS $A_{#2}$\Hyph \VectorSet of \ColsWS $V_{#3}$
into \SideWS $A_{#5}$\Hyph \VectorSet of \ColsWS $V_{#6}$
has presentation
\ShowText{g:A1->A2, D module(#1#2#3)}
\ShowText{f o (ae)=a o f e (#2#3)}{#1}{#2}{#3}
relative to selected bases.
Here
\begin{itemize}
\ShowText{homomorphism of vector space, algebra(#2)}{#1}{#2}{#3}{#4}{#5}{#6}
\ShowText{homomorphism of vector space, algebra 1}{#1}{#2}{#3}{#4}{#5}{#6}
\end{itemize}
\ShowText{matrix of homomorphism relative bases #2#3}{#1}{#2}{#3}{#5}
}

\DefText[4]{matrix of homomorphism relative bases 11}
{
The set of matrices $(g,f)$ is unique and is called
{\bf coordinates of homomorphism}
\eqRef{homomorphism A module #1#2#3}{Vector \SideNS-\Cols}
relative bases
\ShowEq{bases eA eV}{#2}1,
\ShowEq{bases eA eV}{#4}2.
}

\DefText[4]{matrix of homomorphism relative bases 1}
{
The matrix $f$ is unique and is called
{\bf matrix of homomorphism}
$\Vector f$ relative bases \eV[1][,] \eV[2][.]
}

\DefText[6]{Proof, homomorphism of vector space(1)}
{
According to definitions
\ShowRef{Proof, homomorphism of vector space}{#1}{#2}{#3}%
the map
\DrawEq[{\Vector g}A]{homomorphism D algebra #1#2}{\SideNS-\Cols}
\begin{sloppypar}
\noindent
is homomorphism
of $D_{#1}$\Hyph algebra $A_{#2}$
into $D_{#4}$\Hyph algebra $A_{#5}$.
Therefore, equalities
\ShowRef{Proof, homomorphism of vector space 1}{#1}{#2}{#3}%
follow from equalities
\ShowRef{Proof, homomorphism of vector space 2}{#1}{#2}{#3}{#4}%
From the theorem
\refTheorem{linear map of D module, coordinates}{#1#2\Cols},
it follows that the matrix $g$ is unique.
\end{sloppypar}

}

\DefText[6]{Proof, homomorphism of vector space()}
{%
}

\DefProof[6]{homomorphism A module}
{
\ShowText{Proof, homomorphism of vector space(#2)}{#1}{#2}{#3}{#4}{#5}{#6}%
Vector
\ShowEq{vv in V}v{V_1}
has expansion
\DrawEq[v{V_1}]{\SideWS homomorphism \Product, 1}{v#1#2#3}
relative to the basis \eV[V_1][.]
Vector
\ShowEq{vv in V}w{V_2}
has expansion
\DrawEq[w{V_2}]{\SideWS homomorphism \Product, 1}{w#1#2#3}
relative to the basis \eV[V_2][.]
Since $\Vector f$ is a homomorphism,
then the equality
\DrawEq{vv=v*eV \SideWS \Product(#2)}{(#1)}
follows from equalities
\eqRef{homomorphism, f v+w=}{f(#1#2)\SideWS A module},
\eqRef{\SideWS homomorphism, f av=}{(#1#2)\SideWS A module}.
$V_2$\Hyph number
\ShowEq{Vector f(e) \Cols}
has expansion
\DrawEq{f o ei=fij ej \SideWS \Product}{#1#2#3}
relative to basis \eV[V_2][.]
The equality
\DrawEq{vb=a cr f cr e \SideWS \Product(#2)}{(#1)}
follows from equalities
\eqRef{vv=v*eV \SideWS \Product(#2)}{(#1)},
\eqRef{f o ei=fij ej \SideWS \Product}{#1#2#3}.
The equality
\eqRef{v1=v2*a \SideNS-\Cols}{\SideNS(#1#2#3)}
follows from comparison of
\eqRef{\SideWS homomorphism \Product, 1}{w#1#2#3}
and \eqRef{vb=a cr f cr e \SideWS \Product(#2)}{(#1)} and
the theorem
\refTheorem{coordinates of vector}{\SideNS-\Cols}.
From the equality
\eqRef{f o ei=fij ej \SideWS \Product}{#1#2#3}
and from the theorem
\refTheorem{coordinates of vector}{\SideNS-\Cols},
it follows that the matrix $f$ is unique.
}

\DefText{notation for homomorphism of module}
{
We will use notation
\DrawEq{f circ a =}{}
for image of homomorphism $f$.
}

\DefText{be division algebra}
{
Let $D$\Hyph algebra $A$ be division algebra.
}

\DefTheorem{division algebra has unit}
{
\ShowText{be division algebra}
$D$\Hyph algebra $A$ has unit.
}

\DefTheorem{division D algebra, ring D is the field}
{
\ShowText{be division algebra}
The ring $D$ is the field
and subset of the center of $D$\Hyph algebra $A$.
}

\DefText[3]{homomorphism of vector space 2}
{
On the basis of theorems
\refTheorem{homomorphism A module}{\SideNS-\Cols(#1#2#3)},
\refTheorem{matrix generates A module homomorphism}{\SideNS-\Cols(#1#2#3)}
we identify the homomorphism
\newline
\FrameEqRef[{\Vector g}{\Vector f}{}]{homomorphism A module #1#2#3}{Vector \SideNS-\Cols}
\newline
of \SideWS $A$\Hyph vector space $V$ of \ColN s
and coordinates of its presentation
\ShowText{define homomorphism A module by matrix(#1#2#3)}
}

\DefLabeledDefinition{similarity transformation}{\SideNS}
{
The endomorphism
\ShowEq{\SideWS map a En}a{}
of \SideWS $A$\Hyph vector space $V$
is called
\AddIndex{similarity transformation}{similarity transformation}
with respect to the basis $\Basis e$.
}

\AddEq{passive transformation for two bases}
{
\ShowEq{Let V be vector space of and}
\ShowEq{basis e1 e2}
be bases of \SideWS $A$\Hyph vector space $V$.
Let $g$
be passive transformation
of basis $\Basis e_1$ into basis $\Basis e_2$
\DrawEq[12g{}]{e2i=aij e1j \Product-\Cols}{}
}

\DefLabeledTheorem{covariance of eigenvalue}{\SideNS-\Cols}
{
\ShowEq{Let V be vector space of and}
\ShowEq{basis e1 e2}
be bases of \SideWS $A$\Hyph vector space $V$.

\begin{Statement}
\labelStatement{covariance of eigenvalue \SideNS-\Cols}
Eigenvalue $b$
of the endomorphism $\Vector f$
with respect to the basis $\Basis e_1$
does not depend on the choice of the basis $\Basis e_2$.
\hfill\(\odot\)
\end{Statement}

\begin{Statement}
\labelStatement{covariance of eigenvector \SideNS-\Cols}
Eigenvector $\Vector v$
of the endomorphism $\Vector f$
corresponding to eigenvalue $b$
does not depend on the choice of the basis $\Basis e_2$.
\hfill\(\odot\)
\end{Statement}
}

\DefProof{covariance of eigenvalue}
{
Let
\ShowEq{Ai, i=}f
be the matrix of endomorphism $f$ with respect to the basis $\Basis e_i$.
Let $b$ be eigenvalue
of the endomorphism $\Vector f$
with respect to the basis $\Basis e_1$
and
$\Vector v$ be eigenvector
of the endomorphism $\Vector f$
corresponding to eigenvalue $b$.
Let
\ShowEq{Ai, i=}v
be the matrix of vector $\Vector v$ with respect to the basis $\Basis e_i$.

Suppose first that
\ShowEq{basis e1=e2}
Then $\aD En$
is passive transformation
of basis $\Basis e_1$ into basis $\Basis e_2$
\DrawEq[12{\aD En^{}{}}{}]{e2i=aij e1j \Product-\Cols}{}
According to the theorem
\refTheorem{similarity transformation, change of basis}{\SideNS-\Cols},
the similarity transformation
\ShowEq{\SideWS map a En}b1
has the matrix
\ShowEq{ga*g- \SideNS-En}
with respect to the basis $\Basis e_2$.
Therefore, according to the theorem
\refTheorem{eigenvalue of endomorphism, matrix is singular}{\SideNS-\Cols},
we proved following statements.

\begin{ShadedLemma}
\labelLemma{matrix f-bEn singular \SideNS-\Cols}
The matrix
\ShowEq{f-bEn \SideNS}1{}
is \ProductType singular matrix.
\hfill\(\odot\)
\end{ShadedLemma}

\begin{ShadedLemma}
\labelLemma{vector v1 f-bEn \SideNS-\Cols}
The \ColWS vector $v_1$ satisfies to the system of linear equations
\ShowEq{(f-bEn)v \Product-\Cols}
\hfill\(\odot\)
\end{ShadedLemma}

Suppose that
\ShowEq{basis e1 ne e2}
According to the theorem
\refTheorem{exists representation, commuting with active}{\SideNS-\Cols},
there exists unique passive transformation
\DrawEq[12g{}]{e2i=aij e1j \Product-\Cols}{}
and $g$ is \ProductType non\Hyph singular matrix.
According to the theorem
\refTheorem{passive transformation and endomorphism}{\SideNS-\Cols},
\DrawEq[fg{\SideNS}{\Cols}]{f2=a f1 a-}{covariance \SideNS-\Cols}
\begin{sloppypar}
\noindent
According to the theorem
\refTheorem{eigenvalue of endomorphism, matrix is singular}{\SideNS-\Cols},
if $b$ is eigenvalue of endomorphism $\Vector f$,
then we have to prove the lemma
\RefLemma{matrix f-ae singular \SideNS-\Cols}.
\end{sloppypar}

\begin{ShadedLemma}
\labelLemma{matrix f-ae singular \SideNS-\Cols}
The matrix
\[\ShowEq{(f-ae)v matrix}f2bg\]
is \ProductType singular matrix.
\hfill\(\odot\)
\end{ShadedLemma}

\begin{sloppypar}
Since the equality
\ShowEq{f-ae->f-be}
follows from the equality
\eqRef{f2=a f1 a-}{covariance \SideNS-\Cols},
then the lemma
\RefLemma{matrix f-ae singular \SideNS-\Cols}
follows from the lemma
\RefLemma{matrix f-bEn singular \SideNS-\Cols}.
Therefore, we proved the statement
\RefStatement{covariance of eigenvalue \SideNS-\Cols}.
\end{sloppypar}

According to the theorem
\refTheorem{eigenvector coordinates}{\SideNS-\Cols},
if $\Vector v$ is eigenvector of endomorphism $\Vector f$,
then we have to prove the lemma
\RefLemma{vector v2 f-bEn \SideNS-\Cols}.

\begin{ShadedLemma}
\labelLemma{vector v2 f-bEn \SideNS-\Cols}
The \ColWS vector $v_2$ satisfies to the system of linear equations
\ShowEq{(f-bEn)v2 \Product-\Cols}
\end{ShadedLemma}

\begin{sloppypar}
{\sc Proof.}
According to the theorem
\refTheorem{passive transformation of vector space}{\SideNS-\Cols},
\DrawEq[v2v1g{}{\SideNS}{\Cols}]{v2=v1g-}{covariance \SideNS-\Cols}
Since the equality
\ShowEq{f-ae->f-be \SideNS-\Cols}
follows from equalities
\eqRef{f2=a f1 a-}{covariance \SideNS-\Cols},
\eqRef{v2=v1g-}{covariance \SideNS-\Cols},
then the lemma
\RefLemma{vector v2 f-bEn \SideNS-\Cols}
follows from the lemma
\RefLemma{vector v1 f-bEn \SideNS-\Cols}.
\hfill\(\odot\)
\end{sloppypar}

Therefore, we proved the statement
\RefStatement{covariance of eigenvector \SideNS-\Cols}.
}

\DefLabeledTheorem{eigenvalue of endomorphism, matrix is singular}{\SideNS-\Cols}
{
\ShowEq{passive transformation for two bases}
Let $f$ be the matrix of endomorphism $\Vector f$
with respect to the basis $\Basis e_2$.
$A$\Hyph number $b$ is eigenvalue
of endomorphism $\Vector f$ iff the matrix
\DrawEq[f{}bg]{(f-ae)v matrix}{\SideNS-\Cols}
is \ProductType singular.
}

\DefProof{eigenvalue of endomorphism, matrix is singular}
{
The theorem folows from theorems
\RefTheorem{star rows system of linear equations},
\refTheorem{eigenvector coordinates}{\SideNS-\Cols}.
}

\DefLabeledTheorem{eigenvector coordinates}{\SideNS-\Cols}
{
\ShowEq{passive transformation for two bases}
Let $f$ be the matrix of endomorphism $\Vector f$
and $v$ be the matrix of vector $\Vector v$
with respect to the basis $\Basis e_2$.
The vector $\Vector v$ is
eigenvector
of the endomorphism $\Vector f$
with respect to the basis $\Basis e_1$
iff
there exists $A$\Hyph number $b$ such that
the system of linear equations
\DrawEq[fgbv]{(f-ae)v \Product-\Cols}{fgbv}
has non\Hyph trivial solution.
}

\DefProof{eigenvector coordinates}
{
According to theorems
\refTheorem{homomorphism A module}{\SideNS-\Cols(1)},
\refTheorem{similarity transformation, change of basis}{\SideNS-\Cols},
the equality
\ShowEq{e*f*v=e*ae*v \SideNS-\Cols}
follows from the equality
\eqRef{fov=b e o v}{\SideNS}.
The equality
\eqRef{(f-ae)v \Product-\Cols}{fgbv}
follows from the equality
\EqRef{e*f*v=e*ae*v \SideNS-\Cols}
and from the theorem
\refTheorem{coordinates of vector}{\SideNS-\Cols}.
Since the case $\Vector v=0$ is not interesting for us,
the system of linear equations
\eqRef{(f-ae)v \Product-\Cols}{fgbv}
should have non\Hyph trivial solution.
}

\AddEq{remark: eigenvector coordinates}
{
The definition
\refDefinition{Eigenvalue of Endomorphism}{\SideNS}
does not depend on structure of coordinates
of \SideWS $A$\Hyph vector space $V$.
As soon as we choose the basis $\Basis e$
of \SideWS $A$\Hyph vector space $V$,
we see additional structure
which is reflected in theorems
\refTheorem{eigenvector coordinates}{\SideNS-cols},
\refTheorem{eigenvector coordinates}{\SideNS-rows}.
}

\AddEq{remark: eigenvector of matrix}
{
The same matrix may correspond to different endomorphisms
depending on whether we consider left or right
$A$\Hyph vector space.
Therefore, we must consider all definitions
of eigenvectors which may correspond to given matrix.
From theorems
\refTheorem{eigenvector coordinates}{left-cols},
\refTheorem{eigenvector coordinates}{right-cols},
\refTheorem{eigenvector coordinates}{left-rows},
\refTheorem{eigenvector coordinates}{right-rows},
it follows that for given matrix there exist
four types of eigenvectors and
four types of eigenvalues.
Couple eigenvector and eigenvalue of a certain type
is used in solving the problem
where corresponding type of $A$\Hyph vector space arises.
This statement is reflected in the names
of eigenvector and eigenvalue.

Eigenvector of endomorphism depends on
choice of two bases. However when we
consider eigenvector of matrix,
it is not the choice of bases that is important for us,
but the matrix of passive transformation between these bases.
}

\DefLabeledDefinition{eigenvalues of pair of matrices}{\Product}
{
$A$\Hyph number $b$ is called
{\bf \ProductType eigenvalue}
of the pair of matrices $(f,g)$
where $g$ is the \ProductType non\Hyph singular matrix
if the matrix
\DrawEq[f{}bg]{(f-ae)v matrix}{}
is \ProductType singular matrix.
}

\DefLabeledDefinition{eigenvector of pair of matrices}{\Product-\ColN}
{
Let $A$\Hyph number $b$ be
\ProductType eigenvalue of the pair of matrices
\ShowEq{pair of matrices \Product-\ColN}{}
where $g$ is the \ProductType non\Hyph singular matrix.
A \ColWS $v$ is called
{\bf eigen\ColWS}
of the pair of matrices $(f,g)$ corresponding to \ProductType eigenvalue $b$,
if the following equality is true
\DrawEq{\Product-eigen\Cols(f,g)}{}
}

\AddEq{proof: eigencolumn of matrix}
{
The definition
\refDefinition{eigenvector of pair of matrices}{\Product-\ColN}
follows from the definition
\refDefinition{eigenvalues of pair of matrices}{\Product}
and from the theorem
\refTheorem{eigenvector coordinates}{\SideNS-\Cols}.
}

\AddEq{proof: eigenrow of matrix}
{
Since
\ShowEq{bg.\Product.g-}
then the definition
\refDefinition{eigenvector of pair of matrices}{\Product-\ColN}
follows from the definition
\refDefinition{eigenvalues of pair of matrices}{\Product}
and from the theorem
\refTheorem{eigenvector coordinates}{\SideNS-\Cols}.
}

\AddEq{theorem: representation Ao2->V}
{
\begin{ShadedTheorem}
\labelTheorem{representation Ao2->V \SideNS}
Let $f$, $g(\Basis e)$
be twin representations of $D$\Hyph algebra $A$
in $D$\Hyph module $V$
and $\Basis e$ be the basis of \SideWS $A$\Hyph vector space $V$.
Let product in $D$\Hyph module
\AoxA A
be defined according to rule
\DrawEq[pq]{product in algebra AA}{}
Then representation
\ShowEq{representation Ao2->V}
of $D$\Hyph algebra \AoxA A in $D$\Hyph module $V$
defined by the equality
\ShowEq{representation Ao2->V =}
is left \AoxA A\Hyph module.
\end{ShadedTheorem}
}

\AddEq{theorem: twin representation in vector space}
{
\begin{ShadedTheorem}
\labelTheorem{twin representation in \SideWS vector space}
\ShowEq{Let V be vector space}{}
\ShowEq{f:A->*B}fAV
An effective \OtherSideNS\Hyph side representation
\ShowEq{f:A->*B}{g(\Basis e)}AV
of $D$\Hyph algebra $A$
in \SideWS $A$\Hyph vector space $V$
depends on the basis $\Basis e$ of \SideWS $A$\Hyph vector space $V$.
The folowing diagram
\newline
\DrawEq{twin representations of algebra}{\SideNS}
is commutative for any $a$, $b\in A$.
\eqRef{(av)b=a(vb)}{\SideNS}.
\end{ShadedTheorem}
}

\DefProof{twin representation in vector space}
{
Commutativity of the diagram
\eqRef{twin representations of algebra}{\SideNS}
follows from the equality
}

\DefLabeledDefinition{twin representation in vector space}{\SideNS}
{
We call representations $f$ and $g(\Basis e)$
considered in the theorem
\RefTheorem{twin representation in \SideWS vector space}
{\bf twin representations}
of the $D$\Hyph algebra $A$.
The equality
\eqRef{(av)b=a(vb)}{\SideNS}
represents
{\bf associative law}
for twin representations.
This allows us writing of such expression without using of brackets.
}

\DefLabeledTheorem{set of similarity transformations}{\SideNS}
{
Let $V$ be \SideWS $A$\Hyph vector space of columns
and $\Basis e$ be basis of \SideWS $A$\Hyph vector space $V$.
The set of similarity transformations
\ShowEq{\SideWS map a En}a{}
of \SideWS $A$\Hyph vector space $V$
is \OtherSideNS\Hyph side effective representation
of $D$\Hyph algebra $A$
in \SideWS $A$\Hyph vector space $V$.
}

\DefProof{set of similarity transformations}
{
Consider the set of matrices
\ShowEq{\SideWS a En}
of similarity transformations
\ShowEq{\SideWS map a En}a{}
whith respecpect to the basis $\Basis e$.
Let $V^*$ be the set of coordinates of $V$\Hyph numbers
whith respecpect to the basis $\Basis e$.
According to theorems
\RefTheorem{Geometric objects form A vector space, \SideNS-cols},
\RefTheorem{Geometric objects form A vector space, \SideNS-rows},
the set $V^*$ is \SideWS $A$\Hyph vector space.
According to the theorem
\refTheorem{similarity transformation}{\SideNS-\Cols},
the map
\ShowEq{\SideWS map a En}a{}
is the endomorphism
of \SideWS $A$\Hyph vector space of columns.

\ColRowLemma{map a in A->aE in hom is homomorphism}

\ColRowProof{map a in A->aE in hom is homomorphism}

\ColRowProof{set of similarity transformations 1}
}

\AddEq{proof: set of similarity transformations 1}
{
According to the definition
\refDefinition{side representation of group}{\SideNS}
and the lemma
\refLemma{map a in A->aE in hom is homomorphism}{\SideNS-\Cols},
the map
\eqRef{a in A->aE hom \SideNS}{\Cols}
is \OtherSideNS\Hyph side representation
of $D$\Hyph algebra $A$
in \SideWS $A$\Hyph vector space $V^*$.
According to the theorem
\RefTheorem{division algebra},
the equality
\ShowEq{aEn=bEn \SideNS}
implies $a=b$.
According to the definition
\RefDefinition{effective representation of algebra},
the representation
\eqRef{a in A->aE hom \SideNS}{\Cols}
is effective.

According to theorems
\refTheorem{passive transformation and sum of endomorphisms}{\SideNS-\Cols},
\RefTheorem{passive transformation and product of endomorphisms, vector space, \SideNS-\Cols},
the representation
\eqRef{a in A->aE hom \SideNS}{\Cols}
is covariant with respect to choice of basis 
of \SideWS $A$\Hyph vector space $V$.
According to the theorem
\refTheorem{similarity transformation}{\SideNS-\Cols},
the set of maps
\ShowEq{\SideWS map a En}a{}
is \OtherSideNS\Hyph side representation
of $D$\Hyph algebra $A$
in \SideWS $A$\Hyph vector space $V$.
}

\DefLabeledLemma{map a in A->aE in hom is homomorphism}{\SideNS-\Cols}
{
\ShowEq{Let V be vector space of}.
The map
\DrawEq{a in A->aE hom \SideNS}{\Cols}
is homomorphism of $D$\Hyph algebra $A$.
}

\AddEq{proof: map a in A->aE in hom is homomorphism}
{
{\sc Proof.}
According to the theorem
\refTheorem{homomorphism from A1 to A2, D algebra}1,
the lemma follows from equalities
\ShowEq{(a+b)aE \SideNS}
\ShowEq{(daE)v \SideNS}
\ShowEq{(ab)E \SideNS-\Cols}
\hfill\(\odot\)
}

\DefText{passive transformation and endomorphism}
{
We consider transformation of coordinates of endomorphism
of \SideWS $A$\Hyph vector space of \ColsWS
in the theorem
\refTheorem{passive transformation and endomorphism}{\SideNS-\Cols}.
}

\DefLabeledTheorem{passive transformation and endomorphism}{\SideNS-\Cols}
{
\ShowEq{Let V be vector space of and}
$\Basis e_1$, $\Basis e_2$ be bases in \SideWS $A$\Hyph vector space $V$.
Let $g$ be passive transformation of basis $\Basis e_1$ into basis $\Basis e_2$
\DrawEq[12g{}]{e2i=aij e1j \Product-\Cols}{}
Let $f$ be endomorphism of \SideWS $A$\Hyph vector space $V$.
Let
\ShowEq{Ai, i=}f
be the matrix of endomorphism $f$ with respect to the basis $\Basis e_i$.
Then
\DrawEq[fg{\SideNS}{\Cols}]{f2=a f1 a-}{g \SideNS-\Cols}
}

\AddEq[1]{remark: passive transformation and endomorphism}
{
Let
\ShowEq{Ai, i=}{#1}
be the matrix of vector $#1$ with respect to the basis $\Basis e_i$.
Then
\DrawEq[{#1}1{#1}2g{}]{v1=v2*a \SideNS-\Cols}{#1}
follows from the theorem
\refTheorem{passive transformation of vector space}{\SideNS-\Cols}.
}

\AddEq{proof: passive transformation and endomorphism}
{
{\sc Proof of the Theorem
\refTheorem{passive transformation and endomorphism}{\SideNS-\Cols}.}
Let endomorphism $f$ maps a vector $v$ into a vector $w$
\DrawEq{w=f o v}{\SideNS-\Cols}
\ShowEq{passive transformation and endomorphism}

According to the theorem
\refTheorem{homomorphism A module}{\SideNS-cols(1)},
equalities
\DrawEq[w1v1f1]{v1=v2*a \SideNS-\Cols}{v1}
\DrawEq[w2v2f2]{v1=v2*a \SideNS-\Cols}{v2}
follow from the equality
\eqRef{w=f o v}{\SideNS-\Cols}.
The equality
\DrawEq[g]{w2=...v2 a f1 \SideNS-\Cols}g
follows from equalities
\eqRef{v1=v2*a \SideNS-\Cols}v,
\eqRef{v1=v2*a \SideNS-\Cols}{v1},
\eqRef{v1=v2*a \SideNS-\Cols}w.
Since $g$ is regular matrix,
then the equality
\DrawEq[g]{w2=...v2 a f1 a- \SideNS-\Cols}g
follows from the equality
\eqRef{w2=...v2 a f1 \SideNS-\Cols}g.
The equality
\DrawEq[g]{v2 f2=v2 a f1 a- \SideNS-\Cols}g
follows from equalities
\eqRef{v1=v2*a \SideNS-\Cols}{v2},
\eqRef{w2=...v2 a f1 a- \SideNS-\Cols}g.
The equality
\eqRef{f2=a f1 a-}{g \SideNS-\Cols}
follows from the equality
\eqRef{v2 f2=v2 a f1 a- \SideNS-\Cols}g.
\qed
}

\DefText[2]{coordinates of vector with respect to basis (1,2) cols}
{
Let
\DrawEq[{#1_i}n]{a=(a1.n col)}{#2,\SideNS-\Cols}
be matrix of coordinates of the vector $\Vector{#1}$
with respect to the basis
\ShowEq{Basis e i=1#2}
}

\DefText[2]{coordinates of vector with respect to basis (1,2) rows}
{
Let
\ShowEq{phantom a**b}
\DrawEq[{#1_i^{}{}}n]{a=(a1.n row)}{}
\ShowEq{phantom a**b}
be matrix of coordinates of the vector $\Vector{#1}$
with respect to the basis
\ShowEq{Basis e i=1#2}
}

\AddEq[2]{convention: passive transformation maps the basis 1->2}
{%
Let passive transformation $g\in #1$ map the basis \eV[1]
into the basis \eV[2]
\DrawEq[12g{}]{e2i=aij e1j \Product-\Cols}{12 #2}
}

\DefLabeledTheorem{passive transformation of vector space}{\SideNS-\Cols}
{
\ShowEq{Let V be vector space of}.
\ShowConvention{passive transformation maps the basis 1->2}{G(V_*)}1
\ShowText{coordinates of vector with respect to basis (1,2) \Cols}v2
Coordinate transformations
\DrawEq[v1v2g{}]{v1=v2*a \SideNS-\Cols}{v21g}
\DrawEq[v2v1g{}{\SideNS}{\Cols}]{v2=v1g-}{12 \SideNS-\Cols}
do not depend on vector $\Vector v$ or basis $\Basis e$, but is
defined only by coordinates of vector $\Vector v$
relative to basis $\Basis e$.
}

\DefProof{passive transformation of vector space}
{
According to the theorem
\refTheorem{coordinate matrix of vector}{\SideNS-\Cols},
\ShowEq{v=vi ei 12, \SideNS-\Cols                                                                }
The equality
\ShowEq{v1 cr e1=v2 cr g cr e1, \SideNS-\Cols}
follows from equalities
\eqRef{e2i=aij e1j \Product-\Cols}{12 1},
\EqRef{v=vi ei 12, \SideNS-\Cols}.
According to the theorem
\refTheorem{coordinates of vector}{\SideNS-\Cols},
the coordinate transformation
\eqRef{v1=v2*a \SideNS-\Cols}{v21g}
follows from the equality
\EqRef{v1 cr e1=v2 cr g cr e1, \SideNS-\Cols}.
Because $g$ is \ProductType nonsingular matrix,
coordinate transformation
\eqRef{v2=v1g-}{12 \SideNS-\Cols}
follows from the equality
\eqRef{v1=v2*a \SideNS-\Cols}{v21g}.
}

\DefLabeledTheorem{passive transformation and sum of endomorphisms}{\SideNS-\Cols}
{
\ShowEq{Let V be vector space of}.
Let endomorphism $\Vector f$ of \SideWS $A$\Hyph vector space $V$ is sum of endomorphisms
$\Vector f_1$, $\Vector f_2$.
The matrix $f$ of endomorphism $\Vector f$ equal to the sum of matrices $f_1$, $f_2$ of endomorphisms
$\Vector f_1$, $\Vector f_2$
and this equality does not depend on the choice of a basis.
}

\AddEq{two endomorphisms and two bases}
{
Let $\Basis e$, $\Basis e'$ be bases of \SideWS $A$\Hyph vector space $V$
and $g$ be passive transformation of basis $\Basis e$ into basis $\Basis e'$
\ShowEq{\SideWS e'=e*a}{}
Let $f$
be the matrix of endomorphism $\Vector f$ with respect to the basis $\Basis e$.
Let $f'$
be the matrix of endomorphism $\Vector f$ with respect to the basis $\Basis e'$.
Let
\ShowEq{Ai, i=}f
be the matrix of endomorphism $\Vector f_i$ with respect to the basis $\Basis e$.
Let
\ShowEq{Ai, i=}{f'}
be the matrix of endomorphism $\Vector f_i$ with respect to the basis $\Basis e'$.
}

\DefProof{passive transformation and sum of endomorphisms}
{
\ShowEq{two endomorphisms and two bases}
According to the theorem
\refTheorem{sum of homomorphisms, A vector space}{\SideNS-\Cols},
\DrawEq[{}{}]{f=f1+f2 endo}{\SideNS-\Cols}
\DrawEq[{'}{}]{f=f1+f2 endo}{\SideNS-\Cols'}
According to the theorem
\refTheorem{passive transformation and endomorphism}{\SideNS-\Cols},
\DrawEq[{}{}]{f'=g f g- \SideNS-\Cols}f
\DrawEq[1]{f'=g f g- \SideNS-\Cols}1
\DrawEq[2]{f'=g f g- \SideNS-\Cols}2
The equality
\ShowEq{f1+f2=g.g- \SideNS-\Cols}
follows from equalities
\ShowEq{f1+f2=g.g-}
From equalities
\eqRef{f=f1+f2 endo}{\SideNS-\Cols'},
\EqRef{f1+f2=g.g- \SideNS-\Cols}
it follows that the equality
\eqRef{f=f1+f2 endo}{\SideNS-\Cols}
is covariant with respect to choice of basis 
of \SideWS $A$\Hyph vector space $V$.
}

\AddEq{theorem: passive transformation and product of endomorphisms}
{
\begin{ShadedTheorem}
\labelTheorem{passive transformation and product of endomorphisms, vector space, \SideNS-\Cols}
\ShowEq{Let V be vector space of}.
Let endomorphism $\Vector f$ of \SideWS $A$\Hyph vector space $V$
is product of endomorphisms
$\Vector f_1$, $\Vector f_2$.
The matrix $f$ of endomorphism $\Vector f$ equal
to the \ProductType product of matrices
\ShowEq{f1 f2 \SideNS}
of endomorphisms
\ShowEq{vf1 vf2 \SideNS}
and this equality does not depend on the choice of a basis.
\end{ShadedTheorem}
}

\DefLabeledDefinition{coordinates of geometric object, vector space}{\SideNS-\Cols}
{
Orbit
\ShowEq{def coordinates of geometric object, \SideNS-\Cols}
of representation $F_1$
is called
{\bf coordinate manifold of geometric object}
in \SideWS $A$\Hyph vector space $V$ of \ColN s.
For any basis
\newline
\FrameEqRef[{V1}{V2}g{}]{e2i=aij e1j \Product-\Cols}{GOeV}
\newline
corresponding point
\newline
\FrameEqRef{coordinate transformation, vector space W, \SideNS-\Cols}1
\newline
of orbit defines
{\bf coordinates of geometric object}
in coordinate \SideWS $A$\Hyph vector space
relative basis \eV[V2].
}

\DefLabeledDefinition{Left and Right Eigenvalues}{\SideNS-\Cols}
{
Let $a_2$ be \nTimes matrix
which is \ProductType similar to diagonal matrix $a_1$
\ShowEq{D diagonal matrix}{a_1}n
Thus, there exist non\Hyph \ProductType singular matrix $u_2$ such that
\DrawEq{U-*A*U=D \Product-\SideNS}1
\DrawEq[2{2}{a_1}]{A \ProductS U=U \ProductS D \SideNS}1
The \ColWS
\ShowEq{ARow u2i \Cols}u2i
of the matrix $u_2$ satisfies to the equality
\DrawEq[2u]{A*Ui=Ui di \Product-\SideNS}u
The $A$\Hyph number $b(\gi i)$ is called \SideWS \ProductType {\bf eigenvalue}
and \ColWS vector
\ShowEq{ARow u2i \Cols}u2i
is called
\AddIndex{eigen\ColWS}{eigen\ColN}
for \SideWS \ProductType eigenvalue $b(\gi i)$.
}

\DefLabeledTheorem{Eigenvector does not depend on basis}{\SideNS-\Cols}
{
Eigenvector
\ShowEq{ARow u2i \Cols}u2i
of the matrix $a_2$
for \SideWS \ProductType eigenvalue $b(\gi i)$
does not depend on the choice of the basis $\Basis e_2$
\ShowEq{e1i=u2i*e2 \SideNS-\Cols}
}

\AddEq{remark: Eigenvector does not depend on basis}
{
The equality
\EqRef{e1i=u2i*e2 \SideNS-\Cols}
follows from the equality
\eqRef{e2i=aij e1j \Product-\Cols}{21u2}.
However, these equalities have different meaning.
The equality
\eqRef{e2i=aij e1j \Product-\Cols}{21u2}
is representation of passive transformation.
In the equality
\EqRef{e1i=u2i*e2 \SideNS-\Cols},
we see
the expansion of the eigenvector
for given \SideWS \ProductType eigenvalue with respect to selected basis.
The goal of the proof of the theorem
\refTheorem{Eigenvector does not depend on basis}{\SideNS-\Cols}
is to explain the meaning of the equality
\EqRef{e1i=u2i*e2 \SideNS-\Cols}.
}

\DefProof{Eigenvector does not depend on basis}
{
According to the theorem
\refTheorem{matrix generates A module homomorphism}{\SideNS-\Cols(1)},
we can consider the matrix $a_2$ as matrix of automorphism $a$
of \SideWS $A$\Hyph vector space with respect to the basis $\Basis e_2$.
According to the theorem
\refTheorem{passive transformation and endomorphism}{\SideNS-\Cols},
the matrix $u_2$ is the matrix is the matrix of passive transformation
\DrawEq[21u2]{e2i=aij e1j \Product-\Cols}{21u2}
and the matrix $a_1$ is the matrix of automorphism $a$
of \SideWS $A$\Hyph vector space with respect to the basis $\Basis e_1$.

The coordinate matrix of the basis $\Basis e_1$ with respect to the basis $\Basis e_1$
is identity matrix $\aD En$.
Therefore, the vector
\ShowEq{ARow u2i \Cols}e1i
is
eigen\ColWS for \SideWS \ProductType eigenvalue $b(\gi i)$
\DrawEq[1e]{A*Ui=Ui di \Product-\SideNS}e
The equality
\ShowEq{U-*A*U*E=E*D \Product-\SideNS}
follows from equalities
\eqRef{U-*A*U=D \Product-\SideNS}1,
\eqRef{A*Ui=Ui di \Product-\SideNS}e.
The equality
\ShowEq{A*U*E=U*E*D \Product-\SideNS}
follows from the equality
\EqRef{U-*A*U*E=E*D \Product-\SideNS}.
The equality
\ShowEq{U*Ei=Ui \Product-\SideNS}
follows from equalities
\eqRef{A*Ui=Ui di \Product-\SideNS}u,
\EqRef{A*U*E=U*E*D \Product-\SideNS}.
The theorem follows from equalities
\eqRef{e2i=aij e1j \Product-\Cols}{21u2},
\EqRef{U*Ei=Ui \Product-\SideNS}.
}

\DefLabeledTheorem{right eigenvalue is eigenvalue}{\SideNS-\Cols}
{
Any \SideWS \ProductType eigenvalue $b(\gi i)$
is \ProductType eigenvalue.
Eigen\ColWS $\ARow ui$
for \SideWS \ProductType eigenvalue $b(\gi i)$
is eigen\ColWS
for \ProductType eigenvalue $b(\gi i)$.
}

\DefProof{right eigenvalue is eigenvalue}
{
All entries of the row of the matrix
\DrawEq{a1-di En}{\Product-\SideNS}
are equal to zero.
Therefore, the matrix
\eqRef{a1-di En}{\Product-\SideNS}
is \ProductType singular
and \SideWS \ProductType eigenvalue $b(\gii)$
is \ProductType eigenvalue.
Since
\DrawEq{e1ij=01 \Cols}{\SideNS}
then eigenvector
\ShowEq{ARow u2i \Cols}e1i
for \SideWS \ProductType eigenvalue $b(\gii)$
is eigenvector for \ProductType eigenvalue $b(\gii)$.
The equality
\ShowEq{a1*e1=di e1 \Product-\SideNS}
follows from the definition
\refDefinition{eigenvector of matrix}{\Product-\Cols}.

From the equality
\eqRef{U-*A*U=D \Product-\SideNS}1,
it follows that we can present the matrix
\eqRef{a1-di En}{\Product-\SideNS}
as
\ShowEq{U-*A*U-di En \Product-\SideNS}
Since the matrix
\eqRef{a1-di En}{\Product-\SideNS}
is \ProductType singular,
then the matrix
\ShowEq{A-U-*di*U \Product-\SideNS}
is \ProductType singular also.
Therefore,
and \SideWS \ProductType eigenvalue $b(\gii)$
is \ProductType eigenvalue
of the pair of matrices $(a_2,u_2)$.

Equalities
\ShowEq{U-*A*U*e1=di e1 \Product-\SideNS}
\ShowEq{A*U*e1=... \Product-\SideNS}
follow from equalities
\eqRef{U-*A*U=D \Product-\SideNS}1,
\EqRef{a1*e1=di e1 \Product-\SideNS}.
From equalities
\EqRef{e1i=u2i*e2 \SideNS-\Cols},
\EqRef{A*U*e1=... \Product-\SideNS}
and from the theorem
\refTheorem{Eigenvector does not depend on basis}{\SideNS-\Cols},
it follows that
eigen\ColWS
\ShowEq{ARow u2i \Cols}u2i
of the matrix $a_2$
for \SideWS \ProductType eigenvalue $b(\gi i)$
is
eigen\ColWS
\ShowEq{ARow u2i \Cols}u2i
of the pair of matrices $(a_2,u_2)$ corresponding to \ProductType eigenvalue $b(\gii)$.
}

\DefLabeledRemark{right eigenvalue is eigenvalue}{\SideNS-\Cols}
{
The equality
\ShowEq{di*ui=ui*di \Product-\SideNS}
follows from the theorem
\refTheorem{right eigenvalue is eigenvalue}{\SideNS-\Cols}.
The equality
\EqRef{di*ui=ui*di \Product-\SideNS}
seems unusual.
However, if we slightly change it
\ShowEq{di*ui=ui*di 1 \Product-\SideNS}
then the equality
\EqRef{di*ui=ui*di 1 \Product-\SideNS}
follows from the equality
\EqRef{e1i=u2i*e2 \SideNS-\Cols}.
}

\DefLabeledTheorem[1]{Eigenvalue a ox a, c in AAA[b]}{\SideNS-\Cols}
{
Let the \ColWS vector $v$ be
eigen\ColWS for \SideWS \ProductTypeO eigenvalue $b$
of the matrix $a$.
Let $A$\Hyph number $b$ satisfy the condition
\DrawEq[b{}{s#1}]{b in AAA[a]}{bs#1-\Cols}
Then
for any polynomial
\DrawEq{p(c) p in D}{\SideNS-\Cols}
the \ColWS vector
\ShowEq{eigenvector cv \SideNS-\Cols}
\ShowEq{c=p(b)}
is also
eigen\ColWS for \SideWS \ProductTypeO eigenvalue $b$.
}

\DefProof[2]{Eigenvalue a ox a, c in AAA[b]}
{
According to the theorem
\RefTheorem{c in Za=> p(c) in Za},
the equality
\ShowEq{as ox as o cv = ..., \SideNS-\Cols}
follows from the condition
\eqRef{b in AAA[a]}{bs#1-\Cols}.
According to the definition
\refDefinition{Eigenvalue of Matrix of Linear Map}{\SideNS-\Cols},
the equality
\EqRef{as ox as o cv = ..., \SideNS-\Cols}
implies that
the \ColWS vector $#2$ is
eigen\ColWS for \SideWS \ProductTypeO eigenvalue $b$
of the matrix $a$.
}

\DefLabeledDefinition{Eigenvalue of Matrix of Linear Map}{\SideNS-\Cols}
{
$A$\Hyph number $b$ is called
\SideWS \ProductTypeO {\bf eigenvalue}
of the matrix $a$
if there exists \ColWS vector $v$ such that
\DrawEq[v]{Eigenvalue of Linear Map \SideNS-\Cols}d
The \ColWS vector $v$ is called
{\bf eigen\ColWS} for \SideWS \ProductTypeO eigenvalue $b$.
}

\DefLabeledTheorem[1]{Eigenvalue of Linear Map a ox a, 1}{\SideNS-\Cols}
{
Let entries of the matrix $a$ satisfy the equality
\DrawEq{as ox as=a1 ox a2, #1}{\SideNS-\Cols}
Then \SideWS \ProductTypeO eigenvalue $b$
is \SideWS \ProductType eigenvalue of the matrix $a_{#1}$.
}

\DefProofRef[1]{Eigenvalue of Linear Map a ox a, 1}{\SideNS-\Cols}
{
The equality
\DrawEq{as ox as o v = v a12, #1, \Cols}{1-\SideNS}
follows from the equality
\FrameEqRef{as ox as=a1 ox a2, #1}{\SideNS-\Cols}.
\newline
The equality
\ShowEq{as ox as o v = v a12, matrix, #1 1 \SideNS-\Cols}
follows from the equality
\eqRef{as ox as o v = v a12, #1, \Cols}{1-\SideNS}.
The equality
\ShowEq{as ox as o v = v a12, matrix, #1 2 \SideNS-\Cols}
follows from the definition
\refDefinition{Eigenvalue of Matrix of Linear Map}{\SideNS-\Cols}
of \SideWS \ProductTypeO eigenvalue
of the matrix $a$.
The equality
\ShowEq{as ox as o v = v a12, matrix, #1 \SideNS-\Cols}
follows from equalities
\EqRef{as ox as o v = v a12, matrix, #1 1 \SideNS-\Cols},
\EqRef{as ox as o v = v a12, matrix, #1 2 \SideNS-\Cols}.
The theorem follows from the equality
\EqRef{as ox as o v = v a12, matrix, #1 \SideNS-\Cols}
and from the definition
\refDefinition{Left and Right Eigenvalues}{\SideNS-\Cols}
of \SideWS \ProductType eigenvalue
of the matrix $a_{#1}$.
}

\DefLabeledTheorem[1]{Eigenvalue of Linear Map a ox a, 2}{\SideNS-\Cols}
{
Let entries of the matrix $a$ satisfy the equality
\DrawEq{as ox as=a1 ox a2, #1}{\SideNS-\Cols}
Then \SideWS \ProductTypeO eigenvalue $b$
is \ProductTypeA eigenvalue of the matrix $a_{#1}$.
}

\DefProofRef[1]{Eigenvalue of Linear Map a ox a, 2}{\SideNS-\Cols}
{
The equality
\DrawEq{as ox as o v = v a12, #1, \Cols}{2-\SideNS}
follows from the equality
\FrameEqRef{as ox as=a1 ox a2, #1}{\SideNS-\Cols}.
\newline
The equality
\ShowEq{as ox as o v = v a12, matrix, #1 1 \SideNS-\Cols}
follows from the equality
\eqRef{as ox as o v = v a12, #1, \Cols}{2-\SideNS}.
The equality
\ShowEq{as ox as o v = v a12, matrix, #1 2 \SideNS-\Cols}
follows from the definition
\refDefinition{Eigenvalue of Matrix of Linear Map}{\SideNS-\Cols}
of \SideWS \ProductTypeO eigenvalue
of the matrix $a$.
The equality
\ShowEq{as ox as o v = v a12, matrix, #1 \SideNS-\Cols}
follows from equalities
\EqRef{as ox as o v = v a12, matrix, #1 1 \SideNS-\Cols},
\EqRef{as ox as o v = v a12, matrix, #1 2 \SideNS-\Cols}.
The theorem follows from the equality
\EqRef{as ox as o v = v a12, matrix, #1 \SideNS-\Cols},
from the definition
\refDefinition{eigenvalue of matrix}{\ProductA}
of \ProductType eigenvalue
of the matrix $a_{#1}$
and from the definition
\refDefinition{eigenvector of matrix}{\ProductA-\Cols}
of corresponding eigenvector $v$.
}

\DefLabeledTheoremNote[1]{n eigenvectors}{#1}
{
Let $\Base$ be \Algebra.
Let $V$ be $\Base$\Hyph vector space.
If endomorphism has $\gik$ different
eigenvalues\,\footnotemark
\ShowEq{di 1n}bk,
then eigenvectors
\ShowEq{di 1n}vk{}
which correspond to different eigenvalues
are linearly independent.
}
{
See also the theorem on the page
\citeBib{Kurosh: High Algebra}\Hyph 203.
}

\AddEq{remark: n eigenvectors}
{
However, the theorem
\refTheorem{n eigenvectors}{algebra}
is not correct.
We will consider the proof separately
for left and right $A$\Hyph vector spaces.
}

\DefProof{n eigenvectors}
{
We will prove the theorem by induction with respect to $\gik$.

Since eigenvector is non\Hyph zero, the theorem is true for $\gik=\gi 1$.

\begin{Statement}
\labelStatement{theorem is true for k=m-1 \SideNS}
Let the theorem be true for $\gik=\gim-\gi 1$.
\hfill\(\odot\)
\end{Statement}

\begin{sloppypar}
According to theorems
\ShowEq{ref 1 for n eigenvectors}
we can assume that
\ShowEq{e2=e1}
Let
\ShowEq{di 1n+1}bm{}
be different eigenvalues and
\ShowEq{di 1n+1}vm{}
corresponding eigenvectors
\ShowEq{fovi=di vi \SideNS}
\ShowEq{di ne dj}
Let the statement
\RefStatement{ai vi=0 \SideNS}
be true.
\end{sloppypar}

\begin{Statement}
\labelStatement{ai vi=0 \SideNS}
There exists linear dependence
\ShowEq{ai vi=0 \SideNS}
where
\ShowEq{a1 ne 0}
\hfill\(\odot\)
\end{Statement}

The equality
\ShowEq{ai f vi=0 \SideNS}
follows from the equality
\EqRef{ai vi=0 \SideNS}
and the theorem
\ShowEq{ref 2 for n eigenvectors}
The equality
\ShowEq{ai di vi=0 \SideNS}
follows from equalities
\EqRef{fovi=di vi \SideNS},
\EqRef{ai f vi=0 \SideNS}.
Since the product is non\Hyph commutative,
we cannot reduce expression
\EqRef{ai di vi=0 \SideNS}
to linear dependence.
}

\AddEq{remark: Eigenvalue and conjugate class}
{
According to theorems
\refTheorem{Eigenvalue and conjugate class}{\SideNS-\Cols},
\refTheorem{Coefficient comutes with matrix}{\SideNS-\Cols},
the set of eigen\ColsWS
for given \SideWS \ProductType eigenvalue
is not $A$\Hyph vector space.
}

\DefLabeledTheorem{Eigenvalue and conjugate class}{\SideNS-\Cols}
{
Let $A$\Hyph number $b$ be \SideWS \ProductType eigenvalue of the matrix
$a_2$.\refFootnote{Eigenvalue and conjugate class}{\Product}
Then any $A$\Hyph number which $A$\Hyph conjugated with $A$\Hyph number $b$,
is \SideWS \ProductType eigenvalue of the matrix $a_2$.
}

\DefLabeledFootnote{Eigenvalue and conjugate class}{\Product}
{
See the similar statement and proof on the page
\citeBib{Cohn: Skew Fields}\Hyph 375.
}

\DefProof{Eigenvalue and conjugate class}
{
Let $v$ be
eigen\ColWS for \SideWS \ProductType eigenvalue $b$.
Let $c\ne 0$ be $A$\Hyph number.
The theorem follows from the equality
\ShowEq{conjugate eigenvalue \SideNS-\Cols}
}

\DefLabeledTheorem[2]{Coefficient comutes with matrix}{\SideNS-\Cols}
{
Let $v$ be eigen\ColWS for \SideWS \ProductType eigenvalue $b$.
The vector $#1#2$
is eigen\ColWS for \SideWS \ProductType eigenvalue $b$
iff
$A$\Hyph number $c$ commutes with all entries of the matrix $a_2$.
}

\DefProof[2]{Coefficient comutes with matrix}
{
The vector $#1#2$
is eigen\ColWS for \SideWS \ProductType eigenvalue $b$
iff
the following equality is true
\ShowEq{Coefficient comutes with matrix \SideNS-\Cols}
The equality
\EqRef{Coefficient comutes with matrix \SideNS-\Cols}
holds iff
$A$\Hyph number $c$ commutes with all entries of the matrix $a_2$.
}

\DefLabeledDefinition{geometric object, vector space}{\SideNS-\Cols}
{
Let us say the coordinates $w_1$ of vector $\Vector w$
with respect to the basis \eV[W1] are given.
The set of vectors
\DrawEq[w2{w_2}{W2}]{geometric object representative, \SideNS-\Cols}{w=we}
is called
{\bf geometric object}
defined in \SideWS $A$\Hyph vector space $V$ of \ColN s.
For any basis \eV[W2][,]
corresponding point
\newline
\FrameEqRef{coordinate transformation, vector space W, \SideNS-\Cols}1
\newline
of coordinate manifold defines the vector
\DrawEq[w2{w_2}{W2}]{geometric object representative, \SideNS-\Cols}{w=we 1}
which is called
{\bf representative of geometric object}
in \SideWS $A$\Hyph vector space $V$ in basis \eV[V2][.]
}

\AddEq{remark: geometric object, vector space}
{
Since a geometric object is an orbit of representation, we see that
according to theorem
\RefTheorem{proper definition of orbit}
the definition of the geometric object is a proper definition.

We also say that $\Vector w$ is
a \AddIndex{geometric object of type $F$}{type of geometric object}.

Definitions
\refDefinition{coordinates of geometric object, vector space}{\SideNS-cols},
\refDefinition{coordinates of geometric object, vector space}{\SideNS-rows}
introduce a geometric object in coordinate space.
We assume in definitions
\refDefinition{geometric object, vector space}{\SideNS-cols},
\refDefinition{geometric object, vector space}{\SideNS-rows}
that we selected
a basis of vector space $W$.
This allows using a representative of the geometric object
instead of its coordinates.
}

\DefText{Geometric Object of Vector Space}
{
\begin{sloppypar}
Let $V$, $W$ be \SideWS $A$\Hyph vector spaces of \ColN s 
and $G(V_*)$ be symmetry group
of \SideWS $A$\Hyph vector space $V$.
Homomorphism
\DrawEq[F{G(V_*)}{GL(W_*)}{}]{f: A->B}{FG \SideNS-\Cols}
maps passive transformation $g\in G(V_*)$
\DrawEq[{V1}{V2}g{}]{e2i=aij e1j \Product-\Cols}{GOeV}
of \SideWS $A$\Hyph vector spaces $V$
into passive transformation
\ShowEq{Fg in GL}
\DrawEq[{W1}{W2}{F(g)}{}]{e2i=aij e1j \Product-\Cols}{GOeW}
of \SideWS $A$\Hyph vector spaces $W$.
\end{sloppypar}

Then coordinate transformation in
\SideWS $A$\Hyph vector space $W$ gets form
\DrawEq{coordinate transformation, vector space W, \SideNS-\Cols}1
Therefore, the map
\DrawEq{F1(g)=F(g)**-\Product}{\SideNS}
is \OtherSideNS\Hyph side representation
\ShowEq{f:A->*B}{F_1}{G(V_*)}{W_*}
of group $G(V_*)$ in the set $W_*$.
}

\DefProof{passive transformation and product of endomorphisms}
{
\ShowEq{two endomorphisms and two bases}
According to the theorem
\refTheorem{product of homomorphisms, A vector space}{\SideNS-\ColN},
\DrawEq[{}{}]{f=f1*f2 endo \SideNS-\Cols}{\Product}
\DrawEq[{'}{}]{f=f1*f2 endo \SideNS-\Cols}{\Product'}
According to the theorem
\refTheorem{passive transformation and endomorphism}{\SideNS-\Cols},
\DrawEq[{}{}]{f'=g f g- \SideNS-\Cols}{f*}
\DrawEq[1]{f'=g f g- \SideNS-\Cols}{1*}
\DrawEq[2]{f'=g f g- \SideNS-\Cols}{2*}
The equality
\ShowEq{f1*f2=g.g- \SideNS-\Cols}
follows from equalities
\ShowEq{f1*f2=g.g-}
From equalities
\eqRef{f=f1*f2 endo \SideNS-\Cols}{\Product'},
\EqRef{f1*f2=g.g- \SideNS-\Cols}
it follows that the equality
\eqRef{f=f1*f2 endo \SideNS-\Cols}{\Product}
is covariant with respect to choice of basis 
of \SideWS $A$\Hyph vector space $V$.
}

\AddEq{theorem: finite dimension->representation is free}
{
\begin{ShadedTheorem}
\labelTheorem{finite dimension->representation is free, \SideNS}
If \SideWS $A$\Hyph vector space $V$
has finite dimension,
then corresponding representation is free.
\end{ShadedTheorem}
}

\DefProof{finite dimension->representation is free}
{
The theorem follows from the definition
\RefDefinition{free representation of algebra}
and from theorems
\refTheorem{av=bv=>a=b}{\SideNS-cols},
\refTheorem{av=bv=>a=b}{\SideNS-rows}.
}

\AddEq[1]{Let V be vector space}
{
Let $V$ be a \SideWS $A$\Hyph vector space#1
}

\AddEq[1]{Let V be vector space of}
{
Let $V$ be a \SideWS $A$\Hyph vector space of \ColN s#1
}

\AddEq{Let V be vector space of and}
{
Let $V$ be a \SideWS $A$\Hyph vector space of \ColsWS and
}

\AddEq{Let V be vector space of and basis}
{
\ShowEq{Let V be vector space of and}
$\Basis e$ be basis of \SideWS $A$\Hyph vector space $V$.
}

\AddEq{Let V be vector space and basis}
{
\ShowEq{Let V be vector space}{}
and $\Basis e$ be basis of \SideWS $A$\Hyph vector space $V$.
}

\AddEq{remark: automorphism of vector space, group}
{
\ShowEq{Let V be vector space of}.
We proved in the theorem
\refTheorem{automorphism of vector space}{\Product-\Cols}
that, if we select a basis
of \SideWS $A$\Hyph vector space $V$,
then we can identify any automorphism $\Vector f$
of \SideWS $A$\Hyph vector space $V$
with \ProductType nonsingular matrix $f$.
Corresponding transformation of coordinates of vector
\DrawEq[{v'}{}v{}f{}]{v1=v2*a \SideNS-\Cols}{}
is called
\AddIndex{linear transformation}{linear transformation}.
}

\DefLabeledDefinition{Symmetry Group}{\SideNS-\Cols}
{
Normal subgroup $G(V)$ of the group $GL(V)$ such that subgroup $G(V)$ generates
automorphisms which hold properties of the selected structure
is called
{\bf symmetry group}.

Without loss of generality we identify element $g$ of group $G(V)$
with corresponding transformation of representation
and write its action on vector $v\in V$ as
$v\ProductVal g$.
}

\DefText{definition: linear G* representation}
{
The \OtherSideNS\Hyph side representation
of group $G$ in \SideWS $A$\Hyph vector space is called
\AddIndex{linear $G$\Hyph representation}{linear G* representation}.
}

\DefText{linear G* representation}
{
Let us define an additional structure on \SideWS $A$\Hyph vector space $V$.
Then not every automorphism keeps properties of the selected structure.
For imstance, if we introduce norm in
\SideWS $A$\Hyph vector space $V$,
then we are interested in automorphisms which preserve the norm of the vector.
}

\DefLabeledDefinition{active G-representation}{\SideNS}
{
The \OtherSideNS\Hyph side representation
of group $G(V)$ in the set of bases of \SideWS $A$\Hyph vector space $V$ is called
\AddIndex{active \SideNS\Hyph side $G$\Hyph representation}
{active representation, vector space}.
}

\DefText{active G-representation}
{
If the basis \eV is given, then
we can identify the automorphism
of \SideWS $A$\Hyph vector space $V$
and its coordinates with respect to the basis \eV[][.]
The set $G(V_*)$ of coordinates of automorphisms
with respect to the basis \eV is group
isomorphic to the group $G(V)$.

Not every two bases can be mapped by a transformation
from the symmetry group
because not every nonsingular linear transformation belongs to
the representation of group $G(V)$.
Therefore, we can represent
the set of bases as a union of orbits of group $G(V)$.
In particular, if the basis $\eV\in G(V)$,
then the group $G(V)$ is orbit of the basis \eV[][.]
}

\DefLabeledTheorem{representation on basis manifold, vector space}{\SideNS}
{
Active \OtherSideNS\Hyph side $G(V)$\Hyph representation on basis manifold
is single transitive representation.
}

\DefProof{representation on basis manifold, vector space}
{
The theorem follows from the theorem
\refTheorem{active representation is single transitive}{\SideNS-cols}
and the definition
\refDefinition{basis manifold of vector space}{\SideNS-cols},
as well the theorem follows from the theorem
\refTheorem{active representation is single transitive}{\SideNS-rows}
and the definition
\refDefinition{basis manifold of vector space}{\SideNS-rows}.
}

\DefLabeledDefinition{basis manifold of vector space}{\SideNS-\Cols}
{
We call orbit
\ShowEq{basis manifold of vector space, \SideNS-\Cols}
of the selected basis $\Basis e$
the {\bf basis manifold}
of \SideWS $A$\Hyph vector space $V$ of \ColN s.
}

\DefLabeledTheorem{active transformations, vector space}{\SideNS-\Cols}
{
\ShowEq{Let V be vector space of and basis}
Automorphisms of \SideWS $A$\Hyph vector space of \ColsWS form
a \OtherSideNS\Hyph side linear
effective \Group nA\Hyph representation.
}

\DefProof{active transformations, vector space}
{
Let $a$, $b$ be matrices of automorphisms $\Vector a$ and $\Vector b$
with respect to basis $\Basis e$.
According to the theorem
\refTheorem{homomorphism A module}{\SideNS-\Cols(1)},
coordinate transformation has following form
\DrawEq[{v'}{}v{}a{}]{v1=v2*a \SideNS-\Cols}{v'}
\DrawEq[{v''}{}{v'}{}b{}]{v1=v2*a \SideNS-\Cols}{v''}
The equality
\ShowEq{a*b, \SideNS-\Cols}
\DrawEq[{v''}{}v{}{\TheProduct}{}]{v1=v2*a \SideNS-\Cols}{}
\begin{sloppypar}
\noindent
follows from equalities
\eqRef{v1=v2*a \SideNS-\Cols}{v'},
\eqRef{v1=v2*a \SideNS-\Cols}{v''}.
According to the theorem
\refTheorem{product of homomorphisms, A vector space}{\SideNS-\ColN},
the product of automorphisms $\Vector a$ and $\Vector b$
has matrix \TheProduct.
Therefore, automorphisms of \SideWS $A$\Hyph vector space of \ColsWS form
a \OtherSideNS\Hyph side linear
\Group nA\Hyph representation.
\end{sloppypar}

It remains to  prove that
the kernel of inefficiency consists only of identity.
Identity transformation
satisfies to equation
\ShowEq{active transformations, vector space, 2, \SideNS-\Cols}
Choosing values of coordinates as
\ShowEq{ai=delta ik, \Cols}
where we selected $\gik$ we get
\ShowEq{identity transformation, vector space, \SideNS-\Cols}
From \EqRef{identity transformation, vector space, \SideNS-\Cols} it follows
\ShowEq{identity transformation, vector space, 1, \Cols}
Since $\gik$ is arbitrary, we get the conclusion $a=\delta$.
}

\DefLabeledTheorem{active transformations, basis, vector space}{\SideNS-\Cols}
{
Automorphism $a$ acting on each vector of basis of
\SideWS $A$\Hyph vector space of \ColN s
maps a basis into another basis.
}

\DefProof{active transformations, basis, vector space}
{
Let $\Basis e$ be basis of \SideWS $A$\Hyph vector space $V$ of \ColN s.
According to theorem
\refTheorem{automorphism of vector space}{\Product-\Cols},
vector
\ShowEq{vector of basis \Cols}e
maps into a vector
\ShowEq{vector of basis \Cols}{e'}
\DrawEq{automorphism, vector space, 1, \SideNS-\Cols}1
\begin{sloppypar}
\noindent
Let vectors
\ShowEq{vector of basis \Cols}{e'}
be linearly dependent.
Then $\lambda\ne 0$
in the equality
\ShowEq{automorphism, vector space, 2, \SideNS-\Cols}
\end{sloppypar}
\noindent
From equations \eqRef{automorphism, vector space, 1, \SideNS-\Cols}1
and \EqRef{automorphism, vector space, 2, \SideNS-\Cols} it follows that
\ShowEq{automorphism, vector space, 3, \SideNS-\Cols}
and $\lambda\ne 0$. This
contradicts to the statement that vectors
\ShowEq{vector of basis \Cols}e
are linearly independent.
Therefore vectors
\ShowEq{vector of basis \Cols}{e'}
are linearly independent
and form basis.
}

\DefText{extend linear representation to set of bases 1}
{
Thus we can extend
a \OtherSideNS\Hyph side linear $GL(V_*)$\Hyph representation
in \SideWS $A$\Hyph vector space $V_*$
to the set of bases
of \SideWS $A$\Hyph vector space $V$.
}

\DefText{extend linear representation to set of bases 2}
{
Transformation of this \OtherSideNS-side representation
on the set of bases
of \SideWS $A$\Hyph vector space $V$ is called
\AddIndex{active transformation}{active transformation, vector space}
because the homomorphism of the \SideWS $A$\Hyph vector space induced this transformation
(\BlueText{See also definition in the section
\citeBib{Korn}\Hyph 14.1\Hyph 3
as well the definition on the page
\citeBib{Rashevsky}\Hyph 214}).

\ePrints{2022.01.05,2022.11.26}%
\ifx\Semafor\ValueOn%
\FrameCiteBib{Korn}

\FrameCiteBib{Rashevsky}
\fi
}

\DefText{extend linear representation to set of bases 3}
{
\begin{sloppypar}
According to definition we write the action
of the transformation $a\in GL(V_*)$ on the basis $\Basis e$ as
\ShowEq{active transformation, vector space, \SideNS-\Cols}.
Consider the equality
\DrawEq{active transformation ae x=a ex, \SideNS-\Cols}1
The expression
\ShowEq{active transformation, vector space, \SideNS-\Cols}{}
on the left side of the equality
\eqRef{active transformation ae x=a ex, \SideNS-\Cols}1
is image of basis \eV with respect to active transformation $a$.
The expression
\ShowEq{active transformation ae x=a ex 1, \SideNS-\Cols}
on the right side of the equality
\eqRef{active transformation ae x=a ex, \SideNS-\Cols}1
is expansion of vector $\Vector v$ with respect to basis \eV[][.]
Therefore,
the expression on the right side of the equality
\eqRef{active transformation ae x=a ex, \SideNS-\Cols}1
is image of vector $\Vector v$ with respect to endomorphism $a$
and the expression on the left side of the equality
\newline
\FrameEqRef{active transformation ae x=a ex, \SideNS-\Cols}1
\newline
is expansion of image of vector $\Vector v$
with respect to image of basis \eV[][.]
Therefore, from the equality
\eqRef{active transformation ae x=a ex, \SideNS-\Cols}1
it follows that endomorphism $a$ of \SideWS $A$\Hyph vector space and
corresponding active transformation $a$ act synchronously
and coordinates $a\circ \Vector v$ of image of the vector $\Vector v$ with respect to
the image
\ShowEq{active transformation, vector space, \SideNS-\Cols}{}
of the basis \eV
are the same as coordinates of the vector $\Vector v$ with respect to the basis $\Basis e$.
\end{sloppypar}
}

\DefLabeledTheorem{active representation is single transitive}{\SideNS-\Cols}
{
Active \OtherSideNS\Hyph side \Group nA\Hyph representation on the set of bases
is single transitive representation.
The set of bases identified with tragectory
\ShowEq{basis manifold of V, \SideNS-\Cols}e{\Group nA}{}
of active \OtherSideNS\Hyph side \Group nA\Hyph representation
is called
the {\bf basis manifold}
of \SideWS $A$\Hyph vector space $V$.
}

\DefLabeledTheorem{exists representation, commuting with active}{\SideNS-\Cols}
{
On the basis manifold
\ShowEq{basis manifold of V, \SideNS-\Cols}eG{}
of \SideWS $A$\Hyph vector space of \ColN s,
there exists single transitive
\SideNS\Hyph side $G(V)$\Hyph representation, commuting with active.
}

\AddEq{proof: exists representation, commuting with active}
{
\begin{proof}
Theorem
\refTheorem{representation on basis manifold, vector space}{\SideNS}
means that the basis manifold
\ShowEq{basis manifold of V, \SideNS-\Cols}eG{}
is a homogenous space of group $G$.
The theorem follows from the theorem
\RefTheorem{two representations of group}.
\end{proof}
}

\AddEq{remark: passive transformation, vector space}
{
\ePrints{2022.11.26,2306.00880}%
\ifx\Semafor\ValueOn%
Transformation of
\else
As we see from remark
\RefRemark{one representation of group}
transformation of
\fi
\SideNS\Hyph side $G(V)$\Hyph representation is different from an active transformation
and cannot be reduced to
transformation of space $V$.
}

\DefLabeledDefinition{passive transformation, vector space}{\SideNS-\Cols}
{
A transformation of
\SideNS\Hyph side $G(V)$\Hyph representation is called
{\bf passive transformation}
of basis manifold
\ShowEq{basis manifold of V, \SideNS-\Cols}eG{}
of \SideWS $A$\Hyph vector space of \ColN s,
and the \SideNS\Hyph side $G(V)$\Hyph representation is called
{\bf passive \SideNS\Hyph side $G(V)$\Hyph representation}.
According to the definition
we write the passive transformation of basis $\Basis e$
defined by element $a\in G(V)$ as
\ShowEq{passive transformation symbol, \SideNS-\Cols}
}

\AddEq[1]{table: Passive and Active Representations}
{
\ShowEq{def #1}\SideWS $A$\Hyph vector space&\ShowEq{def #1}\ShowEq{\DefRow}\Group nA&
\ShowEq{def #1}\OtherSideNS-side&\ShowEq{def #1}\SideNS-side\\
of rows&&&\\[10pt]
\hline
\ShowEq{def #1}\SideWS $A$\Hyph vector space&\ShowEq{def #1}\ShowEq{\DefCol}\Group nA&
\ShowEq{def #1}\OtherSideNS-side&\ShowEq{def #1}\SideNS-side\\
of columns&&&\\
\hline
}

\DefLabeledTheorem{principle of covariance}{\SideNS}
{
{\bf(Principle of covariance).}
Representative of geometric object does not depend on selection
of basis \eV[V2][.]
}

\DefProof{principle of covariance}
{
To define representative of geometric object,
we need to select basis $\Basis e_V$,
basis $\Basis e_W$
and coordinates of geometric object $w^\alpha$.
Corresponding representative of geometric object
has form
\DrawEq[w{}wW]{geometric object representative, \SideNS-\Cols}{}
Suppose we map basis $\Basis e_V$
to basis $\Basis e'_V$
by passive transformation
\DrawEq[V]{passive transformation e->e', \SideNS-\Cols}{}
According construction this generates passive transformation
\newline
\FrameEqRef[{W1}{W2}{F(g)}{}]{e2i=aij e1j \Product-\Cols}{GOeW}
\newline
and coordinate transformation
\newline
\FrameEqRef{coordinate transformation, vector space W, \SideNS-\Cols}1
\newline
Corresponding representative of geometric object has form
\ShowEq{invariance principle 3, \SideNS-\Cols}
Therefore representative of geometric object
is invariant relative selection of basis.
}

\AddEq{definition: sum of geometric objects}
{
\begin{ShadedDefinition}
\labelDefinition{sum of geometric objects, \SideNS-\Cols}
\ShowEq{Let V be vector space of and}
\ShowEq{sum of geometric objects, 1}
be geometric objects of the same type
defined in \SideWS $A$\Hyph vector space $V$.
Geometric object
\ShowEq{sum of geometric objects, 2}
is called \AddIndex{sum
\ShowEq{sum of geometric objects, 3}
of geometric objects}{sum of geometric objects}
$\Vector w_1$ and $\Vector w_2$.
\end{ShadedDefinition}
}

\AddEq{definition: product of geometric object and constant}
{
\begin{ShadedDefinition}
\labelDefinition{product of geometric object and constant, \SideNS-\Cols}
\ShowEq{Let V be vector space of and}
\DrawEq[w1{w_1}W]{geometric object representative, \SideNS-\Cols}{}
be geometric object
defined in \SideWS $A$\Hyph vector space $V$.
Geometric object
\ShowEq{product of geometric object and constant, 2, \SideNS}
is called \AddIndex{product
\ShowEq{product of geometric object and constant, 3, \SideNS}
of geometric object $\Vector w_1$ and constant $k\in A$}
{product of geometric object and constant}.
\end{ShadedDefinition}
}

\AddEq{theorem: Geometric objects form A vector space}
{
\begin{ShadedTheorem}
\labelTheorem{Geometric objects form A vector space, \SideNS-\Cols}
Geometric objects of type $F$
defined in \SideWS $A$\Hyph vector space $V$ of \ColN s
form \SideWS $A$\Hyph vector space of \ColN s.
\end{ShadedTheorem}
}

\DefProof{Geometric objects form A vector space}
{
The statement of theorems
\RefTheorem{Geometric objects form A vector space, \SideNS-cols},
\RefTheorem{Geometric objects form A vector space, \SideNS-rows}
follows from immediate verification
of the properties of vector space.
}

\DefLabeledTheorem{coordinate matrix of basis and passive transformation}{\SideNS-\Cols}
{
The coordinate matrix of basis $\Basis e'$ relative basis $\Basis e$
of \SideWS $A$\Hyph vector space $V$ of \ColsWS
is identical with the matrix of passive transformation mapping
basis $\Basis e$ to basis $\Basis e'$.
}

\DefProof{coordinate matrix of basis and passive transformation}
{
According to the theorem
\refTheorem{coordinate matrix of vector}{\SideNS-\Cols},
the coordinate matrix of basis $\Basis e'$ relative basis $\Basis e$
consist from \Rows which are coordinate matrices of vectors
$\aD{\Vector e'}i$ relative the basis $\Basis e$. Therefore,
\ShowEq{coordinate matrix of basis and passive transformation, \SideNS-\Cols}{e'}
At the same time the passive transformation $a$ mapping one basis to another has a form
\ShowEq{coordinate matrix of basis and passive transformation, \SideNS-\Cols}a
According to the theorem
\refTheorem{coordinates of vector}{\SideNS-\Cols},
\ShowEq{coordinate matrix of basis and passive transformation, \Cols}
for any $\gii$. This proves the theorem.
}

\DefLabeledSloppyRemark{identify basis and matrix of coordinates}{\SideNS-\Cols}
{
According to theorems
\RefTheorem{single transitive representation of group},
\refTheorem{coordinate matrix of basis and passive transformation}{\SideNS-\Cols},
we can identify a basis $\Basis e'$ 
of the basis manifold
\ShowEq{basis manifold of V, \SideNS-\Cols}eG{}
of \SideWS $A$\Hyph vector space of columns
and the matrix $g$ of coordinates of the basis $\Basis e'$
with respect to the basis $\Basis e$
\DrawEq[{}{}]{passive transformation e->e', \SideNS-\Cols}1
Since $\Basis e'=\Basis e$, then the equality
\eqRef{passive transformation e->e', \SideNS-\Cols}1
gets form
\ShowEq{passive transformation e->e, \SideNS-\Cols}
Based on the equality
\EqRef{passive transformation e->e, \SideNS-\Cols},
we will use notation
\ShowEq{basis manifold of V, \SideNS-\Cols}{\aD En}G{}
for basis manifold
\ShowEq{basis manifold of V, \SideNS-\Cols}eG.
}

\DefText{Passive Transformation}
{
An active transformation changes bases and vectors uniformly
and coordinates of vector relative basis do not change.
A passive transformation changes only the basis and it leads to transformation
of coordinates of vector relative to basis.
}

\DefText{transformation of coordinates of vector}
{
We consider transformation of coordinates of vector in the theorem
\refTheorem{passive transformation of vector space}{\SideNS-\Cols}.
}

\AddEq{remark: active representation is single transitive}
{
To prove theorems
\refTheorem{active representation is single transitive}{\SideNS-cols},
\refTheorem{active representation is single transitive}{\SideNS-rows},
it is sufficient to show that
at least one transformation of \OtherSideNS\Hyph side representation is defined for any two bases
and this transformation is unique.
}

\AddEq{left active representation is single transitive}
{
coordinate matrix of original basis over matrix of automorphism
\ShowEq{left e'=e*a}{}
}

\AddEq{right active representation is single transitive}
{
matrix of automorphism over coordinate matrix of original basis
\ShowEq{right e'=e*a}{}
}

\DefProofSloppy{active representation is single transitive}
{
Homomorphism $a$ operating on basis $\Basis e$ has form
\DrawEq{automorphism, vector space, 1, \SideNS-\Cols}{}
where
\ShowEq{vector of basis \Cols}{e'}
is coordinate matrix of vector
\ShowEq{vector of basis \Cols}{\Vector e'}
relative basis $\Basis h$ and
\ShowEq{vector of basis \Cols}e
is coordinate matrix of vector
\ShowEq{vector of basis \Cols}{\Vector e}
relative basis $\Basis h$.
Therefore, coordinate matrix of image of basis equal to
\ProductType product of
\ShowEq{\SideWS active representation is single transitive}
According to the theorem
\refTheorem{coordinate matrix of basis}{\Product-\Cols},
matrices $g$ and $e$ are nonsingular. Therefore, matrix
\ShowEq{homomorphism on A basis, \SideNS-\Cols}
is the matrix of automorphism mapping basis $\Basis e$
to basis $\Basis e'$.

Suppose elements $g_1$, $g_2$ of group $G$ and basis $\Basis e$ satisfy equation
\ShowEq{two transformations on basis manifold, \SideNS-\Cols}
According to the theorems
\refTheorem{coordinate matrix of basis}{\Product-\Cols}
and
\RefTheorem{two cr-products equal},
we get $g_1=g_2$.
This proves statement of theorem.
}

\DefLabeledTheorem{representation of homomorphism relative different bases}{\SideNS-\Cols}
{
Let
\ShowEq{Bases e12}V{}
be bases of \SideWS $A$\Hyph vector space $V$ of columns.
Let
\ShowEq{Bases e12}W{}
be bases of \SideWS $A$\Hyph vector space $W$ of columns.
Let $a_1$ be matrix of homomorphism
\DrawEq{A:V->W}{different bases \SideNS-\Cols}
relative to bases
\ShowEq{Bases eVW}VW1{}
and $a_2$ be matrix of homomorphism
\eqRef{A:V->W}{different bases \SideNS-\Cols}
relative to bases
\ShowEq{Bases eVW}VW2.
Suppose the basis $\Basis e_{V1}$ has coordinate matrix $b$ relative
the basis $\Basis e_{V2}$
\DrawEq[Vb]{coordinate matrix, f, g, \SideNS-\Cols}{}
and $\Basis e_{W1}$ has coordinate matrix $c$ relative
the basis $\Basis e_{W2}$
\DrawEq[Wc]{coordinate matrix, f, g, \SideNS-\Cols}C
Then there is relationship between matrices $a_1$ and $a_2$
\DrawEq[c]{representation of homomorphism relative different bases, \SideNS-\Cols}{}
}

\DefProof{representation of homomorphism relative different bases}
{
Vector $\Vector v\in V$ has expansion
\ShowEq{expansion of vector v, \SideNS-\Cols}
relative to bases
\ShowEq{Bases e12}V.
Since $a$ is homomorphism, we can write it as
\ShowEq{representation of homomorphism relative different bases, 1, \SideNS-\Cols}
relative to bases
\ShowEq{Bases eVW}VW1{}
and as
\ShowEq{representation of homomorphism relative different bases, 2, \SideNS-\Cols}
relative to bases
\ShowEq{Bases eVW}VW2.
According to the theorem
\refTheorem{coordinate matrix of basis}{\Product-\Cols},
matrix $c$ has \ProductType inverse and from equation
\eqRef{coordinate matrix, f, g, \SideNS-\Cols}C it follows that
\ShowEq{coordinate matrix, W2, W1, \SideNS-\Cols}
The equality
\ShowEq{representation of homomorphism relative different bases, 3, \SideNS-\Cols}
\begin{sloppypar}
\noindent
follows from equalities
\EqRef{representation of homomorphism relative different bases, 2, \SideNS-\Cols},
\EqRef{coordinate matrix, W2, W1, \SideNS-\Cols}.
From the theorem
\refTheorem{coordinates of vector}{\SideNS-\Cols}
and comparison of equations
\EqRef{representation of homomorphism relative different bases, 1, \SideNS-\Cols}
and \EqRef{representation of homomorphism relative different bases, 3, \SideNS-\Cols} it follows that
\ShowEq{representation of homomorphism relative different bases, 4, \SideNS-\Cols}
Since vector $a$ is arbitrary vector,
the statement of theorem follows
from the theorem
\refTheorem{active transformations, vector space}{\SideNS-\Cols}
and equation
\EqRef{representation of homomorphism relative different bases, 4, \SideNS-\Cols}.
\end{sloppypar}
}

\DefLabeledTheorem{coordinate transformations form representation, vector space}{\SideNS-\Cols}
{
Coordinate transformations
\newline
\FrameEqRef[v2v1g{}{\SideNS}{\Cols}]{v2=v1g-}{12 \SideNS-\Cols}
\newline
form effective linear
\OtherSideNS\Hyph side contravariant $G$\Hyph representation which is called
{\bf coordinate representation in \SideWS $A$\Hyph vector space}
of \ColN s.
}

\DefProof{coordinate transformations form representation, vector space}
{
\ShowText{passive representation of left vector space}
According to the definition
\RefDefinition{Left side contravariant representation},
coordinate transformations
form linear right\Hyph side contravariant $G$\Hyph representation.
Suppose coordinate transformation does not change vectors $\delta_k$.
Then unit of group $G$ corresponds to it because representation
is single transitive. Therefore,
coordinate representation is effective.
}

\AddEq{theorem: representation of automorphism relative different bases}
{
\begin{ShadedTheorem}
\labelTheorem{representation of automorphism relative different bases, \SideNS-\Cols}
Let $\Vector a$ be automorphism
of left $A$\Hyph vector space $V$ of columns.
Let $a_1$ be matrix of this automorphism defined
relative to basis $\Basis e_1$ and
$a_2$ be matrix of the same automorphism defined
relative to basis $\Basis e_2$.
Suppose the basis $\Basis e_1$ has coordinate matrix $b$ relative
the basis $\Basis e_2$
\DrawEq[{}b]{coordinate matrix, f, g, \SideNS-\Cols}{}
Then there is relationship between matrices $a_1$ and $a_2$
\DrawEq[c]{representation of homomorphism relative different bases, \SideNS-\Cols}{}
\end{ShadedTheorem}
}

\DefProof{representation of automorphism relative different bases}
{
Statement follows from theorem
\refTheorem{representation of homomorphism relative different bases}{\SideNS-\Cols},
because in this case $c=b$.
}

\DefLabeledTheorem[4]{basis of vector space}{\SideNS-\Cols}
{
\ShowEq{Let V be vector space of}.
\ShowEq{Let be basis of vector space}{#2}iIA{#1}V{#2}{\Cols}-
\ShowEq{Let be basis of vector space}{#4}jJA{#3}V{#4}{\Cols}-
If $|I|$ and $|J|$ are finite numbers then
\ShowEq{|I|=|J|}
}

\DefProof{basis of vector space}
{
Let
\ShowEq{|I|=m}Im{}
and
\ShowEq{|I|=m}Jn.
Let
\DrawEq{m<n gi}{1 \SideNS-\Cols}
Because $\Basis e_1$ is a basis, any vector
\ShowEq{basis e2 of V \Cols}
has expansion
\ShowEq{basis e2 relative e1 \SideNS-\Cols}
Because $\Basis e_2$ is a basis,
\DrawEq{basis e2 of V lambda}{\SideNS-\Cols}
should follow from
\ShowEq{basis e2 relative e1, lambda, \SideNS-\Cols}
Because $\Basis e_1$ is a basis we get
\ShowEq{basis e2 relative e1, lambda=0, \SideNS-\Cols}
According to
\eqRef{m<n gi}{1 \SideNS-\Cols},
\ShowEq{Rank a<m \Product}
and system of linear equations
\EqRef{basis e2 relative e1, lambda=0, \SideNS-\Cols}
has more variables then equations. According to the theorem
\RefTheorem{star rows system of linear equations, solution},
we get the statement $\lambda\ne 0$ which contradicts
statement \eqRef{basis e2 of V lambda}{\SideNS-\Cols}. 
Therefore, statement
\DrawEq{m<n gi}-
is not valid.
In the same manner we can prove that the statement
\ShowEq{n<m gi}
is not valid.
This completes the proof of the theorem.
}

\AddEq{definition: linearly independent, AoxA module}
{
\begin{ShadedDefinition}
\labelDefinition{linearly independent, AoxA module \Cols}
Let $V$ be left \AoxA A\Hyph module of \ColN s.
The set of vectors
\ShowEq{Vector A \Cols}
of left \AoxA A\Hyph module $V$ is
{\bf linearly independent}
if
\ShowEq{AoxA linearly independent 1 \Cols},
follows from the equaility
\ShowEq{AoxA linearly independent \Cols}
Otherwise the set of vectors $\aD ai$
is {\bf linearly dependent}.
\end{ShadedDefinition}
}

\AddEq{definition: basis, AoxA module}
{
\begin{ShadedDefinition}
\labelDefinition{basis, AoxA module \Cols}
We call set of vectors
\ShowEq{basis, AoxA module \Cols}
\AddIndex{basis}{basis}
for \AoxA A\Hyph module $V$
if the set of vectors
\ShowEq{basis vector \Cols}
is linearly independent and adding to this set any other vector
we get a new set which is linearly dependent.
\end{ShadedDefinition}
}

\DefLabeledDefinition{Eigenvalue of Endomorphism}{\SideNS}
{
\ShowEq{Let V be vector space and basis}
The vector $v\in V$ is called
{\bf eigenvector}
of the endomorphism
\ShowEq{f:A->B}{\Vector f}VV
with respect to the basis $\Basis e$,
if there exists $b\in A$ such that
\DrawEq[bv]{fov=b e o v}{\SideNS}
$A$\Hyph number $b$ is called
{\bf eigenvalue}
of the endomorphism $f$
with respect to the basis $\Basis e$.
}

\DefLabeledTheorem{no endomorphism fov=bv}{\SideNS}
{
\ShowEq{Let V be vector space and basis}
If
\ShowEq{b not in ZA},
then there is no endomorphism $\Vector f$ such that
\ShowEq{fov=bv \SideNS}
}

\DefProofSloppy{no endomorphism fov=bv}
{
Let there exist endomorphism $\Vector f$
such that the equality
\EqRef{fov=bv \SideNS}
is true.
The equality
\ShowEq{b(av)=a(bv) \SideNS}
follows from equalities
\ShowRef{b(av)=a(bv)}
and from the equality
\EqRef{fov=bv \SideNS}.
According to the theorem
\refTheorem{av=bv=>a=b}{\SideNS-\Cols},
the equality
\DrawEq{ab=ba}{\SideNS}
for any $a\in A$
follows from the equality
\EqRef{b(av)=a(bv) \SideNS}.
The equality
\eqRef{ab=ba}{\SideNS}
for any $a\in A$
contradicts to the statement
\ShowEq{b not in ZA}.
Therefore, considered endomorphism $\Vector f$
does not exist.
}

\AddEq{remark: Eigenvalue of Endomorphism}
{
\begin{sloppypar}
According to the theorem
\refTheorem{no endomorphism fov=bv}{\SideNS},
as opposed to the \OtherSideNS\Hyph side product
of vector over scalar,
\SideNS\Hyph side product
of vector over scalar
is not endomorphism
of \SideWS $A$\Hyph vector space.
Comparing the definition
\refDefinition{similarity transformation}{\SideNS}
and the theorem
\refTheorem{no endomorphism fov=bv}{\SideNS},
we see that similarity transformation
in non\Hyph commutative algebra
and similarity transformation
in commutative algebra
play the same role.
\end{sloppypar}
}

\DefLabeledTheorem{similarity transformation, change of basis}{\SideNS-\Cols}
{
\ShowEq{Let V be vector space of and}
\ShowEq{basis e1 e2}
be bases of \SideWS $A$\Hyph vector space $V$.
Let $g$
be passive transformation
of basis $\Basis e_1$ into basis $\Basis e_2$
\DrawEq[12g{}]{e2i=aij e1j \Product-\Cols}{}
The similarity transformation
\ShowEq{\SideWS map a En}a1
has the matrix
\DrawEq[ag]{ga*g- \SideNS-\Cols}a
with respect to the basis $\Basis e_2$.
}

\DefProof{similarity transformation, change of basis}
{
According to the theorem
\refTheorem{passive transformation and endomorphism}{\SideNS-\Cols},
the similarity transformation
\ShowEq{\SideWS map a En}a1
has the matrix
\ShowEq{ga*g- 1 \SideNS-\Cols}
with respect to the basis $\Basis e_2$.
The theorem follows from the expression
\EqRef{ga*g- 1 \SideNS-\Cols}.
}

\DefLabeledTheorem{similarity transformation}{\SideNS-\Cols}
{
\ShowEq{Let V be vector space of and basis}
Let
\ShowEq{dim V=n}
and $\aD En$ be \nTimes identity matrix.
For any $A$\Hyph number $a$, there exist endomorphism
\ShowEq{\SideWS map a En}a{}
of \SideWS $A$\Hyph vector space $V$
which has the matrix
\ShowEq{\SideWS a En}
whith respecpect to the basis $\Basis e$.
}

\DefProof{similarity transformation}
{
Let $\Basis e$ be the basis
of \SideWS $A$\Hyph vector space of columns.
According to the theorem
\refTheorem{matrix generates A module homomorphism}{\SideNS-\Cols(1)}
and the definition
\refDefinition{end aut morphism module}{\SideNS},
the map
\ShowEq{v in V -> av in V \SideNS-\Cols}
is endomorphism
of \SideWS $A$\Hyph vector space $V$.
}

\DefLabeledDefinition[6]{isomorphism module}{\SideNS(#1#2#3)}
{
Homomorphism\,\refFootnote{iso end aut morphism}{\SideNS(#1#2#3)}
\DrawEq[gf]{homomorphism A module #1#2#3}{}
is called
\AddIndex{isomorphism}{isomorphism}
between \SideWS $A_{#2}$\Hyph \VectorSet $V_{#3}$
and \SideWS $A_{#5}$\Hyph \VectorSet $V_{#6}$,
if there exists the map
\DrawEq[gf]{homomorphism- A module #1#2#3}{}
which is homomorphism.
}

\DefLabeledDefinition{end aut morphism module}{\SideNS}
{
A homomorphism\,\refFootnote{iso end aut morphism}{\SideNS(1)}
\DrawEq[fVV{}]{f: A->B}{}
in which source and target are the same
\SideWS $\Base$\Hyph\VectorSet
is called
\AddIndex{endomorphism}{endomorphism}.
Endomorphism
\DrawEq[fVV{}]{f: A->B}{}
of \SideWS $\Base$\Hyph\VectorSet $V$
is called
\AddIndex{automorphism}{automorphism},
if there exists the map $f^{-1}$
which is endomorphism.
}

\DefLabeledDefinition{iso end aut morphism vector space}{\SideNS}
{
Homomorphism
\ShowEq{f:A->B}fVW
is called\,\refFootnote{iso end aut morphism}1
\AddIndex{isomorphism}{isomorphism}
between \SideWS $A$\Hyph vector spaces $V$ and $W$,
if correspondence $f^{-1}$ is homomorphism.
A homomorphism
\ShowEq{f:A->B}fVV
in which source and target are the same
\SideWS $A$\Hyph vector space
is called
\AddIndex{endomorphism}{endomorphism}.
Endomorphism
\ShowEq{f:A->B}fVV
of \SideWS $A$\Hyph vector space $V$
is called
\AddIndex{automorphism}{automorphism},
if correspondence $f^{-1}$ is endomorphism.
}

\DefText{Basis Manifold}
{
\section{Dimension of \SideWSC \texorpdfstring{$A$}{A}\Hyph Vector Space}

\ShowText{Dimension of Vector Space}

\section{Basis Manifold for \SideWSC \texorpdfstring{$A$}{A}-Vector Space}

\ShowEq{Basis Manifold for Vector Space}

\section{Passive Transformation in \SideWSC \texorpdfstring{$A$}{A}-Vector Space}

\ShowEq{Passive Transformation}
}

\DefLabeledTheorem[3]{rank of matrix}{\Product-\Cols}
{
Let $a$ be a matrix,
\ShowEq{\Product-rank a=k<m}{#1}
and $\SA T$ is \ProductType major submatrix.
Then \ColWS \ARow a{#2} is
a \SideWS linear composition of \ColsWS \ARow a{#3}
\ShowEq{rank of matrix, \Product-\Cols}
\ShowEq{rank of matrix, 1, \Product-\Cols}
\ShowEq{rank of matrix, 2, \Product-\Cols}
}

\DefProof[5]{rank of matrix}
{
If the matrix $a$ has $\gik$ \RowsNS, then assuming that \ColWS \ARow a{#2}
is a \SideWS linear combination
\EqRef{rank of matrix, 1, \Product-\Cols}
of \ColsWS \ARow a{#3} with coefficients $\Coef$,
we get system of \SideWS linear equations
\EqRef{rank of matrix, 2, \Product-\Cols}.
We assume that $\ACol x{#4}=\Coef$ are unknown variables
in the system of \SideWS linear equations
\EqRef{rank of matrix, 2, \Product-\Cols}.
According to the theorem
\refTheorem{nonsingular system of linear equations}{\SideNS-\Cols},
the system of \SideWS linear equations
\EqRef{rank of matrix, 2, \Product-\Cols}
has a unique solution
and this solution is nontrivial because all \ \, \ProductType quasideterminants are
different from $0$.

It remains to prove this statement in case when a number of \RowsRWSA of the matrix $a$
is more than $\gik$.
I get \ColWS \ARow a{#2} and \RowNWS $\ACol a{#5}$.
According to assumption, submatrix
$\SATpr$ is a \RC singular matrix
and its \ProductType quasideterminant
\DrawEq{singular matrix and quasideterminant}{\SideNS-\Cols}
According to the equality
\eqRef{j i quasideterminant =}{\Product}
the equality
\eqRef{singular matrix and quasideterminant}{\SideNS-\Cols}
has form
\ShowEq{rank of matrix, 4, \Product}
Matrix
\ShowEq{rank of matrix, 5, \Product-\Cols}
does not depend on $\gi{#5}$, Therefore, for any
\ShowEq{\rIn}
\ShowEq{rank of matrix, 6, \Product-\Cols}
To prove the equality
\EqRef{rank of matrix, 2, \Product-\Cols},
it is necessary to prove the equality
\ShowEq{rank of matrix, 9, \Product-\Cols}
From the equality
\ShowEq{rank of matrix, 7, \Product-\Cols}
it follows that
\ShowEq{rank of matrix, 8, \Product-\Cols}
The equality
\EqRef{rank of matrix, 9, \Product-\Cols}
follows from equalities
\EqRef{rank of matrix, 5, \Product-\Cols},
\EqRef{rank of matrix, 8, \Product-\Cols}.
The equality
\EqRef{rank of matrix, 2, \Product-\Cols}
follows from equalities
\EqRef{rank of matrix, 6, \Product-\Cols},
\EqRef{rank of matrix, 9, \Product-\Cols}.
}

\DefLabeledTheorem{set of eigenvectors is vector space}{\Product-\Cols}
{
Let $A$\Hyph number $b$ be
\ProductType eigenvalue
of the matrix $a$.
The set of eigen\ColsWS
of matrix $a$ corresponding to \ProductType eigenvalue $b$
is \SideWS $A$\Hyph vector space of \ColN s.
}

\DefProof{set of eigenvectors is vector space}
{
According to the definitions
\refDefinition{eigenvalue of matrix}{\Product},
\refDefinition{eigenvector of matrix}{\Product-\Cols},
coordinates of eigen\ColWS
are solution of the system of linear equations
\ShowEq{\Product-eigen\Cols=0}
with \ProductType singular matrix \ShowEq{f-bEn \SideNS}{}.
According to the theorem
\RefTheorem{Solutions homogenous system of linear equations},
the set of solutions of the system of linear equations
\EqRef{\Product-eigen\Cols=0}
is \SideWS $A$\Hyph vector space of \ColN s.
}

\DefLabeledTheorem{set of eigenvectors fg is vector space}{\Product-\Cols}
{
Let $A$\Hyph number $b$ be
\ProductType eigenvalue of the pair of matrices
\ShowEq{pair of matrices \Product-\ColN}{}
where $g$ is the \ProductType non\Hyph singular matrix.
The set of eigen\ColsWS
of the pair of matrices $(f,g)$ corresponding to \ProductType eigenvalue $b$
is \SideWS $A$\Hyph vector space of \ColN s.
}

\DefProof{set of eigenvectors fg is vector space}
{
According to the definitions
\refDefinition{eigenvalues of pair of matrices}{\Product},
\refDefinition{eigenvector of pair of matrices}{\Product-\ColN},
coordinates of eigen\ColWS
are solution of the system of linear equations
\ShowEq{\Product-eigen\Cols=0 fg}
with \ProductType singular matrix
\[f-\ShowEq{ga*g- \SideNS-\Cols}bg\]
According to the theorem
\RefTheorem{Solutions homogenous system of linear equations},
the set of solutions of the system of linear equations
\EqRef{\Product-eigen\Cols=0 fg}
is \SideWS space of \ColN s.
}

\DefRemark{simplify definitions of eigenvalue}
{
Let
\ShowEq{gij in ZA}
Then
\ShowEq{gb rc g-=b}
\ShowEq{gb cr g-=b}
From equalities
\EqRef{gb rc g-=b},
\EqRef{gb cr g-=b},
it follows that,
if the condition
\EqRef{gij in ZA},
is satisfied,
then we can simplify definitions
of eigenvalue and eigenvector of matrix
and get definitions
which correspond to definitions
in commutative algebra.
}

\DefLabeledTheorem{columns of matrix are linearly dependent}{\Product-\Cols}
{
Let matrix $a$ have $\gi{\Lbls}$ \ColN s.
If
\ShowEq{\Product-rank a=k<m}{\Lbls}
then \ColsWS of the matrix are \SideWS linearly dependent
\DrawEq[{\lambda}a]{columns of matrix are linearly dependent}{}
}

\DefProof[1]{columns of matrix are linearly dependent}
{
Let \ColWS \ARow a{#1} be a \SideWS linear composition of \Rows
\EqRef{rank of matrix, 1, \Product-\Cols}. We assume
\ShowEq{\Cols-of-matrix linearly dependent, 1}
and the rest
\ShowEq{\Cols-of-matrix linearly dependent, 2}
}

%% file: Vector.Space.2020.Stmt.Eq.tex

\def\spanb{\text{span}(\ARow{\Vector a}i,\iIg)}

\newcommand\ARow[2]{\aU {#1}{#2}}
\newcommand\ACol[2]{\aD {#1}{#2}}
\newcommand\EBase[2]{\ECol {##1}{##2}}%

\newcommand\ColRowTheorem[1]
{
\TwoColText
{
\ShowEq{\DefCol}
\ShowTheorem{#1}
}
{
\ShowEq{\DefRow}
\ShowTheorem{#1}
}
}

\newcommand\ColRowLemma[1]
{
\TwoColText
{
\ShowEq{\DefCol}
\ShowLemma{#1}
}
{
\ShowEq{\DefRow}
\ShowLemma{#1}
}
}

\newcommand\ColRowProof[1]
{
\TwoColText
{
\ShowEq{\DefCol}
\ShowProof{#1}
}
{
\ShowEq{\DefRow}
\ShowProof{#1}
}
}

\newcommand\ProveColRowTheorem[1]
{
\ColRowTheorem{#1}
\ColRowProof{#1}
}

\newcommand\ColRowRemark[1]
{
\TwoColText
{
\ShowEq{\DefCol}
\ShowRemark{#1}
}
{
\ShowEq{\DefRow}
\ShowRemark{#1}
}
}

\newcommand\ColRowText[1]
{
\TwoColText
{
\ShowEq{\DefCol}
\ShowText{#1}
}
{
\ShowEq{\DefRow}
\ShowText{#1}
}
}

\newcommand\ColRowDefinition[1]
{
\TwoColText
{
\ShowEq{\DefCol}
\ShowDefinition{#1}
}
{
\ShowEq{\DefRow}
\ShowDefinition{#1}
}
}

\newcommand\AUD[3]{\ACol{\ARow{#1}{#3}}{#2}}%

\AddEq{def row}
{%
\def\ProductTypeO{\CRo}%
\def\Lbls{m}%
\def\Cols{rows}%
\def\Coef{\pRs}%
\def\rIn{r in N-T}%
\renewcommand\ARow[2]{\aU {##1}{##2}}%
\renewcommand\ACol[2]{\aD {##1}{##2}}%
\renewcommand\EBase[2]{\ERow {##1}{##2}}%
\def\ArMatrix{\ColMatrix}
\def\AcMatrix{\RowMatrix}
\def\ColN{row}%
\def\ColNS{row}%
\def\ColNWS{row }%
\def\RowLbWS{row }%
\def\RowN{col}%
\def\Rows{cols}%
\def\RowNWS{column }%
\ShowEq{def row text}%
}

\AddEq{def col}
{%
\def\ProductTypeO{\RCo}%
\def\Lbls{n}%
\def\Cols{cols}%
\def\Coef{\tRr}%
\def\rIn{p in M-S}%
\renewcommand\ARow[2]{\aD {##1}{##2}}%
\renewcommand\ACol[2]{\aU {##1}{##2}}%
\renewcommand\EBase[2]{\ECol {##1}{##2}}%
\def\ArMatrix{\RowMatrix}
\def\AcMatrix{\ColMatrix}
\def\ColN{column}%
\def\ColNWS{column }%
\def\ColNS{col}%
\def\RowN{row}%
\def\Rows{rows}%
\def\RowNWS{row }%
\ShowEq{def col text}%
}

\AddEq{rc-rows}
{%
\ShowEq{def rc}%
\ShowEq{def row}%
}

\AddEq{rc-cols}
{%
\ShowEq{def rc}%
\ShowEq{def col}%
}

\AddEq{cr-rows}
{%
\ShowEq{def cr}%
\ShowEq{def row}%
}

\AddEq{cr-cols}
{%
\ShowEq{def cr}%
\ShowEq{def col}%
}

\AddEq{-rows}
{%
\ShowEq{def opRow}%
\ShowEq{def row}%
}

\AddEq{-cols}
{%
\ShowEq{def opCol}%
\ShowEq{def col}%
}

\AddEq[6]{prolog homomorphism of vector space(1)}
{
\ShowEq{Let be basis of vector space}{#3}iIA{#2}V{#3}{\Cols}-
\ShowEq{Let be basis of vector space}{#6}jJA{#5}V{#6}{\Cols}-
}

\AddEq[6]{prolog homomorphism of vector space(11)}
{
\ShowEq{Let be basis of algebra and C}{A_{#2}}{#2}kKD{#1}A{#2}-
\ShowEq{Let be basis of algebra and C}{A_{#5}}{#5}lLD{#4}A{#5}-
\ShowEq{Let be basis of vector space}{{V_{#3}}}iIA{#2}V{#3}{\Cols}-
\ShowEq{Let be basis of vector space}{{V_{#6}}}jJA{#5}V{#6}{\Cols}-
}

\AddEq[6]{prolog homomorphism of vector space(111)}
{
\ShowEq{Let be basis of algebra and C}{A_{#2}}{#2}kKD{#1}A{#2}-
\ShowEq{Let be basis of algebra and C}{A_{#5}}{#5}lLD{#4}A{#5}-
\ShowEq{Let be basis of vector space}{{V_{#3}}}iIA{#2}V{#3}{\Cols}-
\ShowEq{Let be basis of vector space}{{V_{#6}}}jJA{#5}V{#6}{\Cols}-
}

\DefRef{define homomorphism A module}
{
\ShowRef{Morphism of Diagram of Representations}%
\refDefinition{module over associative algebra}{\SideWS \VectorSetNS},
}

\DefRef[3]{Proof, homomorphism of vector space}
{%
\ShowRef{Morphism of Diagram of Representations}%
\refDefinition{homomorphism from A1 to A2, D algebra}{#1#2},
\refDefinition{homomorphism A module}{\SideNS(#1#2#3)},
}

\DefRef[3]{Proof, homomorphism of vector space 1}
{
\eqRef{f:V1->V2, D module \Cols}{#1#2#3\SideNS},
\eqRef{f o ea=efa i (#1#2)}{(\Cols-\SideNS)algebra},
\eqRef{f o ea=efa (#1#2)(\Cols)}{#1#2#3\SideNS}
}

\DefRef[4]{Proof, homomorphism of vector space 2}
{
\def\TempA{}%
\def\TempB{#1}%
\ifx\TempA\TempB%
\def\Parm{a}%
\else%
\def\Parm{h(a)}%
\fi%
\newline
\FrameEqRef[f{\Parm}{}]{f:V1->V2, D module \Cols}{algebra(#1#2)}
\newline
\FrameEqRef[hf{#2}{#4}a]{f o ea=efa i (#1#2)}{(\Cols)algebra}
\newline
\FrameEqRef[hf12]{f o ea=efa (#1#2)(\Cols)}{algebra}
\newline
}

\DefText{g:A1->A2, D module(111)}
{
\DrawEq[g{h(a)}{}]{f:V1->V2, D module \Cols}{111\SideNS}
\DrawEq[hg{A_1}{A_2}a]{f o ea=efa i (11)}{(\Cols-\SideNS)algebra}
\DrawEq[hg{A_1}{A_2}]{f o ea=efa (11)(\Cols)}{111\SideNS}
}

\DefText{g:A1->A2, D module(11)}
{
\DrawEq[ga]{f:V1->V2, D module \Cols}{11\SideNS}
\DrawEq[hg{A_1}{A_2}a]{f o ea=efa i (1)}{(\Cols-\SideNS)algebra}
\DrawEq[hg{A_1}{A_2}]{f o ea=efa (1)(\Cols)}{11\SideNS}
}

\DefText{g:A1->A2, D module(1)}
{
}

\DefText{Vector Space Type}
{
\,

\ShowDefinition{module type}

\ProveTheorem{linear combination in module type}

\ProveTheorem{coordinate matrix of vector}

\ShowTheorem{coordinates of vector}
\ShowProof{coordinates of vector *}

\ShowTheorem{Two bases of module}
\ShowLemma{Two bases of module 1}
\ShowProof{Two bases of module 1}
\ShowLemma{Two bases of module 2}
\ShowProof{Two bases of module 2}
\ShowProof{Two bases of module}
}

\AddEq[6]{extended basis}
{
\eV[#1_{#2}#3_{#4}]=(\EBase {#1_{#2}#3_{#4}}{#5#6}
=\Multiply{\EBase{#3_{#4}}#6}{\EBase {#1_{#2}}#5})
}

\AddEq{Two bases of module 1 a}
{
a=\Multiply{\ACol ai}{\EBase 1i}
}

\AddEq{Two bases of module 1 ai}
{
\ACol ai=\Multiply{\ACol a{ik}}{\EBase{}k}
}

\AddEq{Two bases of module 1 a=ai}
{
a=\Multiply{\Multiply{\ACol a{ik}}{\EBase{}k}}{\EBase 1i}
=\Multiply{\ACol a{ik}}{\EBase 2{ik}}
}

\AddEq{extended basis 1}
{
\Multiply{\ACol a{ik}}{\EBase 2{ik}}=0
}

\AddEq{extended basis 2}
{
\Multiply{\ACol a{ik}}{\Multiply{\EBase{}k}{\EBase 1i}}=0
}

\AddEq{extended basis 3}
{
\Multiply{\ACol a{ik}}{\EBase{}k}=0\ \ \ \jJg kK
}

\AddEq{extended basis 4}
{
\ACol a{ik}=0\ \ \ \jJg kK\ \ \ \jJg iI
}

\AddEq{linear span, 0}
{
\[
\Vector a=\ArMatrix{\Vector a}m
=(\ARow{\Vector a}i,\iIg)
\]
}

\AddEq{linear span, 6}
{
$\Vector b$, $\ARow{\Vector a}i$
}

\AddEq{e*v=e*w,}
{
v\ProductVal e=w\ProductVal e\ \ \ \ \ACol vi\ARow ei=\ACol wi\ARow ei
}

\AddEq{e*v=e*w,left}
{
v\ProductVal e=w\ProductVal e\ \ \ \ \ACol vi\ARow ei=\ACol wi\ARow ei
}

\AddEq{e*v=e*w,right}
{
e\ProductVal v=e\ProductVal w\ \ \ \ \ARow ei\ACol vi=\ARow ei\ACol wi
}

\AddEq{matrix a = set ai cols}
{
a=
\begin{pmatrix}
\AcMatrix{\ARow a1}n&...&\AcMatrix{\ARow am}n
\end{pmatrix}
=\PMatrix amn
}

\AddEq{matrix a = set ai rows}
{
a=
\begin{pmatrix}
\AcMatrix{\ARow a1}n\\...\\ \AcMatrix{\ARow am}n
\end{pmatrix}
=\PMatrix anm
}

\AddEq{=> v=w}
{
\[
v=w\ \ \ \ \ACol vi=\ACol wi
\]
}

\AddEq[1]{linear span, b left}
{
\Vector{#1}=#1\ProductVal e=\ACol{#1}j\ARow ej
\ \ \ \ \ #1=\AcMatrix {#1}n
}

\AddEq[1]{linear span, b right}
{
\Vector{#1}=e\ProductVal #1=\ARow ej\ACol{#1}j
\ \ \ \ \ #1=\AcMatrix {#1}n
}

\AddEq{e*b=e*a*x left}
{
\ACol bj\ARow ej=\ACol xi(\ACol{\ARow ai}j\ARow ej)
=(\ACol xi\ACol{\ARow ai}j)\ARow ej
}

\AddEq{e*b=e*a*x right}
{
\ARow ej\ACol bj=(\ARow ej\ACol{\ARow ai}j)\ACol xi
=\ARow ej(\ACol{\ARow ai}j\ACol xi)
}

\AddEq{b=a*x left}
{
\ACol bj=\ACol xi\ACol{\ARow ai}j
}

\AddEq{b=a*x right}
{
\ACol bj=\ACol{\ARow ai}j\ACol xi
}

\AddEq{a*x=b left}
{
\begin{split}
\ACol x1\ACol{\ARow a1}1+...+\ACol xm\ACol{\ARow am}1&=\ACol b1
\\...&...\\
\ACol x1\ACol{\ARow a1}n+...+\ACol xm\ACol{\ARow am}n&=\ACol bn
\end{split}
}

\AddEq{a*x=b right}
{
\begin{split}
\ACol{\ARow a1}1\ACol x1+...+\ACol{\ARow am}1\ACol xm&=\ACol b1
\\...&...\\
\ACol{\ARow a1}n\ACol x1+...+\ACol{\ARow am}n\ACol xm&=\ACol bn
\end{split}
}

\DefText{linear span in vector space}
{
\,

\ShowDefinition{linear span, vector space}

\ProveTheorem{vector in linear span}

\ProveTheorem{linear span is vector space}

\ProveTheorem{linear span and system of equations}

\ProveTheorem{matrix and system of linear equations}

\ShowDefinition{nonsingular system of linear equations}

\ProveTheorem{nonsingular system of linear equations}
}

\AddEq{x=a-*b, matrix right}
{
x=a^{\InverseVal}\ProductVal b
}

\AddEq{x=a-*b, matrix left}
{
x=b\ProductVal a^{\InverseVal}
}

\AddEq{x=a-*b, quasideterminant right}
{
x=\mathcal H\DetVal a\ProductVal b
}

\AddEq{x=a-*b, quasideterminant left}
{
x=b\ProductVal \mathcal H\DetVal a
}

\AddEquation{a*x=b left-cols}
{
\ColMatrix xm\CRstar\PMatrix amn=\ColMatrix bn
}

\AddEquation{a*x=b right-cols}
{
\PMatrix amn\RCstar\ColMatrix xm=\ColMatrix bn
}

\AddEquation{a*x=b left-rows}
{
\RowMatrix xm\RCstar\PMatrix anm=\ColMatrix bn
}

\AddEquation{a*x=b right-rows}
{
\PMatrix anm\CRstar\RowMatrix xm=\RowMatrix bn
}

\AddEq{a*x=b 1 left-cols}
{
x\CRstar a=b
}

\AddEq{a*x=b 1 right-cols}
{
a\RCstar x=b
}

\AddEq{a*x=b 1 left-rows}
{
x\RCstar a=b
}

\AddEq{a*x=b 1 right-rows}
{
a\CRstar x=b
}

\AddEq{linear span, 1}
{
\[
\Vector b\in\spanb
\]
}

\AddEq{linear span, 3 right}
{
\begin{align*}
\Vector b+\Vector c
=\Vector a\ProductVal b+\Vector a\ProductVal c
=\Vector a\ProductVal(b+c)
&\in\spanb\\
\Vector bk=(\Vector a\ProductVal b)k=\Vector a\ProductVal(bk)
&\in\spanb
\end{align*}
}

\AddEq{linear span, 3 left}
{
\begin{align*}
\Vector b+\Vector c
=a\ProductVal b+a\ProductVal c
=a\ProductVal(b+c)
&\in\spanb\\
\Vector bk=(a\ProductVal b)k=a\ProductVal(bk)
&\in\spanb
\end{align*}
}

\AddEq[2]{linear span, 2 right}
{
\Vector{#1}=\Vector a\ProductVal #2=\ARow{\Vector a}i\ACol{#2}i
}

\AddEq[2]{linear span, 2 left}
{
\Vector{#1}=#2\ProductVal\Vector a=\ACol{#2}i\ARow{\Vector a}i
}

\AddEq{linear span, vector space, 2}
{
\begin{align*}
\Vector b&=a\RCstar b\\
\Vector c&=a\RCstar c
\end{align*}
}

\AddEq[1]{b not in ZA}
{
$b\not\in Z(A)$#1
}

\AddEq[2]{D diagonal matrix}
{
\[
#1=\mathrm{diag}(b(\gi 1),...,b(\gi {#2}))
\]
}

\AddEq[3]{A rc U=U rc D right}
{
a_{#1}\RCstar u_{#2}=u_{#2}\RCstar #3
}

\AddEq[3]{A rc U=U rc D left}
{
u_{#2}\RCstar a_{#1}=#3\RCstar u_{#2}
}

\AddEq[3]{A cr U=U cr D left}
{
u_{#2}\CRstar a_{#1}=#3\CRstar u_{#2}
}

\AddEq[3]{A cr U=U cr D right}
{
a_{#1}\CRstar u_{#2}=u_{#2}\CRstar #3
}

\AddEq{U-*A*U=D rc-right}
{
u_2^{\RCInverse}\RCstar a_2\RCstar u_2=a_1
}

\AddEq{U-*A*U=D rc-left}
{
u_2\RCstar a_2\RCstar u_2^{\RCInverse}=a_1
}

\AddEq{U-*A*U=D cr-left}
{
u_2\CRstar a_2\CRstar u_2^{\CRInverse}=a_1
}

\AddEq{U-*A*U=D cr-right}
{
u_2^{\CRInverse}\CRstar a_2\CRstar u_2=a_1
}

\AddEq{U-*A*U-di En rc-right}
{
\begin{align*}
&\,u_2^{\RCInverse}\RCstar a_2\RCstar u_2-b(\gii)\aD En
\\ = &\, 
u_2^{\RCInverse}\RCstar( a_2-u_2b(\gii)\RCstar u_2^{\RCInverse})\RCstar u_2
\end{align*}
}

\AddEq{U-*A*U-di En rc-left}
{
\begin{align*}
&\,u_2\RCstar a_2\RCstar u_2^{\RCInverse}-b(\gii)\aD En
\\ = &\, 
u_2\RCstar( a_2-u_2^{\RCInverse}\RCstar b(\gii)u_2)\RCstar u_2^{\RCInverse}
\end{align*}
}

\AddEq{U-*A*U-di En cr-left}
{
\begin{align*}
&\,u_2\CRstar a_2\CRstar u_2^{\CRInverse}-b(\gii)\aD En
\\ = &\, 
u_2\CRstar( a_2-u_2^{\CRInverse}\CRstar b(\gii)u_2)\CRstar u_2^{\CRInverse}
\end{align*}
}

\AddEq{U-*A*U-di En cr-right}
{
\begin{align*}
&\,u_2^{\CRInverse}\CRstar a_2\CRstar u_2-b(\gii)\aD En
\\ = &\, 
u_2^{\CRInverse}\CRstar( a_2-u_2b(\gii)\CRstar u_2^{\CRInverse})\CRstar u_2
\end{align*}
}

\AddEquation{A-U-*di*U rc-right}
{
a_2-u_2b(\gii)\RCstar u_2^{\RCInverse}
}

\AddEquation{A-U-*di*U rc-left}
{
a_2-u_2^{\RCInverse}\RCstar b(\gii)u_2
}

\AddEquation{A-U-*di*U cr-left}
{
a_2-u_2^{\CRInverse}\CRstar b(\gii)u_2
}

\AddEquation{A-U-*di*U cr-right}
{
a_2-u_2b(\gii)\CRstar u_2^{\CRInverse}
}

\AddEquation{a1*e1=di e1 rc-left}
{
\ShowEq{ARow u2i \Cols}e1i
\RCstar a_1
=\ShowEq{ARow u2i \Cols}e1i
b(\gii)
}

\AddEq[3]{linear combination, -rows}
{
\[
\sum_{\gii=\gi 1}^{\gi{#3}}#1(\gii)#2(\gii)=\aD{#1}i\aU{#2}i
\]
}

\AddEq[3]{linear combination, left-rows}
{
\[
\sum_{\gii=\gi 1}^{\gi{#3}}#1(\gii)#2(\gii)=\aD{#1}i\aU{#2}i
\]
}

\AddEq[3]{linear combination, right-cols}
{
\[
\sum_{\gii=\gi 1}^{\gi{#3}}#2(\gii)#1(\gii)=\aD{#2}i\aU{#1}i
\]
}

\AddEq[3]{linear combination, left-cols}
{
\[
\sum_{\gii=\gi 1}^{\gi{#3}}#1(\gii)#2(\gii)=\aU{#1}i\aD{#2}i
\]
}

\AddEq[3]{linear combination, right-rows}
{
\[
\sum_{\gii=\gi 1}^{\gi{#3}}#2(\gii)#1(\gii)=\aU{#2}i\aD{#1}i
\]
}

\AddEq[3]{linear combination, -cols}
{
\[
\sum_{\gii=\gi 1}^{\gi{#3}}#1(\gii)#2(\gii)=\aU{#1}i\aD{#2}i
\]
}

\AddEquation{a1*e1=di e1 rc-right}
{
a_1\RCstar
\ShowEq{ARow u2i \Cols}e1i
=b(\gii)
\ShowEq{ARow u2i \Cols}e1i
}

\AddEquation{a1*e1=di e1 cr-left}
{
\ShowEq{ARow u2i \Cols}e1i
\CRstar a_1
=\ShowEq{ARow u2i \Cols}e1i
b(\gii)
}

\AddEquation{a1*e1=di e1 cr-right}
{
a_1\CRstar
\ShowEq{ARow u2i \Cols}e1i
=b(\gii)
\ShowEq{ARow u2i \Cols}e1i
}

\AddEquation{U-*A*U*e1=di e1 rc-left}
{
\ShowEq{ARow u2i \Cols}e1i
\RCstar u_2\RCstar a_2\RCstar u_2^{\RCInverse}
=\ShowEq{ARow u2i \Cols}e1i
b(\gii)
}

\AddEquation{U-*A*U*e1=di e1 rc-right}
{
u_2^{\RCInverse}\RCstar a_2\RCstar u_2\RCstar
\ShowEq{ARow u2i \Cols}e1i
=b(\gii)
\ShowEq{ARow u2i \Cols}e1i
}

\AddEquation{U-*A*U*e1=di e1 cr-left}
{
\ShowEq{ARow u2i \Cols}e1i
\CRstar u_2\CRstar a_2\CRstar u_2^{\CRInverse}
=\ShowEq{ARow u2i \Cols}e1i
b(\gii)
}

\AddEquation{U-*A*U*e1=di e1 cr-right}
{
u_2^{\CRInverse}\CRstar a_2\CRstar u_2\CRstar
\ShowEq{ARow u2i \Cols}e1i
=b(\gii)
\ShowEq{ARow u2i \Cols}e1i
}

\AddEquation{A*U*e1=... rc-left}
{
\begin{aligned}
&\,\ShowEq{ARow u2i \Cols}e1i
\RCstar u_2\RCstar a_2
\\ =&\,\ShowEq{ARow u2i \Cols}e1i
b(\gii)\RCstar u_2
\\ =&\,\ShowEq{ARow u2i \Cols}e1i
\RCstar u_2\RCstar u_2^{\RCInverse}\RCstar
b(\gii)u_2
\end{aligned}
}

\AddEquation{A*U*e1=... rc-right}
{
\begin{aligned}
&\,a_2\RCstar u_2\RCstar
\ShowEq{ARow u2i \Cols}e1i
\\ =&\,u_2\RCstar b(\gii)
\ShowEq{ARow u2i \Cols}e1i
\\ =&\,u_2b(\gii)\RCstar 
u_2^{\RCInverse}\RCstar u_2\RCstar
\ShowEq{ARow u2i \Cols}e1i
\end{aligned}
}

\AddEquation{A*U*e1=... cr-left}
{
\begin{aligned}
&\,\ShowEq{ARow u2i \Cols}e1i
\CRstar u_2\CRstar a_2
\\ =&\,\ShowEq{ARow u2i \Cols}e1i
b(\gii)\CRstar u_2
\\ =&\,\ShowEq{ARow u2i \Cols}e1i
\CRstar u_2\CRstar u_2^{\CRInverse}\CRstar
b(\gii)u_2
\end{aligned}
}

\AddEquation{A*U*e1=... cr-right}
{
\begin{aligned}
&\,a_2\CRstar u_2\CRstar
\ShowEq{ARow u2i \Cols}e1i
\\ =&\,u_2\CRstar b(\gii)
\ShowEq{ARow u2i \Cols}e1i
\\ =&\,u_2b(\gii)\CRstar 
u_2^{\CRInverse}\CRstar u_2\CRstar
\ShowEq{ARow u2i \Cols}e1i
\end{aligned}
}

\AddEq{expansion relative basis, vector space, 0}
{
$\Vector v$, $\aD ei$, \iIg,
}

\AddEq{expansion relative basis, vector space, (0)}
{
$\Vector v$, $e(\gii)$, \iIg,
}

\AddEq{expansion relative basis, vector space, 1}
{
b\Vector v+c\CRstar e=0
}

\AddEquation{expansion relative basis, vector space, 1 left}
{
b\Vector v+c(\gii) e(\gii)=0
}

\AddEquation{expansion relative basis, vector space, 1 right}
{
\Vector vb+e(\gii)c(\gii) =0
}

\AddEq{expansion relative basis, vector space, 2}
{
\Vector v=(-cb^{-1})\CRstar e
}

\AddEquation{expansion relative basis, vector space, 2 left}
{
\Vector v=(-b^{-1}c(\gii)) e(\gii)
}

\AddEquation{expansion relative basis, vector space, 2 right}
{
\Vector v= e(\gii)(-c(\gii)b^{-1})
}

\AddEq{expansion relative basis, vector space, 3}
{
\Vector v=v'\RCstar e
}

\AddEquation{expansion relative basis, vector space, 3 left}
{
\Vector v=v'(\gii) e(\gii)
}

\AddEquation{expansion relative basis, vector space, 3 right}
{
\Vector v=e(\gii)v'(\gii)
}

\AddEq{expansion relative basis, vector space, 4}
{
0=(v'-v)\CRstar e
}

\AddEquation{expansion relative basis, vector space, 4 left}
{
0=(v'(\gii)-v(\gii)) e(\gii)
}

\AddEquation{expansion relative basis, vector space, 4 right}
{
0= e(\gii)(v'(\gii)-v(\gii))
}

\AddEq{expansion relative basis, vector space, 5}
{
v'-v=0
}

\AddEq{expansion relative basis, vector space, 5()}
{
v'(\gii)-v(\gii)=0\ \ \ \,\iIg
}

\AddEq{expansion relative basis, vector space}
{
\Vector v=v\CRstar e
}

\AddEquation{expansion relative basis, vector space left}
{
\Vector v=v(\gii) e(\gii)\ \ \ \,v(\gii)=-b^{-1}c(\gii)
}

\AddEquation{expansion relative basis, vector space right}
{
\Vector v=e(\gii)v(\gii)\ \ \ \,v(\gii)=-c(\gii)b^{-1}
}

\AddEq{linear span, set}
{
$\{\ARow {\Vector a}i\in V, \iIg\}$
}

\AddEq{linear span, vector space}
{
\symb{\spanb}{linear span, vector space}1
}

\AddEq{e2=e1}
{
$\Basis e_2=\Basis e_1$.
}

\AddEquation{di*ui=ui*di rc-left}
{
\ShowEq{ARow u2i \Cols}u1i
\RCstar u_2^{\RCInverse}\RCstar
b(\gii)u_2=b(\gii)
\ShowEq{ARow u2i \Cols}u1i
}

\AddEquation{di*ui=ui*di rc-right}
{
u_2b(\gii)\RCstar 
u_2^{\RCInverse}\RCstar
\ShowEq{ARow u2i \Cols}u1i
=
\ShowEq{ARow u2i \Cols}u1i
b(\gii)
}

\AddEquation{di*ui=ui*di cr-left}
{
\ShowEq{ARow u2i \Cols}u1i
\CRstar u_2^{\CRInverse}\CRstar
b(\gii)u_2=b(\gii)
\ShowEq{ARow u2i \Cols}u1i
}

\AddEquation{di*ui=ui*di cr-right}
{
u_2b(\gii)\CRstar 
u_2^{\CRInverse}\CRstar
\ShowEq{ARow u2i \Cols}u1i
=
\ShowEq{ARow u2i \Cols}u1i
b(\gii)
}

\AddEquation{di*ui=ui*di 1 rc-left}
{
\begin{aligned}
&\,\ShowEq{ARow u2i \Cols}u1i
\RCstar u_2^{\RCInverse}
b(\gii)\\ =&\,b(\gii)
\ShowEq{ARow u2i \Cols}u1i
\RCstar u_2^{\RCInverse}
\end{aligned}
}

\AddEquation{di*ui=ui*di 1 rc-right}
{
\begin{aligned}
&\,b(\gii) 
u_2^{\RCInverse}\RCstar
\ShowEq{ARow u2i \Cols}u1i
\\ =&\,
u_2^{\RCInverse}\RCstar
\ShowEq{ARow u2i \Cols}u1i
b(\gii)
\end{aligned}
}

\AddEquation{di*ui=ui*di 1 cr-left}
{
\begin{aligned}
&\,\ShowEq{ARow u2i \Cols}u1i
\CRstar u_2^{\CRInverse}
b(\gii)\\ =&\,b(\gii)
\ShowEq{ARow u2i \Cols}u1i
\CRstar u_2^{\CRInverse}
\end{aligned}
}

\AddEquation{di*ui=ui*di 1 cr-right}
{
\begin{aligned}
&\,b(\gii) 
u_2^{\CRInverse}\CRstar
\ShowEq{ARow u2i \Cols}u1i
\\ =&\,
u_2^{\CRInverse}\CRstar
\ShowEq{ARow u2i \Cols}u1i
b(\gii)
\end{aligned}
}

\AddEq[3]{di 1n}
{
$#1(\gi 1)$, ..., $#1(\gi{#2})$#3
}

\AddEq[3]{di 1n+1}
{
$#1(\gi 1)$, ..., $#1(\gi{#2}-\gi 1)$, $#1(\gi{#2})$#3
}

\AddEq{di ne dj}
{
\[\gii\ne\gij\Rightarrow b(\gii)\ne b(\gij)\]
}

\AddEq{ref 1 for n eigenvectors}
{
\ePrints{348311191}%
\ifx\Semafor\ValueOn%
\refTheorem[\RefDiffEq]{covariance of eigenvalue}{\SideNS-cols},
\refTheorem[\RefDiffEq]{covariance of eigenvalue}{\SideNS-rows},
\else%
\refTheorem{covariance of eigenvalue}{\SideNS-cols},
\refTheorem{covariance of eigenvalue}{\SideNS-rows},
\fi%
}

\AddEq{ref 2 for n eigenvectors}
{
\ePrints{348311191}%
\ifx\Semafor\ValueOn%
\refTheorem[\RefDiffEq]{define homomorphism of vector space}{\SideNS}.
\else%
\refTheorem{define homomorphism A module}{\SideNS(1)}.
\fi%
}

\DefText{Left and Right Eigenvalues}
{
\AddIndex{}{eigenvalue}
\TwoColText
{
\ShowEq{def left}
\ShowEq{\LeftText}
\ShowDefinition{Left and Right Eigenvalues}
}
{
\ShowEq{def right}
\ShowEq{\RightText}
\ShowDefinition{Left and Right Eigenvalues}
}
}

\DefText[6]{homomorphism of vector space, algebra(1)}
{
\ShowText{coordinates of the linear map 1}aA{A_1}{(#1)(#2)}{}
\ShowText{coordinates of the linear map 2(#1)}ahkK{D_{#4}}
\ShowText{coordinates of the linear map 3}aA{A_2}g{}b{(#1)(#2)}
\ShowText{coordinates of the linear map 4}gA{A_{#2}}{A_{#5}}kK
\ShowText{Morphism of D algebra}{#1}{#2}g
}

\DefText[6]{homomorphism of vector space, algebra()}
{
}

\DefText[6]{homomorphism of vector space, algebra 1}
{
\ShowText{coordinates of the linear map 1}vV{V_1}{(#1)(#2)(#3)}{\SideNS}
\ShowText{coordinates of the linear map 2(#2)}vgiI{A_{#5}}
\ShowText{coordinates of the linear map 3}vV{V_2}f{\SideNS}w{(#1)(#2)}
\ShowText{coordinates of the linear map 4}fV{V_{#3}}{V_{#6}}iI
}

\AddEq[3]{bases eA eV}
{
$(\eV[A_{#1}][,]\eV[V_{#2}][)]$#3
}

\AddEq[1]{Left and Right Eigenvalues 1}
{
\ShowText{Left and Right Eigenvalues}

\ShowEq{def #1}
\ShowDefinition{spectrum of matrix}

\TwoColText
{
\ShowEq{def left}
\ShowEq{\LeftText}
\ShowTheorem{Eigenvalue and conjugate class}
}
{
\ShowEq{def right}
\ShowEq{\RightText}
\ShowTheorem{Eigenvalue and conjugate class}
}

\ShowFootnote{Eigenvalue and conjugate class}

\TwoColText
{
\ShowEq{def left}
\ShowEq{\LeftText}
\ShowProof{Eigenvalue and conjugate class}
}
{
\ShowEq{def right}
\ShowEq{\RightText}
\ShowProof{Eigenvalue and conjugate class}
}

\TwoColText
{
\ShowEq{def left}
\ShowEq{\LeftText}
\ShowTheorem{Coefficient comutes with matrix}vc
\ShowProof{Coefficient comutes with matrix}vc
}
{
\ShowEq{def right}
\ShowEq{\RightText}
\ShowTheorem{Coefficient comutes with matrix}cv
\ShowProof{Coefficient comutes with matrix}cv
}

\TwoColText
{
\ShowEq{def left}
\ShowEq{\LeftText}
\ShowRemark{Eigenvalue and conjugate class}
}
{
\ShowEq{def right}
\ShowEq{\RightText}
\ShowRemark{Eigenvalue and conjugate class}
}

\TwoColText
{
\ShowEq{def left}
\ShowEq{\LeftText}
\ShowTheorem{Eigenvector does not depend on basis}
\begin{sloppypar}
\ShowProof{Eigenvector does not depend on basis}
\end{sloppypar}
}
{
\ShowEq{def right}
\ShowEq{\RightText}
\ShowTheorem{Eigenvector does not depend on basis}
\begin{sloppypar}
\ShowProof{Eigenvector does not depend on basis}
\end{sloppypar}
}

\TwoColText
{
\ShowEq{def left}
\ShowEq{\LeftText}
\ShowRemark{Eigenvector does not depend on basis}
}
{
\ShowEq{def right}
\ShowEq{\RightText}
\ShowRemark{Eigenvector does not depend on basis}
}

\TwoColText
{
\ShowEq{def left}
\ShowEq{\LeftText}
\begin{sloppypar}
\ProveTheorem{right eigenvalue is eigenvalue}
\end{sloppypar}
}
{
\ShowEq{def right}
\ShowEq{\RightText}
\begin{sloppypar}
\ProveTheorem{right eigenvalue is eigenvalue}
\end{sloppypar}
}

\TwoColText
{
\ShowEq{def left}
\ShowEq{\LeftText}
\begin{sloppypar}
\ShowRemark{right eigenvalue is eigenvalue}
\end{sloppypar}
}
{
\ShowEq{def right}
\ShowEq{\RightText}
\begin{sloppypar}
\ShowRemark{right eigenvalue is eigenvalue}
\end{sloppypar}
}

}

\AddEq{prove n eigenvectors}
{
\ShowEq{def AVector}
\ShowTheorem{n eigenvectors}{algebra}

\ShowRemark{n eigenvectors}

\TwoColText
{
\ShowEq{def left}
\ShowProof{n eigenvectors}
}
{
\ShowEq{def right}
\ShowProof{n eigenvectors}
}
}

\AddEq{gik<gin}
{
$\gik$, $\gik<\gin$,
}

\AddEq{rank rc u2 =k}
{
\[
\RCRank u_2=\gik
\]
}

\AddEq{rank rc u2 <k}
{
\[
\RCRank u_2<\gik
\]
}

\AddEquation{w rc u2=Ek left}
{
u_2\RCstar w=\aD Ek
}

\AddEquation{w rc u2=Ek right}
{
w\RCstar u_2=\aD Ek
}

\AddEquation{w rc a2 rc u2=a1 left}
{
u_2\RCstar a_2\RCstar w=a_1
}

\AddEquation{w rc a2 rc u2=a1 right}
{
w\RCstar a_2\RCstar u_2=a_1
}

\AddEq{gin x gik left}
{
$\gin\times\gik$
}

\AddEq{gin x gik right}
{
$\gik\times\gin$
}

\AddEq{as ox as=a1 ox a2, 1}
{
\Entry[s0]a{}ij\otimes\Entry[s1]a{}ij
=1\otimes\Entry a1ij
}

\AddEq{as ox as=a1 ox a2, 0}
{
\Entry[s0]a{}ij\otimes\Entry[s1]a{}ij
=\Entry a0ij\otimes 1
}

\AddEquation{as ox as o cv = ..., left-cols}
{
\begin{split}
&\,(\Entry[s0]a{}ij\otimes\Entry[s1]a{}ij)\circ (c\aU vj)
\\=&\,\Entry[s0]a{}ijc\aU vj\Entry[s1]a{}ij
\\=&\,c(\Entry[s0]a{}ij\aU vj\Entry[s1]a{}ij)
\\=&\,cb\aU vi=b(c\aU vi)
\end{split}
}

\AddEquation{as ox as o cv = ..., right-cols}
{
\begin{split}
&\,(\Entry[s0]a{}ij\otimes\Entry[s1]a{}ij)\circ (\aU vjc)
\\=&\,\Entry[s0]a{}ij\aU vjc\Entry[s1]a{}ij
\\=&\,(\Entry[s0]a{}ij\aU vj\Entry[s1]a{}ij)c
\\=&\,\aU vibc=(\aU vic)b
\end{split}
}

\AddEquation{as ox as o cv = ..., left-rows}
{
\begin{split}
&\,(\Entry[s0]a{}ij\otimes\Entry[s1]a{}ij)\circ (c\aD vi)
\\=&\,\Entry[s0]a{}ijc\aD vi\Entry[s1]a{}ij
\\=&\,c(\Entry[s0]a{}ij\aD vi\Entry[s1]a{}ij)
\\=&\,cb\aD vj=b(c\aD vj)
\end{split}
}

\AddEquation{as ox as o cv = ..., right-rows}
{
\begin{split}
&\,(\Entry[s0]a{}ij\otimes\Entry[s1]a{}ij)\circ (\aD vic)
\\=&\,\Entry[s0]a{}ij\aD vic\Entry[s1]a{}ij
\\=&\,(\Entry[s0]a{}ij\aD vi\Entry[s1]a{}ij)c
\\=&\,\aD vjbc=(\aD vfc)b
\end{split}
}

\AddEq{as ox as o v = v a12, 1, cols}
{
(\Entry[s0]a{}ij\otimes\Entry[s1]a{}ij)\circ \aU vj
=(1\otimes\Entry a1ij)\circ \aU vj
=\aU vj\Entry a1ij
}

\AddEquation{as ox as o v = v a12, matrix, 0 left-cols}
{
\PMatrix{a_0^{}{}}nn\RCstar\ColMatrix vn=
b\ColMatrix vn
}

\AddEquation{as ox as o v = v a12, matrix, 1 left-cols}
{
\ColMatrix vn\CRstar\PMatrix{a_1^{}{}}nn=
b\ColMatrix vn
}

\AddEquation{as ox as o v = v a12, matrix, 0 right-cols}
{
\PMatrix{a_0^{}{}}nn\RCstar\ColMatrix vn=
\ColMatrix vnb
}

\AddEquation{as ox as o v = v a12, matrix, 1 right-cols}
{
\ColMatrix vn\CRstar\PMatrix{a_1^{}{}}nn=
\ColMatrix vnb
}

\AddEquation{as ox as o v = v a12, matrix, 0 right-rows}
{
\PMatrix{a_0^{}{}}nn\CRstar\RowMatrix vn=
\RowMatrix vnb
}

\AddEquation{as ox as o v = v a12, matrix, 1 right-rows}
{
\RowMatrix vn\RCstar\PMatrix{a_1^{}{}}nn=
\RowMatrix vnb
}

\AddEquation{as ox as o v = v a12, matrix, 0 left-rows}
{
\PMatrix{a_0^{}{}}nn\CRstar\RowMatrix vn=
b\RowMatrix vn
}

\AddEquation{as ox as o v = v a12, matrix, 1 left-rows}
{
\RowMatrix vn\RCstar\PMatrix{a_1^{}{}}nn=
b\RowMatrix vn
}

\AddEquation{as ox as o v = v a12, matrix, 0 1 left-cols}
{
\PMatrix{a_0^{}{}}nn\RCstar\ColMatrix vn=
\ShowEq{matrix of maps n x n}
\RCcirc\ColMatrix vn
}

\AddEquation{as ox as o v = v a12, matrix, 1 1 left-cols}
{
\ColMatrix vn\CRstar\PMatrix{a_1^{}{}}nn=
\ShowEq{matrix of maps n x n}
\RCcirc\ColMatrix vn
}

\AddEquation{as ox as o v = v a12, matrix, 0 1 right-cols}
{
\PMatrix{a_0^{}{}}nn\RCstar\ColMatrix vn=
\ShowEq{matrix of maps n x n}
\RCcirc\ColMatrix vn
}

\AddEquation{as ox as o v = v a12, matrix, 1 1 right-cols}
{
\ColMatrix vn\CRstar\PMatrix{a_1^{}{}}nn=
\ShowEq{matrix of maps n x n}
\RCcirc\ColMatrix vn
}

\AddEquation{as ox as o v = v a12, matrix, 0 1 right-rows}
{
\PMatrix{a_0^{}{}}nn\CRstar\RowMatrix vn=
\ShowEq{matrix of maps n x n}
\CRcirc\RowMatrix vn
}

\AddEquation{as ox as o v = v a12, matrix, 1 1 right-rows}
{
\RowMatrix vn\RCstar\PMatrix{a_1^{}{}}nn=
\ShowEq{matrix of maps n x n}
\CRcirc\RowMatrix vn
}

\AddEquation{as ox as o v = v a12, matrix, 0 1 left-rows}
{
\PMatrix{a_0^{}{}}nn\CRstar\RowMatrix vn=
\ShowEq{matrix of maps n x n}
\CRcirc\RowMatrix vn
}

\AddEquation{as ox as o v = v a12, matrix, 1 1 left-rows}
{
\RowMatrix vn\RCstar\PMatrix{a_1^{}{}}nn=
\ShowEq{matrix of maps n x n}
\CRcirc\RowMatrix vn
}

\AddEquation{as ox as o v = v a12, matrix, 0 2 left-cols}
{
\ShowEq{matrix of maps n x n}
\RCcirc\ColMatrix vn
= b\ColMatrix vn
}

\AddEquation{as ox as o v = v a12, matrix, 1 2 left-cols}
{
\ShowEq{matrix of maps n x n}
\RCcirc\ColMatrix vn
= b\ColMatrix vn
}

\AddEquation{as ox as o v = v a12, matrix, 0 2 right-cols}
{
\ShowEq{matrix of maps n x n}
\RCcirc\ColMatrix vn
= \ColMatrix vnb
}

\AddEquation{as ox as o v = v a12, matrix, 1 2 right-cols}
{
\ShowEq{matrix of maps n x n}
\RCcirc\ColMatrix vn
= \ColMatrix vnb
}

\AddEquation{as ox as o v = v a12, matrix, 0 2 right-rows}
{
\ShowEq{matrix of maps n x n}
\CRcirc\RowMatrix vn
= \RowMatrix vnb
}

\AddEquation{as ox as o v = v a12, matrix, 1 2 right-rows}
{
\ShowEq{matrix of maps n x n}
\CRcirc\RowMatrix vn
= \RowMatrix vnb
}

\AddEquation{as ox as o v = v a12, matrix, 0 2 left-rows}
{
\ShowEq{matrix of maps n x n}
\CRcirc\RowMatrix vn
= b\RowMatrix vn
}

\AddEquation{as ox as o v = v a12, matrix, 1 2 left-rows}
{
\ShowEq{matrix of maps n x n}
\CRcirc\RowMatrix vn
= b\RowMatrix vn
}

\AddEq{as ox as o v = v a12, 0, cols}
{
(\Entry[s0]a{}ij\otimes\Entry[s1]a{}ij)\circ \aU vj
=(\Entry a0ij\otimes 1)\circ \aU vj
=\Entry a0ij\aU vj
}

\AddEq{as ox as o v = v a12, 1, rows}
{
(\Entry[s0]a{}ij\otimes\Entry[s1]a{}ij)\circ \aD vi
=(1\otimes\Entry a1ij)\circ \aD vi
=\aD vi\Entry a1ij
}

\AddEq{as ox as o v = v a12, 0, rows}
{
(\Entry[s0]a{}ij\otimes\Entry[s1]a{}ij)\circ \aD vi
=(\Entry a0ij\otimes 1)\circ \aD vi
=\Entry a0ij\aD vi
}

\AddEq{Eigenvalue of Matrix of Linear Map}
{
\AddIndex{}{eigenvalue}
\TwoColText
{
\ShowEq{def left}
\ShowEq{\DefMatrix}
\ShowDefinition{Eigenvalue of Matrix of Linear Map}
}
{
\ShowEq{def right}
\ShowEq{\DefMatrix}
\ShowDefinition{Eigenvalue of Matrix of Linear Map}
}

\TwoColText
{
\ShowEq{def left}
\ShowEq{\DefMatrix}
\ShowTheorem{Eigenvalue of Linear Map a ox a, 1}1
}
{
\ShowEq{def right}
\ShowEq{\DefMatrix}
\ShowTheorem{Eigenvalue of Linear Map a ox a, 1}0
}

\ShowEq{def left}
\ShowEq{\DefMatrix}
\ShowProof{Eigenvalue of Linear Map a ox a, 1}1

\ShowEq{def right}
\ShowEq{\DefMatrix}
\ShowProof{Eigenvalue of Linear Map a ox a, 1}0

\TwoColText
{
\ShowEq{def left}
\ShowEq{\DefMatrix}
\ShowTheorem{Eigenvalue of Linear Map a ox a, 2}0
}
{
\ShowEq{def right}
\ShowEq{\DefMatrix}
\ShowTheorem{Eigenvalue of Linear Map a ox a, 2}1
}

\ShowEq{def left}
\ShowEq{\DefMatrix}
\ShowProof{Eigenvalue of Linear Map a ox a, 2}0

\ShowEq{def right}
\ShowEq{\DefMatrix}
\ShowProof{Eigenvalue of Linear Map a ox a, 2}1

\TwoColText
{
\ShowEq{def left}
\ShowEq{\DefMatrix}
\ShowTheorem{Eigenvalue a ox a, c in AAA[b]}0
\ShowProof{Eigenvalue a ox a, c in AAA[b]}0{cv}
}
{
\ShowEq{def right}
\ShowEq{\DefMatrix}
\ShowTheorem{Eigenvalue a ox a, c in AAA[b]}1
\ShowProof{Eigenvalue a ox a, c in AAA[b]}1{vc}
}
}

\AddEq{eigenvector cv left-cols}
{
\[
cv=\ColMatrix{cv}n
\]
}

\AddEq{eigenvector cv right-cols}
{
\[
vc=
\begin{pmatrix}
\aU v1c\\...\\ \aU vnc
\end{pmatrix}
\]
}

\AddEq{eigenvector cv left-rows}
{
\[
cv=\RowMatrix{cv}n
\]
}

\AddEq{eigenvector cv right-rows}
{
\[
vc=
\begin{pmatrix}
\aD v1c&...& \aD vnc
\end{pmatrix}
\]
}

\AddEq[1]{Eigenvalue of Linear Map left-cols}
{
a\RCcirc #1=b#1
}

\AddEq[1]{Eigenvalue of Linear Map right-cols}
{
a\RCcirc #1=#1b
}

\AddEq[1]{Eigenvalue of Linear Map left-rows}
{
a\CRcirc #1=b#1
}

\AddEq[1]{Eigenvalue of Linear Map right-rows}
{
a\CRcirc #1=#1b
}

\AddEq{a=matrix of maps n x n}
{
a=
\ShowEq{matrix of maps n x n}
}

\AddEq{matrix of maps n x n}
{
\begin{pmatrix}
\Entry[s0]a{}11\otimes\Entry[s1]a{}11&...&\Entry[s0]a{}1n\otimes\Entry[s1]a{}1n\\
...&...&...\\
\Entry[s0]a{}n1\otimes\Entry[s1]a{}n1&...&\Entry[s0]a{}nn\otimes\Entry[s1]a{}nn
\end{pmatrix}
}

\AddEq{Left and Right Eigenvalues not equal}
{
\ShowRemark{Left and Right Eigenvalues not equal 0}

\TwoColText
{
\ShowEq{def left}
\ShowEq{\DefRow}
\ShowRemark{Left and Right Eigenvalues not equal}
}
{
\ShowEq{def right}
\ShowEq{\DefCol}
\ShowRemark{Left and Right Eigenvalues not equal}
}
}

\AddEquation{ai vi=0 left}
{
\begin{aligned}
&\,a(\gi 1)v(\gi 1)+ ...\\ +&\,a(\gi m-\gi 1)v(\gi m-\gi 1)\\ +&\,a(\gim)v(\gim)=0
\end{aligned}
}

\AddEquation{ai vi=0}
{
a(\gi 1)v(\gi 1)+ ...+a(\gi m-\gi 1)v(\gi m-\gi 1)+a(\gim)v(\gim)=0
}

\AddEquation{ai vi=0 right}
{
\begin{aligned}
&\,v(\gi 1)a(\gi 1)+ ...\\ +&\,v(\gi m-\gi 1)a(\gi m-\gi 1)\\ +&\,v(\gim)a(\gim)=0
\end{aligned}
}

\AddEquation{ai f vi=0}
{
a(\gi 1)(f\circ v(\gi 1)) + ... +a(\gi m-\gi 1)(f\circ v(\gi m-\gi 1))
+a(\gim)(f\circ v(\gim))=0
}

\AddEquation{ai f vi=0 left}
{
\begin{aligned}
&\,a(\gi 1)(f\circ v(\gi 1)) + ...\\ +&\,a(\gi m-\gi 1)(f\circ v(\gi m-\gi 1))
\\ +&\,a(\gim)(f\circ v(\gim))=0
\end{aligned}
}

\AddEquation{ai f vi=0 right}
{
\begin{aligned}
&\,(f\circ v(\gi 1))a(\gi 1)+ ...\\ +&\,(f\circ v(\gi m-\gi 1))a(\gi m-\gi 1)
\\ +&\,(f\circ v(\gim))a(\gim)=0
\end{aligned}
}

\AddEquation{ai di vi=0 left}
{
\begin{aligned}
&\,a(\gi 1)v(\gi 1)b(\gi 1) + ...
\\ +&\,a(\gi m-\gi 1)\\ *&\,v(\gi m-\gi 1)b(\gi m-\gi 1)
\\ +&\,a(\gim)v(\gi m)b(\gi m) =0
\end{aligned}
}

\AddEquation{ai di vi=0}
{
\begin{aligned}
&\,a(\gi 1)b(\gi 1)v(\gi 1) + ...
\\ +&\,a(\gi m-\gi 1)b(\gi m-\gi 1)v(\gi m-\gi 1)
+a(\gim)b(\gi m)v(\gi m) =0
\end{aligned}
}

\AddEquation{ai di vi=0 right}
{
\begin{aligned}
&\,b(\gi 1)v(\gi 1)a(\gi 1)+ ...
\\ +&\,b(\gi m-\gi 1)v(\gi m-\gi 1)\\ *&\,a(\gi m-\gi 1)
\\ +&\,b(\gi m)v(\gi m)a(\gim) =0
\end{aligned}
}

\AddEq{a1 ne 0}
{
$a(\gi 1)\ne 0$.
}

\AddEq{Corollary e c-ac=}
{
\begin{\CorollaryStyle}
\labelCorollary{e c-ac=c- eat c}
\ShowEq{ce c-ac c-=eat}
\ShowEq{e c-ac=c- eat c}
\end{\CorollaryStyle}
}

\AddEquation{U-*A*U*E=E*D rc-right}
{
\begin{aligned}
&\,u_2^{\RCInverse}\RCstar a_2\RCstar u_2\RCstar\aD[1]ei
\\ =&\,\aD[1]eib(\gii)
\end{aligned}
}

\AddEquation{U-*A*U*E=E*D rc-left}
{
\aU{e_1}i\RCstar u_2\RCstar a_2\RCstar u_2^{\RCInverse}=b(\gii)\aU{e_1}i
}

\AddEquation{U-*A*U*E=E*D cr-left}
{
\aD[1]ei\CRstar u_2\CRstar a_2\CRstar u_2^{\CRInverse}=b(\gii)\aD[1]ei
}

\AddEquation{U-*A*U*E=E*D cr-right}
{
u_2^{\CRInverse}\CRstar a_2\CRstar u_2\CRstar\aU{e_1}i=\aU{e_1}ib(\gii)
}

\AddEquation{A*U*E=U*E*D rc-right}
{
a_2\RCstar u_2\RCstar\aD[1]ei=u_2\RCstar\aD[1]eib(\gii)
}

\AddEquation{A*U*E=U*E*D rc-left}
{
\aU{e_1}i\RCstar u_2\RCstar a_2=b(\gii)\aU{e_1}i\RCstar u_2
}

\AddEquation{A*U*E=U*E*D cr-left}
{
\aD[1]ei\CRstar u_2\CRstar a_2=b(\gii)\aD[1]ei\CRstar u_2
}

\AddEquation{A*U*E=U*E*D cr-right}
{
a_2\CRstar u_2\CRstar\aU{e_1}i=u_2\CRstar\aU{e_1}ib(\gii)
}

\AddEquation{U*Ei=Ui rc-right}
{
u_2\RCstar\aD[1]ei=\aD[2]ui
}

\AddEquation{U*Ei=Ui rc-left}
{
\aU{e_1}i\RCstar u_2=\aU{u_2}i
}

\AddEquation{U*Ei=Ui cr-left}
{
\aD[1]ei\CRstar u_2=\aD[2]ui
}

\AddEquation{U*Ei=Ui cr-right}
{
u_2\CRstar\aU{e_1}i=\aU{u_2}i
}

\AddEquation{e1i=u2i*e2 left-cols}
{
\ShowEq{ARow u2i \Cols}e1i
=
\ShowEq{ARow u2i \Cols}u2i
\CRstar\Basis e_2
}

\AddEquation{e1i=u2i*e2 right-cols}
{
\ShowEq{ARow u2i \Cols}e1i
=
\Basis e_2\RCstar
\ShowEq{ARow u2i \Cols}u2i
}

\AddEquation{e1i=u2i*e2 left-rows}
{
\ShowEq{ARow u2i \Cols}e1i
=
\ShowEq{ARow u2i \Cols}u2i
\RCstar\Basis e_2
}

\AddEquation{e1i=u2i*e2 right-rows}
{
\ShowEq{ARow u2i \Cols}e1i
=
\Basis e_2\CRstar
\ShowEq{ARow u2i \Cols}u2i
}

\AddEq{a1-di En}
{
a_1-b(\gii)\aD En
}

\AddEq{e1ij=01 cols}
{
\Entry e1ji=\aUD{\delta}ji
}

\AddEq{e1ij=01 rows}
{
\Entry e1ij=\aUD{\delta}ij
}

\AddEq[2]{A*Ui=Ui di rc-right}
{
a_{#1}\RCstar\aD[{#1}{}]{#2}i=\aD[{#1}{}]{#2}ib(\gii)
}

\AddEq[2]{A*Ui=Ui di rc-left}
{
\aU{#2_{#1}}i\RCstar a_{#1}=b(\gii)\aU{#2_{#1}}i
}

\AddEq[2]{A*Ui=Ui di cr-right}
{
a_{#1}\CRstar\aU{#2_{#1}}i=\aU{#2_{#1}}ib(\gii)
}

\AddEq[2]{A*Ui=Ui di cr-left}
{
\aD[{#1}{}]{#2}i\CRstar a_{#1}=b(\gii)\aD[{#1}{}]{#2}i
}

\AddEq{symb coordinates of geometric object}
{
\symb{\mathcal O(V,W,\Basis e_V,w)}{coordinates of geometric object}{}
}

\AddEq{conjugate eigenvalue left-cols}
{
\[
cv\CRstar a_2=cb v=cb c^{-1}\,cv
\]
}

\AddEq{conjugate eigenvalue right-cols}
{
\[
a_2\RCstar vc=vb c=vc\,c^{-1}b c
\]
}

\AddEq{conjugate eigenvalue left-rows}
{
\[
cv\RCstar a_2=cb v=cb c^{-1}\,cv
\]
}

\AddEq{conjugate eigenvalue right-rows}
{
\[
a_2\CRstar vc=vb c=vc\,c^{-1}b c
\]
}

\AddEquation{Coefficient comutes with matrix left-cols}
{
b\aU vic=\aU vja_2^{}\aUD{}ijc=\aU vjca_2^{}\aUD{}ij
}

\AddEquation{Coefficient comutes with matrix right-cols}
{
c\aU vib=ca_2^{}\aUD{}ij\aU vj=\aUD aijc\aU vj
}

\AddEquation{Coefficient comutes with matrix left-rows}
{
b\aD vic=\aD vja_2^{}\aUD{}jic=\aD vjca_2^{}\aUD{}ji
}

\AddEquation{Coefficient comutes with matrix right-rows}
{
c\aD vib=ca_2^{}\aUD{}ji\aD vj=a_2^{}\aUD{}jic\aD vj
}

\AddEq{Vector A cols}
{
$\aD ai$, \iIg,
}

\AddEq{Vector A rows}
{
$\aU ai$, \iIg,
}

\AddEq[1]{AoxA linearly independent 1 cols}
{
$\aU bi\aU ci=0$, $\gii=\gi 1$, ..., $\gin$#1
}

\AddEq{AoxA linearly independent cols}
{
\[\aU bi\aD ai\aU ci=0\]
}

\AddEq[1]{AoxA linearly independent 1 rows}
{
$\aD bi\aD ci=0$, $\gii=\gi 1$, ..., $\gin$#1
}

\AddEq{AoxA linearly independent rows}
{
\[\aD bi\aU ai\aD ci=0\]
}

\AddEq{basis of AoxA module}
{
\symb{\Basis{e}}{Basis}{}
}

\AddEq{basis, AoxA module cols}
{
\[
\ShowSymbol{Basis}{}
=(\aD ei,\iIg)
\]
}

\AddEq{basis, AoxA module rows}
{
\[
\ShowSymbol{Basis}{}
=(\aU ei,\iIg)
\]
}

\AddEquation{def coordinates of geometric object, left-cols}
{
\begin{aligned}
&\,\ShowSymbol{coordinates of geometric object}{}
\\=&\,(w\CRstar F(G)^{\CRInverse},G\CRstar \eV[V])
\end{aligned}
}

\AddEquation{def coordinates of geometric object, right-cols}
{
\begin{aligned}
&\,\ShowSymbol{coordinates of geometric object}{}
\\=&\,(F(G)^{\RCInverse}\RCstar w,\eV[V]\RCstar G)
\end{aligned}
}

\AddEquation{def coordinates of geometric object, left-rows}
{
\begin{aligned}
&\,\ShowSymbol{coordinates of geometric object}{}
\\=&\,(w\RCstar F(G)^{\RCInverse},G\RCstar \eV[V])
\end{aligned}
}

\AddEquation{def coordinates of geometric object, right-rows}
{
\begin{aligned}
&\,\ShowSymbol{coordinates of geometric object}{}
\\=&\,(F(G)^{\CRInverse}\CRstar w,\eV[V]\CRstar G)
\end{aligned}
}

\AddEq[4]{geometric object representative, left-cols}
{
\Vector #1_{#2}=#3\CRstar e_{#4}
}

\AddEq[4]{geometric object representative, right-cols}
{
\Vector #1_{#2}=e_{#4}\RCstar #3
}

\AddEq[4]{geometric object representative, left-rows}
{
\Vector #1_{#2}=#3\RCstar e_{#4}
}

\AddEq[4]{geometric object representative, right-rows}
{
\Vector #1_{#2}=e_{#4}\CRstar #3
}

\AddEq{invariance principle 3, left-cols}
{
\[
\begin{aligned}
\Vector w'&=w'\CRstar e'_W\\
&=w\CRstar F(a)^{\CRInverse}\CRstar F(a)\CRstar e_W\\
&=w\CRstar e_W=\Vector w
\end{aligned}
\]
}

\AddEq{invariance principle 3, right-cols}
{
\[
\begin{aligned}
\Vector w'&=e'_W\RCstar w'\\
&=e_W\RCstar F(a)\RCstar F(a)^{\RCInverse}\RCstar w\\
&=e_W\RCstar w=\Vector w
\end{aligned}
\]
}

\AddEq{invariance principle 3, left-rows}
{
\[
\begin{aligned}
\Vector w'&=w'\RCstar e'_W\\
&=w\RCstar F(a)^{\RCInverse}\RCstar F(a)\RCstar e_W\\
&=w\RCstar e_W=\Vector w
\end{aligned}
\]
}

\AddEq{invariance principle 3, right-rows}
{
\[
\begin{aligned}
\Vector w'&=e'_W\CRstar w'\\
&=e_W\CRstar F(a)\CRstar F(a)^{\CRInverse}\CRstar w\\
&=e_W\CRstar w=\Vector w
\end{aligned}
\]
}

\AddEquation{coordinate matrix, W2, W1, left-cols}
{
\Basis e_{W2}=c^{\CRInverse}\CRstar \Basis e_{W1}
}

\AddEquation{coordinate matrix, W2, W1, right-cols}
{
\Basis e_{W2}=e_{W1}\RCstar \Basis c^{\RCInverse}
}

\AddEquation{coordinate matrix, W2, W1, left-rows}
{
\Basis e_{W2}=c^{\RCInverse}\RCstar \Basis e_{W1}
}

\AddEquation{coordinate matrix, W2, W1, right-rows}
{
\Basis e_{W2}=e_{W1}\CRstar \Basis c^{\CRInverse}
}

\AddEq{sum of geometric objects, 1}
{
\DrawEq[w1{w_1}W]{geometric object representative, \SideNS-\Cols}{}
\DrawEq[w2{w_2}W]{geometric object representative, \SideNS-\Cols}{}
}

\AddEq{sum of geometric objects, 2}
{
\DrawEq[w{}{(w_1+w_2)}W]{geometric object representative, \SideNS-\Cols}{}
}

\AddEq{product of geometric object and constant, 2, left}
{
\DrawEq[w2{(kw_1)}W]{geometric object representative, \SideNS-\Cols}{}
}

\AddEq{product of geometric object and constant, 2, right}
{
\DrawEq[w2{(w_1k)}W]{geometric object representative, \SideNS-\Cols}{}
}

\AddEq{product of geometric object and constant, 3, left}
{
\[\Vector w_2=k\Vector w_1\]
}

\AddEq{product of geometric object and constant, 3, right}
{
\[\Vector w_2=\Vector w_1k\]
}

\AddEq{sum of geometric objects, 3}
{
\[\Vector w=\Vector w_1+\Vector w_2\]
}

\AddEquation{representation of homomorphism relative different bases, 3, left-cols}
{
\Vector w=v\CRstar b\CRstar a_2\CRstar c^{\CRInverse}\CRstar e_{W1}
}

\AddEquation{representation of homomorphism relative different bases, 4, left-cols}
{
v\CRstar a_1=v\CRstar b\CRstar a_2\CRstar c^{\CRInverse}
}

\AddEquation{representation of homomorphism relative different bases, 3, right-cols}
{
\Vector w=e_{W1}\RCstar c^{\RCInverse}\RCstar a_2\RCstar b\RCstar v
}

\AddEquation{representation of homomorphism relative different bases, 4, right-cols}
{
v\CRstar a_1=c^{\RCInverse}\RCstar a_2\RCstar b\RCstar v
}

\AddEquation{representation of homomorphism relative different bases, 3, left-rows}
{
\Vector w=v\RCstar b\RCstar a_2\RCstar c^{\RCInverse}\RCstar e_{W1}
}

\AddEquation{representation of homomorphism relative different bases, 4, left-rows}
{
v\CRstar a_1=v\RCstar b\RCstar a_2\RCstar c^{\RCInverse}
}

\AddEquation{representation of homomorphism relative different bases, 3, right-rows}
{
\Vector w=e_{W1}\CRstar c^{\CRInverse}\CRstar a_2\CRstar b\CRstar v
}

\AddEquation{representation of homomorphism relative different bases, 4, right-rows}
{
v\CRstar a_1=c^{\CRInverse}\CRstar a_2\CRstar b\CRstar v
}

\AddEquation{representation of homomorphism relative different bases, 1, left-cols}
{
\Vector w=v\CRstar a_1\CRstar e_{W1}
}

\AddEquation{representation of homomorphism relative different bases, 2, left-cols}
{
\Vector w=v\RCstar b\CRstar a_2\CRstar e_{W2}
}

\AddEquation{representation of homomorphism relative different bases, 1, right-cols}
{
\Vector w=e_{W1}\RCstar a_1\RCstar v
}

\AddEquation{representation of homomorphism relative different bases, 2, right-cols}
{
\Vector w=e_{W2}\RCstar a_2\RCstar b\RCstar v
}

\AddEquation{representation of homomorphism relative different bases, 1, left-rows}
{
\Vector w=v\RCstar a_1\RCstar e_{W1}
}

\AddEquation{representation of homomorphism relative different bases, 2, left-rows}
{
\Vector w=v\CRstar b\RCstar a_2\RCstar e_{W2}
}

\AddEquation{representation of homomorphism relative different bases, 1, right-rows}
{
\Vector w=e_{W1}\CRstar a_1\CRstar v
}

\AddEquation{representation of homomorphism relative different bases, 2, right-rows}
{
\Vector w=e_{W2}\RCstar a_2\CRstar b\CRstar v
}

\AddEq{expansion of vector v, left-cols}
{
\[
\Vector v=v\CRstar e_{V1}=v\CRstar b\CRstar e_{V2}
\]
}

\AddEq{expansion of vector v, right-cols}
{
\[
\Vector v=e_{V1}\CRstar v=e_{V2}\CRstar b\CRstar v
\]
}

\AddEq{expansion of vector v, left-rows}
{
\[
\Vector v=v\RCstar e_{V1}=v\RCstar b\RCstar e_{V2}
\]
}

\AddEq{expansion of vector v, right-rows}
{
\[
\Vector v=e_{V1}\RCstar v=e_{V2}\RCstar b\RCstar v
\]
}

\AddEq[1]{representation of homomorphism relative different bases, left-cols}
{
a_1=b\CRstar a_2\CRstar #1^{\CRInverse}
}

\AddEq[1]{representation of homomorphism relative different bases, right-cols}
{
a_1=#1^{\RCInverse}\RCstar a_2\RCstar b
}

\AddEq[1]{representation of homomorphism relative different bases, left-rows}
{
a_1=b\RCstar a_2\RCstar #1^{\RCInverse}
}

\AddEq[1]{representation of homomorphism relative different bases, right-rows}
{
a_1=#1^{\CRInverse}\CRstar a_2\CRstar b
}

\AddEq[2]{coordinate matrix, f, g, left-cols}
{
\Basis e_{#1 1}=#2\CRstar\Basis e_{#1 2}
}

\AddEq[2]{coordinate matrix, f, g, right-cols}
{
\Basis e_{#1 1}=\Basis e_{#1 2}\RCstar #2
}

\AddEq[2]{coordinate matrix, f, g, left-rows}
{
\Basis e_{#1 1}=#2\RCstar\Basis e_{#1 2}
}

\AddEq[2]{coordinate matrix, f, g, right-rows}
{
\Basis e_{#1 1}=\Basis e_{#1 2}\CRstar #2
}

\AddEq{a*b, left-cols}
{
\gdef\TheProduct{\ensuremath{a\CRstar b}}%
}

\AddEq{a*b, right-cols}
{
\gdef\TheProduct{\ensuremath{b\RCstar a}}%
}

\AddEq{a*b, left-rows}
{
\gdef\TheProduct{\ensuremath{a\RCstar b}}%
}

\AddEq{a*b, right-rows}
{
\gdef\TheProduct{\ensuremath{b\CRstar a}}%
}

\AddEq{two transformations on basis manifold, left-cols}
{
\[
e\CRstar g_1=e\CRstar g_2
\]
}

\AddEq{two transformations on basis manifold, right-cols}
{
\[
g_1\RCstar e=g_2\RCstar e
\]
}

\AddEq{two transformations on basis manifold, left-rows}
{
\[
e\RCstar g_1=e\RCstar g_2
\]
}

\AddEq{two transformations on basis manifold, right-rows}
{
\[
g_1\CRstar e=g_2\CRstar e
\]
}

\AddEq[3]{basis manifold of V, left-cols}
{
\ensuremath{\Basis {#1}\CRstar #2(V)}#3
}

\AddEq[3]{basis manifold of V, right-cols}
{
\ensuremath{#2(V)\RCstar \Basis {#1}}#3
}

\AddEq[3]{basis manifold of V, left-rows}
{
\ensuremath{\Basis {#1}\RCstar #2(V)}#3
}

\AddEq[3]{basis manifold of V, right-rows}
{
\ensuremath{#2(V)\CRstar \Basis {#1}}#3
}

\AddEq{basis manifold of vector space, left-cols}
{
\symb{\Basis e\CRstar G(V)}{basis manifold}1
}

\AddEq{basis manifold of vector space, right-cols}
{
\symb{G(V)\RCstar \Basis e}{basis manifold}1
}

\AddEq{basis manifold of vector space, left-rows}
{
\symb{\Basis e\RCstar G(V)}{basis manifold}1
}

\AddEq{basis manifold of vector space, right-rows}
{
\symb{G(V)\CRstar \Basis e}{basis manifold}1
}

\AddEq{Fg in GL}
{
$F(g)\in GL(W_*)$
}

\AddEq[1]{coordinate matrix of basis and passive transformation, left-cols}
{
\[
\aD{\Vector e'}i=\aD {#1}i\CRstar\Vector e
\]
}

\AddEq[1]{coordinate matrix of basis and passive transformation, right-cols}
{
\[
\aD{\Vector e'}i=\Vector e\RCstar\aD {#1}i
\]
}

\AddEq[1]{coordinate matrix of basis and passive transformation, left-rows}
{
\[
\aU{\Vector e'}i=\aU {#1}i\RCstar\Vector e
\]
}

\AddEq[1]{coordinate matrix of basis and passive transformation, right-rows}
{
\[
\aU{\Vector e'}i=\Vector e\CRstar\aU {#1}i
\]
}

\AddEq{coordinate matrix of basis and passive transformation, cols}
{
\[
\aD {e'}i=\aD ai
\]
}

\AddEq{coordinate matrix of basis and passive transformation, rows}
{
\[
\aU {e'}i=\aU ai
\]
}

\AddEq[1]{passive transformation e->e', left-cols}
{
\Basis e'_{#1}=g\CRstar \Basis e_{#1}
}

\AddEq[1]{passive transformation e->e', left-rows}
{
\Basis e'_{#1}=g\RCstar \Basis e_{#1}
}

\AddEq[1]{passive transformation e->e', right-cols}
{
\Basis e'_{#1}=\Basis e_{#1}\RCstar g
}

\AddEq[1]{passive transformation e->e', right-rows}
{
\Basis e'_{#1}=\Basis e_{#1}\CRstar g
}

\AddEquation{passive transformation e->e, left-cols}
{
\Basis e=\aD En\CRstar \Basis e
}

\AddEquation{passive transformation e->e, left-rows}
{
\Basis e=\aD En\RCstar \Basis e
}

\AddEquation{passive transformation e->e, right-cols}
{
\Basis e=\Basis e\RCstar \aD En
}

\AddEquation{passive transformation e->e, right-rows}
{
\Basis e=\Basis e\CRstar \aD En
}

\AddEq{homomorphism on A basis, left-cols}
{
\[
a=e^{\CRInverse}\CRstar e'
\]
}

\AddEq{homomorphism on A basis, right-cols}
{
\[
a=e'\RCstar e^{\RCInverse}
\]
}

\AddEq{homomorphism on A basis, left-rows}
{
\[
a=e^{\RCInverse}\RCstar e'
\]
}

\AddEq{homomorphism on A basis, right-rows}
{
\[
a=e'\CRstar e^{\CRInverse}
\]
}

\AddEq{passive transformation symbol, left-cols}
{
$a\CRstar \Basis e$.
}

\AddEq{passive transformation symbol, right-cols}
{
$\Basis e\RCstar a$.
}

\AddEq{passive transformation symbol, left-rows}
{
$a\RCstar \Basis e$.
}

\AddEq{passive transformation symbol, right-rows}
{
$\Basis e\CRstar a$.
}

\AddEq[1]{active transformation, vector space, left-cols}
{
$\Basis e\CRstar a$#1
}

\AddEq[1]{active transformation, vector space, right-cols}
{
$a\RCstar\Basis e$#1
}

\AddEq[1]{active transformation, vector space, left-rows}
{
$\Basis e\RCstar a$#1
}

\AddEq[1]{active transformation, vector space, right-rows}
{
$a\CRstar\Basis e$#1
}

\AddEq{active transformation ae x=a ex, left-cols}
{
v\CRstar\EqText{\Basis e\CRstar a}=\EqText{v\CRstar\Basis e}\CRstar a
}

\AddEq{active transformation ae x=a ex, right-cols}
{
\EqText{a\RCstar\Basis e}\RCstar v=a\RCstar\EqText{\Basis e\RCstar v}
}

\AddEq{active transformation ae x=a ex, left-rows}
{
v\RCstar\EqText{\Basis e\RCstar a}=\EqText{v\RCstar\Basis e}\RCstar a
}

\AddEq{active transformation ae x=a ex, right-rows}
{
\EqText{a\CRstar\Basis e }\CRstar v=a\CRstar\EqText{\Basis e \CRstar v}
}

\AddEq{active transformation ae x=a ex 1, left-cols}
{
$v\CRstar\Basis e$
}

\AddEq{active transformation ae x=a ex 1, right-cols}
{
$\Basis e\RCstar v$
}

\AddEq{active transformation ae x=a ex 1, left-rows}
{
$v\RCstar\Basis e$
}

\AddEq{active transformation ae x=a ex 1, right-rows}
{
$\Basis e \CRstar v$
}

\AddEq{v*g, left-cols}
{
$v\CRstar g$.
}

\AddEq{v*g, right-cols}
{
$g\RCstar v$.
}

\AddEq{v*g, left-rows}
{
$v\RCstar g$.
}

\AddEq{v*g, right-rows}
{
$g\CRstar v$.
}

\AddEq{active transformations, vector space, 2, left-cols}
{
\[
\aU vi=\aU vj\aUD aij
\]
}

\AddEq{active transformations, vector space, 2, right-cols}
{
\[
\aU vi=\aUD aij\aU vj
\]
}

\AddEq{active transformations, vector space, 2, left-rows}
{
\[
\aD vi=\aD vj\aUD aji
\]
}

\AddEq{active transformations, vector space, 2, right-rows}
{
\[
\aD vi=\aUD aji\aD vj
\]
}

\AddEq{ai=delta ik, cols}
{
$\aU vi=\aUD{\delta}ik$
}

\AddEq{ai=delta ik, rows}
{
$\aD vi=\aUD{\delta}ki$
}

\AddEq{identity transformation, vector space, 1, cols}
{
\[
\aUD{\delta}ik=\aUD aik
\]
}

\AddEq{identity transformation, vector space, 1, rows}
{
\[
\aUD{\delta}ki=\aUD aki
\]
}

\AddEquation{identity transformation, vector space, left-cols}
{
\aUD{\delta}ik=\aUD{\delta}jk\aUD aij
}

\AddEquation{identity transformation, vector space, right-cols}
{
\aUD{\delta}ik=\aUD aij\aUD{\delta}jk
}

\AddEquation{identity transformation, vector space, left-rows}
{
\aUD{\delta}ki=\aUD{\delta}kj\aUD aji
}

\AddEquation{identity transformation, vector space, right-rows}
{
\aUD{\delta}ik=\aUD aji\aUD{\delta}kj
}

\AddEq{automorphism, vector space, 1, left-cols}
{
\aD{e'}i=\aD ei\CRstar a
}

\AddEq{automorphism, vector space, 1, right-cols}
{
\aD{e'}i=a\RCstar \aD ei
}

\AddEq{automorphism, vector space, 1, left-rows}
{
\aU{e'}i=\aU ei\RCstar a
}

\AddEq{automorphism, vector space, 1, right-rows}
{
\aU{e'}i=a\CRstar \aU ei
}

\AddEquation{b(av)=a(bv) left}
{
\begin{aligned}
\RedText{(ba)v}&=\BlueText{b(av)=\Vector f\circ(av)}
\\ &=\BlueText{a(\Vector f\circ v)=a(bv)}
\\ &=\RedText{(ab)v}
\end{aligned}
}

\AddEquation{b(av)=a(bv) right}
{
\begin{aligned}
\RedText{v(ab)}&=\BlueText{(va)b=\Vector f\circ(va)}
\\ &=\BlueText{(\Vector f\circ v)a=(vb)a}
\\ &=\RedText{v(ba)}
\end{aligned}
}

\AddEquation{fov=bv left}
{
\Vector f\circ v=bv
}

\AddEquation{fov=bv right}
{
\Vector f\circ v=vb
}

\AddEquation{fovi=di vi left}
{
f\circ v(\gi i)=v(\gii)b(\gii)
}

\AddEquation{fovi=di vi right}
{
f\circ v(\gi i)=b(\gii)v(\gii)
}

\AddEquation{fovi=di vi}
{
f\circ v(\gi i)=b(\gii)v(\gii)
}

\AddEq{basis vector cols}
{
$\aD ei$
}

\AddEq{basis vector rows}
{
$\aU ei$
}

\AddEquation{automorphism, vector space, 2, left-cols}
{
\lambda\CRstar e'=0
}

\AddEquation{automorphism, vector space, 2, right-cols}
{
e'\RCstar\lambda =0
}

\AddEquation{automorphism, vector space, 2, left-rows}
{
\lambda\RCstar e'=0
}

\AddEquation{automorphism, vector space, 2, right-rows}
{
e'\CRstar\lambda =0
}

\AddEq{automorphism, vector space, 3, left-cols}
{
\[
\lambda\CRstar e'\CRstar a^{\CRInverse}
=\lambda\CRstar e=0
\]
}

\AddEq{automorphism, vector space, 3, right-cols}
{
\[
a^{\CRInverse}\RCstar e' \RCstar\lambda
=e\RCstar\lambda =0
\]
}

\AddEq{automorphism, vector space, 3, left-rows}
{
\[
\lambda\RCstar e'\RCstar a^{\CRInverse}
=\lambda\RCstar e=0
\]
}

\AddEq{automorphism, vector space, 3, right-rows}
{
\[
a^{\CRInverse}\CRstar e' \CRstar\lambda
=e\CRstar\lambda =0
\]
}

\AddEq[1]{vector of basis cols}
{
\ensuremath{\aD {#1}i}
}

\AddEq[1]{vector of basis rows}
{
\ensuremath{\aU {#1}i}
}

\AddEq[3]{ARow u2i cols}
{
\ensuremath{\aD[#2]{#1}{#3}}
}

\AddEq[3]{ARow u2i rows}
{
\ensuremath{\aU {#1_{#2}}{#3}}
}

\AddEq{Automorphisms of vector space, cr-rows}
{
\[
\begin{aligned}
v'&=f\CRstar v\\
v''=g\CRstar v'&=g\CRstar f\CRstar v
\end{aligned}
\]
}

\AddEq{Automorphisms of vector space, rc-rows}
{
\[
\begin{aligned}
v'&=v\RCstar f\\
v''=v'\RCstar g&=v\RCstar f\RCstar g
\end{aligned}
\]
}

\AddEq{Automorphisms of vector space, cr-cols}
{
\[
\begin{aligned}
v'&=v\CRstar f\\
v''=v'\CRstar g&=v\CRstar f\CRstar g
\end{aligned}
\]
}

\AddEq{Automorphisms of vector space, rc-cols}
{
\[
\begin{aligned}
v'&=f\RCstar v\\
v''=g\RCstar v'&=g\RCstar f\RCstar v
\end{aligned}
\]
}

\AddEq{basis e2 of V cols}
{
\ECol 2j, $\gi j\in\gi J$,
}

\AddEq{basis e2 of V rows}
{
\ERow 2j, $\gi j\in\gi J$,
}

\AddEq{basis e2 of V lambda}
{
\lambda=0
}

\AddEq{basis e2 relative e1 left-cols}
{
\[
\ECol 2j=\aD aj\CRstar e_1
\]
}

\AddEq{basis e2 relative e1 right-cols}
{
\[
\ECol 2j=e_1\RCstar \aD aj
\]
}

\AddEq{basis e2 relative e1 right-rows}
{
\[
\ERow 2j=e_1\CRstar \aU aj
\]
}

\AddEq{basis e2 relative e1 left-rows}
{
\[
\ERow 2j=\aU aj\RCstar e_1
\]
}

\AddEquation{basis e2 relative e1, lambda, right-cols}
{
e_2\RCstar \lambda
= e_1\RCstar a\RCstar\lambda=0
}

\AddEquation{basis e2 relative e1, lambda, left-cols}
{
\lambda\CRstar e_2
=\lambda\CRstar a\CRstar e_1=0
}

\AddEquation{basis e2 relative e1, lambda, right-rows}
{
e_2\CRstar \lambda
= e_1\CRstar a\CRstar\lambda=0
}

\AddEquation{basis e2 relative e1, lambda, left-rows}
{
\lambda\RCstar e_2
=\lambda\RCstar a\RCstar e_1=0
}

\AddEquation{basis e2 relative e1, lambda=0, right-cols}
{
a\RCstar\lambda=0
}

\AddEquation{basis e2 relative e1, lambda=0, left-cols}
{
\lambda\CRstar a=0
}

\AddEquation{basis e2 relative e1, lambda=0, right-rows}
{
a\CRstar\lambda=0
}

\AddEquation{basis e2 relative e1, lambda=0, left-rows}
{
\lambda\RCstar a=0
}

\AddEquation{avi=bvi left-cols}
{
a\aU vi=b\aU vi
}

\AddEquation{avi=bvi right-cols}
{
\aU via=\aU vib
}

\AddEquation{avi=bvi left-rows}
{
a\aD vi=b\aD vi
}

\AddEquation{avi=bvi right-rows}
{
\aD via=\aD vib
}

\AddEq{Rank a<m cr}
{
$\CRRank a\le\gi m$
}

\AddEq{Rank a<m rc}
{
$\RCRank a\le\gi m$
}

\AddEq{av=bv, v left}
{
\forall v\in V,av=bv\ \ \ \ a,b\in A
}

\AddEq{av=bv, v right}
{
\forall v\in V,va=vb\ \ \ \ a,b\in A
}

\AddEq[1]{w=v.rc.f}
{
w_{#1}=f_{#1}\RCstar v_{#1}
}

\AddEq[1]{w=v.cr.f}
{
w_{#1}=v_{#1}\CRstar f_{#1}
}

\AddEq[1]{v2 f2=v2 a f1 a- left-cols}
{
v_2\CRstar f_2=v_2\CRstar #1\CRstar f_1\CRstar #1^{\CRInverse}
}

\AddEq[1]{v2 f2=v2 a f1 a- right-cols}
{
f_2\RCstar v_2=#1^{\RCInverse}\RCstar f_1\RCstar #1\RCstar v_2
}

\AddEq[1]{v2 f2=v2 a f1 a- left-rows}
{
v_2\RCstar f_2=v_2\RCstar #1\RCstar f_1\RCstar #1^{\RCInverse}
}

\AddEq[1]{v2 f2=v2 a f1 a- right-rows}
{
f_2\CRstar v_2=#1^{\CRInverse}\CRstar f_1\CRstar #1\CRstar v_2
}

\AddEq[1]{w2=...v2 a f1 a- left-cols}
{
w_2=v_1\CRstar f_1=v_2\CRstar #1\CRstar f_1\CRstar #1^{\CRInverse}
}

\AddEq[1]{w2=...v2 a f1 a- right-cols}
{
w_2=f_1\RCstar v_1=#1^{\RCInverse}\RCstar f_1\RCstar #1\RCstar v_2
}

\AddEq[1]{w2=...v2 a f1 a- left-rows}
{
w_2=v_1\RCstar f_1=v_2\RCstar #1\RCstar f_1\RCstar #1^{\RCInverse}
}

\AddEq[1]{w2=...v2 a f1 a- right-rows}
{
w_2=f_1\CRstar v_1=#1^{\CRInverse}\CRstar f_1\CRstar #1\CRstar v_2
}

\AddEq[1]{w2=...v2 a f1 left-cols}
{
w_2\CRstar #1=v_1\CRstar f_1=v_2\CRstar #1\CRstar f_1
}

\AddEq[1]{w2=...v2 a f1 right-cols}
{
#1\RCstar w_2=f_1\RCstar v_1=f_1\RCstar #1\RCstar v_2
}

\AddEq[1]{w2=...v2 a f1 left-rows}
{
w_2\RCstar #1=v_1\RCstar f_1=v_2\RCstar #1\RCstar f_1
}

\AddEq[1]{w2=...v2 a f1 right-rows}
{
#1\CRstar w_2=f_1\CRstar v_1=f_1\CRstar #1\CRstar v_2
}

\AddEq[4]{f2=a f1 a-}
{
#1_2=
\ShowEq{a f1 a- #3-#4}{#1}{#2}
}

\AddEq[2]{a f1 a- left-cols}
{
#2\CRstar #1_1\CRstar #2^{\CRInverse}
}

\AddEq[2]{a f1 a- right-cols}
{
#2^{\RCInverse}\RCstar #1_1\RCstar #2
}

\AddEq[2]{a f1 a- left-rows}
{
#2\RCstar #1_1\RCstar #2^{\RCInverse}
}

\AddEq[2]{a f1 a- right-rows}
{
#2^{\CRInverse}\CRstar #1_1\CRstar #2
}

\AddEq[1]{f'=g f g- left-cols}
{
f'_{#1}=g\CRstar f_{#1}\CRstar g^{\CRInverse}
}

\AddEq[1]{f'=g f g- right-cols}
{
f'_{#1}=g^{\RCInverse}\RCstar f_{#1}\RCstar g
}

\AddEq[1]{f'=g f g- left-rows}
{
f'_{#1}=g\RCstar f_{#1}\RCstar g^{\RCInverse}
}

\AddEq[1]{f'=g f g- right-rows}
{
f'_{#1}=g^{\CRInverse}\CRstar f_{#1}\CRstar g
}

\AddEquation{f1+f2=g.g- left-cols}
{
\begin{aligned}
f'&=g\CRstar f\CRstar g^{\CRInverse}
\\ &=g\CRstar (f_1+f_2)\CRstar g^{\CRInverse}
\\ &=g\CRstar f_1\CRstar g^{\CRInverse}\\ &+g\CRstar f_2\CRstar g^{\CRInverse}
\\ &=f'_1+f'_2
\end{aligned}
}

\AddEquation{f1+f2=g.g- right-cols}
{
\begin{aligned}
f'&=g^{\RCInverse}\RCstar f\RCstar g
\\ &=g^{\RCInverse}\RCstar (f_1+f_2)\RCstar g
\\ &=g^{\RCInverse}\RCstar f_1\RCstar g\\ &+g^{\RCInverse}\RCstar f_2\RCstar g
\\ &=f'_1+f'_2
\end{aligned}
}

\AddEquation{f1+f2=g.g- left-rows}
{
\begin{aligned}
f'&=g\RCstar f\RCstar g^{\RCInverse}
\\ &=g\RCstar (f_1+f_2)\RCstar g^{\RCInverse}
\\ &=g\RCstar f_1\RCstar g^{\RCInverse}\\ &+g\RCstar f_2\RCstar g^{\RCInverse}
\\ &=f'_1+f'_2
\end{aligned}
}

\AddEquation{f1+f2=g.g- right-rows}
{
\begin{aligned}
f'&=g^{\CRInverse}\CRstar f\CRstar g
\\ &=g^{\CRInverse}\CRstar (f_1+f_2)\CRstar g
\\ &=g^{\CRInverse}\CRstar f_1\CRstar g\\ &+g^{\CRInverse}\CRstar f_2\CRstar g
\\ &=f'_1+f'_2
\end{aligned}
}

\AddEquation{f1*f2=g.g- left-cols}
{
\begin{aligned}
f'&=g\CRstar f\CRstar g^{\CRInverse}
\\ &=g\CRstar f_1\CRstar f_2\CRstar g^{\CRInverse}
\\ &=g\CRstar f_1\CRstar g^{\CRInverse}\\ &\CRstar g\CRstar f_2\CRstar g^{\CRInverse}
\\ &=f'_1\CRstar f'_2
\end{aligned}
}

\AddEquation{f1*f2=g.g- right-cols}
{
\begin{aligned}
f'&=g^{\RCInverse}\RCstar f\RCstar g
\\ &=g^{\RCInverse}\RCstar f_1\RCstar f_2\RCstar g
\\ &=g^{\RCInverse}\RCstar f_1\RCstar g\\ &\RCstar g^{\RCInverse}\RCstar f_2\RCstar g
\\ &=f'_1\RCstar f'_2
\end{aligned}
}

\AddEquation{f1*f2=g.g- left-rows}
{
\begin{aligned}
f'&=g\RCstar f\RCstar g^{\RCInverse}
\\ &=g\RCstar f_1\RCstar f_2\RCstar g^{\RCInverse}
\\ &=g\RCstar f_1\RCstar g^{\RCInverse}\\ &\RCstar g\RCstar f_2\RCstar g^{\RCInverse}
\\ &=f'_1\RCstar f'_2
\end{aligned}
}

\AddEquation{f1*f2=g.g- right-rows}
{
\begin{aligned}
f'&=g^{\CRInverse}\CRstar f\CRstar g
\\ &=g^{\CRInverse}\CRstar f_1\CRstar f_2\CRstar g
\\ &=g^{\CRInverse}\CRstar f_1\CRstar g\\ &\CRstar g^{\CRInverse}\CRstar f_2\CRstar g
\\ &=f'_1\CRstar f'_2
\end{aligned}
}

\AddEq{f1*f2=g.g-}
{
\eqRef{f=f1*f2 endo \SideNS-\Cols}{\Product},
\eqRef{f'=g f g- \SideNS-\Cols}{f*},
\eqRef{f'=g f g- \SideNS-\Cols}{1*},
\eqRef{f'=g f g- \SideNS-\Cols}{2*}.
}

\AddEq{f1+f2=g.g-}
{
\eqRef{f=f1+f2 endo}{\SideNS-\Cols},
\eqRef{f'=g f g- \SideNS-\Cols}f,
\eqRef{f'=g f g- \SideNS-\Cols}1,
\eqRef{f'=g f g- \SideNS-\Cols}2.
}

\AddEq[2]{Bases e12}
{
\eV[#1 1], \eV[#1 2]#2
}

\AddEq{A:V->W}
{
a:V\rightarrow W
}

\AddEq[4]{(f-ae)v cr-cols}
{
{#4}\CRstar(\ShowEq{(f-ae)v matrix}{#1}{}{#3}{#2})=0
}

\AddEq[4]{(f-ae)v rc-cols}
{
(\ShowEq{(f-ae)v matrix}{#1}{}{#3}{#2})\RCstar {#4}=0
}

\AddEq[4]{(f-ae)v cr-rows}
{
(\ShowEq{(f-ae)v matrix}{#1}{}{#3}{#2})\CRstar {#4}=0
}

\AddEq[4]{(f-ae)v rc-rows}
{
{#4}\RCstar(\ShowEq{(f-ae)v matrix}{#1}{}{#3}{#2})=0
}

\AddEquation{(f-bEn)v cr-cols}
{
v_1\CRstar(\ShowEq{f-bEn left}1)=0
}

\AddEquation{(f-bEn)v rc-cols}
{
(\ShowEq{f-bEn right}1)\RCstar v_1=0
}

\AddEquation{(f-bEn)v cr-rows}
{
(\ShowEq{f-bEn right}1)\CRstar v_1=0
}

\AddEquation{(f-bEn)v rc-rows}
{
v_1\RCstar(\ShowEq{f-bEn left}1)=0
}

\AddEquation{(f-bEn)v2 cr-cols}
{
v_2\CRstar(f_2-\ShowEq{ga*g- \SideNS-\Cols}bg)=0
}

\AddEquation{(f-bEn)v2 rc-cols}
{
(f_2-\ShowEq{ga*g- \SideNS-\Cols}bg)\RCstar v_2=0
}

\AddEquation{(f-bEn)v2 cr-rows}
{
(f_2-\ShowEq{ga*g- \SideNS-\Cols}bg)\CRstar v_2=0
}

\AddEquation{(f-bEn)v2 rc-rows}
{
v_2\RCstar(f_2-\ShowEq{ga*g- \SideNS-\Cols}bg)=0
}

\AddEq[4]{(f-ae)v matrix}
{
\ensuremath{#1_{#2}-\ShowEq{ga*g- \SideNS-\Cols}{#3}{#4}}
}

\AddEq{bg.rc.g-}
{
\[
\begin{aligned}
\ShowEq{ga*g- left-rows}bg
&=g\RCstar (bg^{\RCInverse})
\\ &=(g^{\RCInverse})^{\RCInverse}\RCstar (bg^{\RCInverse})
\end{aligned}
\]
}

\AddEquation{gb rc g-=b}
{
\ShowEq{ga*g- right-cols}bg
=(g^{\RCInverse}\RCstar  g)b=\aD Enb
}

\AddEquation{gb cr g-=b}
{
\ShowEq{ga*g- left-cols}bg
=b(g\CRstar g^{\CRInverse})=b\aD En
}

\AddEquation{gij in ZA}
{
\aUD gij\in Z(A)\ \ \ \,\gii=\gi 1, ..., \gin\ \ \ \,\gij=\gi 1, ..., \gin
}

\AddEq{bg.cr.g-}
{
\[
\begin{aligned}
\ShowEq{ga*g- right-rows}bg
&=(g^{\CRInverse} b)\CRstar g
\\ &=(g^{\CRInverse} b)\CRstar (g^{\CRInverse})^{\CRInverse}
\end{aligned}
\]
}

\AddEq{rc-eigencols(f,g)}
{
f\RCstar v=\ShowEq{ga*g- \SideNS-\Cols}bg\RCstar v
}

\AddEq{cr-eigencols(f,g)}
{
v\CRstar f=v\CRstar \ShowEq{ga*g- \SideNS-\Cols}bg
}

\AddEq{rc-eigenrows(f,g)}
{
v\RCstar f=v\RCstar \ShowEq{ga*g- \SideNS-\Cols}bg
}

\AddEq{cr-eigenrows(f,g)}
{
f\CRstar v=\ShowEq{ga*g- \SideNS-\Cols}bg\CRstar v
}

\AddEq[1]{pair of matrices rc-column}
{
$(f,g)$#1
}

\AddEq[1]{pair of matrices cr-column}
{
$(f,g)$#1
}

\AddEq[1]{pair of matrices rc-row}
{
$(f,g^{\RCInverse})$#1
}

\AddEq[1]{pair of matrices cr-row}
{
$(f,g^{\CRInverse})$#1
}

\AddEq{rc-eigencols}
{
f\RCstar v=bv
}

\AddEq{rc-eigenrows}
{
v\RCstar f=vb
}

\AddEq{cr-eigencols}
{
v\CRstar f=vb
}

\AddEq{cr-eigenrows}
{
f\CRstar v=bv
}

\AddEq{gi in I, gj in J}
{
\[
\gii\in\giI,\ \gij\in\giJ
\]
}

\AddEquation{rc-det ij aIJ=0}
{
\RCdet ij\,\aUD{(a-bE)}IJ=0
}

\AddEquation{cr-det ij aIJ=0}
{
\CRdet ij\,\aUD{(a-bE)}IJ=0
}

\AddEq{rc-det ij a=0}
{
\RCdet ij\,(a-bE)=0
}

\AddEq{cr-det ij a=0}
{
\CRdet ij\,(a-bE)=0
}

\AddEq{f-ae->f-be}
{
\[
\begin{aligned}
&\, \ShowEq{(f-ae)v matrix}f2bg
\\ =&\,\ShowEq{a f1 a- \SideNS-\Cols}fg
-\ShowEq{ga*g- \SideNS-\Cols}bg
\\ =&\,
\ShowEq{g(f-bEn)g \SideNS-\Cols}
\end{aligned}
\]
}

\AddEq{f-ae->f-be left-cols}
{
\[
\begin{aligned}
&\, v_2\CRstar(\ShowEq{(f-ae)v matrix}f2bg)
\\ =&\,
v_2\CRstar(\ShowEq{a f1 a- \SideNS-\Cols}fg
-\ShowEq{ga*g- \SideNS-\Cols}bg)
\\ =&\,
\ShowEq{v1g- \SideNS-\Cols}v1g{}\CRstar
\ShowEq{g(f-bEn)g \SideNS-\Cols}
\end{aligned}
\]
}

\AddEq{f-ae->f-be right-cols}
{
\[
\begin{aligned}
&\,(\ShowEq{(f-ae)v matrix}f2bg)\RCstar v_2
\\ =&\,
(\ShowEq{a f1 a- \SideNS-\Cols}fg
-\ShowEq{ga*g- \SideNS-\Cols}bg)\RCstar v_2
\\ =&\,
\ShowEq{g(f-bEn)g \SideNS-\Cols}\RCstar
\ShowEq{v1g- \SideNS-\Cols}v1g{}
\end{aligned}
\]
}

\AddEq{f-ae->f-be left-rows}
{
\[
\begin{aligned}
&\, v_2\RCstar(\ShowEq{(f-ae)v matrix}f2bg)
\\ =&\,
v_2\RCstar(\ShowEq{a f1 a- \SideNS-\Cols}fg
-\ShowEq{ga*g- \SideNS-\Cols}bg)
\\ =&\,
\ShowEq{v1g- \SideNS-\Cols}v1g{}\RCstar
\ShowEq{g(f-bEn)g \SideNS-\Cols}
\end{aligned}
\]
}

\AddEq{f-ae->f-be right-rows}
{
\[
\begin{aligned}
&\,(\ShowEq{(f-ae)v matrix}f2bg)\CRstar v_2
\\ =&\,
(\ShowEq{a f1 a- \SideNS-\Cols}fg
-\ShowEq{ga*g- \SideNS-\Cols}bg)\CRstar v_2
\\ =&\,
\ShowEq{g(f-bEn)g \SideNS-\Cols}\CRstar
\ShowEq{v1g- \SideNS-\Cols}v1g{}
\end{aligned}
\]
}

\AddEq{rc-det ij f-bg IJ=0}
{
\[
\RCdet ij\,\aUD{(\ShowEq{(f-ae)v matrix}f{}bg)}IJ=0
\]
\[
\gii\in\giI,\ \gij\in\giJ
\]
}

\AddEq{cr-det ij f-bg IJ=0}
{
\[
\CRdet ij\,\aUD{(\ShowEq{(f-ae)v matrix}f{}bg)}IJ=0
\]
\[
\gii\in\giI,\ \gij\in\giJ
\]
}

\AddEq{rc-det ij f-bg =0}
{
\RCdet ij\,(\ShowEq{(f-ae)v matrix}f{}bg)=0
}

\AddEq{cr-det ij f-bg =0}
{
\CRdet ij\,(\ShowEq{(f-ae)v matrix}f{}bg)=0
}

\AddEq{rc-rank a=k<n}
{
\[\RCstar\rank a=\gik<\gin\]
}

\AddEq{cr-rank a=k<n}
{
\[\CRstar\rank a=\gik<\gin\]
}

\AddEquation{rc-eigencols=0 fg}
{
(f-\ShowEq{ga*g- \SideNS-\Cols}bg)\RCstar v=0
}

\AddEquation{rc-eigenrows=0 fg}
{
v\RCstar(f-\ShowEq{ga*g- \SideNS-\Cols}bg)=0
}

\AddEquation{cr-eigencols=0 fg}
{
v\CRstar(f-\ShowEq{ga*g- \SideNS-\Cols}bg)=0
}

\AddEquation{cr-eigenrows=0 fg}
{
(f-\ShowEq{ga*g- \SideNS-\Cols}bg)\CRstar v=0
}

\AddEquation{rc-eigencols=0}
{
(\ShowEq{f-bEn right}{})\RCstar v=0
}

\AddEquation{rc-eigenrows=0}
{
v\RCstar(\ShowEq{f-bEn left}{})=0
}

\AddEquation{cr-eigencols=0}
{
v\CRstar(\ShowEq{f-bEn left}{})=0
}

\AddEquation{cr-eigenrows=0}
{
(\ShowEq{f-bEn right}{})\CRstar v=0
}

\AddEq{basis e1=e2}
{
$\Basis e_1=\Basis e_2$.
}

\AddEq{basis e1 ne e2}
{
$\Basis e_1\ne\Basis e_2$.
}

\AddEq[2]{f-bEn left}
{
\ensuremath{f_{#1}-\aD Enb}#2
}

\AddEq[2]{f-bEn right}
{
\ensuremath{f_{#1}-b\aD En}#2
}

\AddEq[4]{Bases eVW}
{
\eV[#1 #3][,] \eV[#2 #3][#4]
}

\AddEq{f=g+h}
{
\[f=g+h\]
}

\AddEquation{fov=(g+h)ov left-cols}
{
\begin{aligned}
&\,\Vector f\circ (v\CRstar e_U)
\\=&\,v\CRstar f\CRstar e_V
\\=&\,v\CRstar(g+h)\CRstar e_V
\\=&\,v\CRstar g\CRstar e_V+v\CRstar h\CRstar e_V
\\=&\,\Vector g\circ (v\CRstar e_U)+\Vector h\circ (v\CRstar e_U)
\end{aligned}
}

\AddEquation{fov=(g+h)ov right-cols}
{
\begin{aligned}
&\,\Vector f\circ (e_U\RCstar v)
\\=&\,e_V\RCstar f\RCstar v
\\=&\,e_V\RCstar(g+h)\RCstar v
\\=&\,e_V\RCstar g\RCstar v+e_V\RCstar h\RCstar v
\\=&\,\Vector g\circ (e_U\RCstar v)+\Vector h\circ (e_U\RCstar v)
\end{aligned}
}

\AddEquation{fov=(g+h)ov left-rows}
{
\begin{aligned}
&\,\Vector f\circ (v\RCstar e_U)
\\=&\,v\RCstar f\RCstar e_V
\\=&\,v\RCstar(g+h)\RCstar e_V
\\=&\,v\RCstar g\RCstar e_V+v\RCstar h\RCstar e_V
\\=&\,\Vector g\circ (v\RCstar e_U)+\Vector h\circ (v\CRstar e_U)
\end{aligned}
}

\AddEquation{fov=(g+h)ov right-rows}
{
\begin{aligned}
&\,\Vector f\circ (e_U\CRstar v)
\\=&\,e_V\CRstar f\CRstar v
\\=&\,e_V\CRstar(g+h)\CRstar v
\\=&\,e_V\CRstar g\CRstar v+e_V\CRstar h\CRstar v
\\=&\,\Vector g\circ (e_U\CRstar v)+\Vector h\circ (e_U\RCstar v)
\end{aligned}
}

\AddEquation{fov=(gh)ov left-cols}
{
\begin{aligned}
&\,\Vector f\circ (u\CRstar e_U)
=u\CRstar f\CRstar e_W
\\=&\,u\CRstar g\CRstar h\CRstar e_W
\\=&\,\Vector h\circ(u\CRstar g\CRstar e_V)
\\=&\,\Vector h\circ(\Vector g\circ(u\CRstar e_U))
\\=&\,(\Vector h\circ\Vector g)\circ(u\CRstar e_U)
\end{aligned}
}

\AddEquation{fov=(gh)ov right-cols}
{
\begin{aligned}
&\,\Vector f\circ (e_U\RCstar u)
=e_W\RCstar f\RCstar u
\\=&\,e_W\RCstar h\RCstar g\RCstar u
\\=&\,\Vector h\circ(e_V\RCstar g\RCstar u)
\\=&\,\Vector h\circ(\Vector g\circ(e_U\RCstar u))
\\=&\,(\Vector h\circ\Vector g)\circ(e_U\RCstar u)
\end{aligned}
}

\AddEquation{fov=(gh)ov left-rows}
{
\begin{aligned}
&\,\Vector f\circ (u\RCstar e_U)
=u\RCstar f\RCstar e_W
\\=&\,u\RCstar g\RCstar h\RCstar e_W
\\=&\,\Vector h\circ(u\RCstar g\RCstar e_V)
\\=&\,\Vector h\circ(\Vector g\circ(u\RCstar e_U))
\\=&\,(\Vector h\circ\Vector g)\circ(u\RCstar e_U)
\end{aligned}
}

\AddEquation{fov=(gh)ov right-rows}
{
\begin{aligned}
&\,\Vector f\circ (e_U\CRstar u)
=e_W\CRstar f\CRstar u
\\=&\,e_W\CRstar h\CRstar g\CRstar u
\\=&\,\Vector h\circ(e_V\CRstar g\CRstar u)
\\=&\,\Vector h\circ(\Vector g\circ(e_U\CRstar u))
\\=&\,(\Vector h\circ\Vector g)\circ(e_U\CRstar u)
\end{aligned}
}

\AddEq{w=f o v}
{
w=f\circ v
}

\AddEquation{fov=gov+hov left-cols}
{
\begin{aligned}
&\,\Vector f\circ(v\CRstar e_U)
\\=&\,\Vector g\circ(v\CRstar e_U)+\Vector h\circ(v\CRstar e_U)
\\=&\,v\CRstar g\CRstar e_V+v\CRstar h\CRstar e_V
\\=&\,v\CRstar(g+h)\CRstar e_V
\end{aligned}
}

\AddEquation{fov=gov+hov right-cols}
{
\begin{aligned}
&\,\Vector f\circ(e_U\RCstar v)
\\=&\,\Vector g\circ(e_U\RCstar v)+\Vector h\circ(e_U\RCstar v)
\\=&\,e_V\RCstar g\RCstar v+e_V\RCstar h\RCstar v
\\=&\,e_V\RCstar (g+h)\RCstar v
\end{aligned}
}

\AddEquation{fov=gov+hov left-rows}
{
\begin{aligned}
&\,\Vector f\circ(v\RCstar e_U)
\\=&\,\Vector g\circ(v\RCstar e_U)+\Vector h\circ(v\RCstar e_U)
\\=&\,v\RCstar g\RCstar e_V+v\RCstar h\RCstar e_V
\\=&\,v\RCstar(g+h)\CRstar e_V
\end{aligned}
}

\AddEquation{fov=gov+hov right-rows}
{
\begin{aligned}
&\,\Vector f\circ(e_U\CRstar v)
\\=&\,\Vector g\circ(e_U\CRstar v)+\Vector h\circ(e_U\CRstar v)
\\=&\,e_V\CRstar g\CRstar v+e_V\CRstar h\CRstar v
\\=&\,e_V\CRstar (g+h)\CRstar v
\end{aligned}
}

\AddEquation{f o v=g o h o v left-cols}
{
\begin{aligned}
&\,\Vector f\circ(u\CRstar e_U)
\\=&\,\Vector h\circ(\Vector g\circ(u\CRstar e_U))
\\=&\,\Vector h\circ(u\CRstar g\CRstar e_V)
\\=&\,u\CRstar g\CRstar h\CRstar e_W
\end{aligned}
}

\AddEquation{f o v=g o h o v right-cols}
{
\begin{aligned}
&\,\Vector f\circ(e_U\RCstar u)
\\=&\,\Vector h\circ(\Vector g\circ(e_U\RCstar u))
\\=&\,\Vector h\circ(u\RCstar g\RCstar e_V)
\\=&\,e_W\RCstar h\RCstar g\RCstar u
\end{aligned}
}

\AddEquation{f o v=g o h o v left-rows}
{
\begin{aligned}
&\,\Vector f\circ(u\RCstar e_U)
\\=&\,\Vector h\circ(\Vector g\circ(u\RCstar e_U))
\\=&\,\Vector h\circ(u\RCstar g\RCstar e_V)
\\=&\,u\RCstar g\RCstar h\RCstar e_W
\end{aligned}
}

\AddEquation{f o v=g o h o v right-rows}
{
\begin{aligned}
&\,\Vector f\circ(e_U\CRstar u)
\\=&\,\Vector h\circ(\Vector g\circ(e_U\CRstar u))
\\=&\,\Vector h\circ(u\CRstar g\CRstar e_V)
\\=&\,u\CRstar g\CRstar h\CRstar e_W
\end{aligned}
}

\AddEq{sum of homomorphisms f o v=}
{
\forall v\in V: f\circ v=g\circ v+h\circ v
}

\AddEq{commutative law}
{
v+w=w+v
}

\AddEq{(g+h)o(u+v)=}
{
\begin{aligned}
&\,f\circ (u+v)
\\=&\,g\circ(u+v)+h\circ(u+v)
\\=&\,g\circ u+g\circ v
\\+&\,h\circ u+h\circ u
\\=&\,f\circ u+f\circ v
\end{aligned}
}

\AddEq{(h o g)o(u+v)=}
{
\begin{aligned}
&\,f\circ (u+v)
\\=&\,h\circ (g\circ(u+v))
\\=&\,h\circ (g\circ u+g\circ v))
\\=&\,h\circ (g\circ u)+h\circ (g\circ v)
\\=&\,f\circ u+f\circ v
\end{aligned}
}

\AddEquation{(g+h)o(av)= left}
{
\begin{aligned}
&\,f\circ (av)
\\=&\,g\circ(av)+h\circ(av)
\\=&\,a(g\circ v)+a(h\circ v)
\\=&\,a(g\circ v+h\circ v)
\\=&\,a(f\circ v)
\end{aligned}
}

\AddEquation{(h o g)o(av)= right}
{
\begin{aligned}
f\circ (va)
&=h \circ (g\circ(va))
\\ &=h\circ ((g\circ v)a)
\\ &=(h\circ (g\circ v))a
\\ &=(f\circ v)a
\end{aligned}
}

\AddEquation{(h o g)o(av)= left}
{
\begin{aligned}
f\circ (av)
&=h \circ (g\circ(av))
\\ &=h\circ (a(g\circ v))
\\ &=a(h\circ (g\circ v))
\\ &=a(f\circ v)
\end{aligned}
}

\AddEquation{(g+h)o(av)= right}
{
\begin{aligned}
&\,f\circ (va)
\\=&\,g\circ(va)+h\circ(va)
\\=&\,(g\circ v)a+(h\circ v)a
\\=&\,(g\circ v+h\circ v)a
\\=&\,(f\circ v)a
\end{aligned}
}

\AddEq{f=h o g}
{
$f=h\circ g$
}

\AddEq{f1 f2 left}
{
$f_1$, $f_2$
}

\AddEq{f1 f2 right}
{
$f_2$, $f_1$
}

\AddEq{vf1 vf2 left}
{
$\Vector f_1$, $\Vector f_2$
}

\AddEq{vf1 vf2 right}
{
$\Vector f_2$, $\Vector f_1$
}

\AddEq{diagram product of homomorphisms, A vector space}
{
\xymatrix{
U\ar[rr]^f\ar[dr]^g & & W\\
&V\ar[ur]^h &
}
}

\AddEq{f o u=h o g o u}
{
\begin{aligned}
\forall u\in U:
f\circ u&=(h\circ g)\circ u
\\ &=h\circ(g\circ u)
\end{aligned}
}

\AddEq{Vector f(e) cols}
{
\ensuremath{\Vector f\circ \aD[V_1]ei}
}

\AddEq{Vector f(e) rows}
{
\ensuremath{\Vector f\circ \aU{e_{V_1}}i}
}

\AddEq[4]{f=g*h}
{
#1=#2 #3 #4
}

\DefText{define homomorphism A module by matrix(111)}
{
\newline
\FrameEqRef[hg{A_1}{A_2}]{f o ea=efa (11)(\Cols)}{111\SideNS}
\newline
\FrameEqRef[vf{V_1}{V_2}]{f o (ae)=ga o f e, vector space \Product-\Cols}{\SideNS(111)}
\newline
}

\DefText{define homomorphism A module by matrix(11)}
{
\newline
\FrameEqRef[hf{A_1}{A_2}]{f o ea=efa (1)(\Cols)}{11\SideNS}
\newline
\FrameEqRef[vf{V_1}{V_2}]{f o (ae)=ga o f e, vector space \Product-\Cols}{\SideNS(11)}
\newline
}

\DefText{define homomorphism A module by matrix(1)}
{
\newline
\FrameEqRef[vf12]{f o (ae)=a o f e, vector space \Product-\Cols}{\SideNS(1)}
\newline
}

\DefText{define homomorphism D module by matrix(1)}
{
\newline
\FrameEqRef[hf12{}]{f o ea=efa (11)(\Cols)}{module}
\newline
}

\DefText{define homomorphism D module by matrix()}
{
\newline
\FrameEqRef[hf12{}]{f o ea=efa (1)(\Cols)}{module}
\newline
}

\DefRef{b(av)=a(bv)}
{
\newline
\FrameEqRef{associative law, \SideWS module}1
\newline
\FrameEqRef[faav]{\SideWS homomorphism, f av=}{()\SideWS A module}
\newline
}

\DefRef{homomorphism, f av=}
{
\newline
\FrameEqRef[faav]{\SideWS homomorphism, f av=}{()\SideWS A module}
\newline
}

\DefRef{homomorphism, f v+w=}
{
\newline
\FrameEqRef[fuv]{homomorphism, f v+w=}{f()\SideWS A module}
\newline
}

\DefText[3]{f o (ae)=a o f e (11)}
{
\DrawEq[w{}{(\Vector g\circ v)}{}f{}]{v1=v2*a \SideNS-\Cols}{\SideNS(#1#2#3)}
\DrawEq[gf{V_1}{V_2}{}]{f o ev=efv i (11)}{(#1)(\Cols)algebra, \SideNS-module}
\DrawEq[vf{V_1}{V_2}]{f o (ae)=ga o f e, vector space \Product-\Cols}{\SideNS(#1#2#3)}
}

\DefText[3]{f o (ae)=a o f e (1)}
{
\DrawEq[w{}v{}f{}]{v1=v2*a \SideNS-\Cols}{\SideNS(#1#2#3)}
\DrawEq[vf{V_1}{V_2}{}]{f o ev=efv i (1)}{()(\Cols)algebra, \SideNS-module}
\DrawEq[vf12]{f o (ae)=a o f e, vector space \Product-\Cols}{\SideNS(#1#2#3)}
}

\AddEq[4]{f o (ae)=a o f e, vector space cr-cols}
{
\Vector{#2}\circ(#1\CRstar e_{#3})=#1\CRstar #2\CRstar e_{#4}
}

\AddEq[4]{f o (ae)=a o f e, vector space rc-cols}
{
\Vector{#2}\circ(e_{#3}\RCstar #1)=e_{#4}\RCstar #2\RCstar #1
}

\AddEq[4]{f o (ae)=a o f e, vector space cr-rows}
{
\Vector{#2}\circ(e_{#3}\CRstar #1)=e_{#4}\CRstar #2\CRstar #1
}

\AddEq[4]{f o (ae)=a o f e, vector space rc-rows}
{
\Vector{#2}\circ(#1\RCstar e_{#3})=#1\RCstar #2\RCstar e_{#4}
}

\AddEq[4]{f o (ae)=ga o f e, vector space cr-cols}
{
\Vector{#2}\circ(#1\CRstar e_{#3})=(\Vector g\circ #1)\CRstar #2\CRstar e_{#4}
}

\AddEq[4]{f o (ae)=ga o f e, vector space rc-cols}
{
\Vector{#2}\circ(e_{#3}\RCstar #1)=e_{#4}\RCstar #2\RCstar (\Vector g\circ #1)
}

\AddEq[4]{f o (ae)=ga o f e, vector space cr-rows}
{
\Vector{#2}\circ(e_{#3}\CRstar #1)=e_{#4}\CRstar #2\CRstar (\Vector g\circ #1)
}

\AddEq[4]{f o (ae)=ga o f e, vector space rc-rows}
{
\Vector{#2}\circ(#1\RCstar e_{#3})=(\Vector g\circ #1)\RCstar #2\RCstar e_{#4}
}

\AddEq[2]{left homomorphism cr, 1}
{
\Vector {#1}=#1\CRstar e_{#2}
}

\AddEq[2]{left homomorphism rc, 1}
{
\Vector {#1}=#1\RCstar e_{#2}
}

\AddEq[2]{right homomorphism cr, 1}
{
\Vector {#1}=e_{#2}\CRstar #1
}

\AddEq[2]{right homomorphism rc, 1}
{
\Vector {#1}=e_{#2}\RCstar #1
}

\AddEq{vb=a cr f cr e left cr()}
{
\Vector w=v\CRstar f\CRstar e_{V_2}
}

\AddEq{vb=a cr f cr e left rc()}
{
\Vector w=v\RCstar f\RCstar e_{V_2}
}

\AddEq{vb=a cr f cr e right rc()}
{
\Vector w=e_{V_2}\RCstar f\RCstar v
}

\AddEq{vb=a cr f cr e right cr()}
{
\Vector w=e_{V_2}\CRstar f\CRstar v
}

\AddEq{vb=a cr f cr e left cr(1)}
{
\Vector w=(\Vector g\circ v)\CRstar f\CRstar e_{V_2}
}

\AddEq{vb=a cr f cr e left rc(1)}
{
\Vector w=(\Vector g\circ v)\RCstar f\RCstar e_{V_2}
}

\AddEq{vb=a cr f cr e right rc(1)}
{
\Vector w=e_{V_2}\RCstar f\RCstar (\Vector g\circ v)
}

\AddEq{vb=a cr f cr e right cr(1)}
{
\Vector w=e_{V_2}\CRstar f\CRstar (\Vector g\circ v)
}

\AddEq{vv=v*eV left cr()}
{
\begin{aligned}
\Vector w&=\Vector f\circ \Vector v=\Vector f\circ (v\CRstar e_{V_1})
\\ &=v\CRstar (\Vector f\circ e_{V_1})
\end{aligned}
}

\AddEq{vv=v*eV left rc()}
{
\begin{aligned}
\Vector w&=\Vector f\circ \Vector v=\Vector f\circ (v\RCstar e_{V_1})
\\ &=v\RCstar (\Vector f\circ e_{V_1})
\end{aligned}
}

\AddEq{vv=v*eV right rc()}
{
\begin{aligned}
\Vector w&=\Vector f\circ \Vector v=\Vector f\circ (e_{V_1}\RCstar v)
\\ &=(\Vector f\circ e_{V_1})\RCstar v
\end{aligned}
}

\AddEq{vv=v*eV right cr()}
{
\begin{aligned}
\Vector w&=\Vector f\circ \Vector v=\Vector f\circ (e_{V_1}\CRstar v)
\\ &=(\Vector f\circ e_{V_1})\CRstar v
\end{aligned}
}

\AddEq{vv=v*eV left cr(1)}
{
\begin{aligned}
\Vector w&=\Vector f\circ \Vector v=\Vector f\circ (v\CRstar e_{V_1})
\\ &=(\Vector g\circ v)\CRstar (\Vector f\circ e_{V_1})
\end{aligned}
}

\AddEq{vv=v*eV left rc(1)}
{
\begin{aligned}
\Vector w&=\Vector f\circ \Vector v=\Vector f\circ (v\RCstar e_{V_1})
\\ &=(\Vector g\circ v)\RCstar (\Vector f\circ e_{V_1})
\end{aligned}
}

\AddEq{vv=v*eV right rc(1)}
{
\begin{aligned}
\Vector w&=\Vector f\circ \Vector v=\Vector f\circ (e_{V_1}\RCstar v)
\\ &=(\Vector f\circ e_{V_1})\RCstar (\Vector g\circ v)
\end{aligned}
}

\AddEq{vv=v*eV right cr(1)}
{
\begin{aligned}
\Vector w&=\Vector f\circ \Vector v=\Vector f\circ (e_{V_1}\CRstar v)
\\ &=(\Vector f\circ e_{V_1})\CRstar (\Vector g\circ v)
\end{aligned}
}

\AddEq{f o ei=fij ej left cr}
{
\ShowEq{Vector f(e) cols}
=\aD fi\CRstar e_{V_2}
=\aUD fji \aD[V_2]ej
}

\AddEq{f o ei=fij ej left rc}
{
\ShowEq{Vector f(e) rows}
=\aU fi\RCstar e_{V_2}
=\aUD fij \aU{e_{V_2}}j
}

\AddEq{f o ei=fij ej right cr}
{
\ShowEq{Vector f(e) rows}
=e_{V_2}\CRstar \aU fi
=\aD[V_2]ej\aUD fij
}

\AddEq{f o ei=fij ej right rc}
{
\ShowEq{Vector f(e) cols}
=e_{V_2}\RCstar \aD fi
=\aU{e_{V_2}}j\aUD fji
}

\AddEq{ga*g- left-En}
{
\[
\ShowEq{ga*g- \SideNS-\Cols}b{\aD En}
=\aD Enb
\]
}

\AddEq{ga*g- right-En}
{
\[
\ShowEq{ga*g- \SideNS-\Cols}b{\aD En}
=b\aD En
\]
}

\AddEq{g(f-bEn)g left-cols}
{
g\CRstar(
\ShowEq{f-bEn left}1{}
)\CRstar g^{\CRInverse}
}

\AddEq{g(f-bEn)g right-cols}
{
g^{\RCInverse}\RCstar(
\ShowEq{f-bEn right}1{}
)\RCstar g
}

\AddEq{g(f-bEn)g left-rows}
{
g\RCstar(
\ShowEq{f-bEn left}1{}
)\RCstar g^{\RCInverse}
}

\AddEq{g(f-bEn)g right-rows}
{
g^{\CRInverse}\CRstar(
\ShowEq{f-bEn right}1{}
)\CRstar g
}

\AddEq[2]{ga*g- left-cols}
{
\ensuremath{#2#1\CRstar #2^{\CRInverse}}
}

\AddEq[2]{ga*g- right-cols}
{
\ensuremath{#2^{\RCInverse}\RCstar #1 #2}
}

\AddEq[2]{ga*g- left-rows}
{
\ensuremath{#2#1\RCstar #2^{\RCInverse}}
}

\AddEq[2]{ga*g- right-rows}
{
\ensuremath{#2^{\CRInverse}\CRstar #1 #2}
}

\AddEquation{ga*g- 1 left-cols}
{
g\CRstar \aD Ena\CRstar g^{\CRInverse}
}

\AddEquation{ga*g- 1 right-cols}
{
g^{\RCInverse}\RCstar a\aD En\RCstar g
}

\AddEquation{ga*g- 1 left-rows}
{
g\RCstar \aD Ena\RCstar g^{\RCInverse}
}

\AddEquation{ga*g- 1 right-rows}
{
g^{\CRInverse}\CRstar a\aD En\CRstar g
}

\AddEq{(a+b)aE left}
{
\[
\aD En(a+b)=\aD Ena+\aD Enb 
\]
}

\AddEq{(a+b)aE right}
{
\[
(a+b)\aD En=a\aD En+b\aD En
\]
}

\AddEq{(daE)v left}
{
\[
\aD En(ba)=b(\aD Ena)
\]
}

\AddEq{(daE)v right}
{
\[
(ba)\aD En=b(a\aD En)
\]
}

\AddEq{aEn=bEn left}
{
\[
\aD Ena=\aD Enb
\]
}

\AddEq{aEn=bEn right}
{
\[
a\aD En=b\aD En
\]
}

\AddEq[2]{e in basis manifold}
{
$\Basis e_{#1}\in
\ShowEq{basis manifold of V, \SideNS-\Cols}{\aD En}{#2}{}$
}

\AddEq{F1(g)=F(g)**-rc}
{
F_1(g)=F(g)^{\RCInverse}
}

\AddEq{F1(g)=F(g)**-cr}
{
F_1(g)=F(g)^{\CRInverse}
}

\AddEq{coordinate transformation, vector space W, left-cols}
{
w_2=w_1\CRstar F(g)^{\CRInverse}
}

\AddEq{coordinate transformation, vector space W, right-cols}
{
w_2=F(g)^{\RCInverse}\RCstar w_1
}

\AddEq{coordinate transformation, vector space W, left-rows}
{
w_2=w_1\RCstar F(g)^{\RCInverse}
}

\AddEq{coordinate transformation, vector space W, right-rows}
{
w_2=F(g)^{\CRInverse}\CRstar w_1
}

\AddEq[1]{f=f1+f2 endo}
{
f#1=f#1_1+f#1_2
}

\AddEq[1]{left e'=e*a}
{
\DrawEq[{e'_{#1}}{e_{#1}}{\ProductVal}g]{f=g*h}{}
}

\AddEq[1]{right e'=e*a}
{
\DrawEq[{e'_{#1}}g{\ProductVal}{e_{#1}}]{f=g*h}{}
}

\AddEq[1]{f=f1*f2 endo left-cols}
{
f#1=f#1_1\RCstar f#1_2
}

\AddEq[1]{f=f1*f2 endo right-cols}
{
f#1=f#1_1\CRstar f#1_2
}

\AddEq[1]{f=f1*f2 endo left-rows}
{
f#1=f#1_1\CRstar f#1_2
}

\AddEq[1]{f=f1*f2 endo right-rows}
{
f#1=f#1_1\RCstar f#1_2
}

\AddEq[5]{matrix fIJ}
{
\[
#1=(\AUD{#1}{#2}{#4},\jJg{#2}{#3},\jJg{#4}{#5})
\]
}

\DefText[4]{matrices of numbers(11)}
{
\ShowText{matrix of numbers and C}D{#1}gkKlL{#3}{#4}
\ShowText{matrices of numbers(1)}{#1}{#2}{}{}
}

\DefText[4]{matrices of numbers(1)}
{
\ShowText{matrix of numbers}A{#2}fiIjJ
}

\AddEq{passive transformation and endomorphism}
{
\,\newline
\TwoColText
{
\ShowRemark{passive transformation and endomorphism}v
}
{
\ShowRemark{passive transformation and endomorphism}w
}
}

\AddEq{fov=gov cols()}
{
\[
\Vector f\circ(e_1v)=e_2fv =\Vector g\circ(e_1v)
\]
}

\AddEq{fov=gov rows()}
{
\[
\Vector f\circ(ve_1)=vfe_2=\Vector g\circ(ve_1)
\]
}

\AddEq{fov=gov cols(1)}
{
\[
\Vector f\circ(e_1v)=e_2fh(v) =\Vector g\circ(e_1v)
\]
}

\AddEq{fov=gov rows(1)}
{
\[
\Vector f\circ(ve_1)=h(v)fe_2=\Vector g\circ(ve_1)
\]
}

\AddEq{fov=gov left-cols}
{
\[
\Vector f\circ(v\CRstar e_V)=v\CRstar f\CRstar e_W=\Vector g\circ(v\CRstar e_V)
\]
}

\AddEq{fov=gov left-rows}
{
\[
\Vector f\circ(v\RCstar e_V)=v\RCstar f\RCstar e_W=\Vector g\circ(v\RCstar e_V)
\]
}

\AddEq{fov=gov right-cols}
{
\[
\Vector f\circ(e_V\RCstar v)=e_W\RCstar f\RCstar v=\Vector g\circ(e_V\RCstar v)
\]
}

\AddEq{fov=gov right-rows}
{
\[
\Vector f\circ(e_V\CRstar v)=e_W\CRstar f\CRstar v=\Vector g\circ(e_V\CRstar v)
\]
}

\AddEquation{fo(v+w) cols(1)}
{
\begin{aligned}
&\,\Vector f\circ(e_1(v+w))
\\=&\,e_2fh(v+w) 
\\=&\,e_2f(h(v)+h(w)) 
\\=&\,e_2fh(v)+e_2fh(w) 
\\=&\,\Vector f\circ(e_1v)+\Vector f\circ(e_1w)
\end{aligned}
}

\AddEquation{fo(v+w) rows(1)}
{
\begin{aligned}
&\,\Vector f\circ((v+w)e_1)
\\=&\,h(v+w)f e_2 
\\=&\,(h(v)+h(w)) fe_2
\\=&\,h(v)fe_2+h(w)f e_2
\\=&\,\Vector f\circ(ve_1)+\Vector f\circ(we_1)
\end{aligned}
}

\AddEquation{fo(v+w) cols()}
{
\begin{aligned}
&\,\Vector f\circ(e_1(v+w))
\\=&\,e_2f(v+w) 
\\=&\,e_2fv+e_2fw 
\\=&\,\Vector f\circ(e_1v)+\Vector f\circ(e_1w)
\end{aligned}
}

\AddEquation{fo(v+w) rows()}
{
\begin{aligned}
&\,\Vector f\circ((v+w)e_1)
\\=&\,(v+w)f e_2 
\\=&\,vfe_2+wf e_2
\\=&\,\Vector f\circ(ve_1)+\Vector f\circ(we_1)
\end{aligned}
}

\AddEq{fo(v+w) ()left-cols}
{
\begin{aligned}
&\,\Vector f\circ((v+w)\CRstar e_{V_1})
\\=&\,(v+w)\CRstar f\CRstar e_{V_2} 
\\=&\,v\CRstar f\CRstar e_{V_2}+w\CRstar f\CRstar e_{V_2}
\\=&\,\Vector f\circ(v\CRstar e_{V_1})+\Vector f\circ(w\CRstar e_{V_1})
\end{aligned}
}

\AddEq{fo(v+w) ()left-rows}
{
\begin{aligned}
&\,\Vector f\circ((v+w)\RCstar e_{V_1})
\\=&\,(v+w)\RCstar f\RCstar e_{V_2} 
\\=&\,v\RCstar f\RCstar e_{V_2}+w\RCstar f\RCstar e_{V_2}
\\=&\,\Vector f\circ(v\RCstar e_{V_1})+\Vector f\circ(w\RCstar e_{V_1})
\end{aligned}
}

\AddEq{fo(v+w) ()right-cols}
{
\begin{aligned}
&\,\Vector f\circ(e_{V_1}\RCstar(v+w))\\=&\,e_{V_2}\RCstar f\RCstar(v+w)
\\=&\,e_{V_2}\RCstar f\RCstar v+e_{V_2}\RCstar f\RCstar w
\\=&\,\Vector f\circ(e_{V_1}\RCstar v)+\Vector f\circ(e_{V_1}\RCstar w)
\end{aligned}
}

\AddEq{fo(v+w) ()right-rows}
{
\begin{aligned}
&\,\Vector f\circ(e_{V_1}\CRstar(v+w))\\=&\,e_{V_2}\CRstar f\CRstar(v+w)
\\=&\,e_{V_2}\CRstar f\CRstar v+e_{V_2}\CRstar f\CRstar w
\\=&\,\Vector f\circ(e_{V_1}\CRstar v)+\Vector f\circ(e_{V_1}\CRstar w)
\end{aligned}
}

\AddEq{fo(v+w) (1)left-cols}
{
\begin{aligned}
&\,\Vector f\circ((v+w)\CRstar e_{V_1})
\\=&\,(g\circ(v+w))\CRstar f\CRstar e_{V_2} 
\\=&\,(g\circ v)\CRstar f\CRstar e_{V_2}
\\+&\,(g\circ w)\CRstar f\CRstar e_{V_2}
\\=&\,\Vector f\circ(v\CRstar e_{V_1})+\Vector f\circ(w\CRstar e_{V_1})
\end{aligned}
}

\AddEq{fo(v+w) (1)left-rows}
{
\begin{aligned}
&\,\Vector f\circ((v+w)\RCstar e_{V_1})
\\=&\,(g\circ(v+w))\RCstar f\RCstar e_{V_2} 
\\=&\,(g\circ v)\RCstar f\RCstar e_{V_2}
\\+&\,(g\circ w)\RCstar f\RCstar e_{V_2}
\\=&\,\Vector f\circ(v\RCstar e_{V_1})+\Vector f\circ(w\RCstar e_{V_1})
\end{aligned}
}

\AddEq{fo(v+w) (1)right-cols}
{
\begin{aligned}
&\,\Vector f\circ(e_{V_1}\RCstar(v+w))\\=&\,e_{V_2}\RCstar f\RCstar(g\circ(v+w))
\\=&\,e_{V_2}\RCstar f\RCstar(g\circ v)
\\+&\,e_{V_2}\RCstar f\RCstar(g\circ w)
\\=&\,\Vector f\circ(e_{V_1}\RCstar v)+\Vector f\circ(e_{V_1}\RCstar w)
\end{aligned}
}

\AddEq{fo(v+w) (1)right-rows}
{
\begin{aligned}
&\,\Vector f\circ(e_{V_1}\CRstar(v+w))\\=&\,e_{V_2}\CRstar f\CRstar(g\circ(v+w))
\\=&\,e_{V_2}\CRstar f\CRstar(g\circ v)
\\+&\,e_{V_2}\CRstar f\CRstar(g\circ w)
\\=&\,\Vector f\circ(e_{V_1}\CRstar v)+\Vector f\circ(e_{V_1}\CRstar w)
\end{aligned}
}

\AddEquation{fo(va) cols(1)}
{
\begin{aligned}
&\,\Vector f\circ(e_1(av) )\\=&\, e_2f h(av) 
\\=&\,h(a)(e_2f h(v))
\\=&\,h(a)(\Vector f\circ(e_1 v))
\end{aligned}
}

\AddEquation{fo(va) rows(1)}
{
\begin{aligned}
&\,\Vector f\circ((av) e_1)\\=&\,h(av) f e_2
\\=&\,h(a)( h(v)fe_2)
\\=&\,h(a)(\Vector f\circ( ve_1))
\end{aligned}
}

\AddEquation{fo(va) cols()}
{
\begin{aligned}
&\,\Vector f\circ(e_1(av) )\\=&\,e_2f (av)
\\=&\,a(e_2f v)
\\=&\,a(\Vector f\circ(e_1 v))
\end{aligned}
}

\AddEquation{fo(va) rows()}
{
\begin{aligned}
&\,\Vector f\circ((av)e_1)\\=&\,(av) f e_2
\\=&\,a( vfe_2)
\\=&\,a(\Vector f\circ( ve_1))
\end{aligned}
}

\AddEq{fo(va) ()left-cols}
{
\begin{aligned}
&\,\Vector f\circ((av)\CRstar e_{V_1})\\=&\,(av)\CRstar f\CRstar e_{V_2} 
\\=&\,a(v\CRstar f\CRstar e_{V_2})
\\=&\,a(\Vector f\circ(v\CRstar e_{V_1}))
\end{aligned}
}

\AddEq{fo(va) ()left-rows}
{
\begin{aligned}
&\,\Vector f\circ((av)\RCstar e_{V_1})\\=&\,(av)\RCstar f\CRstar e_{V_2} 
\\=&\,a(v\RCstar f\RCstar e_{V_2})
\\=&\,a(\Vector f\circ(v\RCstar e_{V_1}))
\end{aligned}
}

\AddEq{fo(va) ()right-cols}
{
\begin{aligned}
&\,\Vector f\circ(e_{V_1}\RCstar(va))\\=&\,e_{V_2}\RCstar f\RCstar(va)
\\=&\,(e_{V_2}\RCstar f\RCstar v)a
\\=&\,(\Vector f\circ(e_{V_1}\RCstar v))a
\end{aligned}
}

\AddEq{fo(va) ()right-rows}
{
\begin{aligned}
&\,\Vector f\circ(e_{V_1}\CRstar(va))\\=&\,e_{V_2}\CRstar f\CRstar(va)
\\=&\,(e_{V_2}\CRstar f\CRstar v)a
\\=&\,(\Vector f\circ(e_{V_1}\CRstar v))a
\end{aligned}
}

\AddEq{fo(va) (1)left-cols}
{
\begin{aligned}
&\,\Vector f\circ((av)\CRstar e_{V_1})\\=&\,(\Vector g\circ(av))\CRstar f\CRstar e_{V_2} 
\\=&\,(\Vector g\circ a)((\Vector g\circ v)\CRstar f\CRstar e_{V_2})
\\=&\,(\Vector g\circ a)(\Vector f\circ(v\CRstar e_{V_1}))
\end{aligned}
}

\AddEq{fo(va) (1)left-rows}
{
\begin{aligned}
&\,\Vector f\circ((av)\RCstar e_{V_1})\\=&\,(\Vector g\circ(av))\RCstar f\CRstar e_{V_2} 
\\=&\,(\Vector g\circ a)((\Vector g\circ v)\RCstar f\RCstar e_{V_2})
\\=&\,(\Vector g\circ a)(\Vector f\circ(e_{V_1}\CRstar v))
\end{aligned}
}

\AddEq{fo(va) (1)right-cols}
{
\begin{aligned}
&\,\Vector f\circ(e_{V_1}\RCstar(va))\\=&\,e_{V_2}\RCstar f\RCstar(\Vector g\circ(va))
\\=&\,(e_{V_2}\RCstar f\RCstar(\Vector g\circ v))(\Vector g\circ a)
\\=&\,(\Vector f\circ(e_{V_1}\RCstar v))(\Vector g\circ a)
\end{aligned}
}

\AddEq{fo(va) (1)right-rows}
{
\begin{aligned}
&\,\Vector f\circ(e_{V_1}\CRstar(va))\\=&\,e_{V_2}\CRstar f\CRstar(\Vector g\circ(va))
\\=&\,(e_{V_2}\CRstar f\CRstar(\Vector g\circ v))(\Vector g\circ a)
\\=&\,(\Vector f\circ(e_{V_1}\CRstar v))(\Vector g\circ a)
\end{aligned}
}

\AddEquation{e*f*v=e*ae*v left-cols}
{
\begin{aligned}
&\,\RedText{v\CRstar f\CRstar e}
=\Vector f\circ\Vector v
=\Vector{b\Basis e}\circ\Vector v
\\ =&\,\RedText{v\CRstar
\ShowEq{ga*g- \SideNS-\Cols}bg
\CRstar e}
\end{aligned}
}

\AddEquation{e*f*v=e*ae*v right-cols}
{
\begin{aligned}
&\,\RedText{e\RCstar f\RCstar v}
=\Vector f\circ\Vector v
=\Vector{b\Basis e}\circ\Vector v
\\ =&\,\RedText{e\RCstar
\ShowEq{ga*g- \SideNS-\Cols}bg
\RCstar v}
\end{aligned}
}

\AddEquation{e*f*v=e*ae*v left-rows}
{
\begin{aligned}
&\,\RedText{v\RCstar f\RCstar e}
=\Vector f\circ\Vector v
=\Vector{b\Basis e}\circ\Vector v
\\ =&\,\RedText{v\RCstar
\ShowEq{ga*g- \SideNS-\Cols}bg
\RCstar e}
\end{aligned}
}

\AddEquation{e*f*v=e*ae*v right-rows}
{
\begin{aligned}
&\,\RedText{e\CRstar f\CRstar v}
=\Vector f\circ\Vector v
=\Vector{\Basis eb}\circ\Vector v
\\ =&\,\RedText{e\CRstar
\ShowEq{ga*g- \SideNS-\Cols}bg
\CRstar v}
\end{aligned}
}

\AddEq{phantom a**b}
{
\[\vphantom{a^b}\]
}

\AddEquation{twin to left product}
{
(v,a)\rightarrow
\ShowEq{left map a En}a{}
\circ v
}

\AddEquation{twin to right product}
{
(a,v)\rightarrow
\ShowEq{right map a En}a{}
\circ v
}

\AddEq{av left}
{
$va$
}

\AddEq{av right}
{
$av$
}

\AddEquation{(av)b= left}
{
\begin{aligned}
(av)b
&=
\ShowEq{left map a En}b{}
\circ(av)\\ &=a(
\ShowEq{left map a En}b{}
\circ v)=a(vb)
\end{aligned}
}

\AddEquation{(av)b= right}
{
\begin{aligned}
a(vb)&=
\ShowEq{right map a En}a{}
\circ(vb)\\ &=(
\ShowEq{right map a En}a{}
\circ v)b=(av)b
\end{aligned}
}

\AddEq{twin representations of algebra}
{
\xymatrix{
V\ar[rr]^{f(a)}\ar[d]^{g(\Basis e)(b)} & & V\ar[d]^{g(\Basis e)(b)}\\
V\ar[rr]_{f(a)}& &V
}
}

\AddEq{representation Ao2->V}
{
\ShowEq{f:A->*B}h{\AoxA A}V
}

\AddEq{representation Ao2->V =}
{
\[
h(a\otimes b)\circ v=avb
\]
}

\AddEq{(av)b=a(vb)}
{
(av)b=a(vb)
}

\AddEq{ab in A,v in V}
{
$\forall a$, $b\in A$, $v\in V$
}

\AddEq{be=e cr En right}
{
\[
\Basis e=e\CRstar\Basis{\aD En}
\]
}

\AddEq{be=e cr En left}
{
\[
\Basis e=\Basis{\aD En}\RCstar e
\]
}

\AddEq{a in A->aE hom left}
{
a\in A\rightarrow \aD Ena\in\End(A,V^*)
}

\AddEq{a in A->aE hom right}
{
a\in A\rightarrow a\aD En\in\End(A,V^*)
}

\AddEq{(ab)E left-cols}
{
\[
\aD En(ab)=(\aD Ena)\CRstar(\aD Enb)
\]
}

\AddEq{(ab)E right-cols}
{
\[
((ab)\aD En)=(a\aD En)\RCstar(b\aD En)
\]
}

\AddEq{(ab)E left-rows}
{
\[
\aD En(ab)=(\aD Ena)\RCstar(\aD Enb)
\]
}

\AddEq{(ab)E right-rows}
{
\[
((ab)\aD En)=(a\aD En)\CRstar(b\aD En)
\]
}

\AddEq{dim V=n}
{
$\dim V=\gin$
}

\AddEq{left a En}
{
$\aD Ena$
}

\AddEq{right a En}
{
$a\aD En$
}

\AddEq{v in V -> av in V right-cols}
{
\[
\Vector{a\Basis e}:e\RCstar v\in V\rightarrow
e\RCstar (a\aD En)\RCstar v\in V
\]
}

\AddEq{v in V -> av in V left-cols}
{
\[
\Vector{\Basis ea}:v\CRstar e\in V\rightarrow
v\CRstar (\aD Ena)\CRstar e\in V
\]
}

\AddEq{v in V -> av in V right-rows}
{
\[
\Vector{a\Basis e}:e\CRstar v\in V\rightarrow
e\CRstar (a\aD En)\CRstar v\in V
\]
}

\AddEq{v in V -> av in V left-rows}
{
\[
\Vector{\Basis ea}:v\RCstar e\in V\rightarrow
v\RCstar (\aD Ena)\RCstar e\in V
\]
}

\AddEq{basis e1 e2}
{
$\Basis e_1$, $\Basis e_2$
}

\AddEq{a1n= b= column}
{
\[
\aD a1=
\ColMatrix{\aD a1}m
\ \ \ ...\ \ \ \,
\aD an=
\ColMatrix{\aD an}m
\ \ \ \,
b=
\ColMatrix bm
\]
}

\AddEquation{star rows system of linear equations 1}
{
\begin{pmatrix}
\aUD a11 &...&\aUD a1n\\
... & ... & ... \\
\aUD am1 &...&\aUD amn
\end{pmatrix}
\RCstar
\begin{pmatrix}
\aU x1\\...\\ \aU xn
\end{pmatrix}
=
\begin{pmatrix}
\aU b1\\...\\ \aU bm
\end{pmatrix}
}

\AddEq{star rows system of linear equations 2}
{
\begin{array}{cccc}
\aUD a11\aU x1 &+...&+\aUD a1n\aU xn&=\aU b1\\
... & ... & ...& ... \\
\aUD am1\aU x1 &+...&+\aUD amn\aU xn&=\aU bm
\end{array}
}

\AddEquation{star rows system of linear equations, rank}
{
\RCRank(\aUD aji)=\RCRank
\begin{pmatrix}
\aUD aji&\aU bj
\end{pmatrix}
}

\AddEquation{star rows system of linear equations}
{
\aUD aji\aU xi=\aU bj
}

\AddEq{aD1n}
{
$\aD a1$, ..., $\aD an$.
}

\AddEq[1]{rc-rank a=k<m}
{
\[\RCRank a=\gik<\gi{#1}\]
}

\AddEq[1]{cr-rank a=k<m}
{
\symb{\CRRank  a}{cr-rank of matrix}{}
$$\ShowSymbol{cr-rank of matrix}{}=\gik<\gi{#1}$$
}

\AddEq{rc-rank of matrix}
{
\symb{\RCRank a}{rc-rank of matrix}1,
}

\AddEq{cr-rank of matrix}
{
\symb{\CRRank a}{cr-rank of matrix}1,
}

\AddEq[2]{left a rc b}
{
#1\ProductVal #2
}

\AddEq[2]{right a rc b}
{
#2\ProductVal #1
}

\AddEq{rows-of-matrix linearly dependent, 1}
{
$\aD {\lambda}p=-1$, $\aD {\lambda}s=\pRs$
}

\AddEq{rows-of-matrix linearly dependent, 2}
{
$\aD {\lambda}c=0$.
}

\AddEq[2]{columns of matrix are linearly dependent}
{
\ShowEq{\SideWS a rc b}{#1}{#2}=0
}

\AddEq{cols-of-matrix linearly dependent, 1}
{
$\aU {\lambda}r=-1$, $\aU {\lambda}t=\aUD Rtr$
}

\AddEq{cols-of-matrix linearly dependent, 2}
{
$\aU {\lambda}c=0$.
}

\AddEq{representation left vector space of maps B->A}
{
\[
\xymatrix@C=30pt
{
f:A\ar[r]|-{*}&A^B
}
\ \ \ f(a):g\in A^B\rightarrow ag\in A^B
\ \ \ (ag)(b)=ag(b)
\]
}

\AddEq{representation right vector space of maps B->A}
{
\[
\xymatrix@C=30pt
{
f:A\ar[r]|-{*}&A^B
}
\ \ \ f(a):g\in A^B\rightarrow ga\in A^B
\ \ \ (ga)(b)=g(b)a
\]
}

\AddEq{representation left vector space of algebra A}
{
\[
\xymatrix@C=30pt
{
f:A\ar[r]|-{*}&A
}
\ \ \ f(a):b\in A\rightarrow ab\in A
\]
}

\AddEq{representation right vector space of algebra A}
{
\[
\xymatrix@C=30pt
{
f:A\ar[r]|-{*}&A
}
\ \ \ f(a):b\in A\rightarrow ba\in A
\]
}

\AddEquation{left a->ab}
{
f(a):b\in A\rightarrow ab\in A
}

\AddEquation{right a->ab}
{
f(a):b\in A\rightarrow ba\in A
}

\AddEq{left vector space of maps B->A n}
{
\symb{\mathcal L^nA^B}{left vector space of maps B->A}1
}

\AddEq{right vector space of maps B->A n}
{
\symb{\mathcal R^nA^B}{right vector space of maps B->A}1
}

\AddEq{left vector space of algebra A n}
{
\symb{\mathcal L^nA}{left vector space of algebra A}1
}

\AddEq{right vector space of algebra A n}
{
\symb{\mathcal R^nA}{right vector space of algebra A}1
}

\AddEq{left vector space of maps B->A k}
{
\begin{align*}
\mathcal L^1A^B&=\mathcal LA^B\\
\mathcal L^kA^B&=\mathcal L^{k-1}A^B\oplus\mathcal LA^B
\end{align*}
}

\AddEq{right vector space of maps B->A k}
{
\begin{align*}
\mathcal R^1A^B&=\mathcal RA^B\\
\mathcal R^kA^B&=\mathcal R^{k-1}A^B\oplus\mathcal RA^B
\end{align*}
}

\AddEq{left vector space of algebra A k}
{
\begin{align*}
\mathcal L^1A&=\mathcal LA\\
\mathcal L^kA&=\mathcal L^{k-1}A\oplus\mathcal LA
\end{align*}
}

\AddEq{right vector space of algebra A k}
{
\begin{align*}
\mathcal R^1A&=\mathcal RA\\
\mathcal R^kA&=\mathcal R^{k-1}A\oplus\mathcal RA
\end{align*}
}

\AddEq{left vector space of maps B->A}
{
\symb{\mathcal LA^B}{left vector space of maps B->A}1
}

\AddEq{right vector space of maps B->A}
{
\symb{\mathcal RA^B}{right vector space of maps B->A}1
}

\AddEq{left vector space of algebra A}
{
\symb{\mathcal LA}{left vector space of algebra A}1
}

\AddEq{right vector space of algebra A}
{
\symb{\mathcal RA}{right vector space of algebra A}1
}

\AddEq{left a(h+g)=}
{
\[
a(h(b)+g(b))=ah(b)+ag(b)
\]
}

\AddEq{right a(h+g)=}
{
\[
(h(b)+g(b))a=h(b)a+g(b)a
\]
}

\AddEq{left a(b+c)=}
{
\[
a(b+c)=ab+ac
\]
}

\AddEq{right a(b+c)=}
{
\[
(b+c)a=ba+ca
\]
}

\AddEq{left (a1+a2)h=}
{
\[
(a_1+a_2)h(b)=a_1h(b)+a_2h(b)
\]
}

\AddEq{right (a1+a2)h=}
{
\[
h(b)(a_1+a_2)=h(b)a_1+h(b)a_2
\]
}

\AddEq{left a1a2h=}
{
\[
(a_2a_1)h(b)=a_2(a_1h(b))
\]
}

\AddEq{right a1a2h=}
{
\[
h(b)(a_1a_2)=(h(b)a_1)a_2
\]
}

\AddEq{left (a1+a2)b=}
{
\[
(a_1+a_2)b=a_1b+a_2b
\]
}

\AddEq{right (a1+a2)b=}
{
\[
b(a_1+a_2)=ba_1+ba_2
\]
}

\AddEq{left a1a2b=}
{
\[
(a_2a_1)b=a_2(a_1b)
\]
}

\AddEq{right a1a2b=}
{
\[
b(a_1a_2)=(ba_1)a_2
\]
}

\AddEquation{rank of matrix, rc-cols}
{
a_{\gi N\setminus\gi T}=\aD aT\RCstar R
}

\AddEquation{rank of matrix, 1, rc-cols}
{
\aD ar=\aD aT\RCstar \aD Rr
}

\AddEquation{rank of matrix, 2, rc-cols}
{
\aUD aar=\aUD aat\ \tRr
}

\AddEquation{rank of matrix, rc-rows}
{
a^{\gi M\setminus\gi S}=R\RCstar \aU aS
}

\AddEquation{rank of matrix, 1, rc-rows}
{
\aU ap=\aU Rp\RCstar\aU aS
}

\AddEquation{rank of matrix, 2, rc-rows}
{
\aUD apb=\aUD Rps\aUD asb
}

\AddEquation{rank of matrix, cr-rows}
{
a^{\gi M\setminus\gi S}=\aU aS\CRstar R
}

\AddEquation{rank of matrix, 1, cr-rows}
{
\aU ap=\aU aS\CRstar\aU Rp
}

\AddEquation{rank of matrix, 2, cr-rows}
{
\aUD abp=\aUD abs\aUD Rsp
}

\AddEquation{rank of matrix, cr-cols}
{
a_{\gi N\setminus\gi T}=R\CRstar\aD aT
}

\AddEquation{rank of matrix, 1, cr-cols}
{
\aD ar=\aD Rr\CRstar\aD aT
}

\AddEquation{rank of matrix, 2, cr-cols}
{
\aUD air=\aUD Rtr\aUD ait
}

\AddEq{rank of matrix, 4, rc}
{
\[
\pA r
-\pA T\RCstar
\SATm\RCstar \SA r=0
\]
}

\AddEq{rank of matrix, 4, cr}
{
\[
\pA r
-\SA r\CRstar
\SATm\CRstar\pA T =0
\]
}

\AddEquation{rank of matrix, 5, rc-cols}
{
\aD Rr
=\SATm\RCstar \SA r
}

\AddEquation{rank of matrix, 5, rc-rows}
{
\aU Rp
=\pA T\RCstar\SATm
}

\AddEquation{rank of matrix, 5, cr-cols}
{
\aD Rr
=\SA r\CRstar\SATm
}

\AddEquation{rank of matrix, 5, cr-rows}
{
\aU Rp
=\SATm\CRstar\pA T
}

\AddEquation{rank of matrix, 6, rc-cols}
{
\pA r=\pA T\RCstar\aD Rr
}

\AddEquation{rank of matrix, 6, rc-rows}
{
\pA r=\aU Rp\RCstar \SA r
}

\AddEquation{rank of matrix, 6, cr-cols}
{
\pA r=\aD Rr\CRstar\pA T
}

\AddEquation{rank of matrix, 6, cr-rows}
{
\pA r=\SA r\CRstar \aU Rp
}

\AddEq{rank of matrix, 7, rc-cols}
{
\[
\aUD akT\RCstar\SATm
=\aUD {\delta}kS
\]
}

\AddEq{rank of matrix, 7, rc-rows}
{
\[
\SATm\RCstar \SA l
=\Tdl
\]
}

\AddEq{rank of matrix, 7, cr-cols}
{
\[
\SATm\CRstar\aUD akT
=\aUD {\delta}kS
\]
}

\AddEq{rank of matrix, 7, cr-rows}
{
\[
\SA l\CRstar\SATm
=\Tdl
\]
}

\AddEquation{rank of matrix, 8, rc-cols}
{
\begin{split}
\aUD akr&=
\aUD {\delta}kS\RCstar\SA r
\\ &=
\aUD akT\RCstar\SATm
\RCstar
\SA r
\end{split}
}

\AddEquation{rank of matrix, 8, rc-rows}
{
\begin{split}
\pA l&=
\pA T\RCstar \Tdl
\\ &=
\pA T\RCstar
\SATm\RCstar \SA l
\end{split}
}

\AddEquation{rank of matrix, 8, cr-cols}
{
\begin{split}
\aUD akr&=
\SA r\CRstar\aUD {\delta}kS
\\ &=
\SA r
\CRstar
\SATm\CRstar\aUD akT
\end{split}
}

\AddEquation{rank of matrix, 8, cr-rows}
{
\begin{split}
\pA l&=
\Tdl\CRstar\pA T
\\ &=
\SA l\CRstar\SATm
\CRstar
\pA T
\end{split}
}

\AddEquation{rank of matrix, 9, rc-cols}
{
\aUD akr=
\aUD akT\RCstar \aD Rr
}

\AddEquation{rank of matrix, 9, rc-rows}
{
\pA l=
\aU Rp\RCstar \SA l
}

\AddEquation{rank of matrix, 9, cr-cols}
{
\aUD akr=
\aD Rr\CRstar \aUD akT
}

\AddEquation{rank of matrix, 9, cr-rows}
{
\pA l=
\SA l\CRstar \aU Rp
}

%% file: Biring.Stmt.English.tex
\input{Biring.Stmt.Eq}

\DefTheorem{transpose of identity}
{
Transpose of identity is identity
\DrawEq{transpose of identity}E
}

\DefProof{transpose of identity}
{
The theorem follows from definitions
\RefDefinition{transpose of matrix},
\RefDefinition{biring matrix}.
}

\DefDefinition{biring matrix}
{
The set of matrices $\mathcal A$ is a \AddIndex{biring}{biring}
if we defined on $\mathcal A$ an unary operation, say transpose,
and three binary operations,
say \RC product, \CR product and sum, such that
\begin{itemize}
\item \RC product and sum define structure of ring on $\mathcal A$
\item \CR product and sum define structure of ring on $\mathcal A$
\item both products have common identity
\ShowEq{delta n=Matrix}
\end{itemize}
}

\DefLabeledTheorem{product of matrices is distributive over addition}{\Product}
{
The \ProductType product in $\Omega$\Hyph ring
{\bf is distributive}
over addition
\DrawEq{a.\Product.(b1+b2)=}1
\DrawEq{(b1+b2).\Product.a=}1
}

\DefProof{product of matrices is distributive over addition}
{
The equality
\ShowEq{a.\Product.(b1+b2)= 1}
follows from equalities
\ShowEq{a.Product.(b1+b2)= 2}
The equality
\eqRef{a.\Product.(b1+b2)=}1
follows from the equality
\EqRef{a.\Product.(b1+b2)= 1}.
The equality
\ShowEq{(b1+b2).\Product.a= 1}
follows from equalities
\ShowEq{(b1+b2).Product.a= 2}
The equality
\eqRef{(b1+b2).\Product.a=}1
follows from the equality
\EqRef{(b1+b2).\Product.a= 1}.
}

\DefSummary{power matrix}
{
For each product we introduce power
\ShowEq{show power matrix}
as well as the inverse matrix
\ShowEq{show inverse matrix}
}

\DefLabeledDefinition{power matrix}{\ProductS}
{
We introduce \AddIndex{\ProductType power}{\ProductS power} of
the \nTimes matrix $a$
using recursive definition
\ShowEq{\ProductS power}
\DrawEq{\ProductS power, 0}1
\DrawEq{\Product-power, n}1
}

\AddEq{remark: generalized index structure}
{
We will use the symbol $\cdot$ in
front of a generalized index when we need to describe its structure.
we put the sign $'-'$ in place of the index whose 
position was changed.
For instance, if an original term was $\aD a{ij}$
I will use notation $a_{\gii\cdot}^{}{}^{\gij}_-$ instead of notation $\aUD aji$.
}

\DefText{range of index}
{
Even though the structure of a generalized index is arbitrary we assume that
there exists a one-to-one map of the interval of positive integers
\ShowEq{gi 1...n}
to the range of index. Let $\gi I$ be the range
of the index $\gii$. We denote the power of this set by symbol
$|\gi I|$ and assume that $|\gi I|=\gin$. If we want to
enumerate elements $\aD ai$ we use notation
\ShowEq{a gi 1...n}
}

\DefText{submatrix}
{
Representation of coordinates of a vector as a matrix allows
making a notation more compact. The question of the presentation
of vector as a row or a column of the matrix is just a question of convention.
We extend the concept of generalized index to entries of the matrix.
A matrix is a two dimensional table, the rows and columns
of which are enumerated by generalized indeices.
Since we use generalized index,
we cannot tell whether index $a$ of
matrix enumerates rows or columns
until we know the structure of index.
However as we can see
below the form of presentation of matrix is not important for us.
To make sure that notation offered below
is consistent with the traditional
we will assume that the matrix is presented in the form
\DrawEq[anm]{a=(a11.nm matrix)}{}
The upper index enumerates
rows and the lower index enumerates columns.
}

\DefDefinition{notation for submatrix}
{
I use the following names and notation for different
\AddIndex{submatrices}{submatrix}
of the matrix $a$
\begin{itemize}
\item
\AddIndex{column of matrix}{column of matrix}
with the index $\gii$
\ShowEq{column of matrix}
The upper index enumerates entries
of the column and the lower index enumerates columns.
\begin{description}
\item
[\ShowEq{A from columns T}\ \ ]
the submatrix obtained from the matrix $a$ by selecting
columns with an index from the set $\gi T$
\item [\ShowEq{A without column a}\ \ ]
the submatrix obtained from the matrix $a$ by deleting
column $\aD ai$
\item [\ShowEq{A without columns T}\ \ ]
the submatrix obtained from the matrix $a$ by deleting
columns with an index from the set $\gi T$
\end{description}
\item
\AddIndex{row of matrix}{row of matrix}
with the index $\gij$
\ShowEq{row of matrix}
The lower index enumerates entries
of the row and the upper index enumerates rows.
\begin{description}
\item [\ShowEq{A from rows S}\ \ ]
the submatrix obtained from the matrix $a$ by selecting rows with an index
from the set $\gi S$
\item [\ShowEq{A without row b}\ \ ]
the submatrix obtained from the matrix $a$ by deleting
row $\aU aj$
\item [\ShowEq{A without rows S}\ \ ]
the submatrix obtained from the matrix $a$ by deleting rows with an index
from the set $\gi S$
\end{description}
\end{itemize}
}

\DefText{interval [kl]}
{
Let $\gik$, $\gil$, $\gik\le\gil$, be integers.
We will consider the set
\ShowEq{interval [kl]}
}

\DefText{Kronecker symbol}
{
Let
\ShowEq{I,|I|=n}
be a set of indices.
We introduce the
\AddIndex{Kronecker symbol}{Kronecker symbol}
\ShowEq{Kronecker symbol}
\ShowEq{Kronecker symbol=}
}

\DefRemark{combine the notation of indices}
{
We will combine the notation of indices. 
Thus \ShowEq{A from b a} is \nTimes[1] submatrix.
The same time this is the notation for a matrix entry.
This allows an identifying of \nTimes[1] matrix and its entry.
The index $a$ is number of
column of matrix and
the index $b$ is number of row of matrix.
}

\DefLabeledTheorem{singular matrix and quasideterminant}{\Product}
{
Let matrix $a$ be \ProductType singular matrix and submatrix $\SA T$
be \ProductType major submatrix. Then
\DrawEq{singular matrix and quasideterminant}{}
}

\DefText{proof of singular matrix and quasideterminant}
{
To understand why submatrix
\DrawEq{singular matrix and quasideterminant, 0}{\Product}
does not have \ProductType inverse matrix,\,\footnote{It is natural
to expect relationship between \ProductType singularity of the matrix
and its \ProductType quasideterminant
similar to relationship which is known in commutative case.
However \ProductType quasideterminant is defined not always.
For instance, it is not defined when \ProductType inverse matrix
has too much elements equal $0$.
As it follows from this theorem, the \ProductType quasideterminant
is undefined also in case when \ProductType rank of the matrix is less then $n-1$.}
we assume that there exists \ProductType inverse matrix
\ShowEq{B rc-1}
We write down the system of linear equations
\newline
\FrameEqRef{row [J] * col J =0}{\Product}
\newline
\FrameEqRef{row J * col J =E}{\Product}
\newline
in case
\ShowEq{I= J= \Product}
(then
\ShowEq{[I]= [J]= \Product})
\DrawEq{singular matrix and quasideterminant, 1}{\Product}
\DrawEq{singular matrix and quasideterminant, 2}{\Product}
and will try to solve this system.
We multiply
\eqRef{singular matrix and quasideterminant, 1}{\Product}
by
$\Bf$
\DrawEq{singular matrix and quasideterminant, 3}{\Product}
Now we can substitute
\eqRef{singular matrix and quasideterminant, 3}{\Product}
into
\eqRef{singular matrix and quasideterminant, 2}{\Product}
\DrawEq{singular matrix and quasideterminant, 4}{\Product}
From
\eqRef{singular matrix and quasideterminant, 4}{\Product}
it follows that
\DrawEq{singular matrix and quasideterminant, 5}{\Product}
The equality
\DrawEq{singular matrix and quasideterminant, 6}{\Product}
follows from the equality
\newline
\FrameEqRef{j i quasideterminant =}{\Product}
\newline
and from the equality
\eqRef{singular matrix and quasideterminant, 5}{\Product}.
Thus we proved that quasideterminant $\Bq$
is defined and its equation to $0$ is necessary and sufficient condition
that the matrix $b$ is singular.
}

\DefLabeledTheorem{rank of coordinate matrix}{\Product-\Cols}
{
Let $V$ be \SideWS $A$\Hyph vector space of \ColN s.
Let
\ShowEq{set of m vectors}
be set of vectors
linearly independent from left. Then \ProductType rank of their
coordinate matrix equal $\gim$.
}

\DefProof{rank of coordinate matrix}
{
Let \(\Basis e\) be the basis of \SideWS $A$\Hyph vector space $V$.
According to the theorem
\refTheorem{coordinate matrix of vector}{\SideNS-\Cols},
the coordinate matrix of set of vectors $(\aU {\Vector a}i)$ relative basis $\Basis e$
consists from \ColN s which are coordinate matrices of vectors
$\aU {\Vector a}i$ relative the basis $\Basis e$. Therefore \ProductType rank of
this matrix cannot be more then $\gim$.

Let \ProductType                                                                                                                                    rank of
the coordinate matrix be less then $\gim$.
According to the theorem
\refTheorem{columns of matrix are linearly dependent}{\Product-\Cols},
\ColN s of matrix are \SideWS linear dependent
\DrawEq[{\lambda}a]{columns of matrix are linearly dependent}{\Product-\Cols}
Suppose
\DrawEq[{\lambda}ac]{c=a rc b}-.
From the equation
\eqRef{columns of matrix are linearly dependent}{\Product-\Cols}
it follows that \SideWS linear composition of rows
\DrawEq[c{\Vector e}]{columns of matrix are linearly dependent}{}
of vectors of basis equal $0$. This contradicts to statement that vectors $\Vector e$
form basis. We proved statement of theorem.
}

\DefLabeledDefinition{Inverse matrix nm}{\Product-\SideNS}
{
Let $a$ be \nmTimes nm matrix
\DrawEq[anm]{a=(a11.nm matrix)}{}
and
\ShowEq{gi k=min nm}
The matrix $a$ is called
\SideWS \ProductType regular
if there exists \nmTimes mn matrix $b$
\DrawEq[bmn]{a=(a11.nm matrix)}{}
such that
\ProductType product
\DrawEq[bac]{c=a rc b}{Inverse \Product-\SideNS}
has submatrix
\ShowEq{cIJ=Ek}
and other elements of the matrix $c$ are $0$.
The matrix $b$ is called
\SideWS \ProductType inverse
for the matrix $a$.
}

\DefLabeledTheorem{Inverse matrix rank}{\Product-\SideNS}
{
Let
\ShowEq{gi k=min nm}
The matrix $a$ has
\SideWS \ProductType inverse matrix $b$
iff
\DrawEq[k]{rank a=gi n}{k \Product-\SideNS}
}

\DefProof[5]{Inverse matrix rank}
{
Let
\ShowEq{I=[1k]}{#3}
\begin{sloppypar}
Let $\gi{#2}$ be set of indices
of \SideWS linear independent \RowsRWSA of matrix $a$.
If $|\gi{#2}|<\gik$,
then we will complement the set $\gi{#2}$
with indices of the remaining \RowsRWSA of matrix $a$
in such a way as to get $|\gi{#2}|=\gik$.
According to the definition
\ShowEq{ref definition product \Product},
the equality
\end{sloppypar}
\DrawEq[{#2}{#3}]{aI rc bJ=cIJ}{\Product-\SideNS}
follows from the equality
\newline
\FrameEqRef[bac]{c=a rc b}{Inverse \Product-\SideNS}
\newline
According to definitions
\ShowEq{ref definition product \Product},
\refDefinition{Inverse matrix nm}{\Product-\SideNS},
For each index \jJg {#4}{#3},
the system of linear equations
\DrawEq[{#1}{#5}{#4}]{linear equations rank \SideNS}{\Product}
follows from the equality
\eqRef{aI rc bJ=cIJ}{\Product-\SideNS}.
We assume
that unknown variables are
\ShowEq{xj=bij}{#4}{#5}

If we consider all systems of linear equations
\eqRef{linear equations rank \SideNS}{\Product},
then their augmented matrix will have the form
\ShowEq{linear equations for Inverse augmented \Cols}{#1}{#5}
Rank of augmented matrix
is $\gik$.
According to the theorem
\RefTheorem{star rows system of linear equations},
system of linear equations
\eqRef{linear equations rank \SideNS}{\Product}
have solution iff condition
\eqRef{rank a=gi n}{k \Product-\SideNS}.
is met.
In particular, the set of \RowsRWSA
\ShowEq{set vi Acols}a{#5}{#2}
is \SideWS linear independent.
}

\DefLabeledDefinition[1]{Inverse Matrix set}{\Product-\SideNS}
{
The set
\ShowEq{set of inverse matrices}{#1}
is the set of
\SideWS \ProductType inverse matrices
for the matrix $a$.
}

\DefLabeledTheorem{a in left = b in right}{\Product}
{
\ShowEq{a in left}ab{\SideNS}{}
iff
\ShowEq{a in left}ba{\OtherSideNS}.
}

\DefProof{a in left = b in right}
{
The theorem follows from definitions
\ShowEq{ref for a in left}
}

\DefLabeledDefinition[4]{divisor of matrix}{\Product-\SideNS}
{
Let $a$ be \nmTimes nm matrix
\DrawEq[anm]{a=(a11.nm matrix)}{}
and $c$ be \nmTimes {#1}{#2} matrix
\DrawEq[c{#1}{#2}]{a=(a11.nm matrix)}{}
if there exists \nmTimes {#3}{#4} matrix $b$
\DrawEq[b{#3}{#4}]{a=(a11.nm matrix)}{}
such that
\DrawEq[abc]{c=a rc b}{}
and other elements are $0$.
}

\DefLabeledTheorem{quasideterminant and inverse}{\Product}
{
Expression for \ProductType inverse matrix has form
\DrawEq{quasideterminant and inverse}{\Product}
\ShowEq{quasideterminant definition= H}
}

\DefProof{quasideterminant and inverse}
{
The equality
\eqRef{quasideterminant and inverse}{\Product}
follows from the equality
\DrawEq{j i quasideterminant H}{\Product}
}

\DefProof{inverse product of matrix over scalar}
{
To prove the equality
\eqRef{inverse product of matrix over scalar, 1}{\Product}
we proceed by induction on size of the matrix.

Since
\ShowEq{inverse product of matrix over scalar, 10}
the statement is evident for \nTimes[1] matrix.

Let the statement holds for \nTimes[(n-1)] matrix. Then from the equality
\eqRef{inverse minor}{\Product},
it follows that
\ShowEq{inverse product of matrix over scalar, 1, 2}
\DrawEq{inverse product of matrix over scalar, 1, 1}{\Product}
The equality
\eqRef{inverse product of matrix over scalar, 1}{\Product}
follows from the equality
\eqRef{inverse product of matrix over scalar, 1, 1}{\Product}.
In the  same manner we prove the equality
\eqRef{inverse product of matrix over scalar, 2}{\Product}.
}

\DefProof{transpose power}
{
We proceed by induction on $n$.

\begin{sloppypar}
For $n=0$ the statement immediately follows from equalities
\newline
\FrameEqRef{rc power, 0}1
\newline
\FrameEqRef{cr power, 0}1
\newline
\FrameEqRef{transpose of identity}E
\end{sloppypar}

Suppose the statement of theorem holds when $n=k-1$
\DrawEq[{k-1}{}]{transpose power}{k-1 \Product}
It follows from
\newline
\FrameEqRef{\Product-power, n}1
\newline
that
\DrawEq{transpose power, 1}{\Product}
It follows from
\eqRef{transpose power, 1}{\Product}
and
\eqRef{transpose power}{k-1 \Product}
that
\DrawEq{transpose power, 2}{\Product}
It follows from
\eqRef{transpose power, 2}{\Product}
and
\eqRef{cr transpose}{Theorem}
that
\DrawEq{transpose power, 3}{\Product}
\eqRef{transpose power}{\Product}
follows
from
\eqRef{transpose power, 1}{\Product}
and
\eqRef{\ProductA-power, n}1.
}

\DefLabeledTheorem{transpose inverse}{\Product}
{
Let \nTimes matrix $a$ have \ProductType inverse matrix.
Then transpose
matrix $a^T$ has \ProductTypeA inverse matrix and these
matrices satisfy the equality
\DrawEq{transpose inverse}{\Product}
}

\DefProof{transpose inverse}
{
If we get transpose of both side
\eqRef{\Product-inverse matrix}1
and apply
\newline
\FrameEqRef{transpose of identity}E
\newline
we get
\ShowEq{transpose inverse, 0}
Applying
\newline
\FrameEqRef{\ProductS transpose}{Theorem}
\newline
we get
\DrawEq{transpose inverse, 1}{\Product}
\eqRef{transpose inverse}{\Product} follows from comparison
\eqRef{cr-inverse matrix}1 and
\eqRef{transpose inverse, 1}{\Product}.
}

\DefLabeledTheorem{inverse unit matrix order 2}{\Product}
{
Let
\DrawEq{unit matrix order 2}{\Product}
Then
\DrawEq{inverse matrix, unit}{\Product}
}

\DefProof{inverse unit matrix order 2}
{
It is clear from
\eqRef{quasideterminant, matrix 2x2 11}{\Product},
and
\eqRef{quasideterminant, matrix 2x2 22}{\Product}
that
\ShowEq{inverse matrix, unit 1 1}
and
\ShowEq{inverse matrix, unit 2 2}
However expression for
\ShowEq{inverse matrix, unit 1 2}
and
\ShowEq{inverse matrix, unit 2 1}
cannot be defined from
\eqRef{quasideterminant, matrix 2x2 12}{\Product}
and
\eqRef{quasideterminant, matrix 2x2 21}{\Product}
since
\ShowEq{inverse matrix, unit 01}
We can transform these expressions. For instance
\ShowEq{inverse matrix, unit 02}
It follows immediately that
\ShowEq{inverse matrix, unit 03}
In the same manner we can find that
\ShowEq{inverse matrix, unit 04}
This completes the proof of
\eqRef{inverse matrix, unit}{\Product}.
}

\DefLabeledTheorem{quasideterminant, expression}{\Product}
{
Expression for
\QDef\hyph \ProductType quasideterminant
has the following form
\DrawEq{quasideterminant, expression}{\Product}
}

\DefProof{quasideterminant, expression}
{
The equality
\eqRef{quasideterminant, expression}{\Product}
follows from the equalities
\eqRef{inverse matrix, Hadamard}{\Product},
\eqRef{j i quasideterminant =}{\Product}.
}

\DefLabeledTheorem{rdet=a--}{\Product}
{
\DrawEq{rdet=a--}{\Product}
\DrawEq{rdet-=a-}{\Product}
}

\DefProof{rdet=a--}
{
The equality
\eqRef{rdet=a--}{\Product}
follows from the equality
\eqRef{j i quasideterminant =}{\Product}.
The equality
\eqRef{rdet-=a-}{\Product}
follows from the equality
\eqRef{rdet=a--}{\Product}.
}

\DefLabeledTheorem{quasideterminant, matrix 2x2}{\Product}
{
Consider matrix
\ShowEq{inverse matrix 2x2, 0}
Then
\DrawEq{quasideterminant, matrix 2x2}{\Product}
\DrawEq{inverse matrix 2x2}{\Product}
}

\DefProof{quasideterminant, matrix 2x2}
{
According to the equality
\eqRef{j i quasideterminant =}{\Product}
\DrawEq{quasideterminant, matrix 2x2 11}{\Product}
\DrawEq{quasideterminant, matrix 2x2 21}{\Product}
\DrawEq{quasideterminant, matrix 2x2 12}{\Product}
\DrawEq{quasideterminant, matrix 2x2 22}{\Product}
The equality
\eqRef{quasideterminant, matrix 2x2}{\Product}
follows from equalities
\ShowEq{ref quasideterminant matrix}
}

\DefLabeledDefinitionNote{quasideterminant}{\Product}
{
\AddIndex{\QDef\hyph \ProductType quasideterminant}{j i \Product-quasideterminant}
of \nTimes matrix $a$
is formal expression\,\footnotemark
\ShowEq{j i quasideterminant definition}
\DrawEq{j i quasideterminant =}{\Product}
We consider \QDef\hyph \ProductType quasideterminant
as an entry of the matrix
\ShowEq{quasideterminant definition}
\DrawEq{quasideterminant definition=}{\Product}
which is called
\AddIndex{\ProductType quasideterminant}{\Product-quasideterminant}.
}
{
In the definition
\refDefinition{quasideterminant}{\Product},
I follow the definition
\citeBib{math.QA-0208146}-1.2.2
on page=9.
}

\DefLabeledTheorem{inverse matrix, Hadamard}{\Product}
{
Let \nTimes matrix $a$ have
a \RC inverse matrix.
Then entries of the \RC inverse matrix satisfy to the equality
\DrawEq{inverse matrix}{\Product}
\DrawEq{inverse matrix, Hadamard}{\Product}
}

\DefProof{inverse matrix, Hadamard}
{
The equality
\eqRef{inverse matrix}{\Product}
follows from the equality
\eqRef{inverse minor}{\Product}
when we assume
\ShowEq{I=i J=j}
The equality
\eqRef{inverse matrix, Hadamard}{\Product}
follows from equalities
\EqRef{Hadamard inverse of matrix 2 entry},
\eqRef{inverse matrix}{\Product}.
}

\DefLabeledTheoremNote{inverse minor}{\Product}
{
Let \nTimes matrix $a$ have
\ProductType inverse matrix.\,\footnotemark
Then \nTimes[k] submatrix of \ProductType inverse matrix
satisfies to the equality
\DrawEq{inverse minor}{\Product}
\ShowEq{|I|=|J|=k}
}
{
This statement and its proof is based
on statement 1.2.1 on page
\citeBib{math.QA-0208146}\Hyph 8
for matrix over free division ring.
}

\DefProof{inverse minor}
{
We can represent the definition
\newline
\FrameEqRef{\Product-inverse matrix}1
\newline
of \ProductType inverse matrix
the following way
\ShowEq{a*a- submatrix, \Product}
The equality
\ShowEq{a[J]*aJ=0 \Product}
\DrawEq{row [J] * col J =0}{\Product}
and the equality
\ShowEq{aJ*aJ=E \Product}
\DrawEq{row J * col J =E}{\Product}
follow from the equality
\EqRef{a*a- submatrix, \Product}
We multiply the equality
\eqRef{row [J] * col J =0}{\Product}
by
\ShowEq{minor [I][J] inverse}
\DrawEq{row [J] * col J =0 1}{\Product}
The equality
\DrawEq{row J * col J =E 1}{\Product}
follows from equalities
\eqRef{row [J] * col J =0 1}{\Product},
\eqRef{row J * col J =E}{\Product}.
The equality
\eqRef{inverse minor}{\Product} follows from
\eqRef{row J * col J =E 1}{\Product}
when we both sides of the equality
\eqRef{row J * col J =E 1}{\Product}
multiply by
\ShowEq{minor IJ inverse}
}

\DefLabeledDefinition{singular matrix}{\Product}
{
If $\gin\times\gin$ matrix $a$ has \ProductType inverse matrix
we call such matrix
\AddIndex{\ProductType nonsingular matrix}{\Product nonsingular matrix}.
Otherwise, we call such matrix
\AddIndex{\ProductType singular matrix}{\Product singular matrix}.
}

\DefLabeledDefinition{inverse matrix}{\Product}
{
Let $a$ be \nTimes matrix.
If there exists \nTimes matrix $b$
such that
\ShowEq{\Product-inverse matrix 1}
then the matrix
\ShowEq{\Product-inverse element}
is called
\AddIndex{\ProductType inverse element}{\Product-inverse element}
of the matrix $a$
\DrawEq{\Product-inverse matrix}1
\ePrints{2020.06.01}
\ifx\Semafor\ValueOn
Matrix $a$ is called
\RC regular,
if there exists \RC inverse matrix.
\fi
}

\AddEq{theorem: two products equal}
{
\begin{ShadedTheorem}
\labelTheorem{two \Product-products equal}
Let matrix $a$ have \ProductType inverse matrix. Then for any matrices
$b$ and $c$ the equality
\DrawEq{two products equal, 2}{\Product}
follows from the equality
\ShowEq{two \Product-products equal, 1}
\end{ShadedTheorem}
}

\DefProof{two products equal}
{
The equality \eqRef{two products equal, 2}{\Product} follows from the equality
\EqRef{two \Product-products equal, 1} if we multiply both parts of the equality
\EqRef{two \Product-products equal, 1} over
\ShowEq{a-\Product}
}

\DefLabeledTheorem{inverse matrix of product}{\Product}
{
Let \nTimes matrices $a$, $b$ have \ProductType inverse matrix. Then
the matrix
\ShowEq{a.\Product.b}
also has \ProductType inverse matrix and
\ShowEq{\Product-inverse matrix of product}
}

\DefProof{inverse matrix of product}
{
The equality
\ShowEq{\Product.ab b-a-}
follows from the equality
\newline
\FrameEqRef{\Product-inverse matrix}1
\newline
and from the definition
\RefDefinition{biring matrix}.
According to the definition
\refDefinition{inverse matrix}{\Product}
and the equality
\EqRef{\Product.ab b-a-},
the matrix
\ShowEq{a.\Product.b}
has \ProductType inverse matrix
\ShowEq{\Product.ab ba-}
The equality
\EqRef{\Product-inverse matrix of product}
follows from equalities
\EqRef{\Product.ab b-a-},
\EqRef{\Product.ab ba-}.
}

%% file: Biring.Stmt.Eq.tex

\AddEq{def minor rc}
{
\def\InvMatrix{(a^{\RCInverse})}
\def\MinorA{\aUD a{[J]}{[I]}}
\def\MinorB{\aUD a{[J]}I}
\def\MinorC{\aUD aJ{[I]}}
\def\MinorD{\aUD aJI}
\def\MinorE{\aUD{\InvMatrix}{[I]}J}
\def\MinorF{\aUD{\InvMatrix}IJ}

\def\mMinorA{\aUD{(ma)}{[J]}{[I]}}
\def\mMinorB{\aUD{(ma)}{[J]}I}
\def\mMinorC{\aUD{(ma)}J{[I]}}
\def\mMinorD{\aUD{(ma)}JI}
\def\mMinorF{(((ma)^{\RCInverse})\aUD{}IJ)^{\RCInverse}}

\def\minorA{\aUD a{[j]}{[i]}}
\def\minorB{\aUD a{[j]}i}
\def\minorC{\aUD aj{[i]}}
\def\minorD{\aUD aji}
\def\minorF{\aUD{\InvMatrix}ij}

\def\QDef{$(\aUD{}ji)$}
\def\QEntry{\detVal jia}
\def\QEntryA{\aUD{(\mathcal{H}a^{\InverseVal})}ji}
\def\QEntryB{(\aUD{(a^{\InverseVal})}ij)^{-1}}
\def\RCDetA{\mathcal{H}\RCdet{[j]}{[i]}a}

\def\BAA{\aUD a11-\aUD a12\aBB\aUD a21}
\def\BAB{\aUD a12-\aUD a11\aBA\aUD a22}
\def\BBA{\aUD a21-\aUD a22\aAB\aUD a11}
\def\BBB{\aUD a22-\aUD a21\aAA\aUD a12}
}

\AddEq{def minor cr}
{
\def\InvMatrix{(a^{\CRInverse})}
\def\MinorA{\aUD a{[I]}{[J]}}
\def\MinorB{\aUD aI{[J]}}
\def\MinorC{\aUD a{[I]}J}
\def\MinorD{\aUD aIJ}
\def\MinorE{\aUD{\InvMatrix}J{[I]}}
\def\MinorF{\aUD{\InvMatrix}JI}

\def\mMinorA{\aUD{(ma)}{[I]}{[J]}}
\def\mMinorB{\aUD{(ma)}I{[J]}}
\def\mMinorC{\aUD{(ma)}{[I]}J}
\def\mMinorD{\aUD{(ma)}IJ}
\def\mMinorF{(((ma)^{\RCInverse})\aUD{}JI)^{\RCInverse}}

\def\minorA{\aUD a{[i]}{[j]}}
\def\minorB{\aUD ai{[j]}}
\def\minorC{\aUD a{[i]}j}
\def\minorD{\aUD aij}
\def\minorF{\aUD{\InvMatrix}ji}

\def\QDef{$(\aUD{}ij)$}
\def\QEntry{\detVal ija}
\def\QEntryA{\aUD{(\mathcal{H}a^{\InverseVal})}ij}
\def\QEntryB{(\aUD{(a^{\InverseVal})}ji)^{-1}}
\def\RCDetA{\mathcal{H}\RCdet{[i]}{[j]}a}

\def\BAA{\aUD a11-\aUD a21\aBB\aUD a12}
\def\BAB{\aUD a12-\aUD a22\aBA\aUD a11}
\def\BBA{\aUD a21-\aUD a11\aAB\aUD a22}
\def\BBB{\aUD a22-\aUD a12\aAA\aUD a21}
}
\def\RCDetDT{\RCdet ija^T}

\newcommand{\SA}[1]{\aUD aS{#1}}
\def\Bb{\aUD bpr}
\def\Bf{(\aUD bST)^{\InverseVal}}
\def\SATm{((\SA T)^{\InverseVal})}
\def\Bq{\detVal prb}
\def\tRr{\aUD Rtr}
\def\pRs{\aUD Rps}

\AddEq{def rank rc}
{
\def\Ba{\aUD bpT}
\def\Bc{\aUD {b^{\RCInverse}_{}{}}rp}
\def\Bd{\aUD {b^{\RCInverse}_{}{}}Tp}
\def\Be{\aUD bSr}
}

\AddEq{def rank cr}
{
\def\Ba{\aUD bSr}
\def\Bc{\aUD {b^{\CRInverse}_{}{}}rp}
\def\Bd{\aUD {b^{\CRInverse}_{}{}}rS}
\def\Be{\aUD bpT}
}

\AddEq{r in N-T}
{
$\gi r\in\gi N\setminus\gi T$
}

\AddEq{p in M-S}
{
$\gi p\in\gi M\setminus\gi S$
}

\AddEq{quasideterminant, matrix 2x2}
{
\DetVal a=
\begin{pmatrix}
\BAA&\BAB
\\
\BBA&\BBB
\end{pmatrix}
\ifx\texFuture\Defined
\RCDet a=
\aUD a11-\aUD a12(\aUD a22)^{-1}\aUD a21&\aUD a12-\aUD a11(\aUD a21)^{-1}\aUD a22\\
\aUD a21-\aUD a22(\aUD a12)^{-1}\aUD a11&\aUD a22-\aUD a21(\aUD a11)^{-1}\aUD a12
\fi
}

\AddEq{inverse matrix 2x2}
{
a^{\InverseVal}=
\begin{pmatrix}
(\BAA)^{-1}&(\BBA)^{-1}
\\
(\BAB)^{-1}&(\BBB)^{-1}
\end{pmatrix}
}

\AddEq[1]{rank a=gi n}
{
\RankVal a=\gi{#1}
}

\AddEq[3]{c=a rc b}
{
#3=
\ShowEq{\SideWS a rc b}{#1}{#2}
}

\AddEq[3]{a rc b=c}
{
\ShowEq{\SideWS a rc b}{#1}{#2}
=#3
}

\AddEq{gi k=min nm}
{
\[
\gik=\min(\gin,\gim)
\]
}

\AddEq{gi n<=m left}
{
$\gin\le\gim$.
}

\AddEq{gi n<=m right}
{
$\gim\le\gin$.
}

\AddEq[2]{linear equations for Inverse augmented cols}
{
\[
\begin{pmatrix}
\aUD a{#2_1}1&...&\aUD a{#2_1}{#1}&\aUD{\delta}11&...&\aUD{\delta}1k
\\...&...&...&...&...&...\\
\aUD a{#2_k}1&...&\aUD a{#2_k}{#1}&\aUD{\delta}k1&...&\aUD{\delta}kk
\end{pmatrix}
\]
}

\AddEq[1]{set of inverse matrices}
{
\symb[#1-0]{\SideNS(a^{\InverseVal})}{set of inverse matrices}1
}

\AddEq[4]{a in left}
{
$#1\in #3(#2^{\InverseVal})$#4
}

\AddEq{ref for a in left}
{
\refDefinition{Inverse matrix nm}{\Product-\SideNS},
\refDefinition{Inverse matrix nm}{\Product-\OtherSideNS},
\refDefinition{Inverse Matrix set}{\Product-\SideNS},
\refDefinition{Inverse Matrix set}{\Product-\OtherSideNS}.
}

\AddEq[2]{linear equations for Inverse augmented rows}
{
\[
\begin{pmatrix}
\aUD a1{#2_1}&...&\aUD a1{#2_k}
\\...&...&...\\
\aUD a{#1}{#2_1}&...&\aUD a{#1}{#2_k}
\\
\aUD{\delta}11&...&\aUD{\delta}1k
\\...&...&...\\
\aUD{\delta}k1&...&\aUD{\delta}kk
\end{pmatrix}
\]
}

\AddEq[3]{linear equations rank left}
{
\begin{split}
\AUD b1{#3}\AUD a{#2_1}1+...+\AUD b{#1}{#3}\AUD a{#2_1}{#1}&=\AUD{\delta}1{#3}
\\ ...&... \\
\AUD b1{#3}\AUD a{#2_k}1+...+\AUD b{#1}{#3}\AUD a{#2_k}{#1}&=\AUD{\delta}k{#3}
\end{split}
}

\AddEq[3]{linear equations rank right}
{
\begin{split}
\AUD a{#2_1}1\AUD b1{#3}+...+\AUD a{#2_1}{#1}\AUD b{#1}{#3}&=\AUD{\delta}1{#3}
\\ ...&... \\
\AUD a{#2_k}1\AUD b1{#3}+...+\AUD a{#2_k}{#1}\AUD b{#1}{#3}&=\AUD{\delta}k{#3}
\end{split}
}

\AddEq{set of m vectors}
{
$(\aU{\Vector a}i,\gii\in\gi M,|\gi M|=\gim)$
}

\AddEq[2]{xj=bij}
{
$\ACol x{#2}=\AUD b{#2}{#1}$.
}

\DefText{Notation for Matrix Entry}
{
\ShowText{submatrix}

\ShowDefinition{notation for submatrix}

\ShowRemark{combine the notation of indices}

\ShowText{Kronecker symbol}

\ShowText{interval [kl]}
}

\AddEq[1]{I=[1k]}
{
\[
\gi{#1}=[\gi 1,\gik]
\]
}

\AddEq[2]{aI rc bJ=cIJ}
{
\ShowEq{c=a rc b}{\ARow b{#2}}{\ACol a{#1}}{\AUD c{#1}{#2}}
}

\AddEq{cIJ=Ek}
{
\[\aUD cIJ=\aD Ek\]
}

\AddEq{interval [kl]}
{
\[
[\gik,\gil]=\{\gii:\gik\le\gii\le\gil\}
\]
}

\AddEq[1]{v(i)=vi col}
{
\[
#1(\gii)=\aU{#1}i
\]
}

\AddEq[1]{v(i)=vi row}
{
\[
#1(\gii)=\aD{#1}i
\]
}

\AddEq{ref quasideterminant matrix}
{
\eqRef{quasideterminant, matrix 2x2 11}{\Product},
\eqRef{quasideterminant, matrix 2x2 21}{\Product},
\eqRef{quasideterminant, matrix 2x2 12}{\Product},
\eqRef{quasideterminant, matrix 2x2 22}{\Product}.
}

\AddEq{quasideterminant, matrix 2x2 11}
{
\detVal 11a=\BAA
}

\AddEq{quasideterminant, matrix 2x2 21}
{
\detVal 21a=\BAB
}

\AddEq{quasideterminant, matrix 2x2 12}
{
\detVal  12a=\BBA
}

\AddEq{quasideterminant, matrix 2x2 22}
{
\detVal 22a=\BBB
}

\def\aAA{(\aUD a11)^{-1}}
\def\aBA{(\aUD a21)^{-1}}
\def\aAB{(\aUD a12)^{-1}}
\def\aBB{(\aUD a22)^{-1}}
\DefText{Quasideterminant}
{
\section{\texorpdfstring{\ProductType}{\Product}-\SectionTitle}

\ShowEq{def minor \Product}
\ProveTheorem{inverse minor}

\ProveTheorem{inverse matrix, Hadamard}

\ShowDefinition{quasideterminant}

\ProveTheorem{quasideterminant and inverse}

\ProveTheorem{quasideterminant, expression}

\ProveTheorem{rdet=a--}

\ProveTheorem{quasideterminant, matrix 2x2}

\ProveTheorem{inverse product of matrix over scalar}

\ProveTheorem{inverse unit matrix order 2}
}

\AddEq{inverse matrix, unit}
{
a^{\InverseVal}=\UnitMatrix
}

\AddEq{unit matrix order 2}
{
a=\UnitMatrix
}

\AddEq{inverse product of matrix over scalar, 2}
{
(am)^{\InverseVal}=m^{-1}a^{\InverseVal}
}

\AddEq{inverse product of matrix over scalar, 1}
{
(ma)^{\InverseVal}=a^{\InverseVal}m^{-1}
}

\DefLabeledTheorem{inverse product of matrix over scalar}{\Product}
{
\DrawEq{inverse product of matrix over scalar, 1}{\Product}
\DrawEq{inverse product of matrix over scalar, 2}{\Product}
}

\AddEq{inverse product of matrix over scalar, 10}
{
\[
\left(ma\right)^{\InverseVal}=\left((ma)^{-1}\right)=\left(a^{-1}m^{-1}\right)
=\left(a^{-1}\right)m^{-1}=a^{\InverseVal}m^{-1}
\]
}

\AddEq{inverse product of matrix over scalar, 1, 2}
{
\begin{align*}
\mMinorF
&=\mMinorD
-\mMinorC\ProductVal
\left(\mMinorA\right)^{\InverseVal}\ProductVal
\mMinorB\\
&=m\MinorD
-m\ \MinorC\ProductVal
\left(\MinorA\right)^{\InverseVal}\ m^{-1}\ProductVal
m\ \MinorB\\
&=m\ \MinorD
-m\ \MinorC\ProductVal
\left(\MinorA\right)^{\RCInverse}\ProductVal
\MinorB
\end{align*}
}

\AddEq{inverse product of matrix over scalar, 1, 1}
{
\begin{split}
\mMinorF
=m\ \MinorF
\end{split}
}

\AddEq{rdet=a--}
{
\PMatrix[a]{\DetVal}nn
=
\begin{pmatrix}
(\aUD{(a^{\InverseVal})}11)^{-1}&...&(\aUD{(a^{\InverseVal})}n1)^{-1}\\
...&...&...\\
(\aUD{(a^{\InverseVal})}1n)^{-1}&...&(\aUD{(a^{\InverseVal})}nn)^{-1}
\end{pmatrix}
}

\AddEq{rdet-=a-}
{
\begin{pmatrix}
(\detVal 11a)^{-1}&...&(\detVal n1a)^{-1}\\
...&...&...\\
(\detVal 1na)^{-1}&...&(\detVal nna)^{-1}
\end{pmatrix}
=
\PMatrix{(a^{\InverseVal})}nn
}

\AddEq{quasideterminant, expression}
{
\QEntry=\minorD
-\minorC\ProductVal
\RCDetA\ProductVal
\minorB
}

\DefProofRef{singular matrix and quasideterminant}{\Product}
{
\ShowText{proof of singular matrix and quasideterminant}
}

\DefProof{singular matrix and quasideterminant *}
{
\ShowText{proof of singular matrix and quasideterminant}
}

\AddEq{[I]= [J]= cr}
{
\(\gi{[I]}=\gi S\), \(\gi{[J]}=\gi T\)
}

\AddEq{I= J= cr}
{
$\gi I=\{\gi p\}$, $\gi J=\{\gi r\}$
}

\AddEq{[I]= [J]= rc}
{
\(\gi{[I]}=\gi T\), \(\gi{[J]}=\gi S\)
}

\AddEq{I= J= rc}
{
$\gi I=\{\gi r\}$, $\gi J=\{\gi p\}$
}

\AddEq{singular matrix and quasideterminant, 1}
{
\Be\ \Bc
+
\aUD bST\ProductVal\Bd
=0
}

\AddEq{singular matrix and quasideterminant, 2}
{
\Bb\ \Bc
+
\Ba\ProductVal\Bd
=1
}

\AddEq{singular matrix and quasideterminant, 3}
{
\Bf\ProductVal\Be\ \Bc
+
\Bd
=0
}

\AddEq{singular matrix and quasideterminant, 4}
{
\Bb \ \Bc
-\Ba\ProductVal
\Bf\ProductVal
\Be
\ \Bc
=1
}

\AddEq{singular matrix and quasideterminant, 5}
{
(\Bb-\Ba\ProductVal
\Bf\ProductVal
\Be)
\ \Bc=1
}

\AddEq{singular matrix and quasideterminant, 6}
{
(\Bq)\Bc=1
}

\AddEq{quasideterminant and inverse}
{
a^{\InverseVal}=\mathcal{H}\DetVal a
}

\AddEq{B rc-1}
{
\(b^{\InverseVal}\).
}

\AddEq{singular matrix and quasideterminant}
{
\detVal pr\SATpr=0
}

\AddEq{singular matrix and quasideterminant, 0}
{
b=\SATpr
}

\AddEq{j i quasideterminant H}
{
\QEntry=\QEntryA=\QEntryB
}

\AddEq{quasideterminant definition= H}
{
\[
\DetVal\,a=\mathcal{H}a^{\InverseVal}
\]
}

\AddEq{j i quasideterminant definition}
{
\symb{\QEntry}{j i \Product-quasideterminant definition}{}
}

\AddEq{j i quasideterminant =}
{
\ShowSymbol{j i \Product-quasideterminant definition}{}
=\minorD
-\minorC\ProductVal
\left(\minorA\right)^{\InverseVal}\ProductVal
\minorB
=(\aUD{(a^{\InverseVal})}ij)^{-1}
}

\AddEq{quasideterminant definition}
{
\symb{\DetVal\,a}{\Product-quasideterminant definition}{}
}

\AddEq{quasideterminant definition=}
{
\begin{split}
\ShowSymbol{\Product-quasideterminant definition}{}&=
\PMatrix[a]{\DetVal}nn
\\ &=
\begin{pmatrix}
(\aUD{(a^{\InverseVal})}11)^{-1}&...&(\aUD{(a^{\InverseVal})}n1)^{-1}\\
...&...&...\\
(\aUD{(a^{\InverseVal})}1n)^{-1}&...&(\aUD{(a^{\InverseVal})}nn)^{-1}
\end{pmatrix}
\end{split}
}

\AddEq{inverse matrix}
{
\minorF
=\left(
\minorD
-\minorC\ProductVal
\left(\minorA\right)^{\InverseVal}\ProductVal
\minorB
\right)^{-1}
}

\AddEq{inverse matrix, Hadamard}
{
\aUD{(\mathcal{H}a^{\InverseVal})}ji
=\minorD
-\minorC\ProductVal
\left(\minorA\right)^{\InverseVal}\ProductVal
\minorB
}

\AddEq{inverse minor}
{
\left(\MinorF\right)^{\InverseVal}
=\MinorD
-\MinorC\ProductVal
\left(\MinorA\right)^{\InverseVal}\ProductVal
\MinorB
}

\AddEq[1]{minor matrix a- cols J}
{
$\aD{(a^{\InverseVal})}J$#1
}

\AddEq[3]{submatrix IJ}
{
\begin{pmatrix}
\aUD {#1}{#2}{#3}&\aUD {#1}{#2}{[#3]}\\
\aUD {#1}{[#2]}{#3}&\aUD {#1}{[#2]}{[#3]}
\end{pmatrix}
}

\AddEq{submatrix E k n}
{
=
\begin{pmatrix}
\aD Ek&0\\0&\aD E{n-k}
\end{pmatrix}
}

\AddEquation{a*a- submatrix, rc}
{
\ShowEq{submatrix IJ}aJI
\RCstar
\ShowEq{submatrix IJ}{(a^{\RCInverse})}IJ
\ShowEq{submatrix E k n}
}

\AddEquation{a*a- submatrix, cr}
{
\ShowEq{submatrix IJ}aIJ
\CRstar
\ShowEq{submatrix IJ}{(a^{\CRInverse})}JI
\ShowEq{submatrix E k n}
}

\AddEq{a[J]*aJ=0 rc}
{
($\aU a{[J]}\RCstar\aD{(a^{\RCInverse})}J=0$)
}

\AddEq{a[J]*aJ=0 cr}
{
($\aD a{[J]}\CRstar\aU{(a^{\CRInverse})}J=0$)
}

\AddEq{aJ*aJ=E rc}
{
($\aU aJ\RCstar\aD{(a^{\RCInverse})}J=\aD Ek$)
}

\AddEq{aJ*aJ=E cr}
{
($\aD aJ\CRstar\aU{(a^{\CRInverse})}J=\aD Ek$)
}

\AddEq{row [J] * col J =0}
{
\MinorB\ProductVal\MinorF
+
\MinorA\ProductVal\MinorE
=0
}

\AddEq{row J * col J =E}
{
\MinorD\ProductVal\MinorF
+
\MinorC\ProductVal\MinorE
=\aD Ek
}

\AddEq{row [J] * col J =0 1}
{
\left(\MinorA\right)^{\InverseVal}\ProductVal
\MinorB\ProductVal\MinorF
+\MinorE=0
}

\AddEq{row J * col J =E 1}
{
\MinorD\ProductVal\MinorF
-\MinorC\ProductVal
\left(\MinorA\right)^{\InverseVal}\ProductVal
\MinorB\ProductVal
\MinorF
=\aD Ek
}

\AddEq{minor [I][J] inverse}
{
$\left(\MinorA\right)^{\InverseVal}$
}

\AddEq{minor IJ inverse}
{
$\left(\MinorF\right)^{\InverseVal}$.
}

\AddEq{|I|=|J|=k}
{
$|\gi I|=|\gi J|=\gik<\gin$.
}

\AddEq{delta n=Matrix}
{
\[
\UnitNMatrix=(\aUD {\delta}ij)\ \ \ \gi 1\le\gii,\gij\le\gin
\]
}

\AddEquation{b'i=. fb}
{
\aU{b'}i=
b'^\cdot_{}{}^-_{\gik}=f^\cdot{}{}^-_{\gik}{}_\cdot^{}{}^{\gil}_- b^\cdot_{}{}^-_{\gil}
=\aUD fij\aU bj
}

\AddEq{column of matrix}
{
\[
\symb{\aD ai}{minor matrix}3
=
\begin{pmatrix}
\aUD a1i\\...\\ \aUD ami
\end{pmatrix}
\]
}

\AddEq{A from columns T}
{
\symb{\aD aT}{minor matrix}1
}

\AddEq{A without column a}
{
\symb{a_{[\gii]}}{minor matrix}1
}

\AddEq{A without columns T}
{
\symb{a_{[\gi T]}}{minor matrix}1
}

\AddEq{row of matrix}
{
\[
\symb{\aU aj}{minor matrix}3
=
\begin{pmatrix}
\aUD aj1&...&\aUD ajn
\end{pmatrix}
\]
}

\AddEq{A from rows S}
{
\symb{\aU aS}{minor matrix}1
}

\AddEq{A without row b}
{
\symb{a^{[\gij]}}{minor matrix}1
}

\AddEq{A without rows S}
{
\symb{a^{[\gi S]}}{minor matrix}1
}

\AddEq{A from b a}
{
\symb{\aUD aji}{minor matrix}1
}

\AddEq{I,|I|=n}
{
$\gi I$, $|\gi I|=\gin$,
}

\AddEq{Kronecker symbol}
{
\symb{\aUD{\delta}ij}{Kronecker symbol}{}
}

\AddEq{Kronecker symbol=}
{
\[
\ShowSymbol{Kronecker symbol}{}=\left\{\begin{array}{ccc}
1\ \ \ \gii=\gij\\
0\ \ \ \gii\ne \gij
\end{array}\right.\ \ \ \gii,\gij\in\gi I
\]
}

\AddEq{a gi 1...n}
{
$\aD a1$, ..., $\aD an$.
}

\AddEq{gi 1...n}
{
$\gi 1$, ..., $\gin$
}

\AddEq{show power matrix}
{

\FrameEqRef{rc power, 0}1
\FrameEqRef{rc-power, n}1

\FrameEqRef{cr power, 0}1
\FrameEqRef{cr-power, n}1

\noindent
}

\AddEq{show inverse matrix}
{

\FrameEqRef{rc-inverse matrix}1
\FrameEqRef{cr-inverse matrix}1
}

\AddEq{power matrix}
{
\TwoColText
{
\ShowEq{def rc}
\ShowDefinition{power matrix}
}
{
\ShowEq{def cr}
\ShowDefinition{power matrix}
}

\TwoColText
{
\ShowTheorem{rc-power, 1}
\begin{sloppypar}
\ShowProof{rc-power, 1}
\end{sloppypar}
}
{
\ShowTheorem{cr-power, 1}
\begin{sloppypar}
\ShowProof{cr-power, 1}
\end{sloppypar}
}

\TwoColText
{
\ShowEq{def rc}
\ProveTheorem{transpose power}
}
{
\ShowEq{def cr}
\ProveTheorem{transpose power}
}

\TwoColText
{
\ShowEq{def rc}
\ShowDefinition{inverse matrix}

\ShowDefinition{singular matrix}
}
{
\ShowEq{def cr}
\ShowDefinition{inverse matrix}

\ShowDefinition{singular matrix}
}

\TwoColText
{
\ShowEq{def rc}
\ShowTheorem{transpose inverse}
\begin{sloppypar}
\ShowProof{transpose inverse}
\end{sloppypar}
}
{
\ShowEq{def cr}
\ShowTheorem{transpose inverse}
\begin{sloppypar}
\ShowProof{transpose inverse}
\end{sloppypar}
}
}

\AddEq{transpose inverse, 1}
{
\aD En=a^T\ProductAVal( a^{\InverseVal})^T
}

\AddEq{transpose inverse, 0}
{
\[
(a\ProductVal a^{\InverseVal})^T=\aD En^T=\aD En
\]
}

\AddEq{transpose inverse}
{
(a^T)^{\InverseAVal}=(a^{\InverseVal})^T
}

\DefLabeledTheorem{transpose power}{\Product}
{
\DrawEq[n]{transpose power}{\Product}
}

\AddEq[1]{transpose power}
{
(a^T)^{\PowerVal {#1}}=(a^{\PowerAVal {#1}})^T
}

\AddEq{transpose power, 1}
{
(a^T)^{\PowerVal k}=(a^T)^{\PowerVal {k-1}}\ProductVal a^T
}

\AddEq{transpose power, 2}
{
(a^T)^{\PowerVal k}=(a^{\PowerAVal {k-1}})^T\ProductVal a^T
}

\AddEq{transpose power, 3}
{
(a^T)^{\PowerVal k}=(a^{\PowerAVal {k-1}}\ProductAVal a)^T
}

\AddEq{Quasideterminant}
{
\ShowRemark{we do not have definition of determinant for division algebra}

\ShowEq{def minor rc}
\ShowDefinition{RC-quasideterminant}

\ShowTheorem{quasideterminant rc cr, matrix 2x2}
}

\AddEquation{cr-inverse matrix 1}
{
a\CRstar b=\UnitNMatrix
}

\AddEq{cr-inverse element}
{
\symb{a^{\CRInverse}}{cr-inverse matrix}{}
$b=\ShowSymbol{cr-inverse matrix}{}$
}

\AddEquation{rc-inverse matrix 1}
{
a\RCstar b=\UnitNMatrix
}

\AddEq{rc-inverse element}
{
\symb{a^{\RCInverse}}{rc-inverse matrix}{}
$b=\ShowSymbol{rc-inverse matrix}{}$
}

\AddEquation{rc-inverse matrix of product}
{
(a\RCstar b)^{\RCInverse}=b^{\RCInverse}\RCstar a^{\RCInverse}
}

\AddEquation{cr-inverse matrix of product}
{
(a\CRstar b)^{\CRInverse}=b^{\CRInverse}\CRstar a^{\CRInverse}
}

\AddEquation{rc.ab b-a-}
{
\begin{split}
&\,(a\RCstar b)\RCstar\BlueText{(b^{\RCInverse}\RCstar a^{\RCInverse})}
\\=&\,
a\RCstar (b\RCstar b^{\RCInverse})\RCstar a^{\RCInverse}
\\=&\,
a\RCstar \aD En\RCstar a^{\RCInverse}
=
a\RCstar a^{\RCInverse}
\\=&\,
\aD En
\end{split}
}

\AddEquation{cr.ab b-a-}
{
\begin{split}
&\,(a\CRstar b)\CRstar\BlueText{(b^{\CRInverse}\CRstar a^{\CRInverse})}
\\=&\,
a\CRstar (b\CRstar b^{\CRInverse})\CRstar a^{\CRInverse}
\\=&\,
a\CRstar \aD En\CRstar a^{\CRInverse}
=
a\CRstar a^{\CRInverse}
\\=&\,
\aD En
\end{split}
}

\AddEquation{rc.ab ba-}
{
(a\RCstar b)\RCstar\BlueText{(a\RCstar b)^{\RCInverse}}=\aD En
}

\AddEquation{cr.ab ba-}
{
(a\CRstar b)\CRstar\BlueText{(a\CRstar b)^{\CRInverse}}=\aD En
}

\AddEq{a.rc.b}
{
$a\RCstar b$
}

\AddEq{a.cr.b}
{
$a\CRstar b$
}

\AddEquation{two rc-products equal, 1}
{
b\RCstar a=c\RCstar a
}

\AddEquation{two cr-products equal, 1}
{
b\CRstar a=c\CRstar a
}

\AddEq{two products equal, 2}
{
b=c
}

\AddEq{a-rc}
{
$a^{\RCInverse}$.
}

\AddEq{a-cr}
{
$a^{\CRInverse}$.
}

\AddIndex{}{distributive law}%
\AddEq{product of matrices is distributive over addition}
{
\TwoColText
{
\ShowEq{def rc}
\ShowTheorem{product of matrices is distributive over addition}
}
{
\ShowEq{def cr}
\ShowTheorem{product of matrices is distributive over addition}
}

\TwoColText
{
\ShowEq{def rc}
\begin{sloppypar}
\ShowProof{product of matrices is distributive over addition}
\end{sloppypar}
}
{
\ShowEq{def cr}
\begin{sloppypar}
\ShowProof{product of matrices is distributive over addition}
\end{sloppypar}
}
}

\AddEq{a.rc.(b1+b2)=}
{
a\RCstar(b_1+b_2)=a\RCstar b_1+a\RCstar b_2
}

\AddEq{a.cr.(b1+b2)=}
{
a\CRstar(b_1+b_2)=a\CRstar b_1+a\CRstar b_2
}

\AddEq{(b1+b2).rc.a=}
{
(b_1+b_2)\RCstar a=b_1\RCstar a+b_2\RCstar a
}

\AddEq{(b1+b2).cr.a=}
{
(b_1+b_2)\CRstar a=b_1\CRstar a+b_2\CRstar a
}

\AddEquation{a.rc.(b1+b2)= 1}
{
\begin{aligned}
&\,\aUD{(a\RCstar(b_1+b_2))}ij
\\=&\,\aUD aik\aUD{(b_1+b_2)}kj
\\=&\,\aUD aik(\aUD{b_1}kj+\aUD{b_2}kj)
\\=&\,\aUD aik\aUD{b_1}kj+\aUD aik\aUD{b_2}kj
\\=&\,\aUD{(a\RCstar b_1)}ij+\aUD{(a\RCstar b_2)}ij
\end{aligned}
}

\AddEquation{a.cr.(b1+b2)= 1}
{
\begin{aligned}
&\,\aUD{(a\CRstar(b_1+b_2))}ij
\\=&\,\aUD akj\aUD{(b_1+b_2)}ik
\\=&\,\aUD akj(\aUD{b_1}ik+\aUD{b_2}ik)
\\=&\,\aUD akj\aUD{b_1}ik+\aUD akj\aUD{b_2}ik
\\=&\,\aUD{(a\CRstar b_1)}ij+\aUD{(a\CRstar b_2)}ji
\end{aligned}
}
	
\AddEq{a.Product.(b1+b2)= 2}
{
\EqRef{a(b+c)=.},
\EqRef{\Product-product of matrices},
\EqRef{(a+b)=}.
}

\AddEquation{(b1+b2).rc.a= 1}
{
\begin{aligned}
&\,\aUD{((b_1+b_2)\RCstar a)}ij
\\=&\,\aUD{(b_1+b_2)}ik\aUD akj
\\=&\,(\aUD{b_1}ik+\aUD{b_2}ik)\aUD akj
\\=&\,\aUD{b_1}ik\aUD akj+\aUD{b_2}ik\aUD akj
\\=&\,\aUD{(b_1\RCstar a)}ij+\aUD{(b_2\RCstar a)}ij
\end{aligned}
}

\AddEquation{(b1+b2).cr.a= 1}
{
\begin{aligned}
&\,\aUD{((b_1+b_2)\CRstar a)}ij
\\=&\,\aUD{(b_1+b_2)}kj\aUD aik
\\=&\,(\aUD{b_1}kj+\aUD{b_2}kj)\aUD aik
\\=&\,\aUD{b_1}kj\aUD aik+\aUD{b_2}kj\aUD aik
\\=&\,\aUD{(b_1\CRstar a)}ij+\aUD{(b_2\CRstar a)}ij
\end{aligned}
}
`	
\AddEq{(b1+b2).Product.a= 2}
{
\EqRef{(a+b)c=.},
\EqRef{\Product-product of matrices},
\EqRef{(a+b)=}.
}

%% file: Stmt.Omega.English.tex
\input{Stmt.Omega.Eq}

\DefDefinition{polyadditive map}
{
A map
\ShowEq{f:A->B}g{A^n}A
is called
\AddIndex{polyadditive map}{polyadditive map}
if for any
\ShowEq{fi(a+b)=}
}

\DefTheorem{shifts on group commuting}
{
Left and right shifts on multiplicative $\Omega$\Hyph group $A_1$ are commuting.
}

\DefProof{shifts on group commuting}
{
The theorem follows from the associativity of product on multiplicative $\Omega$\Hyph group $A_1$
\ShowEq{L(a)oR(b)=}
}

\DefProof{two representations of group}
{
We use group coordinates for $A_2$\Hyph numbers $a_2$.
Then according to theorem \RefTheorem{single transitive representation of group}
we can write the left shift $L(a_1)$ instead of the transformation $f(a_1)$.

Let
\ShowEq{a b in A}2.
Then we can find
one and only one
\ShowEq{a in A}1{}
such that
\ShowEq{b=aRa}
We assume
\ShowEq{ha=Ra}
For some
\ShowEq{b in A}1,
we have
\ShowEq{c=Lba}
According to the theorem \RefTheorem{shifts on group commuting}, the diagram
\ShowEq{Diagram: two representations of group}
is commutative.

Changing $b_1$ we get that $c_2$ is an arbitrary $A_2$\Hyph number.

We see from the diagram that if $a_2=b_2$ then $c_2=d_2$ and therefore $h(e)=\delta$.
On other hand if $a_2\neq b_2$ then $c_2\neq d_2$ because
the left\Hyph side $A_1$\Hyph representation $f$ is single transitive.
Therefore the right\Hyph side $A_1$\Hyph representation $h$ is effective.

In the same way we can show that for given $c_2$ we can find $a_1$
such that
\ShowEq{d2=c2 h}
Therefore the right\Hyph side $A_1$\Hyph representation $h$ is single transitive.

In general the product of transformations of
the left\Hyph side $A_1$\Hyph representation $f$ is not commutative and therefore
the right\Hyph side $A_1$\Hyph representation $h$ is different from the left\Hyph side $A_1$\Hyph representation $f$.
In the same way we can create a left\Hyph side $A_1$\Hyph representation $f$
using the right\Hyph side $A_1$\Hyph representation $h$.
}

\DefDefinition{semiring of sets}
{
A nonempty system of sets $\mathcal S$ is called
\AddIndex{semiring of sets}{semiring of sets},\,\footnote{
See also the definition
\citeBib{Kolmogorov Fomin}\Hyph 2,  page 32.}
if
\ifx\setCACAA\undefined
\StartLabelItem
\else
\StartLabelItem[definition]
\fi
\begin{enumerate}
\item
\ShowEq{empty in S}
\item
If
\ShowEq{A, B in S}
then
\ShowEq{AB in S}
\item
If
\ShowEq{A, A1 in S},
then the set \(A\) can be represented as
\DrawEq{A=+Ai}{def}
where
\ShowEq{Ai*Aj=0}
\labelItem{A=+Ai}
\end{enumerate}
The representation
\eqRef{A=+Ai}{def}
of the set \(A\) is called
\AddIndex{finite expansion of set}{finite expansion of set}
\(A\).
}

\DefDefinition{norm on Omega group}
{
\AddIndex{Norm on $\Omega$\Hyph group}
{norm on Omega group} $A$\,\footnote{
I made definition according\,to definition
from\,\citeBib{Bourbaki: General Topology: Chapter 5 - 10},
IX,\,\S 3.2 and definition \citeBib{Arnautov Glavatsky Mikhalev}-1.1.12,
p. 23.} is a map
\ShowEq{d->|d|}
which satisfies the following axioms
\ifx\setCACAA\undefined
\StartLabelItem
\else
\StartLabelItem[definition]
\fi
\begin{enumerate}
\ShowEq{|a|>=0}
\ShowEq{|a|=0}
if, and only if, $a=0$
\ShowEq{|a+b|<=|a|+|b|}
\ShowEq{|-a|=|a|}
\end{enumerate}

The $\Omega$\Hyph group $A$,
endowed with the structure defined by a given norm on $A$, is called
\AddIndex{normed $\Omega$\Hyph group}{normed Omega group}.
}

\DefTheorem{h=f1...fn omega, converges uniformly}
{
Let $A$ be complete $\Omega$\Hyph group.
Let
\ShowEq{omega in Omega}{}{}
be $n$\Hyph ari operation.
Let sequence of maps
\ShowEq{fim M(X,A)}
into complete $\Omega$\Hyph group $A$
converge uniformly
to the map $f_i$.
Let range of the map
$f_i$
be compact set.
Then sequence of maps
\ShowEq{hm=f1m...fnm omega}
into complete $\Omega$\Hyph group $A$
converges uniformly
to the map
\DrawEq{h=f1...fn omega}{converges uniformly}
}

\DefTheorem{|int h|<int|omega||f1n|}
{
Let
\ShowEq{omega in Omega}{}{}
be $n$\Hyph ary operation in
Abelian \(\Omega\)\Hyph group \(A\).
Let
\ShowEq{fi:X->A}
be $\mu$\Hyph measurable map
with compact range.
Since map \(f_i\), \(i=1\), ..., \(n\),
is integrable map, then map
\DrawEq{h=f1...fn omega}{}
is integrable map and
\DrawEq{|int h|<int|omega||f1n|}{measurable}
}

\DefTheorem{h=fX g1 g2 representation}
{
Let $\mu$ be a \(\sigma\)\Hyph additive measure defined on the set $X$.
Let
\ShowEq{f:A->*B}g{A_1}{A_2}
be representation
of $\Omega_1$\Hyph group $A_1$ with norm $\|x\|_1$
in $\Omega_2$\Hyph group $A_2$ with norm $\|x\|_2$.
Let
\ShowEq{gi:X->A 12}
be integrable map
with compact range.
Then map
\ShowEq{h=fX g1 g2}
is integrable map and
\DrawEq{|int h|<int|g1||g2|}{measurable}
}

\DefTheorem{polymorphism of representations}
{
Let the map
\ShowEq{map r1n,R}
be polymorphism of
representations $f_1$, ..., $f_n$ into representation $f$.
For any
\ShowEq{k,1n}
the map
\ShowEq{map r,R}{r_{1k}}{r_2}{}
satisfies to the equality
\ShowEq{polymorphism of representation, ak}
Let $\omega_1\in\Omega_1(p)$.
For any
\ShowEq{k,1n}
the map $r_{1k}$ satisfies to the equality
\ShowEq{polymorphism of representation, omega1}
Let $\omega_2\in\Omega_2(p)$.
For any
\ShowEq{k,1n}
the map $r_2$ satisfies to the equality
\ShowEq{polymorphism of representation, omega2}
}

\DefDefinition{polymorphism of representations}
{
Let
\ShowEq{A1n A}
be $\Omega_1$\Hyph algebras.
Let
\ShowEq{B1n}
be $\Omega_2$\Hyph algebras.
Let, for any
\ShowEq{k,1n}
\ShowEq{f:A->*B}{f_k}{A_k}{B_k}
be representation of $\Omega_1$\Hyph algebra $A_k$
in $\Omega_2$\Hyph algebra $B_k$.
Let
\ShowEq{f:A->*B}fAB
be representation of $\Omega_1$\Hyph algebra $A$
in $\Omega_2$\Hyph algebra $B$.
The map
\ShowEq{polymorphism of representation}
is called
\AddIndex{polymorphism of representations}
{polymorphism of representations}
$f_1$, ..., $f_n$ into representation $f$,
if, for any
\ShowEq{k,1n}
provided that all variables except variables
\ShowEq{ak bk in}
have given value, the map
\ShowEq{map r,R}{r_{1k}}{r_2}{}
is a morphism of representation $f_k$ into representation $f$.

If $f_1=...=f_n$, then we say that the map
\ShowEq{map r1n,R}
is polymorphism of representation $f_1$ into representation $f$.

If $f_1=...=f_n=f$, then we say that the map
\ShowEq{map r1n,R}
is polymorphism of representation $f$.
}

\DefDefinition{Omega group}
{
Let sum which is not necessarily commutative
be defined
in $\Omega_1$\Hyph algebra $A$.
We use the symbol $+$ to denote sum.
Let
\ShowEq{Omega=...-+}
If $\Omega_1$\Hyph algebra $A$ is group
relative to sum
and any operation $\omega\in\Omega$
is polyadditive map,
then $\Omega_1$\Hyph algebra $A$ is called
\AddIndex{$\Omega$\Hyph group}{Omega group}.
If $\Omega$\Hyph group $A$ is associative group
relative to sum,
then $\Omega_1$\Hyph algebra\,$A$\,is\,called
\AddIndex{associative $\Omega$\Hyph group}{associative Omega group}.
If $\Omega$\Hyph group $A$ is Abelian group
relative\,to sum,
then $\Omega_1$\Hyph algebra\,$A$ is\,called
\AddIndex{Abelian $\Omega$\Hyph group}{Abelian Omega group}.%
}

\DefTheorem{operation is distributive over addition}
{
Let
\ShowEq{omega n ari}{}{}n{}
be polyadditive map.
The operation $\omega$
\AddIndex{is distributive}{distributive law}
over addition
\ShowEq{omega is distributive over addition}
}

\DefDefinition{multiplicative Omega group}
{
Let product
\DrawEq{c=ab}{}
be operation of $\Omega_1$\Hyph algebra $A$.
Let
\ShowEq{Omega=...-o}
If $\Omega_1$\Hyph algebra $A$ is group
with respect to product
and, for any operation
\ShowEq{omega n ari}{}{}n,
the product is distributive over the operation $\omega$
\DrawEq{a omega=omega a left}{}
\DrawEq{a omega=omega a right}{}
then $\Omega_1$\Hyph algebra $A$ is called
\AddIndex{multiplicative $\Omega$\Hyph group}
{multiplicative Omega group}.
}

\DefDefinition{Abelian multiplicative Omega group}
{
If
\ShowEq{Abelian multiplicative Omega group}
then multiplicative $\Omega$\Hyph group is called
\AddIndex{Abelian}{Abelian multiplicative Omega group}.
}

\DefDefinition{Omega ring}
{
Let sum
\ShowEq{c1=a+b}
\ePrints{8428-0408,2207.06506}%
\ifx\Semafor\ValueOff%
which is not necessarily commutative
\fi%
and product
\DrawEq{c=ab}{}
be operations of $\Omega_1$\Hyph algebra $A$.
Let
\ShowEq{Omega=...-o+}
If $\Omega_1$\Hyph algebra $A$ is
\ePrints{8428-0408,2207.06506}
\ifx\Semafor\ValueOn
Abelian
\fi
$\Omega\cup\{*\}$\Hyph group
and
multiplicative $\Omega\cup\{+\}$\Hyph group,
then $\Omega_1$\Hyph algebra $A$ is called
\AddIndex{$\Omega$\Hyph ring}{Omega ring}.
}

\DefTheorem{product in Omega ring is distributive over addition}
{
The product in $\Omega$\Hyph ring
\AddIndex{is distributive}{distributive law}
over addition
\ShowEq{a*(b1+b2)=}
}

\DefProof{product in Omega ring is distributive over addition}
{
The theorem follows from the definitions
\RefDefinition{Omega group},
\RefDefinition{multiplicative Omega group},
\RefDefinition{Omega ring}.
}

\DefDefinition{matrix over Omega ring}
{
Let $A$ be $\Omega$\Hyph ring. The
\AddIndex{matrix}{matrix}
over $\Omega$\Hyph ring $A$ is a table of $A$\Hyph numbers $a^i_j$
where the index $i$ is the number of row
and the index $j$ is the number of column.
}

%% file: Stmt.Omega.Eq.tex

\AddEq{f(a+b)=}
{
\[
f(a+b)=f(a)+f(b)
\]
}

\AddEq{b=aRa}
{
\[b_2=a_2* a_1=a_2\circ R(a_1)\]
}

\AddEq{ha=Ra}
{
\[h(a)=R(a)\]
}

\AddEq[2]{b in A}
{
$b_{#1}\in A_{#1}$#2 
}

\AddEq{c=Lba}
{
\[c_2=f(b_1)\Act a_2=L(b_1)\circ a_2\ \ \ \ d_2=f(b_1)\Act b_2=L(b_1)\circ b_2\]
}

\AddEquation{Diagram: two representations of group}
{
\xymatrix{
a_2\ar[rrr]^{h(a_1)=R(a_1)}\ar[d]^{f(b_1)=L(b_1)}
&&& b_2\ar[d]^{f(b_1)=L(b_1)}\\
c_2\ar[rrr]_{h(a_1)=R(a_1)}&&&d_2
}
}

\AddEq{d2=c2 h}
{
$d_2=c_2\Act h(a_1)$.
}

\AddEq[2]{a b in A}
{
$a_{#1}$, $b_{#1}\in A_{#1}$#2
}

\AddEq{L(a)oR(b)=}
{
\[(L(a)\circ c)\circ R(b)=(a*c)*b = a*(c*b)=L(a)\circ (c\circ R(b))\]
}

\AddEq{fi(a+b)=}
{
$i$, $i=1$, ..., $n$,
\[
f(a_1,...,a_i+b_i,...,a_n)=f(a_1,...,a_i,...,a_n)+f(a_1,...,b_i,...,a_n)
\]
}

\AddEq{polymorphism of representation}
{
\[
\begin{pmatrix}
r_{1k}:A_k\rightarrow A
&
k=1,...,n
&
r_2:B_1\times...\times B_n\rightarrow B
\end{pmatrix}
\]
}

\AddEq{map r1n,R}
{
$((r_{1,k}, k=1,...,n)\ \ r_2)$
}

\AddEquation{polymorphism of representation, ak}
{
r_2(m_1,...,\BlueText{f_k(a_k)(m_k)},...,m_n)
=f(\RedText{r_{1k}(a_k)})(\BlueText{r_2(m_1,...,m_n)})
}

\AddEquation{reduced polymorphism of representation, ak}
{
r_2(m_1,...,\BlueText{f_k(a)(m_k)},...,m_n)
=f(a)(\BlueText{r_2(m_1,...,m_n)})
}

\AddEquation{polymorphism of representation, omega1}
{
r_{1k}(a_{k\cdot 1}...a_{k\cdot p}\omega_1)
=
r_{1k}(a_{k\cdot 1})
...
r_{1k}(a_{k\cdot p})\omega_1
}

\AddEquation{polymorphism of representation, omega2}
{
\begin{array}{r@{\,}l}
&\,r_2(m_1,...,m_{k\cdot 1}...m_{k\cdot p}\omega_2,...,m_n)
\\
=&\,
r_2(m_1,...,m_{k\cdot 1},...,m_n)
...
r_2(m_1,...,m_{k\cdot p},...,m_n)
\omega_2
\end{array}
}

\AddEq{gi:X->A 12}
{
\[
\begin{matrix}
g_i:X\rightarrow A_i&i=1,2
\end{matrix}
\]
}

\AddEq{Ai*Aj=0}
{
\(i\ne j=>A_i\cap A_j=\emptyset\)
}

\AddEq{A=+Ai}
{
A=\bigcup_{i=1}^nA_i\ \ \ A_i\in\mathcal S
}

\AddEq{A, A1 in S}
{
\(A\), \(A_1\in\mathcal S\),
\(A_1\subset A\)%
}

\AddEq{AB in S}
{
\(A\cap B\in\mathcal S\)
\labelItem{AB in S}
}

\AddEq{A, B in S}
{
\(A\), \(B\in\mathcal S\),
}

\AddEq{empty in S}
{
\(\emptyset\in\mathcal S\)
\labelItem{empty in S}
}

\AddEq{|-a|=|a|}
{
\item $\|-a\|=\|a\|$
\labelItem{|-a|=|a|}
}

\AddEq{|a+b|<=|a|+|b|}
{
\item $\|a+b\|\le\|a\|+\|b\|$
\labelItem{|a+b|<=|a|+|b|}
}

\AddEq{d->|d|}
{
\[d\in A\rightarrow \|d\|\in R\]
}

\AddEq{|a|=0}
{
\item $\|a\|=0$
\labelItem{|a|=0}
}

\AddEq{|a|>=0}
{
\item $\|a\|\ge 0$
\labelItem{|a|>=0}
}

\AddEq{fim M(X,A)}
{
$f_{i\cdot m}\in M(X,A)$, $i=1$, ..., $n$, $m=1$, ...,
}

\AddEq{|int h|<int|omega||f1n|}
{
\left\|\int_Xd\mu(x) h(x)\right\|\le
\int_Xd\mu(x)(\|\omega\|\|f_1(x)\|...\|f_n(x)\|)
}

\AddEq{fi:X->A}
{
\[
\begin{matrix}
f_i:X\rightarrow A&i=1,...,n
\end{matrix}
\]
}

\AddEq{|int h|<int|g1||g2|}
{
\left\|\int_Xd\mu(x) h(x)\right\|_2\le
\int_Xd\mu(x)(\|f\|\|g_1(x)\|_1\|g_2(x)\|_2)
}

\AddEq{h=fX g1 g2}
{
\[h=f_X(g_1)(g_2)\]
}

\AddEq{ak bk in}
{
$a_k\in A_k$, $b_k\in B_k$
}

\AddEq{B1n}
{
$B_1$, ..., $B_n$, $B$
}

\AddEq{A1n A}
{
$A_1$, ..., $A_n$, $A$
}

\AddEq{omega is distributive over addition}
{
\[
a_1...(a_i+b_i)...a_n\omega=a_1...a_i...a_n\omega+a_1...b_i...a_n\omega
\ \ \ i=1, ..., n
\]
}

\AddEq{Omega=...-o+}
{
$\Omega=\Omega_1\setminus\{+,*\}$.
}

\AddEq{c=ab}
{
c_1=a_1* b_1
}

\AddEq{a*(b1+b2)=}
{
\begin{align*}
a*(b_1+b_2)&=a*b_1+a*b_2\\
(b_1+b_2)*a&=b_1*a+b_2*a
\end{align*}
}

\AddEq{f(ab)=f(a)f(b)}
{
f(a* b)=f(a)\circ f(b)
}

\AddEq{c1=a+b}
{
\[c_1=a_1+ b_1\]
}

\AddEquation{Abelian multiplicative Omega group}
{
a* b=b* a
}

\AddEq{a * b + =}
{
a_{11}a_{21}+a_{12}a_{22}
=
(a_{11}+a_{12})(a_{21}+a_{22})
}

\DefEq
{
\begin{align*}
(a_{11}+a_{12})(a_{21}+a_{22})
&=(a_{11}+a_{12})a_{21}+(a_{11}+a_{12})a_{22}
\\
&=a_{11}a_{21}+a_{12}a_{21}
+a_{11}a_{22}+a_{12}a_{22}
\end{align*}
}
{a+a a+a}

\AddEq{a omega=omega a right}
{
(b_1...b_n\omega)* a=(b_1* a)...(b_n* a)\omega
}

\AddEq{a omega=omega a left}
{
a*(b_1...b_n\omega)=(a* b_1)...(a* b_n)\omega
}

\AddEq{Omega=...-o}
{
$\Omega=\Omega_1\setminus\{*\}$.
}

\AddEq{Omega=...-+}
{
\[\Omega=\Omega_1\setminus\{+\}\]
}

%% file: Prolog.Ref.tex

\ifx\UseRussian\Defined
\def\TheoremFollows{Теорема является следствием теоремы }
\else
\def\TheoremFollows{The theorem follows from the theorem }
\fi

\ePrints{4975-6381,9835-2163,0767-8264,7287-9339,9856-6693,5284-0163}
\ifx\Semafor\ValueOn%
\def\RefLinearMapA{4993-2400}%
\xRefDef{4993-2400}
\def\RefLinearMap{5114-6019}%
\xRefDef{5114-6019}
\fi
\ePrints{1601.03259,1506.00061}%
\Items{309618526,CACAA.06.121,1801.01628,CACAA.07,2207.06506,2021.11.26,2112.00613}%
\ifx\Semafor\ValueOn%
\def\RefLinearMap{1502.04063}%
\xRefDef{1502.04063}
\ePrints{2207.06506,1601.03259,2021.11.26,2112.00613,1506.00061}%
\ifx\Semafor\ValueOff%
\def\RefLinearMapA{0701.238}%
\xRefDef{0701.238}
\fi
\fi
\ePrints{348311191}%
\ifx\Semafor\ValueOn%
\def\RefLinearMap{2021.01.06}%
\xRefDef{2021.01.06}
\fi
\ePrints{1211.6965}%
\ifx\Semafor\ValueOn%
\def\RefLinearMap{1003.1544}%
\def\RefLinearMapA{0701.238}%
\xRefDef{1003.1544}
\xRefDef{0701.238}
\fi
\ePrints{MAlgebra4,MAlgebra2,2306.00880,2022.11.26,2023.01.07,357606033}%
\ifx\Semafor\ValueOn%
\def\RefLinearMap{2207.06506}%
\def\RefLinearMapA{2207.06506}%
\def\SectionTitle{}
\def\ProductType{}
\def\ColN{}
\ShowEq{def right}
\xRefDef{2207.06506}
\fi
\ePrints{5148-4632}%
\ifx\Semafor\ValueOn%
\def\RefLinearMap{8428-0408}%
\def\RefLinearMapA{8428-0408}%
\def\SectionTitle{}
\def\ProductType{}
\def\ColN{}
\ShowEq{def right}
\xRefDef{8428-0408}
\fi

\ePrints{330532713,1801.01628,2021.11.26,2112.00613}
\ifx\Semafor\ValueOn
\def\RefQuadratic{1506.00061}%
\xRefDef{1506.00061}
\fi%
\ePrints{0921-2370,0767-8264,5284-0163}
\ifx\Semafor\ValueOn
\def\RefQuadratic{7287-9339}%
\xRefDef{7287-9339}
\fi

\ePrints{330532713,2207.06506,2204.06320,2022.01.05,348311191}
\Items{2023.01.07}%
\ifx\Semafor\ValueOn
\def\PartA{}
\def\PartB{}
\def\SideRu{}
\ShowEq{def left}
\def\RefDiffEq{1801.01628}%
\xRefDef{1801.01628}
\fi
\ePrints{8428-0408}
\ifx\Semafor\ValueOn
\def\RefDiffEq{0767-8264}%
\xRefDef{0767-8264}
\fi

\ePrints{309618526,1801.01628,322019412,323236904,CACAA.07,330532713,323966352,2207.06506}
\Items{2021.01.06}
\ifx\Semafor\ValueOn
\def\RefCalculus{1601.03259}%
\xRefDef{1601.03259}
\fi%
\ePrints{9835-2163,0767-8264,9856-6693,7287-9339,0921-2370,8428-0408,5284-0163}
\ifx\Semafor\ValueOn
\def\RefCalculus{4975-6381}%
\xRefDef{4975-6381}
\fi

\ePrints{2020.06.01}
\ifx\Semafor\ValueOn
\def\RefCountable{1211.6965}%
\xRefDef{1211.6965}
\fi

\ePrints{1908.04418,323966352,324329551,0921-2370}
\ifx\Semafor\ValueOn
\def\RefGravity{0803.3276}%
\xRefDef{0803.3276}
\fi%
\ePrints{6860-2955}
\ifx\Semafor\ValueOn
\def\RefGravity{0803.3276}%
\xRefDef{0803.3276}
\fi

\AddEq{SetupRefRepresentation}
{
\def\RefRepresentation{}%
\ePrints{5114-6019,1502.04063,1305.4547}
\ifx\Semafor\ValueOn
\def\RefRepresentation{0912.3315}%
\xRefDef{0912.3315}
\fi
\ePrints{0767-8264,367250917}
\ifx\Semafor\ValueOn
\def\RefRepresentation{6860-2955}%
\xRefDef{6860-2955}
\fi
\ePrints{MAlgebra2,2204.06320,2022.01.05,357606033}
\ifx\Semafor\ValueOn
\def\RefRepresentation{1908.04418}%
\xRefDef{1908.04418}
\fi
}

\AddEq{SetupRefPolymorphism}
{
\def\RefPolymorphism{}%
\ePrints{1801.01628,CACAA.07}%
\ifx\Semafor\ValueOn
\def\RefPolymorphism{1502.04063}%
\fi
\ePrints{MAlgebra2,323236904}%
\ifx\Semafor\ValueOn
\def\RefPolymorphism{1502.04063}%
\xRefDef{1502.04063}
\fi
\ePrints{9856-6693}%
\ifx\Semafor\ValueOn
\def\RefPolymorphism{5114-6019}%
\fi
}

\AddEq{SetupRefOmegaNorm}
{
\def\RefTheoremOmegaNorm{}%
\ePrints{1102.5168,1502.04063}%
\ifx\Semafor\ValueOn%
\def\RefTheoremOmegaNorm{1305.4547}
\xRefDef{1305.4547}
\fi
\ePrints{CACAA.06.121}%
\ifx\Semafor\ValueOn%
\def\RefTheoremOmegaNorm{CACAA.04.001}
\xRefDef{CACAA.04.001}
\fi
}

\AddEq{SetupRefMeasure}
{
\def\RefMeasure{}%
\ePrints{9856-6693}
\ifx\Semafor\ValueOn
\def\RefMeasure{5410-9916}
\xRefDef{5410-9916}
\fi
\ePrints{323236904}
\ifx\Semafor\ValueOn
\def\RefMeasure{1310.5591}
\xRefDef{1310.5591}
\fi
}

\AddEq{SetupRef}
{
\ShowEq{SetupRefRepresentation}
\ShowEq{SetupRefOmegaNorm}
\ShowEq{SetupRefMeasure}
\ShowEq{SetupRefPolymorphism}
}

\AddEq{PreliminaryRefOn}
{
\ePrints{1506.00061,1601.03259}%
\ifx\Semafor\ValueOn%
\def\RefRepresentation{1502.04063}%
\fi
\ePrints{1801.01628,2207.06506}
\ifx\Semafor\ValueOn
\def\RefRepresentation{1908.04418}%
\xRefDef{1908.04418}
\fi
\ePrints{5148-4632,8428-0408,5284-0163}
\ifx\Semafor\ValueOn
\def\RefRepresentation{6860-2955}%
\xRefDef{6860-2955}
\fi
\ePrints{4975-6381}
\ifx\Semafor\ValueOn
\def\RefRepresentation{5114-6019}%
\fi

\ePrints{9835-2163,9856-6693,0767-8264,5284-0163}%
\ifx\Semafor\ValueOn%
\def\RefPolymorphism{5114-6019}%
\fi

\ShowEq{PreliminaryRefOmegaNormOn}
\ShowEq{PreliminaryRefMeasureOn}
}

\AddEq{PreliminaryRefOmegaNormOn}
{
\ePrints{5410-9916,4975-6381}%
\ifx\Semafor\ValueOn%
\def\RefTheoremOmegaNorm{5059-9176}%
\xRefDef{5059-9176}%
\fi
}

\AddEq{PreliminaryRefOmegaNormOff}
{
\def\RefTheoremOmegaNorm{}%
}

\AddEq{PreliminaryRefMeasureOff}
{
\def\RefMeasure{}%
}

\AddEq{PreliminaryRefMeasureOn}
{
\ePrints{1601.03259,309618526}
\ifx\Semafor\ValueOn
\def\RefMeasure{1310.5591}
\xRefDef{1310.5591}
\fi
\ePrints{4975-6381}
\ifx\Semafor\ValueOn
\def\RefMeasure{5410-9916}
\xRefDef{5410-9916}
\fi
\ePrints{CACAA.06.121}%
\ifx\Semafor\ValueOn%
\def\RefMeasure{CACAA.04.001}
\fi
}

\AddEq{PreliminaryRefOff}
{
\def\RefRepresentation{}%
\def\RefPolymorphism{}%
\ShowEq{PreliminaryRefOmegaNormOff}%
\ShowEq{PreliminaryRefMeasureOff}%
}

\ShowEq{SetupRef}%

\newcommand\ProofTheorem[2]
{
\begin{proof}
\TheoremFollows
\RefTheorem[#1]{#2}.
\end{proof}
}%

\newcommand\proofTheorem[3]
{
\begin{proof}
\TheoremFollows
\refTheorem[#1]{#2}{#3}.
\end{proof}
}%


\AddEq[1]{ref definition product rc}
{
\RefDefinition{row over column product}#1
}

\AddEq[1]{ref definition product cr}
{
\RefDefinition{column over row product}#1
}

\DefEq
{
\ePrints{1601.03259,4975-6381}%
\ifx\Semafor\ValueOff%
\EqRef[\RefCalculus]{derivative of Order n, algebra},
\else%
\EqRef{derivative of Order n, algebra},
\fi%
}
{ref derivative of Order n, algebra}

\AddEq{ref product in category}
{
\ePrints{8428-0408,2207.06506,5284-0163,1801.01628}
\ifx\Semafor\ValueOn
\RefDefinition{product in category, 2020},
\else
\RefDefinition{product in category},
\fi
}

\AddEq[1]{ref definition: left-side representation of algebra}
{
\ePrints{1310.5591,5059-9176,1305.4547}%
\ifx\Semafor\ValueOff%
\ePrints{5410-9916}%
\ifx\Semafor\ValueOn%
\RefDefinition{representation of algebra}#1
\else
\RefDefinition{left-side representation of algebra}#1
\fi
\else
\RefDefinition[\RefRepresentation]{left-side representation of algebra}#1
\fi
}

\DefEq
{
\ePrints{1310.5591}
\ifx\Semafor\ValueOff
\RefDefinition{morphism of representations of universal algebra}.
\else
\RefDefinition[0912.3315]{morphism of representations of universal algebra}.
\fi
}
{ref definition: morphism of representations of universal algebra}

\DefEq
{
\ePrints{1310.5591}%
\ifx\Semafor\ValueOff%
\RefDefinition{closed ball}
\else%
\RefDefinition[1305.4547]{closed ball}
\fi%
}
{ref definition: closed ball}

\DefEq
{
\ePrints{1310.5591}%
\ifx\Semafor\ValueOff%
\RefDefinition{open ball}\Pt
\else%
\RefDefinition[1305.4547]{open ball}\Pt
\fi%
}
{ref definition: open ball}

\DefEq
{
\ePrints{1310.5591}%
\ifx\Semafor\ValueOff%
\RefDefinition{open set},
\else%
\RefDefinition[1305.4547]{open set},
\fi%
}
{ref definition: open set}

\AddEq[1]{ref definition: fundamental sequence}
{
\ePrints{1310.5591}%
\ifx\Semafor\ValueOff%
\refDefinition{fundamental sequence}{normed Omega group}#1
\else%
\RefDefinition[1305.4547]{fundamental sequence, normed Omega group}#1
\fi%
}

\DefEq
{
\ePrints{1310.5591}%
\ifx\Semafor\ValueOff%
\RefDefinition{sequence converges uniformly}\Pt
\else%
\RefDefinition[1305.4547]{sequence converges uniformly}\Pt
\fi%
}
{ref definition: sequence converges uniformly}

\DefEq
{
\ePrints{1310.5591}%
\ifx\Semafor\ValueOff%
\RefDefinition{compact set},
\else%
\RefDefinition[1305.4547]{compact set},
\fi%
}
{ref definition: compact set}

\DefEq
{
\ePrints{1310.5591}%
\ifx\Semafor\ValueOff%
\RefDefinition{Omega group},
\else%
\RefDefinition[1305.4547]{Omega group},
\fi%
}
{ref definition: Omega group}

\DefEq
{
\ePrints{1310.5591}%
\ifx\Semafor\ValueOff%
\RefItem{|a|>=0}
\else%
\RefItem[1305.4547]{|a|>=0}
\fi%
}
{ref item |a|>=0}

\DefEq
{
\ePrints{1310.5591}%
\ifx\Semafor\ValueOff%
\RefItem{|a|=0}.
\else%
\RefItem[1305.4547]{|a|=0}.
\fi%
}
{ref item |a|=0}

\DefEq
{
\ePrints{1310.5591}%
\ifx\Semafor\ValueOff%
\RefItem{|a+b|<=|a|+|b|}\Pt
\else%
\RefItem[1305.4547]{|a+b|<=|a|+|b|}\Pt
\fi%
}
{ref item |a+b|<=|a|+|b|}

\DefEq
{
\ePrints{1310.5591}%
\ifx\Semafor\ValueOff%
\EqRef{|a omega|<|omega||a|1n}.
\else%
\EqRef[1305.4547]{|a omega|<|omega||a|1n}.
\fi%
}
{ref EqRef |a omega|<|omega||a|1n}

\DefEq
{
\ePrints{1310.5591}%
\ifx\Semafor\ValueOff%
\EqRef{|a-b|>|a|-|b|},
\else%
\EqRef[1305.4547]{|a-b|>|a|-|b|},
\fi%
}
{ref EqRef |a-b|>|a|-|b|}

\DefEq
{
\ePrints{1310.5591}%
\ifx\Semafor\ValueOff%
\EqRef{|fab|<|f||a||b|}.
\else%
\EqRef[1305.4547]{|fab|<|f||a||b|}.
\fi%
}
{ref |fab|<|f||a||b|}

\AddEq{ref tensor product and polylinear map}
{
\RefTheorem[\RefLinearMap]{there exists tensor product of modules},
\RefTheorem[\RefLinearMap]{tensor product and polylinear map}.
}

\DefEq
{
\ePrints{1502.04063,5114-6019}
\ifx\Semafor\ValueOff
\RefDefinition{reduced morphism of representations},
\else
\xRef[0912.3315]{remark: reduced morphism of representations},
\fi
}
{ref reduced morphism of representations}

\DefEq
{
\RefTheorem[1202.6021]{expand linear mapping, quaternion}.
}
{ref expand linear mapping, quaternion}

\DefEq
{
\RefTheorem[0912.4061]{ax-xa=1 quaternion algebra}.
}
{ref ax-xa=1 quaternion algebra}

\AddEq{ref theorem derivative, representation in algebra}
{
\ePrints{9835-2163,0767-8264}
\ifx\Semafor\ValueOn
\RefTheorem{derivative, representation in algebra},
\else
\RefTheorem[\RefCalculus]{derivative, representation in algebra},
\fi
}

\AddEq{ref partial derivatives equal}
{
\ePrints{9835-2163,0767-8264}
\ifx\Semafor\ValueOn
\EqRef{partial derivatives equal, Direct Sum of Banach Modules},
\else
\ePrints{CACAA.07}
\ifx\Semafor\ValueOff
\EqRef[0812.4763]{Gateaux partial derivatives equal, D vector space},
\else
\EqRef[CACAA.05.001]{partial derivatives equal, D vector space},
\fi
\fi
\EqRef{dM/dy=dN/dx yx},
\EqRef{d int f=f}.
}

\AddEq{ref aE, quaternion, Jacobian matrix}
{
\ePrints{CACAA.07}
\ifx\Semafor\ValueOn
\RefTheorem[CACAA.01.109]{aE, quaternion, Jacobian matrix},
\else
\RefTheorem[1107.1139]{aE, quaternion, Jacobian matrix},
\fi
}

\AddEq{ref Ea, quaternion, Jacobian matrix}
{
\ePrints{CACAA.07}
\ifx\Semafor\ValueOn
\RefTheorem[CACAA.01.109]{Ea, quaternion, Jacobian matrix},
\else
\RefTheorem[1107.1139]{Ea, quaternion, Jacobian matrix},
\fi
}

\DefEq
{
\RefTheorem[1003.1544]{quaternion conjugation}.
}
{ref quaternion conjugation}

\DefEq
{
\RefTheorem[1003.1544]{Quaternion over real field, matrix}.
}
{ref Quaternion over real field, matrix}

\DefEq
{
\RefTheorem[1601.03259]{representation of derivative, algebra A->B}.
}
{ref differentiable map A->B}

\DefEq
{
\RefTheorem[1601.03259]{map is continuous, derivative}.
}
{ref map is continuous, derivative}

\DefEq
{
\RefTheorem[1302.7204]{monomial of power k}.
}
{ref monomial of power k}

\DefEq
{
\RefTheorem[1302.7204]{map A(k+1)->pk otimes is polylinear map}.
}
{ref map A(k+1)->pk otimes is polylinear map}

\DefEq
{
\RefTheorem[1302.7204]{r=+q circ p}.
}
{ref r=+q circ p}

\DefEq
{
\RefTheorem[1302.7204]{r=+q circ p 1}.
}
{ref r=+q circ p 1}

\DefEq
{
\RefTheorem[1105.4307]{Re Im A->F}.
}
{ref Re Im A->F}

\DefEq
{
\RefTheorem[1105.4307]{conjugation in algebra}.
}
{ref conjugation in algebra}

\DefEq
{
\RefTheorem[1105.4307]{structural constant of algebra with unit}.
}
{ref structural constant of algebra with unit}

%% file: Biblio.English.tex
\OpenBiblio


\BiblioItem{Doctor Ouch}
{
Kornei Chukovsky. Doctor Ouch.
\\
Translator and illustrator Jan Seabaugh.
\\
Viveca Smith Publishing, 2004, ISBN-10: 0974055107.
}%

\BiblioItem{Einstein: Electrodynamics of Moving Bodies}
{
Albert Einstein,
On the Electrodynamics of Moving Bodies, 1905,
\\
The Principle of Relativity: A Collection of Original
Memoirs on the Special and General Theory of Relativity , 37 - 65,
\\
Courier Dover Publications, 1952; ISBN-13: 978-0486600819
\\
Zur Elektrodynamik der bewegter K\"orper. Ann. Phys., 1905, 17, 891-921. 
}%

\BiblioItem{Einstein: On the Relativity Principle}
{
Albert Einstein,
On the Relativity Principle and the Conclusions Drawn from It, 1907,
\\
The Collected Papers of Albert Einstein, Volume 2:
The Swiss Years: Writings, 1900-1909. English translation. 252 - 311.
\\
Anna Beck, translator; Peter Havas, consultant.
Princeton University Press, 1989; ISBN-13: 9780691085494
\\
\"Uber das Relativit\"atsprinzip und die aus demselben gezogenen Folgerungen. 
Jahrb. d. Radioaktivit\"at u. Elektronik, 1907, 4, 411-462. 
}%

\BiblioItem{Einstein: Foundations of general relativity}
{
Albert Einstein,
Die Grundlage der allgemeinen Relativit\"atstheorie,
Ann. Phys., 1916, {\bf 49}, 769 - 822,\\
Einstein's Annalen Papers: The Complete Collection 1901-1922,
edited by J\"urgen Renn, 517 - 571,\\
Wiley-VCH Verlag GmbH \& Co. KGaA, 2005
}%

\BiblioItem{Einstein: Geometry and Experience}
{
Albert Einstein, Geometry and Experience, (1921)\\
Albert Einstein, Sidelights on Relativity, 25 - 56,\\
Courier Dover Publications, 1983
}%

\BiblioItem{Einstein: Main problems of general relativity}
{
Albert Einstein,
Grundgedanken und Probleme der Relativit\"atstheorie, (1923),\\
Nobelstiftelsen, Les Prix Nobel en 1921 - 1922,
Imprimerie Royale, Stockholm, 1923
}%

\BiblioItem{Einstein: Noneuclidean Geometry and Physics}
{
Albert Einstein,
Nichtenklidische Geometrie in der Physik Neue Rundschan, (1925)
Berlin, S. 16 - 20
}%

\BiblioItem{Einstein: Isaak Newton}
{
Albert Einstein,
Isaak Newton, 1927,
Out of My Later Years, 
Citadel Press, 1995, 219 - 223
}%

\BiblioItem{Einstein: On Science}
{
Albert Einstein,
On Science, 
Cosmic Religion, with Other Opinions and Aphorisms,142 - 146,
New York, 1931, 97 - 103
}%

\BiblioItem{Einstein: Autobiographical Notes}
{
Albert Einstein,
Autobiographical Notes, 1949,\\
Paul A. Schilpp, editor, Albert Einstein: Philosopher-Scientist,
Evanston, 
Illinois, The Library of Living Philosophers, 1949, 1 - 95
}%

\BiblioItem{Feynman 1}
{
Richard Phillips Feynman, Robert B. Leighton, Matthew Linzee Sands.
The Feynman lectures on physics: Volume 1.
Mainly Mechanics, Radiation, and Heat.
Addison\Hyph Wesley, 1965.
}%

\BiblioItem{0538731877}
{
James Shipman, Jerry D. Wilson and Aaron Todd.
Introduction to Physical Science.
Cengage Learning, 2009; ISBN 0538731877.
}%

\BiblioItem{Cite: 104}
{
Cite 104, Source unknown
}%

\BiblioItem{Ghez}
{
Ghez et al.,
The First Measurement of Spectral Lines in a Short-Period Star Bound to the Galaxy's Central Black Hole: A Paradox of Youth,
\href{http://www.journals.uchicago.edu/ApJ/journal/issues/ApJL/v586n2/16990/brief/16990.abstract.html}{ApJL, 586, L127} (2003),
eprint \href{http://arxiv.org/abs/astro-ph/0302299}{arXiv:astro-ph/0302299} (2003)
}%

\BiblioItem{Schodel}
{
R. Sch\"odel et al.,
A star in a 15.2-year orbit around the supermassive black hole at the centre of the Milky Way,
\href{http://www.nature.com/cgi-taf/DynaPage.taf?file=/nature/journal/v419/n6908/abs/nature01121_fs.html}{Nature 419, 694} (2002)
}%

\BiblioItem{Mielke}
{
Eckehard W. Mielke, Affine generalization of the Komar complex of general relativity,
\href{http://prola.aps.org/searchabstract/PRD/v63/i4/e044018}{Phys. Rev. D 63, 044018} (2001)
}%

\BiblioItem{Obukhov}
{
Yu. N. Obukhov and J. G. Pereira, Metric\hyph affine approach to teleparallel gravity,
\href{http://scitation.aip.org/getabs/servlet/GetabsServlet?prog=normal&id=PRVDAQ000067000004044016000001&idtype=cvips&gifs=Yes}
{Phys. Rev. D 67, 044016} (2003),
eprint \href{https://arxiv.org/abs/gr-qc/0212080}{arXiv:gr-qc/0212080} (2002)
}%

\BiblioItem{2110.04354}
{
Dimitri Gurevich, Varvara Petrova, Pavel Saponov,
q-Casimir and q-cut-and-join operators related to Reflection Equation Algebras,
eprint \href{https://arxiv.org/abs/2110.04354}{arXiv:gr-qc/2110.04354} (2021)
}%

\BiblioItem{Sardanashvily}
{
Giovanni Giachetta, Gennadi Sardanashvily, Dirac Equation in Gauge and Affine-Metric Gravitation Theories,
eprint \href{http://arxiv.org/abs/gr-qc/9511035}{arXiv:gr-qc/9511035} (1995)
}%

\BiblioItem{Gauge}
{
Frank Gronwald and Friedrich W. Hehl, On the Gauge Aspects of Gravity, eprint
\href{http://arxiv.org/abs/gr-qc/9602013}{arXiv:gr-qc/9602013} (1996)
}%

\BiblioItem{Neeman}
{
Yuval Neeman, Friedrich W. Hehl, Test Matter in a Spacetime with Nonmetricity, eprint
\href{http://arxiv.org/abs/gr-qc/9604047}{arXiv:gr-qc/9604047} (1996)
}%

\BiblioItem{torsion}
{
F. W. Hehl, P. von der Heyde, G. D. Kerlick, and J. M. Nester,
General relativity with spin and torsion: Foundations and prospects,\\
\href{http://prola.aps.org/abstract/RMP/v48/i3/p393_1}{Rev. Mod. Phys. 48, 393} (1976)
}%

\BiblioItem{Megged}
{
O. Megged, Post-Riemannian Merger of Yang-Mills Interactions with Gravity,
eprint \href{http://arxiv.org/abs/hep-th/0008135}{arXiv:hep-th/0008135} (2001)
}%


\BiblioItem{gr-qc-9604027}
{
Yu.N. Obukhov, E.J. Vlachynsky, W. Esser, R. Tresguerres and F.W. Hehl,
An exact solution of the metric\hyph affine gauge theory with dilation, shear, and spin charges,
eprint \href{http://arxiv.org/abs/gr-qc/9604027}{arXiv:gr-qc/9604027} (1996)
}%

\BiblioItem{4419-7514}
{
Mari\'an Fabian, Petr Habala, Petr H\'ajek, Vicente Montesinos, V\'aclav Zizler.
Banach Space Theory: The Basis for Linear and Nonlinear Analysis.
\\
Springer; New York, 2010; ISBN-13: 978-1441975140
}%

\BiblioItem{Weinberg I}
{
Steven Weinberg.
The Quantum Theory of Fields. Volume I. Foundations.
Cambridge university press, 1995
}%

\BiblioItem{Weinberg II}
{
Steven Weinberg.
The Quantum Theory of Fields. Volume II. Modern applications.
Cambridge university press, 1996
}%

\BiblioItem{Reinhardt}
{
Walter Greiner, Joachim Reinhardt. Field Quantization. Springer.
}%

\BiblioItem{978-3540875604}
{
Walter Greiner, Joachim Reinhardt. Quantum Electrodynamics. Springer, 2009.
}%

\BiblioItem{978-1898563020}
{
H. Robert Mills. Practical Astronomy. Woodhead Publishing, 1994. ISBN-13: 978-1898563020.
}%

\BiblioItem{Landau I}
{
L. D. Landau, E. M. Lifshich.
Course of theoretical physics, volume 1.
Mechanics.
\\
Translated from the Russian by J. B. Sykes and J. S. Bell.
Pergamon Press, 1969
}%

\BiblioItem{Landau}
{
L. D. Landau, E. M. Lifshich, The classical theory of fields.
\\
Translated from the Russian by Morton Hamermesh.
Pergamon Press, 1971
}%

\BiblioItem{Landau III}
{
L. D. Landau, E. M. Lifshich,
Course of Theoretical Physics, Volume 3.
Quantum Mechanics Non-Relativistic Theory, Third Edition.
\\
Translated from the Russian by J. B. Sykes and J. S. Bell.
Butterworth-Heinemann, 1981, ISBN 978-0750635394.
}%

\BiblioItem{Wheeler}
{
Ignazio Ciufolini, John Wheeler. Gravitation and Inertia.
Princeton university press.
}%

\BiblioItem{Gravitation MTW}
{
Charles W. Misner, Kip S. Thorne, John Archibald Wheeler.
Gravitation.
W. H. Freeman and Company, San Francisco, 1973.
}%

\BiblioItem{Anderson98}
{
J. D. Anderson, P. A. Laing, E. L. Lau, A. S. Liu, M. M. Nieto, and S. G. Turyshev,
Indication, from Pioneer 10/11, Galileo, and Ulysses Data, of an Apparent Anomalous, Weak, Long-Range Acceleration,
\href{http://prola.aps.org/abstract/PRL/v81/i14/p2858_1}{Phys. Rev. Lett. 81, 2858}, (1998),
eprint \href{http://arxiv.org/abs/gr-qc/9808081}{arXiv:gr-qc/9808081} (1998)
}%

\BiblioItem{Anderson02}
{
J. D. Anderson, P. A. Laing, E. L. Lau, A. S. Liu, M. M. Nieto, and S. G. Turyshev,
Study of the anomalous acceleration of Pioneer 10 and 11,
\href{http://prola.aps.org/searchabstract/PRD/v65/i8/e082004}{Phys. Rev. D 65, 082004, 50 pp.}, (2002),
eprint \href{http://arxiv.org/abs/gr-qc/0104064}{arXiv:gr-qc/0104064} (2001)
}%


\BiblioItem{H. Aslaksen}
{
H. Aslaksen.  Quaternionic determinants \textit{Math.
Intelligencer} {\bf 18}(3), pp.57-65, (1996).
}%

\BiblioItem{L. Chen: Definition of determinant}
{
L. Chen, Definition of determinant and Cramer solutions over
quaternion field, \textit{Acta Math. Sinica (N.S.)} {\bf 7},
pp.171-180, (1991).
}%

\BiblioItem{L. Chen: Inverse matrix}
{
L. Chen,
Inverse matrix and properties of double determinant over quaternion
field, \textit{Sci. China, Ser. A} {\bf 34}, pp.528-540, (1991).
}%

\BiblioItem{N. Cohen S. De Leo}
{
N. Cohen, S. De Leo, The quaternionic determinant, \textit{The Electronic Journal Linear
Algebra} {\bf 7}, pp.100-111, (2000).
}%

\BiblioItem{Dyson: Quaternion determinants}
{
F. J. Dyson, Quaternion determinants, \textit{Helvetica Phys.
Acta} {\bf 45}, pp. 289-302, (1972).
}%

\BiblioItem{Melvin Hausner}
{
Melvin Hausner,
A Vector Space Approach to Geometry,
Dover Publications, 1998
}%

\BiblioItem{Serge Lang}
{
Serge Lang,
Algebra, Springer, 2002
}%

\BiblioItem{9780534423230}
{
Charles Lanski.
Concepts In Abstract Algebra.
American Mathematical Soc., 2005, ISBN 978-0534423230
}%

\BiblioItem{Burris Sankappanavar}
{
S. Burris, H.P. Sankappanavar,
A Course in Universal Algebra, Springer-Verlag (March, 1982),
\\eprint
\href{http://www.math.uwaterloo.ca/~snburris/htdocs/ualg.html}
{http://www.math.uwaterloo.ca/~snburris/htdocs/ualg.html}
\\(The Millennium Edition)
}%

\BiblioItem{Shilov single 12}
{
G. E. Shilov,
Calculus, Single Variable Functions, Parts 1 - 2,
Moscow, Nauka, 1969
}%

\BiblioItem{Shilov single 3}
{
G. E. Shilov,
Calculus, Single Variable Functions, Part 3,
Moscow, Nauka, 1970
}%

\BiblioItem{Shilov}
{
G. E. Shilov,
Calculus, Multivariable Functions,
Moscow, Nauka, 1972
}%

\BiblioItem{Kolmogorov Fomin}
{
A. N. Kolmogorov and S. V. Fomin.
Introductory Real Analysis.
\\
Translated and edited by Richard A. Silverman.
\\
Dover Publication, 1975, ISBN-13: 978-0486612263
}%

\BiblioItem{Lebedev Vorovich}
{
L. P. Lebedev, I. I. Vorovich,
Functional Analysis in Mechanics,
Springer, 2002
}%

\BiblioItem{8176-4374}
{
Mariano Giaquinta, Giuseppe Modica,
Mathematical Analysis: Linear and Metric Structures and Continuity.
\\
Springer, 2007, ISBN-13: 978-0-8176-4374-4
}%

\BiblioItem{319-09203}
{
Ralf Schiffler,
Quiver Representations.
\\
Springer, 2014, ISBN-13: 978-3-319-09203-4
}%

\BiblioItem{511-34545}
{
Ibrahim Assem, Daniel Simson, Andrzej Skowro\'nski,\\
Elements of the Representation Theory
of Associative Algebras.
\\
Volume 1: Techniques of Representation Theory.
\\
Springer, 2014, ISBN-13: 978-3-319-09203-4
}%

\DefBiblioItem{Rashevsky}

\BiblioItem
{Kurosh: High Algebra}
{
A. G. Kurosh, Higher Algebra,
\\
George Yankovsky translator,
\\
Mir Publishers, 1988, ISBN: 978-5030001319
}%

\BiblioItem
{Kurosh: General Algebra}
{
A. G. Kurosh, Lectures on General Algebra,
Chelsea Pub Co, 1965 
}%

\BiblioItem{Sabinin: Smooth Quasigroups}
{
Lev V. Sabinin, Smooth Quasigroups and Loops,
Kluwer Academic Publisher, 1999 
}%

\BiblioItem{978-0-8176-8384-9}
{
Garret Sobczyk, New Foundations in Mathematics: The Geometric Concept of Number,
\\
Springer, 2013, ISBN: 978-0-8176-8384-9
}%

\BiblioItem{978-0538497817}
{
James Stewart, Calculus,
\\
Cengage Learning, 2012, ISBN: 978-0-538-49781-7
}%

\BiblioItem{470-38334}
{
William E. Boyce, Richard C. DiPrima,
Elementary Differential Equations and Boundary Value Problems,
\\
John Wiley \& Sons, Inc., 2009, ISBN 978-0-470-38334-6
}%

\BiblioItem{Dubrovin Fomenko Novikov part 1}
{
B. A. Dubrovin, A. T. Fomenko, S. P. Novikov,
Modern Geometry - Methods and Applications,\\
Part I, The Geometry of Surfaces, Transformation Groups, and Fields,\\
Translated by Robert G. Burns,\\
Springer - New York, 1992
}%

\BiblioItem{Dubrovin Fomenko Novikov part 2}
{
B. A. Dubrovin, A. T. Fomenko, S. P. Novikov,
Modern Geometry - Methods and Applications,
Part II: The Geometry and Topology of Manifolds,\\
Translated by Robert G. Burns,\\
Springer - New York, 1985
}%

\BiblioItem{Kobayashi Nomizu vol 1}
{
Kobayashi S, Nomizu K,
Foundations of Differential Geometry, volume I,\\
Interscience Publishers, 1963
}%

\BiblioItem{Lichnerowicz}
{
Andre Lichnerowicz,
Global Theory of Connections and Holonomy Groups,\\
Kluwer Academic Publishers, 1976, ISBN-13: 978-9028604964
}%

\DefBiblioItem{Korn}%

\BiblioItem{Hocking Young Topology}
{
John G. Hocking, Gail S. Young,
Topology,\\
Courier Dover Publications, 1988
}%

\BiblioItem{Olver: Lie groups to differential equations}
{
Peter J. Olver,
Applications of Lie groups to differential equations,\\
Springer, 2000
}%

\BiblioItem{1708.01190}
{
Nathan BeDell,
Doing Algebra over an Associative Algebra,
\\
eprint \href{https://arxiv.org/abs/1708.01190}{arXiv:1708.01190} (2017)
}%

\BiblioItem{Karthika Viji 2021}
{
S. Karthika, M. Viji,
Unit elements in the path algebra of an acyclic quiver,
Indian J Pure Appl Math\\
Published online: 10 June 2021
}%

\BiblioItem{Tartaglia}
{
Angelo Tartaglia and Matteo Luca Ruggiero,
Angular Momentum Effects in Michelson\Hyph Morley Type Experiments,
Gen.Rel.Grav. 34, 1371-1382 (2002),\\
eprint \href{http://arxiv.org/abs/gr-qc/0110015}{arXiv:gr-qc/0110015} (2001)
}%

\BiblioItem{Tomozawa}
{
Yukio Tomozawa, Speed of Light in Gravitational Fields, eprint
\href{http://arxiv.org/abs/astro-ph/0303047}{arXiv:astro-ph/0303047} (2004)
}%

\BiblioItem{Magueijo}
{
Joao Magueijo,
Covariant and locally Lorentz-invariant varying speed of light theories,
\href{http://prola.aps.org/abstract/PRD/v62/i10/e103521}{Phys. Rev. D 62, 103521} (2000),
eprint \href{http://arxiv.org/abs/gr-qc/0007036}{arXiv:gr-qc/0007036} (2000)
}%

\BiblioItem{Bassett}
{
Bruce A. Bassett, Stefano Liberati, Carmen Molina-Paris, and Matt Visser,
Geometrodynamics of variable-speed-of-light cosmologies,
\href{http://prola.aps.org/abstract/PRD/v62/i10/e103518}{Phys. Rev. D 62}, 103518 (2000),
eprint \href{http://arxiv.org/abs/astro-ph/0001441}{arXiv:astro-ph/0001441} (2000)
}%

\BiblioItem{C.A. Deavours The Quaternion Calculus}
{
C.A. Deavours, The Quaternion Calculus, 
American Mathematical Monthly, {\bf 80} (1973), pp. 995 - 1008
}%

\BiblioItem{Straumann}
{
Lochlainn O'Raifeartaigh and Norbert Straumann,
Gauge theory: Historical origins and some modern developments,
\href{http://prola.aps.org/abstract/RMP/v72/i1/p1_1}{Rev. Mod. Phys. 72, 1} (2000)
}%

\BiblioItem{Lammerzahl}
{
Claus L\"ammerzahl, Mark P. Haugan,
On the interpretation of Michelson\Hyph Morley experiments,
{Phys. Lett. A282 223-229} (2001),\\
eprint \href{http://arxiv.org/abs/gr-qc/0103052}{arXiv:gr-qc/0103052} (2001)
}%

\BiblioItem{0305117}
{
Holger Mueller, Sven Herrmann, Claus Braxmaier, Stephan Schiller, Achim Peters.
Modern Michelson-Morley Experiment using Cryogenic Optical Resonators.
eprint \href{http://arxiv.org/abs/physics/0305117}{arXiv:physics/0305117} (2003)
\\
Phys. Rev. Lett. 91:020401, 2003
}%

\BiblioItem{0706.2031}
{
Holger Mueller, Paul Louis Stanwix, Michael Edmund Tobar,
Eugene Ivanov, Peter Wolf, Sven Herrmann, Alexander Senger,
Evgeny Kovalchuk, Achim Peters.
Relativity tests by complementary rotating Michelson-Morley experiments.
eprint \href{http://arxiv.org/abs/0706.2031}{arXiv:0706.2031 [physics.class-ph]} (2006)
\\
Phys. Rev. Lett. 99:050401, 2007
}%

\BiblioItem{1008.1205}
{
M. Nagel, K. M\"ohle, K. D\"oringshoff, S. Herrmann, A. Senger, E.V. Kovalchuk, A. Peters.
Testing Lorentz Invariance by Comparing Light Propagation in Vacuum and Matter.
eprint \href{http://arxiv.org/abs/1008.1205}{arXiv:1008.1205 [physics.ins-det]} (2010)
}%

\BiblioItem{1109.4897}
{
The OPERA Collaboration.
Measurement of the neutrino velocity with the OPERA detector in the CNGS beam.
eprint \href{http://arxiv.org/abs/1109.4897}{arXiv:1109.4897 [hep-ex]} (2011)
}%

\BiblioItem{Ranada}
{
Antonio F. Ranada,
Pioneer acceleration and variation of light speed: experimental situation,
eprint \href{http://arxiv.org/abs/gr-qc/0402120}{arXiv:gr-qc/0402120} (2004)
}%

\BiblioItem{Gelfand Minlos: rotation and Lorentz groups}
{
Izrail Moiseevich Gelfand, Robert Adolfovich Minlos,
Representations of the rotation and Lorentz groups and their applications;\\
Engl. transl. ed. H. K. Farahat; Transl. by G. Cummins and T. Boddongton;\\
Pergamon Press, 1963
}%

\DefBiblioItem{math.QA-0208146}%

\DefBiblioItem{q-alg-9705026}

\BiblioItem{Gelfand Retakh 1991}
{
I. Gelfand and V. Retakh, Determinants of Matrices over Noncommutative Rings, Funct.
Anal. Appl. 25 (1991), no. 2, 91-102
}%

\BiblioItem{Gelfand Retakh 1992}
{
I. Gelfand and V. Retakh, A Theory of Noncommutative Determinants and Characteristic
Functions of Graphs, Funct. Anal. Appl. 26 (1992), no. 4, 1-20
}%

\BiblioItem{hep-th-9407124}
{
I. M. Gelfand, D. Krob, A. Lascoux, B. Leclerc, V.S. Retakh and J.-Y. Thibon,
Noncommutative symmetric functions,\\
eprint \href{http://arxiv.org/abs/hep-th/9407124}{arXiv:hep-th/9407124} (1994)
}%

\BiblioItem{0911.4454}
{
Vladimir Retakh,
From factorizations of noncommutative polynomials to combinatorial topology,\\
eprint \href{http://arxiv.org/abs/0911.4454}{arXiv:0911.4454} (2009)
}%

\BiblioItem{Naimark Shtern: Theory of group representations}
{
Mark Aronovich Naimark, Aleksandr Isaakovich Shtern,
Theory of group representations;\\
Heidelberg, 1982
}%

\BiblioItem{Barut Raczka: Theory of group representations}
{
Asim Orhan Barut; Ryszard R\c{a}czka;
Theory of group representations and applications;\\
World Scientific Publishing Co. Pre. Ltd., 1986
}%

\BiblioItem{Mihalev Pilz: concise handbook of algebra}
{
Aleksandr Vasilevich Mikhalev; G\"{u}nter Pilz;
The concise handbook of algebra;\\
Kluwer Academic Publishers, 2002
}%

\BiblioItem{McCrimmon: Jordan Algebras}
{
Kevin McCrimmon;
A Taste of Jordan Algebras;\\
Springer, 2004
}%

\BiblioItem{Zharinov: foundation of mathematical physics}
{
V. V. Zharinov,
Algebraic and geometric foundation of mathematical physics,\\
Lecture courses of the scientific and educational center, 9, Steklov Math. Institute of RAS,\\
Moscow, 2008
}%

\BiblioItem{Shafarevich: Basic notions of algebra}
{
I. R. Shafarevich,
Basic notions of algebra,\\
Translated from the Russian by M. Reid,\\
Springer, 2005
}%

\BiblioItem{Coppel: Number Theory}
{
W.A. Coppel,
Number Theory: An Introduction to Mathematics,\\
Springer, 2009
}%

\BiblioItem{978-0486497952}
{
Michael J. Field,
Differential Calculus and Its Applications,\\
Dover Publications, 2012; ISBN-13: 978-0486497952
}%

\BiblioItem{Elsgolts: Differential Equations}
{
Lev Elsgolts,
Differential Equations and the Calculus of Variations,\\
Translated from the Russian by George Yankovsky,\\
MIR Publishers, Moscow, 1977
}%

\BiblioItem{Baez Huerta: algebra of grand unified theories}
{
John Baez; John Huerta;
The algebra of grand unified theories;\\
Bull. Amer. Math. Soc. {\bf 47} (2010), 483-552
}%

\BiblioItem{J. Fan: Determinants}
{
J. Fan, Determinants and multiplicative functionals
on quaternion matrices, \textit{Linear Algebra and Its
Applications} {\bf 369}, pp. 193-201, (2003).
}%

\BiblioItem{Carl Faith 1}
{
Carl Faith, Algebra: Rings, Modules and Categories I,
Springer - Verlag, Berlin - Heidelberg - New York, 1973
}%

\BiblioItem{Gilson Nimmo Ohta}
{
 C.R.Gilson, J.J.C.Nimmo, Y.Ohta, Quasideterminant solutions of a non-Abelian Hirota-Miwa
 equation, \textit{Journal of Physics A: Mathematical and Theoretical} {\bf 40}(42), pp.
 12607-12617,(2007).
}%

\BiblioItem{Haider Hassan}
{
B. Haider, M. Hassan, Quasideterminant solutions of an integrable chiral model in two
 dimensions, \textit{Journal of Physics A: Mathematical and Theoretical} {\bf 42} (35), art. no.
 355211, (2009).
}%



\BiblioItem{0702447}
{
I.I. Kyrchei, Cramer's rule for quaternion systems of linear equations,
\textit{Journal of Mathematical Sciences} {\bf 155}(6), 839-858, (2008).
 Translated from  \textit{Fundamental and Appl. Math.}
 {\bf 13}(4), pp.67-94, (2007). (in Russian)\\
eprint
\href{http://arxiv.org/abs/math/0702447}{arXiv:math.RA/0702447}
(2007)
}%

\BiblioItem{1004.4380}
{
I.I. Kyrchei, Cramer's rule for some quaternion matrix
    equations,  \textit{Applied Mathematics and Computation} {\bf 217}(5), pp.2024-2030, (2010).\\eprint
\href{http://arxiv.org/abs/1004.4380
}{arXiv:math.RA/arXiv:1004.4380 } (2010)
}%

\BiblioItem{1005.0736}
{
I.I. Kyrchei,Determinantal representations of the Moore-Penrose inverse
 over the quaternion skew field and corresponding Cramer's rules,
 \\
eprint
\href{http://arxiv.org/abs/1005.0736}{arXiv:math.RA/1005.0736}
(2010)
}%

\BiblioItem{0412.391}
{
Aleks Kleyn,
Basis Manifold,
eprint \href{http://arxiv.org/abs/math.DG/0412391}{arXiv:math.DG/0412391} (2007)
}%

\BiblioItem{0405.027}
{
Aleks Kleyn,
Reference Frame in General Relativity,\\
eprint \href{http://arxiv.org/abs/gr-qc/0405027}{arXiv:gr-qc/0405027} (2008)
}%

\BiblioItem{0405.028}
{
Aleks Kleyn, Metric\hyph Affine Manifold,\\
eprint \href{http://arxiv.org/abs/gr-qc/0405028}{arXiv:gr-qc/0405028} (2008)
}%

\BiblioItem{0612.111}
{
Aleks Kleyn,
Biring of Matrices,\\
eprint \href{http://arxiv.org/abs/math.OA/0612111}{arXiv:math.OA/0612111} (2007)
}%

\BiblioItem{0701.238}
{
Aleks Kleyn,
Lectures on Linear Algebra over Division Ring,\\
eprint \href{http://arxiv.org/abs/math.GM/0701238}{arXiv:math.GM/0701238} (2010)
}%

\BiblioItem{0702.561}
{
Aleks Kleyn,
Fibered Universal Algebra,\\
eprint \href{http://arxiv.org/abs/math.DG/0702561}{arXiv:math.DG/0702561} (2007)
}%

\BiblioItem{math.RA-0501237}
{
Aleks Kleyn,
Vector Space Over Division Ring,\\
eprint \href{http://arxiv.org/abs/math.RA/0412391}{arXiv:math.RA/0501237} (2007)
}%

\BiblioItem{math.RA-0501237v1}
{
Aleks Kleyn,
Module Over Division Ring, version 1,\\
eprint \href{http://arxiv.org/abs/math/0501237v1}{arXiv:math.RA/0501237v1} (2005)
}%

\BiblioItem{0707.2246}
{
Aleks Kleyn,
Fibered Correspondence,\\
eprint \href{http://arxiv.org/abs/0707.2246}{arXiv:0707.2246} (2007)
}%

\BiblioItem{0803.2620}
{
Aleks Kleyn,
Morphism of \Ts{T}Representations,\\
eprint \href{http://arxiv.org/abs/0803.2620}{arXiv:0803.2620} (2008)
}%

\BiblioItem{0803.3276}
{
Aleks Kleyn,
Lorentz Transformation and General Covariance Principle,\\
eprint \href{http://arxiv.org/abs/0803.3276}{arXiv:0803.3276} (2009)
}%

\DefBiblioItem{0812.4763}%

\BiblioItem{0906.0135}
{
Aleks Kleyn,
Introduction into Geometry over Division Ring,\\
eprint \href{http://arxiv.org/abs/0906.0135}{arXiv:0906.0135} (2010)
}%

\BiblioItem{0909.0855}
{
Aleks Kleyn,
Quaternion Rhapsody,\\
eprint \href{http://arxiv.org/abs/0909.0855}{arXiv:0909.0855} (2010)
}%

\BiblioItem{0912.3315}
{
Aleks Kleyn,
Representation of Universal Algebra,\\
eprint \href{http://arxiv.org/abs/0912.3315}{arXiv:0912.3315} (2009)
}%

\BiblioItem{0912.4061}
{
Aleks Kleyn,
Linear Equation in Finite Dimensional Algebra,\\
eprint \href{http://arxiv.org/abs/0912.4061}{arXiv:0912.4061} (2010)
}%

\BiblioItem{1001.4852}
{
Aleks Kleyn,
The Matrix of Linear Maps,\\
eprint \href{http://arxiv.org/abs/1001.4852}{arXiv:1001.4852} (2010)
}%

\BiblioItem{1003.1544}
{
Aleks Kleyn,
Linear Maps of Free Algebra,\\
eprint \href{http://arxiv.org/abs/1003.1544}{arXiv:1003.1544} (2010)
}%

\BiblioItem{1006.2597}
{
Aleks Kleyn,
The G\^ateaux Derivative and Integral over Banach Algebra,\\
eprint \href{http://arxiv.org/abs/1006.2597}{arXiv:1006.2597} (2010)
}%

\BiblioItem{1011.3102}
{
Aleks Kleyn,
Polylinear Map of Free Algebra,\\
eprint \href{http://arxiv.org/abs/1011.3102}{arXiv:1011.3102} (2010)
}%

\BiblioItem{1102.1776}
{
Aleks Kleyn, Ivan Kyrchei,
Correspondence between Row\Hyph Column Determinants
and Quasideterminants of Matrices over Quaternion Algebra,\\
eprint \href{http://arxiv.org/abs/1102.1776}{arXiv:1102.1776} (2011)
}%

\BiblioItem{1104.5197}
{
Aleks Kleyn,
$C^*$-Rhapsody,\\
eprint \href{http://arxiv.org/abs/1104.5197}{arXiv:1104.5197} (2011)
}%

\BiblioItem{1105.4307}
{
Aleks Kleyn,
Algebra with Conjugation,\\
eprint \href{http://arxiv.org/abs/1105.4307}{arXiv:1105.4307} (2011)
}%

\BiblioItem{1107.1139}
{
Aleks Kleyn,
Linear Maps of Quaternion Algebra,\\
eprint \href{http://arxiv.org/abs/1107.1139}{arXiv:1107.1139} (2011)
}%

\BiblioItem{1107.5037}
{
Aleks Kleyn,
Orthogonal Basis and Motion in Finsler Geometry,\\
eprint \href{http://arxiv.org/abs/1107.5037}{arXiv:1107.5037} (2011)
}%

\BiblioItem{1111.6035}
{
Aleks Kleyn,
Basis of Representation of Universal Algebra,\\
eprint \href{http://arxiv.org/abs/1111.6035}{arXiv:1111.6035} (2011)
}%

\BiblioItem{1201.4158}
{
Aleks Kleyn, Alexandre Laugier,
Orthonormal Basis in Minkowski Space,\\
eprint \href{http://arxiv.org/abs/1201.4158}{arXiv:1201.4158} (2012)
}%

\BiblioItem{1202.6021}
{
Aleks Kleyn,
Maps of Conjugation of Quaternion Algebra,\\
eprint \href{http://arxiv.org/abs/1202.6021}{arXiv:1202.6021} (2012)
}%

\BiblioItem{1206.0200}
{
Aleks Kleyn,
Algebra of Fractions of Algebra with Conjugation,\\
eprint \href{http://arxiv.org/abs/1206.0200}{arXiv:1206.0200} (2012)
}%

\BiblioItem{1211.6965}
{
Aleks Kleyn,
Free Algebra with Countable Basis,\\
eprint \href{http://arxiv.org/abs/1211.6965}{arXiv:1211.6965} (2012)
}%

\DefBiblioItem{1302.7204}%

\BiblioItem{1305.4547}
{
Aleks Kleyn,
Normed $\Omega$-Group,\\
eprint \href{http://arxiv.org/abs/1305.4547}{arXiv:1305.4547} (2013)
}%

\BiblioItem{1310.5591}
{
Aleks Kleyn,
Integral of Map into Abelian $\Omega$\Hyph group,\\
eprint \href{http://arxiv.org/abs/1310.5591}{arXiv:1310.5591} (2013)
}%

\BiblioItem{1412.5425}
{
Aleks Kleyn,
Division in Associative $D$-Algebra,\\
eprint \href{http://arxiv.org/abs/1412.5425}{arXiv:1412.5425} (2014)
}%

\DefBiblioItem{1502.04063}%

\BiblioItem{1505.03625}
{
Aleks Kleyn,
Derivative of Map of Banach algebra,\\
eprint \href{http://arxiv.org/abs/1505.03625}{arXiv:1505.03625} (2015)
}%

\BiblioItem{1506.00061}
{
Aleks Kleyn,
Quadratic Equation over Associative $D$\Hyph Algebra,\\
eprint \href{http://arxiv.org/abs/1506.00061}{arXiv:1506.00061} (2015)
}%

\DefBiblioItem{1601.03259}%

\DefBiblioItem{1801.01628}

\DefBiblioItem{1908.04418}

\BiblioItem{2020.06.01}
{
Aleks Kleyn,
System of Differential Equations over Quaternion Algebra,\\
Geometry of differential equations seminar 2020\\
eprint \href{https://gdeq.org/files/Aleks_Kleyn-2020.06.01.English.pdf}{talk:2020.06.01} (2020)
}%

\BiblioItem{2021.01.06}
{
Aleks Kleyn,
Calculus over Quaternion Algebra,\\
Joint Mathematics Meetings January 2020\\
eprint \href{http://arxiv.org/abs/1908.04418}{arXiv:1908.04418} (2019)
}%

\BiblioItem{2112.00613}
{
Aleks Kleyn,
Polynomial in Non-Commutative Algebra,\\
eprint \href{http://arxiv.org/abs/2112.00613}{arXiv:2112.00613} (2021)
}%

\BiblioItem{322019412}
{
Aleks Kleyn,
Research Diary, 2017\\
eprint \href{https://www.researchgate.net/publication/322019412}{RG:322019412} (2017)
}%

\BiblioItem{323966352}
{
Aleks Kleyn,
Research Diary, 2018\\
eprint \href{https://www.researchgate.net/publication/323966352}{RG:323966352} (2018)
}%

\BiblioItem{348311191}
{
Aleks Kleyn,
Research Diary, 2021\\
eprint \href{https://www.researchgate.net/publication/348311191}{RG:348311191} (2021)
}%

\BiblioItem{323236904}
{
Aleks Kleyn,
Crash Course in Calculus over  Banach Algebra\\
eprint \href{https://www.researchgate.net/publication/323236904}{RG:323236904} (2018)
}%

\DefBiblioItem{2207.06506}

\BiblioItem{MAlgebra2}
{
Aleks Kleyn,
Introduction into Noncommutative Algebra,
Volume 2, Module over Algebra\\
eprint \href{http://arxiv.org/abs/MAlgebra2}{arXiv:MAlgebra2} (2019)
}%

\BiblioItem{8433-5163}
{
Aleks Kleyn,
Linear Maps of Free Algebra: First Steps in Noncommutative Linear Algebra,\\
Lambert Academic Publishing, 2010
}%

\BiblioItem{8443-0072}
{
Aleks Kleyn,
Representation Theory: Representation of Universal Algebra,\\
Lambert Academic Publishing, 2011
}%

\BiblioItem{4776-3181}
{
Aleks Kleyn.\\
Linear Algebra over Division Ring: System of Linear Equations.\\
CreateSpace Independent Publishing Platform, 2012;\\
ISBN-13: 978-1-4776-3181-2
}%

\BiblioItem{4975-6381}
{
Aleks Kleyn.\\
Single Variable Calculus: Noncomutative Banach Algebra.\\
CreateSpace Independent Publishing Platform, 2014;\\
ISBN-13: 978-1-4975-6381-0
}%

\BiblioItem{4993-2400}
{
Aleks Kleyn.\\
Linear Algebra over Division Ring: Vector Space.\\
CreateSpace Independent Publishing Platform, 2014;\\
ISBN-13: 978-1-4993-2400-6
}%

\BiblioItem{5059-9176}
{
Aleks Kleyn.\\
Normed \(\Omega\)-Group.\\
CreateSpace Independent Publishing Platform, 2015;\\
ISBN-13: 978-1-5059-9176-5
}%

\BiblioItem{5114-6019}
{
Aleks Kleyn.\\
Representation of Universal Algebra: Polymorphism.\\
CreateSpace Independent Publishing Platform, 2015;\\
ISBN-13: 978-1-5114-6019-4
}%

\BiblioItem{5148-4632}
{
Aleks Kleyn.\\
Noncommutative Algebra: Introduction.\\
CreateSpace Independent Publishing Platform, 2018;\\
ISBN-13: 978-1-5114-6019-4
}%

\BiblioItem{5410-9916}
{
Aleks Kleyn.\\
Lebesgue Integral in Abelian $\Omega$-Group.\\
CreateSpace Independent Publishing Platform, 2016;\\
ISBN-13: 978-1-5410-9916-6
}%

\BiblioItem{9856-6693}
{
Aleks Kleyn,
Crash Course in Calculus over  Banach Algebra\\
CreateSpace Independent Publishing Platform, 2018;\\
ISBN-13: 978-1-9856-6693-1
}%

\BiblioItem{7287-9339}
{
Aleks Kleyn,
Quadratic Equation over Associative $D$\Hyph Algebra\\
Kindle Direct Publishing, 2018;\\
ISBN-13: 978-1-7287-9339-9
}%

\BiblioItem{9835-2163}
{
Aleks Kleyn,
Differential Equation over Banach Algebra\\
Kindle Direct Publishing, 2018;\\
ISBN-13: 978-1-9835-2163-8
}%

\BiblioItem{6860-2955}
{
Aleks Kleyn,
Diagram of Representations of Universal Algebras,\\
Kindle Direct Publishing, 2019;\\
ISBN-13: 978-1-6860-2955-4
}%

\BiblioItem{0767-8264}
{
Aleks Kleyn,
System of Differential Equations over Banach Algebra: Exponent,\\
Kindle Direct Publishing, 2019;\\
ISBN-13: 978-1-0767-8264-9
}%

\BiblioItem{8428-0408}
{
Aleks Kleyn,
Introduction into Noncommutative Algebra, Volume 1: Division Algebra,\\
Kindle Direct Publishing, 2022;\\
ISBN-13: 979-8-8428-0408-5
}%

\BiblioItem{CACAA.01.109}
{
Aleks Kleyn,
Mappings of Conjugation of Quaternion Algebra.\\
Clifford Analysis, Clifford Algebras and their applications,
volume 1, Issue 1, pages 109 - 121, 2012
}%

\BiblioItem{CACAA.01.291}
{
Aleks Kleyn,
Introduction into Calculus over Division Ring.\\
Clifford Analysis, Clifford Algebras and their applications,
volume 1, Issue 4, pages 291 - 355, 2012
}%

\BiblioItem{CACAA.02.097}
{
Aleks Kleyn,
Polynomial over Associative $D$-Algebra.\\
Clifford Analysis, Clifford Algebras and their applications,
volume 2, Issue 2, pages 97 - 115, 2013
}%

\BiblioItem{CACAA.04.001}
{
Aleks Kleyn,
Integral of Map into Abelian $\Omega$-group.\\
Clifford Analysis, Clifford Algebras and their applications,
volume 4, Issue 1, pages 1 - 68, 2013
}%

\BiblioItem{CACAA.05.001}
{
Aleks Kleyn,
Introduction into Calculus over Division Ring.\\
Clifford Analysis, Clifford Algebras and their applications,
volume 5, issue 1, pages 1 - 68, 2016 
}%

\BiblioItem{CACAA.06.121}
{
Aleks Kleyn,
Differential Forms in Banach Algebra.\\
Clifford Analysis, Clifford Algebras and their applications,
volume 6, issue 2, pages 121 - 214, 2017 
}%

\BiblioItem{GJSFRA.13.1.39}
{
Aleks Kleyn,
Reference frame and Lorentz transformation,\\
Global Journals of Science Frontier Research A,
volume 13, issue 1, pages 39 - 55, 2013 
}%

\BiblioItem{1807.05583}
{
V. Sokolov, T. Wolf,
Non\Hyph commutative generalization of integrable quadratic ODE systems,\\
eprint \href{http://arxiv.org/abs/1807.05583}{arXiv:1807.05583} (2019)
}%

\BiblioItem{1506.05848}
{
Rida T. Farouki, Graziano Gentili, Carlotta Giannelli, Alessandra Sestini,
Caterina Stoppato,\\
Solution of a quadratic quaternion equation with mixed coefficients,\\
eprint \href{http://arxiv.org/abs/1506.05848}{arXiv:1506.05848} (2015)
}%

\BiblioItem{2003.05263}
{
Florian-Horia Vasilescu,
Spectrum and Analytic Functional Calculus in Real and Quaternionic Frameworks,\\
eprint \href{http://arxiv.org/abs/2003.05263}{arXiv:2003.05263} (2003)
}%

\BiblioItem{1702.04935}
{
M. Irene Falc\~{a}o, Fernando Miranda, Ricardo Severino, M. Joana Soares,
Weierstrass method for quaternionic polynomial root-finding,\\
eprint \href{http://arxiv.org/abs/1702.04935}{arXiv:1702.04935} (2017)
}%

\DefBiblioItem{1812.03397}%

\BiblioItem{Lauve: Quantum coordinates}
{
A. Lauve, Quantum- and quasi-Plucker coordinates,
\textit{Journal of Algebra} {\bf 296}(2), pp.440-461,
(2006).
}%

\BiblioItem{Lewis D. W. Quaternion algebras}
{
Lewis D. W. Quaternion algebras and the algebraic legacy
of Hamilton's quaternions, \textit{Irish Math. Soc. Bulletin} {\bf
57}, pp. 41-64, (2006).
}%

\BiblioItem{0812.2865}
{
Jos\'e Miguel Figueroa-O'Farrill,
Three lectures on 3-algebras,
eprint \href{http://arxiv.org/abs/0812.2865}{arXiv:0812.2865} (2008)
}%

\BiblioItem{1202.0951}
{
Daniel Edward Clark,
Deconvolution of point processes,
eprint \href{http://arxiv.org/abs/1202.0951}{arXiv:1202.0951} (2012)
}%

\BiblioItem{1202.4546}
{
Ming-Liang Hu,
Disentanglement, Bell-nonlocality violation
and teleportation capacity of the decaying tripartite states,
eprint \href{http://arxiv.org/abs/1202.4546}{arXiv:1202.4546} (2012)
}%

\BiblioItem{1203.1629}
{
Borivoje Dakic, Yannick Ole Lipp, Xiaosong Ma, Martin Ringbauer,
Sebastian Kropatschek, Stefanie Barz, Tomasz Paterek, Vlatko Vedral,
Anton Zeilinger, Caslav Brukner, Philip Walther,
Quantum Discord as Optimal Resource for Quantum Communication,
eprint \href{http://arxiv.org/abs/1203.1629}{arXiv:1203.1629} (2012)
}%

\BiblioItem{Li Nimmo: Darboux transformations}
{
C.X.Li, J.J.C. Nimmo, Darboux transformations for a twisted
derivation and quasideterminant solutions to the super KdV
equation, \textit{Proceedings of the Royal Society A:
Mathematical, Physical and Engineering Sciences} {\bf 466} (2120),
pp. 2471-2493, (2010).
}%

\BiblioItem{Schiebold: Cauchy-type determinants}
{
C. Schiebold, Cauchy-type determinants and integrable
systems, \textit{Linear Algebra and Its Applications} {\bf 433}
(2), pp. 447-475, (2010)
}%

\BiblioItem{Suzuki: Noncommutative spectral decomposition}
{
T. Suzuki, Noncommutative
spectral decomposition with qua\-si\-de\-ter\-mi\-nant, \textit{Advances in
Mathematics} {\bf 217}(5), pp. 2141-2158, (2008).
}%

\BiblioItem{1105.3456}
{
C. W. F. Everitt, D. B. DeBra, B. W. Parkinson, J. P. Turneaure, J. W. Conklin,
M. I. Heifetz, G. M. Keiser, A. S. Silbergleit, T. Holmes, J. Kolodziejczak,
M. Al-Meshari, J. C. Mester, B. Muhlfelder, V. Solomonik, K. Stahl, P. Worden,
W. Bencze, S. Buchman, B. Clarke, A. Al-Jadaan, H. Al-Jibreen, J. Li, J. A. Lipa,
J. M. Lockhart, B. Al-Suwaidan, M. Taber, S. Wang,\\
Gravity Probe B: Final Results of a Space Experiment to Test General Relativity,\\
eprint \href{http://arxiv.org/abs/1105.3456}{arXiv:1105.3456[gr-qc]} (2011)
}%

\BiblioItem{0009305}
{
G. S. Asanov.
Can Neutrinos and High-Energy Particles Test Finsler Metric of Space-Time?\\
eprint \href{http://arxiv.org/abs/hep-ph/0009305}{arXiv:hep-ph/0009305} (2000)
}%

\BiblioItem{Asanov 2004}
{
G. S. Asanov.
Finsleroid - space supplemented by angle and scalar product.\\
Hypercomplex Numbers in Geometry and Physics, {\bf 1}, 2004, p. 40 - 62
}%

\BiblioItem{1004.3007}
{
Sergiu I. Vacaru,
Principles of Einstein-Finsler Gravity and Perspectives in Modern Cosmology,\\
eprint \href{http://arxiv.org/abs/1004.3007}{arXiv:1004.3007[math-ph]} (2010)
}%

\BiblioItem{1012.4148}
{
Sergiu I. Vacaru.
Principles of Einstein-Finsler Gravity and Cosmology.\\
eprint \href{http://arxiv.org/abs/1012.4148}{arXiv:1012.4148[physics.gen-ph]} (2010)
}%

\BiblioItem{1112.5641}
{
Christian Pfeifer, Mattias N.R. Wohlfarth.
Finsler geometric extension of Einstein gravity.\\
eprint \href{http://arxiv.org/abs/1112.5641}{arXiv:1112.5641[gr-qc]} (2011)
}%

\BiblioItem{0711.0056}
{
Zhe Chang, Xin Li.
Lorentz Invariance Violation and Symmetry in Randers\Hyph Finsler Spaces.\\
eprint \href{http://arxiv.org/abs/0711.0056}{arXiv:0711.0056[hep-th]} (2011)
}%

\DefBiblioItem{1510.02224}%

\BiblioItem{1902.09800}
{
Dong Cheng, Kit Ian Kou, Yong Hui Xia.
Floquet Theory for Quaternion-valued Differential Equations.\\
eprint \href{http://arxiv.org/abs/1902.09800}{arXiv:1902.09800} (2019)
}%

\BiblioItem{Zharinov Kursy NOC, 9}
{
В. В. Жаринов,
Алгебро-геометрические основы математической физики.
\\
Лекц. курсы НОЦ, 9, МИАН, М., 2008, 3–209
}%

\BiblioItem{Rund Finsler geometry}
{
Hanno Rund,
The differential geometry of Finsler spaces.
\\
Springer - Verlag, Berlin - G\"ottingen - Heidelberg, 1959
}%

\BiblioItem{Smirnov vol 1}
{
V. I. Smirnov,
A Course of Higher Mathematics, volume I.
\\
Translated by D. E. Brown.
\\
Translation, edited and additions made by I. N. Sneddon.
\\
Pergamon Press, Addison-Wesley Publishing Company, 1964
}%

\BiblioItem{Beem Dostoglou Ehrlich}
{
John K. Beem, Stamatis A. Dostoglou, Paul E. Ehrlich,
Advances in differential geometry and general relativity.
\\
American Mathematical Society, 2004
}%

\BiblioItem{978-0719033414}
{
Malcolm Pemberton, Nicholas Rau,
Mathematics for economists: an introductory textbook.
\\
Manchester University Press, November 2001; ISBN-13: 978-0719033414
}%

\BiblioItem{0 521 59180 5}
{
Cyrus D. Cantrell,
Modern mathematical methods for physicists and engineers.
\\
Cambridge University Press, 2000
}%

\BiblioItem{Arveson spectral theory}
{
William Arveson,
A short course on spectral theory.
\\
Springer - Verlag, New York, 2002
}%

\BiblioItem{Robert Hermann}
{
Robert Hermann,
Topics in the mathematics of quantum mechanics.
\\
Math Sci Press, 1973
}%

\BiblioItem{9705.009}
{
John C. Baez,
An Introduction to n-Categories,\\
eprint \href{http://arxiv.org/abs/q-alg/9705009}{arXiv:q-alg/9705009} (1997)
}%

\BiblioItem{0105.155}
{
John C. Baez,
The Octonions,\\
eprint \href{http://arxiv.org/abs/math.RA/0105155}{arXiv:math.RA/0105155} (2002)
}%

\BiblioItem{John Baez: Math Blogs}
{
John C. Baez,
What do mathematicians need to know about blogging?,\\
Notices of the American Mathematical Society,
(2010), 3, {\bf 57}, 333,\\
\url{http://www.ams.org/notices/201003/rtx100300333p.pdf}
}%

\BiblioItem{Tolstoi about Anna Karenina}
{
Tolstoi about Anna Karenina,
in book A Karenina Companion, by C. J. G. Turner,
published by Wilfrid Laurier University Press (August 1993)
}%

\BiblioItem
{Cohn: Universal Algebra}
{
Paul M. Cohn,
Universal Algebra,
Springer, 1981
}%

\BiblioItem
{Cohn: Algebra 1}
{
Paul M. Cohn,
Algebra, Volume 1,
John Wiley \& Sons, 1982
}%

\BiblioItem
{Cohn: Algebra 3}
{
Paul M. Cohn,
Algebra, Volume 3,
John Wiley \& Sons, 1991
}%

\DefBiblioItem{Cohn: Skew Fields}%

\BiblioItem
{Lam: Noncommutative Rings}
{
T. Y. Lam,
A First Course in
Noncommutative Rings,
Springer-Verlag, 1991
}%

\BiblioItem
{Maunder: Algebraic Topology}
{
C. R. F. Maunder,
Algebraic Topology,
Dover Publications, Inc, Mineola, New York, 1996
}%

\BiblioItem{Pommaret: Partial Differential Equations}
{
J.-F. Pommaret,
Partial Differential Equations and Group Theory,
Springer, 1994
}%

\BiblioItem{Bourbaki: Set Theory}
{
N. Bourbaki,
Theory of sets,
Springer, 2004
}%

\BiblioItem{Bourbaki: Algebra 1}
{
N. Bourbaki,
Algebra 1,
Springer, 2004
}%

\BiblioItem{Bourbaki: Algebra 2}
{
N. Bourbaki,
Algebra II, Chapters 4 - 7,//
Translated by P. M. Cohn & J. Howie,//
Springer, 2004
}%

\BiblioItem
{Bourbaki: General Topology 1}
{
N. Bourbaki,
General Topology, Chapters 1 - 4,
Springer, 1989
}

\BiblioItem{Bourbaki: General Topology: Chapter 5 - 10}
{
N. Bourbaki,
General Topology, Chapters 5 - 10,
Springer, 1989
}

\BiblioItem{Bourbaki: Topological Vector Space}
{
N. Bourbaki,
Topological Vector Spaces, Chapters 1 - 5,
Transl. by H. G. Eggleston $\&$ S. Madan,
Springer, 2003
}

\BiblioItem{Bourbaki: Group Lie}
{
N. Bourbaki,
Lie Groups and Lie Algebras, Chapters 1 - 3,
Springer, 1989
}

\BiblioItem{Bourbaki: Coxeter Group Lie}
{
N. Bourbaki,
Lie Groups and Lie Algebras, Chapters 4 - 6,
Translator Andrew Pressley,
Springer, 2002
}

\BiblioItem{Bourbaki: Real Group Lie}
{
N. Bourbaki,
Lie Groups and Lie Algebras, Chapters 7 - 9,
Translator Andrew Pressley,
Springer, 2005
}

\BiblioItem{Shabat: Complex Analysis}
{
Shabat B. V.,
Introduction to Complex Analysis,
Moscow, Nauka, 1969
}

\BiblioItem{Pontryagin: Topological Group}
{
L. S. Pontryagin,
Selected Works, Volume Two, Topological Groups,
Gordon and Breach Science Publishers, 1986
}

\BiblioItem
{Eisenhart: Riemannian Geometry}
{
Eisenhart,
Riemannian Geometry,
Princeton University Press, Princeton, 1949
}

\BiblioItem
{Eisenhart: Continuous Groups of Transformations}
{
Eisenhart,
Continuous Groups of Transformations,
Dover Publications, New York, 1961
}

\BiblioItem
{Condon Odabasi}
{
Edward Uhler Condon, Halis Odabasi,
Atomic Structure,
CUP Archive, 1980
}

\BiblioItem{Postnikov: Differential Geometry}
{
Postnikov M. M.,
Geometry IV: Differential geometry,
Moscow, Nauka, 1983
}

\BiblioItem{Fikhtengolts: Calculus volume 1}
{
Fikhtengolts G. M.,
Differential and Integral Calculus Course, volume 1,
Moscow, Nauka, 1969
}

\BiblioItem{Fikhtengolts: Calculus volume 2}
{
Fikhtengolts G. M.,
Differential and Integral Calculus Course, volume 2,
Moscow, Nauka, 1969
}

\BiblioItem{Fikhtengolts: Calculus volume 3}
{
Fikhtengolts G. M.,
Differential and Integral Calculus Course, volume 3,
Moscow, Nauka, 1969
}

\BiblioItem{Hatcher: Algebraic Topology}
{
Allen Hatcher,
Algebraic Topology,
Cambridge University Press, 2002
}

\BiblioItem{geometry of differential equations}
{
Krasil'shchik I. S., Lychagin V. V., Vinogradov A. M.,
Geometry of Jet Spaces and Nonlinear Partial Differential Equations,
\\
Translated from the Russian by A. B. Sosinskii,
\\
Gordon and Breach Science Publishers, 1985
}

\BiblioItem{Basic Concepts of Differential Geometry}
{
Alekseyevskii D. V., Vinogradov A. M., Lychagin V. V.,
Basic Concepts of Differential Geometry
\\
VINITI Summary 28
\\
Moscow. VINITI, 1988
}

\BiblioItem{cohomological analysis}
{
A. M. Vinogradov,
Cohomological Analysis of Partial Differential Equations
and Secondary Calculus,
American Mathematical Society, 2001
}

\BiblioItem{0801.1734}
{
Brandon S. DiNunno, Richard A. Matzner,
The Volume Inside a Black Hole,\\
eprint \href{http://arxiv.org/abs/0801.1734v1}{arXiv:0801.1734v1} (2008)
}

\BiblioItem{0702.447}
{
Ivan Kyrchei,
Cramer's rule for some quaternion matrix equations,\\
eprint \href{http://arxiv.org/abs/math/0702447}{arXiv:math.RA/0702447} (2007)
}

\BiblioItem{Izrail M. Gelfand: Quaternion Groups}
{
I. M. Gelfand, M. I. Graev,
Representation of Quaternion Groups over Localy Compact and
Functional Fields,\\
Funct. Anal. Appl. {\bf 2} (1968) 19 - 33;\\
Izrail Moiseevich Gelfand, Semen Grigorevich Gindikin,\\
Izrail M. Gelfand: Collected Papers, volume II, 435 - 449,\\
Springer, 1989
}

\BiblioItem{Richard D. Schafer}
{
Richard D. Schafer,
An Introduction to Nonassociative Algebras,
Dover Publications, Inc., New York, 1995
}

\BiblioItem{Bamberg Sternberg}
{
Paul Bamberg, Shlomo Sternberg,
A course in mathematics for students of physics,
Cambridge University Press, 1991
}

\BiblioItem{Conway Smith}
{
John Horton Conway, Derek Alan Smith,
On quaternions and octonions: their geometry, arithmetic, and symmetry,
A K Peters, Natick, Massachussets, 2003
}

\BiblioItem{Fueter}
{
Fueter, R.
Die Funktionentheorie der Differentialgleichungen $\Delta u = 0$ und
$\Delta \Delta u = 0$ mit vier reellen Variablen.
Comment. Math. Helv. {\bf 7} (1935), 307-330
}

\DefBiblioItem{Sudbery Quaternionic Analysis}%

\BiblioItem{Sudbery 2657821}
{
A. Sudbery,
Quaternionic Analysis,\\
eprint \href{https://www.researchgate.net/publication/2657821}{ResearchGate:2657821} (1977)
}

\BiblioItem{0902.4771}
{
Fabrizio Colombo, Graziano Gentili, Irene Sabadini,
A Cauchy kernel for slice regular functions,\\
eprint \href{http://arxiv.org/abs/0902.4771v1}{arXiv:0902.4771v1} (2009)
}

\BiblioItem{Vadim Komkov}
{
Vadim Komkov,
Variational Principles of Continuum Mechanics with Engineering Applications: Critical Points Theory,\\
Springer, 1986
}

\BiblioItem{Alain Connes 1994}
{
Alain Connes,
Noncommutative Geometry,\\
Academic Press, 1994
}

\BiblioItem{Hamilton papers 3}
{
Sir William Rowan Hamilton,
The Mathematical Papers, Vol. III, Algebra,\\
Cambridge at the University Press, 1967
}

\BiblioItem{Hamilton Elements of Quaternions 1}
{
Sir William Rowan Hamilton,
Elements of Quaternions, Volume I,\\
Longmans, Green, and Co., London, New York, and Bombay, 1899
}

\BiblioItem{Cartan geometry in reper}
{
Elie Cartan, Vladislav V. Goldberg, Serge\u{i} Pavlovich Finikov,\\
Riemannian geometry in an orthogonal frame:
from lectures delivered by Elie Cartan at the Sorbonne in 1926-1927,\\
translated by Vladislav V. Goldberg,\\
World Scientific, 2001
}

\BiblioItem{Cartan differential form}
{
Henri Cartan.
Differential forms.\\
Kershaw Publishing Company Limited, London, 1971
}

\BiblioItem{Arnautov Glavatsky Mikhalev}
{
V. I. Arnautov, S. T. Glavatsky, A. V. Mikhalev,\\
Introduction to the theory of topological rings and modules,
Volume 1995,\\
Marcel Dekker, Inc, 1996
}

\BiblioItem{Moore Yaqub}
{
Hal G. Moore, Adil Yaqub,
A first course in linear algebra with applications,
Edition 3, Academic Press, 1998 
}

\BiblioItem{math.CV-0405471}
{
S. V. Ludkovsky,
Differentiable functions of Cayley-Dickson numbers,\\
eprint \href{http://arxiv.org/abs/math.CV/0405471}{arXiv:math.CV/0405471} (2004)
}%

\BiblioItem{W.Bertram H.Glockner K.Neeb}
{
W.Bertram, H.Glockner, K.Neeb,
Differential Calculus over General Base Fields and Rings,
Expositiones Mathematicae (2004), Volume 22, Issue 3, Pages 213-282
}

\CloseBiblio

%% file: Index.English.tex
\OpenIndex
\SetIndexSpace%
\Index
   {$1$\Hyph form}%
   {1-form}%
\SetIndexSpace%
\Index
   {$2$\Hyph ary fibered relation}%
   {2 ary fibered relation}%
\SetIndexSpace%
\Index
   {$A$\Hyph algebra of polynomials over $D$\Hyph algebra $A$}%
   {algebra of polynomials over algebra}%
\Index
   {$A$\Hyph number}%
   {A number}%
\Index
   {$\mathcal A(A)$\Hyph map}%
   {A(A) map}%
\Index
   {$A*$\Hyph module}%
   {A*-module}%
\Index
   {$A*$\Hyph vector space}%
   {A*-vector space}%
\Index
   {$A$\Hyph module}%
   {module over algebra}%
\Index
   {$A$\Hyph valued function}%
   {A valued function}%
\Index
   {$A$\Hyph representation in $\Omega$\Hyph algebra}%
   {A representation of algebra}%
\Index
   {$A$\hyph vector space}%
   {A vector space}%
\Index
   {Abelian multiplicative $\Omega$\Hyph group}%
   {Abelian multiplicative Omega group}%
\Index
   {Abelian $\Omega$\Hyph group}%
   {Abelian Omega group}%
\Index
   {Abelian semigroup}%
   {Abelian semigroup}%
\Index
   {absolute value}%
   {absolute value}%
\Index
   {active \sT{G}representation}%
   {active representation, vector space}%
\Index
   {active representation}%
   {active representation}%
\Index
   {active representation in basis manifold}%
   {active representation in basis manifold}%
\Index
   {active representation of group $G(\Vector f)$ in basis manifold of tower of representations}%
   {active representation in basis manifold, tower of representations}%
\Index
   {active transformation of basis manifold}%
   {active transformation of basis}%
\Index
   {active transformation on basis manifold}%
   {active transformation}%
\Index
   {active transformation on the set of \rcd bases}%
   {active transformation, vector space}%
\Index
   {additive map}%
   {additive map}%
\Index
   {affine basis}%
   {Affine Basis}%
\Index
   {affine functional}%
   {affine functional}%
\Index
   {affine representation of Lie group}%
   {affine representation of Lie group}%
\Index
   {affine space}%
   {affine space}%
\Index
   {affine structure on set}%
   {affine structure on set}%
\Index
   {affine transformation}%
   {affine transformation}%
\Index
   {affine transformation group}%
   {affine transformation group}%
\Index
   {affine transformation group}%
   {affine transformation group}%
\Index
   {affine transformation on basis manifold}%
   {affine transformation}%
\Index
   {algebra of fractions of algebra with conjugation}%
   {algebra of fractions of algebra with conjugation}%
\Index
   {algebra of polynomials over $D$\Hyph algebra}%
   {algebra of polynomials over D algebra}%
\Index
   {algebra of rational mappings of algebra}%
   {algebra of rational mappings of algebra}%
\Index
   {algebra of sets}%
   {algebra of sets}%
\Index
   {algebra over ring}%
   {algebra over ring}%
\Index
   {algebra with conjugation}%
   {algebra with conjugation}%
\Index
   {alternation of polylinear map}%
   {alternation of polylinear map}%
\Index
   {alternative representation of matrix}%
   {Alternative representation}%
\Index
   {anholonomic coordinate}%
   {anholonomic coordinate}%
\Index
   {anholonomic coordinates of connection}%
   {anholonomic coordinates of connection}%
\Index
   {anholonomic coordinates of vector}%
   {vector anholonomic coordinates}%
\Index
   {anholonomic coordinates on manifold}%
   {anholonomic coordinates on manifold}%
\Index
   {anholonomity object}%
   {anholonomity object}%
\Index
   {antilinear homomorphism}%
   {antilinear homomorphism}%
\Index
   {antilinear map}%
   {antilinear map}%
\Index
   {antisymmetric $2$\Hyph ary fibered relation}%
   {antisymmetric 2 ary fibered relation}%
\Index
   {$A\RCstar$\Hyph basis for vector space}%
   {Arc basis, vector space}%
\Index
   {arity}%
   {arity}%
\Index
   {arity of operation}%
   {arity of operation}%
\Index
   {associative $D$\Hyph algebra}%
   {associative D algebra}%
\Index
   {associative law}%
   {associative law}%
\Index
   {associative multiplicative $\Omega$\Hyph group}%
   {associative multiplicative Omega group}%
\Index
   {associative $\Omega$\Hyph group}%
   {associative Omega group}%
\Index
   {associative operation}%
   {associative operation}%
\Index
   {associator of $D$\Hyph algebra}%
   {associator of algebra}%
\Index
   {augmented matrix}%
   {augmented matrix}%
\Index
   {auto parallel line}%
   {auto parallel line}%
\Index
   {automorphism}%
   {automorphism}%
\Index
   {automorphism of diagram of representations}%
   {automorphism of diagram of representations}%
\Index
   {automorphism of representation of $\Omega$\Hyph algebra}%
   {automorphism of representation}%
\Index
   {automorphism of tower of representations}%
   {automorphism of tower of representations}%
\Index
   {automorphism of vector space}%
   {automorphism of vector space}%
\Index
   {$(^j_i)$\hyph \CR quasideterminant}%
   {j i cr-quasideterminant}%
\Index
   {norm of quaternion}%
   {norm of quaternion}%
\SetIndexSpace%
\Index
   {$B$\Hyph set}%
   {B set}%
\Index
   {Banach $D$\Hyph algebra}%
   {Banach algebra}%
\Index
   {Banach $D$\Hyph module}%
   {Banach module}%
\Index
   {base of fibered correspondence}%
   {base of fibered correspondence}%
\Index
   {base of mapping}%
   {base of map}%
\Index
   {basis}%
   {Basis}%
\Index
   {basis dual to basis}%
   {basis dual to basis}%
\Index
   {basis dual to basis}%
   {dual basis}%
\Index
   {basis for \crd vector space}%
   {basis, crd vector space}%
\Index
   {basis for \dcr vector space}%
   {basis, dcr vector space}%
\Index
   {basis for \drc vector space}%
   {basis, drc vector space}%
\Index
   {basis for module}%
   {basis, module}%
\Index
   {basis for \rcd vector space}%
   {basis, rcd vector space}%
\Index
   {basis for vector space}%
   {basis, vector space}%
\Index
   {basis manifold}%
   {basis manifold}%
\Index
   {basis manifold of affine space}%
   {Basis Manifold, Affine Space}%
\Index
   {basis manifold of central affine space}%
   {Basis Manifold, Central Affine Space}%
\Index
   {basis manifold of Euclid space}%
   {Basis Manifold, Euclid Space}%
\Index
   {basis manifold of Euclid space}%
   {Basis Manifold, Euclid Space, division ring}%
\Index
   {basis manifold of tower of representations}%
   {basis manifold tower of representations}%
\Index
   {basis of algebra $\mathcal L(A;A)$}%
   {basis of algebra L(A,A)}%
\Index
   {basis of diagram of representations}%
   {basis of diagram of representations}%
\Index
   {basis of representation}%
   {basis of representation}%
\Index
   {basis of tower of representations}%
   {basis of tower of representations}%
\Index
   {basis vector of representation of Lie group over algebra $A$}%
   {basis vector of representation of Lie group over algebra A}%
\Index
   {biring}%
   {biring}%
\Index
   {Borel algebra}%
   {Borel algebra}%
\Index
   {Borel set}%
   {Borel set}%
\Index
   {Borel\Hyph measurable map}%
   {Borel-measurable map}%
\Index
   {bundle of level $2$}%
   {bundle of level 2}%
\Index
   {bundle of level $n$}%
   {bundle of level n}%
\SetIndexSpace%
\Index
   {\subs row of matrix}%
   {c row}%
\Index
   {$c$\hyph row of matrix}%
   {c-row}%
\Index
   {can be embeded}%
   {can be embeded}%
\Index
   {cancellation law}%
   {cancellation law}%
\Index
   {canonical map}%
   {canonical map}%
\Index
   {canonical map}%
   {canonical map}%
\Index
   {canonical remainder of the division}%
   {canonical remainder of the division}%
\Index
   {canonical representation of division with remainder}%
   {canonical representation of division with remainder}%
\Index
   {carrier of $\Omega$\Hyph algebra}%
   {carrier of Omega-algebra}%
\Index
   {Cartan connection}%
   {Cartan connection}%
\Index
   {Cartan curvature}%
   {Cartan curvature}%
\Index
   {Cartan derivative}%
   {Cartan derivative}%
\Index
   {Cartan equation}%
   {Cartan equation}%
\Index
   {Cartan symbol}%
   {Cartan symbol}%
\Index
   {Cartan transport}%
   {Cartan transport}%
\Index
   {Cartesian power}%
   {Cartesian power}%
\Index
   {Cartesian power $\Bundle A$ of bundle $\Bundle B$}%
   {Cartesian power A of bundle B}%
\Index
   {Cartesian power $A$ of set $B$}%
   {Cartesian power of set}%
\Index
   {Cartesian power $n$ of bundle $\Bundle E$}%
   {Cartesian power n of bundle E}%
\Index
   {Cartesian power $n$ of $\mathfrak{H}$\Hyph algebra}%
   {Cartesian power of algebra}%
\Index
   {Cartesian power of systems of subsets}%
   {Cartesian power of systems of subsets}%
\Index
   {Cartesian product of groups}%
   {Cartesian product of groups}%
\Index
   {Cartesian product of measures}%
   {Cartesian product of measures}%
\Index
   {Cartesian product of \(\Omega\)\Hyph algebras}%
   {Cartesian product of Omega algebras}%
\Index
   {Cartesian product of systems of subsets}%
   {Cartesian product of systems of subsets}%
\Index
   {category of \drc vector spaces}%
   {category of drc vector spaces}%
\Index
   {category of fibered correspondences over diagonal}%
   {category of fibered correspondences over diagonal}%
\Index
   {category of left-side representations}%
   {category of left-side representations}%
\Index
   {category of left-side representations of $\Omega_1$\Hyph algebra $A$}%
   {category of left-side representations of Omega1 algebra}%
\Index
   {category of reduced fibered correspondences}%
   {category of reduced fibered correspondences}%
\Index
   {category of representations}%
   {category of representations}%
\Index
   {Cauchy sequence}%
   {Cauchy sequence}%
\Index
   {center of $A$\Hyph number}%
   {center of A number}%
\Index
   {center of $D$\Hyph algebra $A$}%
   {center of algebra}%
\Index
   {center of ring $D$}%
   {center of ring}%
\Index
   {central affine basis}%
   {Central Affine Basis}%
\Index
   {closed ball}%
   {closed ball}%
\Index
   {closure of set}%
   {closure of set}%
\Index
   {coefficient of polynomial}%
   {coefficient of polynomial}%
\Index
   {colinear vectors}%
   {colinear vectors}%
\Index
   {column $D*$\Hyph vector}%
   {column D* vector}%
\Index
   {column determinant}%
   {column determinant}%
\Index
   {column vector}%
   {column vector}%
\Index
   {common factor}%
   {common factor}%
\Index
   {commutative $D$\Hyph algebra}%
   {commutative D algebra}%
\Index
   {commutative diagram of correspondences}%
   {commutative diagram of correspondences}%
\Index
   {commutative diagram of representations of universal algebras}%
   {commutative diagram of representations}%
\Index
   {commutative law}%
   {commutative law}%
\Index
   {commutative operation}%
   {commutative operation}%
\Index
   {commutativity of representations}%
   {commutativity of representations}%
\Index
   {commutator of $D$\Hyph algebra}%
   {commutator of algebra}%
\Index
   {compact set}%
   {compact set}%
\Index
   {compact\hyph open topology}%
   {compact open topology}%
\Index
   {complete division ring}%
   {complete division ring}%
\Index
   {complete measure}%
   {complete measure}%
\Index
   {complete normed $\Omega$\Hyph group}%
   {complete Omega group}%
\Index
   {complete ring}%
   {complete ring}%
\Index
   {complete system of linear partial differential equations}%
   {Complete System of Linear Partial Differential Equations}%
\Index
   {completely integrable system}%
   {completely integrable system}%
\Index
   {completion of normed $\Omega$\Hyph group}%
   {completion of normed Omega group}%
\Index
   {completion of representation}%
   {completion of representation}%
\Index
   {component of derivative}%
   {component of derivative}%
\Index
   {component of derivative of Second Order}%
   {component of derivative of Second Order}%
\Index
   {component of linear map}%
   {component of linear map}%
\Index
   {component of polylinear map}%
   {component of polylinear map}%
\Index
   {component of the G\^ateaux derivative}%
   {component of Gateaux derivative}%
\Index
   {component of the G\^ateaux derivative of second order}%
   {component of Gateaux derivative of Second Order}%
\Index
   {composition of fibered correspondences}%
   {composition of fibered correspondences}%
\Index
   {composition of reduced fibered correspondences}%
   {composition of reduced fibered correspondences}%
\Index
   {condition of reducibility of products}%
   {condition of reducibility of products}%
\Index
   {congruence}%
   {congruence}%
\Index
   {conjugate of quaternion $x$}%
   {conjugate of quaternion}%
\Index
   {conjugated affine space}%
   {conjugated affine space}%
\Index
   {conjugated $D$\Hyph  module}%
   {conjugated D module}%
\Index
   {conjugated vector space}%
   {conjugated vector space}%
\Index
   {conjugation in algebra}%
   {conjugation in algebra}%
\Index
   {conjugation in ring}%
   {conjugation in ring}%
\Index
   {conjugation transformation}%
   {conjugation transformation}%
\Index
   {connected set}%
   {connected set}%
\Index
   {connection coefficients in affine space}%
   {connection coefficients, affine space}%
\Index
   {connection in principal fibre bundle}%
   {connection in principal bundle}%
\Index
   {contact point of set}%
   {contact point of set}%
\Index
   {continues basis}%
   {continues basis}%
\Index
   {continuous correspondence}%
   {continuous correspondence}%
\Index
   {continuous map}%
   {continuous map}%
\Index
   {continuous multivariable map}%
   {continuous multivariable map}%
\Index
   {continuous Schauder basis}%
   {continuous Schauder basis}%
\Index
   {contravariant representation}%
   {contravariant representation}%
\Index
   {convex set}%
   {convex set}%
\Index
   {coordinate isomorphism}%
   {coordinate isomorphism}%
\Index
   {coordinate matrix of set of vectors}%
   {coordinate matrix of set of vectors}%
\Index
   {coordinate matrix of vector}%
   {coordinate matrix of vector}%
\Index
   {coordinate matrix of vector field in \rcD basis}%
   {coordinate matrix of vector field in drc basis}%
\Index
   {coordinate \rcd vector space}%
   {coordinate rcd vector space}%
\Index
   {coordinate reference frame}%
   {coordinate reference frame}%
\Index
   {coordinate representation}%
   {coordinate representation}%
\Index
   {coordinate representation in tuple of $\VX\Omega$\Hyph algebras}%
   {coordinate tower of representations, Omega algebra}%
\Index
   {coordinate representation of group in vector space}%
   {coordinate representation, vector space}%
\Index
   {coordinate representation of vector}%
   {coordinate representation of vector}%
\Index
   {coordinate vector bundle}%
   {coordinate vector bundle}%
\Index
   {coordinate vector space}%
   {coordinate vector space}%
\Index
   {coordinates}%
   {coordinates}%
\Index
   {coordinates of $A_2$\Hyph number $m$ relative to set $X$}%
   {coordinates relative to set}%
\Index
   {coordinates of associator}%
   {coordinates of associator}%
\Index
   {coordinates of basis}%
   {coordinates of basis}%
\Index
   {coordinates of basis of representation}%
   {coordinates of basis relative to basis, representation}%
\Index
   {coordinates of endomorphism of representation}%
   {coordinates of endomorphism, representation}%
\Index
   {coordinates of endomorphism of tower of representations}%
   {coordinates of endomorphism, tower of representations}%
\Index
   {coordinates of geometric object}%
   {coordinates of geometric object}%
\Index
   {coordinates of homomorphism}%
   {coordinates of homomorphism}%
\Index
   {coordinates of morphism of diagram of representations}%
   {coordinates of morphism, diagram of representations}%
\Index
   {coordinates of point $A$ of affine space $\overset{\circ}{A}$ relative to basis $(O,\Basis e)$}%
   {coordinates in affine space}%
\Index
   {coordinates of reduced morphism of representation}%
   {coordinates of reduced morphism of representation}%
\Index
   {coordinates of representation}%
   {coordinates of representation, drc vector space}%
\Index
   {coordinates of representation}%
   {coordinates of representation}%
\Index
   {coordinates of set of vectors}%
   {coordinates of set of vectors}%
\Index
   {coordinates of vector}%
   {coordinates of vector}%
\Index
   {coordinates of vector field in \Drc basis}%
   {coordinates of vector field in drc basis}%
\Index
   {coordinates of vector relative to Hamel basis}%
   {coordinates of vector, Hamel basis}%
\Index
   {coordinates of vector relative to Schauder basis}%
   {coordinates of vector, Schauder basis}%
\Index
   {coproduct of objects in category}%
   {coproduct in category}%
\Index
   {correspondence continuous on the set}%
   {correspondence continuous on the set}%
\Index
   {correspondence of homomorphism}%
   {correspondence of homomorphism}%
\Index
   {cosine}%
   {cosine}%
\Index
   {covariant representation}%
   {covariant representation}%
\Index
   {\CR exponent}%
   {CR exponent}%
\Index
   {\CR inverse element of biring}%
   {cr-inverse element}%
\Index
   {\CR matrix group}%
   {cr-matrix group}%
\Index
   {\CR nonsingular matrix}%
   {cr nonsingular matrix}%
\Index
   {\CR power}%
   {cr power}%
\Index
   {\CR product (product column over row)}%
   {cr-product}%
\Index
   {$\CRcirc$\Hyph product of matrices of maps}%
   {cr product of matrices of maps}%
\Index
   {\CR singular matrix}%
   {cr singular matrix}%
\Index
   {\CR inverse matrix}%
   {cr-inverse matrix}%
\Index
   {\CR quasideterminant}%
   {cr-quasideterminant}%
\Index
   {\crd vector}%
   {crd vector}%
\Index
   {\crd vector space}%
   {crd vector space}%
\Index
   {$C^*$\Hyph algebra}%
   {Cstar-algebra}%
\Index
   {curvilinear coordinates of point in affine space}%
   {curvilinear coordinates of point in affine space}%
\SetIndexSpace%
\Index
   {$D$\Hyph linear functional}%
   {D linear functional}%
\Index
   {$D*$\hyph matrices vector space}%
   {matrices vector space}%
\Index
   {$D*$\hyph  vector space}%
   {D* vector space}%
\Index
   {$D*$\Hyph module}%
   {D*-module}%
\Index
   {$D$\Hyph affine connection on manifold with affine connections}%
   {D affine connection, affine manifold}%
\Index
   {$D$\Hyph algebra}%
   {D algebra}%
\Index
   {$D$\Hyph module}%
   {D-module}%
\Index
   {$D$\Hyph module}%
   {D module}%
\Index
   {$D$\Hyph valued variable}%
   {D valued variable}%
\Index
   {$D$\Hyph vector function}%
   {d vector function}%
\Index
   {$D$\Hyph affine connection coefficients on manifold}%
   {D affine connection coefficients, manifold}%
\Index
   {\dcr vector}%
   {dcr vector}%
\Index
   {\dcr vector space}%
   {dcr vector space}%
\Index
   {definite integral}%
   {definite integral}%
\Index
   {derivative of map}%
   {derivative of map}%
\Index
   {derivative of order $n$}%
   {derivative of Order n}%
\Index
   {derivative of second order}%
   {derivative of Second Order}%
\Index
   {determinant of matrix}%
   {determinant}%
\Index
   {deviation of trajectories}%
   {deviation of trajectories}%
\Index
   {diagonal in bundle}%
   {diagonal in bundle}%
\Index
   {diagram of correspondences}%
   {diagram of correspondences}%
\Index
   {diagram of representations}%
   {diagram of representations}%
\Index
   {diagram of representations of universal algebras}%
   {diagram of representations of algebras}%
\Index
   {differentiable map}%
   {differentiable map}%
\Index
   {differential equation with separated variables}%
   {differential equation with separated variables}%
\Index
   {differential form of degree $p$}%
   {differential form of degree p}%
\Index
   {differential of independent variable}%
   {differential of independent variable}%
\Index
   {differential of map}%
   {differential of map}%
\Index
   {differential $p$\Hyph form}%
   {differential p form}%
\Index
   {differential separable equation}%
   {differential separable equation}%
\Index
   {dimension of \rcd vector space}%
   {dimension of vector space}%
\Index
   {direct product of bundles}%
   {Cartesian product of bundles}%
\Index
   {direct product of $D$\Hyph vector spaces}%
   {direct product of D vector spaces}%
\Index
   {direct product of division rings}%
   {direct product of division rings}%
\Index
   {direct product of \Ts{G}representations}%
   {direct product of G* representations}%
\Index
   {direct product of \(\Omega\)\Hyph algebras}%
   {direct product of Omega algebras}%
\Index
   {direct product of \rcd vector spaces}%
   {direct product, rcd vector space}%
\Index
   {direct product of representations of fibered group}%
   {direct product of representations of fibered group}%
\Index
   {direct product of representations of group}%
   {direct product of representations of group}%
\Index
   {direct product of total spaces}%
   {Cartesian product of total spaces}%
\Index
   {direct sum}%
   {direct sum}%
\Index
   {direct sum of representations}%
   {direct sum of representations}%
\Index
   {direction over commutative ring}%
   {direction over commutative ring}%
\Index
   {distributive law}%
   {distributive law}%
\Index
   {division algebra}%
   {division algebra}%
\Index
   {division with remainder}%
   {division with remainder}%
\Index
   {division without remainder}%
   {division without remainder}%
\Index
   {divisor of polynomial}%
   {divisor of polynomial}%
\Index
   {double determinant}%
   {double determinant}%
\Index
   {\Drc linear map of vector bundles}%
   {drc linear map of vector bundles}%
\Index
   {\drc vector}%
   {drc vector}%
\Index
   {\drc vector space}%
   {drc vector space}%
\Index
   {$D\star$\Hyph antilinear homomorphism}%
   {Dstar antilinear homomorphism}%
\Index
   {$\mathcal D\star$\Hyph vector bundle}%
   {Dstar vector bundle}%
\Index
   {$\mathcal D\star$\Hyph vector field}%
   {Dstar vector field}%
\Index
   {$\mathcal D\star$\hyph linear composition of vector fields}%
   {linear composition of vector fields}%
\Index
   {$\mathcal D\star$\hyph product of vector field over scalar}%
   {Dstar product of vector field over scalar, vector space}%
\Index
   {dual space of \rcd vector space}%
   {dual space of rcd vector space}%
\Index
   {duality principle for biring}%
   {duality principle for biring}%
\Index
   {duality principle for biring of matrices}%
   {duality principle for biring of matrices}%
\SetIndexSpace%
\Index
   {effective \Ts{G}representation}%
   {effective G* representation}%
\Index
   {effective representation}%
   {effective representation}%
\Index
   {effective representation of division ring}%
   {effective representation of division ring}%
\Index
   {effective representation of fibered $\Omega$\Hyph algebra}%
   {effective representation of fibered Omega-algebra}%
\Index
   {effective representation of group}%
   {effective representation of group}%
\Index
   {effective representation of ring}%
   {effective representation of ring}%
\Index
   {effective \Ts representation of fibered division ring}%
   {effective representation of fibered division ring}%
\Index
   {effective \Ts representation of fibered group}%
   {effective representation of fibered group}%
\Index
   {eigencolumn}%
   {eigencolumn}%
\Index
   {eigenrow}%
   {eigenrow}%
\Index
   {eigenvalue}%
   {eigenvalue}%
\Index
   {eigenvector}%
   {eigenvector}%
\Index
   {Einstein equation}%
   {Einstein equation}%
\Index
   {endomorphism}%
   {endomorphism}%
\Index
   {endomorphism of diagram of representations}%
   {endomorphism of diagram of representations}%
\Index
   {endomorphism of representation of $\Omega$\Hyph algebra}%
   {endomorphism of representation}%
\Index
   {endomorphism of representation regular on generating set $X$}%
   {endomorphism of representation, regular on set}%
\Index
   {endomorphism of representation singular on generating set $X$}%
   {endomorphism of representation, singular on set}%
\Index
   {endomorphism of tower of representations}%
   {endomorphism of tower of representations}%
\Index
   {endomorphism of tower of representations regular on tuple of generating sets}%
   {endomorphism of representation, regular on tuple}%
\Index
   {endomorphism of tower of representations singular on tuple of generating sets}%
   {endomorphism of representation, singular on tuple}%
\Index
   {enhanced Lie group}%
   {enhanced Lie group}%
\Index
   {epimorphism}%
   {epimorphism}%
\Index
   {equivalence}%
   {equivalence}%
\Index
   {equivalence generated by representation $f$}%
   {equivalence of representation}%
\Index
   {equivalent norms}%
   {equivalent norms}%
\Index
   {essential parameters in a set of functions}%
   {essential parameters}%
\Index
   {Euclidean metric on division ring}%
   {Euclidean metric on division ring}%
\Index
   {Euclidean scalar product in $D$\Hyph vector space}%
   {Euclidean scalar product, vector space}%
\Index
   {Euclidean scalar product on division ring}%
   {Euclidean scalar product on division ring}%
\Index
   {everywhere dense subset}%
   {everywhere dense subset}%
\Index
   {exact differential equation}%
   {exact differential equation}%
\Index
   {exact sequence}%
   {exact sequence}%
\Index
   {expansion of vector relative to basis converges}%
   {expansion converges}%
\Index
   {expansion of vector relative to basis converges normally}%
   {expansion converges normally}%
\Index
   {exponent}%
   {exponent}%
\Index
   {extension of correspondence}%
   {extension of correspondence}%
\Index
   {extension of measure}%
   {extension of measure}%
\Index
   {exterior differential}%
   {exterior differential}%
\Index
   {exterior product}%
   {exterior product}%
\Index
   {extreme line}%
   {extreme line}%
\SetIndexSpace%
\Index
   {factor group}%
   {factor group}%
\Index
   {fibered coordinate isomorphism}%
   {fibered coordinate isomorphism}%
\Index
   {fibered correspondence from $\Bundle A$ to $\Bundle B$}%
   {fibered correspondence from A to B}%
\Index
   {fibered correspondence in $\Bundle{A}$}%
   {fibered correspondence in A}%
\Index
   {fibered correspondence of homomorphism}%
   {fibered correspondence of homomorphism}%
\Index
   {fibered equivalence}%
   {fibered equivalence}%
\Index
   {fibered group}%
   {fibered group}%
\Index
   {fibered identification morphism}%
   {fibered identification morphism}%
\Index
   {fibered little group}%
   {fibered little group}%
\Index
   {fibered morphism from bundle $\Bundle A$ into $\Bundle B$}%
   {fibered morphism from A into B}%
\Index
   {fibered natural morphism}%
   {fibered natural morphism}%
\Index
   {fibered $\Omega$\Hyph algebra}%
   {fibered Omega-algebra}%
\Index
   {fibered $\Omega$\Hyph subalgebra}%
   {fibered Omega-subalgebra}%
\Index
   {fibered ordering}%
   {fibered ordering}%
\Index
   {fibered preordering}%
   {fibered preordering}%
\Index
   {fibered ring}%
   {fibered ring}%
\Index
   {fibered stability group}%
   {fibered stability group}%
\Index
   {fibered subset}%
   {fibered subset}%
\Index
   {field equation}%
   {field equation}%
\Index
   {field-strength tensor}%
   {field-strength tensor}%
\Index
   {filter $\mathfrak{F}$ converges to $A$}%
   {filter converges}%
\Index
   {finite expansion of set}%
   {finite expansion of set}%
\Index
   {Finsler metric}%
   {Finsler metric}%
\Index
   {Finsler space}%
   {Finsler space}%
\Index
   {Finsler structure}%
   {Finsler structure}%
\Index
   {first integral}%
   {first integral}%
\Index
   {first Newton law}%
   {First Newton law}%
\Index
   {frame\Hyph dragging effect}%
   {frame dragging effect}%
\Index
   {free $A$\Hyph module}%
   {free A module}%
\Index
   {free Abelian group}%
   {free Abelian group}%
\Index
   {free algebra over ring}%
   {free algebra over ring}%
\Index
   {free module}%
   {free module}%
\Index
   {free representation}%
   {free representation}%
\Index
   {free representation of group}%
   {free representation of group}%
\Index
   {free \Ts representation of fibered group}%
   {free representation of fibered group}%
\Index
   {Frenet transport}%
   {Frenet transport}%
\Index
   {function homogeneous of degree $k$}%
   {function homogeneous}%
\Index
   {function of division ring \Ds differentiable in the Fr\'echet sense}%
   {function Dstar differentiable in Frechet sense, division ring}%
\Index
   {fundamental sequence}%
   {fundamental sequence}%
\SetIndexSpace%
\Index
   {$G$\Hyph reference frame}%
   {G reference frame}%
\Index
   {$G$\Hyph basis of vector space}%
   {G-basis}%
\Index
   {$G$\Hyph coordinates of basis}%
   {G-coordinates}%
\Index
   {$G$\Hyph space}%
   {GSpace}%
\Index
   {the G\^ateaux \dcr derivative of map $f$ of $D$\Hyph vector space $V$ to $D$\Hyph vector space $W$}%
   {Gateaux dcr derivative of map, D vector space}%
\Index
   {the G\^ateaux derivative of map}%
   {Gateaux derivative of map}%
\Index
   {the G\^ateaux derivative of order $n$}%
   {Gateaux derivative of Order n}%
\Index
   {the G\^ateaux derivative of second order}%
   {Gateaux derivative of Second Order}%
\Index
   {the G\^ateaux \Ds derivative of map $f$ of division ring $D$}%
   {Gateaux Dstar derivative of map, division ring}%
\Index
   {the G\^ateaux mixed partial derivative}%
   {Gateaux partial derivative of Second Order}%
\Index
   {the G\^ateaux partial \dcr derivative of map $f^{\gi b}$ with respect to variable $x^{\gi a}$}%
   {Gateaux partial dcr derivative of map with respect to variable, D vector space}%
\Index
   {the G\^ateaux partial derivative}%
   {Gateaux partial derivative}%
\Index
   {the G\^ateaux partial \rcd derivative of map $f^{\gi b}$ with respect to variable $x^{\gi a}$}%
   {Gateaux partial rcd derivative of map with respect to variable, D vector space}%
\Index
   {the G\^ateaux \rcd derivative of map $f$ of $D$\hyph vector space $V$ to $D$\hyph vector space $W$}%
   {Gateaux rcd derivative of map, D vector space}%
\Index
   {the G\^ateaux \sD derivative of map $f$ of division ring $D$}%
   {Gateaux starD derivative of map, division ring}%
\Index
   {generating set}%
   {generating set}%
\Index
   {generator of linear map}%
   {generator of linear map}%
\Index
   {geodetic effect}%
   {geodetic effect}%
\Index
   {geometric object}%
   {geometric object}%
\Index
   {group algebra}%
   {group algebra}%
\Index
   {group of automorphisms of representation}%
   {group of automorphisms of representation}%
\SetIndexSpace%
\Index
   {Hadamard inverse of matrix}%
   {Hadamard inverse of matrix}%
\Index
   {Hamel basis}%
   {Hamel basis}%
\Index
   {hermitian conjugated vector}%
   {hermitian conjugated vector}%
\Index
   {hermitian conjugation in division ring}%
   {hermitian conjugation, division ring}%
\Index
   {hermitian matrix}%
   {hermitian matrix}%
\Index
   {hermitian metric on division ring}%
   {hermitian metric on division ring}%
\Index
   {hermitian scalar product in $D$\Hyph vector space}%
   {hermitian scalar product, vector space}%
\Index
   {hermitian scalar product on division ring}%
   {hermitian scalar product on division ring}%
\Index
   {highest common factor}%
   {highest common factor}%
\Index
   {holomorphic map}%
   {holomorphic map}%
\Index
   {holonomic coordinates of connection}%
   {holonomic coordinates of connection}%
\Index
   {holonomic coordinates of vector}%
   {vector holonomic coordinates}%
\Index
   {homogeneous bundle of fibered group}%
   {homogeneous bundle of fibered group}%
\Index
   {homogeneous linear geometric object}%
   {homogeneous linear geometric object}%
\Index
   {homogeneous map of degree $k$ over field $F$}%
   {homogeneous map of degree over field, D vector space}%
\Index
   {homogeneous polynomial}%
   {homogeneous polynomial}%
\Index
   {homogeneous space}%
   {homogeneous space}%
\Index
   {homomorphic image}%
   {homomorphic image}%
\Index
   {homomorphism}%
   {homomorphism}%
\Index
   {homomorphism of fibered groups}%
   {homomorphism of fibered groups}%
\Index
   {homomorphism of fibered universal algebras}%
   {homomorphism of fibered universal algebras}%
\Index
   {homomorphism of vector space}%
   {homomorphism of vector space}%
\Index
   {horizontal component of vector}%
   {horizontal component of vector}%
\Index
   {horizontal subspace}%
   {horizontal subspace}%
\Index
   {horizontal vector}%
   {horizontal vector}%
\Index
   {hyperbolic cosine}%
   {hyperbolic cosine}%
\Index
   {hyperbolic sine}%
   {hyperbolic sine}%
\SetIndexSpace%
\Index
   {ideal of algebra}%
   {ideal of algebra}%
\Index
   {image of map}%
   {image of map}%
\Index
   {indefinite integral}%
   {indefinite integral}%
\Index
   {independent points}%
   {independent points}%
\Index
   {induction over diagram of representations}%
   {induction over diagram of representations}%
\Index
   {infinitesimal generator of representation}%
   {infinitesimal generator}%
\Index
   {infinitesimal generators of group Lie}%
   {infinitesimal generators of group Lie}%
\Index
   {integrable differential equation}%
   {integrable differential equation}%
\Index
   {integrable differential form}%
   {integrable differential form}%
\Index
   {integrable form}%
   {integrable form}%
\Index
   {integrable map}%
   {integrable map}%
\Index
   {integral of differential $1$\Hyph form along path}%
   {integral of differential 1 form along path}%
\Index
   {invariance principle in tower of representations of universal algebras}%
   {invariance principle, tower of representations g}%
\Index
   {inverse fibered correspondence}%
   {inverse fibered correspondence}%
\Index
   {inverse reduced fibered correspondence}%
   {inverse reduced fibered correspondence}%
\Index
   {involution in quaternion algebra}%
   {involution, quaternion algebra}%
\Index
   {isomorphism}%
   {isomorphism}%
\Index
   {isomorphism of fibered $\Omega$\Hyph algebras}%
   {isomorphism of fibered Omega-algebras}%
\Index
   {isomorphism of repesentations of $\Omega$\Hyph algebra}%
   {isomorphism of repesentations of Omega algebra}%
\Index
   {isomorphism of vector spaces}%
   {isomorphism of vector spaces}%
\Index
   {isotropic vector}%
   {isotropic vector}%
\Index
   {Lebesgue integral}%
   {Lebesgue integral}%
\SetIndexSpace%
\Index
   {$(^j_i)$\hyph $\RCcirc$\Hyph quasideterminant}%
   {j i RCcirc-quasideterminant}%
\Index
   {the Jacobi matrix of map}%
   {Jacobi matrix of map}%
\Index
   {Jacobian complete system of differential equations}%
   {Jacobian complete system of differential equations}%
\Index
   {Jacobian complete system of \drv differential equations}%
   {Jacobian complete system of drc differential equations}%
\Index
   {$(ji)$\hyph quasideterminant}%
   {j i quasideterminant}%
\Index
   {the Jacobi\Hyph G\^ateaux matrix of map}%
   {Jacobi Gateaux matrix of map}%
\SetIndexSpace%
\Index
   {kernel of group\Hyph homomorphism}%
   {ker group-homomorphism}%
\Index
   {kernel of homomorphism}%
   {kernel of homomorphism}%
\Index
   {kernel of inefficiency of \Ts{G}representation}%
   {kernel of inefficiency of G* representation}%
\Index
   {kernel of inefficiency of representation of fibered group}%
   {kernel of inefficiency of representation of fibered group}%
\Index
   {kernel of inefficiency of representation of group}%
   {kernel of inefficiency of representation of group}%
\Index
   {kernel of linear map}%
   {kernel of linear map}%
\Index
   {kernel of map}%
   {kernel of map}%
\Index
   {Killing equation}%
   {Killing equation}%
\Index
   {Killing equation of second type}%
   {Killing equation second type}%
\Index
   {Killing vector of second type}%
   {Killing vector second type}%
\Index
   {Kronecker symbol}%
   {Kronecker symbol}%
\SetIndexSpace%
\Index
   {latitude}%
   {latitude}%
\Index
   {leading coefficient of polynomial}%
   {leading coefficient of polynomial}%
\Index
   {Lebesgue extension of measure}%
   {Lebesgue extension of measure}%
\Index
   {Lebesgue measurable set}%
   {Lebesgue measurable}%
\Index
   {Lebesgue measure}%
   {Lebesgue measure}%
\Index
   {left $A$\Hyph module}%
   {left A module}%
\Index
   {left $A$\Hyph vector space}%
   {left A vector space}%
\Index
   {left $A$\Hyph column space}%
   {left A-column space}%
\Index
   {left $A$\Hyph row space}%
   {left A-row space}%
\Index
   {left cofactor of entry of matrix}%
   {left cofactor, matrix}%
\Index
   {left coset}%
   {left coset}%
\Index
   {left $D$\hyph vector space of columns}%
   {left vector space of columns}%
\Index
   {left $D$\hyph vector space of rows}%
   {left vector space of rows}%
\Index
   {left defined Lie algebra of Lie group}%
   {left defined Lie algebra}%
\Index
   {left double cofactor of entry of matrix}%
   {left double cofactor}%
\Index
   {left fraction}%
   {left fraction}%
\Index
   {left ideal of algebra}%
   {left ideal of algebra}%
\Index
   {left invariant vector field}%
   {left invariant vector}%
\Index
   {left linear combination}%
   {left linear combination}%
\Index
   {left linear dependent}%
   {left linear dependent}%
\Index
   {left module}%
   {left module}%
\Index
   {left principal ideal}%
   {left principal ideal}%
\Index
   {left shift of module}%
   {left shift of module}%
\Index
   {left shift on fibered group}%
   {left shift, fibered group}%
\Index
   {left shift on group}%
   {left shift}%
\Index
   {left shift on group}%
   {left shift, group}%
\Index
   {left structural constant of Lie algebra}%
   {left structural constant of Lie algebra}%
\Index
   {left vector space}%
   {left vector space}%
\Index
   {left zero divisor}%
   {left zero divisor}%
\Index
   {left-ordered cycle notation of permutation}%
   {left-ordered cycle notation of permutation}%
\Index
   {left\Hyph side $A_1$\Hyph representation}%
   {left-side A representation}%
\Index
   {left\Hyph side product}%
   {left-side product}%
\Index
   {left-side product of map over scalar}%
   {left-side product of map over scalar}%
\Index
   {left\Hyph side product of vector over scalar}%
   {left-side product of vector over scalar}%
\Index
   {left-side representation}%
   {left-side representation}%
\Index
   {left-side representation of fibered $\Omega$\Hyph algebra}%
   {left-side representation of fibered Omega-algebra}%
\Index
   {left-side representation of $\Omega_1$\Hyph algebra $A$ in $\Omega_2$\Hyph algebra $M$}%
   {left-side representation of algebra}%
\Index
   {left-side transformation}%
   {left-side transformation}%
\Index
   {left-side transformation on bundle}%
   {left-side transformation of bundle}%
\Index
   {Lie algebra of Lie group}%
   {algebra Lie group Lie}%
\Index
   {Lie derivative}%
   {Lie derivative}%
\Index
   {Lie derivative of connection}%
   {Lie derivative of connection}%
\Index
   {Lie derivative of metric}%
   {Lie derivative of metric}%
\Index
   {Lie group basic operators}%
   {Lie group basic operators}%
\Index
   {lift of correspondence}%
   {lift of correspondence}%
\Index
   {lift of mapping}%
   {lift of map}%
\Index
   {limit of correspondence with respect to the filter}%
   {limit of correspondence with respect to the filter}%
\Index
   {limit of filter}%
   {limit of filter}%
\Index
   {limit of sequence}%
   {limit of sequence}%
\Index
   {limit set of filter}%
   {limit set of filter}%
\Index
   {linear combination}%
   {linear combination}%
\Index
   {linear functional}%
   {linear functional}%
\Index
   {linear \Ts{G}representation}%
   {linear G* representation}%
\Index
   {linear geometric object}%
   {linear geometric object}%
\Index
   {linear homogeneous equation}%
   {linear homogeneous equation}%
\Index
   {linear homomorphism}%
   {linear homomorphism}%
\Index
   {linear map}%
   {linear map}%
\Index
   {linear map generated by map}%
   {linear map generated by map}%
\Index
   {linear map of division ring}%
   {linear map of division ring}%
\Index
   {linear representation of group}%
   {linear representation of group}%
\Index
   {linear representation of Lie group}%
   {linear representation of Lie group}%
\Index
   {linear span}%
   {linear span, vector space}%
\Index
   {linear transformation group}%
   {linear transformation group}%
\Index
   {linear transformation of affine space}%
   {linear transformation, affine space}%
\Index
   {linearly dependent}%
   {linearly dependent}%
\Index
   {linearly dependent set}%
   {linearly dependent set}%
\Index
   {linearly dependent vector fields}%
   {linearly dependent vector fields}%
\Index
   {linearly independent set}%
   {linearly independent set}%
\Index
   {little group}%
   {little group}%
\Index
   {local reference frame}%
   {local reference frame}%
\Index
   {locally compact at point $p$ space}%
   {locally compact at point space}%
\Index
   {locally compact space}%
   {locally compact space}%
\Index
   {longitude}%
   {longitude}%
\Index
   {Lorentz transformation}%
   {Lorentz transformation}%
\SetIndexSpace%
\Index
   {$m$\Hyph dimensional parallelepiped}%
   {m dimensional parallelepiped}%
\Index
   {$m$\Hyph vector}%
   {m-vector}%
\Index
   {major submatrix}%
   {major submatrix}%
\Index
   {manifold with $D$\Hyph affine connections}%
   {manifold with D- affine connections}%
\Index
   {map continuous with respect to set of arguments}%
   {map continuous with respect to set of arguments}%
\Index
   {map differentiable in the G\^ateaux sense}%
   {map differentiable in Gateaux sense}%
\Index
   {map is compatible with operation}%
   {map is compatible with operation}%
\Index
   {map of conjugation}%
   {map of conjugation}%
\Index
   {map of $\gi n$ $D$\Hyph valued variables}%
   {map of n D valued variables}%
\Index
   {map of type $G$ on manifold}%
   {map of type G on manifold}%
\Index
   {map polylinear over finite dimensional algebras}%
   {map polylinear over finite dimensional algebras}%
\Index
   {map projective over commutative ring}%
   {map projective over commutative ring}%
\Index
   {mapping of rings polylinear over commutative ring}%
   {map polylinear over commutative ring, ring}%
\Index
   {mapping space}%
   {mapping space}%
\Index
   {matrix}%
   {matrix}%
\Index
   {matrix of antilinear homomorphism}%
   {matrix of antilinear homomorphism}%
\Index
   {matrix of bilinear function}%
   {matrix of bilinear function}%
\Index
   {matrix of endomorphisms of $\Omega$\Hyph algebra}%
   {matrix of endomorphisms of Omega algebra}%
\Index
   {matrix of fibered \Drc linear map}%
   {matrix of fibered drc linear map}%
\Index
   {matrix of homomorphism}%
   {matrix of homomorphism}%
\Index
   {matrix of linear homomorphism}%
   {matrix of linear homomorphism}%
\Index
   {matrix of linear map}%
   {matrix of linear map}%
\Index
   {matrix of linear maps}%
   {matrix of linear maps}%
\Index
   {matrix of maps}%
   {matrix of maps}%
\Index
   {matrix of quadratic map}%
   {matrix of quadratic map, division ring}%
\Index
   {Maxwell equation}%
   {Maxwell equation}%
\Index
   {measurable map}%
   {measurable map}%
\Index
   {measure}%
   {measure}%
\Index
   {method of successive differentiation}%
   {method of successive differentiation}%
\Index
   {metric tensor in Minkowski space}%
   {metric tensor, Minkowski space}%
\Index
   {metric-affine manifold}%
   {metric-affine manifold}%
\Index
   {Minkowski space}%
   {Minkowski space, Finsler}%
\Index
   {minor matrix}%
   {minor matrix}%
\Index
   {module over ring}%
   {module over ring}%
\Index
   {monomial of power $k$}%
   {monomial of power}%
\Index
   {monomorphism}%
   {monomorphism}%
\Index
   {morphism from diagram of representations into diagram of representations}%
   {morphism from diagram of representations into diagram of representations}%
\Index
   {morphism from tower of representations into tower of representations}%
   {morphism from tower of representations into tower of representations}%
\Index
   {morphism of fibered \Ts representations from $\Bundle F$ into $\Bundle G$}%
   {morphism of fibered representations from f into g}%
\Index
   {morphism of representation $f$}%
   {morphism of representation f}%
\Index
   {morphism of representations from $f$ into $g$}%
   {morphism of representations from f into g}%
\Index
   {morphism of representations of $\Omega_1$\Hyph algebra in $\Omega_2$\Hyph algebra}%
   {morphism of representations of Omega1 algebra in Omega2 algebra}%
\Index
   {morphism of \Ts representations of fibered $\Omega$\Hyph algebra}%
   {morphism of representations of fibered Omega algebra}%
\Index
   {motion of Minkowski space}%
   {motion, Minkowski space}%
\Index
   {movement on basis manifold}%
   {movement transformation}%
\Index
   {multiplicative map}%
   {multiplicative map}%
\Index
   {multiplicative $\Omega$\Hyph group}%
   {multiplicative Omega group}%
\SetIndexSpace%
\Index
   {$n$\Hyph ary fibered relation}%
   {fibered relation}%
\Index
   {$n$\Hyph ary operation on set}%
   {n-ary operation on set}%
\Index
   {natural homomorphism}%
   {natural homomorphism}%
\Index
   {neutral element of operation}%
   {neutral element of operation}%
\Index
   {nonmetricity}%
   {nonmetricity}%
\Index
   {nonsingular bilinear function}%
   {nonsingular bilinear function}%
\Index
   {nonsingular system of linear equations}%
   {nonsingular system of linear equations}%
\Index
   {nonsingular system of right $A^*$\Hyph linear equations}%
   {nonsingular system of right AU* linear equations}%
\Index
   {nonsingular tensor}%
   {nonsingular tensor}%
\Index
   {nonsingular transformation}%
   {nonsingular transformation}%
\Index
   {norm in quaternion algebra}%
   {norm, quaternion algebra}%
\Index
   {norm of functional}%
   {norm of functional}%
\Index
   {norm of map}%
   {norm of map}%
\Index
   {norm of octonion}%
   {norm of octonion}%
\Index
   {norm of operation}%
   {norm of operation}%
\Index
   {norm of polylinear map}%
   {norm of polymap}%
\Index
   {norm of representation}%
   {norm of representation}%
\Index
   {norm on $D$\Hyph algebra}%
   {norm on D algebra}%
\Index
   {norm on $D$\Hyph vector space}%
   {norm on D vector space}%
\Index
   {norm on $D$\Hyph module}%
   {norm on D module}%
\Index
   {norm on $\Omega$\Hyph group}%
   {norm on Omega group}%
\Index
   {norm on ring}%
   {norm on ring}%
\Index
   {normal basis}%
   {normal basis}%
\Index
   {normal subgroup}%
   {normal subgroup}%
\Index
   {normed $D$\Hyph algebra}%
   {normed D algebra}%
\Index
   {normed $D$\Hyph module}%
   {normed D module}%
\Index
   {normed $D$\Hyph vector space}%
   {normed D vector space}%
\Index
   {normed $\Omega$\Hyph group}%
   {normed Omega group}%
\Index
   {normed ring}%
   {normed ring}%
\Index
   {not complete group}%
   {not complete group}%
\Index
   {not complete $\Omega$\Hyph algebra}%
   {not complete Omega algebra}%
\Index
   {nucleus of $D$\Hyph algebra $A$}%
   {nucleus of algebra}%
\SetIndexSpace%
\Index
   {octonion algebra}%
   {octonion algebra}%
\Index
   {open ball}%
   {open ball}%
\Index
   {open set}%
   {open set}%
\Index
   {operation on bundle}%
   {operation on bundle}%
\Index
   {operation on set}%
   {operation on set}%
\Index
   {operator domain}%
   {operator domain}%
\Index
   {opposite algebra to algebra $P$}%
   {opposite algebra}%
\Index
   {opposite fibered preordering}%
   {opposite fibered preordering}%
\Index
   {orbit of linear map}%
   {orbit of linear map}%
\Index
   {orbit of representation}%
   {orbit of representation}%
\Index
   {orbit of representation of fibered group}%
   {orbit of representation of fibered group}%
\Index
   {orbit of representation of group}%
   {orbit of representation of group}%
\Index
   {origin of coordinate system of affine space}%
   {origin of coordinate system of affine space}%
\Index
   {orthogonal basis in Minkowski space}%
   {orthogonal basis, Minkowski space}%
\Index
   {orthogonality in Minkowski space}%
   {Minkowski orthogonality}%
\Index
   {orthonormal basis}%
   {Orthonormal Basis, division ring}%
\Index
   {orthonormal basis in Minkowski space}%
   {orthonormal basis, Minkowski space}%
\Index
   {orthonornal basis}%
   {Orthonornal Basis}%
\Index
   {outer measure}%
   {outer measure}%
\SetIndexSpace%
\Index
   {parallel shift of affine space}%
   {parallel shift, affine space}%
\Index
   {parallelogram}%
   {parallelogram}%
\Index
   {parity of permutation}%
   {parity of permutation}%
\Index
   {partial derivative}%
   {partial derivative}%
\Index
   {partial derivative of second order}%
   {partial derivative of second order}%
\Index
   {partial linear map}%
   {partial linear map}%
\Index
   {passive $G$\Hyph representation}%
   {passive G representation}%
\Index
   {passive representation}%
   {passive representation}%
\Index
   {passive representation in basis manifold}%
   {passive representation in basis manifold}%
\Index
   {passive representation of group $G(\Vector f)$ in basis manifold of tower of representations}%
   {passive representation in basis manifold, tower of representations}%
\Index
   {passive transformation of basis manifold}%
   {passive transformation, vector space}%
\Index
   {passive transformation of basis manifold}%
   {passive transformation of basis}%
\Index
   {passive transformation of the basis manifold of tower of representations}%
   {passive transformation of basis, tower of representations}%
\Index
   {passive transformation on basis manifold}%
   {passive transformation}%
\Index
   {permutability property of trace}%
   {permutability property of trace}%
\Index
   {permutation}%
   {permutation}%
\Index
   {pfaffian derivative}%
   {pfaffian derivative}%
\Index
   {polyadditive map}%
   {polyadditive map}%
\Index
   {polylinear map}%
   {polylinear map}%
\Index
   {polylinear skew symmetric map}%
   {polylinear map skew symmetric}%
\Index
   {polylinear symmetric map}%
   {polylinear map symmetric}%
\Index
   {polymorphism of representations}%
   {polymorphism of representations}%
\Index
   {polynomial}%
   {polynomial}%
\Index
   {polyvector}%
   {polyvector}%
\Index
   {potential energy}%
   {potential energy}%
\Index
   {power of measure}%
   {power of measure}%
\Index
   {prime $A$\Hyph number}%
   {prime A number}%
\Index
   {principal ideal}%
   {principal ideal}%
\Index
   {principle of covariance}%
   {principle of covariance}%
\Index
   {product in category}%
   {product in category}%
\Index
   {product of geometric object and constant}%
   {product of geometric object and constant}%
\Index
   {product of geometric object and constant in vector space}%
   {product of geometric object and constant, vector space}%
\Index
   {product of homomorphisms}%
   {product of homomorphisms}%
\Index
   {product of measures}%
   {product of measures}%
\Index
   {product of morphisms of diagram of representations}%
   {product of morphisms of diagram of representations}%
\Index
   {product of morphisms of representations of universal algebra}%
   {product of morphisms of representations of universal algebra}%
\Index
   {product of morphisms of tower of representations}%
   {product of morphisms of tower of representations}%
\Index
   {product of morphisms of \Ts representations of fibered $\Omega$\Hyph algebra}%
   {product of morphisms of representations of fibered Omega algebra}%
\Index
   {product of polynomials}%
   {product of polynomials}%
\Index
   {product of rings of sets}%
   {product of rings of sets}%
\Index
   {projection of bundle $\Bundle E$ along fiber $E$}%
   {projection of bundle along fiber}%
\Index
   {projective map is continuous in direction over field}%
   {projective map is continuous in direction over field}%
\Index
   {pseudo\Hyph Euclidean metric on division ring}%
   {pseudo-Euclidean metric on division ring}%
\Index
   {pseudo\Hyph Euclidean scalar product in $D$\Hyph vector space}%
   {pseudo-Euclidean scalar product, vector space}%
\Index
   {pseudo-Euclidean scalar product on division ring}%
   {pseudo-Euclidean scalar product on division ring}%
\SetIndexSpace%
\Index
   {quadratic equation}%
   {quadratic equation}%
\Index
   {quadratic form in division ring}%
   {quadratic form, division ring}%
\Index
   {quadratic map of division ring}%
   {Quadratic Map of Division Ring}%
\Index
   {quasi affine transformation on basis manifold}%
   {quasi affine transformation}%
\Index
   {quasi affine transformation on basis manifold}%
   {quasi affine drc transformation}%
\Index
   {quasi movement on basis manifold}%
   {quasi movement, division ring}%
\Index
   {quasi movement on basis manifold}%
   {quasi movement}%
\Index
   {quasibasis}%
   {quasibasis}%
\Index
   {quasiclosed ring of maps}%
   {quasiclosed ring of maps}%
\Index
   {quasideterminant}%
   {quasideterminant definition}%
\Index
   {quasiexponent}%
   {quasiexponent}%
\Index
   {quasimotion of Minkowski space}%
   {Quasimotion, Minkowski space}%
\Index
   {quaternion algebra}%
   {quaternion algebra}%
\Index
   {quaternion algebra $E$ over the field $F$}%
   {quaternion algebra over the field}%
\Index
   {quotient}%
   {quotient divided by}%
\Index
   {quotient bundle}%
   {quotient bundle}%
\SetIndexSpace%
\Index
   {$(\aUD{}ji)$\hyph \RC quasideterminant}%
   {j i rc-quasideterminant}%
\Index
   {\sups row of matrix}%
   {r row}%
\Index
   {$R$\Hyph module}%
   {R- module}%
\Index
   {$r$\hyph row of matrix}%
   {r-row}%
\Index
   {rank of Hermitian matrix by principal minors}%
   {rank of Hermitian matrix by principal minors}%
\Index
   {rank of matrix}%
   {rank of matrix}%
\Index
   {rank of quadratic map of division ring}%
   {rank of quadratic map, division ring}%
\Index
   {\RC exponent}%
   {RC exponent}%
\Index
   {\RC inverse element of biring}%
   {rc-inverse element}%
\Index
   {\RC matrix group}%
   {rc-matrix group}%
\Index
   {\RC nonsingular matrix}%
   {r? nonsingular matrix}%
\Index
   {\RC power}%
   {rc power}%
\Index
   {\RC product (product of row over column)}%
   {rc-product}%
\Index
   {$\RCcirc$\Hyph product of matrices of maps}%
   {rc product of matrices of maps}%
\Index
   {\RC quasideterminant}%
   {rc-quasideterminant}%
\Index
   {\RC singular matrix}%
   {rc singular matrix}%
\Index
   {\RC inverse matrix}%
   {rc-inverse matrix}%
\Index
   {$\RCcirc$\Hyph nonsingular matrix}%
   {RCcirc nonsingular matrix}%
\Index
   {$\RCcirc$\Hyph nonsingular system of additive equations}%
   {RCcirc nonsingular system of additive equations}%
\Index
   {$\RCcirc$\Hyph quasideterminant}%
   {RCcirc-quasideterminant definition}%
\Index
   {$\RCcirc$\Hyph singular matrix}%
   {RCcirc singular matrix}%
\Index
   {\rcd affine plane}%
   {rcd affine plane}%
\Index
   {\rcd affine space}%
   {rcd affine space}%
\Index
   {\rcd vector}%
   {rcd vector}%
\Index
   {\rcd vector space}%
   {rcd vector space}%
\Index
   {reduced Cartesian product of bundles}%
   {reduced Cartesian product of bundles}%
\Index
   {reduced Cartesian product of total spaces}%
   {reduced Cartesian product of total spaces}%
\Index
   {reduced fibered correspondence from $\Bundle{A}$ to $\Bundle B$}%
   {reduced fibered correspondence from A to B}%
\Index
   {reduced fibered correspondence in $\Bundle{A}$}%
   {reduced fibered correspondence in A}%
\Index
   {reduced morphism of representations}%
   {reduced morphism of representations}%
\Index
   {reduced polymorphism of representations}%
   {reduced polymorphism of representations}%
\Index
   {reduced quadratic equation}%
   {reduced quadratic equation}%
\Index
   {reducible biring}%
   {reducible biring}%
\Index
   {reference frame in event space}%
   {reference frame in event space}%
\Index
   {reference frame manifold}%
   {reference frame manifold}%
\Index
   {reflexive $2$\Hyph ary fibered relation}%
   {reflexive 2 ary fibered relation}%
\Index
   {reflexive correspondence}%
   {reflexive correspondence}%
\Index
   {regular endomorphism}%
   {regular endomorphism}%
\Index
   {regular endomorphism of tower of representations}%
   {regular endomorphism of tower of representations}%
\Index
   {regular quadratic map in division ring}%
   {regular quadratic map, division ring}%
\Index
   {relatively prime $A$\Hyph numbers}%
   {relatively prime A numbers}%
\Index
   {remainder of the division}%
   {remainder of the division}%
\Index
   {representation conjugated to representation}%
   {representation conjugated to representation}%
\Index
   {\Ts{A}representation in $\Omega_2$\Hyph algebra}%
   {A* representation of algebra}%
\Index
   {representation of group}%
   {representation of group}%
\Index
   {representation of $\Omega$\Hyph algebra in representation}%
   {representation of Omega algebra in representation}%
\Index
   {representation of $\Omega$\Hyph algebra in tower of representations}%
   {representation of Omega algebra in tower of representations}%
\Index
   {representation of $\Omega$\Hyph algebra $A$ in category $\mathcal B$}%
   {representation of Omega algebra in category}%
\Index
   {\sT{A}representation of $\Omega_1$\Hyph algebra $A$ in $\Omega_2$\Hyph algebra}%
   {*A representation of algebra}%
\Index
   {representation of $\Omega_1$\Hyph algebra $A$ in $\Omega_2$\Hyph algebra $M$}%
   {representation of algebra}%
\Index
   {representative of geometric object}%
   {representative of geometric object}%
\Index
   {restriction of correspondence $\Phi$ to set $C$}%
   {restriction of correspondence}%
\Index
   {right $A_*$\Hyph vector space}%
   {right A subs vector space}%
\Index
   {right $A$\Hyph vector space}%
   {right A vector space}%
\Index
   {right $A$\Hyph column space}%
   {right A-column space}%
\Index
   {right $A$\Hyph row space}%
   {right A-row space}%
\Index
   {right cofactor of entry of matrix}%
   {right cofactor, matrix}%
\Index
   {right coset}%
   {right coset}%
\Index
   {right $D$\Hyph module}%
   {right D module}%
\Index
   {right $D$\hyph vector space of columns}%
   {right vector space of columns}%
\Index
   {right $D$\hyph vector space of rows}%
   {right vector space of rows}%
\Index
   {right defined Lie algebra of Lie group}%
   {right defined Lie algebra}%
\Index
   {right double cofactor of entry of matrix}%
   {right double cofactor}%
\Index
   {right fraction}%
   {right fraction}%
\Index
   {right ideal of algebra}%
   {right ideal of algebra}%
\Index
   {right invariant vector field}%
   {right invariant vector}%
\Index
   {right linear combination}%
   {right linear combination}%
\Index
   {right module}%
   {right module}%
\Index
   {right module over $D$\Hyph algebra $A$}%
   {right module over algebra}%
\Index
   {right principal ideal}%
   {right principal ideal}%
\Index
   {right shift on group}%
   {right shift}%
\Index
   {right shift on group}%
   {right shift, group}%
\Index
   {right structural constant of Lie algebra}%
   {right structural constant of Lie algebra}%
\Index
   {right vector space}%
   {right vector space}%
\Index
   {right zero divisor}%
   {right zero divisor}%
\Index
   {right-ordered cycle notation of permutation}%
   {right-ordered cycle notation of permutation}%
\Index
   {right\Hyph side $A_1$\Hyph representation}%
   {right-side A representation}%
\Index
   {right\Hyph side product}%
   {right-side product}%
\Index
   {right\Hyph side product of vector over scalar}%
   {right-side product of vector over scalar}%
\Index
   {right-side representation}%
   {right-side representation}%
\Index
   {right-side representation of fibered $\Omega$\Hyph algebra}%
   {right-side representation of fibered Omega-algebra}%
\Index
   {right-side representation of $\Omega_1$\Hyph algebra $A$ in $\Omega_2$\Hyph algebra $M$}%
   {right-side representation of algebra}%
\Index
   {right-side transformation}%
   {right-side transformation}%
\Index
   {ring has characteristic $0$}%
   {ring has characteristic 0}%
\Index
   {ring has characteristic $p$}%
   {ring has characteristic p}%
\Index
   {ring of sets}%
   {ring of sets}%
\Index
   {ring of sets generated by semiring of sets}%
   {ring of sets generated by semiring}%
\Index
   {ring with conjugation}%
   {ring with conjugation}%
\Index
   {root of polynomial}%
   {root of polynomial}%
\Index
   {row $*D$\Hyph vector}%
   {row *D vector}%
\Index
   {row $D*$\Hyph vector}%
   {row D* vector}%
\Index
   {row determinant}%
   {row determinant}%
\Index
   {row vector}%
   {row vector}%
\SetIndexSpace%
\Index
   {$\star A$\Hyph module}%
   {starA-module}%
\Index
   {scalar algebra of algebra}%
   {scalar algebra of algebra}%
\Index
   {scalar algebra of ring}%
   {scalar algebra of ring}%
\Index
   {scalar of element of algebra}%
   {scalar of algebra}%
\Index
   {scalar of element of ring}%
   {scalar of ring}%
\Index
   {scalar potential}%
   {scalar potential}%
\Index
   {Schauder basis}%
   {Schauder basis}%
\Index
   {second axiom of countability}%
   {second axiom of countability}%
\Index
   {second Newton law}%
   {Second Newton law}%
\Index
   {section of bundle}%
   {section of bundle}%
\Index
   {semigroup}%
   {semigroup}%
\Index
   {semiring of sets}%
   {semiring of sets}%
\Index
   {sequence converges}%
   {sequence converges}%
\Index
   {sequence converges almost everywhere}%
   {converges almost everywhere}%
\Index
   {sequence converges uniformly}%
   {sequence converges uniformly}%
\Index
   {series converges normally}%
   {series converges normally}%
\Index
   {set admits operation}%
   {set admits operation}%
\Index
   {set is closed with respect to operation}%
   {set is closed with respect to operation}%
\Index
   {set is dense in set}%
   {dense in set}%
\Index
   {set of coordinates of representation}%
   {coordinate set of representation}%
\Index
   {set of invertible elements of algebra}%
   {set of invertible elements of algebra}%
\Index
   {set of $\Omega_2$\Hyph words of representation}%
   {word set of representation}%
\Index
   {set of tuples of coordinates of diagram of representations}%
   {coordinate set of diagram of representations}%
\Index
   {set of tuples of coordinates of tower of representations}%
   {coordinate set of tower of representations}%
\Index
   {set of tuples of $\Omega$\Hyph words}%
   {set of tuples of Omega words}%
\Index
   {set of tuples of $\Vector\Omega$\Hyph words of tower of representations}%
   {word set of tower of representations}%
\Index
   {set of zeros of algebra}%
   {set of zeros of algebra}%
\Index
   {similarity transformation}%
   {similarity transformation}%
\Index
   {simple map}%
   {simple map}%
\Index
   {simple polyvector}%
   {simple polyvector}%
\Index
   {simplex}%
   {simplex}%
\Index
   {sine}%
   {sine}%
\Index
   {single transitive representation of fibered $\Omega$\Hyph algebra}%
   {single transitive representation of fibered Omega-algebra}%
\Index
   {single transitive representation of group}%
   {single transitive representation of group}%
\Index
   {single transitive representation of $\Omega$\Hyph algebra $A$}%
   {single transitive representation of algebra}%
\Index
   {singular endomorphism}%
   {singular endomorphism}%
\Index
   {singular linear map}%
   {singular linear map}%
\Index
   {skew product of vectors}%
   {skew product of vectors}%
\Index
   {skew symmetric polylinear map}%
   {skew symmetric polylinear map}%
\Index
   {space of orbits of \Ts{G}representation}%
   {space of orbits of G* representation}%
\Index
   {space of orbits of left\Hyph side representation}%
   {space of orbits of left side representation}%
\Index
   {spacelike vector}%
   {spacelike vector}%
\Index
   {speed of deviation}%
   {speed of deviation}%
\Index
   {spherical coordinates}%
   {spherical coordinates}%
\Index
   {square root}%
   {square root}%
\Index
   {$(\mathcal S\RCstar,\mathcal T\RCstar)$\Hyph linear map of vector bundles}%
   {src trc linear map of vector bundles}%
\Index
   {($S\star$, $\star T$)\hyph bimodule}%
   {(Sstar,starT)-bimodule}%
\Index
   {stability group}%
   {stability group}%
\Index
   {stable set of representation}%
   {stable set of representation}%
\Index
   {standard component of derivative}%
   {standard component of derivative}%
\Index
   {standard component of the G\^ateaux derivative}%
   {standard component of Gateaux derivative}%
\Index
   {standard component of linear map}%
   {standard component of linear map}%
\Index
   {standard component of polylinear map}%
   {standard component of polylinear map}%
\Index
   {standard component of tensor}%
   {standard component of tensor}%
\Index
   {standard component over field $F$ of bilitnear map $f$}%
   {standard component of bilinear map, division ring}%
\Index
   {standard coordinates of basis}%
   {standard coordinates of basis}%
\Index
   {standard coordinates of basis}%
   {standard coordinates of basis}%
\Index
   {standard representation of the derivative}%
   {derivative, standard representation}%
\Index
   {standard representation of the G\^ateaux derivative}%
   {Gateaux derivative, standard representation}%
\Index
   {standard representation of linear map}%
   {linear map, standard representation}%
\Index
   {standard representation of matrix}%
   {Standard representation}%
\Index
   {standard representation of polylinear map}%
   {polylinear map, standard representation}%
\Index
   {standard representation of quadratic map of division ring over field $F$}%
   {quadratic map, standard representation, division ring}%
\Index
   {standard representation over field $F$ of bilinear map of division ring}%
   {bilinear map, standard representation, division ring}%
\Index
   {$\star R$\hyph module}%
   {starR-module}%
\Index
   {$\star D$\hyph product of vector over scalar}%
   {starD product of vector over scalar, vector space}%
\Index
   {starlike set}%
   {starlike set}%
\Index
   {\sT representation of fibered group}%
   {starT representation of fibered group}%
\Index
   {\sT representation of fibered group}%
   {starT representation of fibered group}%
\Index
   {\sT representation of fibered $\Omega$\Hyph algebra}%
   {starT representation of fibered Omega-algebra}%
\Index
   {\sT shift on fibered group}%
   {starT shift, fibered group}%
\Index
   {\sT transformation on bundle}%
   {starT transformation of bundle}%
\Index
   {structural constants}%
   {structural constants}%
\Index
   {subalgebra of $\Omega$\Hyph algebra}%
   {subalgebra of Omega-algebra}%
\Index
   {subbundle}%
   {subbundle}%
\Index
   {subbundle of $\mathcal D\star$\hyph vector space}%
   {subbundle of Dstar vector bundle}%
\Index
   {subgroup of $\Omega$\Hyph group}%
   {subgroup of Omega group}%
\Index
   {submodule}%
   {submodule}%
\Index
   {submodule generated by set}%
   {submodule generated by set}%
\Index
   {subrepresentation}%
   {subrepresentation}%
\Index
   {subrepresentation generated by set $X$}%
   {subrepresentation generated by set}%
\Index
   {subrepresentation of representation}%
   {subrepresentation of representation}%
\Index
   {sum of geometric objects in vector space}%
   {sum of geometric objects, vector space}%
\Index
   {sum of geometric objects}%
   {sum of geometric objects}%
\Index
   {sum of maps}%
   {sum of maps}%
\Index
   {sum of polynomials}%
   {sum of polynomials}%
\Index
   {superposition of coordinates}%
   {superposition of coordinates,}%
\Index
   {superposition of coordinates of the tower of representations $\Vector f$ and the element $\VX a$}%
   {superposition of coordinates, tower of representations}%
\Index
   {symmetric $2$\Hyph ary fibered relation}%
   {symmetric 2 ary fibered relation}%
\Index
   {symmetric bilinear map of $D$\Hyph vector space to division ring}%
   {symmetric bilinear map, vector space to division ring}%
\Index
   {symmetric correspondence}%
   {symmetric correspondence}%
\Index
   {symmetric polylinear map}%
   {symmetric polylinear map}%
\Index
   {symmetric polylinear mapping into associative algebra}%
   {polylinear map symmetric, associative algebra}%
\Index
   {symmetrization of polylinear map}%
   {symmetrization of polylinear map}%
\Index
   {symmetry group}%
   {symmetry group}%
\Index
   {symmetry group}%
   {SymmetryGroup}%
\Index
   {synchronization of reference frame}%
   {synchronization of reference frame}%
\Index
   {system of additive equations}%
   {system of additive equations}%
\Index
   {system of \drc linear equations}%
   {system of drc linear equations}%
\Index
   {system of left $A_*$\Hyph linear equations}%
   {system of left AD* linear equations}%
\Index
   {system of linear equations}%
   {system of linear equations}%
\Index
   {system of \rcd linear equations}%
   {system of rcd linear equations}%
\Index
   {system of right $A^*$\Hyph linear equations}%
   {system of right AU* linear equations}%
\SetIndexSpace%
\Index
   {$T_1$\Hyph space}%
   {T1 space}%
\Index
   {Taylor polynomial}%
   {Taylor polynomial, division ring}%
\Index
   {Taylor series}%
   {Taylor series, division ring}%
\Index
   {tensor inverse to tensor}%
   {inverse tensor}%
\Index
   {tensor power}%
   {tensor power}%
\Index
   {tensor product}%
   {tensor product}%
\Index
   {the Fr\'echet \Ds derivative of map $f$ of division ring $D$ at point $x$}%
   {Frechet Dstar derivative of map, division ring}%
\Index
   {timelike vector}%
   {timelike vector}%
\Index
   {topological $D$\Hyph vector space}%
   {topological D vector space}%
\Index
   {topological $D$\Hyph algebra}%
   {topological D algebra}%
\Index
   {topological division ring}%
   {topological division ring}%
\Index
   {topological ring}%
   {topological ring}%
\Index
   {torsion form}%
   {torsion form}%
\Index
   {torsion tensor}%
   {torsion tensor}%
\Index
   {tower of bundles}%
   {tower of bundles}%
\Index
   {tower of effective representations}%
   {tower of effective representations}%
\Index
   {tower of representations of $\Omega$\Hyph algebras}%
   {tower of representations of algebras}%
\Index
   {tower of subrepresentations}%
   {tower of subrepresentations}%
\Index
   {tower of subrepresentations of tower of representations $\Vector f$ generated by tuple of sets $\VX X$}%
   {subrepresentation generated by tuple of sets}%
\Index
   {trace of quaternion}%
   {trace, quaternion algebra}%
\Index
   {transformation coordinated with equivalence}%
   {transformation coordinated with equivalence}%
\Index
   {transformation of universal algebra}%
   {transformation of universal algebra}%
\Index
   {transformation on bundle}%
   {transformation of bundle}%
\Index
   {transitive $2$\Hyph ary fibered relation}%
   {transitive 2 ary fibered relation}%
\Index
   {transitive correspondence}%
   {transitive correspondence}%
\Index
   {transitive representation of fibered $\Omega$\Hyph algebra}%
   {transitive representation of fibered Omega-algebra}%
\Index
   {transitive representation of group}%
   {transitive representation of group}%
\Index
   {transitive representation of $\Omega$\Hyph algebra $A$}%
   {transitive representation of algebra}%
\Index
   {trivial kernel of homomorphism}%
   {trivial kernel}%
\Index
   {\Ts representation of fibered group}%
   {Tstar representation of fibered group}%
\Index
   {\Ts representation of fibered $\Omega$\Hyph algebra}%
   {Tstar representation of fibered Omega-algebra}%
\Index
   {tuple of equivalence generated by tower of representations $\Vector f$}%
   {tuple of equivalence of tower of representations}%
\Index
   {tuple of generating sets of tower of representations}%
   {tuple of generating sets of tower of representations}%
\Index
   {tuple of $\Omega$\Hyph words}%
   {tuple of Omega words}%
\Index
   {tuple of $\Vector{\Omega}$\Hyph words of element of tower of representations relative to tuple of generating sets}%
   {tuple of words relative to tuple of generating sets, tower of representations}%
\Index
   {tuple of stable sets of diagram of representations}%
   {tuple of stable sets of diagram of representations}%
\Index
   {tuple of stable sets of tower of representation}%
   {tuple of stable sets of tower of representations}%
\Index
   {twin representations}%
   {twin representations}%
\Index
   {twin representations of division ring}%
   {twin representations of division ring}%
\Index
   {twin representations of fibered group}%
   {twin representations of fibered group}%
\Index
   {twin representations of group}%
   {twin representations of group}%
\Index
   {type of geometric object}%
   {type of geometric object}%
\SetIndexSpace%
\Index
   {unit interval}%
   {unit interval}%
\Index
   {unit of ring of sets}%
   {unit of ring of sets}%
\Index
   {unit sphere in $D$\Hyph algebra}%
   {unit sphere in algebra}%
\Index
   {unit sphere in division ring}%
   {unit sphere in division ring}%
\Index
   {unit vector}%
   {unit vector}%
\Index
   {unital algebra}%
   {unital algebra}%
\Index
   {unital extension}%
   {unital extension}%
\Index
   {unital ring}%
   {unital ring}%
\Index
   {unitarity law}%
   {unitarity law}%
\Index
   {universal algebra}%
   {universal algebra}%
\Index
   {universally attracting object of category}%
   {universally attracting}%
\Index
   {universally repelling  object of category}%
   {universally repelling}%
\SetIndexSpace%
\Index
   {basis for vector  bundle}%
   {basis, vector bundle}%
\Index
   {valued division ring}%
   {valued division ring}%
\Index
   {vector}%
   {vector}%
\Index
   {vector $*A$\Hyph space}%
   {*A-vector space}%
\Index
   {vector bundle}%
   {vector bundle}%
\Index
   {vector module of algebra}%
   {vector module of algebra}%
\Index
   {vector module of ring}%
   {vector module of ring}%
\Index
   {vector of element of algebra}%
   {vector of algebra}%
\Index
   {vector of element of ring}%
   {vector of ring}%
\Index
   {vector potential}%
   {vector potential}%
\Index
   {vector space}%
   {vector space}%
\Index
   {vector space type}%
   {vector space type}%
\Index
   {vertical component of vector}%
   {vertical component of vector}%
\Index
   {vertical subspace}%
   {vertical subspace}%
\Index
   {vertical vector}%
   {vertical vector}%
\SetIndexSpace%
\Index
   {Wronskian determinant}%
   {Wronskian determinant}%
\Index
   {Wronskian matrix}%
   {Wronskian matrix}%
\SetIndexSpace%
\Index
   {zero divisor}%
   {zero divisor}%
\SetIndexSpace%
\Index
   {$\mu$\Hyph measurable map}%
   {mu measurable map}%
\SetIndexSpace%
\Index
   {$\Omega$\Hyph algebra}%
   {Omega-algebra}%
\Index
   {$\Omega$\Hyph group}%
   {Omega group}%
\Index
   {$\Omega$\Hyph groupoid}%
   {Omega groupoid}%
\Index
   {$\Omega$\Hyph linear mapping}%
   {Omega linear map}%
\Index
   {\(\Omega\)\Hyph ring}%
   {Omega ring}%
\Index
   {$\Omega_2$\Hyph word of element of representation relative to generating set}%
   {word of element relative to generating set, representation}%
\SetIndexSpace%
\Index
   {$\sigma$\Hyph algebra of sets}%
   {sigma algebra of sets}%
\Index
   {$\sigma$\Hyph ring of sets}%
   {sigma ring of sets}%
\Index
   {\(\sigma\)\Hyph additive measure}%
   {sigma-additive measure}%

\CloseIndex

%% file: Symbol.English.tex
\def\indexname{Special Symbols and Notations}
\OpenIndex

\SetIndexSpace
\Symb%
   {direct sum}%
   {direct sum}%
   {0}{0}%
\Symb%
   {unit interval}%
   {unit interval}%
   {0}{0}%

\SetIndexSpace
\Symb%
   {set of vectors whose expansion relative to the basis $\Basis e$ converges normally}%
   {A plus Schauder}%
   {A}{0}%
\Symb%
   {active representation in basis manifold}%
   {active representation in basis manifold}%
   {A}{0}%
\Symb%
   {$A$\Hyph algebra of polynomials over $D$\Hyph algebra $A$}%
   {algebra of polynomials over algebra}%
   {A}{0}%
\Symb%
   {algebra of polynomials over $D$\Hyph algebra $A$}%
   {algebra of polynomials over D algebra}%
   {A}{0}%
\Symb%
   {algebra of rational mappings of algebra $A$}%
   {algebra of rational mappings of algebra}%
   {A}{0}%
\Symb%
   {affine space}%
   {An}%
   {A}{0}%
\Symb%
   {associator of $D$\Hyph algebra}%
   {associator of algebra}%
   {A}{0}%
\Symb%
   {category of left-side representations of $\Omega_1$\Hyph algebra $A$}%
   {category of left-side representations of Omega1 algebra}%
   {A}{0}%
\Symb%
   {category of representations}%
   {category of representations}%
   {A}{0}%
\Symb%
   {commutator of $D$\Hyph algebra}%
   {commutator of algebra}%
   {A}{0}%
\Symb%
   {component of linear map}%
   {component of linear map, vector}%
   {A}{0}%
\Symb%
   {component $p$ of polylinear mapping $\Vector A$}%
   {component of polyadditive map, D vector space}%
   {A}{0}%
\Symb%
   {component of polylinear map}%
   {component of polylinear map, vector}%
   {A}{0}%
\Symb%
   {conjugated $D$\Hyph  module}%
   {conjugated D module}%
   {A}{0}%
\Symb%
   {coordinates of associator}%
   {coordinates of associator}%
   {A}{0}%
\Symb%
   {\CR power of element $A$ of biring}%
   {cr power}%
   {A}{0}%
\Symb%
   {\crd vector}%
   {crd vector}%
   {A}{0}%
\Symb%
   {\CR inverse matrix}%
   {cr-inverse matrix}%
   {A}{0}%
\Symb%
   {\CR product}%
   {cr-product}%
   {A}{0}%
\Symb%
   {\dcr vector}%
   {dcr vector}%
   {A}{0}%
\Symb%
   {derivative of left shift}%
   {derivative of left shift}%
   {A}{0}%
\Symb%
   {derivative of left shift in $1$\Hyph parameter Lie group}%
   {derivative of left shift, 1-Parameter Group}%
   {A}{0}%
\Symb%
   {derivative of left shift in $1$\Hyph parameter Lie D group}%
   {derivative of left shift, 1-Parameter Group, algebra}%
   {A}{0}%
\Symb%
   {derivative of right shift}%
   {derivative of right shift}%
   {A}{0}%
\Symb%
   {derivative of right shift in $1$\Hyph parameter Lie group}%
   {derivative of right shift, 1-Parameter Group}%
   {A}{0}%
\Symb%
   {derivative of right shift in $1$\Hyph parameter Lie D group}%
   {derivative of right shift, 1-Parameter Group, algebra}%
   {A}{0}%
\Symb%
   {derivative of left shift}%
   {derivative of Tstar shift}%
   {A}{0}%
\Symb%
   {\drc vector}%
   {drc vector}%
   {A}{0}%
\Symb%
   {coordinates of vector $a$ relative to Hamel basis}%
   {Hamel basis, coordinates}%
   {A}{0}%
\Symb%
   {hermitian conjugation in division ring}%
   {hermitian conjugation, division ring}%
   {A}{0}%
\Symb%
   {tensor inverse to tensor $a$}%
   {inverse tensor}%
   {A}{0}%
\Symb%
   {isomorphic}%
   {isomorphic}%
   {A}{0}%
\Symb%
   {$(^j_i)$\hyph\CR quasideterminant}%
   {j i CR quasideterminant definition}%
   {A}{0}%
\Symb%
   {$(ji)$\hyph quasideterminant of matrix $\bfA$}%
   {j i quasideterminant definition}%
   {A}{0}%
\Symb%
   {$(^j_i)$\hyph $\RCcirc$\Hyph quasideterminant}%
   {j i RCcirc-quasideterminant definition}%
   {A}{0}%
\Symb%
   {$(^j_i)$\hyph \RC quasideterminant}%
   {j i RC-quasideterminant definition}%
   {A}{0}%
\Symb%
   {left fraction}%
   {left fraction}%
   {A}{0}%
\Symb%
   {left principal ideal}%
   {left principal ideal}%
   {A}{0}%
\Symb%
   {left shift in $D$\Hyph algebra}%
   {left shift, D algebra}%
   {A}{0}%
\Symb%
   {linear combination}%
   {linear combination}%
   {A}{0}%
\Symb%
   {little group}%
   {little group}%
   {A}{0}%
\Symb%
   {transformation of matrix}%
   {matrix, replacing its column}%
   {A}{0}%
\Symb%
   {transformation of matrix}%
   {matrix, replacing its row}%
   {A}{0}%
\Symb%
   {minor matrix}%
   {minor matrix}%
   {A}{0}%
\Symb%
   {norm on $D$\Hyph module}%
   {norm on D module}%
   {A}{0}%
\Symb%
   {$\Omega$\Hyph algebra}%
   {Omega-algebra}%
   {A}{0}%
\Symb%
   {opposite algebra to algebra $A$}%
   {opposite algebra}%
   {A}{0}%
\Symb%
   {orbit of linear map}%
   {orbit of linear map}%
   {A}{0}%
\Symb%
   {derivative}%
   {overline nabla_l, definition 2}%
   {A}{0}%
\Symb%
   {partial linear map}%
   {partial linear map}%
   {A}{0}%
\Symb%
   {principal ideal}%
   {principal ideal}%
   {A}{0}%
\Symb%
   {quasideterminant of matrix $\bfA$}%
   {quasideterminant definition}%
   {A}{0}%
\Symb%
   {\RC power of element $A$ of biring}%
   {rc power}%
   {A}{0}%
\Symb%
   {$\RCcirc$\Hyph quasideterminant}%
   {RCcirc-quasideterminant definition}%
   {A}{0}%
\Symb%
   {\rcd vector}%
   {rcd vector}%
   {A}{0}%
\Symb%
   {\RC inverse matrix}%
   {rc-inverse matrix}%
   {A}{0}%
\Symb%
   {\RC product}%
   {rc-product}%
   {A}{0}%
\Symb%
   {\RC quasideterminant}%
   {RC-quasideterminant definition}%
   {A}{0}%
\Symb%
   {right principal ideal}%
   {right principal ideal}%
   {A}{0}%
\Symb%
   {right shift in $D$\Hyph algebra}%
   {right shift, D algebra}%
   {A}{0}%
\Symb%
   {coordinates of vector $a$ relative to Schauder basis}%
   {Schauder basis, coordinates}%
   {A}{0}%
\Symb%
   {set of additive maps}%
   {set additive maps}%
   {A}{0}%
\Symb%
   {set of homogeneous polynomials}%
   {set of homogeneous polynomials}%
   {A}{0}%
\Symb%
   {set of invertible elements of algebra $A$}%
   {set of invertible elements of algebra}%
   {A}{0}%
\Symb%
   {set of vectors generated by vector $v$}%
   {set of vectors generated by vector}%
   {A}{0}%
\Symb%
   {set of zeros of algebra $A$}%
   {set of zeros of algebra}%
   {A}{0}%
\Symb%
   {set of polylinear maps of rings $R_1$, ..., $R_n$ into module $S$}%
   {set polylinear maps, ring}%
   {A}{0}%
\Symb%
   {simplex}%
   {simplex}%
   {A}{0}%
\Symb%
   {skew product of vectors $\Vector a_1$, ..., $\Vector a_m$}%
   {skew product of vectors}%
   {A}{0}%
\Symb%
   {space of orbits of left\Hyph side representation}%
   {space of orbits of left side representation}%
   {A}{0}%
\Symb%
   {space of orbits of representation}%
   {space of orbits of representation}%
   {A}{0}%
\Symb%
   {space of orbits of effective right\Hyph side representation}%
   {space of orbits of right-side representation}%
   {A}{0}%
\Symb%
   {square root}%
   {square root}%
   {A}{0}%
\Symb%
   {stability group}%
   {stability group}%
   {A}{0}%
\Symb%
   {\sT shift}%
   {starT shift, fibered group}%
   {A}{0}%
\Symb%
   {tensor power of algebra $A$}%
   {tensor power of algebra}%
   {A}{0}%
\Symb%
   {anholonomic coordinates of vector}%
   {vector anholonomic coordinates}%
   {A}{0}%
\Symb%
   {holonomic coordinates of vector}%
   {vector holonomic coordinates}%
   {A}{0}%

\SetIndexSpace
\Symb%
   {basis manifold of tower of representations $\Vector f$}%
   {basis manifold tower of representations}%
   {B}{0}%
\Symb%
   {basis manifold of affine space}%
   {Basis Manifold, Affine Space}%
   {B}{0}%
\Symb%
   {basis manifold of central affine space}%
   {BCAn}%
   {B}{0}%
\Symb%
   {basis manifold of Euclid space}%
   {BEn}%
   {B}{0}%
\Symb%
   {Borel algebra}%
   {Borel algebra}%
   {B}{0}%
\Symb%
   {canonical remainder of the division}%
   {canonical remainder of the division}%
   {B}{0}%
\Symb%
   {Cartesian power}%
   {Cartesian power}%
   {B}{0}%
\Symb%
   {Cartesian power $\Bundle A$ of bundle $\Bundle B$}%
   {Cartesian power A of bundle B}%
   {B}{0}%
\Symb%
   {Cartesian power $A$ of set $B$}%
   {Cartesian power of set}%
   {B}{0}%
\Symb%
   {closed ball}%
   {closed ball}%
   {B}{0}%
\Symb%
   {closure of set}%
   {closure of set}%
   {B}{0}%
\Symb%
   {coproduct in category}%
   {coproduct in category}%
   {B}{0}%
\Symb%
   {basis manifold of central affine space}%
   {FCAn}%
   {B}{0}%
\Symb%
   {basis manifold of Euclid space}%
   {FEn}%
   {B}{0}%
\Symb%
   {lattice of subrepresentations}%
   {lattice of subrepresentations}%
   {B}{0}%
\Symb%
   {lattice of towers of subrepresentations of tower of representations $\Vector f$}%
   {lattice of subrepresentations, tower of representations}%
   {B}{0}%
\Symb%
   {open ball}%
   {open ball}%
   {B}{0}%
\Symb%
   {product in category}%
   {product in category}%
   {B}{0}%
\Symb%
   {right fraction}%
   {right fraction}%
   {B}{0}%
\Symb%
   {tensor power of representation}%
   {tensor power of representation}%
   {B}{0}%

\SetIndexSpace
\Symb%
   {$\sigma$\Hyph algebra of sets measurable with respect to measure $\mu$}%
   {algebra of sets measurable with respect to measure}%
   {C}{0}%
\Symb%
   {central affine space}%
   {CAn}%
   {C}{0}%
\Symb%
   {central affine space}%
   {central affine space}%
   {C}{0}%
\Symb%
   {continuity class}%
   {class Cn}%
   {C}{0}%
\Symb%
   {$j$th column determinant of matrix $\bfA$}%
   {column determinant}%
   {C}{0}%
\Symb%
   {cosine}%
   {cosine}%
   {C}{0}%
\Symb%
   {$\CRcirc$\Hyph product of matrices of maps}%
   {cr product of matrices of maps}%
   {C}{0}%
\Symb%
   {hyperbolic cosine}%
   {hyperbolic cosine}%
   {C}{0}%
\Symb%
   {left structural constant of Lie algebra}%
   {left structural constant of Lie algebra}%
   {C}{0}%
\Symb%
   {right structural constant of Lie algebra}%
   {right structural constant of Lie algebra}%
   {C}{0}%
\Symb%
   {set of continuous multivariable maps}%
   {set continuous multivariable maps}%
   {C}{0}%
\Symb%
   {structural constants}%
   {structural constants}%
   {C}{0}%

\SetIndexSpace
\Symb%
   {basis vector of representation of Lie group over algebra $A$}%
   {basis vector of representation of Lie group over algebra A}%
   {D}{0}%
\Symb%
   {coordinates of basis vector of representation of Lie group over algebra $A$}%
   {basis vector of representation of Lie group over algebra A, coordinates}%
   {D}{0}%
\Symb%
   {component of derivative of map $f(x)$}%
   {component of derivative}%
   {D}{0}%
\Symb%
   {component of derivative of second order of map $f(x)$}%
   {component of derivative of Second Order}%
   {D}{0}%
\Symb%
   {component of the G\^ateaux derivative of map $f(x)$}%
   {component of Gateaux derivative}%
   {D}{0}%
\Symb%
   {component of the G\^ateaux derivative of map $f(x)$}%
   {component of Gateaux derivative of map, D vector space, short}%
   {D}{0}%
\Symb%
   {component of the G\^ateaux derivative of second order of map $f(x)$}%
   {component of Gateaux derivative of Second Order}%
   {D}{0}%
\Symb%
   {component of the G\^ateaux derivative of second order of map $f(x)$}%
   {component of Gateaux derivative of Second Order, D vector space}%
   {D}{0}%
\Symb%
   {component of the G\^ateaux derivative of map $f(x)$}%
   {component of Gateaux derivative, vector space}%
   {D}{0}%
\Symb%
   {conjugation in algebra}%
   {conjugation in algebra}%
   {D}{0}%
\Symb%
   {conjugation in ring}%
   {conjugation in ring}%
   {D}{0}%
\Symb%
   {coordinate \rcd vector space}%
   {coordinate rcd vector space}%
   {D}{0}%
\Symb%
   {coordinate reference frame}%
   {coordinate reference frame, extensive definition}%
   {D}{0}%
\Symb%
   {coordinate vector bundle}%
   {coordinate vector bundle}%
   {D}{0}%
\Symb%
   {derivative of map $f$}%
   {derivative of map}%
   {D}{0}%
\Symb%
   {derivative of map $f$}%
   {derivative of map inline}%
   {D}{0}%
\Symb%
   {derivative of order $n$}%
   {derivative of Order n}%
   {D}{0}%
\Symb%
   {derivative of order $n$}%
   {derivative of Order n inline}%
   {D}{0}%
\Symb%
   {derivative of second order}%
   {derivative of Second Order}%
   {D}{0}%
\Symb%
   {derivative of second order}%
   {derivative of Second Order inline}%
   {D}{0}%
\Symb%
   {diagonal in bundle $\Bundle A$}%
   {diagonal in bundle, 1}%
   {D}{0}%
\Symb%
   {differential of independent variable}%
   {differential of independent variable}%
   {D}{0}%
\Symb%
   {differential of map $f$}%
   {differential of map}%
   {D}{0}%
\Symb%
   {direct product of division rings $D_1$, ..., $D_n$}%
   {direct product of division rings, 1 n}%
   {D}{0}%
\Symb%
   {double determinant of matrix $\bfA$}%
   {double determinant}%
   {D}{0}%
\Symb%
   {exterior differential}%
   {exterior differential}%
   {D}{0}%
\Symb%
   {the Fr\'echet \Ds derivative of map $f$ of division ring}%
   {Frechet Dstar derivative of map, division ring}%
   {D}{0}%
\Symb%
   {the G\^ateaux \dcr derivative of map $f$ of $D$\Hyph vector space $V$ to $D$\Hyph vector space $W$}%
   {Gateaux dcr derivative of map, D vector space}%
   {D}{0}%
\Symb%
   {the G\^ateaux derivative of map $f$}%
   {Gateaux derivative of map}%
   {D}{0}%
\Symb%
   {the G\^ateaux derivative of map $f$}%
   {Gateaux derivative of map, fraction}%
   {D}{0}%
\Symb%
   {the G\^ateaux derivative of order $n$}%
   {Gateaux derivative of Order n}%
   {D}{0}%
\Symb%
   {the G\^ateaux derivative of order $n$ of map $f$ of division ring}%
   {Gateaux derivative of Order n, division ring}%
   {D}{0}%
\Symb%
   {the G\^ateaux derivative of order $n$ of map $f$ of algebra}%
   {Gateaux derivative of Order n, fraction, algebra}%
   {D}{0}%
\Symb%
   {the G\^ateaux derivative of order $n$ of map $f$ of division ring}%
   {Gateaux derivative of Order n, fraction, division ring}%
   {D}{0}%
\Symb%
   {the G\^ateaux derivative of second order}%
   {Gateaux derivative of Second Order}%
   {D}{0}%
\Symb%
   {the G\^ateaux derivative of second order of mapping $f$ of algebra}%
   {Gateaux derivative of Second Order, fraction, algebra}%
   {D}{0}%
\Symb%
   {the G\^ateaux derivative of second order of map $f$ of division ring}%
   {Gateaux derivative of Second Order, fraction, division ring}%
   {D}{0}%
\Symb%
   {the G\^ateaux differential of map $f$}%
   {Gateaux differential of map, vector}%
   {D}{0}%
\Symb%
   {the G\^ateaux \Ds derivative of map $f$ of division ring $D$}%
   {Gateaux Dstar derivative of map, division ring}%
   {D}{0}%
\Symb%
   {the G\^ateaux Jacobian of map of $D$\Hyph vector space}%
   {Gateaux Jacobian of map, D vector space}%
   {D}{0}%
\Symb%
   {the G\^ateaux partial \dcr derivative of map $f^{\gi b}$ with respect to variable $v^{\gi a}$}%
   {Gateaux partial dcr derivative of map, 1, D vector space}%
   {D}{0}%
\Symb%
   {the G\^ateaux partial \dcr derivative of map $f^{\gi b}$ with respect to variable $v^{\gi a}$}%
   {Gateaux partial dcr derivative of map, 2, D vector space}%
   {D}{0}%
\Symb%
   {the G\^ateaux partial \dcr derivative of map $f^{\gi b}$ with respect to variable $v^{\gi a}$}%
   {Gateaux partial dcr derivative of map, 3, D vector space}%
   {D}{0}%
\Symb%
   {the G\^ateaux partial derivative}%
   {Gateaux partial derivative}%
   {D}{0}%
\Symb%
   {the G\^ateaux mixed partial derivative}%
   {Gateaux partial derivative of Second Order}%
   {D}{0}%
\Symb%
   {the G\^ateaux partial \rcd derivative of map $f^{\gi b}$ with respect to variable $x^{\gi a}$}%
   {Gateaux partial rcd derivative of map, 1, D vector space}%
   {D}{0}%
\Symb%
   {the G\^ateaux partial \rcd derivative of map $f^{\gi b}$ with respect to variable $x^{\gi a}$}%
   {Gateaux partial rcd derivative of map, 2, D vector space}%
   {D}{0}%
\Symb%
   {the G\^ateaux partial \rcd derivative of map $f^{\gi b}$ with respect to variable $x^{\gi a}$}%
   {Gateaux partial rcd derivative of map, 3, D vector space}%
   {D}{0}%
\Symb%
   {the G\^ateaux \rcd derivative of map $f$ of $D$\hyph vector space $V$ to $D$\hyph vector space $W$}%
   {Gateaux rcd derivative of map, D vector space}%
   {D}{0}%
\Symb%
   {the G\^ateaux \sD derivative of map $f$ of division ring $D$}%
   {Gateaux starD derivative of map, division ring}%
   {D}{0}%
\Symb%
   {Jacobi matrix of map}%
   {Jacobi matrix of map}%
   {D}{0}%
\Symb%
   {matrices vector space}%
   {matrices vector space}%
   {D}{0}%
\Symb%
   {Cartan derivative}%
   {overbrace D}%
   {D}{0}%
\Symb%
   {derivative}%
   {overline D}%
   {D}{0}%
\Symb%
   {partial derivative}%
   {partial derivative}%
   {D}{0}%
\Symb%
   {partial derivative of second order}%
   {partial derivative of second order}%
   {D}{0}%
\Symb%
   {derivative $e_{(k)}$}%
   {partial(k)}%
   {D}{0}%
\Symb%
   {product of map over scalar}%
   {product of map over scalar}%
   {D}{0}%
\Symb%
   {set of vectors generated by vector $v$}%
   {set of vectors generated by vector}%
   {D}{0}%
\Symb%
   {speed of deviation}%
   {speed of deviation}%
   {D}{0}%
\Symb%
   {standard component of derivative}%
   {standard component of derivative}%
   {D}{0}%
\Symb%
   {standard component of the G\^ateaux derivative}%
   {standard component of Gateaux derivative}%
   {D}{0}%
\Symb%
   {vector space type}%
   {vector space type}%
   {D}{0}%

\SetIndexSpace
\Symb%
   {affine basis}%
   {Affine Basis}%
   {E}{0}%
\Symb%
   {basis}%
   {Basis}%
   {E}{0}%
\Symb%
   {basis for module}%
   {basis for module}%
   {E}{0}%
\Symb%
   {basis in vector space $\Vector V$}%
   {basis in V}%
   {E}{0}%
\Symb%
   {basis of $D$\Hyph module $\mathcal L(D;A_1;A_2)$}%
   {basis L(A1,A2)}%
   {E}{0}%
\Symb%
   {basis manifold}%
   {basis manifold}%
   {E}{0}%
\Symb%
   {basis for vector space}%
   {basis of vector space}%
   {E}{0}%
\Symb%
   {basis for \crd vector space}%
   {basis, crd vector space}%
   {E}{0}%
\Symb%
   {basis for \dcr vector space}%
   {basis, dcr vector space}%
   {E}{0}%
\Symb%
   {basis for \drc vector space}%
   {basis, drc vector space}%
   {E}{0}%
\Symb%
   {basis for \rcd vector space}%
   {basis, rcd vector space}%
   {E}{0}%
\Symb%
   {basis for vector bundle}%
   {basis, vector bundle}%
   {E}{0}%
\Symb%
   {basis of $(n)$\hyph vector space}%
   {basis,n vector space}%
   {E}{0}%
\Symb%
   {Cartesian power of total spaces}%
   {Cartesian power of total spaces}%
   {E}{0}%
\Symb%
   {Cartesian product of total spaces}%
   {Cartesian product of total spaces, definition 1}%
   {E}{0}%
\Symb%
   {central affine basis}%
   {Central Affine Basis}%
   {E}{0}%
\Symb%
   {\CR exponent}%
   {CR exponent}%
   {E}{0}%
\Symb%
   {form of reference frame}%
   {dual forms, reference frame}%
   {E}{0}%
\Symb%
   {Euclid space}%
   {Euclid space}%
   {E}{0}%
\Symb%
   {Euclid space}%
   {Euclid space, division ring}%
   {E}{0}%
\Symb%
   {exponent}%
   {exponent}%
   {E}{0}%
\Symb%
   {Hamel basis}%
   {Hamel basis}%
   {E}{0}%
\Symb%
   {identical transformation of bundle}%
   {identical transformation of bundle}%
   {E}{0}%
\Symb%
   {map of conjugation}%
   {map of conjugation}%
   {E}{0}%
\Symb%
   {linear automorphism of quaternioin algebra}%
   {mapping E, quaternion}%
   {E}{0}%
\Symb%
   {linear automorphism of quaternioin algebra}%
   {mapping E_1, quaternion}%
   {E}{0}%
\Symb%
   {linear automorphism of quaternioin algebra}%
   {mapping E_2, quaternion}%
   {E}{0}%
\Symb%
   {Jacobian matrix of maps of conjugation}%
   {maps of conjugation, Jacobian matrix}%
   {E}{0}%
\Symb%
   {orthonornal basis}%
   {Orthonornal Basis}%
   {E}{0}%
\Symb%
   {image of basis $\Basis e$ under passive transformation $S$}%
   {passive transformation of basis}%
   {E}{0}%
\Symb%
   {pseudo Euclid space}%
   {pseudo Euclid space}%
   {E}{0}%
\Symb%
   {pseudo Euclid space}%
   {pseudo Euclid space, division ring}%
   {E}{0}%
\Symb%
   {quasiexponent}%
   {quasiexponent}%
   {E}{0}%
\Symb%
   {quaternion algebra over the field $F$}%
   {quaternion algebra over the field}%
   {E}{0}%
\Symb%
   {quaternion division algebra over the field}%
   {quaternion division algebra over the fieldL}%
   {E}{0}%
\Symb%
   {\RC exponent}%
   {RC exponent}%
   {E}{0}%
\Symb%
   {reduced Cartesian product of total spaces}%
   {reduced Cartesian product of total spaces, definition 1}%
   {E}{0}%
\Symb%
   {Schauder basis}%
   {Schauder basis}%
   {E}{0}%
\Symb%
   {set of \CR eigenvalues}%
   {set of cr-eigenvalues}%
   {E}{0}%
\Symb%
   {set of endomorphisms}%
   {set of endomorphisms}%
   {E}{0}%
\Symb%
   {set of \RC eigenvalues}%
   {set of rc-eigenvalues}%
   {E}{0}%
\Symb%
   {set of nonsingular \sT transformations of bundle $\Bundle E$}%
   {set of starT nonsingular transformations of bundle}%
   {E}{0}%
\Symb%
   {set of transformations of universal algebra}%
   {set of transformations}%
   {E}{0}%
\Symb%
   {set of nonsingular \Ts transformations of bundle $\Bundle E$}%
   {set of Tstar nonsingular transformations of bundle}%
   {E}{0}%
\Symb%
   {standard coordinates of basis}%
   {standard coordinates of basis}%
   {E}{0}%
\Symb%
   {standard coordinates of reference frame}%
   {standard coordinates of reference frame}%
   {E}{0}%
\Symb%
   {vector field of reference frame}%
   {vector field of reference frame}%
   {E}{0}%
\Symb%
   {vector of basis}%
   {vector of basis}%
   {E}{0}%

\SetIndexSpace
\Symb%
   {alternation of polylinear map}%
   {alternation of polylinear map}%
   {F}{0}%
\Symb%
   {component of linear map $f$ of division ring}%
   {component of linear map, division ring}%
   {F}{0}%
\Symb%
   {component of polylinear map}%
   {component of polylinear map}%
   {F}{0}%
\Symb%
   {conjugation transformation}%
   {conjugation transformation}%
   {F}{0}%
\Symb%
   {exterior product}%
   {exterior product}%
   {F}{0}%
\Symb%
   {fibered morphism from bundle $\Bundle A$ into $\Bundle B$}%
   {fibered morphism from A into B}%
   {F}{0}%
\Symb%
   {filter $\mathfrak{F}$ converges to set $A$}%
   {filter converges}%
   {F}{0}%
\Symb%
   {homomorphism of fibered universal algebras}%
   {homomorphism of fibered universal algebras}%
   {F}{0}%
\Symb%
   {inverse fibered correspondence}%
   {inverse fibered correspondence, 1}%
   {F}{0}%
\Symb%
   {inverse reduced fibered correspondence}%
   {inverse reduced fibered correspondence, 1}%
   {F}{0}%
\Symb%
   {map to Cartesian product}%
   {map to Cartesian product}%
   {F}{0}%
\Symb%
   {norm of functional}%
   {norm of functional}%
   {F}{0}%
\Symb%
   {norm of map}%
   {norm of map}%
   {F}{0}%
\Symb%
   {norm of polylinear map}%
   {norm of polymap}%
   {F}{0}%
\Symb%
   {norm of representation}%
   {norm of representation}%
   {F}{0}%
\Symb%
   {orbit of representation of the group}%
   {orbit of representation}%
   {F}{0}%
\Symb%
   {orthonormal basis}%
   {Orthonormal Basis, division ring}%
   {F}{0}%
\Symb%
   {quaternion algebra  over field ${\rm {\mathbb{F}}}$}%
   {quaternion algebra F a b}%
   {F}{0}%
\Symb%
   {reference frame}%
   {reference frame}%
   {F}{0}%
\Symb%
   {reference frame, extensive definition}%
   {reference frame, extensive definition}%
   {F}{0}%
\Symb%
   {standard component of biadditive map $f$ over field $F$}%
   {standard component of biadditive map, division ring}%
   {F}{0}%
\Symb%
   {standard component of linear map}%
   {standard component of linear map, G}%
   {F}{0}%
\Symb%
   {standard component of polylinear map}%
   {standard component of polylinear map}%
   {F}{0}%
\Symb%
   {standard component of quadratic map $f$ over field $F$}%
   {standard component of quadratic map, division ring}%
   {F}{0}%
\Symb%
   {standard component of tensor}%
   {standard component of tensor}%
   {F}{0}%
\Symb%
   {sum of maps}%
   {sum of maps}%
   {F}{0}%
\Symb%
   {symmetrization of polylinear map}%
   {symmetrization of polylinear map}%
   {F}{0}%

\SetIndexSpace
\Symb%
   {affine transformation group}%
   {affine transformation group}%
   {G}{0}%
\Symb%
   {affine transformation group}%
   {affine transformation group}%
   {G}{0}%
\Symb%
   {Cartesian product of groups $G_1$, ..., $G_n$}%
   {Cartesian product of groups, 1 n}%
   {G}{0}%
\Symb%
   {\CR matrix group}%
   {cr-matrix group}%
   {G}{0}%
\Symb%
   {fibered little group of section $h$}%
   {fibered little group}%
   {G}{0}%
\Symb%
   {fibered stability group of section $h$}%
   {fibered stability group}%
   {G}{0}%
\Symb%
   {group of automorphisms of representation $f$}%
   {group of automorphisms of representation}%
   {G}{0}%
\Symb%
   {group of homomorphisms of vector space $\Vector V$}%
   {GV}%
   {G}{0}%
\Symb%
   {indefinite integral}%
   {indefinite integral}%
   {G}{0}%
\Symb%
   {left defined Lie algebra of Lie group}%
   {left defined Lie algebra of Lie group}%
   {G}{0}%
\Symb%
   {Lie algebra of Lie group}%
   {Lie algebra of Lie group}%
   {G}{0}%
\Symb%
   {linear transformation group}%
   {linear transformation group}%
   {G}{0}%
\Symb%
   {little group}%
   {little group}%
   {G}{0}%
\Symb%
   {orbit of effective Ts representation of group}%
   {orbit of effective starT representation of fibered group}%
   {G}{0}%
\Symb%
   {orbit of effective \Ts representation of fibered group}%
   {orbit of effective Tstar representation of fibered group}%
   {G}{0}%
\Symb%
   {\RC matrix group}%
   {rc-matrix group}%
   {G}{0}%
\Symb%
   {right defined Lie algebra of Lie group}%
   {right defined Lie algebra}%
   {G}{0}%
\Symb%
   {stability group}%
   {stability group}%
   {G}{0}%

\SetIndexSpace
\Symb%
   {Hadamard inverse of matrix}%
   {Hadamard inverse of matrix}%
   {H}{0}%
\Symb%
   {horizontal component of vector}%
   {horizontal component of vector}%
   {H}{0}%
\Symb%
   {horizontal subspace}%
   {horizontal subspace}%
   {H}{0}%
\Symb%
   {quaternion algebra}%
   {quaternion algebra}%
   {H}{0}%
\Symb%
   {quaternion algebra}%
   {quaternion algebra H a b}%
   {H}{0}%
\Symb%
   {set of homomorphisms}%
   {set of homomorphisms}%
   {H}{0}%

\SetIndexSpace
\Symb%
   {image of map}%
   {image of map}%
   {I}{0}%
\Symb%
   {infinitesimal generator of representation}%
   {infinitesimal generator i of representation}%
   {I}{0}%
\Symb%
   {infinitesimal generator of representation}%
   {infinitesimal generator of representation}%
   {I}{0}%
\Symb%
   {Lie group infinitesimal generators}%
   {Lie group infinitesimal generators}%
   {I}{0}%
\Symb%
   {map of conjugation}%
   {map of conjugation}%
   {I}{0}%
\Symb%
   {Jacobian matrix of maps of conjugation}%
   {maps of conjugation, Jacobian matrix}%
   {I}{0}%
\Symb%
   {vector module of algebra $A$}%
   {vector module of algebra}%
   {I}{0}%
\Symb%
   {vector module of ring $D$}%
   {vector module of ring}%
   {I}{0}%
\Symb%
   {vector of element $d$ of algebra}%
   {vector of algebra}%
   {I}{0}%
\Symb%
   {vector of element $d$ of ring}%
   {vector of ring}%
   {I}{0}%

\SetIndexSpace
\Symb%
   {closure operator of representation $f$}%
   {closure operator, representation}%
   {J}{0}%
\Symb%
   {closure operator of tower of representations $\Vector f$}%
   {closure operator, tower of representations}%
   {J}{0}%
\Symb%
   {Jacobian matrix of right shift}%
   {Ea, quaternion, Jacobian matrix}%
   {J}{0}%
\Symb%
   {map of conjugation}%
   {map of conjugation}%
   {J}{0}%
\Symb%
   {Jacobian matrix of maps of conjugation}%
   {maps of conjugation, Jacobian matrix}%
   {J}{0}%
\Symb%
   {subrepresentation generated by generating set $X$}%
   {subrepresentation generated by set}%
   {J}{0}%
\Symb%
   {tower of subrepresentations of tower of representations $\Vector f$ generated by tuple of sets $\VX X$}%
   {subrepresentation generated by tuple of sets}%
   {J}{0}%

\SetIndexSpace
\Symb%
   {kernel of group\Hyph homomorphism}%
   {ker group-homomorphism}%
   {K}{0}%
\Symb%
   {kernel of homomorphism}%
   {kernel of homomorphism}%
   {K}{0}%
\Symb%
   {kernel of linear map}%
   {kernel of linear map}%
   {K}{0}%
\Symb%
   {kernel of map}%
   {kernel of map}%
   {K}{0}%
\Symb%
   {map of conjugation}%
   {map of conjugation}%
   {K}{0}%
\Symb%
   {Jacobian matrix of maps of conjugation}%
   {maps of conjugation, Jacobian matrix}%
   {K}{0}%

\SetIndexSpace
\Symb%
   {Cartesian power of systems of subsets}%
   {Cartesian power of systems of subsets}%
   {L}{0}%
\Symb%
   {Cartesian product of systems of subsets}%
   {Cartesian product of systems of subsets}%
   {L}{0}%
\Symb%
   {left $ij$th cofactor of entry of matrix}%
   {left cofactor, matrix}%
   {L}{0}%
\Symb%
   {left double $ij$th cofactor of entry of matrix}%
   {left double cofactor}%
   {L}{0}%
\Symb%
   {left shift}%
   {left shift}%
   {L}{0}%
\Symb%
   {left vector space of algebra $A$}%
   {left vector space of algebra A}%
   {L}{0}%
\Symb%
   {left vector space of maps}%
   {left vector space of maps B->A}%
   {L}{0}%
\Symb%
   {Lie derivative of connection}%
   {Lie derivative of connection}%
   {L}{0}%
\Symb%
   {Lie derivative of metric}%
   {Lie derivative of metric}%
   {L}{0}%
\Symb%
   {limit of correspondence $\Phi$ with respect to the filter $\mathfrak{F}$}%
   {limit of correspondence with respect to the filter}%
   {L}{0}%
\Symb%
   {limit of sequence}%
   {limit of sequence}%
   {L}{0}%
\Symb%
   {linear combination}%
   {linear combination}%
   {L}{0}%
\Symb%
   {module of skew symmetric polylinear maps}%
   {module of skew symmetric polylinear maps}%
   {L}{0}%
\Symb%
   {passive transformation}%
   {passive transformation}%
   {L}{0}%
\Symb%
   {$D$\Hyph module of continuous linear mappings of normed $D$\Hyph module $A_1$ into normed $D$\Hyph module $A_2$}%
   {set continuous linear mappings, module}%
   {L}{0}%
\Symb%
   {set of continuous linear maps}%
   {set continuous linear maps, vector}%
   {L}{0}%
\Symb%
   {set of continuous polylinear maps}%
   {set continuous polylinear maps}%
   {L}{0}%
\Symb%
   {set of linear maps}%
   {set linear maps}%
   {L}{0}%
\Symb%
   { inverse matrices}%
   {set of inverse matrices}%
   {L}{0}%
\Symb%
   {set of left-side nonsingular transformations of universal algebra $M$}%
   {set of left-side nonsingular transformations}%
   {L}{0}%
\Symb%
   {set of polylinear maps}%
   {set polylinear maps}%
   {L}{0}%
\Symb%
   {set of $n$\hyph linear maps}%
   {set polylinear maps An}%
   {L}{0}%
\Symb%
   {set of polylinear maps}%
   {set polylinear maps, D vector space}%
   {L}{0}%
\Symb%
   {set of polylinear maps of algebras $A_1$, ..., $A_n$ into algebra $A$}%
   {set polylinear maps, finite dimensional algebra}%
   {L}{0}%

\SetIndexSpace
\Symb%
   {set of left-side transformations of the universal algebra $M$}%
   {set of left-side transformations}%
   {M}{0}%
\Symb%
   {set of maps to $\Omega$\Hyph group $A$}%
   {set of maps to Omega group}%
   {M}{0}%
\Symb%
   {set of right-side transformations of universal algebra $M$}%
   {set of right-side transformations}%
   {M}{0}%
\Symb%
   {space of orbits of \Ts{G}representation}%
   {space of orbits of G* representation}%
   {M}{0}%

\SetIndexSpace
\Symb%
   {norm of quaternion $x$}%
   {norm, quaternion algebra}%
   {N}{0}%
\Symb%
   {nucleus of $D$\Hyph algebra $A$}%
   {nucleus of algebra}%
   {N}{0}%

\SetIndexSpace
\Symb%
   {coordinates of geometric object}%
   {coordinates of geometric object}%
   {O}{0}%
\Symb%
   {geometric object}%
   {geometric object}%
   {O}{0}%
\Symb%
   {octonion algebra}%
   {octonion algebra}%
   {O}{0}%
\Symb%
   {orbit of representation of fibered group $\Bundle G$}%
   {orbit of representation of fibered group}%
   {O}{0}%
\Symb%
   {tensor product}%
   {tensor product}%
   {O}{0}%

\SetIndexSpace
\Symb%
   {bundle}%
   {bundle}%
   {P}{0}%
\Symb%
   {bundle of level $2$}%
   {bundle of level 2}%
   {P}{0}%
\Symb%
   {bundle of level $n$}%
   {bundle of level n}%
   {P}{0}%
\Symb%
   {Cartesian power $n$ of bundle $\bundle{}{p}{E}{}$}%
   {Cartesian power of bundle}%
   {P}{0}%
\Symb%
   {Cartesian product of bundles}%
   {Cartesian product of bundles, definition 1}%
   {P}{0}%
\Symb%
   {passive representation in basis manifold}%
   {passive representation in basis manifold}%
   {P}{0}%
\Symb%
   {reduced Cartesian product of bundles}%
   {reduced Cartesian product of bundles, definition 1}%
   {P}{0}%
\Symb%
   {set of nonsingular \sT transformations of bundle $\bundle{}pE{}$}%
   {set of starT nonsingular transformations of bundle, projection}%
   {P}{0}%
\Symb%
   {set of nonsingular \Ts transformations of bundle $\bundle{}pE{}$}%
   {set of Tstar nonsingular transformations of bundle, projection}%
   {P}{0}%

\SetIndexSpace
\Symb%
   {active transformation}%
   {active transformation}%
   {R}{0}%
\Symb%
   {Cartan curvature}%
   {Cartan curvature}%
   {R}{0}%
\Symb%
   {\CR rank of matrix}%
   {cr-rank of matrix}%
   {R}{0}%
\Symb%
   {diagonal in bundle  $\bundle{}pA{}$}%
   {diagonal in bundle, 2}%
   {R}{0}%
\Symb%
   {diagonal in bundle $\Bundle A$}%
   {diagonal in reduced bundle, 2}%
   {R}{0}%
\Symb%
   {image of $m$ under endomorphism $R$ of effective representation}%
   {endomorphism image, effective representation}%
   {R}{0}%
\Symb%
   {image of tuple $\VX a$ under endomorphism $\VX r$ of tower of effective representations}%
   {endomorphism image, tower of effective representations}%
   {R}{0}%
\Symb%
   {curvature}%
   {GLn curvature_overline}%
   {R}{0}%
\Symb%
   {product of rings of sets}%
   {product of rings of sets}%
   {R}{0}%
\Symb%
   {$\RCcirc$\Hyph product of matrices of maps}%
   {rc product of matrices of maps}%
   {R}{0}%
\Symb%
   {\RC rank of matrix}%
   {rc-rank of matrix}%
   {R}{0}%
\Symb%
   {right $ij$th cofactor of entry of matrix}%
   {right cofactor, matrix}%
   {R}{0}%
\Symb%
   {right double $ij$th cofactor of entry of matrix}%
   {right double cofactor}%
   {R}{0}%
\Symb%
   {right shift}%
   {right shift}%
   {R}{0}%
\Symb%
   {right vector space of algebra $A$}%
   {right vector space of algebra A}%
   {R}{0}%
\Symb%
   {right vector space of maps}%
   {right vector space of maps B->A}%
   {R}{0}%
\Symb%
   {$i$th row determinant of matrix $\bfA$}%
   {row determinant}%
   {R}{0}%
\Symb%
   {scalar algebra of algebra $A$}%
   {scalar algebra of algebra}%
   {R}{0}%
\Symb%
   {scalar algebra of ring $D$}%
   {scalar algebra of ring}%
   {R}{0}%
\Symb%
   {scalar of element $d$ of algebra}%
   {scalar of algebra}%
   {R}{0}%
\Symb%
   {scalar of element $d$ of ring}%
   {scalar of ring}%
   {R}{0}%
\Symb%
   { inverse matrices}%
   {set of inverse matrices}%
   {R}{0}%
\Symb%
   {set of right-side nonsingular transformations of universal algebra $M$}%
   {set of right-side nonsingular transformations}%
   {R}{0}%
\Symb%
   {spherical coordinates}%
   {spherical coordinates}%
   {R}{0}%

\SetIndexSpace
\Symb%
   {composition of fibered correspondences}%
   {composition of fibered correspondences}%
   {S}{0}%
\Symb%
   {hyperbolic sine}%
   {hyperbolic sine}%
   {S}{0}%
\Symb%
   {inverse fibered correspondence}%
   {inverse fibered correspondence, 2}%
   {S}{0}%
\Symb%
   {inverse reduced fibered correspondence}%
   {inverse reduced fibered correspondence, 2}%
   {S}{0}%
\Symb%
   {Lebesgue integral}%
   {Lebesgue integral}%
   {S}{0}%
\Symb%
   {linear span in vector space}%
   {linear span, vector space}%
   {S}{0}%
\Symb%
   {image of basis $\VX  X$ under passive transformation $\VX s$}%
   {passive transformation of basis, tower of representations}%
   {S}{0}%
\Symb%
   {set of permutations}%
   {set of permutations}%
   {S}{0}%
\Symb%
   {set of transpositions}%
   {set of transpositions}%
   {S}{0}%
\Symb%
   {sine}%
   {sine}%
   {S}{0}%
\Symb%
   {symmetric group}%
   {symmetric group}%
   {S}{0}%

\SetIndexSpace
\Symb%
   {category of left-side representations}%
   {category of left-side representations}%
   {T}{0}%
\Symb%
   {tangent plane to Lie group $G$}%
   {tangent plane to Lie group}%
   {T}{0}%
\Symb%
   {trace of quaternion $x$}%
   {trace, quaternion algebra}%
   {T}{0}%

\SetIndexSpace
\Symb%
   {affine space}%
   {affine space}%
   {V}{0}%
\Symb%
   {conjugated affine space}%
   {conjugated affine space}%
   {V}{0}%
\Symb%
   {conjugated vector space}%
   {conjugated vector space}%
   {V}{0}%
\Symb%
   {coordinate vector space}%
   {coordinate vector space}%
   {V}{0}%
\Symb%
   {coordinates in vector space}%
   {coordinates in vector space}%
   {V}{0}%
\Symb%
   {direct product of $\RCstar D_i$\hyph vector spaces $\Vector V_1$, ..., $\Vector V_n$}%
   {direct product, rcd vector space, 1 n}%
   {V}{0}%
\Symb%
   {dual space of \rcd vector space $\Vector V$}%
   {dual space of rcd vector space}%
   {V}{0}%
\Symb%
   {hermitian conjugated vector}%
   {hermitian conjugated vector}%
   {V}{0}%
\Symb%
   {linear composition of vectors}%
   {linear composition of vectors}%
   {V}{0}%
\Symb%
   {set of vectors generated by vector $v$}%
   {set of vectors generated by vector}%
   {V}{0}%
\Symb%
   {vector space}%
   {V}%
   {V}{0}%
\Symb%
   {vector space of initial values of system of differential equations}%
   {vector space of initial values}%
   {V}{0}%
\Symb%
   {vector space of solutions of system of differential equations}%
   {vector space of solutions}%
   {V}{0}%
\Symb%
   {vertical component of vector}%
   {vertical component of vector}%
   {V}{0}%
\Symb%
   {vertical subspace}%
   {vertical subspace}%
   {V}{0}%

\SetIndexSpace
\Symb%
   {set of coordinates of representation $J(f,X)$}%
   {coordinate set of representation}%
   {W}{0}%
\Symb%
   {set of tuples of coordinates of tower of representations $\Vector J(\Vector f,\VX X)$}%
   {coordinate set of tower of representations}%
   {W}{0}%
\Symb%
   {coordinates of basis $X'$ relative to basis $X$ of representation}%
   {coordinates of basis relative to basis, representation}%
   {W}{0}%
\Symb%
   {coordinates of element $m$ relative to set $X$}%
   {coordinates of element relative to set, representation}%
   {W}{0}%
\Symb%
   {tuple of coordinates of element $\Vector a*$ relative to tuple of sets $\VX X$}%
   {coordinates of element, tower of representations}%
   {W}{0}%
\Symb%
   {coordinates of element $m$ of representation $f$ relative to set $X$}%
   {coordinates relative to set}%
   {W}{0}%
\Symb%
   {geometric object}%
   {geometric object}%
   {W}{0}%
\Symb%
   {set of tuples of $\Omega$\Hyph words}%
   {set of tuples of Omega words}%
   {W}{0}%
\Symb%
   {set of coordinates of set $B\subset J(f,X)$}%
   {subset of coordinates of representation}%
   {W}{0}%
\Symb%
   {coordinates of tuple of sets $\VX B$ relative to tuple of sets $\VX X$}%
   {subset of coordinates of tower of representations}%
   {W}{0}%
\Symb%
   {coordinates of set $B_k$ relative to tuple of sets $\VX X$}%
   {subset of coordinates of tower of representations, k}%
   {W}{0}%
\Symb%
   {set of $\Omega_2$\Hyph words representing set $B\subset J(f,X)$}%
   {subset of words of representation}%
   {W}{0}%
\Symb%
   {superposition of coordinates}%
   {superposition of coordinates}%
   {W}{0}%
\Symb%
   {superposition of coordinates of the tower of representations $\Vector f$ and the element $\VX a$}%
   {superposition of coordinates, tower of representations}%
   {W}{0}%
\Symb%
   {tuple of $\Omega$\Hyph words}%
   {tuple of Omega words}%
   {W}{0}%
\Symb%
   {$\Omega_2$\Hyph word representing element $m\in J(f,X)$}%
   {word of element relative to generating set, representation}%
   {W}{0}%
\Symb%
   {set of $\Omega_2$\Hyph words of representation $J(f,X)$}%
   {word set of representation}%
   {W}{0}%
\Symb%
   {set of tuples of $\VX{\Omega}$\Hyph words of tower of representations $\Vector J(\Vector f,\VX X)$}%
   {word set of tower of representations}%
   {W}{0}%
\Symb%
   {tuple of words of element $\Vector a*$ relative to tuple of sets $\VX X$}%
   {words of element, tower of representations}%
   {W}{0}%

\SetIndexSpace
\Symb%
   {conjugate of quaternion $x$}%
   {conjugate of quaternion}%
   {X}{0}%
\Symb%
   {local basis of affine space}%
   {local basis of affine space}%
   {X}{0}%
\Symb%
   {Wronskian matrix}%
   {Wronskian matrix}%
   {X}{0}%
\Symb%
   {anholonomic coordinate}%
   {x(k)}%
   {X}{0}%

\SetIndexSpace
\Symb%
   {center of $A$\Hyph number}%
   {center of A number}%
   {Z}{0}%
\Symb%
   {center of $D$\Hyph algebra $A$}%
   {center of algebra}%
   {Z}{0}%
\Symb%
   {center of ring $D$}%
   {center of ring}%
   {Z}{0}%

\SetIndexSpace
\Symb%
   {deviation of trajectories}%
   {deviation of trajectories}%
   {Delta}{1}%
\Symb%
   {identical transformation}%
   {identical transformation}%
   {Delta}{1}%
\Symb%
   {image of vector $\Vector e_k\in\Basis e$ under isomorphism to coordinate vector space}%
   {image of vector e_k, coordinate vector space}%
   {Delta}{1}%
\Symb%
   {Kronecker symbol}%
   {Kronecker symbol}%
   {Delta}{1}%

\SetIndexSpace
\Symb%
   {anholonomic coordinates of connection}%
   {anholonomic coordinates of connection}%
   {Gamma}{1}%
\Symb%
   {Cartan symbol}%
   {Cartan symbol}%
   {Gamma}{1}%
\Symb%
   {connection}%
   {conection overline}%
   {Gamma}{1}%
\Symb%
   {connection coefficients in $D$\Hyph affine space}%
   {connection coefficients, D affine space}%
   {Gamma}{1}%
\Symb%
   {connection in $D$\Hyph affine manifold}%
   {connection, affine manifold}%
   {Gamma}{1}%
\Symb%
   {$D$\Hyph affine connection coefficients on manifold}%
   {D affine connection coefficients, manifold}%
   {Gamma}{1}%
\Symb%
   {holonomic coordinates of connection}%
   {holonomic coordinates of connection}%
   {Gamma}{1}%
\Symb%
   {Cartan connection}%
   {overbrace Gamma i kl}%
   {Gamma}{1}%
\Symb%
   {set of sections of bundle}%
   {set of sections of bundle}%
   {Gamma}{1}%

\SetIndexSpace
\Symb%
   {inverse operator to operator $\psi_l$}%
   {inverse operator to operator psi l}%
   {Lambda}{1}%
\Symb%
   {inverse operator to operator $\psi_r$}%
   {inverse operator to operator psi r}%
   {Lambda}{1}%

\SetIndexSpace
\Symb%
   {Cartesian product of measures}%
   {Cartesian product of measures}%
   {Mu}{1}%
\Symb%
   {power of measure}%
   {power of measure}%
   {Mu}{1}%
\Symb%
   {product of measures}%
   {product of measures}%
   {Mu}{1}%
\Symb%
   {product of measures}%
   {product of measures, otimes}%
   {Mu}{1}%

\SetIndexSpace
\Symb%
   {anholonomity object}%
   {anholonomity object}%
   {Omega}{1}%
\Symb%
   {definite integral}%
   {definite integral}%
   {Omega}{1}%
\Symb%
   {integral of differential $1$\Hyph form along path}%
   {integral of differential 1 form along path}%
   {Omega}{1}%
\Symb%
   {norm of operation}%
   {norm of operation}%
   {Omega}{1}%
\Symb%
   {operator domain}%
   {operator domain}%
   {Omega}{1}%
\Symb%
   {set of differential $p$\Hyph forms}%
   {set of differential p forms}%
   {Omega}{1}%
\Symb%
   {set of $n$\Hyph ary operations of $\Omega$\Hyph algebra}%
   {set of n-ary operations}%
   {Omega}{1}%
\Symb%
   {set of $n$\Hyph ary operators}%
   {set of n-ary operators}%
   {Omega}{1}%

\SetIndexSpace
\Symb%
   {left basic operator of Lie group over algebra $A$}%
   {L basic operator of Lie group over algebra A}%
   {Psi}{1}%
\Symb%
   {left basic operator of group Lie}%
   {Lie Basic Operator L}%
   {Psi}{1}%
\Symb%
   {left basic operator of Lie 1-parameter group}%
   {Lie Basic Operator L, 1-Parameter Group}%
   {Psi}{1}%
\Symb%
   {left basic operator of Lie 1-parameter group over algebra $A$}%
   {Lie Basic Operator L, 1-Parameter Group, algebra}%
   {Psi}{1}%
\Symb%
   {right basic operator of group Lie}%
   {Lie Basic Operator R}%
   {Psi}{1}%
\Symb%
   {right basic operator of Lie 1-parameter group}%
   {Lie Basic Operator R, 1-Parameter Group}%
   {Psi}{1}%
\Symb%
   {right basic operator of Lie 1-parameter group over algebra $A$}%
   {Lie Basic Operator R, 1-Parameter Group, algebra}%
   {Psi}{1}%
\Symb%
   {right basic operator of Lie group over algebra $A$}%
   {R basic operator of Lie group over algebra A}%
   {Psi}{1}%

\SetIndexSpace
\Symb%
   {fibered subset}%
   {fibered subset}%
   {Sigma}{1}%
\Symb%
   {parity of permutation}%
   {parity of permutation}%
   {Sigma}{1}%
\Symb%
   {subbundle}%
   {subbundle}%
   {Sigma}{1}%

\SetIndexSpace
\Symb%
   {Cartan derivative}%
   {overbrace nabla_l}%
   {Nabla}{2}%
\Symb%
   {derivative}%
   {overline nabla_l, definition 1}%
   {Nabla}{2}%

\SetIndexSpace
\Symb%
   {Lie group composition law}%
   {Lie group composition law}%
   {Phi}{1}%
\Symb%
   {restriction of correspondence $\Phi$ to set $C$}%
   {restriction of correspondence}%
   {Phi}{1}%

\SetIndexSpace
\Symb%
   {Cartesian product of bundles}%
   {Cartesian product of bundles, definition 2}%
   {Pi}{1}%
\Symb%
   {Cartesian product of groups $G_i$, $i\in I$}%
   {Cartesian product of groups}%
   {Pi}{1}%
\Symb%
   {Cartesian product of groups $G_1$, ..., $G_n$}%
   {Cartesian product of groups, i 1 n}%
   {Pi}{1}%
\Symb%
   {Cartesian product of total spaces}%
   {Cartesian product of total spaces, definition 2}%
   {Pi}{1}%
\Symb%
   {coproduct in category}%
   {coproduct in category}%
   {Pi}{1}%
\Symb%
   {direct product of division rings $D_i$, $i\in I$}%
   {direct product of division rings}%
   {Pi}{1}%
\Symb%
   {direct product of division rings $D_1$, ..., $D_n$}%
   {direct product of division rings, i 1 n}%
   {Pi}{1}%
\Symb%
   {direct product of $\RCstar D_i$\hyph vector spaces $\Vector V_i$, $i\in I$}%
   {direct product, rcd vector space}%
   {Pi}{1}%
\Symb%
   {direct product of $\RCstar D_i$\hyph vector spaces}%
   {direct product, rcd vector space, i 1 n}%
   {Pi}{1}%
\Symb%
   {product in category}%
   {product in category}%
   {Pi}{1}%
\Symb%
   {reduced Cartesian product of bundles}%
   {reduced Cartesian product of bundles, definition 2}%
   {Pi}{1}%
\Symb%
   {reduced Cartesian product of total spaces}%
   {reduced Cartesian product of total spaces, definition 2}%
   {Pi}{1}%

\CloseIndex

%% file: Prolog.English.tex
\ifx\setCACAA\Defined
\maketitle
\input{Prolog.Eq}
\input{\FilePrefix Abstract.\BookNumber.\TheLanguage}
\input{\FilePrefix Preface.\BookNumber.\TheLanguage}
\else
\Prolog
\fi

\DefText{Preliminary Definitions}
{
\ifx\PrintBook\undefined
This section
\else
This chapter
\fi
contains definitions and theorems
which are necessary for an understanding of the text of this
\ifx\PrintBook\undefined
paper.
\else
book.
\fi
So the reader may read the statements from this
\ifx\PrintBook\undefined
section
\else
chapter
\fi
in process of reading the main text of the
\ifx\PrintBook\undefined
paper.
\else
book.
\fi
}

\ePrints{4975-6381,6860-2955}%
\Items{5410-9916,9835-2163,7287-9339,309618526,CACAA.06.121,0767-8264}
\ifx\Semafor\ValueOn%
\Chapter{Preliminary Definitions}

\ShowText{Preliminary Definitions}
\fi

%% file: Preliminary.Representation.English.tex
\input{Preliminary.Representation.Eq}

\ePrints{4975-6381,6860-2955}%
\Items{5410-9916,9835-2163,7287-9339}
\ifx\Semafor\ValueOff
\Chapter{Representation of Universal Algebra}

\ShowText{Preliminary Definitions}
\fi

\ePrints{4975-6381,1506.00061,8428-0408,2207.06506,5148-4632,1908.04418,6860-2955}
\Items{1601.03259,5410-9916,9835-2163,7287-9339,5284-0163,1801.01628}
\ifx\Semafor\ValueOn

\ePrints{1908.04418,6860-2955}
\ifx\Semafor\ValueOn
\Section{Equivalence}

\begin{\DefinitionStyle}
\labelDefinition{equivalence}
Correspondence
\ShowEq{Phi in AxA}
is called
\AddIndex{equivalence}{equivalence},
if\,\footnotemark
\StartLabelItem
\begin{enumerate}
\item
correspondence $\Phi$ is
\AddIndex{reflexive}{reflexive correspondence}
\ShowEq{reflexive correspondence}
\item
correspondence $\Phi$ is
\AddIndex{symmetric}{symmetric correspondence}
\ShowEq{symmetric correspondence}
\item
correspondence $\Phi$ is
\AddIndex{transitive}{transitive correspondence}
\ShowEq{transitive correspondence}
\end{enumerate}
\end{\DefinitionStyle}
\footnotetext{\,
See also the definition on page
\citeBib{Cohn: Universal Algebra}\Hyph 14.
}

\begin{\TheoremStyle}
For the map
\ShowEq{f:A->B}fAB
the set
\ShowEq{kernel of map}
\ShowEq{def kernel of map}
is equivalence and
is called
\AddIndex{kernel of map}{kernel of map}.\,\footnotemark
\end{\TheoremStyle}
\footnotetext{\,
See also the definition on page
\citeBib{Cohn: Universal Algebra}\Hyph 16.
}
\begin{proof}
\begin{ShadedLemma}
\labelLemma{ker - reflexive correspondence}
{\it
Correspondence
\ShowEq{ker f}
is reflexive.
}
\end{ShadedLemma}

{\sc Proof.}
From the equality
\ShowEq{fa=fa}
and from the definition
\EqRef{def kernel of map},
it follows that
\DrawEq[aa]{ab in ker}{aa}
The lemma follows from the statement
\eqRef{ab in ker}{aa}
and from the definition
\RefItem{reflexive correspondence}.
\hfill\(\odot\)

\begin{ShadedLemma}
\labelLemma{ker - symmetric correspondence}
{\it
Correspondence
\ShowEq{ker f}
is symmetric.
}
\end{ShadedLemma}

{\sc Proof.}
The equality
\DrawEq[ab]{fa=fb}{ab symmetric}
follows from the statement
\DrawEq[ab]{ab in ker}{}
and from the definition
\EqRef{def kernel of map}.
The equality
\DrawEq[ba]{fa=fb}{ba}
follows from the equality
\eqRef{fa=fb}{ab symmetric}.
The statement
\DrawEq[ba]{ab in ker}{}
follows from the equality
\eqRef{fa=fb}{ba}
and from the definition
\EqRef{def kernel of map}.
Therefore, we proved the statement
\ShowEq{ab in ker -> ba in ker}
The lemma follows from the statement
\EqRef{ab in ker -> ba in ker}
and from the definition
\RefItem{symmetric correspondence}.
\hfill\(\odot\)

\begin{ShadedLemma}
\labelLemma{ker - transitive correspondence}
{\it
Correspondence
\ShowEq{ker f}
is transitive.
}
\end{ShadedLemma}

{\sc Proof.}
The equality
\DrawEq[ab]{fa=fb}{ab transitive}
follows from the statement
\DrawEq[ab]{ab in ker}{}
and from the definition
\EqRef{def kernel of map}.
The equality
\DrawEq[bc]{fa=fb}{bc}
follows from the statement
\DrawEq[bc]{ab in ker}{}
and from the definition
\EqRef{def kernel of map}.
The equality
\DrawEq[ac]{fa=fb}{ac}
follows from equalities
\eqRef{fa=fb}{ab transitive},
\eqRef{fa=fb}{bc}.
The statement
\DrawEq[ac]{ab in ker}{}
follows from the equality
\eqRef{fa=fb}{ac}
and from the definition
\EqRef{def kernel of map}.
Therefore, we proved the statement
\ShowEq{ab,ac in ker -> ac in ker}
The lemma follows from the statement
\EqRef{ab,ac in ker -> ac in ker}
and from the definition
\RefItem{symmetric correspondence}.
\hfill\(\odot\)

The theorem follows from lemmas
\ShowEq{lemmas for ker}
and from the definition
\RefDefinition{equivalence}.
\end{proof}
\else
\Section{Universal Algebra}
\fi

\ePrints{8428-0408,2207.06506}
\ifx\Semafor\ValueOn
\begin{\DefinitionStyle}
For the map
\ShowEq{f:A->B}fAB
the set
\ShowEq{image of map}
\ShowEq{show image of map}
is called
\AddIndex{image of map}{image of map}
$f$.
\end{\DefinitionStyle}
\fi

\ePrints{5410-9916,9835-2163,5284-0163,1801.01628,8428-0408,2207.06506}
\ifx\Semafor\ValueOff
\begin{\TheoremStyle}
\labelTheorem{maps and kernel equivalence}
Let $N$ be equivalence
on the set $A$.
Consider category $\mathcal A$ whose objects are
maps\,\footnotemark
\ShowEq{maps category}
We define morphism $f_1\rightarrow f_2$
to be map $h:S_1\rightarrow S_2$
making following diagram commutative
\ShowEq{maps category, diagram}
The map
\ShowEq{maps category, universal}
is universally repelling in the category $\mathcal A$.\,\footnotemark
\end{\TheoremStyle}
\ShowPrevFootnote{
The statement of the theorem
\RefTheorem{maps and kernel equivalence}
is similar to the statement on p. \citeBib{Serge Lang}-119.
}
\ShowNextFootnote{
See definition of universal object of category in definition
on p. \citeBib{Serge Lang}-57.
}
\begin{proof}
Consider diagram
\ShowEq{maps category, universal, diagram}
\ShowEq{maps category, universal, ker}
From the statement
\EqRef{maps category, universal, ker}
and the equality
\ShowEq{maps category 1}
it follows that
\ShowEq{maps category 2}
Therefore, we can uniquely define the map $h$
using the equality
\ShowEq{maps category, h}
\end{proof}
\fi

\ePrints{1908.04418,6860-2955}
\ifx\Semafor\ValueOn
\Section{Universal Algebra}
\fi

\begin{\DefinitionStyle}
\labelDefinition{Cartesian power}
For any sets\,\footnotemark
$A$, $B$,
\AddIndex{Cartesian power}{Cartesian power}
\ShowEq{Cartesian power}
is the set of maps
\ShowEq{f:A->B}fAB
\end{\DefinitionStyle}
\footnotetext{\,
I follow the definition from the example (iv) on the
page \citeBib{Cohn: Universal Algebra}\Hyph 5.
}

\begin{\DefinitionStyle}
\labelDefinition{operation on set}
For any $n\ge 0$, a map\,\footnotemark
\ShowEq{o:An->A}
is called
\AddIndex{$n$\Hyph ary operation on set}{n-ary operation on set}
$A$ or just
\AddIndex{operation on set}{operation on set}
$A$.
For any
\ShowEq{b1n in B}aA
we use either notation
\ShowEq{a1no=oa1n}
to denote image of map $\omega$.
\end{\DefinitionStyle}
\footnotetext{\,
\ePrints{MAlgebra2,5148-4632,1908.04418,6860-2955}
\ifx\Semafor\ValueOn
Definitions
\RefDefinition{operation on set},
\RefDefinition{set is closed with respect to operation}
follow
\else
Definition
\RefDefinition{operation on set}
follows
\fi
the definition in the example (vi) on the page
page \citeBib{Cohn: Universal Algebra}\Hyph 13.
}

\ePrints{9835-2163,5284-0163,1801.01628}
\ifx\Semafor\ValueOff
\begin{remark}
\labelRemark{o in AAn}
According to definitions
\RefDefinition{Cartesian power},
\RefDefinition{operation on set},
$n$\Hyph ari operation
\ShowEq{o in AAn}
\qed
\end{remark}
\fi

\begin{\DefinitionStyle}
\AddIndex{Operator domain}{operator domain}
is the set of operators\,\footnotemark
\ShowEq{operator domain}
with a map
\ShowEq{f:A->B}a{\Omega}N
If
\ShowEq{omega in Omega}{},
then $a(\omega)$ is called the
\AddIndex{arity}{arity}
of operator $\omega$. If
\ShowEq{a(o)=n}
then operator $\omega$ is called $n$\Hyph ary.
\ShowEq{set of n-ary operators}
We use notation
\ShowEq{set of n-ary operators =}
for the set of $n$\Hyph ary operators.
\end{\DefinitionStyle}
\footnotetext{\,
I follow the definition (1),
page \citeBib{Cohn: Universal Algebra}\Hyph 48.
}

\begin{\DefinitionStyle}
\labelDefinition{Omega-algebra}
Let $A$ be a set. Let $\Omega$ be an operator domain.\,\footnotemark
The family of maps
\ShowEq{O(n)->AAn}
is called $\Omega$\Hyph algebra structure on $A$.
The set $A$ with $\Omega$\Hyph algebra structure is called
\AddIndex{$\Omega$\Hyph algebra}{Omega-algebra}
\ShowEq{Omega-algebra}
or
\AddIndex{universal algebra}{universal algebra}.
\ePrints{9835-2163,0767-8264,5284-0163,1801.01628}
\ifx\Semafor\ValueOff
The set $A$ is called
\AddIndex{carrier of $\Omega$\Hyph algebra}{carrier of Omega-algebra}.
\fi
\end{\DefinitionStyle}
\footnotetext{\,
I follow the definition (2),
page \citeBib{Cohn: Universal Algebra}\Hyph 48.
}

\ePrints{9835-2163,0767-8264,5284-0163,1801.01628}
\ifx\Semafor\ValueOff
The operator domain $\Omega$ describes a set of $\Omega$\Hyph algebras.
An element of the set $\Omega$ is called operator,
because an operation assumes certain set.
According to the remark
\ref{remark: o in AAn}
and the definition
\RefDefinition{Omega-algebra},
for each operator
\ShowEq{omega n ari}{}{}n,
we match $n$\Hyph ary operation $\omega$ on $A$.
\fi

\ePrints{MAlgebra2,5148-4632,1908.04418,6860-2955}
\ifx\Semafor\ValueOn

\begin{\TheoremStyle}
\labelTheorem{Cartesian power is universal algebra}
Let the set $B$ be $\Omega$\Hyph algebra.
Then the set $B^A$ of maps
\ShowEq{f:A->B}fAB
also is $\Omega$\Hyph algebra.
\end{\TheoremStyle}
\begin{proof}
Let
\ShowEq{omega n ari}{}{}n.
For maps
\ShowEq{f1n in B**A}
we define the operation $\omega$ by the equality
\DrawEq{f1n omega=}{}
\end{proof}
\fi

\ePrints{MAlgebra2,5148-4632,1908.04418,6860-2955}
\ifx\Semafor\ValueOn
\begin{\DefinitionStyle}
\labelDefinition{set is closed with respect to operation}
{\it
Let
\ShowEq{B subset A}.
Since, for any
\ShowEq{b1n in B}bB
\ShowEq{b1no in B}
then we say that $B$
\AddIndex{is closed with respect to}{set is closed with respect to operation}
$\omega$ or that $B$
\AddIndex{admits operation}{set admits operation}
$\omega$.
}
\qed
\end{\DefinitionStyle}

\begin{\DefinitionStyle}
$\Omega$\Hyph algebra $B_{\Omega}$ is
\AddIndex{subalgebra}{subalgebra of Omega-algebra}
of $\Omega$\Hyph algebra $A_{\Omega}$
if following statements are true\,\footnotemark
\StartLabelItem
\begin{enumerate}
\item
\ShowEq{B subset A}.
\item
if operator
\ShowEq{omega in Omega}{}{}
defines operations $\omega_A$ on $A$ and $\omega_B$ on $B$, then
\ShowEq{oAB=oB}
\end{enumerate}
\qed
\end{\DefinitionStyle}
\footnotetext{\,
I follow the definition on
page \citeBib{Cohn: Universal Algebra}\Hyph 48.
}
\fi

\ShowFootnote{iso end aut morphism}1

\begin{\DefinitionStyle}
\labelDefinition{homomorphism}
Let $A$, $B$ be $\Omega$\Hyph algebras and
\ShowEq{omega n ari}{}{}n.
The map\,\refFootnote{iso end aut morphism}1
\ShowEq{f:A->B}fAB
\AddIndex{is compatible with operation}{map is compatible with operation}
$\omega$, if, for all
\ShowEq{b1n in B}aA
\ShowEq{afo=aof}
The map $f$ is called
\AddIndex{homomorphism}{homomorphism}
from $\Omega$\Hyph algebra $A$ to $\Omega$\Hyph algebra $B$,
if $f$ is compatible with each
\ShowEq{omega in Omega}{}.
\ePrints{1908.04418,6860-2955}%
\ifx\Semafor\ValueOn%
We use notation
\ShowEq{set of homomorphisms}
for the set of homomorphisms
from $\Omega$\Hyph algebra $A$ to $\Omega$\Hyph algebra $B$.
\fi
\end{\DefinitionStyle}

\ePrints{1908.04418,6860-2955}
\ifx\Semafor\ValueOn
\begin{\TheoremStyle}
\labelTheorem{Hom empty A B=B**A}
Since operator domain is empty,
then a homomorphism
from $\Omega$\Hyph algebra $A$ to $\Omega$\Hyph algebra $B$
is a map
\ShowEq{f:A->B}fAB
Therefore,
\ShowEq{Hom empty A B=B**A}
\end{\TheoremStyle}
\begin{proof}
The theorem follows from definitions
\RefDefinition{Cartesian power},
\RefDefinition{homomorphism}.
\end{proof}
\fi

\ePrints{MAlgebra2,5148-4632,1908.04418,6860-2955,8428-0408,2207.06506}
\ifx\Semafor\ValueOn
\ShowFootnote{iso end aut morphism}2

\begin{\DefinitionStyle}
\labelDefinition{isomorphism}
Homomorphism
\ShowEq{f:A->B}fAB
is called\,\refFootnote{iso end aut morphism}2
\AddIndex{isomorphism}{isomorphism}
between $A$ and $B$, if correspondence $f^{-1}$ is homomorphism.
If there is an isomorphism between $A$ and $B$, then we
say that $A$ and $B$ are isomorphic and write
\ShowEq{isomorphic}
An injective homomorphis is called
\AddIndex{monomorphism}{monomorphism}.
A surjective homomorphis is called
\AddIndex{epimorphism}{epimorphism}.
\end{\DefinitionStyle}
\fi

\begin{\DefinitionStyle}
\labelDefinition{endomorphism}
A homomorphism
\ShowEq{f:A->B}fAA
in which source and target are the same algebra is called
\AddIndex{endomorphism}{endomorphism}.
We use notation
\ShowEq{set of endomorphisms}
for the set of endomorphisms
of $\Omega$\Hyph algebra $A$.
\ePrints{4975-6381,1506.00061,7287-9339}%
\Items{1601.03259,MAlgebra2,5148-4632,9835-2163,0767-8264,5284-0163,1801.01628}%
\ifx\Semafor\ValueOff%
Isomorphism
\ShowEq{f:A->B}fAA
is called
\AddIndex{automorphism}{automorphism}.
\fi
\end{\DefinitionStyle}

\ePrints{1908.04418,6860-2955}
\ifx\Semafor\ValueOn

\begin{\TheoremStyle}
\labelTheorem{End A=Hom AA}
\ShowEq{End A=Hom AA}
\end{\TheoremStyle}
\begin{proof}
The theorem follows from the definitions
\RefDefinition{homomorphism},
\RefDefinition{endomorphism}.
\end{proof}

\begin{\TheoremStyle}
\labelTheorem{End empty A=A**A}
Since operator domain is empty,
then an endomorphism of the set $A$
is a map
\ShowEq{t:A->A}
Therefore,
\ShowEq{End empty A=A**A}
\end{\TheoremStyle}
\begin{proof}
The theorem follows from the theorems
\RefTheorem{Hom empty A B=B**A},
\RefTheorem{End A=Hom AA}.
\end{proof}
\fi

\ePrints{4975-6381,1506.00061,7287-9339}%
\Items{8428-0408,1601.03259,MAlgebra2,5148-4632,5410-9916,9835-2163,0767-8264,5284-0163,1801.01628}%
\ifx\Semafor\ValueOff%
\begin{\DefinitionStyle}
If there is a monomorphism from $\Omega$\Hyph algebra $A$
to $\Omega$\Hyph algebra $B$, then we say that
\AddIndex{$A$ can be embeded in $B$}{can be embeded}.
\end{\DefinitionStyle}

\begin{\DefinitionStyle}
If there is an epimorphism from $A$ to $B$, then $B$ is called
\AddIndex{homomorphic image}{homomorphic image} of algebra $A$.
\end{\DefinitionStyle}
\fi
\fi

\ePrints{MAlgebra2,5148-4632,4975-6381,5410-9916,9835-2163,0767-8264,5284-0163,1801.01628}
\ifx\Semafor\ValueOn
\ShowConvention{A number}{\Omega}{algebra}
\fi

\ePrints{1908.04418,6860-2955,8428-0408,2207.06506,5284-0163,1801.01628}
\ifx\Semafor\ValueOn

\Section{Cartesian Product of Universal Algebras}

\ePrints{8428-0408,2207.06506,5284-0163,1801.01628}
\ifx\Semafor\ValueOn
\begin{\DefinitionStyle}
\labelDefinition{universally attracting object of category}
Let $\mathcal C$ be a category.
An object $P$ of category $\mathcal C$ is called
\AddIndex{universally attracting}{universally attracting}\,\footnotemark
if, for any object $R$ of category $\mathcal C$,
there exists unique morphism
\ShowEq{f:A->B}fRP
\end{\DefinitionStyle}
\footnotetext{\,
See also the definition in \citeBib{Serge Lang}, pages 57.
}

\begin{\DefinitionStyle}
\labelDefinition{universally repelling object of category}
Let $\mathcal C$ be a category.
An object $P$ of category $\mathcal C$ is called
\AddIndex{universally repelling}{universally repelling}\,\footnotemark
if, for any object $R$ of category $\mathcal C$,
there exists unique morphism
\ShowEq{f:A->B}fPR
\end{\DefinitionStyle}
\footnotetext{\,
See also the definition in \citeBib{Serge Lang}, pages 57.
}

\begin{\DefinitionStyle}
\labelDefinition{product in category, 2020}
Let $\mathcal A$ be a category.
Let
\ShowEq{set Bi}B
be the set of objects of $\mathcal A$.
Let
\ShowEq{category A(Bi)}
be a category
whose objects are tuples $(P,f)$
where $P$ is object of category $\mathcal A$
and $f$ is set of morphisms
\ShowEq{set f:A->B}fP{B_i}
Universally attracting object of category
\ShowEq{category A(Bi)}
\ShowEq{product in category}
is called a
\AddIndex{product of set of objects
\ShowEq{set Bi}B
in category $\mathcal A$}
{product in category}.\,\footnotemark

If $|I|=n$, then we also will use notation
\ShowEq{product in category, 1 n}
for product of set of objects
$\{B_i,\iI\}$ in $\mathcal A$.
\end{\DefinitionStyle}
\footnotetext{\,
I made definition according to the definition on page
\citeBib{Serge Lang}\Hyph 58.
}

\begin{\TheoremStyle}
Let $\mathcal A$ be a category.
Let
\ShowEq{set Bi}B
be the set of objects of $\mathcal A$.
The product of set of objects
\ShowEq{set Bi}B
in category $\mathcal A$
is an object
\ShowEq{product in category}
and set of morphisms
\ShowEq{set f:A->B}fP{B_i}
such that
for any object $R$ and set of morphisms
\ShowEq{set f:A->B}gR{B_i}
there exists a unique morphism
\ShowEq{f:A->B}hRP
such that diagram
\ShowEq{product in category diagram}
is commutative for all $i\in I$.
\end{\TheoremStyle}

\begin{proof}
The theorem follows from definitions
\RefDefinition{universally attracting object of category},
\RefDefinition{product in category, 2020}.
\end{proof}
\else
\ShowDefinition{product in category}
\fi

\begin{example}
Let \(\mathcal S\) be the category of sets.\,\footnote{\,
See also the example in
\citeBib{Serge Lang},
page 59.
}
According to the definition
\ShowEq{ref product in category}
Cartesian product
\ShowEq{Cartesian product of sets}
of family of sets
\ShowEq{Ai iI}A{}
and family of projections on the \(i\)\Hyph th factor
\ShowEq{projection on i factor}
are product in the category \(\mathcal S\).
\qed
\end{example}

\begin{\TheoremStyle}
\labelTheorem{product exists in category of Omega algebras}
The product exists in the category \(\mathcal A\) of \(\Omega\)\Hyph algebras.
Let \(\Omega\)\Hyph algebra \(A\)
and family of morphisms
\ShowEq{p:A->Ai i in I}
be product in the category \(\mathcal A\).
Then
\StartLabelItem
\begin{enumerate}
\item
\labelItem{set A is Cartesian product}
The set \(A\) is Cartesian product
of family of sets
\ShowEq{Ai iI}A{}
\item
\labelItem{homomorphism is projection on factor}
The homomorphism of \(\Omega\)\Hyph algebra
\ShowEq{projection on i factor}
is projection on \(i\)\Hyph th factor.
\item
We can represent any \(A\)\Hyph number $a$
as tuple
\ShowEq{tuple represent A number}
of \(A_i\)\Hyph numbers.
\labelItem{tuple represent A number}
\item
Let
\ShowEq{omega in Omega}{}{}
be n\Hyph ary operation.
Then operation $\omega$ is defined componentwise
\ShowEq{operation is defined componentwise}
where
\ShowEq{a=ai 1n}.
\labelItem{operation is defined componentwise}
\end{enumerate}
\end{\TheoremStyle}
\begin{proof}
Let
\ShowEq{Cartesian product of sets}
be Cartesian product
of family of sets
\ShowEq{Ai iI}A{}
and, for each \iI, the map
\ShowEq{projection on i factor}
be projection on the \(i\)\Hyph th factor.
Consider the diagram of morphisms in category of sets $\mathcal S$
\ShowEq{operation is defined componentwise, diagram}
where the map $g_i$ is defined by the equality
\ShowEq{gi()=}
According to the definition
\ShowEq{ref product in category}
the map $\omega$ is defined uniquely from the set of diagrams
\EqRef{operation is defined componentwise, diagram}
\ShowEq{omega(ai)=(omega ai)}
The equality
\EqRef{operation is defined componentwise}
follows from the equality
\EqRef{omega(ai)=(omega ai)}.
\end{proof}

\begin{\DefinitionStyle}
If \(\Omega\)\Hyph algebra \(A\)
and family of morphisms
\ShowEq{p:A->Ai i in I}
is product in the category \(\mathcal A\),
then \(\Omega\)\Hyph algebra \(A\) is called
\AddIndex{direct}{direct product of Omega algebras}
or
\AddIndex{Cartesian product of \(\Omega\)\Hyph algebras}
{Cartesian product of Omega algebras}
\ShowEq{Ai iI}A.
\end{\DefinitionStyle}

\begin{\TheoremStyle}
\labelTheorem{map from product into product}
Let set \(A\) be
Cartesian product of sets
\ShowEq{Ai iI}A{}
and set \(B\) be
Cartesian product of sets
\ShowEq{Ai iI}B.
For each \iI, let
\ShowEq{f:A->B i}
be the map from the set $A_i$ into the set $B_i$.
For each \iI, consider commutative diagram
\ShowEq{homomorphism of Cartesian product of Omega algebras diagram}
where maps
\ShowEq{pi p'i}
are projection on the \(i\)\Hyph th factor.
The set of commutative diagrams
\EqRef{homomorphism of Cartesian product of Omega algebras diagram}
uniquely defines map
\ShowEq{f:A->B}fAB
\DrawEq{f:A->B=}{}
\end{\TheoremStyle}
\begin{proof}
For each \iI, consider commutative diagram
\ShowEq{homomorphism of Cartesian product of Omega algebras}
Let \(a\in A\).
According to the statement
\RefItem{tuple represent A number},
we can represent \(A\)\Hyph number \(a\)
as tuple of \(A_i\)\Hyph numbers
\ShowEq{a=p(a)i}
Let
\ShowEq{b=f(a)}
According to the statement
\RefItem{tuple represent A number},
we can represent \(B\)\Hyph number \(b\)
as tuple of \(B_i\)\Hyph numbers
\ShowEq{b=p(b)i}
From commutativity of diagram $(1)$
and from equalities
\EqRef{b=f(a)},
\EqRef{b=p(b)i},
it follows that
\ShowEq{b=g(a)i}
From commutativity of diagram $(2)$
and from the equality
\EqRef{a=p(a)i},
it follows that
\ShowEq{b=f(a)i}
\end{proof}

\begin{\TheoremStyle}
\labelTheorem{homomorphism of Cartesian product of Omega algebras}
Let \(\Omega\)\Hyph algebra \(A\) be
Cartesian product of \(\Omega\)\Hyph algebras
\ShowEq{Ai iI}A{}
and \(\Omega\)\Hyph algebra \(B\) be
Cartesian product of \(\Omega\)\Hyph algebras
\ShowEq{Ai iI}B.
For each \iI,
let the map
\ShowEq{f:A->B i}
be homomorphism of \(\Omega\)\Hyph algebra.
Then the map
\ShowEq{f:A->B}fAB
defined by the equality
\DrawEq{f:A->B=}{homomorphism}
is homomorphism of \(\Omega\)\Hyph algebra.
\end{\TheoremStyle}
\begin{proof}
Let
\ShowEq{omega in Omega}{}{}
be n\Hyph ary operation.
Let
\ShowEq{a=ai 1n},
\ShowEq{b=bi 1n}.
From equalities
\EqRef{operation is defined componentwise},
\eqRef{f:A->B=}{homomorphism},
it follows that
\ShowEq{f:A->B omega}
\end{proof}

\begin{\DefinitionStyle}
Equivalence on $\Omega$\Hyph algebra $A$,
which is subalgebra of $\Omega$\Hyph algebra $A^2$,
is called
\AddIndex{congruence}{congruence}\,\footnotemark
on $A$.
\end{\DefinitionStyle}
\footnotetext{\,
I follow the definition on
page \citeBib{Cohn: Universal Algebra}\Hyph 57.
}

\begin{\TheoremStyle}[first isomorphism theorem]
\labelTheorem{first isomorphism theorem}
Let
\ShowEq{f:A->B}fAB
be homomorphism of $\Omega$\Hyph algebras with kernel $s$.
Then there is decomposition
\ShowEq{decomposition of map f}
\StartLabelItem
\begin{enumerate}
\item
The \AddIndex{kernel of homomorphism}{kernel of homomorphism}
\ShowEq{kernel of homomorphism}
is a congruence on $\Omega$\Hyph algebra $A$.
\labelItem{kernel of homomorphism}
\item
The set
\ShowEq{A/ker f}
is $\Omega$\Hyph algebra.
\labelItem{A/ker f is Omega-algebra}
\item
The map
\ShowEq{p:A->/ker}
is epimorphism and is called
\AddIndex{natural homomorphism}{natural homomorphism}.
\labelItem{natural homomorphism}
\item
The map 
\ShowEq{q:A/ker->f(A)}
is the isomorphism
\labelItem{q:A/ker->f(A) isomorphism}
\item
The map 
\ShowEq{r:f(A)->B}
is the monomorphism
\labelItem{r:f(A)->B monomorphism}
\end{enumerate}
\end{\TheoremStyle}
\begin{proof}
The statement
\RefItem{kernel of homomorphism}
follows from the proposition II.3.4
(\citeBib{Cohn: Universal Algebra}, page 58).
Statements
\RefItem{A/ker f is Omega-algebra},
\RefItem{natural homomorphism}
follow from the theorem II.3.5
(\citeBib{Cohn: Universal Algebra}, page 58)
and from the following definition.
Statements
\RefItem{q:A/ker->f(A) isomorphism},
\RefItem{r:f(A)->B monomorphism}
follow from the theorem II.3.7
(\citeBib{Cohn: Universal Algebra}, page 60).
\end{proof}

\ePrints{8428-0408,2207.06506,5284-0163,1801.01628}
\ifx\Semafor\ValueOn
\ShowDefinition{coproduct in category}

\ShowTheorem{coproduct in category}
\ShowProof{coproduct in category}
\fi
\fi

\ePrints{1908.04418,6860-2955}
\ifx\Semafor\ValueOn
\Section{Semigroup}

Usually the operation
\ShowEq{omega n ari}{}{}2{}
is called product
\ShowEq{abo=ab}
or sum
\ShowEq{abo=a+b}

\begin{\DefinitionStyle}
{\it
Let $A$ be $\Omega$\Hyph algebra and
\ShowEq{omega n ari}{}{}2.
$A$\Hyph number $e$ is called
\AddIndex{neutral element of operation}{neutral element of operation}
$\omega$, when for any $A$\Hyph number $a$ following equations are true
\ShowEq{left neutral element}
\ShowEq{right neutral element}
}
\qed
\end{\DefinitionStyle}

\begin{\DefinitionStyle}
{\it
Let $A$ be $\Omega$\Hyph algebra.
The operation
\ShowEq{omega n ari}{}{}2{}
is called
\AddIndex{associative}{associative operation}
if the following equality is true
\ShowEq{associative operation}
}
\qed
\end{\DefinitionStyle}

\begin{\DefinitionStyle}
{\it
Let $A$ be $\Omega$\Hyph algebra.
The operation
\ShowEq{omega n ari}{}{}2{}
is called
\AddIndex{commutative}{commutative operation}
if the following equality is true
\ShowEq{commutative operation}
}
\qed
\end{\DefinitionStyle}

\begin{\DefinitionStyle}
\labelDefinition{semigroup}
{\it
Let
\ShowEq{Omega=omega}
If the operation
\ShowEq{omega n ari}{}{}2{}
is associative, then $\Omega$\Hyph algebra is called
\AddIndex{semigroup}{semigroup}.
If the operation in the semigroup is commutative,
then the semigroup is called
\AddIndex{Abelian semigroup}{Abelian semigroup}.
}
\qed
\end{\DefinitionStyle}
\fi

\ePrints{4975-6381,1506.00061,7287-9339,8428-0408,2207.06506,1601.03259,5148-4632}%
\Items{MAlgebra2,5410-9916,9835-2163,0767-8264,5284-0163,1801.01628}%
\ifx\Semafor\ValueOn%
\Section{Representation of Universal Algebra}
\labelSection{Representation of Universal Algebra}

\ShowDefinition{representation of algebra}

We also use notation
\ShowEq{f:A->*B}f{A_1}{A_2}
to denote the representation of $\Omega_1$\Hyph algebra $A_1$
in $\Omega_2$\Hyph algebra $A_2$.

\begin{\DefinitionStyle}
\labelDefinition{effective representation of algebra}
Let the map
\ShowEq{representation of algebra}
be an isomorphism of the $\Omega_1$\Hyph algebra $A_1$ into
\ShowEq{End O2}{A_2}.
Then the representation
\ShowEq{f:A->*B}f{A_1}{A_2}
of the $\Omega_1$\Hyph algebra $A_1$ is called
\AddIndex{effective}{effective representation}.\,\footnotemark
\end{\DefinitionStyle}
\footnotetext{\,
See similar definition of effective representation of group in
\citeBib{Postnikov: Differential Geometry}, page 16,
\citeBib{Basic Concepts of Differential Geometry}, page 111,
\citeBib{Cohn: Algebra 1}, page 51
(Cohn calls such representation faithful).
}

\ePrints{MAlgebra2,5148-4632}
\ifx\Semafor\ValueOn

\begin{\TheoremStyle}
\labelTheorem{representation is effective}
The representation
\ShowEq{f:A->*B}g{A_1}{A_2}
is effective iff the statement
\ShowEq{a1 ne b1}
implies that there exists
\ShowEq{a in A}2{}
such that\,\footnotemark
\ShowEq{fam ne fbm}
\end{\TheoremStyle}
\footnotetext{\,
In case of group, the theorem
\RefTheorem{representation is effective}
has the following form.
{\it
The representation
\ShowEq{f:A->*B}g{A_1}{A_2}
is effective iff, for any $A_1$\Hyph number $a\ne e$,
there exists
\ShowEq{a in A}2{}
such that
\ShowEq{fam ne m}
}
}
\ProofTheorem{\RefRepresentation}{representation is effective}
\fi

\ePrints{1506.00061,7287-9339,5284-0163,1801.01628,2207.06506,8428-0408}%
\ifx\Semafor\ValueOn%
\begin{\DefinitionStyle}
\labelDefinition{free representation of algebra}
The representation
\ShowEq{f:A->*B}g{A_1}{A_2}
of the $\Omega_1$\Hyph algebra $A_1$ is called
\AddIndex{free}{free representation},\,\footnotemark
if the statement
\ShowEq{faa=fba}
for any
\ShowEq{a in A}2{}
implies that $a=b$.
\end{\DefinitionStyle}
\footnotetext{\,
See similar definition of free representation of group in
\citeBib{Postnikov: Differential Geometry}, page 16.
}

\ePrints{1506.00061,7287-9339}%
\ifx\Semafor\ValueOff%
\begin{\TheoremStyle}
\labelTheorem{Free representation is effective}
Free representation is effective.
\end{\TheoremStyle}
\ProofTheorem{\RefRepresentation}{Free representation is effective}

\begin{remark}
Representation of rotation group in affine space is effective.
However this representation is not free, since origin
is fixed point of every transformation.
\qed
\end{remark}

\begin{\DefinitionStyle}
\labelDefinition{transitive representation of algebra}
The representation
\ShowEq{f:A->*B}g{A_1}{A_2}
of $\Omega_1$\Hyph algebra is called
\AddIndex{transitive}{transitive representation of algebra}\,\footnotemark
if, for any
\ShowEq{ab in A}2,
there exists such $g$ that
\[a=f(g)(b)\]
The representation of $\Omega_1$\Hyph algebra is called
\AddIndex{single transitive}{single transitive representation of algebra}
if it is transitive and free.
\end{\DefinitionStyle}
\footnotetext{\,
See similar definition of transitive representation of group in
\citeBib{Basic Concepts of Differential Geometry}, page 110,
\citeBib{Cohn: Algebra 1}, page 51.
}

\ShowTheorem{Representation is single transitive iff}
\ShowProof{Representation is single transitive iff}

\ShowTheorem{single transitive representation generates algebra}
\ProofTheorem{\RefRepresentation}{single transitive representation generates algebra}
\fi
\fi

\ShowDefinition{morphism of representations of universal algebra}

\ShowRemark{morphism of representations of universal algebra as map}

\ePrints{5148-4632,MAlgebra2,8428-0408,2207.06506,5284-0163,1801.01628}
\ifx\Semafor\ValueOn

\ShowRemark{notation for morphism of representations}

\ShowTheorem{Tuple of maps is morphism of representations iff}
\ProofTheorem{\RefRepresentation}{Tuple of maps is morphism of representations iff}

\ShowRemark{morphism of representations of universal algebra}
\fi

\ShowDefinition{morphism of representation f}

\ePrints{aa}
\ifx\Semafor\ValueOn

\ShowTheorem{unique morphism of representations of universal algebra}
\ProofTheorem{\RefRepresentation}{unique morphism of representations of universal algebra}

\ShowDefinition{transformation coordinated with equivalence}

\ShowTheorem{transformation correlated with equivalence}
\ProofTheorem{\RefRepresentation}{transformation correlated with equivalence}

\ShowTheorem{decompositions of morphism of representations}
\ProofTheorem{\RefRepresentation}{decompositions of morphism of representations}
\fi

\ShowDefinition{reduced morphism of representations}

\ShowTheorem{map is reduced morphism of representations iff}
\ProofTheorem{\RefRepresentation}{map is reduced morphism of representations iff}
\fi

\ePrints{8428-0408,2207.06506}
\ifx\Semafor\ValueOn

\ShowDefinition{category of representations A1(mA2)}

\ShowTheorem{product of effective representations}
\ProofTheorem{\RefRepresentation}{product of effective representations}
\fi

\ePrints{5148-4632,MAlgebra2,5410-9916,8428-0408,2207.06506,5284-0163,1801.01628}
\ifx\Semafor\ValueOn
\input{\FilePrefix Preliminary.Omega.\TheLanguage}
\fi

\ePrints{MAlgebra2,5148-4632,8428-0408,2207.06506,5284-0163,1801.01628}
\ifx\Semafor\ValueOn

\Section{Basis of Representation of Universal Algebra}

\ShowDefinition{stable set of representation}

\ShowTheorem{subrepresentation of representation}
\ProofTheorem{\RefRepresentation}{subrepresentation of representation}

\ShowTheorem{lattice of subrepresentations}
\ProofTheorem{\RefRepresentation}{lattice of subrepresentations}

\ShowRemark{closure operator, representation}

\ShowDefinition{generating set of representation}

\ShowTheorem{structure of subrepresentations}
\ProofTheorem{\RefRepresentation}{structure of subrepresentations}

\ePrints{MAlgebra2,5148-4632,8428-0408,2207.06506,5284-0163,1801.01628}
\ifx\Semafor\ValueOn
\ShowDefinition{word of element relative to set, representation}

\ShowTheorem{map of words of representation}
\ProofTheorem{\RefRepresentation}{map of words of representation}

\ShowDefinition{quasibasis of representation}

\ShowTheorem{X is quasibasis of representation}
\ProofTheorem{\RefRepresentation}{X is quasibasis of representation}

\ShowRemark{X is quasibasis of representation}

\ShowRemark{representation in form of Omega2-word is ambiguous}

\ShowRemark{three reasons of ambiguity in Omega2-word}

\ShowDefinition{equivalence on the set of Omega2-words}

\ShowTheorem{equivalence generated by basis}
\ProofTheorem{\RefRepresentation}{equivalence generated by basis}

\ShowDefinition{basis of representation}
\else
\begin{\DefinitionStyle}
\labelDefinition{basis of representation}
If the set $X\subset A_2$ is generating set of representation
$f$, then any set $Y$, $X\subset Y\subset A_2$
also is generating set of representation $f$.
If there exists minimal set $X$ generating
the representation $f$, then the set $X$ is called
\AddIndex{basis of representation}{basis of representation} $f$.
\qed
\end{\DefinitionStyle}

\begin{\TheoremStyle}
\labelTheorem{X is basis of representation}
The generating set $X$ of representation $f$ is basis
iff for any $m\in X$
the set $X\setminus\{m\}$ is not
generating set of representation $f$.
\end{\TheoremStyle}
\ProofTheorem{\RefRepresentation}{X is basis of representation}

\begin{remark}
\labelRemark{X is basis of representation}
The proof of the theorem
\RefTheorem{X is basis of representation}
gives us effective method for constructing the basis of the representation $f$.
Choosing an arbitrary generating set, step by step, we
remove from set those elements which have coordinates
relative to other elements of the set.
If the generating set of the representation is infinite,
then this construction may not have the last step.
If the representation has finite generating set,
then we need a finite number of steps to construct a basis of this representation.

As noted by Paul Cohn in
\citeBib{Cohn: Universal Algebra}, p. 82, 83,
the representation may have inequivalent bases.
For instance, the cyclic group of order six has
bases $\{a\}$ and $\{a^2,a^3\}$ which we cannot map
one into another by endomorphism of the representation.
\qed
\end{remark}
\fi

\ShowTheorem{automorphism uniquely defined by image of basis}
\ProofTheorem{\RefRepresentation}{automorphism uniquely defined by image of basis}
\fi

\ePrints{MAlgebra2,5148-4632,8428-0408,2207.06506}
\ifx\Semafor\ValueOn

\Section{Diagram of Representations of Universal Algebras}

\ShowDefinition{diagram of representations}

\ShowRemark{diagram of representations}

\ShowDefinition{commutative diagram of representations}

Consider the theorem
\RefTheorem{diagram of representations, define map fik}
for the purposes of illustration of the definition
\RefDefinition{diagram of representations}.

\ShowTheorem{diagram of representations, define map fik}
\ProofTheorem{\RefRepresentation}{diagram of representations, define map fik}

\ShowDefinition{Morphism of Diagram of Representations}

\ShowRemark{Morphism of Diagram of Representations}
\fi

\ePrints{MAlgebra2,5148-4632,4975-6381,1601.03259,309618526,CACAA.06.121}
\ifx\Semafor\ValueOn
\Section{Permutation}

\begin{\DefinitionStyle}
{\it
An injective map of finite set into itself is called
\AddIndex{permutation}{permutation}.\,\footnotemark
}
\qed
\end{\DefinitionStyle}
\footnotetext{\,
You can see definition and properties of permutation in
\citeBib{Kurosh: High Algebra}, pages 27 - 32,
\citeBib{Cohn: Algebra 1}, pages 58, 59.
}

Usually we write a permutation $\sigma$ as a matrix
\ShowEq{permutation as matrix}
The notation
\EqRef{permutation as matrix}
is equivalent to the statement
\ShowEq{a->sigma a}
So the order of columns in the notation
\EqRef{permutation as matrix}
is not essential.

Since there is an order
on the set
\ShowEq{set a1n}{}
(for instance, we assume, that $a_i$ precedes $a_j$
when $i<j$),
then we may assume that elements of first row
are written according to the intended order
and we will identify permutation with second row
\ShowEq{permutation as matrix 2}

\ShowDefinition{parity of permutation}
\fi

\ShowEq{PreliminaryRefOff}

%% file: Preliminary.Representation.Eq.tex
\ShowEq{PreliminaryRefOn}

\AddEq{Phi in AxA}
{
$\Phi\in A\times A$
}

\AddEq{reflexive correspondence}
{
\labelItem{reflexive correspondence}
\[
(a,a)\in\Phi
\]
}

\AddEq{symmetric correspondence}
{
\labelItem{symmetric correspondence}
\[
(a,b)\in\Phi\Rightarrow (b,a)\in\Phi
\]
}

\AddEq{transitive correspondence}
{
\labelItem{transitive correspondence}
\[
(a,b),(b,c)\in\Phi\Rightarrow (a,c)\in\Phi
\]
}

\AddEq{kernel of map}
{
\symb{\mathrm{ker}\,f}{kernel of map}{}
}

\AddEquation{def kernel of map}
{
\ShowSymbol{kernel of map}{}
=\{
(a,b):a,b\in A,f(a)=f(b)
\}
}

\AddEq{ker f}
{
$\mathrm{ker}\,f$
}

\AddEq{fa=fa}
{
\[f(a)=f(a)\]
}

\AddEq [2]{ab in ker}
{
(#1,#2)\in\mathrm{ker}\,f
}

\AddEq [2]{fa=fb}
{
f(#1)=f(#2)
}

\AddEquation{ab in ker -> ba in ker}
{
(a,b)\in\mathrm{ker}\,f\Rightarrow(b,a)\in\mathrm{ker}\,f
}

\AddEquation{ab,ac in ker -> ac in ker}
{
(a,b),(b,c)\in\mathrm{ker}\,f\Rightarrow(a,c)\in\mathrm{ker}\,f
}

\AddEq{lemmas for ker}
{
\RefLemma{ker - reflexive correspondence},
\RefLemma{ker - symmetric correspondence},
\RefLemma{ker - transitive correspondence}
}

\AddEq{image of map}
{
\symb{\im f}{image of map}{}
}

\AddEq{show image of map}
{
\[
\ShowSymbol{image of map}{}
=\{f(a):a\in A\}
\]
}

\AddEq{maps category}
{
\[
\begin{matrix}
f_1:A\rightarrow S_1&\mathrm{ker}\,f_1\supseteq N
\\
f_2:A\rightarrow S_2&\mathrm{ker}\,f_2\supseteq N
\end{matrix}
\]
}

\AddEq{maps category, diagram}
{
\[
\xymatrix{
&S_1\ar[dd]^h
\\
A
\ar[ru]^{f_1}\ar[rd]_{f_2}
\\
&S_2
}
\]
}

\AddEq{maps category, universal}
{
\[
\mathrm{nat}\,N:A\rightarrow A/N
\]
}

\AddEq{maps category, universal, diagram}
{
\[
\xymatrix
{
&A/N\ar[dd]^h
\\
A\ar[ur]^{j=\mathrm{nat}\,N}\ar[dr]_f
\\
&S
}
\]
}

\AddEq{maps category 1}
{
\[j(a_1)=j(a_2)\]
}

\AddEq{maps category 2}
{
\[f(a_1)=f(a_2)\]
}

\AddEq{maps category, h}
{
\[h(\BlueText{j(b)})=f(b)\]
}

\AddEq{Cartesian power}
{
\symb{B^A}{Cartesian power}1
}

\AddEq{a1no=oa1n}
{
$\omega(a_1,...,a_n)$, $a_1...a_n\omega$
}

\AddEq{o in AAn}
{
$\omega\in A^{A^n}$.
}

\AddEq{isomorphic}
{
\symb{A\cong B}{isomorphic}1.
}

\AddEq{t:A->A}
{
\[t:A\rightarrow A\]
}

\AddEq[1]{set Bi}
{
$\{#1_i,\iI\}$
}

\AddEq{product in category, 1 n}
{
\symb[Pi-1]{\prod_{i=1}^nB_i}{product in category}{i 1 n}
\symb[B-0]{B_1\times...\times B_n}{product in category}{1 n}
\[
P=\ShowSymbol{product in category}{i 1 n}
=\ShowSymbol{product in category}{1 n}
\]
}

\AddEq{product in category diagram}
{
\[
\xymatrix{
P\ar[r]^{f_i}&B_i&f_i\circ h=g_i
\\
R\ar[ur]_{g_i}\ar[u]^h
}
\]
}

\AddEq{gi()=}
{
\[
g_i(a_1, ..., a_n)=p_i(a_1)...p_i(a_n)\omega
\]
}

\AddEquation{omega(ai)=(omega ai)}
{
a_1...a_n\omega=(p_i(a_1)...p_i(a_n)\omega,\iI)
}

\AddEquation{homomorphism of Cartesian product of Omega algebras diagram}
{
\xymatrix
{
B\ar[r]^{p'_i}&B_i
\\
A\ar[u]^f\ar[r]_{p_i}&A_i\ar[u]_{f_i}
}
}

\AddEquation{homomorphism of Cartesian product of Omega algebras}
{
\xymatrix
{
B\ar[rrr]^{p'_i}\ar@{}[dr]^(.6){(1)}&&&B_i
\\
&&
\\
A\ar[uu]^f\ar[uurrr]^{g_i}\ar[rrr]_{p_i}&&&A_i\ar[uu]_{f_i}\ar@{}[ul]^(.8){(2)}
}
}

\AddEquation{a=p(a)i}
{
\begin{matrix}
a=(a_i,\iI)&a_i=p_i(a)\in A_i
\end{matrix}
}

\AddEquation{b=f(a)}
{
b=f(a)\in B
}

\AddEquation{b=p(b)i}
{
\begin{matrix}
b=(b_i,\iI)&b_i=p'_i(b)\in B_i
\end{matrix}
}

\AddEq{decomposition of map f}
{
\[
\xymatrix
{
A/\mathrm{ker}\,f\ar[r]^q&f(A)\ar[d]_r
\\
A\ar[u]^p\ar[r]^f&B
}
\ \ \ f=p\circ q\circ r
\]
}

\AddEq{kernel of homomorphism}
{
\symb{\mathrm{ker}\,f}{kernel of homomorphism}0
$\ShowSymbol{kernel of homomorphism}0
=f\circ f^{-1}$
}

\AddEq{A/ker f}
{
$A/\mathrm{ker}\,f$
}

\AddEq{p:A->/ker}
{
\[
p:a\in A\rightarrow a^{\mathrm{ker}\,f}\in A/\mathrm{ker}\,f
\]
}

\AddEq{q:A/ker->f(A)}
{
\[
q:p(a)\in A/\mathrm{ker}\,f\rightarrow f(a)\in f(A)
\]
}

\AddEq{commutative operation}
{
\[ab\omega=ba\omega\]
}

\AddEq{associative operation}
{
\[
a(bc\omega)\omega=(ab\omega)c\omega
\]
}

\AddEq{representation of algebra}
{
\[f:A_1\rightarrow\End(\Omega_2;A_2)\]
}

\AddEq [2]{a in A}
{
$a_{#1}\in A_{#1}$#2
}

\AddEq{fam ne m}
{
\[f(a_1)(a_2)\ne a_2\]
}

\AddEq[2]{ab in A}
{
$a$, $b\in #1$#2
}

\AddEq{ri:A->B}
{
\[
r_i:A_i\rightarrow B_i
\]
}

\AddEquation{morphism of representations of universal algebra, definition, 2}
{
r_2\circ\BlueText{f(a)}=g(\RedText{r_1(a)})\circ r_2
}

\AddEq{F:A1+A2->B1+B2}
{
\[
F:A_1\cup A_2\rightarrow B_1\cup B_2
\]
}

\AddEq{F:A1+A2->B1+B2 1}
{
\[
F(A_1)=B_1 \ \ \ \ F(A_2)=B_2
\]
}

\AddEq{faa=fba}
{
\[f(a_1)(a_2)=f(b_1)(a_2)\]
}

\AddEq{fam ne fbm}
{
\[f(a)(m)\ne f(b)(m)\]
}

\AddEq{a1 ne b1}
{
$a_1\ne b_1$, $a_1$, $b_1\in A_1$,
}

\AddEquation{permutation as matrix}
{
\sigma=
\begin{pmatrix}
a_1&...&a_n
\\
\sigma(a_1)&...&\sigma(a_n)
\end{pmatrix}
}

\AddEquation{permutation as matrix 2}
{
\sigma=
\begin{pmatrix}
\sigma(a_1)&...&\sigma(a_n)
\end{pmatrix}
}

\AddEq{a->sigma a}
{
\[
\sigma:a\in A\rightarrow \sigma(a)\in A\ \ \ A=\{a_1,...,a_n\}
\]
}

\AddEq{Omega=omega}
{
$\Omega=\{\omega\}$.
}

\AddEquation{left neutral element}
{
ea\omega=a
}

\AddEquation{right neutral element}
{
ae\omega=a
}

\AddEq{abo=ab}
{
\[ab\omega=ab\]
}

\AddEq{abo=a+b}
{
\[ab\omega=a+b\]
}

\AddEq{r:f(A)->B}
{
\[
r:f(a)\in f(A)\rightarrow f(a)\in B
\]
}

\AddEq{f:A->B omega}
{
\begin{align*}
f(a_1...a_n\omega)&=f(a_{1i}...a_{ni}\omega,\iI)
\\&=(f_i(a_{1i}...a_{ni}\omega),\iI)
\\&=((f_i(a_{1i}))...(f_i(a_{ni})),\iI)
\\&=(b_{1i}...b_{ni}\omega,\iI)
\end{align*}
\[
f(a_1)...f(a_n)\omega=b_1...b_n\omega
=(b_{1i}...b_{ni}\omega,\iI)
\]
}

\AddEq{b=bi 1n}
{
\(b_1=(b_{1i},\iI)\), ..., \(b_n=(b_{ni},\iI)\)
}

\AddEquation{b=g(a)i}
{
b_i=g_i(b)
}

\AddEq{f:A->B=}
{
f(a_i,\iI)=(f_i(a_i),\iI)
}

\AddEq{pi p'i}
{
\(p_i\), \(p'_i\)
}

\AddEq{f:A->B i}
{
\[f_i:A_i\rightarrow B_i\]
}

\AddEquation{operation is defined componentwise, diagram}
{
\xymatrix{
A\ar[r]^{p_i}&A_i&p_i\circ \omega=g_i
\\
A^n\ar[ur]_{g_i}\ar[u]^{\omega}
}
}

\AddEq{Cartesian product of sets}
{
\[A=\prod_{\iI}A_i\]
}

\AddEq{projection on i factor}
{
\[p_i:A\rightarrow A_i\]
}

\AddEq{tuple represent A number}
{
\((p_i(a),\iI)\)
}

\AddEq{p:A->Ai i in I}
{
\[p_i:A\rightarrow A_i\ \ \ \iI\]
}

\AddEq[2]{Ai iI}
{
\((#1_i,\iI)\)#2
}

\AddEq[3]{set f:A->B}
{
\[
\{#1_i:#2\rightarrow #3,i\in I\}
\]
}

\AddEq{End empty A=A**A}
{
$\End(\emptyset;A)=A^A$.
}

\AddEq{End A=Hom AA}
{
$\End(\Omega;A)=\Hom(\Omega;A\rightarrow A)$
}

\AddEq{set of endomorphisms}
{
\symb{\End(\Omega;A)}{set of endomorphisms}1
}

\AddEq{Hom empty A B=B**A}
{
$\Hom(\emptyset;A\rightarrow B)=B^A$.
}

\AddEq{set of homomorphisms}
{
\symb{\Hom(\Omega;A\rightarrow B)}{set of homomorphisms}1
}

\AddEq[2]{omega in Omega}
{
\(\omega\in\Omega_{#1}\)#2
}

\AddEq{a=ai 1n}
{
\(a_1=(a_{1i},\iI)\), ..., \(a_n=(a_{ni},\iI)\)
}

\AddEquation{operation is defined componentwise}
{
a_1...a_n\omega=(a_{1i}...a_{ni}\omega,\iI)
}

\AddEq{oAB=oB}
{
\[\omega_A|B=\omega_B\]
}

\AddEq{a(o)=n}
{
$a(\omega)=n$,
}

\AddEq{set of n-ary operators}
{
\symb{\Omega(n)}{set of n-ary operators}{}
}

\AddEq{Omega-algebra}
{
\symb{A_{\Omega}}{Omega-algebra}1
}

\AddEq{f1n in B**A}
{
$f_1$, ..., $f_n\in B^A$,
}

\AddEq [1]{B subset A}
{
$B\subseteq A$#1
}

\AddEq{b1no in B}
{
$b_1...b_n\omega\in B$,
}

\AddEq{f1n omega=}
{
(f_1...f_n\omega)(x)=f_1(x)...f_n(x)\omega
}

\AddEq[4]{omega n ari}
{
$\omega_{#1}\in\Omega_{#2}(#3)$#4
}

\AddEquation{afo=aof}
{
f(a_1)...f(a_n)\omega=f(a_1...a_n\omega)
}

\AddEq{O(n)->AAn}
{
\[\Omega(n)\rightarrow A^{A^n}\ \ \ n\in N\]
}

\AddEq{set of n-ary operators =}
{
\[
\ShowSymbol{set of n-ary operators}{}=
\{\omega\in\Omega:a(\omega)=n\}
\]
}

\AddEq{operator domain}
{
\symb{\Omega}{operator domain}1
}

\AddEq[2]{b1n in B}
{
$#1_1$, ..., $#1_n\in #2$,
}

\AddEq{o:An->A}
{
\[\omega:A^n\rightarrow A\]
}

\AddEquation{maps category, universal, ker}
{
\mathrm{ker}\,f\supseteq N
}

%% file: Preliminary.Omega.English.tex
\input{Preliminary.Omega.Eq}

\ePrints{4975-6381,1601.03259,309618526,CACAA.06.121}
\ifx\Semafor\ValueOn
\Chapter{Lebesgue Integral in Abelian \texorpdfstring{$\Omega$}{Omega}-Group}
\fi

\ePrints{5148-4632,MAlgebra2,5410-9916,4975-6381,1601.03259,309618526,CACAA.06.121}
\Items{8428-0408,2207.06506,5284-0163,1801.01628}
\ifx\Semafor\ValueOn
\ePrints{5284-0163,1801.01628}
\ifx\Semafor\ValueOff
\Section{\texorpdfstring{$\Omega$}{Omega}-Group}

\begin{ShadedDefinition}
\labelDefinition{additive map}
Let sum be defined
in $\Omega$\Hyph algebra $A$.
A map
\ShowEq{f:A->B}fAA
of $\Omega_1$\Hyph algebra $A$
is called
\AddIndex{additive map}{additive map}
if
\ShowEq{f(a+b)=}
\end{ShadedDefinition}

\ShowDefinition{polyadditive map}

\ShowDefinition{Omega group}

\ShowTheorem{operation is distributive over addition}
\begin{proof}
The theorem follows from definitions
\RefDefinition{polyadditive map},
\RefDefinition{Omega group}.
\end{proof}

\ShowEq{contents: norm on Omega group}
\fi

\ePrints{5284-0163,1801.01628}
\ifx\Semafor\ValueOn
\Section{Representation of Multiplicative \texorpdfstring{$\Omega$}{Omega}-Group}
\fi
\ePrints{5410-9916,4975-6381,1601.03259}
\ifx\Semafor\ValueOff
\ShowDefinition{multiplicative Omega group}
\fi

\ifx\texFuture\Defined
\begin{ShadedDefinition}
\labelDefinition{multiplicative map}
{\it
Let $A$, $B$ be
multiplicative $\Omega$\Hyph groups.
The map
\ShowEq{f:A->B}fAB
is called
\AddIndex{multiplicative}{multiplicative map},
if
\DrawEq{f(ab)=f(a)f(b)}{}
}
\end{ShadedDefinition}
\fi

\ePrints{5148-4632,MAlgebra2,8428-0408,2207.06506}
\ifx\Semafor\ValueOn
\ShowDefinition{Abelian multiplicative Omega group}

\ShowDefinition{Omega ring}

\ShowTheorem{product in Omega ring is distributive over addition}
\ShowProof{product in Omega ring is distributive over addition}

\ePrints{8428-0408,2207.06506}
\ifx\Semafor\ValueOff
\ShowDefinition{matrix over Omega ring}

\ShowConvention{Einstein summation}

\ShowText{matrix operations}

\ShowDefinition{biring matrix}
 
\ShowTheorem{transpose of identity}
\ShowProof{transpose of identity}

\ShowTheorem{duality principle for biring}

\ShowTheorem{duality principle for biring of matrices}

\ShowEq{remark: reducible biring}
\fi
\fi
\fi

\ePrints{aaa}
\ifx\Semafor\ValueOn
\begin{ShadedDefinition}
\labelDefinition{orbit of left-side representation}
Let the map
\ShowEq{f:A->*B}g{A_1}{A_2}
be the left\Hyph side representation
of multiplicative $\Omega$\Hyph group $A_1$
in $\Omega_2$\Hyph algebra $A_2$.
For any
\ShowEq{a in A}2,
we define
\AddIndex{orbit of representation}{orbit of representation}
of the multiplicative $\Omega$\Hyph group $A_1$ as set
\ShowEq{orbit of representation}
}
\end{ShadedDefinition}
\fi

\ePrints{5148-4632,MAlgebra2,5284-0163,1801.01628,8428-0408,2207.06506}
\ifx\Semafor\ValueOn
\ePrints{5284-0163,1801.01628}
\ifx\Semafor\ValueOff
\Section{Representation of Multiplicative \texorpdfstring{$\Omega$}{Omega}-Group}
\fi

We assume that transformations of representation
\ShowEq{f:A->*B}g{A_1}{A_2}
of multiplicative $\Omega$\Hyph group $A_1$
in $\Omega_2$\Hyph algebra $A_2$
may act on $A_2$\Hyph numbers either on the left or on the right.
We must consider an endomorphism of $\Omega_2$\Hyph algebra $A_2$ as operator.
In such case, the following conditions must be satisfied.

\ShowText{representation of group}

\ShowDefinition{Left side covariant representation}
\ShowDefinition{Left side contravariant representation}

\ShowTheorem{left right shift}
\ProofTheorem{\RefRepresentation}{left right shift}

\ShowEq{def left}
\ShowDefinition{orbit of representation}

\ShowEq{def right}
\ShowDefinition{orbit of representation}

\ShowTheorem{proper definition of orbit}
\ProofTheorem{\RefRepresentation}{proper definition of orbit}

\ShowTheorem{duality principle, algebra representation}

\ShowTheorem{single transitive representation of group}
\ProofTheorem{\RefRepresentation}{single transitive representation of group}

\ShowTheorem{shifts on group commuting}
\ShowProof{shifts on group commuting}

\ShowTheorem{two representations of group}
\ShowProof{two representations of group}

\ShowDefinition{twin representations}

\ShowRemark{one representation of group}
\fi

\ePrints{MAlgebra2,5148-4632,8428-0408,2207.06506}
\ifx\Semafor\ValueOn

\Section{Tensor Product of Representations}

\ShowDefinition{reduced polymorphism of representations}

\ShowDefinition{tensor product of representations}

\ShowTheorem{tensor product of representations}
\ProofTheorem{\RefRepresentation}{tensor product of representations}

\ShowTheorem{representation, tensor product}
\ProofTheorem{\RefRepresentation}{representation, tensor product}

\ShowTheorem{tensor product and polymorphism}
\ProofTheorem{\RefRepresentation}{tensor product and polymorphism}

\ShowTheorem{B times->B otimes}
\ProofTheorem{\RefRepresentation}{B times->B otimes}
\fi

\ePrints{4975-6381,1601.03259,309618526,CACAA.06.121}
\ifx\Semafor\ValueOn

\Section{Algebra of Sets}

\ShowEq{contents: Algebra of Sets}

\Section{Lebesgue Integral}

\ShowEq{contents: Lebesgue Integral}

\ShowEq{PreliminaryRefOmegaNormOff}
\ShowEq{PreliminaryRefMeasureOff}
\fi

%% file: Preliminary.Omega.Eq.tex

\ePrints{5148-4632,MAlgebra2}
\ifx\Semafor\ValueOn
\input{\FilePrefix Biring.Stmt.\TheLanguage}
\fi

\ePrints{4975-6381,1601.03259,309618526,CACAA.06.121}
\ifx\Semafor\ValueOn
\ePrints{4975-6381,1601.03259}
\ifx\Semafor\ValueOn
\ShowEq{PreliminaryRefOmegaNormOn}
\fi
\ShowEq{PreliminaryRefMeasureOn}
\fi

\DefEq
{
\ePrints{5410-9916,4975-6381,309618526,CACAA.06.121,1601.03259}
\ifx\Semafor\ValueOn
\ShowDefinition{norm on Omega group}

\ShowTheorem{|a-b|>|a|-|b|}
\ePrints{CACAA.06.121}%
\ifx\Semafor\ValueOff%
\ProofTheorem{\RefTheoremOmegaNorm}{|a-b|>|a|-|b|}
\fi

\ShowEq{=Omega}

\ePrints{4975-6381,1601.03259}
\ifx\Semafor\ValueOff
\ShowDefinition{open ball}

\ShowDefinition{closed ball}
\fi

\ShowDefinition{open set}

\ePrints{5410-9916}
\ifx\Semafor\ValueOn
\ShowDefinition{continuous map, Omega group}

\ShowTheorem{continuous map, open set}
\ProofTheorem{\RefTheoremOmegaNorm}{continuous map, open set}

\ShowTheorem{image of interval is interval, real field}
\ProofTheorem{\RefTheoremOmegaNorm}{image of interval is interval, real field}

\ShowDefinition{norm of operation}

\ShowTheorem{|fx|<|f||x|1n}
\ProofTheorem{\RefTheoremOmegaNorm}{|fx|<|f||x|1n}

\ShowDefinition{norm of representation}

\ShowTheorem{|fab|<|f||a||b|}
\ProofTheorem{\RefTheoremOmegaNorm}{|fab|<|f||a||b|}

\ShowDefinition{compact set}

\ShowTheorem{norm of compact set is bounded}
\ProofTheorem{\RefTheoremOmegaNorm}{norm of compact set is bounded}
\fi

\ShowDefinition{limit of sequence}

\ShowDefinition{fundamental sequence}

\ePrints{5410-9916}
\ifx\Semafor\ValueOn
\ShowTheorem{lim a=lim b}
\ProofTheorem{\RefTheoremOmegaNorm}{lim a=lim b}
\fi

\ShowDefinition{complete Omega group}

\ePrints{5410-9916}
\ifx\Semafor\ValueOn
\ShowTheorem{c1+c2 in B}
\ProofTheorem{\RefTheoremOmegaNorm}{c1+c2 in B}

\ShowTheorem{set of maps to Omega group}
\ProofTheorem{\RefTheoremOmegaNorm}{set of maps to Omega group}

\ShowEq{remark: set of maps to Omega group}

\ShowDefinition{limit of sequence, map to Omega group}

\ShowDefinition{sequence converges uniformly}

\ShowTheorem{sequence converges uniformly, fn-fm}
\ProofTheorem{\RefTheoremOmegaNorm}{sequence converges uniformly, fn-fm}

\ShowTheorem{h=f+g, converges uniformly}
\ProofTheorem{\RefTheoremOmegaNorm}{h=f+g, converges uniformly}

\ShowTheorem{h=f1...fn omega, converges uniformly}
\ProofTheorem{\RefTheoremOmegaNorm}{h=f1...fn omega, converges uniformly}

\ShowTheorem{representation f generates representation fX}
\ProofTheorem{\RefTheoremOmegaNorm}{representation f generates representation fX}

\ShowTheorem{fX(g1)(g2), converges uniformly}
\ProofTheorem{\RefTheoremOmegaNorm}{fX(g1)(g2), converges uniformly}
\fi

\ePrints{5410-9916}
\ifx\Semafor\ValueOff
\ShowDefinition{series converges normally}
\fi

\fi
}
{contents: norm on Omega group}

\DefEq
{
\ShowDefinition{semiring of sets}

\ShowDefinition{algebra of sets}

\ShowEq{remark: A cup/minus B}

\ShowDefinition{sigma algebra of sets}

\ShowDefinition{Borel algebra}
}
{contents: Algebra of Sets}

\AddEq{contents: Lebesgue Integral}
{
\ShowDefinition{C_X C_Y measurable}

\ShowEq{example: measurable map}



\ShowDefinition{sigma-additive measure}

\ShowDefinition{simple map}

\ShowTheorem{measurable simple map}
\ProofTheorem{\RefMeasure}{measurable simple map}

\ShowEq{remark - measure - effective representation of real field}

\ShowDefinition{integrable map}

\ShowDefinition{Integral of Map over Set of Finite Measure}

\ShowTheorem{|int f|<int|f|}
\ProofTheorem{\RefMeasure}{|int f|<int|f|}

\ShowTheorem{int |f|<M mu}
\ProofTheorem{\RefMeasure}{int |f|<M mu}

\ShowTheorem{int f+g X}
\ProofTheorem{\RefMeasure}{int f+g X}

\ShowTheorem{|int h|<int|omega||f1n|}
\ProofTheorem{\RefMeasure}{|int h|<int|omega||f1n|}

\ShowTheorem{h=fX g1 g2 representation}
\ProofTheorem{\RefMeasure}{h=fX g1 g2 representation}
}

%% file: Preliminary.D.Module.English.tex
\input{Preliminary.D.Module.Eq}

\Chapter{Linear Algebra}

\ShowText{Preliminary Definitions}

\ePrints{4975-6381,1506.00061,7287-9339}%
\Items{1601.03259,309618526,CACAA.06.121,9835-2163,0767-8264,5284-0163,1801.01628}%
\Items{MAlgebra2,5148-4632}%
\ifx\Semafor\ValueOn%
\Section{Module over Ring}

\ShowText{Preliminary, Module over Ring}

\section{\texorpdfstring{$D$}{D}-Module Type}

\ShowText{Preliminary, D-module type}

\section{Linear Map of \texorpdfstring{$D$}{D}-Module}

\subsection{General Definition}

\ShowText{Preliminary, linear map}1122{h(p)}{h(a)}

\subsection{Linear Map When Rings
\ShowEq{A1=A2 pdf}D}

\ShowText{Preliminary, linear map}{}1{}2pa

\ShowDefinition{D module, endomorphism}

\ShowTheorem{set of endomorphisms - D module}
\ProofTheorem{\RefLinearMap}{set of endomorphisms - D module}

\ePrints{1601.03259,4975-6381,9835-2163,5284-0163,1801.01628}
\ifx\Semafor\ValueOn
\begin{ShadedDefinition}
\labelDefinition{tensor product of algebras}
Let $A_1$, ..., $A_n$ be
free modules over commutative ring $D$.\,\footnotemark
Consider category $\mathcal A_1$ whose objects are
polylinear maps
\ShowEq{polylinear maps category}
where $S_1$, $S_2$ are modules over ring $D$,
We define morphism
\ShowEq{f->g}fg
to be linear map
\ShowEq{f->g 1}
making following diagram commutative
\ShowEq{polylinear maps category, diagram}
Universal object
\ShowEq{tensor product of algebras}
of category $\mathcal A_1$ is called
\AddIndex{tensor product}{tensor product}
of modules $A_1$, ..., $A_n$.
\qed
\end{ShadedDefinition}
\footnotetext{\,
I give definition
of tensor product of $D$\Hyph modules
following to definition in \citeBib{Serge Lang}, p. 601 - 603.
}
\fi

\ePrints{1601.03259,4975-6381,5284-0163,1801.01628}
\ifx\Semafor\ValueOn
\ShowTheorem{there exists tensor product of modules}
\ProofTheorem{\RefLinearMap}{there exists tensor product of modules}

\begin{ShadedTheorem}
\labelTheorem{Tensor product is distributive over sum}
Let $D$ be the commutative ring.
Let $A_1$, ..., $A_n$ be $D$\Hyph modules.
Tensor product is distributive over sum
\ShowEq{tensors 1, tensor product}
The representation of the ring $D$
in tensor product is defined by equation
\ShowEq{tensors 2, tensor product}
\end{ShadedTheorem}
\ProofTheorem{\RefLinearMap}{Tensor product is distributive over sum}

\ShowTheorem{tensor product and polylinear map}
\begin{proof}
See the proof of theorems
\ShowEq{ref tensor product and polylinear map}
\end{proof}

\begin{convention}
\labelConvention{isomorphic representations S1=S2}
Algebras $S_1$, $S_2$ may be different sets.
However they are indistinguishable for us when we consider them
as isomorphic representations.
In such case, we write the statement $S_1=S_2$.
\qed
\end{convention}

\begin{ShadedTheorem}
\labelTheorem{tensor product is associative}
\ShowEq{A1xA2xA3}
\end{ShadedTheorem}
\ProofTheorem{\RefLinearMap}{tensor product is associative}

\begin{\DefinitionStyle}
\labelDefinition{tensor power of algebra}
Tensor product
\ShowEq{tensor power of algebra}
is called
\AddIndex{tensor power}{tensor power} of module $A$.
\end{\DefinitionStyle}

\begin{ShadedTheorem}
\labelTheorem{V times->V otimes}
The map
\ShowEq{V times->V otimes}
is polylinear map.
\end{ShadedTheorem}
\ProofTheorem{\RefLinearMap}{V times->V otimes}
\fi

\ePrints{4975-6381,8428-0408,1601.03259,5284-0163,1801.01628}
\ifx\Semafor\ValueOn
\ShowTheorem{standard component of tensor, module}
\ProofTheorem{\RefLinearMap}{standard component of tensor, algebra}
\fi
\fi

\input{\FilePrefix Preliminary.D.Algebra.\TheLanguage}

%% file: Preliminary.D.Module.Eq.tex

\DefText{Preliminary, Module over Ring}
{
\ShowEq{def DModule}
\ShowEq{def ab=ba}
\ShowEq{def universal}

\ShowText{sum of transformations of Abelian group}

\ShowTheorem{Representation of ring f(0)=v0}
\ProofTheorem{\RefLinearMap}{Representation of ring f(0)=v0}

\ShowTheorem{effective representation of the ring}
\ProofTheorem{\RefLinearMap}{effective representation of the ring}

\ShowDefinition{module over commutative ring}

\ShowTheorem{module over algebra}
\proofTheorem{\RefLinearMap}{module over algebra}{\SideNS}

\def\TheoremStyle{ParacolTheorem}%
\ShowTheorem{definition of A module}
\def\TheoremStyle{ShadedTheorem}%
\proofTheorem{\RefLinearMap}{definition of A module}{\SideWS \VectorSetNS}

\ShowTheorem{definition of A module, property}
\proofTheorem{\RefLinearMap}{definition of A module, property}{\SideNS}

\ShowTheorem{set of vectors generated by set of vectors}
\proofTheorem{\RefLinearMap}{set of vectors generated by set of vectors}{\SideWS module}

\ShowDefinition{linear combination of vectors}

\ShowConvention{sum av() convention}

\ShowDefinition{generating set of module}

\ShowDefinition{basis of module}

\ShowDefinition{coordinates of vector}

\ShowTheorem{quasibasis of module}
\proofTheorem{\RefLinearMap}{quasibasis of module}{\SideNS}

\ShowTheorem{quasibasis of module is basis}
\proofTheorem{\RefLinearMap}{quasibasis of module is basis}{\SideWS \VectorSetNS}

\ShowTheorem{coordinates of vector}
\proofTheorem{\RefLinearMap}{coordinates of vector}{\SideNS-\Cols}

\ShowDefinition{free module over ring}
}

\DefText{Preliminary, D-module type}
{
\ShowText{D-module type}

\ShowEq{\DefCol}
\ShowDefinition{module type}

\ShowTheorem{coordinate matrix of vector}
\proofTheorem{\RefLinearMap}{coordinate matrix of vector}{\SideNS-\Cols}

\ShowTheorem{coordinates of vector}
\proofTheorem{\RefLinearMap}{coordinates of vector}{\SideNS-\Cols}

\ShowEq{\DefRow}
\ShowDefinition{module type}

\ShowTheorem{coordinate matrix of vector}
\proofTheorem{\RefLinearMap}{coordinate matrix of vector}{\SideNS-\Cols}

\ShowTheorem{coordinates of vector}
\proofTheorem{\RefLinearMap}{coordinates of vector}{\SideNS-\Cols}

\ShowText{D-module type 1}
}

\DefText[6]{Preliminary, linear map}
{
\,
\ShowEq{def DModule}
\ShowDefinition{linear map of D module}{#1}{#2}{#3}{#4}

\ShowTheorem{linear map of D module}{#1}{#2}{#3}{#4}{#5}
\proofTheorem{\RefLinearMap}{linear map of D module}{11}

\ShowEq{\DefCol}
\ShowFootnote{homomorphism of D module}{#1}{#2}{#3}{#4}
\ShowEq{\DefRow}
\ShowFootnote{homomorphism of D module}{#1}{#2}{#3}{#4}

\TwoColText
{
\ShowEq{\DefCol}
\ShowTheorem{linear map of D module, coordinates}{#1}{#2}{#3}{#4}{#6}
\proofTheorem{\RefLinearMap}{linear map of D module, coordinates}{#1#2\Cols}
}
{
\ShowEq{\DefRow}
\ShowTheorem{linear map of D module, coordinates}{#1}{#2}{#3}{#4}{#6}
\proofTheorem{\RefLinearMap}{linear map of D module, coordinates}{#1#2\Cols}
}

\TwoColText
{
\ShowEq{\DefCol}
\ShowTheorem{matrix generates D module homomorphism}{#1}{#2}{#3}{#4}
\proofTheorem{\RefLinearMap}{matrix generates D module homomorphism}{\Cols(#1#2)}
}
{
\ShowEq{\DefRow}
\ShowTheorem{matrix generates D module homomorphism}{#1}{#2}{#3}{#4}
\proofTheorem{\RefLinearMap}{matrix generates D module homomorphism}{\Cols(#1#2)}
}

\ePrints{}
\ifx\Semafor\ValueOn
\ShowDefinition{polylinear map of modules}

\ShowTheorem{polylinear map of modules}
\ProofTheorem{\RefLinearMap}{polylinear map of modules}

\ShowTheorem{sum of polylinear maps, module}
\ProofTheorem{\RefLinearMap}{sum of polylinear maps, module}

\ShowCorollary{sum of linear maps, D module}

\ShowTheorem{module of polylinear maps}
\ProofTheorem{\RefLinearMap}{module of polylinear maps}

\ShowCorollary{product of linear map over scalar, D module}
\fi
}

%% file: Preliminary.D.Algebra.English.tex

\input{Preliminary.D.Algebra.Eq}

\ePrints{4975-6381,1506.00061,7287-9339}%
\Items{1601.03259,MAlgebra2,309618526,CACAA.06.121,9835-2163,0767-8264,5284-0163,1801.01628}%
\ifx\Semafor\ValueOn%
\Section{Algebra over Commutative Ring}

\ShowDefinition{algebra over ring}

\ePrints{1601.03259}
\ifx\Semafor\ValueOff
\ShowTheorem{multiplication in algebra is distributive over addition}
\ProofTheorem{\RefLinearMap}{multiplication in algebra is distributive over addition}
\fi

\ePrints{1506.00061,7287-9339}%
\Items{4975-6381,1601.03259,309618526,CACAA.06.121}%
\ifx\Semafor\ValueOn%
\ShowConvention{A number}D{algebra}
\fi

\ShowDefinition{unital algebra}

The multiplication in algebra can be neither commutative
nor associative. Following definitions are based
on definitions given in \citeBib{Richard D. Schafer}, page 13.

\ShowDefinition{commutator of algebra}

\ShowDefinition{associator of algebra}

\ShowDefinition{nucleus of algebra}

\ShowDefinition{center of algebra}

\ePrints{5284-0163,1801.01628}
\ifx\Semafor\ValueOn
\ProveTheorem{unit of algebra and ring}
\fi

\ePrints{5284-0163,1801.01628,0767-8264}
\ifx\Semafor\ValueOn

\ProveTheorem{c in S=>bc=cb submodule of A}

\ProveTheorem{c in Za=> p(c) in Za}

\ProveTheorem{a in ZAb=>b in ZAa}

\ProveTheorem{a in ZAb, c in ZA=>a in ZAbc}

\ShowDefinition{division algebra}

\ProveTheorem{division algebra}

\ShowDefinition{similar A numbers}

\ePrints{0767-8264}
\ifx\Semafor\ValueOff
\ProveTheorem{set of maps B->A is Abelian group}
\fi
\fi

\ShowText{D-Algebra Type}

\ePrints{309618526,1601.03259,CACAA.06.121}
\ifx\Semafor\ValueOff
\ShowText{Preliminary product in algebra}
\fi

\ePrints{4975-6381,1601.03259,309618526,CACAA.06.121,5284-0163,1801.01628}
\Items{5148-4632,MAlgebra2}%
\ifx\Semafor\ValueOn
\ShowTheorem{Free Algebra over Ring}
\begin{proof}
\ePrints{5148-4632,MAlgebra2}
\ifx\Semafor\ValueOn
The theorem follows from theorems
\RefTheorem[\RefLinearMap]{endomorphism of module from product},
\RefTheorem[\RefLinearMap]{Free Algebra over Ring}.
\else
\TheoremFollows
\RefTheorem[\RefLinearMap]{endomorphism of module from product}
and the corollary
\RefCorollary[\RefLinearMap]{Free Algebra over Ring}.
\fi
\end{proof}
\fi

\ShowText{Preliminary homomorphism D algebra}

\Section{Linear Map of \texorpdfstring{$D$}{D}-Algebra}

\ShowEq{=DD}

\ShowDefinition{linear map from A1 to A2, algebra}

\ePrints{1601.03259,4975-6381}
\ifx\Semafor\ValueOn

\begin{ShadedTheorem}
\labelTheorem{module L(A;A) is algebra}
Let $A$ be algebra over commutative ring $D$.
$D$\Hyph module $\mathcal L(D;A;A)$ equiped by product
\ShowEq{module L(A;A) is algebra}
\ShowEq{product of linear map, algebra 1}
is algebra over $D$.
\end{ShadedTheorem}
\ProofTheorem{\RefLinearMap}{module L(A;A) is algebra}
\fi

\begin{ShadedDefinition}
\labelDefinition{polylinear map of algebras}
Let $A_1$, ..., $A_n$, $S$ be $D$\Hyph algebras.
Polylinear map
\ShowEq{f:A->B}f{\Times}S
of $D$\Hyph modules
$A_1$, ..., $A_n$
into $D$\Hyph module $S$
is called
\AddIndex{polylinear map}{polylinear map}
of $D$\Hyph algebras
$A_1$, ..., $A_n$
into $D$\Hyph algebra $S$.
Let us denote
\ShowEq{set polylinear maps}
set of polylinear maps
of $D$\Hyph algebras
$A_1$, ..., $A_n$
into $D$\Hyph algebra
$S$.
Let us denote
\ShowEq{set polylinear maps An}DAS
set of $n$\hyph linear maps
of $D$\Hyph algebra $A_1$ ($A_1=...=A_n=A_1$)
into $D$\Hyph algebra
$S$.
\qed
\end{ShadedDefinition}

\ePrints{9835-2163,5284-0163,1801.01628}
\ifx\Semafor\ValueOff
\ShowTheorem{tensor product of D-algebras is D-algebra}
\ProofTheorem{\RefLinearMap}{tensor product of D-algebras is D-algebra}
\fi


\ShowTheorem{representation of algebra A2 in LA}
\ProofTheorem{\RefLinearMap}{representation of algebra A2 in LA}

\begin{convention}
I assume sum over index $i$
in expression like
\ShowEq{Sum over repeated index}
\qed
\end{convention}


\ePrints{8428-0408,5284-0163,1801.01628}
\ifx\Semafor\ValueOn
\ShowTheorem{standard representation of map A->A, associative algebra}
\ProofTheorem{\RefLinearMap}{standard representation of map A1 A2, associative algebra}

\ShowDefinition{component of linear map}

\ePrints{5284-0163,1801.01628}
Let $A$ be $D$\Hyph module.
Let $B$ be free finite dimensional associative $D$\Hyph algebra.
According to the theorem
\RefTheorem[\RefLinearMap]{linear map in L(A,A), associative algebra},
left \BoxB{B}module
\ShowEq{L(A->B)}DAB
has finite
basis $\Basis F_{AB}$.
However, it is obvious that in the general case the choice of this basis is arbitrary.
At the same time, the choice of basis $\Basis F$ of left \BoxB{B}module
\ShowEq{L(A->B)}DBB{}
is usually associated with a structure of $D$\Hyph algebra $B$.
So the question arises whether there exists relation between bases
$\Basis F_{AB}$ and $\Basis F$.
\ifx\Semafor\ValueOn
\fi

\ShowTheorem{set FoG generates left module L(A;B), n<m}
\ShowProof{set FoG generates left module L(A;B), n<m}

\ShowTheorem{set FoG generates left module L(A;B), n>m}
\ShowProof{set FoG generates left module L(A;B), n>m}

\ShowRemark{set FoG generates left module L(A;B), n>m}
\ePrints{5284-0163,1801.01628}
\ifx\Semafor\ValueOn
It is easy to see that theorem
\RefTheorem{set FoG generates left module L(A;B), n<m}
is particular case of the theorem
\RefTheorem{set FoG generates left module L(A;B), n>m},
because, in the theorem
\RefTheorem{set FoG generates left module L(A;B), n>m},
$\ker G=\emptyset$.

\ShowTheorem{basis D module L(D->A)}
\ShowProof{basis D module L(D->A)}
\fi
\fi

\ePrints{9835-2163,8428-0408}
\ifx\Semafor\ValueOn
\ShowTheorem{set FoG generates left module L(A;B)}
\ePrints{8428-0408}
\ifx\Semafor\ValueOn
\ShowProof{set FoG generates left module L(A;B)}
\else
\ProofTheorem{8428-0408}{set FoG generates left module L(A;B)}
\fi
\fi

\ePrints{8428-0408,9835-2163,5284-0163,1801.01628}
\ifx\Semafor\ValueOn
\ShowTheorem{coordinates of map A1 A2 FoG, algebra}
\ProofTheorem{\RefLinearMap}{coordinates of map A1 A2, algebra}
\fi

\ePrints{4975-6381,1506.00061,7287-9339,1601.03259}
\ifx\Semafor\ValueOn
\ShowTheorem{standard representation of map A->A, associative algebra}
\ProofTheorem{\RefLinearMap}{standard representation of map A1 A2, associative algebra}

\ePrints{4975-6381}
\ifx\Semafor\ValueOff
\ShowDefinition{component of linear map}
\fi

\ShowTheorem{conjugation transformation}
\ProofTheorem{\RefLinearMap}{conjugation transformation}

\ShowTheorem{representation of composition of linear maps}
\ProofTheorem{\RefLinearMap}{representation of composition of linear maps}

\ShowTheorem{representation of composition of linear maps A->A}
\ProofTheorem{\RefLinearMap}{representation of composition of linear maps A->A}

\fi

\ePrints{1601.03259,4975-6381,309618526,CACAA.06.121}
\ifx\Semafor\ValueOn
\ShowTheorem{coordinates of map A1 A2, algebra}
\ProofTheorem{\RefLinearMap}{coordinates of map A1 A2, algebra}

According to the definition
\ShowEq{Basis F=...}
We also consider the set
\ShowEq{F=...}
having in mind that the equality
\ShowEq{Fi=Fj}
is possible when $i\ne j$.
\fi

\ePrints{4975-6381,1506.00061,7287-9339}
\ifx\Semafor\ValueOn
\begin{convention}
In the equation
\ShowEq{n linear map A LA}
as well as in other expressions of polylinear map,
we have convention that map $f_i$ has variable $x_i$ as argument.
\qed
\end{convention}
\fi

\ePrints{1601.03259,4975-6381,1506.00061,7287-9339,309618526,CACAA.06.121}
\ifx\Semafor\ValueOn

\ShowTheorem{L(An;B) is free D module}
\ProofTheorem{\RefLinearMap}{L(An;B) is free D module}

\begin{ShadedTheorem}
\labelTheorem{polylinear map, algebra} 
Let $A$ be associative $D$\Hyph algebra.
Polylinear map
\ShowEq{polylinear map, algebra}
generated by maps
\ShowEq{I1n in L(A;A)}
has form
\ShowEq{polylinear map, algebra, canonical morphism}
where $\sigma_s$ is a transposition of set of variables
\ShowEq{transposition of set of variables, algebra}
\end{ShadedTheorem}
\ProofTheorem{\RefLinearMap}{polylinear map, algebra}

\ShowTheorem{representation of algebra An in LAnA}
\ProofTheorem{\RefLinearMap}{representation of algebra An in LAnA}

\begin{convention}
Since the tensor
\EqParm{a in Aoxn+1}{=z}
has the expansion
\ShowEq{a=oxi}
then set of permutations
\ShowEq{si in Sn}
and tensor $a$
generate the map
\ShowEq{ox:An->A}
defined by rule
\ShowEq{ox circ =}
\qed
\end{convention}
\fi
\fi

\ePrints{5284-0163,1801.01628}
\ifx\Semafor\ValueOn
\begin{ShadedTheorem}
\labelTheorem{exists A-number a fd=da}
Let
\ShowEq{f:A->B}fDA
is linear map of commutative ring $D$
into $D$\Hyph algebra $A$.
Then there exists $A$\Hyph number $a$ such that
\ShowEq{fd=da}
\end{ShadedTheorem}

\begin{proof}
The theorem follows from the equality
\ShowEq{fd=df1}
when we assume $a=f\circ 1$.
\end{proof}

\begin{corollary}
\labelCorollary{L D->A=A}
\ShowEq{L D->A=A}
\end{corollary}
\fi

\ePrints{1506.00061,7287-9339}
\ifx\Semafor\ValueOn
\Section{Algebra with Conjugation}

Let $D$ be commutative ring.
Let $A$ be $D$\hyph algebra with unit $e$, $A\ne D$.

Let there exist subalgebra $F$ of algebra $A$ such that
$F\ne A$, $D\subseteq F\subseteq Z(A)$, and algebra $A$ is a free
module over the ring $F$.
Let $\Basis e$ be the basis of free module $A$ over ring $F$.
We assume that $e_{\gi 0}=1$.

\begin{ShadedTheorem}
\labelTheorem{structural constant of algebra with unit}
Structural constants of
$D$\hyph algebra with unit $e$ satisfy condition
\ShowEq{structural constant re algebra, 1}
\end{ShadedTheorem}
\begin{proof}
\TheoremFollows
\ShowEq{ref structural constant of algebra with unit}
\end{proof}

Consider maps
\ShowEq{Re Im A->A}
defined by equation
\ShowEq{Re Im A->F 0}
\ShowEq{Re Im A->F 1}
The expression
\ShowEq{Re Im A->F 2}
is called
\AddIndex{scalar of element}{scalar of algebra} $d$.
The expression
\ShowEq{Re Im A->F 3}
is called
\AddIndex{vector of element}{vector of algebra} $d$.

According to
\EqRef{Re Im A->F 1}
\ShowEq{Re Im A->F 4}
We will use notation
\ShowEq{scalar algebra of algebra}
to denote
\AddIndex{scalar algebra of algebra}{scalar algebra of algebra} $A$.

\begin{ShadedTheorem}
\labelTheorem{Re Im A->F}
The set
\ShowEq{vector module of algebra}
is $(\re A)$\Hyph module
which is called
\AddIndex{vector module of algebra}{vector module of algebra} $A$.
\ShowEq{A=re+im}
\end{ShadedTheorem}
\begin{proof}
\TheoremFollows
\ShowEq{ref Re Im A->F}
\end{proof}

According to the theorem \RefTheorem{Re Im A->F},
there is unique defined representation
\ShowEq{d Re Im}

\begin{ShadedDefinition}
\labelDefinition{conjugation in algebra}
The map
\ShowEq{conjugation in algebra}
\ShowEq{conjugation in algebra, 0}
is called
\AddIndex{conjugation in algebra}{conjugation in algebra}
provided that this map satisfies
\ShowEq{conjugation in algebra, 1}
$(\re A)$\Hyph algebra $A$
equipped with conjugation is called
\AddIndex{algebra with conjugation}{algebra with conjugation}.
\qed
\end{ShadedDefinition}

\begin{ShadedTheorem}
\labelTheorem{conjugation in algebra}
The $(\re A)$\Hyph algebra $A$ is
algebra with conjugation
iff structural constants of
$(\re A)$\Hyph algebra $A$ satisfy condition
\ShowEq{conjugation in algebra, 5}
\ShowEq{conjugation in algebra, 4}
\end{ShadedTheorem}
\begin{proof}
\TheoremFollows
\ShowEq{ref conjugation in algebra}
\end{proof}
\fi

\ePrints{1601.03259,7287-9339,4975-6381,1506.00061}
\ifx\Semafor\ValueOn
\Section{Polynomial over Associative \texorpdfstring{$D$}{D}-Algebra}

Let $D$ be commutative ring and
$A$ be associative $D$\Hyph algebra with unit.
Let $\Basis F$ be basis of algebra
\ShowEq{L(A->B)}DAA.

\begin{ShadedTheorem}
\labelTheorem{monomial of power k}
Let $p_k(x)$ be
\AddIndex{monomial of power}{monomial of power} $k$
over $D$\Hyph algebra $A$.
Then
\StartLabelItem
\begin{enumerate}
\item
Monomial of power $0$ has form
\ShowEq{p0(x)=a0}
\item
If $k>0$, then
\ShowEq{monomial of power k}
where $a_k\in A$ and $F\in\Basis F$.
\end{enumerate}
\end{ShadedTheorem}
\begin{proof}
\TheoremFollows
\ShowEq{ref monomial of power k}
\end{proof}

In particular, monomial of power $1$ has form
\ShowEq{p1(x)=aFxa}

\begin{ShadedDefinition}
\labelDefinition{Abelian group of homogeneous polynomials over algebra}
We denote
\ShowEq{set of homogeneous polynomials}k
Abelian group
generated by the set of monomials of power $k$.
Element $p_k(x)$ of Abelian group $A_k[x]$ is called
\AddIndex{homogeneous polynomial}{homogeneous polynomial}
of power $k$.
\qed
\end{ShadedDefinition}

\begin{convention}
Let the tensor
\EqParm{a in Aoxn+1}{=.}
Let
\ShowEq{F1n in F}
When
\ShowEq{x1=xn=x}
we assume
\ShowEq{a xn=}
\qed
\end{convention}

\begin{convention}
If we have few tuples of maps $F\in\Basis F$,
then we will use index like $[k]$
to index tuple
\ShowEq{F[k]=...}
\qed
\end{convention}

\begin{ShadedTheorem}
\labelTheorem{homogeneous polynomial a circ xk}
We can present homogeneous polynomial $p(x)$
in the following form
\ShowEq{p(x)=a circ xk}
\end{ShadedTheorem}
\begin{proof}
\TheoremFollows
\ShowEq{ref map A(k+1)->pk otimes is polylinear map}
\end{proof}

\begin{ShadedDefinition}
\labelDefinition{Polynomial over algebra}
We denote
\ShowEq{algebra of polynomials over algebra}
direct sum\,\footnotemark
of $A$\Hyph modules $A_n[x]$.
An element $p(x)$ of $A$\Hyph module $A[x]$ is called
\AddIndex{polynomial}{polynomial}
over $D$\Hyph algebra $A$.
\qed
\end{ShadedDefinition}
\footnotetext{\,
See the definition of direct sum of modules in
\citeBib{Serge Lang},
page 128.
On the same page, Lang proves the existence of direct sum of modules.
}

\ePrints{4975-6381,1601.03259}
\ifx\Semafor\ValueOn
The following definition
\RefDefinition{otimes -}
is based on the definition
\xRef[1302.7204]{definition: otimes -}.

\ShowDefinition{otimes -}

\begin{ShadedDefinition}
\labelDefinition{product of homogeneous polynomials}
Product of homogeneous polynomials
\ShowEq{polynomials p,r}
is defined by the equality\,\footnotemark
\ShowEq{pr=p circ r}
\qed
\end{ShadedDefinition}
\footnotetext{\,
The definition
\RefDefinition{product of homogeneous polynomials}
is based on the definition
\xRef[1302.7204]{definition: product of homogeneous polynomials}.
}

\ePrints{}
\ifx\Semafor\ValueOn
\begin{ShadedTheorem}
\labelTheorem{p(x)=a k-1 k}
Let
\ShowEq{p(x)=a circ xk}
be monomial of power $k>1$.
Then the polynomial $p(x)$ can be represented using
one of the following forms
\ShowEq{p(x)=a k-1 k 1}
\ShowEq{p(x)=a k-1 k 2}
where
\ShowEq{p(x)=a k-1 k 3}
\end{ShadedTheorem}
\begin{proof}
See the proof of the theorem
\RefTheorem[1302.7204]{p(x)=a k-1 k}.
\end{proof}
\fi
\fi

\ePrints{1506.00061,7287-9339}%
\ifx\Semafor\ValueOn%

\begin{ShadedDefinition}
\labelDefinition{divisor of polynomial}
The polynomial $p(x)$ is called
\AddIndex{divisor of polynomial}{divisor of polynomial} $r(x)$,
if we can represent the polynomial $r(x)$ as
\ShowEq{r=qpq(x)}
\qed
\end{ShadedDefinition}

\ShowTheorem{r=+q circ p}
\begin{proof}
\TheoremFollows
\ShowEq{ref r=+q circ p}
\end{proof}

\begin{ShadedTheorem}
\labelTheorem{r=+q circ p 1}
Let
\ShowEq{p=p+p circ x}
be polynomial of power $1$ and $p_1$ be nonsingular tensor.
Let
\ShowEq{r power k}
be polynomial of power $k$.
Then
\ShowEq{r=+q circ p 1}
\end{ShadedTheorem}
\begin{proof}
\TheoremFollows
\ShowEq{ref r=+q circ p 1}
\end{proof}
\fi
\fi

%% file: Preliminary.D.Algebra.Eq.tex

\ShowEq{def DAlgebra}

\DefText[6]{Preliminary 1 homomorphism of D algebra}
{
\ShowDefinition{homomorphism from A1 to A2, D algebra}{#1}{#2}{#3}{#4}

\ShowTheorem{homomorphism from A1 to A2, D algebra}{#1}{#2}{#3}{#4}{#5}
\proofTheorem{\RefLinearMap}{homomorphism from A1 to A2, D algebra}{#1#2}

\ShowEq{\DefCol}
\ShowFootnote{homomorphism of d algebra}{#1}{#2}{#3}{#4}
\ShowEq{\DefRow}
\ShowFootnote{homomorphism of d algebra}{#1}{#2}{#3}{#4}

\TwoColText
{
\ShowEq{\DefCol}
\ShowTheorem{homomorphism of d algebra, coordinates}{#1}{#2}{#3}{#4}{#6}
\proofTheorem{\RefLinearMap}{homomorphism of d algebra, coordinates}{#1#2\Cols}
}
{
\ShowEq{\DefRow}
\ShowTheorem{homomorphism of d algebra, coordinates}{#1}{#2}{#3}{#4}{#6}
\proofTheorem{\RefLinearMap}{homomorphism of d algebra, coordinates}{#1#2\Cols}
}

\TwoColText
{
\ShowEq{\DefCol}
\ShowTheorem{matrix generates D algebra homomorphism}{#1}{#2}{#3}{#4}
\proofTheorem{\RefLinearMap}{matrix generates D algebra homomorphism}{\Cols(#1#2)}
}
{
\ShowEq{\DefRow}
\ShowTheorem{matrix generates D algebra homomorphism}{#1}{#2}{#3}{#4}
\proofTheorem{\RefLinearMap}{matrix generates D algebra homomorphism}{\Cols(#1#2)}
}
}

\DefText{Preliminary homomorphism D algebra}
{
\ePrints{5284-0163,1801.01628,MAlgebra2,5148-4632}%
\ifx\Semafor\ValueOn
\ePrints{5284-0163,1801.01628}%
\ifx\Semafor\ValueOff
\ShowText{Preliminary 1 homomorphism of D algebra}1122{h(p)}{h(a)}
\fi

\ShowText{Preliminary 1 homomorphism of D algebra}{}1{}2pa
\fi
}

%% file: A.Module.English.tex
\ifx\UseAModule\undefined
\input{A.Module.Eq}
\fi

\chapter{\SideWSC Module over \texorpdfstring{$D$}{D}-Algebra}

\section{\SideWSC Module over Associative \texorpdfstring{$D$}{D}-Algebra}

Let $D$ be commutative associative ring with unit.

\ShowText{definition of A module}

\ShowText{linear combination of vectors}

\ShowText{basis of module}

\section{\texorpdfstring{$A$}{A}-Module Type}

\subsection{\SideWSC \texorpdfstring{$A$}{A}-Module of Columns}

\ShowEq{def AModule}
\ShowEq{\DefCol}
\ShowText{Vector Space Type}

\subsection{\SideWSC \texorpdfstring{$A$}{A}-Module of Rows}

\ShowEq{\DefRow}
\ShowText{Vector Space Type}

\section{Submodule of \SideWSC \texorpdfstring{$A$}{A}-Module}

\ShowDefinition{submodule}

\ShowTheorem{submodule}
\ShowProof{submodule}

%% file: A.Module.Eq.tex
\def\UseAModule{}

\AddEquation{division algebra, basis RCA module, 1}
{
e_1c^1+...+e_nc^n=0
}

\AddEquation{division algebra, basis RCA module, 2}
{
e_1=-e_2c^2(c^1)^{-1}-...-e_nc^n(c^1)^{-1}
}

\DefEq
{
\symb{f+g}{sum of maps}{Acr module}
}
{sum of maps,,Acr module}

\DefEquation
{
x\CRcirc(\ShowSymbol{sum of maps}{Acr module})=x\CRcirc f+x\CRcirc g
}
{sum of maps, Acr module}

\DefEq
{
\symb{f+g}{sum of maps}{A module}
}
{sum of maps,,A module}

\DefEquation
{
x\CRcirc(\ShowSymbol{sum of maps}{A module})=x\CRcirc f+x\CRcirc g
}
{sum of maps, A module}

\DefEquation
{
\begin{matrix}
f+g:V\rightarrow W
&f,g\in\LAVW
\end{matrix}
}
{sum of maps, 1, Acr module}

\DefEq
{
\[
\Vector f:\Vector V\rightarrow \Vector W
\]
}
{f:V->W}

\DefEq
{
\[
b=a\CRcirc f
\]
}
{Acr linear map of module, presentation}

\DefEq
{
\[\Vector v\CRcirc\Vector f=v\CRcirc f\CRcirc \Vector e_W\]
}
{Acr linear map 1}

\DefEq
{
\begin{align*}
(a\CRcirc \Vector v)\CRcirc\Vector f&=a\CRcirc v\CRcirc f\CRcirc \Vector e_W \\
&=a\CRcirc (\Vector v\CRcirc\Vector f)
\end{align*}
}
{Acr linear map 2}

\DefEq
{
\begin{align*}
(a\CRstar x)\CRstar f&=a\CRstar (x\CRstar f)
\\
(a\CRstar x)\CRstar g&=a\CRstar (x\CRstar g)
\end{align*}
}
{sum of maps, 2, Acr module}

\DefEq
{
$f=(f^k_i)$, \iI, $k\in K$,
}
{matrix of Acr linear map}

\DefEq
{
\[
(fm)^k_i=f^k_im
\]
}
{matrix fm of Acr linear map}

\DefEq
{
\[
\Vector fm:\Vector V\rightarrow \Vector W
\]
}
{right action on Acr linear map of module}

\DefEq
{
\[
f(m)=\Vector{\delta}m
\]
}
{starT structure of Acr module, 2}

\DefEquation
{
\Vector v\CRcirc[fd]=d(\Vector v\CRcirc f)=(d\Vector v)\CRcirc f
}
{product of Acr map over scalar from ring}

\DefEq
{
\[
\Vector f:\Vector V\rightarrow\Vector W
\]
}
{product of Acr map over scalar from ring, 2}

\DefEquation
{
\Vector v=v\RCcirc\Vector e_V
}
{product of Acr map over scalar from ring, 1}

\DefEq
{
\[v\RCcirc A\RCcirc m\]
}
{Acr linear maps of module, 2}

\DefEq
{
\[\Vector w=\Vector v\RCcirc [\Vector f\RCcirc m]\]
}
{Acr linear maps of module, 3}

\DefEq
{
\[w=v\RCcirc f\RCcirc g\]
}
{Acr linear maps of module, 4}

\DefEquation
{
\Vector v\RCcirc[\Vector fd]=v^i(f_i^kd)\Vector e_{W\cdot k}
=d(v^if_i^k\Vector e_{W\cdot k})
}
{product of Acr map over scalar from ring, 3}

\DefEq
{
\begin{align*}
(a\CRstar x)\CRstar(f+g) 
&=(a\CRstar x)\CRstar f +(a\CRstar x)\CRstar g\\
&=a\CRstar (x\CRstar f)+a\CRstar (x\CRstar g)\\
&=a\CRstar (x\CRstar f+x\CRstar g)\\
&=a\CRstar (x\CRstar (f+g))
\end{align*}
}
{sum of maps, 3, Acr module}

\DefEq
{
\begin{align*}
x\CRcirc(f+g)
&=x\CRcirc f+x\CRcirc g\\
&=x\CRcirc g+x\CRcirc f\\
&=x\CRcirc(g+f)
\end{align*}
}
{sum of maps, 4, Acr module}

\DefEq
{
\[
0=(v'-v)\CRstar e
\]
}
{expansion relative basis, Astar module over algebra, 4}

\DefEquation
{
b=f\CRcirc a
}
{linear map of A module, presentation}

\DefEq
{
\[
(\Vector g\circ \Vector f)\circ\Vector v
=\Vector g\circ(\Vector f\circ\Vector v)
=(\Vector v\RCcirc\Vector f)\RCcirc\Vector g
=\Vector v\RCcirc(\Vector f\RCcirc\Vector g)
\]
}
{Acr map of module, product 2}

\DefEq
{
\[
\Vector g\circ\Vector f=\Vector f\CRcirc\Vector g
\]
}
{Acr map of module, product 1}

\DefEquation
{
f\circ v=v\RCcirc f
}
{presentation of linear map, A module}

\DefEq
{
\[
\xymatrix{
V\ar[rr]^h\ar[dr]^f & & W\\
&U\ar[ur]^g &
}
\]
}
{product of linear maps}

\DefEquation
{
b=a\CRcirc f
}
{product of linear maps, Acr, f}

\DefEquation
{
b=a\circ f
}
{product of linear maps, A, f}

\DefEquation
{
c=b\CRcirc g
}
{product of linear maps, Acr, g}

\DefEquation
{
c=b\CRcirc g
}
{product of linear maps, A, g}

\DefEq
{
\symb{\mathcal L(A;V;W)}{set linear maps}1
}
{set A linear maps, module}

\DefEq
{
\begin{align*}
h\circ(a\CRcirc e_V)
&=(g\CRcirc f)\circ(a\CRcirc e_V)
=g\circ (f\circ(a\CRcirc e_V))
\\
&=g\circ ((a\CRcirc f)\CRcirc\Vector e_U)
=((a\CRcirc f)\CRcirc g)\CRcirc e_W
\\
&=(a\CRcirc (f\CRcirc g))\CRcirc e_W
=(a\CRcirc h)\CRcirc e_W
\end{align*}
}
{product of linear maps, 1, Acr module}

\DefEquation
{
\begin{matrix}
c=a\CRcirc h= a\CRcirc f\CRcirc g
\\
h=f\CRcirc g
\end{matrix}
}
{product of linear maps, Acr, h}

\DefEquation
{
\begin{matrix}
c=a\circ h= a\circ f\circ g
\\
h=f\circ g
\end{matrix}
}
{product of linear maps, A, h}

\DefEquation
{
b\Vector v+c\CRstar e =0
}
{expansion relative basis, Astar module over algebra, 1}

\DefEquation
{
\Vector v=(-b^{-1}c)\CRstar e
}
{expansion relative basis, Astar module over algebra, 2}

\DefEquation
{
\Vector v=v'\CRstar e 
}
{expansion relative basis, Astar module over algebra, 3}

\DefEq
{
$\Vector v$, $e_i$
}
{expansion relative basis, module over algebra, 0}

\DefEquation
{
\Vector v=v\CRstar e
}
{expansion relative basis, Astar module over algebra}

\DefEq
{
$e_1+e_2$.
}
{Vector A subs row, 3}

\DefEq
{
\[
\begin{matrix}
f_{1,3\cdot 3}:d\circ \Vector v=(d\circ e)\circ\Vector v
&d\in D&\Vector v\in\Vector V
\end{matrix}
\]
}
{representation D in V}

\DefEquation
{
\Vector v=v_1\CRstar\Vector e_{3\cdot 1}=
v_1^k\ \Vector e_{3\cdot 1\cdot k}
}
{expansion vector v, basis V over A 1}

\DefEq
{
\[
\begin{matrix}
f'_{1,3\cdot 2}:D\rightarrow\Vector V
&f'_{1,3\cdot 2}(d)\circ\Vector v=f_{2,3\cdot 2}(f_{1,2}(d)\circ e)\circ\Vector v
\end{matrix}
\]
}
{sT representation D A V 2}

%% file: Homomorphism.A.Module.English.tex
\ifx\HomomorphismAModule\undefined
\input{Homomorphism.A.Module.Eq}\fi

\chapter{Homomorphism of \SideWSC \texorpdfstring{$A$}{A}-\VectorSetNSC}

\ShowText{homomorphism A module}

\section{Sum of Homomorphisms of \texorpdfstring{$A$}{A}-\VectorSetNSC}

\ShowEq{sum of homomorphisms, A vector space}

\section{Product of Homomorphisms of \texorpdfstring{$A$}{A}-\VectorSetNSC}

\ShowEq{product of homomorphisms, A vector space}

%% file: Homomorphism.A.Module.Eq.tex
\def\HomomorphismAModule{}

\AddEq{sum of homomorphisms, A vector space}
{
\ShowTheorem{sum of homomorphisms is homomorphism, A vector space}

\ShowProof{sum of homomorphisms is homomorphism, A vector space}

\TwoColText
{
\ShowEq{\DefCol}
\ShowTheorem{sum of homomorphisms, A vector space}
\begin{sloppypar}
\ShowProof{sum of homomorphisms, A vector space}
\end{sloppypar}
}
{
\ShowEq{\DefRow}
\ShowTheorem{sum of homomorphisms, A vector space}
\begin{sloppypar}
\ShowProof{sum of homomorphisms, A vector space}
\end{sloppypar}
}
}

\AddEq{product of homomorphisms, A vector space}
{
\ShowTheorem{product of homomorphisms is homomorphism, A vector space}
\ShowProof{product of homomorphisms is homomorphism, A vector space}

\TwoColText
{
\ShowEq{\DefCol}
\ShowTheorem{product of homomorphisms, A vector space}gUVhVW
}
{
\ShowEq{\DefRow}
\ShowTheorem{product of homomorphisms, A vector space}hVWgUV
}

\TwoColText
{
\ShowEq{\DefCol}
\begin{sloppypar}
\ShowProof{product of homomorphisms, A vector space}gh
\end{sloppypar}
}
{
\ShowEq{\DefRow}
\begin{sloppypar}
\ShowProof{product of homomorphisms, A vector space}hg
\end{sloppypar}
}
}

%% file: Na.A.Module.English.tex
\ifx\UseAModuleNa\undefined
\input{Na.A.Module.Eq}
\fi

\chapter{\SideWSC Module over \texorpdfstring{$D$}{D}-Algebra}

\section{\SideWSC Module over Non-Associative \texorpdfstring{$D$}{D}-Algebra}

Let $D$ be commutative associative ring with unit.

\ShowText{definition of A module Na}

%% file: Na.A.Module.Eq.tex
\def\UseAModuleNa{}

\DefText{definition of A module Na}
{
\ShowDefinition{module over non-associative algebra}

\ProveTheorem{module over non-associative algebra}

\def\TheoremStyle{ParacolTheorem}%
\ShowTheorem{definition of non-associative A module}
\def\TheoremStyle{ShadedTheorem}%
\ShowTheorem{definition of non-associative A module, property}
}